\documentclass[reqno,english]{amsbook}

\input epsf

\usepackage{graphicx}
\usepackage{amsmath}
\usepackage{amscd}
\usepackage{amssymb}
\usepackage{amstext}
\usepackage{amsmath}
\usepackage[all]{xy}
\usepackage{color}

\def\E{\ifmmode{\mathbb E}\else{$\mathbb E$}\fi} 
\def\N{\ifmmode{\mathbb N}\else{$\mathbb N$}\fi} 
\def\R{\ifmmode{\mathbb R}\else{$\mathbb R$}\fi} 
\def\Q{\ifmmode{\mathbb Q}\else{$\mathbb Q$}\fi} 
\def\C{\ifmmode{\mathbb C}\else{$\mathbb C$}\fi} 
\def\H{\ifmmode{\mathbb H}\else{$\mathbb H$}\fi} 
\def\Z{\ifmmode{\mathbb Z}\else{$\mathbb Z$}\fi} 
\def\P{\ifmmode{\mathbb P}\else{$\mathbb P$}\fi} 
\def\T{\ifmmode{\mathbb T}\else{$\mathbb T$}\fi} 
\def\SS{\ifmmode{\mathbb S}\else{$\mathbb S$}\fi} 
\def\DD{\ifmmode{\mathbb D}\else{$\mathbb D$}\fi} 
\def\K{\ifmmode{\mathbb K}\else{$\mathbb K$}\fi}

\newcommand{\del}{\partial}

\theoremstyle{plain} 
\newtheorem{thm}{Theorem}[section]
\newtheorem{cor}[thm]{Corollary}
\newtheorem{lem}[thm]{Lemma}
\newtheorem{sublem}[thm]{Sublemma}
\newtheorem{prop}[thm]{Proposition}

\newtheorem{conj}[thm]{Conjecture}
\newtheorem{clm}[thm]{Claim}

\theoremstyle{definition}
\newtheorem{defn}[thm]{Definition}
\newtheorem{conv}[thm]{Convention}
\newtheorem{rem}[thm]{Remark}

\newtheorem{exm}[thm]{Example}
\newtheorem{conds}[thm]{Condition}
\newtheorem{notation}[thm]{\rm\bfseries{Notation}}
\newtheorem{warn}[thm]{Warning}
\newtheorem{assum}[thm]{Assumption}

\newtheorem*{thm*}{Theorem}

\numberwithin{equation}{section}

\def\R{{\mathbb R}}

\def\E{{\mathbb E}}
\def\Z{{\mathbb Z}}
\def\C{{\mathbb C}}
\def\R{{\mathbb R}}
\def\P{{\mathbb P}}

\def\N{{\mathbb N}}

\def\11{{\mathbb I}}

\def\CM{{\mathcal M}}
\def\CX{{\mathcal X}}

\def\e{\varepsilon} 
\def\CA{{\mathcal A}}
\def\opname#1{\mathop{\kern0pt{\rm #1}}\nolimits}

\def\dim{\opname{dim}}

\def\Aut{\operatorname{Aut}}

\makeindex

\begin{document}
\quad \vskip1.375truein

\title[Floer theory and Mirror symmetry on toric manifolds]
{Lagrangian Floer theory and Mirror symmetry on compact toric manifolds}

\author{Kenji Fukaya}\address{Simons Center for Geometry and Physics,
State University of New York, Stony Brook, NY 11794-3636 U.S.A., \& Center for Geometry and Physics, Institute for Basic Sciences (IBS), Pohang, Korea} \email{kfukaya@scgp.stonybrook.edu}
\author{Yong-Geun Oh}
\address{Center for Geometry and Physics, Institute for Basic Sciences (IBS), Pohang, Korea \& Department of Mathematics,
POSTECH, Pohang, Korea} \email{yongoh1@postech.ac.kr}
\author{Hiroshi Ohta} 
\address{Graduate School of Mathematics,
Nagoya University, Nagoya, Japan} \email{ohta@math.nagoya-u.ac.jp}
\author{Kaoru Ono}
\address{Research Institute for Mathematical Sciences,
Kyoto University, Kyoto, Japan}
\email{ono@kurims.kyoto-u.ac.jp}

\thanks{KF is supported partially by JSPS Grant-in-Aid for Scientific Research \#
19104001, 23224002, Global COE program G08 and NSF Grant No. 1406423, YO by US NSF grant \# 0904197  
and the IBS project IBS-R003-D1, HO by JSPS Grant-in-Aid
for Scientific Research \# 19340017, 23340015, 15H02054 and KO by JSPS Grant-in-Aid for
Scientific Research, \# 18340014, 21244002, 23224001, 26247006.}

\begin{abstract}
In this paper we study
Lagrangian Floer theory on toric manifolds from the point of view of
mirror symmetry.
We
construct a natural isomorphism between the Frobenius manifold
structures of the (big) quantum cohomology of the toric manifold and
of Saito's theory of singularities of the potential function constructed
in \cite{fooo09} via the Floer cohomology deformed by ambient
cycles. Our proof of the isomorphism involves the open-closed
Gromov-Witten theory of one-loop.
\end{abstract}

\date{Mar. 31, 2013}


\maketitle

\tableofcontents
\chapter{Introduction}
\section{Introduction}
\label{sec:introduction}
The purpose of this paper is to prove a version of
mirror symmetry between compact toric A-model and Landau-Ginzburg B-model.
\par
Let $X$ be a compact toric manifold and take 
a toric K\"ahler structure on it.
In this paper $X$ is not necessarily assumed to be Fano.
In \cite{fooo09} we defined a potential function with bulk,
$\mathfrak{PO}_{\frak b}$, which is a family of
functions of $n$ variables $y_1,\dots,y_n$
parameterized by the cohomology class $\frak b \in H(X;\Lambda_0)$.
(More precisely, it is parameterized by $T^n$-invariant cycles. We explain
this point later in Sections \ref{sec:statements} and
\ref{sec:uptovariables}.)
($\Lambda_0$ is defined in Definition \ref{novring}.)
\par
$\mathfrak{PO}_{\frak b}$ is an element of
an appropriate completion of
$\Lambda_0[y_1,y_1^{-1}\dots,y_n,y_n^{-1}]$,
the Laurent polynomial ring over
(universal) Novikov ring $\Lambda_0$.
We denote this completion by $\Lambda\langle\!\langle y,y^{-1}\rangle\!\rangle_0^{\overset \circ{P}}$.
(See Definition \ref{def:compPo} for its definition.)
We put
$$
\text{\rm Jac}(\frak{PO}_{\frak b})
= \frac{ \Lambda\langle\!\langle y,y^{-1}\rangle\!\rangle_0^{\overset{\circ}P} }
{\text{\rm Clos}_{d_{\overset{\circ}P}}\left(y_i \frac{\partial\frak{PO}_{\frak b}}{\partial y_i}
: i=1,\dots,n\right)}.
$$
Differentiation of $\frak{PO}_{\frak b}$ with respect to the ambient cohomology class
$\frak b$ gives rise to a $\Lambda_0$ module homomorphism
\begin{equation}\label{KSdefikntro}
\frak{ks}_{\frak b}:
H(X;\Lambda_0) \to \text{\rm Jac}(\frak{PO}_{\frak b}).
\end{equation}
We define a quantum ring structure $\cup^{\frak b}$ on $H(X;\Lambda_0)$
by deforming the cup product by $\frak b$ using Gromov-Witten theory.
(See Definition \ref{def:deformcup}.)
The main result of this paper can be stated as follows: 
\begin{thm}\label{mtintro}
Equip $H(X;\Lambda_0)$ with the ring structure $\cup^{\frak b}$. Then
\begin{enumerate}
\item The homomorphism $\frak{ks}_{\frak b}$ in
$(\ref{KSdefikntro})$ is a ring isomorphism.
\item If $\frak{PO}_{\frak b}$ is a Morse function, i.e., all its critical
points are nondegenerate, the isomorphism $\frak{ks}_{\frak b}$ above
intertwines the Poincar\'e duality pairing on $H(X;\Lambda)$ and the residue
pairing on $\text{\rm Jac}(\frak{PO}_{\frak b})\otimes_{\Lambda_0} \Lambda$.
\end{enumerate}
\end{thm}
We refer readers to Section \ref{sec:statements} for the definition of various
notions appearing in Theorem \ref{mtintro}.
Especially the definition of residue pairing which we use in this paper is
given in Definition \ref{res2} and Theorem \ref{cliffordZ}.
\par
We now explain how Theorem \ref{mtintro} can be regarded as a mirror symmetry between
toric A-model and Landau-Ginzburg B-model:
\par
\begin{enumerate}
\item[(a)]
First, the surjectivity of (\ref{KSdefikntro})
implies that the family $\frak{PO}_{\frak b}$
is a versal family in the sense of deformation
theory of singularities. (See Theorem \ref{versality}.)
\item[(b)]
Therefore the tangent space $T_{\frak b}H(X;\Lambda_0)$
carries a ring structure pulled-back from the Jacobian ring by (\ref{KSdefikntro}).
Theorem \ref{mtintro}.1 implies that this pull-back coincides with the quantum
ring structure $\cup^{\frak b}$.
\item[(c)]
Theorem \ref{mtintro}.2 implies that the residue paring, which is defined when
$\text{\rm Jac}(\frak{PO}_{\frak b})$ is a Morse function, can
be extended to arbitrary class $\frak b \in H(X;\Lambda_0)$.
In fact, the Poincar\'e duality pairing of
$H(X;\Lambda_0)$ is independent of $\frak b$ and so obviously extended.
\par
This (extended) residue pairing\index{residue pairing} is nondegenerate and
defines a `Riemannian metric' on $H(X;\Lambda_0)$.
(Here we put parenthesis since our metric is
over the field of fractions $\Lambda$ of Novikov ring and is not over $\R$.)
As in the standard Riemannian geometry, it
determines the Levi-Civita connection $\nabla$ on $H(X;\Lambda_0)$.
Theorem \ref{mtintro}.2 then implies that this connection is flat: In fact,
in the A-model side this connection is noting but the canonical affine
connection on the affine space $H(X;\Lambda_0)$.
\par
Let $w_i$ be affine coordinates of $H(X;\Lambda_0)$ i.e.
the coordinates satisfying $\nabla_{\partial/\partial w_i}\partial/\partial w_j = 0$.
Then we have a function $\Phi$ on $H(X;\Lambda_0)$
that satisfies
\begin{equation}\label{GWpotential}
\langle \text{\bf f}_i \cup^{\frak b} \text{\bf f}_j,
\text{\bf f}_k\rangle_{\text{\rm PD}_X}
= \frac{\partial^3\Phi}{\partial w_i\partial w_j\partial w_k}.
\end{equation}
Here $\{\text{\bf f}_i\}$ is the basis of $H(X;\Lambda_0)$
corresponding to the affine coordinates $w_i$ and
$\langle\cdot,\cdot\rangle_{\text{\rm PD}_X}$ is the
Poincar\'e duality pairing.
\par
In fact, $\Phi$ is constructed from Gromov-Witten invariants and
is called Gromov-Witten potential. (See \cite{Manin:qhm}, for example.)
Using the isomorphism given in Theorem \ref{mtintro}.1
we find that the third derivative of $\Phi$ also
gives the structure constants of the Jacobian ring.
\item[(d)]
In  \cite[Section 10]{fooo09}, we defined an Euler vector field $\frak E$ on $H(X;\Lambda_0)$.
We also remark that the Jacobian ring carries a unit that is parallel with
respect to the Levi-Civita connection: This is obvious in the A-model side
and hence the same holds in the B-model side.
\item[(e)]
The discussion above implies that the miniversal family $\frak{PO}_{\frak b}$ determines
the structure of Frobenius manifold on $H(X;\Lambda_0)$.\index{Frobenius manifold}
K. Saito \cite{Sai83} and M. Saito \cite{Msaito} defined a Frobenius manifold structure
on the parameter space of miniversal deformation of a holomorphic function
germ with isolated singularities. Their proof is based on a very deep notion of
primitive forms. On the other hand, our construction
described above uses mirror symmetry which is different from theirs.
Strictly speaking, our situation is also slightly different from Saito's
in that we work over the Novikov ring but not over $\C$, and deals with
a more global version than that of function germs.
See \cite{saba} (especially its Section 4) for some discussion
of a global version of Saito theory\index{Saito theory} over $\C$. It is not clear to the authors
how difficult it is to adapt Saito's argument to our situation and prove
existence of a Frobenius manifold structure using the primitive forms.
\item[(f)]
Dubrovin \cite{dub} showed that quantum cohomology also defines a
Frobenius manifold structure on $H(X;\Lambda_0)$.
Theorem \ref{mtintro} can be regarded as the coincidence of the two
Frobenius manifolds, one defined by quantum cohomology
(= toric $A$-model) and the other defined by Saito's theory
(= Landau-Ginzburg $B$-model).
\end{enumerate}

\begin{rem}\label{bunken}
\begin{enumerate}
\item Saito invented a notion, which he calls \emph{flat structure}. It is
the same structure as
Frobenius manifold structure which Dubrovin found in
Gromov-Witten theory.
\item
For the case of Fano $X$ and with $\frak b=0$, a statement similar
to Theorem \ref{mtintro}.1 was known to
Batyrev \cite{batyrev:qcrtm92} and Givental \cite{givental2}.
(\cite{givental2} contains its proof relying on the virtual
localization formula which was later established in \cite{grapan}.)
For the case of $\C P^n$, a result closely related to Theorem
\ref{mtintro} was known to Takahashi \cite{taka} ($n=1$) and
Baranikov \cite{barani}. Proofs of both Takahashi and Baranikov are
rather computational. On the other hand, our proof does not rely on
computations but uses a natural construction of the map
$\frak{ks}_{\frak b}$ which we prove is an isomorphism. We prove
this isomorphism property using some geometric argument without
carrying out separate calculations on both sides. In addition to
this geometric argument, we use some important partial calculation
of $\frak{PO}_{\frak b}$ carried out by Cho-Oh \cite{cho-oh} in our
proof of Theorem \ref{mtintro}.
\item
In the work \cite{hori-vafa} of Hori-Vafa a function similar to the
potential function $\frak{PO}_0$ appears\footnote{Hori and Vafa
conjectured that the superpotential in their situation coincides
with the potential in Lagrangian Floer theory. A similar conjecture
was also made by Givental around the same time. This conjecture was
one of the motivations of the work \cite{cho-oh}.} and their
discussion is somewhat close to Theorem  \ref{mtintro}. Also there
are various conjectures closely related to Theorem \ref{mtintro} in
the literature such as those stated by Saito-Takahashi
\cite{SaiTaka}, Iritani
\cite{iritani3} and others.
These conjectures are more wide-ranging
than Theorem \ref{mtintro} in their scope in that they also discuss the Hodge
structure, higher residues and others. Especially in \cite{iritani3}
a similar statement as Theorem \ref{mtintro} is announced as a
theorem, whose proof is deferred to subsequent papers.
On the other hand, in \cite{iritani3} and other references,
superpotential in the B-model side are defined `by hand' and the
role of potential function {\it with bulk} $\frak{PO}_{\frak b}$ was
not discussed. In that sense, Theorem \ref{mtintro} does not seem to
have appeared in the previous literature even as a conjecture.
See also \cite{KKP}.
\item
We also remark that our potential function $\frak{PO}_{\frak b}$ in
Lagrangian Floer theory is slightly different from the
superpotential appearing in the literature we mentioned above: The
`superpotential' in the literature coincides with the leading order
potential function $\frak{PO}_0$ in our terminology. $\frak{PO}_0$
coincides with potential function $\frak{PO}_{\frak b}$ for the case
of $X$ Fano and $\frak b = \text{\bf 0}$. There is an example
where $\frak{PO}_0$ is different from $\frak{PO}_{\frak b}$, $\frak b= \text{\bf 0}$
in non-Fano case. (Theorem \ref{F2them}.) However we can prove that
$\frak{PO}_0$ is the same as
$\frak{PO}_{\frak b}$ for some $\frak b$ up to an appropriate
coordinate change (Theorem \ref{versality}).
\item
We believe that using Theorem \ref{mtintro} we can prove a version
of homological mirror symmetry conjecture (\cite{konts:hms}) between
toric A-model and Landau-Ginzburg B-model, namely an equivalence between the 
Fukaya category of $X$ and the category of matrix factorization of
$\text{\rm Jac}(\frak{PO}_{\frak b})$. We plan to explore this point elsewhere
in the near future, with M. Abouzaid.
\item
It seems important to push our approach further and study the mirror
symmetry along the lines of Saito, Givental and Hori-Vafa. For
example, including higher residues, primitive forms, quantum
D-module and hypersurfaces of toric manifolds into our story are
some interesting projects. We hope to pursue this line of research
in the future.
\item
After the main results of this paper were announced,
an interesting result of Gross \cite{gros} appears which studies a similar problem
in the case of $\C P^2$. See also \cite{grospand}.
\end{enumerate}
\end{rem}
\par
We can apply Theorem \ref{mtintro} to the study of Lagrangian Floer
theory of toric manifolds. Let $\frak M(X,\frak b)$ be the set of
pairs $(u,b)$ with
$$
u \in \text{\rm Int} P, \quad b \in H^1(L(u);\Lambda)/
H^1(L(u);2\pi\sqrt{-1}\Z)
$$
such that
$$
HF((L(u),\frak b,b),(L(u),\frak b,b);\Lambda) \ne 0.
$$
Here $HF((L(u),\frak b,b),(L(u),\frak b,b);\Lambda)$
is the Floer cohomology with bulk deformation.
($\Lambda$ is defined in Definition \ref{novring}
(See \cite{fooo06}, \cite{fooo09}.)
\begin{thm}\label{number}
If $\frak{PO}_{\frak b}$ is a Morse function, we have
\begin{equation}\label{numberformula}
\# \frak M(X,\frak b) = \sum_k \text{\rm rank}H_k(X;\Q).
\end{equation}
If $\frak{PO}_{\frak b}$ is not a Morse function, we have
\begin{equation}\label{number2}
0 < \# \frak M(X,\frak b) < \sum_k \text{\rm rank}H_k(X;\Q).
\end{equation}
\end{thm}
\par
Theorem \ref{number} settles a conjecture in the preprint version of \cite{fooo08}. It also
removes the rationality hypothesis we put in the statement of
 \cite[Theorem 1.4]{fooo08}.
Theorem \ref{number} is proved in Section \ref{sec:statements}
right after Theorem \ref{Mirmain}, using Theorem \ref{mtintro}.1.
\par
We have another application of Theorem  \ref{mtintro} which is
proved in Section \ref{sec:c1}.
\begin{thm}\label{c1iscrit}
Assume $\frak b \in H^2(X;\Lambda_0)$.
The set of eigenvalues of the map $x \mapsto c_1(X) \cup^{\frak b} x :
H(X;\Lambda_0) \to  H(X;\Lambda_0)$ coincides with the set of
critical values of $\frak{PO}^{\frak b}$, with multiplicities counted.
\end{thm}
Here $\cup^{\frak b}$ is the quantum cup product associated with the
bulk deformation by $\frak b$. (See Definition \ref{def:deformcup} for its definition.)
\begin{rem}
Theorem \ref{c1iscrit} was conjectured by M. Kontsevich.
See also \cite{Aur07}. Theorem \ref{c1iscrit} is proved in \cite{fooo08}
in case $X$ is Fano and $\frak b=\text{\bf 0}$.
\end{rem}
\par
Outlines of the proof of Theorem \ref{mtintro} and of the contents
of the following sections are in order. The proof of Theorem
\ref{mtintro}.1 goes along the line we indicated in 
\cite[Remark 6.15]{fooo08}. However the details of the proof are much more
involved partially because we have much amplified the statement.
\par
In Section \ref{sec:statements} we define various notions used in
the statement of Theorem \ref{mtintro} and make a precise statement
thereof. In Section \ref{sec:example}, we make explicit computations
for the case of $X = \C P^n$, $\frak b = \text{\bf 0}$ to illustrate Theorem
\ref{mtintro}. Section \ref{sec:valuation} contains several
(technical) issues related to the completion of Laurent polynomial
ring over Novikov ring.
\par
As we already mentioned, the potential function we defined in
\cite{fooo09} is parametrized by the group of $T^n$-invariant
cycles. The set of $T^n$-invariant cycles, which we denote by
$\mathcal A(\Lambda_0)$, is mapped surjectively onto the cohomology
group $H(X;\Lambda_0)$. However this map is not injective. In
Section \ref{sec:uptovariables}, we prove that if $\frak b$ and
$\frak b'$ give the same cohomology class in $H(X;\Lambda_0)$,
$\frak{PO}_{\frak b}$ and $\frak{PO}_{\frak b'}$ coincide up to an
appropriate change of coordinates. (See Theorem \ref{cchangethem}.)
In Section \ref{sec:welldef} we prove that the map $\frak{ks}_{\frak
b}$, which is originally defined on $\mathcal A(\Lambda_0)$,
descends to a map defined on $H(X;\Lambda_0)$. (See Theorem
\ref{indepenceKS}.)
\par
In Section \ref{sec:ringhomo}, we prove that the map
$\frak{ks}_{\frak b}$ is a ring homomorphism. In Section
\ref{sec:QMHOch} we explain its relation to various results, ideas,
and conjectures appearing in the previous literature. We also
explain there how the discussion of Section \ref{sec:ringhomo} can
be generalized beyond the case of toric manifolds.
\par
In Section \ref{sec:surf}, surjectivity of $\frak{ks}_{\frak b}$ is
proved. We first show that its $\C$-reduction
$\overline{\frak{ks}}_{\frak b}$ is surjective and then use it to
prove surjectivity of $\frak{ks}_{\frak b}$ itself.
\par
Sections \ref{sec:versality} - \ref{sec:injgen} are devoted to the
proof of injectivity of $\frak{ks}_{\frak b}$. Since this proof is
lengthy and nontrivial, we give some outline of this proof
here.
\par
We first remark that $\frak{PO}_{\frak b}$ is not a Laurent
polynomial in general but a formal Laurent power series of
$y_1,\dots,y_n$ with elements of $\Lambda_0$ as its coefficients,
which converges in an appropriate $T$-adic topology. We note that
even in the Fano case we need to work with such formal power series
as soon as we turn on bulk deformations. Therefore we first show
that $\frak{PO}_{\frak b}$ can be transformed into a Laurent {\it
polynomial} by an appropriate change of variables. This is proved in
Section \ref{sec:algebra}. This is possible because of the
\emph{versality} of our family which we prove in Section
\ref{sec:versality} using the surjectivity of
$\overline{\frak{ks}}_{\frak b}$.
\par
The second point to take care of is the fact that we actually need
to take a closure of the ideal generated by
$y_i\frac{\partial\frak{PO}_{\frak b}} {\partial y_i}$. This is
especially important when we work over a Laurent polynomial ring and
not over its completion. For example, the equality
(\ref{numberformula}) will not hold in the non Fano case, if we do
not take the closure of the ideal. (See  \cite[Example 8.2]{fooo08}.)
This point is discussed in Section \ref{sec:injgen}.
\par
After these preparations, we proceed with the proof of injectivity
of $\frak{ks}_{\frak b}$ as follows. Using the surjectivity of
$\frak{ks}_{\frak b}$ and the fact that $H(X;\Lambda_0)$ is a free
$\Lambda_0$ module, we observe that the injectivity of
$\frak{ks}_{\frak b}$ will follow from the inequality
\begin{equation}\label{rankeq}
\text{\rm rank}_{\Lambda}(\text{\rm Jac}(\frak{PO}_{\frak b})
\otimes_{\Lambda_0}\Lambda) \ge \text{\rm
rank}_{\Lambda}H(X;\Lambda).
\end{equation}
One possible way of proving the inequality (\ref{rankeq}) is to
explicitly construct the inverse homomorphism from
$\text{\rm Jac}(\frak{PO}_{\frak b})$ to $H(X;\Lambda)$. We have not been able
to directly construct such a map on $\text{\rm Jac}(\frak{PO}_{\frak b})$
itself but will construct a map from a ring similar to
$\text{\rm Jac}(\frak{PO}_{\frak b})$ instead. An outline of this construction
is now in order:
\par
In \cite{mc-tol}, McDuff-Tolman proved that the Seidel's
\cite{seidel:auto} representation
$$
\pi_1({\rm Ham}(X)) \to H(X;\Lambda_0)
$$
provides $m$ elements $z'_i \in H(X;\Lambda_0)$ (Recall $m$ is the
number of $T^n$-invariant cycles (or toric divisors) of $X$ of
(real) codimension $2$.) that satisfy the quantum Stanley-Reisner
relations (see  \cite[Definition 6.4]{fooo08} and Section
\ref{sec:seidel} of this paper) and which are congruent modulo
$\Lambda_+$ to the Poincar\'e dual of the generators of the ring of
$T^n$-invariant cycles. See Theorem \ref{mcduftolman}. We use this
to show that the quantum cohomology ring is expressed as the
quotient of $\Lambda_0[z_1,\dots,z_m]$ by an ideal which is a
closure of the ideal generated by quantum Stanley-Reisner relations
and other relations congruent to the relations
$y_i\frac{\partial\frak{PO}_{\frak b}} {\partial y_i} = 0$. In other
words, quantum cohomology $(H(X;\Lambda),\cup^{\frak b})$  is
represented as the quotient of polynomial ring of $m$ variables over
$\Lambda$ by the relations which are congruent to those used for the
Jacobian ring $\text{\rm Jac}(\frak{PO}_{\frak b})$. Using this congruence, we
find a family of schemes parameterized by one variable $s$ which
interpolates $Spec(\text{\rm Jac}(\frak{PO}_{\frak b})\otimes_{\Lambda_0}
\Lambda)$ and $Spec(H(X;\Lambda),\cup^{\frak b})$. This is a family
of zero dimensional schemes. Therefore if it is \emph{proper} and
\emph{flat} over $Spec(\Lambda[s])$, the orders (counted with
multiplicity) of the fibers, which is nothing but the rank over
$\Lambda$ of the coordinate ring of the fiber, are independent of
$s$ and hence (\ref{rankeq}) follows by the standard fact in
algebraic geometry. The proof of the properness is based on the fact
that $\frak{PO}_{\frak b}$ has no critical point whose valuation
lies in a sufficiently small neighborhood of the boundary of the
moment polytope. The proof of Theorem \ref{mtintro}.1 is
completed in Section \ref{sec:injgen}.
\par
In Section \ref{sec:c1} we prove Theorem \ref{c1iscrit}.
We use it in Section \ref{sec:hirze} to calculate the
potential function of second Hirzebruch surface $F_2$ at $\frak b=\text{\bf 0}$.
(The same result was previously obtained by Auroux \cite{Aur09} by a different
method. See also \cite{fooo10}.) We also calculate residue pairing in that case
and verify Theorem \ref{mtintro}.2 by a direct calculation.
\par
Sections \ref{sec:operatorp} - \ref{sec:PDRes} are devoted to the
proof of coincidence of the Poincar\'e duality and the residue
pairing. We summarize the main idea of the proof now.
\par
Using the Poincar\'e duality applied to both $X$ and $L(u)$, we
dualize the map $\frak{ks}_{\frak b}$ and obtain the adjoint map
$\frak{ks}_{\frak b}^*$. A description of this map in terms of the
operator $\frak p$ introduced in  \cite[Section 3.8]{fooobook} is given
in Section \ref{sec:operatorp}. (We recall that another operator
$\frak q$ was used in \cite{fooo09} in the definition of the
potential function with bulk.) We then calculate the pairing
\begin{equation}\label{PDofKSmap}
\langle \frak{ks}_{\frak b}^*(PD[pt]),\frak{ks}_{\frak
b}^*(PD[pt])\rangle_{\text{PD}_X}
\end{equation}
where $PD[pt]$ is the top dimensional cohomology class in
$H(L(u);\Lambda_0)$ as a Floer cohomology class in
$$
HF((L(u),\frak b,b);(L(u),\frak b,b);\Lambda_0) \cong
H(L(u);\Lambda_0).
$$
\par
By inspecting the definition of $\frak{ks}^*_{\frak b}$, we find
that the paring (\ref{PDofKSmap}) can be calculated by counting the
number of holomorphic maps $w$ with boundary lying on $L(u)$ from
certain bordered semi-stable curve $\Sigma$ which is obtained by
attaching two disks at their centers. We also require $w$ to satisfy
an appropriate constraint related to the class $[pt]$. Note that
$\Sigma$ is a degeneration of the annulus, i.e., of a compact
bordered Riemann surface of genus 0 with 2 boundary components. So
calculation of (\ref{PDofKSmap}) is turned into a counting problem
of appropriate holomorphic annuli. By deforming complex structures
of the annuli, we prove that (\ref{PDofKSmap}) is equal to an
appropriate double trace of the operation $\frak m_2$ on
$HF((L(u),\frak b,b);(L(u),\frak b,b);\Lambda_0)$, which is our
definition of residue pairing. (See Theorem \ref{annulusmain}.)
Section \ref{sec:annuli} contains several constructions related to
open-closed Gromov-Witten theory of one-loop. These results are new.
For example, they did not appear in \cite{fooobook} \cite{fooobook2}.
\par
One important point in carrying out this proof
is that we need to reconstruct a filtered $A_{\infty}$ structure
so that it becomes cyclically symmetric. In other words
we need to take the perturbation of the moduli space
of pseudo-holomorphic disks so that it is invariant
under the cyclic permutation of the boundary marked points.
This is because fiber products of the moduli spaces of pseudo-holomorphic
disks appear at the boundary of the moduli spaces of pseudo-holomorphic
annuli where we need cyclic symmetry for the consistency of the
perturbation.
\par
The perturbation we used to define the operator $\frak q$ in
\cite{fooo09} is not cyclically symmetric. In \cite{fooo091}, the
cyclically symmetric perturbations of the moduli space of
pseudo-holomorphic disks are produced. (In \cite{fooo091} it is not
assumed that $X$ is a toric manifold.) For this purpose, {\it
continuous family of} multisections is used therein. In Section
\ref{sec:cyclic} we adapt this scheme to the present situation where
$X$ is toric and $L$ is a $T^n$-orbit, and construct a cyclically
symmetric versions $\frak q^{\frak c}$, $\frak m^{\frak c}$ of our
operators $\frak q$, $\frak m$. We also discuss their relationship
with the operators $\frak q$ and $\frak m$ in \cite{fooo09}.
\par
The proof of Theorem \ref{mtintro} is completed in Section \ref{sec:PDRes}.
\par
If a primitive form takes a simple form $dy_1\wedge \dots \wedge
dy_n /(y_1\cdots y_n)$, the residue pairing is the inverse of the
determinant of the Hessian matrix of the potential function. (See
Definition \ref{res1}.) We will prove that this is the case if
either (1) $n=2$ or (2) $X$ is nef and $\frak b$ is of degree $2$,
as follows. It was proved by Cho \cite{Cho05II} that the Floer
cohomology ring $(HF((L(u),\frak b,b),(L(u),\frak
b,b);\Lambda),\frak m_2)$ is isomorphic to the Clifford algebra
associated to the Hessian matrix of $\frak{PO}_{\frak b}$. This
isomorphism between Clifford algebra and Lagrangian Floer cohomology
holds in general if $\frak{PO}$ is Morse. However in general we do
not know whether this isomorphism respects Poincar\'e duality of
$L(u)$. We prove this under the condition that (1) or (2) above
holds. (See Theorem \ref{clifford}.) In Section \ref{sec:ResHess} we
use it to prove coincidence of the residue pairing (in our sense)
and the inverse of the Hessian determinant modulo the higher order
term in the general case.
\par
In Section \ref{sec:cyclic cohomology} we calculate cyclic (co)homology of
Clifford algebra using Connes' periodicity exact sequence.
We use this to prove that residue paring in our sense
detects the second of the nontrivial elements of cyclic cohomology.
Note the first nontrivial element of cyclic cohomology is detected by
the critical value, which in turn becomes the first Chern class.
(We also review the relation between cyclic cohomology and deformation
of cyclic $A_{\infty}$ algebra.)
\par
In Sections \ref{sec:ori}-\ref{sec:sign} we discuss orientation and sign.
Especially we study the case of moduli space of pseudo-holomorphic annuli,
which was not discussed in  \cite{fooobook2}.
\par
In Sections \ref{sec:equikuracot}-\ref{sec:compwiseplusTnequiv}  we discuss the
construction of Kuranishi structure of the moduli space of
pseudo-holomorphic disks which is invariant under $T^n$ action and
other symmetries. We also construct Kuranishi structure on the 
moduli space of pseudo-holomorphic annulus there.
In Section \ref{sec:cornerestimate} we prove a lemma which we need to
apply Stokes' theorem on the moduli space in the situation where 
nonsmoothness of the forgetful map requires to 
use certain exponential decay estimate.
Those sections are rather technical but
necessary and important substance for the constructions of this paper.
\par
Hochschild and cyclic (co)homology of Lagrangian Floer cohomology
and of Fukaya category is known to be related to the
quantum cohomology of the ambient symplectic manifold.
We prove and use this relationship in this paper. In Section \ref{sec:QMHOch}
we summarize various known facts and previous works on this relationship.
\par
A part of the results of this paper is announced in \cite[Remark 10.3]{fooo09}
and \cite[Remark 6.15]{fooo08}, as well as several
lectures by the authors. 
There is a survey article \cite{fooosurvey} related to the results of this paper.
\par
The authors would like to thank M. Abouzaid for helpful discussions
especially on those related to Theorem \ref{mtintro}.2, Remark
\ref{bunken}.5 and Proposition \ref{clifcyccalc}.
The authors would like to thank Dongning Wang who pointed out an error in Section 4.2 in the previous version.
The authors would like to thank the referee for careful reading and 
many useful comments and suggestions.
The second named author thanks NIMS in Korea for its financial support and
hospitality during his stay in the fall of 2009 when a large
chunk of the current research was carried out.
The first named-author thanks MSRI for its financial support and
hospitality during his stay in the spring of 2010 when
he was working on Section \ref{sec:cyclicKura} of this paper.
The third and fourth named authors thank KIAS for its financial support and
hospitality during their stay in the summer of 2010, when the submitted version of this paper was prepared.
\par 
\section{Notations and terminologies}
\label{sec:notations}
\par
\begin{defn}\label{novring}
Let $R$ be a commutative ring with unit.
We define the {\it universal Novikov ring }$\Lambda_0(R)$ as the set of all formal sums
\begin{equation}\label{novformula}
\sum_{i=0}^{\infty} a_i T^{\lambda_i}
\end{equation}
where $a_i \in R$ and $\lambda_i \in \R_{\ge 0}$ such that
$\lim_{i\to \infty} \lambda_i = \infty$.
Here $T$ is a formal parameter.
\par
We allow $\lambda_i \in \R$ in (\ref{novformula})
(namely negative $\lambda_i$)
to define $\Lambda(R)$ which we call {\it universal Novikov field}. It is a field of fraction of $\Lambda_0(R)$. If $R$ is an algebraic closed field of 
characteristic zero, $\Lambda(R)$ is also algebraic closed. (See \cite[Lemma A.1]{fooo08}.)
\par
We require $\lambda_i > 0$ in  (\ref{novformula}) to define $\Lambda_+(R)$,
which is the maximal ideal of $\Lambda_0(R)$.
\par
In this paper we mainly use the case $R = \C$. In such a case we omit $\C$ and
write $\Lambda_0$\index{$\Lambda_0$}, $\Lambda$\index{$\Lambda$}, 
$\Lambda_+$\index{$\Lambda_+$}.

\end{defn}
\begin{defn}\label{gapping}
A {\it discrete monoid} $G$ is a submonoid of additive monoid $\R_{\ge 0}$ which is discrete.
An element (\ref{novformula}) is said to be {\it $G$-gapped} if
all the exponents $\lambda_i$ are contained in $G$.
\end{defn}
For $\beta \in H_2(X,L;\Z)$,
we denote by
$\mathcal M_{k+1;\ell}^{\text{\rm main}}(\beta)$  the
compactified moduli space of the genus zero bordered holomorphic maps
in class $\beta  \in H_2(X,L(u);\Z)$  with $k+1$ boundary marked points and $\ell$ interior
marked points.
Namely it is the set of all isomorphism classes of
$(\Sigma,\vec z,\vec z^+,u)$, where
$\Sigma$ is a genus zero bordered Riemann surface with one boundary
component, $u : (\Sigma,\partial\Sigma) \to (X,L)$ is a holomorphic
map of homology class $\beta$,
$\vec z = (z_0,\dots,z_k)$ is an ordered set of $k+1$ distinct boundary marked points,
and $\vec z^+ = (z^+_1,\dots,z^+_{\ell})$ is an ordered set of  $\ell$ distinct interior marked points.
(Sometimes, we denote the interior marked points
by $\vec z^{\text {int}} = (z^{\text {int}}_1,\dots,z^{\text {int}}_{\ell})$,
instead of $\vec z^+ = (z^+_1,\dots,z^+_{\ell})$.)
We require the order of $k+1$ boundary marked points respects
the counter-clockwise cyclic order of the boundary $\partial\Sigma$
and also assume appropriate stability conditions.
\par
This space is a compact Hausdorff space and has a Kuranishi structure with boundary and corners.
(See \cite[Section 7.1]{fooobook2}.)
\begin{defn}\label{valuation}
A ($\R$ valued non-Archimedean) {\it valuation} of a commutative
ring $R$ with unit is by definition a map
$
\frak v : R \to \R \cup \{\infty\}
$
such that
\begin{enumerate}
\item $\frak v(ab) = \frak v(a) + \frak v(b)$.
\item $\frak v(a+b) \ge \min\{\frak v(a), \frak v(b)\}$.
\item $\frak v(0) = \infty$.
\end{enumerate}
It uniquely induces a valuation on its field of fractions $F(R)$.
(See \cite[p 41 Proposition 4]{BGR}.)
We denote the induced valuation on $F(R)$ by the same symbol $\frak v$.
A ring $R$ with a valuation $\frak v$ is called a
{\it valuation ring} if
$$
R=\{x \in F(R) \mid \frak v(x) \ge 0\}.
$$
\end{defn}
\par
In general,
a ring $R$ with a valuation $\frak v$
may not coincide with $\{x \in F(R) \mid \frak v(x) \ge 0\}$.
In fact, we will see in Section \ref{sec:statements}
that the universal Novikov ring $\Lambda_0$
is a valuation ring, but we also treat several rings
with valuations which are not valuation rings
in the sense of above.
\par
When a valuation $\frak v$ is given, it defines a natural non-Archimedean norm on
$R$ given by
$$
|a|_{\frak v} := e^{-\frak v(a)}.
$$
In other words it satisfies
$$|ab|_{\frak v} = |a|_{\frak v}|b|_{\frak v}, \qquad
|a+b|_{\frak v} \le \max\{ |a|_{\frak v},|b|_{\frak v}\},
\qquad 
|ca|_{\frak v} = |a|_{\frak v}
$$
for $a,b \in R$ and $c$ in the ground ring.
Therefore it defines a metric topology on $R$.
We say $(R,\frak v)$ is {\it complete} if it is complete with respect to the
metric induced by the norm $|\cdot|_{\frak v}$.
\par
For a $T^n$ orbit $L(u)$ of our toric manifold we denote by
$\text{\bf e}_i$, $i=1,\dots,n$ the standard basis of $H^1(L(u);\Z)$.
Namely $\text{\bf e}_i$ is the Poincar\'e dual to the fundamental class of 
codimension one submanifold $T^{i-1} \times \{\text{point}\} \times T^{n-i}$. 
We write an element $b \in H^1(L(u);\Lambda_0)$ as
$$
b = \sum_{i=1}^n x_i\text{\bf e}_i.
$$
We put $y_i = e^{x_i}$.
\par 
\section{Statement of the results}
\label{sec:statements}
\par
In this section we describe the main results of this paper. We also
review various results of \cite{fooo09} and of several other
previous works and introduce various notations that we use in later sections.
\par
We first discuss potential function and its Jacobian ring. Let $(X,\omega)$
be a compact toric manifold and $\pi: X \to P$ the moment map.
Let $n=\dim_{\C} X$ and $m =$ the number of irreducible components of $D =
\pi^{-1}(\partial P)$ = the number of codimension one faces of $P$.
Let $B$ be the number of faces of $P$ of arbitrary codimension.
Let
$\{P_a\}_{a=1,\dots,B}$ be the set of all faces of $P$ such that
$P_i =
\partial_iP$ $(i=1,\dots,m)$ are the faces of codimension one. Denote $D_a =
\pi^{-1}(P_a)$ which are $T^n$-invariant cycles.
We denote by
$\mathcal A^k$ for $k\ne 0$ the free $\Z$ module generated by
$D_a$'s with $\text{\rm codim} P_a = k$. Then $P_0 = P$ and $D_0 = X$.
Let $\mathcal A^0$ be the rank one free abelian group with its basis $[D_0]$. (We note this notation
is slightly different from one in \cite{fooo09}.
Namely for $k\ne 0$,
the group $\mathcal A^k$ is the same as the one used in
\cite{fooo09}.
But in \cite{fooo09} we put $\mathcal A^0=\{0\}$
instead.)
We put 
$\mathcal A = \mathcal A (\Z)=\oplus _{k}\mathcal A^k$ and
$\mathcal A(\Lambda_0) = \mathcal A\otimes_{\Z} \Lambda_0$,\index{$\mathcal A(\Lambda_0)$} 
$\mathcal A(\Lambda_+) = \mathcal A \otimes_{\Z} \Lambda_+$\index{$\mathcal A(\Lambda_+)$}.
\par
We put $\text{\bf f}_a = [D_a]$, $a=0,1,\dots,B$.
An element $\frak b \in \mathcal A(\Lambda_0)$ is written as
$$
\frak b = \sum_{a=0}^{B} w_a(\frak b)\text{\bf f}_a
$$
so $w_a$'s are coordinates of $\frak b$. Let $y_1,\dots,y_n$ be
formal variables.\index{$y_i$} For $u = (u_1,\dots,u_n) \in \text{Int} P$ we
define $n$ variables $y_j(u)$ ($j=1,\dots,n$)\index{$y_i(u)$} by
\begin{equation}\label{yy(u)identify}
y_j(u) = T^{-u_j} y_j.
\end{equation}
We consider the Laurent polynomial ring:
$$
\Lambda[y_1,y_1^{-1},\dots,y_n,y_n^{-1}],
$$
which is canonically isomorphic to
$$
\Lambda[y_1(u),y_1(u)^{-1},\dots,y_n(u),y_n(u)^{-1}],
$$
by (\ref{yy(u)identify}).
We write this ring $\Lambda[y(u),y(u)^{-1}]$ or $\Lambda[y,y^{-1}]$.
If $u' = (u'_1,\dots,u'_n) \in \text{Int} P$ is another point,  
we have a canonical isomorphism
$$
\psi_{u',u}: \Lambda[y(u),y(u)^{-1}]
\to \Lambda[y(u'),y(u')^{-1}]
$$
that is
$$
\psi_{u',u}(y_j(u)) = T^{u'_j-u_j} y_j(u').
$$
\par
For $i=1,\dots,m$, we denote
$$
\frak w_i = \exp w_i = \sum_{k=0}^{\infty} \frac{1}{k!}w_i^k.
$$
We consider the Laurent polynomial ring
$$
\Lambda[w_0,\frak w_1,\frak w_1^{-1},\dots,\frak w_m,\frak w_m^{-1},
w_{m+1},\dots,w_B,y,y^{-1}],
$$
and denote it by
$$
\Lambda[\frak w,\frak w^{-1},w,y,y^{-1}].
$$
\par
A non-Archimedean valuation $\frak v_T$\index{$\frak v_T$}
on $\Lambda$ is given by
\begin{equation}\label{defval0}
\frak v_T\left(\sum_{i=1}^{\infty} a_i T^{\lambda_i}\right)
= \inf \{ \lambda_i \mid a_i \ne 0\}, \quad
\frak v_T (0) =+\infty.
\end{equation}
We define a valuation $\frak v_T^u$\index{$\frak v_T^u$} on
$\Lambda[\frak w,\frak w^{-1},w,y,y^{-1}]$ by assigning
\begin{equation}\label{defval1}
\frak v_T^u(\frak w_i) = \frak v_T^u(w_i) = \frak v_T^u(y_i(u)) = 0,
\quad \frak v_T^u(T) = 1.
\end{equation}
It implies
\begin{equation}\label{defval2}
\frak v_T^u(y_j(u')) = u_j - u'_j, \quad \frak v_T^u(y_j) = u_j.
\end{equation}
We also define\index{$\frak v_T^P$}
$$
\frak v_T^P = \inf_{u \in P} \frak v_T^u.
$$
Let $\Lambda\langle\!\langle\frak w,\frak
w^{-1},w,y,y^{-1}\rangle\!\rangle^P$ be the completion of $
\Lambda[\frak w,\frak w^{-1},w,y,y^{-1}]$ with respect to the norm
$\frak v_T^P$. We note that $\frak v_T^P$ is not a valuation but
$e^{-\frak v_T^P}$ still defines a norm. We also put
\index{$\Lambda\langle\langle\frak w,\frak w^{-1},w,y,y^{-1}\rangle\rangle_0^P$}
$$
\Lambda\langle\!\langle\frak w,\frak w^{-1},w,y,y^{-1}\rangle\!\rangle_0^P
= \{ x \in \Lambda\langle\!\langle\frak w,\frak w^{-1},w,y,y^{-1}\rangle\!\rangle^P
\mid \frak v_T^P(x) \ge 0\}.
$$
We define
$\Lambda\langle\!\langle y,y^{-1}\rangle\!\rangle^P$,
\index{$\Lambda\langle\langle y,y^{-1}\rangle\rangle^P$} 
$\Lambda\langle\!\langle y,y^{-1}\rangle\!\rangle_0^P$
\index{$\Lambda\langle\langle y,y^{-1}\rangle\rangle_0^P$} 
in a similar way.
\par
Let $\Lambda\langle\!\langle y,y^{-1}\rangle\!\rangle_0^u$ (resp. $\Lambda\langle\!\langle\frak w,\frak
w^{-1},w,y,y^{-1}\rangle\!\rangle_0^u$) be the completion of the rings
$\Lambda_0[y(u),y(u)^{-1}]$ (resp. $\Lambda_0[\frak w,\frak
w^{-1},w,y(u),y(u)^{-1}]$) with respect to $\frak v_T^u$. We denote by $\Lambda\langle\!\langle y,y^{-1}\rangle\!\rangle_+^u$ the
ideal of $\Lambda\langle\!\langle y,y^{-1}\rangle\!\rangle_0^u$ consisting of
all the elements with
$\frak v_T^u >0$. 
\par
To handle the case of $\frak b \in \mathcal A(\Lambda_0)$ we need to
use slightly different completion.
\begin{defn}\label{def:compPo}
We put
\index{$\Lambda\langle\langle y,y^{-1}\rangle\rangle_0^{\overset{\circ}P}$}
\index{$\Lambda\langle\langle\frak w,\frak w^{-1},w,y,y^{-1}\rangle\rangle_0^{\overset{\circ}P}$}
$$
\aligned
\Lambda\langle\!\langle y,y^{-1}\rangle\!\rangle_0^{\overset{\circ}P}
&= \bigcap_{u \in \text{\rm Int}P}\Lambda\langle\!\langle y,y^{-1}\rangle\!\rangle_0^u, \\
\Lambda\langle\!\langle\frak w,\frak w^{-1},w,y,y^{-1}\rangle\!\rangle_0^{\overset{\circ}P} &=
\bigcap_{u\in\text{\rm Int}P}\Lambda\langle\!\langle\frak w,\frak w^{-1},w,y,y^{-1}\rangle\!\rangle_0^u.
\endaligned
$$
Here we say an element $\alpha \in \Lambda\langle\!\langle y,y^{-1}\rangle\!\rangle_0^{u_0}$ is in 
the intersection 
$$
\bigcap_{u \in \text{\rm Int}P}\Lambda\langle\!\langle y,y^{-1}\rangle\!\rangle_0^u
$$
if there exists a sequence $\alpha_i \in \Lambda[y,y^{-1}]$ which converges in $\frak v_T^{u}$ for 
any $u \in \text{\rm Int} P$ and $\lim_{i\to\infty} \alpha_i = \alpha$ in  $\frak v_T^{u_0}$ topology.
(Here $u_0 \in \text{\rm Int} P$. It is easy to see that this definition is independent of $u_0$.)
See Definition \ref{lambdaPe}.
\end{defn}
\begin{rem}
In \cite{fooo09} we used slightly different notations 
$\Lambda^P\langle\!\langle y,y^{-1}\rangle\!\rangle$\index{$\Lambda^P\langle\langle y,y^{-1}\rangle\rangle$}, 
$\Lambda_0^{P}\langle\!\langle y,y^{-1}\rangle\!\rangle$\index{$\Lambda_0^{P}\langle\langle y,y^{-1}\rangle\rangle$} 
etc instead of 
$\Lambda\langle\!\langle y,y^{-1}\rangle\!\rangle^P$, 
$\Lambda\langle\!\langle y,y^{-1}\rangle\!\rangle_0^{P}$. 
\end{rem}
In Section \ref{sec:valuation}, we will define a metric on these rings in Definition
\ref{toplamdaP0} and show that indeed they are the completions of
the corresponding Laurent polynomial rings with respect to the
metric.
\par
We extend the $\frak v_T$ norm on $\Lambda_0$ to $\Lambda_0[Z_1,\dots,Z_m]$
by setting $\frak v_T(Z_i) = 0$, $i = 1,\ldots, m$.
Let $\Lambda_0\langle\!\langle Z_1,\dots,Z_m\rangle\!\rangle$
\index{$\Lambda_0\langle\langle Z_1,\dots,Z_m\rangle\rangle$} 
be the completion of the
polynomial ring $\Lambda_0[Z_1,\dots,Z_m]$ with respect to this $\frak v_T$ norm thereon.
\par
Here and hereafter we use the convention of \cite[Section 2]{fooo08} 
on toric manifolds with $P$ being its moment polytope. In
\cite{fooo08} we use the affine functions $\ell_j : \R^n \to \R$,
($j = 1,\cdots, m$) defined by $\ell_j(u)  = \langle u,\vec v_j \rangle -
\lambda_j$ such that
\begin{equation}\label{defpolytope}
P = \{ u \mid \ell_j(u) \ge 0\} = \{ u \mid \langle u,\vec v_j\rangle \geq \lambda_j \}.
\end{equation}
Here $\vec v_j = (v_{j,i})_{i=1}^n$ is the gradient vector.
We have $v_{j,i} \in \Z$.
\begin{defn}\label{zidef}
We define $z_j \in \Lambda_0[y(u),y(u)^{-1}]$\index{$z_j$} by
\begin{equation}\label{zjdef}
z_j = T^{\ell_j(u)} \prod_{i=1}^n y_i(u)^{v_{j,i}}
= T^{-\lambda_j}\prod_{i=1}^n y_i^{v_{j,i}}.
\end{equation}
Here $\lambda_j = \langle u, \vec v_j \rangle -\ell_j(u)$
for $j=1,\dots ,m$.
If $u \in P$, then $\ell_j (u) \ge 0$ and hence
$z_j$ defines an element of $\Lambda_0[y(u),y(u)^{-1}]$
from the first expression.
Note that we use
(\ref{yy(u)identify}) to get the second expression,
from which
we can also regard
$z_j \in \Lambda[y,y^{-1}]$ with ${\frak v}_T^u(z_j)=
-\lambda_j +\langle u, \vec  v_j \rangle=\ell_j(u) \ge 0$
by (\ref{defval1}) and (\ref{defval2}), if $u \in P$.
\end{defn}
\begin{lem}\label{surjhomfromz}
We have surjective
continuous ring homomorphisms\index{$\Lambda_0\langle\langle Z_1,\dots,Z_m\rangle\rangle$}\index{$\Lambda_0[[Z_1,\dots,Z_m]]$}
\begin{equation}\label{ztoysurj}
\Lambda_0\langle\!\langle Z_1,\dots,Z_m\rangle\!\rangle  \to \Lambda\langle\!\langle y,y^{-1}\rangle\!\rangle_0^{P} ,
\quad \Lambda_0[[Z_1,\dots,Z_m]] \to \Lambda\langle\!\langle y,y^{-1}\rangle\!\rangle_0^{\overset{\circ}P}
\end{equation}
which send $Z_j$ to $z_j$. We also have surjective continuous ring
homomorphisms
\begin{equation}\label{ztoysurj2}
\aligned
\Lambda_0\langle\!\langle\frak w,\frak w^{-1},w,Z_1,\dots,Z_m\rangle\!\rangle
&\to \Lambda\langle\!\langle\frak w,\frak w^{-1},w,y,y^{-1}\rangle\!\rangle_0^{P} \\
\Lambda_0\langle\!\langle\frak w,\frak w^{-1},w\rangle\!\rangle [[Z_1,\dots,Z_m]]
&\to \Lambda\langle\!\langle\frak w,\frak w^{-1},w,y,y^{-1}\rangle\!\rangle_0^{\overset{\circ}P},
\endaligned
\end{equation}
which send $Z_j$ to $z_j$ and preserve $\frak w$, $w$.
\end{lem}
The proof will be given in Section \ref{sec:valuation}.

\par
\begin{defn}\label{def24}
Let $G$ be a discrete submonoid of $\R_{\ge 0}$.
\index{$\Lambda\langle\langle y,y^{-1}\rangle\rangle_{+}^{P}$}
\index{$\Lambda\langle\langle y,y^{-1}\rangle\rangle_{+}^{\overset{\circ}P}$}
\begin{enumerate}
\item
An element of  $\Lambda_0[[Z_1,\dots,Z_m]]$ is said to be
{\it $G$-gapped} if the exponents of $T$ appearing in all of its coefficients are contained
in $G$.
\item
An element of $\Lambda\langle\!\langle y,y^{-1}\rangle\!\rangle_0^{P}$ or of
$\Lambda\langle\!\langle y,y^{-1}\rangle\!\rangle_0^{\overset{\circ}P}$ is said to be {\it $G$-gapped} if
it is an image of a $G$-gapped element by the surjective
homomorphism (\ref{ztoysurj})
\item
$\Lambda\langle\!\langle y,y^{-1}\rangle\!\rangle_{+}^{\overset{\circ}P}$  denotes the image of
$\Lambda_{+}[[Z_1,\dots,Z_m]]$ under the map (\ref{ztoysurj}) in
$\Lambda\langle\!\langle y,y^{-1}\rangle\!\rangle_{0}^{\overset{\circ}P}$.
\item
We put
$
\Lambda\langle\!\langle y,y^{-1}\rangle\!\rangle_{+}^{P}
=
\Lambda\langle\!\langle y,y^{-1}\rangle\!\rangle_{+}^{\overset{\circ}P}
\cap
\Lambda\langle\!\langle y,y^{-1}\rangle\!\rangle^{P}
$.
\end{enumerate}
\end{defn}
\begin{rem}
Any element of $\Lambda\langle\!\langle y,y^{-1}\rangle\!\rangle_0^{P}$ is $G$-gapped for some
discrete submonoid $G$ of $\Z_{\geq 0}$. However this is not true
for the elements of $\Lambda\langle\!\langle y,y^{-1}\rangle\!\rangle_0^{\overset{\circ}P}$ in
general.
\end{rem}
Let
$\text{\bf e}_i$, $i=1,\dots,n$ be the baisis of
$H^1(L(u);\Z)$ and write an element $b \in H^1(L(u);\Lambda_0)$ as
$
b = \sum_{i=1}^n x_i\text{\bf e}_i.
$
We put $y_i = e^{x_i}$.
\par
In \cite{fooo09}, we defined the potential function with bulk
$\frak{PO}(w_1,\dots,w_B;y_1, \dots,y_n)$, that is, a `generating
function of open-closed Gromov-Witten invariant' given in  \cite[Section 3.8]{fooobook}. We enhanced this function to
$$
\frak{PO}(w_0,w_1,\dots,w_B;y_1, \dots,y_n) = w_0 +
\frak{PO}(w_1,\dots,w_B;y_1, \dots,y_n)
$$
by adding additional variable $w_0$ corresponding to the top
dimensional cycle $P^0 = P$ (and $D^0 = X$). Hereafter in this
paper, $\frak{PO}$ will always denote this enhanced function given
in the left hand side. We also denote
$$
\frak{PO}_{\frak b}(y_1,
\dots,y_n)
=
\frak{PO}(w_0,w_1,\dots,w_B;y_1,
\dots,y_n)
$$
where $\frak b = \sum_{i=0}^B w_i\text{\bf f}_i$. The following
lemma and the remark thereafter justify our consideration of
$\Lambda\langle\!\langle y,y^{-1}\rangle\!\rangle_0^{\overset{\circ}P}$.
\par
We recall
\index{$\frak{PO}_{\frak b}(y_1,\dots,y_n)$}
\begin{equation}\label{PPformula}
\aligned
\frak{PO}_{\frak b}(y_1,\dots,y_n) &=
\sum_{\ell,\beta}\sum_{j_1,\dots,j_{\ell}}
\frac{w_{j_1}\dots w_{j_{\ell}}}{\ell !}
T^{\beta \cap \omega/2\pi} y_1(u)^{\partial \beta
\cap \text{\bf e}_1}\cdots y_n(u)^{\partial \beta \cap \text{\bf e}_n} \\
&\quad\qquad \text{\rm ev}_0[\mathcal M_{1;\ell}(\beta;D_{j_1},\dots, D_{j_{\ell}})] \cap [L]
\\
&=\sum_{\ell,\beta}
\frac{1}{\ell !}
T^{\beta \cap \omega/2\pi} y_1(u)^{\partial \beta
\cap \text{\bf e}_1}\cdots y_n(u)^{\partial \beta \cap \text{\bf e}_n} \\
&\quad\qquad \text{\rm ev}_0[\mathcal M_{1;\ell}(\beta;\frak b,\dots,\frak b)] \cap [L].
\endaligned
\end{equation}
See \cite[(7.2)]{fooo09}.
(The moduli space $\mathcal M_{1;\ell}(\beta;D_{j_1},\dots, D_{j_{\ell}})$
is defined at the beginning of Section \ref{sec:frakqreview}.)
Note that $y_i$ and $y_i(u)$ are related by (\ref{yy(u)identify}).
\begin{lem}\label{PPnotPO}
\begin{enumerate}
\item
If $\frak b \in \mathcal A(\Lambda_+)$, then
$
\frak{PO}_{\frak b} \in \Lambda\langle\!\langle y,y^{-1}\rangle\!\rangle_0^P.
$
\item
If $\frak b \in \mathcal A(\Lambda_0)$, then we have
$
\frak{PO}_{\frak b} \in \Lambda\langle\!\langle y,y^{-1}\rangle\!\rangle_0^{\overset{\circ}P}.
$
\end{enumerate}
\end{lem}
\begin{proof}
Statement 1 is proved in  \cite[Theorem 3.11]{fooo09}.
Since $\frak{PO}_{\frak b}$ converges in
$\frak v_T^u$ norm for any $u
\in \text{\rm Int} P$, Statement 2 follows.
\end{proof}
\begin{rem}
We do not know whether
$\frak{PO}_{\frak b} \notin \Lambda\langle\!\langle y,y^{-1}\rangle\!\rangle_0^P$
or not, in case of general $\frak b \in \mathcal A(\Lambda_0)
\setminus \mathcal A(\Lambda_+)$.
\end{rem}
\begin{defn} \label{def28}
We define the Jacobian ring $\text{\rm Jac}(\frak{PO})$ by
$$
\text{\rm Jac}(\frak{PO}) = \frac{ \Lambda\langle\!\langle\frak w,\frak
w^{-1},w,y,y^{-1}\rangle\!\rangle_0^{\overset{\circ}P}} {\text{\rm Clos}_{d_{\overset{\circ}P}}\left(\left\{y_i
\frac{\partial\frak{PO}}{\partial y_i}\mid
i=1,\dots,n\right\}\right)}.
$$
Here $\text{\rm Clos}_{d_{\overset{\circ}P}}$ stands for the closure with respect to the
topology induced by the metric (\ref{metricoverPnot}).
\end{defn}
\begin{defn}
If $\frak b = \sum w_i(\frak b) \text{\bf f}_i \in \mathcal A(\Lambda_0)$, we define
\index{$\text{\rm Jac}(\frak{PO}_{\frak b})$}
$$
\text{\rm Jac}(\frak{PO}_{\frak b})
= \frac{ \Lambda\langle\!\langle y,y^{-1}\rangle\!\rangle_0^{\overset{\circ}P}}
{\text{\rm Clos}_{d_{\overset{\circ}P}}(\{y_i \frac{\partial\frak{PO}_{\frak b}}{\partial y_i}
 \mid i=1,\dots,n\})}.
$$
In other words,
$$
\text{\rm Jac}(\frak{PO}_{\frak b})
= \text{\rm Jac}(\frak{PO}) \otimes_{\Lambda\langle\!\langle\frak w,\frak w^{-1},w,y,y^{-1}\rangle\!\rangle_0^{\overset{\circ}P}}\Lambda\langle\!\langle y,y^{-1}\rangle\!\rangle_0^{\overset{\circ}P}.
$$
Here we regard $\Lambda\langle\!\langle y,y^{-1}\rangle\!\rangle_0^{\overset{\circ}P}$ as a
$\Lambda\langle\!\langle\frak w,\frak
w^{-1},w,y,y^{-1}\rangle\!\rangle_0^{\overset{\circ}P}$-module via the ring homomorphism which sends
$w_i$ to $w_i(\frak b)$.
\end{defn}
\begin{rem}
We can prove that the ideal generated by
$\{y_i \frac{\partial\frak{PO}_{\frak b}}{\partial y_i}
\mid i=1,\dots,n\}$ is closed if $\frak b \in \mathcal A(\Lambda_+)$. (Lemma \ref{polyapprox}.2.)
So, in that case, $\text{\rm Clos}_{d_{\overset{\circ}P}}$ in the above
definition is superfluous.
\end{rem}
\begin{rem}
If $\frak b \in \mathcal A(\Lambda_+)$, we have
\index{$\text{\rm Jac}(\frak{PO}_{\frak b})$}
\begin{equation}\label{JacPOK}
\text{\rm Jac}(\frak{PO}_{\frak b})
\cong
\frac{ \Lambda\langle\!\langle y,y^{-1}\rangle\!\rangle_0^{P} }
{\text{\rm Clos}_{\frak v_T^P}(\{y_i \frac{\partial\frak{PO}_{\frak b}}{\partial y_i}
 \mid i=1,\dots,n\})},
\end{equation}
where the closure $\text{\rm Clos}_{\frak v_T^P}$ is taken with respect to the
norm $e^{-\frak v_T^P}$. (See Proposition \ref{extraJac}.)
\par
In case $\frak b \in \mathcal A(\Lambda_0)\setminus\mathcal A(\Lambda_+)$,
there is an example where either $\frak P_{\frak b} \notin \Lambda\langle\!\langle y,y^{-1}\rangle\!\rangle_0^{P}$ 
or
the equality (\ref{JacPOK}) does not hold.
See Example \ref{boundarycrit}.
\end{rem}
\par
We now recall that
\begin{equation}\label{Htohomology}
\pi: \mathcal A(\Lambda_0) \to H(X;\Lambda_0); \quad \frak b \mapsto [\frak b]
\end{equation}
is a surjective map but not injective. We can choose a subset
$\{\overline{\text{\bf f}}_i \mid i=0,\dots,B'\}$ of our basis
$\{\text{\bf f}_i\}$ so that $\pi(\overline{\text{\bf f}}_i)$ forms
a basis of $H(X;\Z)$ and $\overline{\text{\bf f}}_0 = {\text{\bf
f}}_0$, and $\deg \overline{\text{\bf f}}_i = 2$ if and only if
$i\le m'$. With respect to such a basis, we identify $H(X;\Lambda_0)$ with
the subspace of $\mathcal A(\Lambda_0)$ generated by $\{\overline{\text{\bf f}}_i\}$.
Let $\overline w_i$ be the corresponding coordinates and
$\overline{\frak w}_i$ their exponentials. We define
$$
\Lambda\langle\!\langle \overline{\frak w},\overline{\frak w}^{-1},\overline
w,y,y^{-1}\rangle\!\rangle_0^P, \quad \Lambda\langle\!\langle\overline{\frak
w},\overline{\frak w}^{-1},\overline w,y,y^{-1}\rangle\!\rangle_0^{\overset{\circ}P}
$$
in the same way as we did for $\Lambda\langle\!\langle \frak w,\frak
w^{-1},w,y,y^{-1}\rangle\!\rangle_0^P$, $\Lambda\langle\!\langle \frak w,\frak
w^{-1},w,y,y^{-1}\rangle\!\rangle_0^{\overset{\circ}P}$. By restriction, we may regard $\frak{PO}$ as
$$
\frak{PO} \in \Lambda\langle\!\langle\overline{\frak w},\overline{\frak w}^{-1},\overline w,y,y^{-1}
\rangle\!\rangle_0^{\overset{\circ}P}.
$$
We then define another quotient ring
\begin{equation}\label{jacprime}
\text{\rm Jac}'(\frak{PO})
= \frac{ \Lambda\langle\!\langle\overline{\frak w},\overline{\frak w}^{-1},\overline w,y,y^{-1}
\rangle\!\rangle_0^{\overset{\circ}P}}
{\text{\rm Clos}_{d_{\overset{\circ}P}}(\{y_i \frac{\partial\frak{PO}}{\partial y_i}
 \mid i=1,\dots,n\})}.
\end{equation}
(The difference between $\text{\rm Jac}'$ and $\text{\rm Jac}$ is the
number of variables we use.)
\par
We next define a residue pairing on $\text{\rm Jac}(\frak{PO}_{\frak b})$. For
further discussion, we need to restrict ourselves to the case where
$\frak{PO}_{\frak b}$ is a Morse function. We first make this notion
precise.
\par
Define the vector of valuations
$$
\vec{\frak v}_T: (\Lambda \setminus \{0\})^n \to \R^n
$$
by
$$
\vec{\frak v}_T (\frak y_1,\dots,\frak y_n)
= (\frak v_T(\frak y_1),\dots,\frak v_T(\frak y_n)).
$$
\begin{defn}
A point $\frak y = (\frak y_1,\dots,\frak y_n)$ is said to be a {\it critical
point} of $\frak{PO}_{\frak b}$\index{critical
point of $\frak{PO}_{\frak b}$} if the following holds:
\begin{equation}\label{criticaleq}
y_i \frac{\partial\frak{PO}_{\frak b}}{\partial y_i}(\frak y) = 0
\quad i =1,\dots,n
\end{equation}
\begin{equation}\label{eq:inthemp}
\vec{\frak v}_T(\frak y) \in \text{\rm Int}P.
\end{equation}
We denote $\text{\rm Crit}(\frak{PO}_{\frak b})$\index{$\text{\rm Crit}(\frak{PO}_{\frak b})$} the set of critical
points of $\frak{PO}_{\frak b}$.
\end{defn}

The following lemma will be proved in Section \ref{sec:valuation}.

\begin{lem}\label{pointP}
\begin{enumerate}
\item[(A)] The following conditions for $\frak y = (\frak y_1,\dots,\frak y_n)
\in (\Lambda \setminus \{0\})^n$ are equivalent:
\begin{enumerate}
\item[1]
$(\ref{eq:inthemp})$ holds.
\item[2]
The ring homomorphism
$
\Lambda [y,y^{-1}] \to \Lambda, \quad y_i \mapsto \frak y_i
$
extends to a continuous homomorphism:
$
\Lambda\langle\!\langle y,y^{-1}\rangle\!\rangle^{\overset{\circ}P}  \to \Lambda, \quad y_i \mapsto \frak y_i.
$
\end{enumerate}
\item[(B)] Similarly, the following conditions for $\frak y = (\frak y_1,\dots,\frak y_n)
\in (\Lambda \setminus \{0\})^n$ are equivalent:
\begin{enumerate}
\item[1]
$\vec{\frak v}_T(\frak y) \in P.$
\item[2]
The ring homomorphism
$
\Lambda[y,y^{-1}] \to \Lambda, \quad y_i \mapsto \frak y_i
$
extends to a continuous homomorphism:
$
\Lambda\langle\!\langle y,y^{-1}\rangle\!\rangle^{P}  \to \Lambda, \quad y_i \mapsto \frak y_i.
$
\end{enumerate}
\end{enumerate}
\end{lem}

\begin{defn}
\begin{enumerate}
\item
A critical point $\frak y$ is said to be
{\it nondegenerate} if
$$
\det \left[ \frac{\partial^2\frak{PO}_{\frak b}}{\partial y_i\partial y_j}
\right]_{i,j=1}^{i,j=n} (\frak y)
\ne 0.
$$
\item
$\frak{PO}_{\frak b}$ is said to be a {\it Morse function} if
all of its critical points are nondegenerate.
\end{enumerate}
\end{defn}
\begin{prop}\label{Morsesplit}
\begin{enumerate}
\item There is a splitting of Jacobian ring
$$
\text{\rm Jac}(\frak{PO}_{\frak b}) \otimes_{\Lambda_0} \Lambda \cong \prod_{\frak y \in \text{\rm Crit}(\frak{PO}_{\frak b})}
\text{\rm Jac}(\frak{PO}_{\frak b};\frak y)
$$
as a direct product of rings.
Each of $\text{\rm Jac}(\frak{PO}_{\frak b};\frak y)$ is a local ring.
\item $\frak y$ is nondegenerate if and only if
$\text{\rm Jac}(\frak{PO}_{\frak b};\frak y) \cong \Lambda$.
\item
In particular,
if
$\frak{PO}_{\frak b}$ is a Morse function, then there exists an isomorphism
\begin{equation}\label{factorizeJac}
\text{\rm Jac}(\frak{PO}_{\frak b}) \otimes_{\Lambda_0} \Lambda \cong \prod_{\frak y \in \text{\rm Crit}(\frak{PO}_{\frak b})}
\Lambda.
\end{equation}
\end{enumerate}
\end{prop}
When $X$ is Fano and $\frak b=\text{\bf 0}$, this reduces to 
\cite[Proposition 7.10]{fooo08}. We will prove Proposition \ref{Morsesplit} in Section
\ref{sec:localJacobi}.
\begin{defn}\label{1ydef}
We denote by $1_{\frak y} \in \text{\rm Jac}(\frak{PO}_{\frak b};\frak y)$ the unit of
the ring $\text{\rm Jac}(\frak{PO}_{\frak b};\frak y)$.
\end{defn}
We are now ready to define the simplified residue pairing:
\begin{defn}\label{res1}
Assume that $\frak{PO}_{\frak b}$ is a Morse function. We then
define a {\it simplified residue pairing}
$$
\langle \cdot,\cdot \rangle_{\text{\rm res}}^0:
(\text{\rm Jac}(\frak{PO}_{\frak b})
\otimes_{\Lambda_0}\Lambda) \otimes (\text{\rm Jac}(\frak{PO}_{\frak b})
\otimes_{\Lambda_0}\Lambda)
\to \Lambda
$$
by
\begin{equation}\label{rpdefine0}
\langle 1_{\frak y},1_{\frak y'} \rangle_{\text{\rm res}}^0
=
\begin{cases} 0  &\text{if $\frak y \ne \frak y'$,} \\
\left(\det \left[ \frak y_i^u\frak y_j^u\frac{\partial^2\frak{PO}_{\frak b}}{\partial y_i\partial y_j}\right]_{i,j=1}^{i,j=n}
(\frak b;\frak y^u)\right)^{-1} &\text{if $\frak y = \frak y'.$}
\end{cases}
\end{equation}
Here $u = \vec{\frak v}_T(\frak y)$, $T^{u_i}\frak y_i^u
= \frak y_i$ and $\frak{PO}_{\frak b}$ is as in (\ref{PPformula}).
\end{defn}
\begin{rem}
\begin{enumerate}
\item We remark $e^{x_i} = y_i$. Hence $y_i \frac{\partial}{\partial
y_i} = \frac{\partial}{\partial x_i}$. Thus the matrix appearing in
(\ref{rpdefine0}) is the Hessian matrix of $\frak{PO}_{\frak b}$
with respect to the variable $x_i$.
\item In that regard, the value of $\langle 1_{\frak y},1_{\frak y} \rangle^0_{\text{\rm res}}$
is the contribution to the principal part
at the critical point $x_i = \log \frak y_i$ of the oscillatory integral
$$
\int e^{\sqrt{-1}\frak{PO}_{\frak b}(x)/\hbar} dx
$$
in the stationary phase approximation as $\hbar \to 0$.
\end{enumerate}
\end{rem}
In the case when one of the following conditions is satisfied,
$\langle \cdot,\cdot \rangle^0_{\text{\rm res}}$ defined above
gives the residue pairing in Theorem \ref{mtintro}.
\begin{conds}\label{fano2jigen}
$X$ is a compact toric manifold and $\frak b \in H(X;\Lambda_0)$.
We assume one of the following two conditions.
\par\begin{enumerate}
\item $\dim_{\C}X = 2$.
\item $X$ is nef and $\deg \frak b = 2$.
\end{enumerate}
\end{conds}
We remark that a compact toric manifold $X$ is nef if and only if for all
holomorphic map $u : S^2 \to X$ we have $\int u^*c_1(X) \ge 0$.
\par
To define our residue pairing in the general case, we need to 
use the relation of Jacobian ring with $A_{\infty}$ algebra 
associated to the corresponding Lagrangian fiber.
In other words, 
the definition of our residue pairing uses both the B-model (Landau-Ginzburg) side 
and the A-model side, unfortunately.

\begin{defn}\label{Def:Frobeniusalg}\index{Frobenius algebra}
Let $C$ be a $\Z_2$ graded finitely generated free $\Lambda$ module.
A structure of {\it unital Frobenius algebra} of dimension $n$ is
$\langle \cdot,\cdot \rangle : C^k \otimes C^{n-k} \to \Lambda$,
$\cup :  C^k \otimes C^{\ell} \to  C^{k+\ell}$,
$1 \in C^0$, such that:
\begin{enumerate}
\item $\langle \cdot,\cdot\rangle$ is a graded symmetric bilinear form which induces an
isomorphism $x \mapsto (y \mapsto \langle x,y\rangle)$,
$C^k \to {\rm Hom}_{\Lambda}(C^{n-k},\Lambda)$.
\item $\cup$ is an associative product on $C$. $1$ is its unit.
\item  $\langle x\cup y,z\rangle = \langle x,y\cup z\rangle$.
\end{enumerate}
\end{defn}
We remark that we do {\it not} assume that $\cup$ is (graded) commutative.
The cohomology group of an oriented closed manifold becomes a unital Frobenius algebra in an obvious way.
\begin{defn}\label{invariantZ}\index{Frobenius algebra!trace}\index{trace}\index{$Z(C)$}
Let $(C,\langle \cdot,\cdot\rangle,\cup,1)$ be a unital Frobenius algebra.
We take a basis $\text{\bf e}_I$, $I \in \frak I$ of $C$ such that $\text{\bf e}_0$ is the unit.
Let $g_{IJ} = \langle \text{\bf e}_I,\text{\bf e}_J\rangle$ and let
$g^{IJ}$ be its inverse matrix. We define an invariant of $C$ by
\begin{equation}\label{defnformulaZ}
\aligned
Z(C) = \sum_{I_1,I_2,I_3 \in \frak I}\sum_{J_1,J_2,J_3 \in \frak I}
&(-1)^*g^{I_1J_1}g^{I_2J_2}g^{I_30}g^{J_30} \\
&\langle \text{\bf e}_{I_1} \cup \text{\bf e}_{I_2},\text{\bf e}_{I_3}\rangle
\langle \text{\bf e}_{J_1} \cup \text{\bf e}_{J_2},\text{\bf e}_{J_3}\rangle
\endaligned
\end{equation}
where $* =  \deg \text{\bf e}_{I_1}
\deg \text{\bf e}_{J_2}
+ \frac{n(n-1)}{2}$.
We call $Z(C)$ the \emph{trace} of unital Frobenius algebra $C$.
\end{defn}
It is straightforward to check that $Z(C)$ is independent of the choices of the basis.
This invariant is an example of $1$-loop partition functions\index{$1$-loop partition function} and
can be described by the following Feynman diagram.
\par
\hskip3cm
\epsfbox{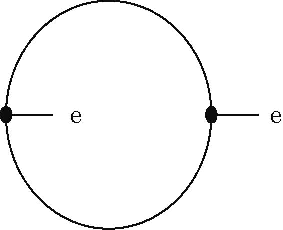}
\par
\centerline{\bf Figure 1.2.1}
\par
We will discuss the invariant $Z$ more in Section \ref{sec:cyclic cohomology}.

\begin{rem}
In \cite{fooo091}, a cyclic unital filtered $A_{\infty}$ algebra is assigned
to an arbitrary relatively spin Lagrangian submanifold $L$ in a symplectic manifold $X$.
So if we fix an element $b$ of the Maurer-Cartan scheme $\mathcal M_{\rm weak}(L;\Lambda_0)$, then
we obtain a unital Frobenius algebra. Therefore $Z$ in Definition \ref{invariantZ}
induces a map
\begin{equation}
Z : \mathcal M_{\rm weak}(L;\Lambda_0) \to \Lambda.
\end{equation}
This map is an invariant of $L$ in that generality.
\par
For the pair $(X,L)$ with $c_1(X) =0$ and $\mu_L = 0$, this
invariant vanishes by the degree reason. There is another invariant
constructed in \cite{fu091} for the case of dimension 3, denoted by $\Psi : \mathcal M(L;\Lambda_0) \to \Lambda_0$,
which is different from $Z$ given here. For example,
the unit did not enter in the definition of $\Psi$ in \cite{fu091}.
\end{rem}
\par
Let $\frak b \in H(X;\Lambda_0)$ and $(u,b) \in \mathfrak M(X;\frak b)$.
In \cite{fooo09} we constructed a (unital) filtered $A_{\infty}$ algebra
$(H(L(u);\Lambda_0),\{\frak m_k^{\frak b,b}\},\text{\rm PD}([L(u)]))$ such that
$\frak m^{\frak b,b}_1 = 0$.
In Section \ref{sec:cyclic} we modify it slightly and obtain another
$(H(L(u);\Lambda_0),\{\frak m_k^{\frak c,\frak b,b}\},\text{\rm PD}([L(u)]))$ which
has the following cyclic property\index{cyclic symmetry} in addition:
\begin{equation}\label{cycsymsign}
\aligned
&\langle \frak m_k^{\frak c,\frak b,b}(x_1,\dots,x_k),x_0\rangle_{{\rm{cyc}}}\\
&= (-1)^{\deg'x_0(\deg'x_1+\dots+\deg'x_k)}\langle\frak m_k^{\frak c,
\frak b,b}(x_0,x_1\dots,x_{k-1}),x_k\rangle_{{\rm{cyc}}}.
\endaligned
\end{equation}
Here we put
\begin{eqnarray}
x \cup^{\frak c,\frak b,b} y &=& (-1)^{\deg x(\deg y + 1)} \frak m^{\frak c,\frak b,b}_2(x,y), \label{deformcup}\\
\langle x, y\rangle_{\text{cyc}} &=& (-1)^{\deg x(\deg y + 1)}\langle x, y\rangle_{\text{\rm PD}_{L(u)}}.
\label{PDandCYC}
\end{eqnarray}
See Remark \ref{rempair} for the precise definition of
$\langle \cdot , \cdot \rangle_{\text{\rm PD}_{L(u)}}$.
Then we can find that $\langle \cdot , \cdot \rangle_{{\rm cyc}}$
satisfies the property (\ref{cycsymsign}),
while the Poincar\'e pairing
$\langle \cdot , \cdot \rangle_{\text{\rm PD}_{L(u)}}$
itself satisfies the property 3 in Definition
\ref{Def:Frobeniusalg}.
Then $(H(L(u);\Lambda),\langle \cdot, \cdot\rangle_{\text{\rm PD}_{L(u)}},
\cup^{\frak c,b,b},\text{\rm PD}[L(u)])$ becomes a unital Frobenius algebra with cyclic
symmetric pairing $\langle \cdot , \cdot \rangle_{{\rm cyc}}$.
We put
\begin{equation}\label{Zbbeq}
Z(\frak b, b) =
Z((H(L(u);\Lambda),\langle \cdot, \cdot\rangle_{\text{\rm PD}_{L(u)}},
\cup^{\frak c,\frak b,b},\text{\rm PD}([L(u)])).
\end{equation}
\begin{defn}\label{res2}
Assume that $\frak{PO}_{\frak b}$ is a Morse function. We then
define a {\it residue pairing}\index{residue pairing}
$$
\langle \cdot,\cdot \rangle_{\text{\rm res}}:
(\text{\rm Jac}(\frak{PO}_{\frak b})
\otimes_{\Lambda_0}\Lambda) \otimes (\text{\rm Jac}(\frak{PO}_{\frak b})
\otimes_{\Lambda_0}\Lambda)
\to \Lambda
$$
by
\begin{equation}\label{rpdefine}
\langle 1_{\frak y},1_{\frak y'} \rangle_{\text{\rm res}}
=
\begin{cases} 0  &\text{if $\frak y \ne \frak y'$,} \\
\left(Z(\frak b, b)\right)^{-1} &\text{if $\frak y = \frak y'$.}
\end{cases}
\end{equation}
\end{defn}
In Section \ref{sec:cyclic}, we also define a unital and filtered $A_{\infty}$
isomorphism $\frak F$
between $(H(L(u);\Lambda_0),\{\frak m_k^{\frak b}\},\text{\rm PD}([L(u)]))$ and
$(H(L(u);\Lambda_0),\{\frak m_k^{\frak c,\frak b}\},\text{\rm PD}([L(u)]))$.
For $b = b_0 + b_+ \in H^1(L(u);\Lambda_0)$,
$b_0 \in  H^1(L(u);\R)$, $b_+ \in H^1(L(u);\Lambda_+)$
we define
$$
b^{\frak c} = b_0 + b_+^{\frak c}$$
by
$$b_+^{\frak c} = \frak F_*(b_+).
$$
We have
$b^{\frak c}_+ \in H^{odd}(L(u);\Lambda_+)$ in general
but $b_+^{\frak c} \in H^1(L(u);\Lambda_+)$ in Case
1 or 2 of Condition \ref{fano2jigen} holds.
We put:
$$
\frak{PO}^{\frak c}_{\frak b} = \frak{PO}_{\frak b} \circ \frak F^{-1}_*.
$$
\par
This is related to Definition \ref{res1} as follows.
\begin{thm}\label{cliffordZ}
\begin{enumerate}
\item
Assume that  $b = \sum x_i \text{\bf e}_i$, $x_i \equiv \overline x_i \mod \Lambda_+$
is a nondegenerate critical point of $\frak {PO}^u_{\frak b}$.
Then
\begin{equation}\label{res0reseq}
Z(\frak b, b)
\equiv \det \left[y_i y_j \frac{\partial^2\frak{PO}_{\frak b}}{\partial y_i\partial y_j}
\right]_{i,j=1}^{i,j=n} (\frak y) \quad \mod T^{\lambda}\Lambda_+.
\end{equation}
Here $\lambda = \frak v_T(Z(\frak b, b))$ and
$\frak y = (e^{x_1},\dots,e^{x_n})$.
\item If Condition \ref{fano2jigen}.1 holds, then
we have
\begin{equation}\label{res0res}
Z(\frak b, b) = \det \left[y_iy_j\frac{\partial^2\frak{PO}^{\frak c}_{\frak b}}{\partial y_i\partial y_j}
\right]_{i,j=1}^{i,j=n} (\frak y).
\end{equation}
\item If Condition \ref{fano2jigen}.2 holds, then
we have
\begin{equation}\label{res1res}
Z(\frak b, b) = \det \left[y_i  y_j\frac{\partial^2\frak{PO}_{\frak b}}{\partial y_i\partial y_j}
\right]_{i,j=1}^{i,j=n} (\frak y).
\end{equation}
\end{enumerate}
\end{thm}
(\ref{res1res}) implies that
\begin{equation}
\langle \cdot,\cdot \rangle_{\text{\rm res}}
= \langle \cdot,\cdot \rangle_{\text{\rm res}}^0
\end{equation}
in case  Condition $\ref{fano2jigen}$.2 holds.
\par
We prove Theorem \ref{cliffordZ} in Sections \ref{sec:clifford} and \ref{sec:ResHess}.
\begin{rem}
The authors do not know a counter example to the equality (\ref{res1res}) in the non-nef case.
But it is unlikely that it holds in the non-nef case.
\end{rem}
\begin{rem}
We define
\begin{equation}
\frak{PO}_{0,\frak b}
= 
\sum_{j=1}^m e^{\frak b \cap D_j} y_1^{v_{j1}} \cdots y_n^{v_{jn}},
\end{equation}
where $\vec v_j = (v_{j1},\dots,v_{jn}) = \text{\rm grad} \,\ell_j$ and 
$D_j$ ($j=1,\dots,m$) is an irreducible component of the toric divisor correspondint to 
$\partial_jP$. The cap  product $\frak b \cap D_j$ is by definition $\frak b_2 \cap D_j$
where $\frak b_2$ is the degree 2 part of $\frak b$.
\par
This function $\frak{PO}_{0,\frak b}$ is the leading order potential function 
(with bulk) that is defined in \cite{fooo08} in the case when $\frak b = 0$.
We remark that $\frak{PO}_{0,\frak b}$ is explicitly calculable.
By definition, (\ref{res0reseq}) implies
\begin{equation}\label{res0reseq2}
Z(\frak b, b)
\equiv \det \left[y_i y_j\frac{\partial^2\frak{PO}_{0,\frak b}}{\partial y_i\partial y_j}
\right]_{i,j=1}^{i,j=n} (\frak y)
\mod T^{\lambda}\Lambda_+.
\end{equation}
\end{rem}
\par
We next turn to the Frobenius manifold structure in the
big quantum cohomology. Since this story is well established
(See \cite{dub}, \cite{Manin:qhm}.), we just briefly recall it
mainly to fix our notation and for readers' convenience.
\par
To simplify the sign issue, we only consider the case where the
cohomology of $X$ consist only of even degrees. (This is certainly the case for
toric manifolds, which is of our main interest in this paper.)
\par
Let $(X,\omega)$ be a symplectic manifold. For
$\alpha \in H_2(X;\Z)$ let
$\mathcal M_{\ell}(\alpha)$ be the moduli space of
stable maps from genus zero semi-stable curves with $\ell$ marked points and of homology class $\alpha$.
There exists an evaluation map
$$
\text{\rm ev}: \mathcal M_{\ell}(\alpha) \to X^{\ell}.
$$
The space $ \mathcal M_{\ell}(\alpha)$ has a virtual fundamental cycle
and hence defines a class
$$
\text{\rm ev}_*[\mathcal M_{\ell}(\alpha)] \in H_{*}(X^{\ell};\Q).
$$
(See \cite{FO}.)
Here
$
* = 2n + 2c_1(X)\cap \alpha + 2\ell - 6
$.
Let $Q_1,\dots,Q_{\ell}$ be cycles such that
\begin{equation}\label{degreecondGW}
\sum \text{\rm codim} \,Q_i = 2n + 2c_1(X)\cap \alpha + 2\ell - 6.
\end{equation}
We define Gromov-Witten invariant by\index{$GW_{\ell}(\alpha;Q_1,\dots,Q_{\ell})$}
$$
GW_{\ell}(\alpha;Q_1,\dots,Q_{\ell}) =
\text{\rm ev}_*[\mathcal M_{\ell}(\alpha)] \cap (Q_1 \times \dots \times Q_{\ell})
\in \Q.
$$
We put $GW_{\ell}(\alpha;Q_1,\dots,Q_{\ell}) = 0$ when (\ref{degreecondGW})
is not satisfied.
\par
We now define\index{$GW_{\ell}(Q_1,\dots,Q_{\ell})$}
\begin{equation}\label{sumGW}
GW_{\ell}(Q_1,\dots,Q_{\ell})
= \sum_{\alpha} T^{\alpha \cap \omega} GW(\alpha;Q_1,\dots,Q_{\ell}).
\end{equation}
(\ref{sumGW}) extends to a $\Lambda_0$ module homomorphism
$$
GW_{\ell}: H(X;\Lambda_0)^{\otimes\ell} \to \Lambda_0.
$$
\begin{defn}\label{def:deformcup}
Let $\frak b \in H(X;\Lambda_0)$ be given. For each given pair $\frak c,\frak d \in H(X;\Lambda_0)$, we define a product
$
\frak c \cup^{\frak b} \frak d \in H(X;\Lambda_0)
$
by the following formula
\begin{equation}\label{defcup}
\langle \frak c \cup^{\frak b} \frak d, \frak e\rangle_{\text{\rm PD}_X}
= \sum_{\ell=0}^{\infty} \frac{1}{\ell!}GW_{\ell+3}(\frak c,\frak d,\frak e,\frak b,\dots,\frak b).
\end{equation}
Here $\langle \cdot,\cdot\rangle_{\text{\rm PD}_X}$ denotes the Poincar\'e duality.
The right hand side converges and defines a graded commutative and associative ring structure on
$H(X;\Lambda_0)$. We call $\cup^{\frak b}$ the {\it deformed quantum cup product}.\index{deformed quantum cup product}
\index{$\cup^{\frak b}$}
\end{defn}
\par
We now discuss convergence of the above deformed quantum cup product
in more detail in an adic topology.
Let $\overline{\text{\bf f}}_i$ ($i=0,\dots,B'$) be a basis of
$H(X;\Lambda_0)$. We assume $\deg \overline{\text{\bf f}}_0 =0$ and
$\deg \overline{\text{\bf f}}_i = 2$ for $i=1,\dots,m'$,
$\deg \overline{\text{\bf f}}_i \ge 4$ for $i=m'+1,\dots,B'$.
Then the divisor axiom of Gromov-Witten invariant (See \cite{Manin:qhm,hori-vafa2}.
For its proof see \cite[Theorem 23.1.4]{FO}  for example.) implies
\begin{equation}\label{divisoraxiom}
\aligned
&GW_{\ell+k}(\alpha;\underbrace{\overline{\text{\bf f}}_i,\dots,
\overline{\text{\bf f}}_i}_{\text{$k$-times}},
Q_{k+1},\dots,Q_{k+\ell}) \\
&= T^{\omega\cap\alpha/2\pi} (\alpha\cap \overline{\text{\bf f}}_i)^{k}
GW_{\ell}(\alpha;Q_{k+1},\dots,Q_{k+\ell})
\endaligned
\end{equation}
for $i=1,\dots,m'$. Now we consider
$
\frak b = \sum_{i=1}^{B'} w_i \overline{\text{\bf f}}_i
$
and put
$$
\frak b_{\rm high} = \sum_{i=m'+1}^{B'} \overline w_i\overline{\text{\bf f}}_i, \qquad
\overline{\frak w}_i = \exp \overline w_i \quad i =1,\dots,m'.
$$
Then \eqref{defcup} can be rewritten as
\begin{equation}\label{defcup2}
\aligned
\langle \frak c \cup^{\frak b} \frak d, \frak e\rangle_{PD_X}
= \sum_{\alpha\in H_2(X;\Z)}\sum_{\ell=0}^{\infty}
\frac{1}{\ell!}&\overline{\frak w}_1^{\alpha\cap \overline{\text{\bf f}}_1}
\dots \overline{\frak w}_{m'}^{\alpha\cap \overline{\text{\bf f}}_{m'}}\\
&GW_{\ell+3}(\alpha;\frak c,\frak d,\frak e,
\underbrace{\frak b_{\rm high},\dots,\frak b_{\rm high}}_{\text{$\ell$-times}}).
\endaligned
\end{equation}
\begin{rem}
Since $\overline{\text{\bf f}}_0$ is a unit, the deformed quantum
product $\cup^{\frak b}$ is independent of $w_0$.
\end{rem}
\par
And as $\ell \to \infty$ in (\ref{defcup2}), we derive
$c_1(X)\cap \alpha \to \infty$ from (\ref{degreecondGW})
and $\alpha \cap \omega\to \infty$ from Gromov's compactness theorem respectively.
Therefore we obtain
\begin{lem}\label{convdefcup}
There exist $\lambda_i\to \infty$, $\lambda_i \ge 0$ and
$$
P_i(\frak c,\frak d,\frak e) \in \Lambda_0[\overline{\frak w}_1,\overline{\frak w}_1^{-1},
\dots,\overline{\frak w}_{m'},\overline{\frak w}_{m'}^{-1},
\overline w_{m'+1},\dots,\overline w_{B'}]
$$
such that
$$
\langle \frak c \cup^{\frak b} \frak d, \frak e\rangle
=
\sum_i T^{\lambda_i} P_i(\frak c,\frak d,\frak e).
$$
\end{lem}
We denote by
$$
\Lambda_0\langle\!\langle\overline{\frak w},\overline{\frak w}^{-1},\overline w\rangle\!\rangle
$$
the completion of
$\Lambda_0[\overline{\frak w}_1,\overline{\frak w}_1^{-1},\dots,\overline{\frak w}_{m'},\overline{\frak w}_{m'}^{-1},
\overline w_{m'+1},\dots,\overline w_{B'}]$
with respect to the non-Archimedean norm on $\Lambda_0$.
Lemma \ref{convdefcup} gives rise to a well-defined map
$$
\cup^{\text{\rm big}}: H(X;\Lambda_0) \otimes_{\Lambda_0}  H(X;\Lambda_0)
\to  H(X;\Lambda_0) \otimes_{\Lambda_0} \Lambda_0\langle\!\langle\overline{\frak w},\overline{\frak w}^{-1},\overline w\rangle\!\rangle
$$
which induces a ring structure on
\begin{equation}\label{bigring}
H(X;\Lambda_0\langle\!\langle\overline{\frak w},\overline{\frak w}^{-1},\overline w\rangle\!\rangle) = H(X;\Lambda_0)
\otimes_{\Lambda_0} \Lambda_0\langle\!\langle\overline{\frak w},\overline{\frak w}^{-1},\overline w\rangle\!\rangle.
\end{equation}
This finishes our description of the family of ring structures $\cup^{\frak
b}$ parameterized by the classes $\frak b \in H(X;\Lambda_0)$.
\par
We now specializes the above to the case of toric manifolds. The following theorem states
that (\ref{bigring}) is isomorphic to the Jacobian ring.
\begin{thm}\label{ringiso}
Let $X$ be a compact toric manifold. Then,
we have a ring isomorphism
$$H(X;\Lambda_0\langle\!\langle\overline{\frak w},\overline{\frak w}^{-1},\overline w
\rangle\!\rangle)
\cong \text{\rm Jac}'(\frak{PO}).
$$
\end{thm}
See (\ref{jacprime}) for $\text{\rm Jac}'(\frak{PO})$.
The isomorphism in Theorem \ref{ringiso} can be described more explicitly.
For this purpose it will be more convenient to fix $\frak b$.
We identify the tangent space
$$
T_{\frak b}H(X;\Lambda_0)
$$
with $H(X;\Lambda_0)$. The coordinate basis of $T_{\frak b}H(X;\Lambda_0)$
of the coordinates $(\overline{\frak w}, \overline{w})$ on $H(X;\Lambda_0)$
is given by the union
$$
\left\{\frac{\partial}{\partial\overline{\frak w}_i}\right\}_{i=1}^{m'} \cup
\left\{\frac{\partial}{\partial\overline w_0}\right\} \cup
\left\{\frac{\partial}{\partial\overline w_j}\right\}_{j=m'+1}^{B'}.
$$
\begin{defn}\label{defn:KSmap}
The {\it Kodaira-Spencer map}\index{Kodaira-Spencer map}
$$
\mathfrak{ks}_{\frak b}: T_{\frak b}H(X;\Lambda_0)
\to \text{\rm Jac}(\frak{PO}_{\frak b})
$$
is defined by
$$\aligned
&\mathfrak{ks}_{\frak b} (\frac{\partial}{\partial\overline{\frak w}_i})
\equiv \frac{\partial}{\partial\overline{\frak w}_i}(\frak{PO}_{\frak b})
\mod \text{\rm Clos}_{d_{\overset{\circ }P}}\left(y_i \frac{\partial\frak{PO}_{\frak b}}{\partial y_i}
 \mid i=1,\dots,n\right),  \\
 &\mathfrak{ks}_{\frak b} (\frac{\partial}{\partial\overline w_i})
\equiv \frac{\partial}{\partial\overline w_i}(\frak{PO}_{\frak b})
\mod \text{\rm Clos}_{d_{\overset{\circ }P}}\left(y_i \frac{\partial\frak{PO}_{\frak b}}
{\partial y_i}
 \mid i=1,\dots,n\right).
\endaligned$$
\end{defn}
\begin{rem} The map
$\mathfrak{ks}$ is the singularity theory analogue
of the Kodaira-Spencer map which is used to study
deformations of complex structures.
\end{rem}

The following is the main theorem of this paper.
\begin{thm}\label{Mirmain}
\begin{enumerate}
\item The Kodaira-Spencer map $\mathfrak{ks}_{\frak b}$ induces a ring isomorphism
$$
\mathfrak{ks}_{\frak b}: T_{\frak b}H(X;\Lambda_0)
\to \text{\rm Jac}(\frak{PO}_{\frak b}).
$$
\item
If $\frak{PO}_{\frak b}$ is a Morse function, then we have
$$
\langle \frak c,\frak d\rangle_{\text{\rm PD}_X}
=
\langle\mathfrak{ks}_{\frak b}\frak c,\mathfrak{ks}_{\frak b}\frak d\rangle_{\text{\rm res}}
$$
for $\frak c,\frak d \in T_{\frak b}H(X;\Lambda_0) = H(X;\Lambda_0)$,
where the left hand side is the Poincar\'e duality.
\end{enumerate}
\end{thm}
\begin{proof}[Proof of Theorem \ref{number}]
Theorem \ref{number}  immediately follows from Theorem \ref{Mirmain}.1,
Proposition \ref{Morsesplit} in this present paper and 
\cite[Theorem 3.16 and Proposition 8.24]{fooo09}.
\end{proof}
The following is immediate from Theorem  \ref{Mirmain}.1 and
Proposition \ref{Morsesplit}.
\begin{cor}
If $\frak{PO}_{\frak b}$ is a Morse function, then
$$
\# \text{\rm Crit}(\frak{PO}_{\frak b})
= \text{\rm rank}_{\Lambda} \text{\rm Jac}(\frak{PO}_{\frak b})
= \text{\rm rank}_{\Q} H(X;\Q).
$$
\end{cor}
Sections \ref{sec:valuation}-\ref{sec:PDRes} of this paper are devoted to the proof of Theorem \ref{Mirmain}.
\par
We note that Theorem  \ref{Mirmain} implies the following mysterious identity:
\begin{cor}
If $\frak{PO}^{\frak b}$ is a Morse function, we have
\begin{equation}\label{sumformula}
\sum_{(u,b) \in \frak M(X,\frak b)} \frac{1}{Z(u,b)} = 0.
\end{equation}
\end{cor}
\begin{proof}
Let $1_X \in H^0(X;\Lambda)$ be the unit. Then
$$
\langle 1_X,1_X \rangle_{\text{\rm PD}_X} = 0.
$$
By Proposition \ref{Morsesplit} we have
$$
1_X = \sum_{\frak y \in \frak M(X,\frak b)} 1_{\frak y}
$$
where $1_{\frak y}$ is the unit of the Jacobian ring $\text{\rm Jac}(\frak{PO}_{\frak b};\frak y)$.
(\ref{sumformula}) now follows from (\ref{rpdefine}) and Theorem \ref{Mirmain}.2.
\end{proof}
\begin{exm}
We take monotone toric blow up of $\C P^2$ at one point.
Its potential function at the (unique) monotone fiber located at $u_0$ is:
\begin{equation}\label{0236}
T^{1/3}(y_1 + y_2 + (y_1y_2)^{-1} + y_1^{-1}).
\end{equation}
See  \cite[Example 8.1]{fooo08}.
(In (\ref{0236}) and in this example, we write $y_i$ instead of 
$y_i(u_0)$.)
The condition for $(y_1,y_2)$ to be critical is
\begin{equation}\label{230}
1 - y_1^{-2}y_2^{-1} - y_1^{-2} = 0, \quad
1 - y_1y_2^2 = 0.
\end{equation}
Therefore $y_1y_2=1/y_2$ and
\begin{equation}\label{231}
y_2^4 + y_2^3 -1 = 0.
\end{equation}
By Theorem \ref{cliffordZ}.3 we have
\begin{equation}
Z(0,(y_1,y_2))
= T^{2/3} \text{\rm det}
\left[
\begin{matrix}
y_1 + (y_1y_2)^{-1} + y_1^{-1}   &   (y_1y_2)^{-1} \\
(y_1y_2)^{-1}           &
y_2 + (y_1y_2)^{-1}
\end{matrix}
\right].
\end{equation}
Using (\ref{230}), (\ref{231}) we have
\begin{equation}
\aligned
T^{-2/3}Z(0,(y_1,y_2))
&= \text{\rm det}
\left[
\begin{matrix}
y_2 + y_2^2 + y_2^{-2}   &   y_2 \\
y_2          &
2y_2\end{matrix}
\right] \\
&= 2y_2^3 + y_2^2 + 2y_2^{-1} =
\frac{4-y_2^3}{y_2}.
\endaligned
\end{equation}
Let $z_i$ ($i=1,2,3,4$) be the 4 solutions of (\ref{231}). Then the left hand side of
(\ref{sumformula}) is:
\begin{equation}\label{hennasiki}
T^{-2/3} \sum_{i=1}^4 \frac{z_i}{4-z_i^3}.
\end{equation}
We have $z^4+z^3-1 = -(z+1)(4-z^3) + 4z +3$. Therefore (\ref{hennasiki}) is
$$
\aligned
T^{-2/3} \sum_{i=1}^4 \frac{z_i(z_i+1)}{4z_i+3}
&= \frac{T^{-2/3}}{16}
\left(
\sum_{i=1}^4 (4z_i+1) - \sum_{i=1}^4\frac{3}{4z_i+3}\right)\\
&=- \frac{3T^{-2/3}}{16} \sum_{i=1}^4\frac{1}{4z_i+3}.
\endaligned
$$
We use $\sum z_i = -1$ in the second equality.
\par
We have $4(z^4+z^3-1) = (z^3+z^2/4-3z/16 +9/64)(4z+3) - 283/64$.
Therefore
$$
\sum_{i=1}^4\frac{1}{4z_i+3}
= \frac{64}{283}\sum_{i=1}^4\left(z_i^3 + \frac{z_i^2}{4}-\frac{3z_i}{16}+\frac{9}{64}\right).
$$
By (\ref{231}) we have $\sum z_i^3 = \sum z_i =-1$, $\sum z_i^2 = 1$. Therefore
$$
\frac{64}{283}\sum_{i=1}^4\left(z_i^3 + \frac{z_i^2}{4}-\frac{3z_i}{16}+\frac{9}{64}\right)
= \frac{64}{283}\left( -1 + \frac{1}{4} + \frac{3}{16} + 4 \times \frac{9}{64}\right)
= 0.
$$
We have thus checked (\ref{sumformula}) directly in the case of monotone toric one point
blow up of $\C P^2$.
\end{exm}
\par

\section{Projective space: an example}
\label{sec:example} Before starting the proof of Theorem \ref{Mirmain}, we illustrate
the theorem for the case $X = \C P^n$ with $\frak b = \text{\bf 0}$
by doing some explicit calculations. In this case the moment polytope
$P$ is
$
\left\{(u_1,\dots,u_n) \mid 0 \le u_i, \sum u_i \le 1\right\}
$
and the potential function is
$$
\frak{PO}_{\text{\bf 0}} = \frak{PO}_{\frak b=\text{\bf 0}} = \sum_{i=1}^n y_i + T
(y_1\cdots y_n)^{-1}.
$$
The critical points are
$\frak y^{(k)} = (T^{\frac{1}{n+1}} e^{\frac{2\pi\sqrt{-1}k}{n+1}}, \dots ,T^{\frac{1}{n+1}} e^{\frac{2\pi\sqrt{-1}k}{n+1}})$, $k=0,\dots,n$
which are all non-degenerate. Therefore we have
$$
\text{\rm Jac}(\frak{PO}_{\text{\bf 0}})\otimes_{\Lambda_0} \Lambda
\cong \prod_{k=0}^{n} \Lambda 1_{\frak y^{(k)}}.
$$
The isomorphism is induced by
$$
\Lambda\langle\!\langle y,y^{-1}\rangle\!\rangle^P  \to \prod_{k=0}^{n} \Lambda 1_{\frak y^{(k)}},
\quad
\mathcal P \mapsto \sum_{k=0}^n\mathcal P(\frak y^{(k)}) 1_{\frak y^{(k)}}.
$$
If we put
$$
\text{\bf f}_k = \pi^{-1}(\{(u_1,\dots,u_n) \in P \mid u_i = 0,
i = n-k+1, \dots, n\}),
$$
$\{\text{\bf f}_0,\dots,\text{\bf f}_n\}$ forms a basis of $H(\C P^n;\Lambda)$.
\par
We derive
$$
\frak{PO}_{w\text{\bf f}_1}(y)
= \frak{PO}_{\text{\bf 0}}(y) + (e^w - 1)y_n
$$
from  \cite[Proposition 4.9]{fooo09} and hence
\begin{equation}\label{ksf1}
\frak{ks}_{\text{\bf 0}}(\text{\bf f}_1)
= [y_n] = T^{\frac{1}{n+1}} \sum_{k=0}^n e^{\frac{2\pi\sqrt{-1}k}{n+1}} 1_{\frak y^{(k)}}
\end{equation}
by definition of $\frak{ks}_{\text{\bf 0}}$.
Using the fact that $\frak{ks}_{\text{\bf 0}}$ is a ring homomorphism,
we have
\begin{equation}\label{ksfm}
\frak{ks}_{\text{\bf 0}}(\text{\bf f}_\ell)
= T^{\frac{\ell}{n+1}} \sum_{k=0}^n e^{\frac{2\pi\sqrt{-1}k\ell}{n+1}} 1_{\frak y^{(k)}}.
\end{equation}
Note this holds for $\ell=0$ also since $\text{\bf f}_0$ is a unit
and $\frak{ks}_{\text{\bf 0}}$ is unital.
We note that
$$
(\frak{ks}_{\text{\bf 0}}(\text{\bf f}_1))^{n+1}
= T \sum_{k=0}^n 1_{\frak y^{(k)}}
= T (\frak{ks}_{\text{\bf 0}}(\text{\bf f}_0)).
$$
Recall that
$$
(\text{\bf f}_1)^{n+1} = T\text{\bf f}_0
$$
is the fundamental relation of the small quantum cohomology ring of
$\C P^n$.
\par
We now calculate the residue pairing.
We remark that Condition $\ref{fano2jigen}$.2 holds
in this case.  The Hessian of $\frak{PO}_{\text{\bf 0}}$
is given by
$$
\text{\rm Hess}_{\frak y^{(k)}} \frak{PO}_{\text{\bf 0}}
= \left[ T^{\frac{1}{n+1}} \frac{\partial^2}{\partial x_i\partial x_j}
\left(e^{x_1} + \dots + e^{x_n} + e^{-(x_1 + \dots + x_n)}\right)\right]_{i,j=1}^{i,j=n}
(\frak x^{(k)})
$$
with $\frak x^{(k)}_i = \frac{2\pi\sqrt{-1}k}{n+1}$. Therefore
$$
\text{\rm Hess}_{\frak y^{(k)}} \frak{PO}_{\text{\bf 0}}
=  T^{\frac{1}{n+1}}e^{\frac{2\pi\sqrt{-1}k}{n+1}} \left[
\delta_{ij} + 1\right]_{i,j=1}^{i,j=n}.
$$
It is easy to see that the determinant of the matrix
$\left[
\delta_{ij} + 1\right]_{i,j=1}^{i,j=n}$ is $n+1$. Therefore
the residue pairing is given by
\begin{equation}\label{CPnresdue}
\langle 1_{\frak y^{(k)}}, 1_{\frak y^{(k')}}\rangle_{\text{\rm res}}
= T^{-\frac{n}{n+1}}e^{-\frac{2\pi\sqrt{-1}kn}{n+1}}\frac{\delta_{kk'}}{1+n}.
\end{equation}
Combining (\ref{ksfm}) and (\ref{CPnresdue}), we obtain
\begin{equation}\label{rescalcmm}
\langle \frak{ks}_{\text{\bf 0}}(\text{\bf f}_\ell),
\frak{ks}_{\text{\bf 0}}(\text{\bf f}_{\ell'})\rangle_{\text{\rm res}}
= \frac{1}{n+1}T^{-\frac{n}{n+1}}\sum_{k=0}^n
e^{-\frac{2\pi\sqrt{-1}kn}{n+1}} T^{\frac{\ell+\ell'}{n+1}}
e^{\frac{2\pi\sqrt{-1}(\ell+\ell')k}{n+1}}.
\end{equation}
It follows that (\ref{rescalcmm}) is $0$ unless $\ell+\ell' = n$ and
$$
\langle \frak{ks}_{\text{\bf 0}}(\text{\bf f}_\ell),
\frak{ks}_{\text{\bf 0}}(\text{\bf f}_{n-\ell})\rangle_{\text{\rm res}}
= 1 =
\langle \text{\bf f}_\ell,
\text{\bf f}_{n-\ell}\rangle_{\text{PD}_{\C P^n}}.
$$
Thus Theorem \ref{Mirmain}.2 holds in this case.
\par
We note that in this case a similar calculation can be found
in the earlier literature \cite{taka},  \cite{barani}, \cite{KKP}.
The main result of the present paper generalizes this
calculation to any compact toric manifold $X$ also with
arbitrary bulk deformation parameter $\frak b \in H(X;\Lambda_0)$.
\begin{rem}
In the above discussion of projective spaces, we derived (\ref{ksfm})
from (\ref{ksf1}) and the fact $\frak{ks}_{\text{\bf 0}}$
is a ring homomorphism. We can also prove (\ref{ksfm}) directly:
Let $\beta_i \in \pi_2(X;L(u))$ as in the end of Section \ref{sec:introduction}.
We consider the moduli space
$\mathcal M_{1;1}(\beta_{n-\ell+1}+\dots+\beta_n)$ and its fiber product
$$
\mathcal M_{1}(\beta_{n-\ell+1}+\dots+\beta_n;\text{\bf f}_{\ell})
: =
\mathcal M_{1;1}(\beta_{n-\ell+1}+\dots+\beta_n)
{}_{\text{\rm ev}^{\text{\rm int}}} \times_{X} \text{\bf f}_{\ell}.
$$
Since the Maslov index of $\beta_{n-\ell+1}+\dots+\beta_n$ is $2\ell$,
it follows that the (virtual) dimension of
$\mathcal M_{1}(\beta_{n-\ell+1}+\dots+\beta_n;\text{\bf f}_{\ell})$
is $n$:
Actually by inspecting the classification of holomorphic disks,
\cite[Theorem 5.2]{cho-oh} (see also  \cite[Theorem 4.4]{Cho04}) we can find that
$\mathcal M_{1;1}(\beta_{n-\ell+1}+\dots+\beta_n)$ is Fredholm regular and
the evaluation map
$$
\text{\rm ev}_{\text{\rm int}}:\mathcal M_{1;1}(\beta_{n-\ell+1}+\dots+\beta_n)
\to X
$$
is transverse to the divisor $\text{\bf f}_\ell$. Therefore
$\mathcal M_{1}(\beta_{n-\ell+1}+\dots+\beta_n;\text{\bf f}_{\ell})$ is
smooth, $n$ dimensional and
$$
\text{\rm ev}: \mathcal M_{1}(\beta_{n-\ell+1}+\dots+\beta_n;\text{\bf f}_{\ell})
\to L(u)
$$
is a diffeomorphism. We can also find that $\mathcal M_{1}(\beta;\text{\bf f}_{\ell})$ is empty
for other $\beta$ with Maslov index $2\ell$.
Therefore by definition and (\ref{yy(u)identify}) we have
$$
\aligned
\frak{PO}(w\text{\bf f}_{\ell};y) -
\frak{PO}(\text{\bf 0};y)
& \equiv wT^{u_{n-\ell+1}+\dots+u_n}y_{n-\ell+1}(u)\cdots y_n(u) \mod w^2 \\
& = w y_{n-\ell+1}\cdots y_n
\mod w^2.
\endaligned
$$
Hence
$$
\frak{ks}_{\text{\bf 0}}(\text{\bf f}_{\ell})
= \sum_{k=0}^n \left.\frac{\partial\frak{PO}_{w\text{\bf f}_{\ell}}}{\partial w}
\right\vert_{w=0}(\frak y^{(k)}) 1_{\frak y^{(k)}}
$$
is given by (\ref{ksfm}).
\end{rem}
\begin{rem}
In \cite{gros} the case of $\C P^2$ is studied including bulk deformation by $\frak b \in H^4(\C P^2)$.
\end{rem}
\par

\chapter{Ring isomorphism}
\section{Norms and completions of polynomial ring over
Novikov ring}
\label{sec:valuation}
In this section we discuss various points related to
the ring $\Lambda\langle\!\langle y,y^{-1}\rangle\!\rangle_0^{\overset{\circ}P} $ which was defined to be
$$
\Lambda\langle\!\langle y,y^{-1}\rangle\!\rangle_0^{\overset{\circ}P}  = \bigcap_{u \in \text{Int} P}
\Lambda\langle\!\langle y,y^{-1}\rangle\!\rangle_0^u
$$
in Section \ref{sec:statements}, and prove Lemma \ref{surjhomfromz}
and Lemma \ref{pointP} whose proofs were postponed therein.
This section is rather technical since we need to deal with completion of
the ring of Laurent polynomials with Novikov ring coefficients. In particular,
we need to make precise the definition of a metric on $\Lambda_0[y,y^{-1}]$
used taking the completion which coincides with $\Lambda\langle\!\langle y,y^{-1}\rangle\!\rangle_0^{\overset{\circ}P} $.
\par
For each sufficiently small $\e > 0$, we define the subset $P_\e$ of $P$ by
$$
P_\e = \{ u \in P \mid \ell_j(u) \geq \e,~ j=1,\dots ,m \}.
$$
Here $\ell_j$ are affine functions appeared
in (\ref{defpolytope}).
\begin{defn}\label{lambdaPe}
We define $\Lambda\langle\!\langle y,y^{-1}\rangle\!\rangle^{P_\e}$\index{$\Lambda\langle\langle y,y^{-1}\rangle\rangle^{P_\e}$} by the
completion of $\Lambda[y,y^{-1}]$ with respect to the \emph{bounded}
metric
$$
d_{P_\e}(x,x'): = \min\{1,\exp(- \frak v_T^{P_\e}(x - x'))\}
$$
on $\Lambda_0[y,y^{-1}]$. And $\Lambda\langle\!\langle y,y^{-1}\rangle\!\rangle^{P_\e}_0$\index{$\Lambda\langle\langle y,y^{-1}\rangle\rangle^{P_\e}_0$} is the
subring of elements $x$ with $\frak v_T^{P_{\epsilon}}(x) \ge 0$.
\end{defn}
\begin{defn}\label{toplamdaP0}
Define a metric $d_{\overset{\circ}P}$ on
$\Lambda[y,y^{-1}]$ by
\begin{equation}\label{metricoverPnot}
d_{\overset{\circ}P}(x,x')
= \sum_{k=1}^{\infty} 2^{-k} d_{P_{1/k}}(x,x').
\end{equation}
The case of $\Lambda\langle\!\langle\frak w,\frak
w^{-1},w,y,y^{-1}\rangle\!\rangle_0^{\overset{\circ}P}$ is similar.
\end{defn}

\begin{prop}\label{metricconf}
$\Lambda\langle\!\langle y,y^{-1}\rangle\!\rangle_0^{\overset{\circ}P}$ is the completion
of a subring\index{$\Lambda\langle\langle y,y^{-1}\rangle\rangle_0^{\overset{\circ}P}$}
$$
\Lambda[y,y^{-1}] \cap \Lambda\langle\!\langle y,y^{-1}\rangle\!\rangle_0^{\overset{\circ}P}
$$
of the
Laurent polynomial ring
with respect to the metric $d_{\overset{\circ}P}$.
\end{prop}
\begin{proof} We note that if $\e> 0$ is sufficiently small $P_\e$ is still
convex and so
\begin{equation}\label{vepsilonwrite}
\frak v_T^{P_\e} = \min \{\frak v_T^{(u(\e,j))} \mid j = 1, \dots,
N \}
\end{equation}
where $u(\e,1), \dots, u(\e,N)$ are the vertices of $P_\e$. Therefore we have
$$
\bigcap_{u \in \text{Int} P} \Lambda\langle\!\langle y,y^{-1}\rangle\!\rangle_0^u  =
\bigcap_{k=1}^\infty \Lambda\langle\!\langle y,y^{-1}\rangle\!\rangle_0^{P_{1/k}}
$$
and hence
$$
\Lambda\langle\!\langle y,y^{-1}\rangle\!\rangle_0^{\overset{\circ}P}  = \bigcap_{k=1}^\infty
\Lambda\langle\!\langle y,y^{-1}\rangle\!\rangle_0^{P_{1/k}} .
$$
But it easily follows from the definition of the metric $d_{\overset{\circ}P}$
that the completion of $\Lambda[y,y^{-1}] \cap \Lambda\langle\!\langle y,y^{-1}\rangle\!\rangle_0^{\overset{\circ}P} $
with respect to $d_{\overset{\circ}P}$ coincides with
$$
\bigcap_{k=1}^\infty \Lambda\langle\!\langle y,y^{-1}\rangle\!\rangle_0^{P_{1/k}}
$$
whose proof is standard and so omitted.
\end{proof}
\par
Note that there exists a continuous homomorphism
$$
\Lambda\langle\!\langle y,y^{-1}\rangle\!\rangle_{0}^{P_{\e_1}}  \to
\Lambda\langle\!\langle y,y^{-1}\rangle\!\rangle_{0}^{P_{\e_2}}
$$
for $\e_1 < \e_2$. The following lemma is an immediate consequence of
Proposition \ref{metricconf} and (\ref{vepsilonwrite}).
\begin{lem}\label{Poisprojlim}
$\Lambda\langle\!\langle y,y^{-1}\rangle\!\rangle_0^{\overset{\circ}P} $ is
isomorphic to the projective limit
$
\varprojlim \Lambda\langle\!\langle y,y^{-1}\rangle\!\rangle_{0}^{P_{\e}} .
$
\end{lem}
The following lemma is a consequence of the proof of 
\cite[Lemma 3.10]{fooo09} whose proof we omit but refer thereto.
\par
\begin{lem}\label{sublemma} Suppose $T^{\mu} y^v$ lies in $\Lambda\langle\!\langle y,y^{-1}\rangle\!\rangle_0^P$.
Then there exists $w(v) \in \Z^m_{\geq 0}$ such that
\begin{equation}\label{yvzwv}
T^{\mu} y^v = T^{\lambda(\mu,v)} z^{w(v)}
\end{equation}
for $z^{w(v)} = z_1^{w(v)^{(1)}}\cdots  z_m^{w(v)^{(m)}}$ and
\begin{equation}\label{lambdalambda'}
\frak v_T^P(T^{\mu} y^v) = \lambda(\mu,v).
\end{equation}
\end{lem}

Now we are ready to give the proofs of Lemmata \ref{surjhomfromz}
and \ref{pointP}.

\begin{proof}[Proof of Lemma \ref{surjhomfromz}]
The case of $\Lambda\langle\!\langle y,y^{-1}\rangle\!\rangle_0^{P} $ immediately follows from
 \cite[Lemma 3.10]{fooo09} and so we will focus on the case of
$\Lambda\langle\!\langle y,y^{-1}\rangle\!\rangle_0^{\overset{\circ}P} $.
\par
By definition of $z_j$'s in \eqref{zjdef}, we have $\frak v_T^u(z_j) = \ell_j(u)$.
In particular, if $u \in \text{\rm Int}  P$, we have $\frak v_T^u(z_j) > 0$.
Therefore any formal power series of $z_j$ with coefficients in $\Lambda_0$
converges in $\frak v_T^{u}$ for any $u \in \text{\rm Int} P$.
Hence the map
\begin{equation}\label{formaltoconv}
\Lambda_0[[Z_1,\dots,Z_m]] \to \Lambda\langle\!\langle y,y^{-1}\rangle\!\rangle_0^{\overset{\circ}P}
\end{equation}
is well defined. It remains to prove its surjectivity.

Now let $x$ be an arbitrary element
$\Lambda\langle\!\langle y,y^{-1}\rangle\!\rangle_0^{\overset{\circ}P} $.
Then for any $u \in \text{\rm Int} P$ we can write $x$ as
$$
x = \sum_{i=1}^{\infty} a_iT^{\lambda_i'} y(u)^{v_i}
$$
of with $0 \leq \lambda'_i \le \lambda'_{i+1}$, $\lim \lambda'_i \to \infty$. (Recall $\frak v_T^u(y_j(u)) = 0$.)
We assume that $(\lambda'_i,v_i) \ne (\lambda'_j,v_j)$, and $a_i \in \C\setminus \{0\}$.
We define
$$
\lambda_i = \lambda(\lambda'_i,u,v_i) = \frak v_T^P(T^{\lambda_i'} y(u)^{v_i}), 
\quad w_i = w(v_i) = (w_i^{(1)},\dots,w_i^{(m)}).
$$
For each $N \in \N$, we define
$$
I_N = \{i \mid w_i^{(1)} + \dots + w_i^{(m)} \le N \},
$$
and
\begin{equation}\label{xNsum}
x_N = \sum_{i \in I_N} a_i T^{\lambda_i'} y(u)^{v_i}.
\end{equation}
A priori, (\ref{xNsum}) may be an infinite sum but the sum converges in $\frak v_T^u$ norm for any $u \in
\text{\rm Int}P$: It is a sub-series of a convergent series.
Therefore the sum converges to an element of
$\Lambda\langle\!\langle y,y^{-1}\rangle\!\rangle_0^{\overset{\circ}P} $ in $d_{P_\e}$ for any $\e > 0$,
since $v_T^u$ is an {\it non-Archimedean} valuation for any $u$.
Therefore it converges in $d_{\overset{\circ}P}$ by the definition \eqref{metricoverPnot}.
\par
Using \eqref{yvzwv}, we can rewrite
$$
\sum_{i \in I_N} a_iT^{\lambda_i'} y(u)^{v_i} = \sum_{i \in I_N}
a_iT^{\lambda_i} z^{w_i}
$$
and so the right hand side also converges in $d_{P_{\epsilon}}$.
On the other hand, if $\vert w_i^{(1)} + \dots + w_i^{(m)} \vert\le
N$, we derive
$$
\vert \frak v_T^{P_{\epsilon}}(a_iT^{\lambda_i} z^{w_i})
- \frak v_T^{P}(a_iT^{\lambda_i} z^{w_i})  \vert
\le \epsilon N
$$
from \eqref{vepsilonwrite} and hence the sum also converges in $d_P$ and
hence $x_N \in \Lambda\langle\!\langle y,y^{-1}\rangle\!\rangle^P $. In particular,
$$
\{ i \mid \lambda_i < C,  w_i^{(1)} + \dots + w_i^{(m)} < N\}
$$
is a finite set for any $C$ and $N$. Therefore we have proved that
$$
\sum_{i \in I_N} a_i T^{\lambda_i} z^{w_i}
$$
is a polynomial with $\lambda_i \geq 0$: Recall
$\lambda_i = \frak v_T^P(T^{\lambda_i'} y(u)^{v_i}) = \frak v_T^P(T^{\lambda_i}z^{w_i}) \geq 0$.
By construction, we have
$$
x = \lim_{N \to \infty} x_N
$$
with respect to the distance $d_{\overset{\circ}P}$ and so $x$ is the image of the formal power series
$$
\sum_{i \in \cup_{N=0}^\infty I_N} a_i T^{\lambda_i}Z^{w_i}.
$$
This implies the surjectivity of (\ref{formaltoconv}).
\par
The proof of the other two cases are similar and is left to the
reader.
\end{proof}
\begin{proof}[Proof of Lemma \ref{pointP}] (A)
1 $\Rightarrow$ 2 is immediate from Lemma \ref{surjhomfromz}.
\par
We next prove 2 $\Rightarrow$ 1. Recall from Definition
\ref{zidef} that
$
z_j = T^{-\lambda_j}\prod_{i=1}^n y_i^{v_{j,i}}
$
($j = 1, \dots, m$).
By the homomorphism $y_i \to \frak y_i$,
$z_j$ is mapped to $\frak z_j =T^{-\lambda_j}\prod_{i=1}^n {\frak y}_i^{v_{j,i}}$.
For each $j = 1, \dots, m$,
$\sum_{k=1}^{\infty} z^k_j$ converges in
$\Lambda\langle\!\langle y,y^{-1}\rangle\!\rangle_0^{\overset{\circ}P}$.
This is because $\frak v_T^u (z_j) = \ell_i(u) > 0$ for any $u \in \text{Int }P$.
By the continuity hypothesis of the evaluation homomorphism substituting $y \mapsto \frak y$,
it follows that $\sum_{k=0}^{\infty} \frak z_j^k$ converges in $\Lambda_0$.
Therefore
$$
\frak v_T(\frak z_j) > 0.
$$
Since we have $\frak z_j = T^{-\lambda_j} \prod_{i=1}^n \frak y_i^{v_{j,i}}$ from
\eqref{zjdef}, we have
$
0 < \frak v_T(\frak z_j) = -\lambda_j + \langle v_j, \vec u \rangle = \ell_j(\vec u)
$
for all $j =1, \dots, m$ with
$
\vec u = \vec{\frak v}_T(\frak y).
$
By the definition of the moment polytope $P$, $\vec u \in \text{Int }P$ which
implies \eqref{eq:inthemp}.
\par
(B) 1 $\Rightarrow$ 2 also follows from Lemma \ref{surjhomfromz}. To prove 2 $\Rightarrow$ 1, we consider
a series $\sum_{k=1}^{\infty} T^{\epsilon k}z^k_j$ for
$\epsilon >0$. Then it converges
in $\Lambda\langle\!\langle y,y^{-1}\rangle\!\rangle_0^{P}$.
Then a similar argument above yields that $\epsilon + \frak v_T(\frak z_j) >0$. Since $\epsilon >0$ can be taken sufficiently small, we have $v_T(\frak z_j) \ge0$.
Then the similar argument shows 1.
\end{proof}
\begin{defn}\label{definitonPcirczero}
We define\index{$\Lambda\langle\langle y,y^{-1}\rangle\rangle^{\overset{\circ}P}$}
$$
\Lambda\langle\!\langle y,y^{-1}\rangle\!\rangle^{\overset{\circ}P}
 =
\Lambda\langle\!\langle y,y^{-1}\rangle\!\rangle_0^{\overset{\circ}P} \otimes_{\Lambda_0}\Lambda.
$$
\end{defn}
\begin{rem}\label{remPcirczero}
We recall that in Section \ref{sec:statements} we define
$\Lambda\langle\!\langle y,y^{-1}\rangle\!\rangle^{P} $
by the completion of $\Lambda[y,y^{-1}]$ with respect to
$\frak v_T^P$.
By considering the homomorphism (\ref{ztoysurj}) over
$\Lambda$ coefficients, we find that
$\Lambda\langle\!\langle y,y^{-1}\rangle\!\rangle^{P}$
coincides with the image of
\begin{equation}\label{formalP}
\Lambda\langle\!\langle Z_1,\dots,Z_m\rangle\!\rangle  \to \Lambda\langle\!\langle y,y^{-1}\rangle\!\rangle^{P}.
\end{equation}
We find that
$$
\Lambda\langle\!\langle y,y^{-1}\rangle\!\rangle^{P}
=
\Lambda\langle\!\langle y,y^{-1}\rangle\!\rangle_0^{P} \otimes_{\Lambda_0}\Lambda.
$$
Similarly, we also find that
$$
\Lambda\langle\!\langle y,y^{-1}\rangle\!\rangle^{P_{\e}}
=
\Lambda\langle\!\langle y,y^{-1}\rangle\!\rangle_0^{P_{\e}} \otimes_{\Lambda_0}\Lambda.
$$
However, we note that (\ref{formalP})
does not extend to
$\Lambda[[Z_1,\dots,Z_m]] \to \Lambda\langle\!\langle y,y^{-1}\rangle\!\rangle^{\overset{\circ}P}$.
\end{rem}
\par
\section{Localization of Jacobian ring at moment polytope}
\label{sec:localJacobi}

In this section, we give the proof of Proposition \ref{Morsesplit}.
We will use the following lemma for the proof.
\begin{lem}\label{Jacfinitedimension}
$\text{\rm Jac}(\frak{PO}_{\frak b}) \otimes_{\Lambda_0} \Lambda$ is a finite dimensional $\Lambda$ vector space.
\end{lem}
Lemma \ref{Jacfinitedimension} is an immediate consequence of Theorem \ref{Mirmain}.1.
(We will use Proposition \ref{Morsesplit} to define residue pairing on Jacobian ring.
The residue pairing will be used to state Theorem \ref{Mirmain}.2. However the proof of
Theorem \ref{Mirmain}.1 will be completed, without using Proposition \ref{Morsesplit},
in Section \ref{sec:surf} when the surjectivity of
$\frak{ks}_{\frak b}$ is established. Thus the logic of arguments used in the proofs
will not be circular.)
\par
\begin{proof}[Proof of Proposition \ref{Morsesplit}]
We first prove Proposition \ref{Morsesplit}.1.
For $y_i \in \Lambda\langle\!\langle y,y^{-1}\rangle\!\rangle^{\overset{\circ}P}$ we define a $\Lambda$ linear map
$$
\widehat{y}_i: \text{\rm Jac}(\frak{PO}_{\frak b}) \otimes_{\Lambda_0} \Lambda \to
\text{\rm Jac}(\frak{PO}_{\frak b}) \otimes_{\Lambda_0} \Lambda
$$
as the multiplication by $y_i$ in the Jacobian ring. Since $y_i^{-1} \in
\text{\rm Jac}(\frak{PO}_{\frak b})$ and the map $y_i \mapsto \widehat y_i$ is
a ring homomorphism from $\Lambda\langle\!\langle y,y^{-1}\rangle\!\rangle^{\overset{\circ}P}$ to
${\rm End}_\Lambda(\text{\rm Jac}(\frak{PO}_{\frak b})
\otimes_{\Lambda_0} \Lambda)$, the homomorphism $\widehat y_i$ is
invertible.

For $\frak y = (\frak y_1,\dots,\frak y_n) \in \Lambda^n$ we put
\begin{equation}
\aligned
\text{\rm Jac}(\frak{PO}_{\frak b};\frak y)  = \{ &\mathcal P \in
\text{\rm Jac}(\frak{PO}_{\frak b}) \otimes_{\Lambda_0} \Lambda \mid\\
&(\widehat{y}_i- \frak y_i)^N \mathcal P = 0,\,  i = 1,
\dots, n  \text{ for sufficiently large $N$}\}.
\endaligned
\end{equation}
Since $\widehat y_i$ is invertible, $\text{\rm Jac}(\frak{PO}_{\frak
b};\frak y) \ne 0$ holds only when 
$\frak y_i \ne 0$. We denote
$$
\frak Y = \{ \frak y \mid \text{\rm Jac}(\frak{PO}_{\frak b};\frak y) \ne 0 \}.
$$
Since $\Lambda$ is an algebraically closed field (See for example \cite[Section A]{fooo08}.)
and $\widehat{y}_i \widehat{y}_j
= \widehat{y}_j \widehat{y}_i$,
the finite dimensional $\Lambda$ vector space 
$\text{\rm Jac}(\frak{PO}_{\frak b}) \otimes_{\Lambda_0} \Lambda$ 
is decomposed to simultaneous (generalized) eigenspaces of $\widehat{y}_j$'s.
Thus we have a decomposition
\begin{equation}\label{decompJac}
\text{\rm Jac}(\frak{PO}_{\frak b}) \otimes_{\Lambda_0} \Lambda \cong
\prod_{\frak y \in \frak Y} \text{\rm Jac}(\frak{PO}_{\frak b};\frak y)
\end{equation}
as a $\Lambda$ vector space. Using the fact $\text{\rm Jac}(\frak{PO}_{\frak b})$ is 
commutative, we can prove that (\ref{decompJac}) is a
direct product decomposition as a ring.
(See the proof of  \cite[Proposition 6.7]{fooo08}.)
\begin{lem}\label{inthemoment}
If $\text{\rm Jac}(\frak{PO}_{\frak b};\frak y) \ne 0$,
then $\vec{\frak v}_T(\frak y) \in \text{\rm Int}P$.
\end{lem}
\begin{proof}
Let $u = \vec{\frak v}_T(\frak y)$. If $u \notin
\text{\rm Int}P$ then there exists $j$ such that
$\ell_j(u) \le 0$.
There exists nonzero $\mathcal P \in \text{\rm Jac}(\frak{PO}_{\frak b})$ such that $y_i \mathcal P = \frak y_i \mathcal P$ for
$i = 1,\dots,n$.
We have
\begin{equation}\label{provePzero}
(z_j)^k \mathcal P = T^{-k\lambda_j} y_1^{kv_{j,1}}\cdots y_n^{kv_{j,n}}\mathcal P
= T^{-k\lambda_j}\frak y_1^{kv_{j,1}}\cdots \frak y_n^{kv_{j,n}}\mathcal P
\end{equation}
Since
$
\lim_{k\to \infty}z_j^k = 0
$
in $\Lambda\langle\!\langle y,y^{-1}\rangle\!\rangle_0^{\overset{\circ}P}$
and
$$
\lim_{k\to \infty}\frak v_T((T^{-k\lambda_j}\frak y_1^{kv_{j,1}}\cdots \frak y_n^{kv_{j,n}})^{-1})
= -\lim_{k\to\infty} k\ell_j(u) \ge 0,
$$
it follows that
$$
(z_j)^k(T^{-k\lambda_j}\frak y_1^{kv_{j,1}}\cdots \frak y_n^{kv_{j,n}})^{-1}
$$
converges to $0$ in $\Lambda\langle\!\langle y,y^{-1}\rangle\!\rangle_0^{\overset{\circ }P}$. Therefore (\ref{provePzero}) implies
$\mathcal P = 0$. This is a contradiction.
\end{proof}
\begin{lem}\label{criequaleigen}
$\frak Y = \text{\rm Crit}(\frak{PO}_{\frak b})$.
\end{lem}
\begin{proof}
Let $\frak y \in \frak Y$. Then there exists
nonzero $\mathcal P \in \text{\rm Jac}(\frak{PO}_{\frak b};\frak y)$
such that $y_i \mathcal P = \frak y_i \mathcal P$ for
$i = 1,\dots,n$. It follows that $y_i \mapsto \frak y_i$ induces a
non-zero homomorphism $\Phi_{\frak y}: \text{\rm Jac}(\frak{PO}_{\frak b}) \to \Lambda$.
Therefore  (\ref{criticaleq}) is satisfied. (\ref{eq:inthemp}) follows from Lemma \ref{inthemoment}.
\par
We next assume $\frak y \in \text{\rm Crit}(\frak{PO}_{\frak b})$.
Then, by Lemma \ref{pointP} and the definition, $y_i \mapsto \frak
y_i$ induces a homomorphism $\Phi_{\frak y}: \text{\rm Jac}(\frak{PO}_{\frak b}) \to \Lambda$.
For given $\frak y' = (\frak y'_1,\dots,\frak y'_n) \in \frak Y$, let $1_{\frak y'}$ be
the unit of $\text{\rm Jac}(\frak{PO}_{\frak b};\frak y')$. By definition, we
have $(y_i - \frak y'_i)^N 1_{\frak y'}=0$ for sufficiently large $N$
and so we have
$0 = \Phi_{\frak y}((y_i - \frak y'_i)^N1_{\frak y'})
= (\frak y_i - \frak y_i')^N\Phi_{\frak y}(1_{\frak y'})$. Therefore if $\frak y \ne \frak
y'$, we obtain $\Phi_{\frak y}(1_{\frak y'}) = 0$. Therefore
the factorization (\ref{decompJac}) implies
$\text{\rm Crit}(\frak{PO}_{\frak b}) \subset \frak Y$
for otherwise we would have $\Phi_{\frak y} \equiv 0$.  
This contradicts to the fact that $\Phi_{\frak y}([1]) = 1$. The proof of Lemma
\ref{criequaleigen} is complete.
\end{proof}
The proof of  Proposition \ref{Morsesplit}.1 is complete 
by the next lemma.
\begin{lem}\label{localityof}
$\text{\rm Jac}(\frak{PO}_{\frak b};\frak y)$ is a local ring.
\end{lem}
\begin{proof}
For any homomorphism $\text{\rm Jac}(\frak{PO}_{\frak b};\frak y) \to \Lambda$,
the element $y_i1_{\frak y}$ must be sent to $\frak y_i$. The lemma then follows
form the fact that $\text{\rm Jac}(\frak{PO}_{\frak b};\frak y)$ is generated
by $y_i1_{\frak y}$ as a $\Lambda$ algebra.
\end{proof}
Finally we prove Proposition \ref{Morsesplit}.2 and 3.
We use the assumption that
$\frak{PO}_{\frak b}$ is a Morse function.
Lemma \ref{criequaleigen} and (\ref{decompJac}) imply that
$\prod_{\frak y \in \text{\rm Crit}(\frak{PO}_{\frak b})}
\Phi_{\frak y}$ is surjective. To prove its injectivity it suffices
to show the following.
\begin{lem}\label{Morsecriterion}
$\frak y$ is a non-degenerate critical point of $\frak{PO}_{\frak b}$ if and only if
$$
\dim_{\Lambda} \text{\rm Jac}(\frak{PO}_{\frak b};\frak y) = 1.
$$
\end{lem}
\begin{proof}
First assume $ \dim_{\Lambda} \text{\rm Jac}(\frak{PO}_{\frak b};\frak y)> 1. $
By an elementary linear algebra we can find $\mathcal P \in
\text{\rm Jac}(\frak{PO}_{\frak b};\frak y)$ and $i_0$ with the following
properties:
\begin{enumerate}
\item $(y_{i_0} - \frak y_{i_0}) \mathcal P \ne 0$.
\item $(y_j - \frak y_j) \mathcal P  = 0$ for $j \ne i_0$.
\item $(y_{i_0} - \frak y_{i_0})^2 \mathcal P = 0$.
\end{enumerate}
We will show by contradiction that $\text{\bf e}_{i_0} =
(0,\dots,0,\overset{\text{$i_0$-th}}{1},0,\dots,0)$ cannot be
contained in the image of the Hessian matrix $\text{\rm Hess}_{\frak
y}\frak{PO}_{\frak b}$.

Suppose to the contrary that there exists a vector $\vec c = (c_1,
\dots,c_n)$ with $c_i \in \Lambda$ such that
\begin{equation}\label{deltaijintheimage}
\text{\bf e}_{i_0} = \left(\frac{\partial^2 \frak{PO}_{\frak
b}}{\partial y_j\partial y_k}(\frak y)\right)\vec c.
\end{equation}
In other words, we have
$$
dy_{i_0}(\frak y) = d\left(\sum_{k=1}^n c_k\frac{\partial \frak{PO}_{\frak b}}{\partial
y_k}\right)(\frak y).
$$
Integrating this, we obtain
\begin{equation}\label{yisekibun}
y_{i_0} - \frak y_{i_0} = \tilde{\mathcal Q} + \sum_{k=1}^n c_k \frac{\partial
\frak{PO}_{\frak b}}{\partial y_k}(y)
\end{equation}
for $\tilde{\mathcal Q} \in \Lambda\langle\!\langle y,y^{-1}\rangle\!\rangle $ such that
\begin{equation}\label{propertiesQ}
\tilde{\mathcal Q}(\frak y) = 0, \quad  d\tilde{\mathcal Q}(\frak y) = 0.
\end{equation}
This implies $\tilde{\mathcal Q} = \sum_{j_1,j_2=1}^n (y_{j_1} - \frak y_{j_1})(y_{j_2} - \frak y_{j_2}) \tilde{\mathcal R}_{j_1,j_2}$ for
some $\tilde{\mathcal R}_{j_1,j_2} \in \Lambda\langle\!\langle y,y^{-1}\rangle\!\rangle^{\overset{\circ}P}$. Therefore Properties 2,
3 of $\mathcal P$ above imply $\mathcal P\mathcal Q = 0$ in
$\text{\rm Jac}(\frak{PO}_{\frak b};\frak y)$. 
(Here $\mathcal Q$ is the equivalence class of 
$\tilde{\mathcal Q}$ in $\text{\rm Jac}(\frak{PO}_{\frak b};\frak y)$.) Thanks to (\ref{yisekibun}), this in turn implies $\mathcal
P(y_{i_0} - \frak y_{i_0}) = 0$. This contradicts to
Property 1 above and so the Hessian matrix $\text{\rm Hess}_{\frak
y}\frak{PO}_{\frak b}$ cannot be non-degenerate.
\par
We next prove the converse. Suppose $ \dim_{\Lambda}
\text{\rm Jac}(\frak{PO}_{\frak b};\frak y) = 1$ for all $\frak y \in \text{\rm Crit}(\frak{PO}_{\frak b})$. 
Using the finite dimensionality of
$\text{\rm Jac}(\frak{PO}_{\frak b})$ over $\Lambda$ and Lemma
\ref{localityof}, there exists a sufficiently large $N$ such that
whenever $\tilde{\mathcal Q}(\frak y') = 0$ for $\tilde{\mathcal Q} \in \Lambda\langle\!\langle y,y^{-1}\rangle\!\rangle^{\overset{\circ}P}$, 
the equality ${\mathcal Q}^N = 0$ holds in 
$\text{\rm Jac}_{\frak b}(\frak{PO};\frak y')$. 
(Here ${\mathcal Q}$ is the equivalence class of $\tilde{\mathcal Q}$ in $\text{\rm Jac}_{\frak b}(\frak{PO};\frak y')$.)
\par
Using this fact, we find an element
$\tilde{\mathcal P} \in \Lambda\langle\!\langle y,y^{-1}\rangle\!\rangle^{\overset{\circ}P}_0$ such that 
its equivalence class $\mathcal P$ is zero in
$\text{\rm Jac}_{\frak b}(\frak{PO};\frak y')$ for all $\frak y' \ne \frak y$
and $\tilde{\mathcal P}(\frak y) \ne 0$. By the hypothesis 
$\dim_{\Lambda}
\text{\rm Jac}(\frak{PO}_{\frak b};\frak y) = 1$, we conclude $(y_i - \frak
y_i)\mathcal P(\frak y) = 0$ for $i = 1, \dots, n$ in
$\text{\rm Jac}(\frak{PO}_{\frak b};\frak y)$.
\par
Therefore it follows from (\ref{decompJac}) that for each $i$ there
exists some vector $\vec{\mathcal G}_i = (\mathcal G_{i,1}, \dots,
\mathcal G_{i,n})$ with $\mathcal G_{i,k} \in 
\Lambda\langle\!\langle y,y^{-1}\rangle\!\rangle^{\overset{\circ}P}$
such that
\begin{equation}\label{PmodJacobi}
(y_i - \frak y_i)\tilde{\mathcal P} = \sum_{k=1}^n \mathcal G_{i,k}
\frac{\partial \frak{PO}_{\frak b}}{\partial y_k}.
\end{equation}
Here $\tilde{\mathcal P}\in \Lambda\langle\!\langle y,y^{-1}\rangle\!\rangle_0^{P_{\e}}$ is an element whose equivalence class in 
$\text{\rm Jac}(\frak{PO}_{\frak b};\frak y)$ is $\mathcal P$ (see \eqref{eq2211}).
By substituting $y = \frak y$ into the $y_i$-derivative 
and $y_j$-derivative of
(\ref{PmodJacobi}), we obtain
$$
\aligned
0\ne\tilde{\mathcal P}(\frak y) &= \sum_{k=1}^n c_k \frac{\partial^2 \frak{PO}_{\frak b}}{\partial y_k\partial y_i}(\frak y), \\
0 &= \sum_{k=1}^n c_k \frac{\partial^2 \frak{PO}_{\frak b}}{\partial
y_k\partial y_j}(\frak y) \qquad \text{$j\ne i$}
\endaligned$$
for $c_k = \mathcal G_{i,k}(\frak y)$. This implies that $\text{\bf
e}_i$ is in the image of the Hessian matrix at $\frak y$. Since $i$
is arbitrary for $1, \dots, n$, the Hessian at $\frak y$ is
nondegenerate, which finishes the proof.
\end{proof}
Now the proof of Proposition \ref{Morsesplit} is complete.
\end{proof}
\begin{prop}\label{extraJac}
We have
\begin{equation}\label{Jacskosisemai1}
\text{\rm Jac}(\frak{PO}_{\frak b}) \cong
\frac{\Lambda\langle\!\langle y,y^{-1}\rangle\!\rangle_0^{P_{\e}}}{\left(y_i\frac{\partial \frak{PO}_{\frak b}}{\partial y_i}: i=1,\dots,n\right)},
\end{equation}
for sufficiently small $\epsilon > 0$.
If $\frak b \in \mathcal A(\Lambda_+)$ then we have
\begin{equation}\label{Jacskosisemai2}
\text{\rm Jac}(\frak{PO}_{\frak b}) \cong
\frac{\Lambda\langle\!\langle y,y^{-1}\rangle\!\rangle_0^{P}}{\left(y_i\frac{\partial \frak{PO}_
{\frak b}}{\partial y_i}: i=1,\dots,n\right)}.
\end{equation}
\end{prop}
\begin{rem}
We remark that $\Lambda\langle\!\langle y,y^{-1}\rangle\!\rangle_0^{P}$
is a $\Lambda_0$ module and is not a $\Lambda$ vector space.
In particular the right hand side of (\ref{Jacskosisemai1}) does not split as 
in Proposition \ref{Morsesplit}.
\end{rem}
\begin{proof}
There is a continuous homomorphism
\begin{equation}\label{toPemap}
\text{\rm Jac}(\frak{PO}_{\frak b})
\to
\frac{\Lambda\langle\!\langle y,y^{-1}\rangle\!\rangle_0^{P_{\e}}}{\text{\rm Clos}_{d^{P_{\epsilon}}}\left(y_i\frac{\partial \frak{PO}_{\frak b}}{\partial y_i}: i=1,\dots,n\right)},
\end{equation}
where the closure $\text{\rm Clos}_{d^{P_{\epsilon}}}$ is taken with respect to the
metric $d^{P_{\epsilon}}$.
Surjectivity of (\ref{toPemap}) follows from
Lemma \ref{polyapprox}.1.
(We may apply this lemma to $P=P_{\epsilon}$.)
To prove the injectivity of
(\ref{toPemap})  it suffices to prove that it is injective after
taking $\otimes_{\Lambda_0} \Lambda$. This is because
$\text{\rm Jac}(\frak{PO}_{\frak b})$ is torsion-free by Theorem \ref{Mirmain}.1.
We use Proposition \ref{Morsesplit}.1. Then
the factor $\text{\rm Jac}(\frak{PO}_{\frak b};\frak y)$
injects to the right hand side of (\ref{toPemap}) since
$\vec{\frak v}_T(\frak y) \in P_{\epsilon}$.
Therefore (\ref{toPemap}) is injective.
\par
We can prove
\begin{equation}\label{eq2211}
\frac{\Lambda\langle\!\langle y,y^{-1}\rangle\!\rangle_0^{P_{\e}}}{\left(y_i\frac{\partial \frak{PO}_{\frak b}}{\partial y_i}: i=1,\dots,n\right)}
\cong
\frac{\Lambda\langle\!\langle y,y^{-1}\rangle\!\rangle_0^{P_{\e}}}{\text{\rm Clos}
_{d_{P^{\epsilon}}}\left(y_i\frac{\partial \frak{PO}_{\frak b}}{\partial y_i}: i=1,\dots,n\right)}
\end{equation}
in the same way as in Lemma \ref{polyapprox}.3.
This implies (\ref{Jacskosisemai1}).
The proof of  (\ref{Jacskosisemai2}) is similar.\end{proof}
\begin{rem}
We used Lemma \ref{polyapprox} in the proof of Proposition \ref{extraJac}.
We do not use Proposition \ref{extraJac} during the proof of
Lemma \ref{polyapprox}.
\end{rem}
\par
\section{Operator $\frak q$: review}
\label{sec:frakqreview}
We first review the construction of operator $\frak q$.\index{operator $\frak q$}
It is described in \cite[Section 7.4]{fooobook2} for the general case and
in \cite[Section 6]{fooo09} for the toric case.
\par
Let
$\mathcal M_{k+1;\ell}^{\text{\rm main}}(\beta)$ be the
compactified moduli space of the genus zero bordered holomorphic maps
$u : (\Sigma,\partial\Sigma) \to (X,L)$,
in class $\beta  \in H_2(X,L(u);\Z)$  with $k+1$ boundary marked points and $\ell$ interior
marked points.
(We require the order of $k+1$ boundary marked points respects
the counter-clockwise cyclic order of the boundary of the
Riemann surface.)
We have the evaluation maps
$$
\text{\rm ev}: \mathcal M_{k+1;\ell}^{\text{\rm main}}(\beta)
\to X^{\ell} \times L(u)^{k+1}
$$
where we put
$$
\text{\rm ev} = (\text{\rm ev}^{\text{int}},\text{\rm ev})
= (\text{\rm ev}^{\text{int}}_1,\dots,\text{\rm ev}^{\text{int}}_{\ell};\text{\rm ev}_0,\dots,\text{\rm ev}_k).
$$
We denote the $\ell$-tuple of $T^n$-invariant cycles by
$$
\text{\bf p} = (\text{\bf p}(1),\dots,\text{\bf p}(\ell)).
$$
Here $\text{\bf p}(i)$ is a $T^n$-invariant cycle in
$\bigoplus_{k>0}\mathcal A^{k}(\Z)$. 
We note that the elements of $\mathcal A^{0}(\Z)$ change
the potential function only by constant and does not
affect the Jacobian ring. 
\par
We define the fiber product\index{$\mathcal M_{k+1;\ell}^{\text{\rm main}}(\beta;\text{\bf p})$}
\begin{equation}
\mathcal M_{k+1;\ell}^{\text{\rm main}}(\beta;\text{\bf p})
= \mathcal M_{k+1;\ell}^{\text{\rm main}}(\beta)
{}_{\text{\rm ev}^{\rm int}}\times_{X^{\ell}}
(\text{\bf p}(1) \times \dots \times \text{\bf p}(\ell)).
\end{equation}
It has a $T^n$-equivariant Kuranishi structure and the evaluation map
\begin{equation}\label{ev0afterfiber}
\text{\rm ev}_0: \mathcal M_{k+1;\ell}^{\text{\rm main}}(\beta;\text{\bf p})
\to L(u)
\end{equation}
is $T^n$-equivariant.
Since the $T^n$ action on $L(u)$ is free,
it follows from $T^n$-equivariance of (\ref{ev0afterfiber}) that
the $T^n$ action on each of the Kuranishi neighborhood of
$\mathcal M_{k+1;\ell}^{\text{\rm main}}(\beta;\text{\bf p})$ is
free.
We use this fact to take  an appropriate multisection of
$\mathcal M_{k+1;\ell}^{\text{\rm main}}(\beta;\text{\bf p})$.
\par
To describe the properties of this multisection, we need to review some notations which
we introduced in  \cite[Section 6]{fooo09}. We denote the set of
shuffles of $\ell$ elements by\index{$\text{\rm Shuff}(\ell)$}
\begin{equation}\label{shuff}
\text{\rm Shuff}(\ell) = \{ (\mathbb L_1,\mathbb L_2) \mid \mathbb
L_1 \cup \mathbb L_2 = \{1,\dots,\ell\}, \mathbb L_1 \cap
\mathbb L_2 = \emptyset \}.
\end{equation}
Let $\underline B = \{1,\dots,B\}$ and $Map(\ell,\underline B)$
the set of all maps $\{1,\dots,\ell\} \to \underline B$.\index{$Map(\ell,\underline B)$}
We will define a map\index{$\text{Split}$}
\begin{equation}\label{split}
\text{Split}: \text{\rm Shuff}(\ell) \times Map(\ell,\underline B)
\longrightarrow \bigcup_{\ell_1+\ell_2 = \ell} Map(\ell_1,\underline
B) \times Map(\ell_2,\underline B)
\end{equation}
as follows: Let $\text{\bf p} \in Map(\ell,\underline B)$ and
$(\mathbb L_1,\mathbb L_2) \in \text{\rm Shuff}(\ell)$. We put
$\ell_j = \# (\mathbb L_j)$ and let $\frak i_j: \{1,\dots,\ell_j\}
\cong \mathbb L_j$ be the order-preserving bijection. We consider
the map $\text{\bf p}_j: \{1,\dots,\ell_j\} \to \underline B$
defined by $ \text{\bf p}_j(i) = \text{\bf p}(\frak i_j(i)) $, and
set
$$
\text{Split}((\mathbb L_1,\mathbb L_2),\text{\bf p})
= (\text{Split}((\mathbb L_1,\mathbb L_2),\text{\bf p})_1,
\text{Split}((\mathbb L_1,\mathbb L_2),\text{\bf p})_2): = (\text{\bf
p}_1,\text{\bf p}_2).
$$
\par
We now define a corresponding gluing map
\begin{equation}\label{gluing}
\aligned
& \text{Glue}^{(\mathbb L_1,\mathbb L_2),\text{\bf
p}}_{\ell_1,\ell_2;k_1,k_2;i;\beta_1,\beta_2} : 
\\
&\quad \mathcal
M^{\text{\rm main}}_{k_1+1;\ell_1}(\beta_1;\text{\bf p}_1)
{}_{\text{\rm ev}_0}\times_{\text{\rm ev}_i} \mathcal M^{\text{\rm
main}}_{k_2+1;\ell_2}(\beta_2;\text{\bf p}_2) \to
\mathcal M^{\text{\rm main}}_{k+1;\ell}(\beta;\text{\bf p})
\endaligned
\end{equation}
below.
Here $k = k_1 + k_2 -1$, $\ell = \ell_1 + \ell_2$,
$\beta = \beta_1 + \beta_2$, and $i = 1,\dots,k_2$.
Put
$$
\mathbb S_j = ((\Sigma_{(j)},\varphi_{(j)},\{z^+_{i,{(j)}}\},
\{z_{i,{(j)}}\}) \in \mathcal M^{\text{\rm
main}}_{k_j+1;\ell_j}(\beta_j;\text{\bf p}_j)
$$
for $j = 1,2$.
\begin{rem}
Here and hereafter, boundary marked points are denoted by
$z_i$ and interior marked points are denoted either by
$z^+_i$ or $z^{\rm int}_i$.
\end{rem}
We glue $z_{0,(1)} \in \partial \Sigma_1$ with $z_{i,(2)}
\in
\partial \Sigma_2$ to obtain
$$
\Sigma = \Sigma_1 \#_{i} \Sigma_2.
$$
Suppose that $(\mathbb S_1,\mathbb S_2)$ is an element of the fiber
product in the left hand side of (\ref{gluing}). Namely we assume
$$
\varphi_{(1)}(z_{0,(1)}) = \varphi_{(2)}(z_{i,(2)}).
$$
This defines a holomorphic map
$$
\varphi = \varphi_{(1)} \#_i \varphi_{(2)}
: \Sigma \to X
$$
by putting $\varphi = \varphi_{(j)}$ on $\Sigma_j$.
\par
Let $m \in \mathbb L_j$. We define $c$ by $\frak i_j(c) = m$, where
$\frak i_j:
\{1,\dots,\ell_j\} \cong \mathbb L_j$ is the order-preserving
bijection. 
We put $z_m^{\text{int}} = z^{\text{int}}_{c;(j)}$ and call it 
the $m$-th interior marked point.
We define the boundary marked points
$(z_0,z_1,\dots,z_k)$ on $\Sigma$ by
$$
(z_0,z_1,\dots,z_k)
= (z_{0,(2)},\dots,z_{i-1,(2)},z_{1,(1)},\dots,
z_{k_1,(1)},z_{i+1,(2)},\dots,z_{k_2,(2)}).
$$
Now we put
$$
\mathbb S = ((\Sigma,\varphi,\{z^+_{i}\}, \{z_{i}\})
$$
and
$$
\text{Glue}^{(\mathbb L_1,\mathbb L_2),\text{\bf
p}}_{\ell_1,\ell_2;k_1,k_2;i;\beta_1,\beta_2} (\mathbb S_1,\mathbb
S_2) = \mathbb S.
$$
\par
The following lemma is proved in \cite[ Section 7.1]{fooobook2}.
\begin{lem}\label{boundaryno2}
We have the following identity as the spaces
with Kuranishi structure:
$$\aligned
\partial \mathcal M_{k+1;\ell}^{\text{\rm main}}(\beta;\text{\bf p})
= \bigcup_{k_1+k_2=k+1 \atop \beta_1+\beta_2=\beta}
\bigcup_{i=1}^{k_2}&\bigcup_{(\mathbb L_1,\mathbb L_2) \in \text{\rm Shuff}(\ell)}
\mathcal M_{k_1+1;\vert\mathbb L_1\vert}
^{\text{\rm main}}(\beta_1;(\text{\rm Split}((\mathbb L_1,\mathbb L_2),\text{\bf p})_1) \\
& {}_{\text{\rm ev}_0}\times_{\text{\rm ev}_i}
\mathcal M_{k_2+1;\vert\mathbb L_2\vert}
^{\text{\rm main}}(\beta_2;(\text{\rm Split}((\mathbb L_1,\mathbb L_2),\text{\bf p})_2).
\endaligned$$
\end{lem}
Let $\sigma: \{1,\dots,\ell\} \to \{1,\dots,\ell\}$ be an
element of $\frak S_{\ell}$, the symmetric group of order $\ell!$.
We have the induced map
$$
\sigma_{*}: \mathcal M_{k+1;\ell}^{\text{\rm main}}(\beta)
\to \mathcal M_{k+1;\ell}^{\text{\rm main}}(\beta)
$$
which permutates the interior marked points.
We define a new $\ell$-tuple $\text{\bf p}_{\sigma}$ by setting
$
\text{\bf p}_{\sigma}(i) = \text{\bf p}({\sigma}(i)).
$
Then $\sigma_{*}$ induces a map
\begin{equation}\label{interiorsyme}
\sigma_{*}: \mathcal M_{k+1;\ell}^{\text{\rm main}}(\beta;\text{\bf p})
\to \mathcal M_{k+1;\ell}^{\text{\rm main}}(\beta;\text{\bf p}_{\sigma}).
\end{equation}
\par
We now consider the following condition for a system of multisections $\frak q$.
(Here we use the same letter $\frak q$ as the operator $\frak q$ to
indicate that this multisection is used to define the operator
$\frak q$.)
We write $\mathcal M_{k+1;\ell}^{\text{\rm main}}(\beta;\text{\bf p})^{\frak q}$ for the
zero set of this multisection.
\begin{conds}\label{perturbforq}
\begin{enumerate}\label{condmultisq}
\item The multisection $\frak q$ is $T^n$-equivariant   and transversal to $0$.
\item The evaluation map
$ \text{\rm ev}_0: \mathcal M_{k+1;\ell}^{\text{\rm main}}(\beta;\text{\bf p})^{\frak q}
\to L(u)$ is a submersion. Here
the left hand side is the zero set of our multisection $\frak q$.
\item It is compatible with other multisections given at the boundary
$$
\mathcal M_{k_1+1;\ell_1}^{\text{\rm main}}(\beta_1;\text{\bf p}_1)^{\frak q}
{}_{\text{\rm ev}_0}\times_{\text{\rm ev}_i}
\mathcal M_{k_2+1;\ell_2}^{\text{\rm main}}(\beta_2;\text{\bf p}_2)^{\frak q}
$$
where $k_1+k_2 = k+1$, $\beta_1 +\beta_2 = \beta$, $i \in \{1,\dots,k_2\}$
and
$$
\text{Split}((\mathbb L_1,\mathbb L_2),\text{\bf p})
= (\text{Split}((\mathbb L_1,\mathbb L_2),\text{\bf p})_1,
\text{Split}((\mathbb L_1,\mathbb L_2),\text{\bf p})_2): = (\text{\bf
p}_1,\text{\bf p}_2).
$$
\item
It is compatible with the forgetful map
$$
\frak{forget}: \mathcal M_{k+1;\ell}^{\text{\rm main}}(\beta;\text{\bf p})^{\frak q}
\to
\mathcal M_{1;\ell}^{\text{\rm main}}(\beta;\text{\bf p})^{\frak q}
$$
which forgets the $1$-st \dots $k$-th boundary marked points.
\item It is invariant under the action of the symmetric group which exchanges the
interior marked points and the factors of $\text{\bf p}$.
\end{enumerate}
\end{conds}
\begin{rem}
Condition \ref{condmultisq}.1 implies
that $ \text{\rm ev}_i: \mathcal M_{k+1;\ell}^{\text{\rm main}}(\beta;\text{\bf p})^{\frak q}
\to L(u)$ is a submersion for each $i = 0,\dots,k+1$.
(So Condition \ref{condmultisq}.2 is a consequence of Condition \ref{condmultisq}.1 actually.)
It is in general impossible to take our perturbation $\frak q$ so that
$$
(\text{\rm ev}_i,\text{\rm ev}_j): \mathcal M_{k+1;\ell}^{\text{\rm main}}(\beta;\text{\bf p})^{\frak q}
\to L(u)^2
$$
is a submersion as far as we consider a \emph{single} multisection.
\end{rem}
The following is proved in \cite{fooo09}.
\begin{lem}[\cite{fooo09} {\rm Lemma $6.5$.
See also Sections \ref {sec:cyclicKura} - \ref{sec:equimulticot} of
this paper.}]\label{lem65toricII}
There exists a system of multisections,
which satisfies Condition  \ref{condmultisq} above.
\end{lem}
\begin{rem}\label{rem235}
Actually we need to fix an energy level $E_0$ and restrict the
construction of the family of multisections for the
moduli space of maps of homology class $\beta$ with $\beta \cap \omega \le E_0$,
by the reason explained in \cite[Subsection 7.2.3]{fooobook2}.
(After then we can use homological algebra to extend construction of
 the operators to all $\beta$'s.)
\par
This process is discussed in 
\cite[Section 7.2]{fooobook2}.
(See also \cite[Remark 11.3]{fooo08}.)
In the current
$T^n$-equivariant case, we can perform the same construction
in a $T^n$-equivariant way.
We will not repeat this kind of remarks in other similar situations in the rest of this chapter.
\end{rem}
\begin{defn}\label{defn66}
Hereafter we call the multisection on
$\mathcal M^{\text{main}}_{k+1;\ell}(\beta;{\text{\bf p}})$
chosen in  \cite[Lemma 6.5]{fooo09} the {\it $\frak q$-multisection}\index{multisection (perturbation)!$\frak q$-multisection}. This is the multisection
we use to define the operator $\frak q$.
\end{defn}
\par
Let $C$ be a graded module.
We define its degree shift $C[1]$ by $C[1]^k = C^{k+1}$. We
put
\begin{equation}
B_kC = \underbrace{C\otimes \dots \otimes C}_{\text{$k$ times}}.
\end{equation}
Let
$$
\widehat BC = \widehat{{\bigoplus_{k=0}^{\infty}}}B_kC
$$
be the completed direct
sum of them.
We also consider a map $\Delta^{k-1}: BC \to (BC)^{\otimes k}$
$$
\Delta^{k-1} = (\Delta \otimes  \underbrace{id \otimes \dots
\otimes id}_{k-2}) \circ (\Delta \otimes  \underbrace{id \otimes
\dots \otimes id}_{k-3}) \circ \dots \circ \Delta.
$$
For an indecomposable element $\text{\bf x} \in BC$, it can be
expressed as
\begin{equation}\label{coproduct:decon}
\Delta^{k-1}(\text{\bf x}) = \sum_c \text{\bf x}^{k;1}_c \otimes
\dots \otimes \text{\bf x}^{k;k}_c
\end{equation}
where $c$ runs over some index set.
\par
Let $\frak S_k$ \index{$\frak S_k$} be the symmetric group of order $k!$.
It acts on $B_kC$ by
\begin{equation}
\sigma \cdot (x_1\otimes \dots \otimes x_k)
= (-1)^* x_{\sigma(1)} \otimes \dots \otimes x_{\sigma(k)}
\end{equation}
where $* = \sum_{i<j : \sigma(i) > \sigma(j)} \deg x_i \deg  x_j$.
\begin{defn}
$E_kC$ is the quotient of $B_kC$ by the
submodule generated by $\sigma\cdot \text{\bf x} - \text{\bf x}$,
where $\sigma \in \frak S_k$, $\text{\bf x} \in B_kC$.
\par
$\widehat EC$ is a completion of the direct product
$\bigoplus_{k=0}^{\infty}E_kC$.
\par
For an element $\text{\bf x} \in B_kC$
we denote its equivalence class in $E_kC$ by $[\text{\bf x}]$.
\end{defn}
\begin{rem}\label{quotientsub}
In \cite{fooobook, fooo09} the authors defined
$E_kC$ as the {\it subset} of elements of $B_kC$ which are {\it invariant} under the
$\frak S_k$ action. Here we take the {\it quotient} instead.
This difference is not essential since we are working on a ring containing $\Q$.
It is slightly more convenient to use the quotient
for the description of this paper.
\par
There are a few changes which occur accordingly.
Especially for an element $\frak b$ of $C$ we define\index{$e^{\frak b}$}
$$
e^{\frak b} = \sum \frac{1}{k!} \frak b^k \in \widehat EC.
$$
Note in case of $\widehat BC[1]$,
$
e^{\frak b} = \sum \frak b^k \in \widehat BC[1] 
$. 
In \cite{fooobook, fooo09}
we put $
e^{\frak b} = \sum \frak b^k
$ also in case of $\widehat EC$ 
because we defined $\widehat EC$ as the subset of $\widehat BC$ and the coalgebra structure 
is the restriction of $\Delta$ in 
\eqref{coproduct:decon}. 
However, when we define $\widehat EC$ as the quotient, 
the coalgebra structure on $\widehat EC$ is given by\index{$\Delta_{\rm shuff}$}
$$\aligned
&\Delta_{\rm shuff}([x_1 \otimes \dots \otimes x_k]) \\
&= \sum_{(\mathbb L_1,\mathbb L_2) \in
\text{\rm Shuffle}(k)}
(-1)^*
[x_{\ell_1(1)} \otimes \dots \otimes x_{\ell_1(k_1)}]
\otimes
[x_{\ell_2(1)} \otimes \dots \otimes x_{\ell_2(k_2)}]
\endaligned$$
where $\mathbb L_1 = \{\ell_1(1),\dots,\ell_1(k_1)\}$,
$\mathbb L_2 = \{\ell_2(1),\dots,\ell_2(k_2)\}$ and
$$
* = \sum_{\ell_2(i) < \ell_1(j)} \deg x_{\ell_2(i)} \deg x_{\ell_1(j)}.
$$
Therefore the equality 
\begin{equation}\label{grouplikeelement}
\Delta_{\rm shuff}\left( \sum_k \frac{1}{k!} \frak b^k\right)
= \left( \sum_k \frac{1}{k!} \frak b^k\right) \otimes
\left( \sum_k \frac{1}{k!} \frak b^k\right)
\end{equation}
holds if $\deg x$ is even. 
Note that when $\widehat BC=\widehat BH(L(u);\Q)$, we 
use the coalgebra structure $\Delta$ in 
\eqref{coproduct:decon} so $e^b=\sum b^k$\index{$e^b$} for 
$b \in H(L(u);\Q)$, while 
we use $\Delta_{\rm shuff}$ and $e^{\frak b}=\sum \frac{\frak b^k}{k!}$ for the case $\widehat EC=\widehat E \mathcal A$ 
in this paper. See also Remark \ref{rem:6.11}.2. 
\end{rem}
We fix a $T^n$-invariant Riemannian metric on $L(u)$. We identify $H(L;\R)$
with the vector space of harmonic forms.
The following lemma can be easily checked and so its proof is omitted.
This characterization of harmonic forms for the torus
turns out to be useful for later discussions.
\begin{lem}\label{invariant-harmonic} Equip $T^n$ with an invariant Riemannian
metric. Then a differential form is harmonic with respect to this metric
if and only if it is $T^n$-invariant.
\end{lem}
We also identify
$H(L(u);\Lambda_0(\R))$ with $H(L(u);\R) \otimes_{\R} \Lambda_0(\R)$.
Hereafter, we write $\Lambda_0$ etc. in place of
$\Lambda_0(\R)$ etc..
\par
We take $h_1,\dots,h_k \in H(L(u);\Lambda_0)$ and identify them with
harmonic forms, that is $T^n$-invariant   forms in our case.
\begin{defn}\label{qedfdef}
We define
\begin{equation}\label{qbetadefform}
\frak q_{\ell;k;\beta}([\text{\bf p}] ; h_1\otimes\dots\otimes h_k)
= \text{\rm ev}_{0*} \text{\rm ev}^* (h_1\times\dots\times h_k),
\end{equation}
where
$$
(\text{\rm ev}_0,\text{\rm ev}): \mathcal M_{k+1;\ell}^{\text{\rm main}}(\beta;\text{\bf p})^{\frak q}
\to L(u) \times L(u)^k.
$$
\begin{rem}
We remark that when we need to use a multisection rather than a single-valued section, the perturbed moduli space
$\mathcal M_{k+1;\ell}^{\text{\rm main}}(\beta;\text{\bf p})^{\frak q}$ is neither a manifold nor an orbifold.
However integration along the fiber and pull-back of differential
forms on such space can be defined  in the way similar to the case of manifold.
(See \cite[Section 12]{fooo09}.)
\par
The same remark applies to several other places. We do not repeat this 
remark in those places.
\end{rem}
Condition $\ref{perturbforq}$ $(5)$ implies that the left hand side of (\ref{qbetadefform})
is independent of the representative $\text{\bf p}$ but depends only
on $[\text{\bf p}] \in \widehat E(\mathcal A[2])$.
We thus obtain:
$$
\frak q_{\ell;k;\beta} :
E_{\ell}(\mathcal A[2]) \otimes B_{k+1}(H(L(u);\Q)[1])  \to H(L(u);\Q)[1] .
$$
\end{defn}
\begin{rem}\label{rem:6.11}
\begin{enumerate}
\item
In our situation, all the forms $h_i$ and the (perturbed) moduli spaces
are $T^n$-invariant. Therefore
$\frak q_{\ell;k;\beta}([\text{\bf p}] ; h_1\otimes\dots\otimes h_k)
$ is also $T^n$-invariant.
So we need not take the procedure to go to the canonical model\index{canonical model} by using
the Green kernel to define $\frak q$.
\item
In  \cite[Section 3.8 (Definition 3.8.67 (3.8.68))]{fooobook} 
there is a factor $1/\ell!$ is in the right hand side.
There is the same factor in \cite[(6.10)]{fooo09}.
We do not put it here, since we take quotient to define $EC$ here.
(See Remark \ref{quotientsub}.)
\end{enumerate}
\end{rem}
\begin{conv}\label{ev0kakikataconv}
Hereafter we write
$$
(\text{\rm ev}_0)_*(\text{\rm ev}^* (h_1\times \dots \times h_k);\mathcal M_{k+1;\ell}^{\text{\rm main}}(\beta;\text{\bf p})^{\frak q} )
$$
for the right hand side of (\ref{qbetadefform}). We will apply the similar
notation in the general context. (See Section \ref{sec:equikuracot} for the usage of
the same notation in the general context.)
\end{conv}
We take $u \in \text{\rm Int}P$.
Let $b_0 = \sum x_i\text{\bf e}_i$ $x_i \in \C$ and associate a representation
$\rho : H_1(L(u);\Z) \to \C^*$ by $\rho(e_i) = e^{x_i}$.
\begin{defn}
Let $h_1,\dots,h_k \in H(L(u);\Lambda_0)$ be $T^n$-invariant forms. We put:
$$\aligned
&\frak q_{\ell;k}^{\rho}([\text{\bf p}] ; h_1\otimes\dots\otimes h_k) \\
&= \sum_{\beta\in H_2(X;L(u);\Z)}
T^{\beta \cap \omega/2\pi}
\rho(\partial \beta)
\frak q_{\ell;k;\beta}([\text{\bf p}] ; h_1\otimes\dots\otimes h_k).
\endaligned$$
\end{defn}
The main property of the operator $\frak q$ above is
 \cite[Theorem 3.8.32]{fooobook} or  \cite[Theorem 2.1]{fooo09}.
Namely we have:
\par
\begin{enumerate}
\item
For each $\beta$ and $\text{\bf x} \in B_k(H(L;\R)[1])$,
$\text{\bf y} \in E_k(H[2])$, we have the following:
\begin{subequations}\label{qproptall}
\begin{equation}\label{qmaineq0}
0 =
\sum_{\beta_1+\beta_2=\beta}\sum_{c_1,c_2}
(-1)^*
\frak q_{\beta_1}(\text{\bf y}^{2;1}_{c_1} ;
\text{\bf x}^{3;1}_{c_2} \otimes
\frak q_{\beta_2}(\text{\bf y}^{2;2}_{c_1} ; \text{\bf x}^{3;2}_{c_2})
\otimes \text{\bf x}^{3;3}_{c_2})
\end{equation}
where
$$
* = \deg'\text{\bf x}^{3;1}_{c_2} +
\deg'\text{\bf x}^{3;1}_{c_2} \deg \text{\bf y}^{2;2}_{c_1}
+\deg \text{\bf y}^{2;1}_{c_1}.
$$
In $(\ref{qmaineq0})$ and hereafter, we write $\frak q_{\beta}(\text{\bf y};\text{\bf x})$ in place of
$\frak q_{\beta;\ell,k}(\text{\bf y};\text{\bf x})$ if
$\text{\bf y} \in E_{\ell}(H[2])$, $\text{\bf x} \in B_{k}(H(L;\R)[1])$.
\item Let $\text{\bf e} = \text{\rm PD}([L])$ be the Poincar\'e dual to the fundamental
class of $L$. Let $\text{\bf x}_i \in B(H(L;\R)[1])$ and we put
$\text{\bf x} = \text{\bf x}_1 \otimes \text{\bf e} \otimes \text{\bf x}_2
\in B(H(L;R)[1])$. Then
\begin{equation}\label{unital0}
\frak q_{\beta}(\text{\bf y};\text{\bf x}) = 0
\end{equation}
except the following case:
\begin{equation}\label{unital20}
\frak q_{\beta_0}(1;\text{\bf e} \otimes x) =
(-1)^{\deg x}\frak q_{\beta_0}(1;x \otimes \text{\bf e}) = x,
\end{equation}
\end{subequations}
where $\beta_0 = 0 \in H_2(X,L;\Z)$ and $x \in H(L;\R)[1]
= B_1(H(L;\R)[1])$.
\end{enumerate}
\par
We can prove them by using  Condition $\ref{perturbforq}$.1-5.
\par
\begin{defn}\label{defmrho}
Let $\frak b \in H^{even}(X;\Lambda_0)$ and $\rho : H_1(L(u);\Z) \to \C^*$. We write
$\frak b = \frak b_0 + \frak b_2 + \frak b_{\text{\rm high}}$ where
$\frak b_0 \in H^{0}(X;\Lambda_0)$, $\frak b_2 \in H^{2}(X;\Lambda_0)$,
$\frak b_{\text{\rm high}} \in H^{2m}(X;\Lambda_0)$ $(m > 1)$.
We define
$
\frak m_{k;\beta}^{\rho} :
B_k(H(L;\R)[1]) \to H(L;\R)[1]
$
by
\begin{equation}\label{mcyclicdef0}
\aligned
&\frak m_{k}^{\frak b,\rho}(h_1,\dots,h_k) \\
&= \sum_{\beta}\exp(\frak b_2\cap \beta)\sum_{\ell=0}^{\infty}
\frac{T^{\beta \cap \omega/2\pi}}{\ell!}
\rho(\partial \beta)
\frak q_{\ell;k;\beta}(\frak b_{\text{\rm high}}^{\ell};h_1,\dots,h_k).
\endaligned\end{equation}
Let $b \in H^{1}(L(u);\Lambda_0)$. We define $b = b_0 + b_+$ where
$b_0 \in H^1(L(u);\C)$ and $b_+ \in H^1(L(u);\Lambda_+)$.
We define $\rho :  H_1(L(u);\Z) \to \C^*$ by $\gamma \mapsto e^{\gamma \cap b_0}$.
We define
\begin{equation}\label{mcyclicdef20}
\frak m_{k}^{\frak b,b}(h_1,\dots,h_k) =
\sum_{l_0,\dots,l_{k}}\frak m_{k+ \sum l_i;\beta}^{\frak b, \rho}(b_+^{l_0},h_1,b_+^{1},\dots,b_+^{l_{k-1}},
h_k,b_+^{l_k}).
\end{equation}
\end{defn}
We can use (\ref{qproptall}) to prove that $\frak m_{k}^{\frak b,b}$ defines a
unital filtered $A_{\infty}$ structure.
\par
In \cite[Definition 8.2]{fooo09}, we defined a chain complex
$
(\Omega(L(u))  \widehat{\otimes}  \Lambda_0; \delta_1^{(\frak b,b)})
$
on the de Rham complex of Lagrangian fiber $L(u)$
(tensored with the universal Novikov ring).
We can descend the (strict and gapped) filtered $A_{\infty}$
algebra on $\Omega(L(u))   \widehat{\otimes}  \Lambda_0$ to
the de Rham cohomology $H(L(u);\Lambda_0)$.
Especially we can define the boundary operator $\frak m_1^{\frak b,b}
= \pm \delta_{\text{\rm can}}^{\frak b,b}$ on $H(L(u);\Lambda_0)$ as above\footnote{Actually it is zero in
case Floer cohomology is isomorphic to the classical cohomology, that
is the case $b$ is a critical point of $\frak{PO}^u_{\frak b}$.}.
This is a general fact proven in  \cite[Section 5.4]{fooobook}.
In the rest of this section we will review its construction
and define the chain homotopy equivalence:
\begin{equation}\label{harmonicproj}
\Pi  :
(\Omega(L(u))\widehat{\otimes}\Lambda_0 ,\delta^{\frak b,b})
\to (H(L(u);\Lambda_0),\delta_{\text{\rm can}}^{\frak b,b}).
\end{equation}
\par
We first recall the definition of $\delta^{\frak b,b}$
in  \cite[Section 8]{fooo09}. \par
We use the correspondence by this perturbed moduli space to define
operators on the de Rham complex $\Omega(L(u))$ of $L(u)$.
Let $h \in \Omega^d(L(u))$ be a degree $d$ smooth differential
form on $L(u)$. Let
\begin{equation}\label{pullbackofh}
\text{\rm ev}_1^{*}(h) \in \Omega^{\deg h}(\mathcal M_{2,\ell}(\frak \beta;\frak b^{\otimes\ell}))
\end{equation}
be the pull-back of $h$. We define
\begin{equation}\label{pushout}
\delta_{\beta;\ell}^{\frak b}(h)
= (\text{\rm ev}_0)_{*}\text{\rm ev}_1^{*}(h)
= \frak q_{\ell;1;\beta}^{\frak b}(h)
\in \Omega^*(L(u))
\end{equation}
where $(\text{\rm ev}_0)_{*}$ is the integration along the fiber.
We note that $T^n$ equivariance of $\text{\rm ev}_0$ and
transitivity of the $T^n$ action on $L(u)$ imply that
$\text{\rm ev}_0$ is a submersion. Therefore integration along the fiber is well-defined
and defines a smooth form on $L$. The degree $*$ is given by
$$
* = \deg h - \mu(\beta) + 1
$$
where $\mu$ is the Maslov index.
\begin{rem}
We remark that $\frak b$ is not in $\mathcal A(\Z)$ but in $\mathcal A(\Lambda_0)$.
So, more precisely we need to define
$\delta_{\beta;\ell}^{\frak b} $ as follows.
We put
$$
\frak b = \sum_{i=1}^B w_i(\frak b) \text{\bf f}_i.
$$
For each multi-index $I = (i_1,\dots,i_{\ell})$,
we define $\text{\bf f}_I = (\text{\bf f}_{i_1},\dots,
\text{\bf f}_{i_{\ell}})$ and consider
$$
\text{\rm ev}^{I} = (\text{\rm ev}^{I}_0,\text{\rm ev}^{I}_1)
:
\mathcal M_{2;\ell}^{\text{\rm main}}(\beta;\text{\bf f}_I)
\to L(u) \times L(u).
$$
Now we put
\begin{equation}\label{pushout2}
\delta_{\beta;\ell}^{\frak b}(h)
= 
\sum_{I \in \{1,\dots,B\}^{\ell}}
w_{i_1}(\frak b)\cdots w_{i_{\ell}}(\frak b)
(\text{\rm ev}^{I}_0)_! (\text{\rm ev}^{I}_1)^*(h).
\end{equation}
Hereafter we write (\ref{pushout}) in place
of (\ref{pushout2}) whenever no confusion can occur.
\end{rem}
We next use $\delta_{\beta;\ell}^{\frak b}$ to define
$\delta^{\frak b,b}$. Let
$$
b = \sum_{i=1}^n \frak x_i \text{\bf e}_i \in H^1(L(u);\Lambda_0).
$$
We put
$
\frak y_i = e^{\frak x_i} \in \Lambda_0 \setminus \Lambda_+,
$
and define
$
\rho^{b}: H_1(L(u);\Z) \to \Lambda_0 \setminus \Lambda_+
$
by
\begin{equation}\label{eq:rhob}
\rho^{b}\left(\sum k_i \text{\bf e}^*_i\right)
= \frak y_1^{k_1} \cdots \frak y_n^{k_n}.
\end{equation}
\begin{defn}\label{bdryformdef}
We define
$
\delta^{\frak b,b}: \Omega(L(u)) \widehat{\otimes} \Lambda_0
\to \Omega(L(u)) \widehat{\otimes} \Lambda_0
$
by
$$
\delta^{\frak b,b} = \frak m_{1,0} + \sum_{\beta,\ell}
\frac{1}{\ell!} 
T^{\beta
\cap \omega/2\pi}
 \rho^{b}(\partial\beta)\delta_{\beta,\ell}^{\frak b},
$$
Here $\frak m_{1,0}$ is defined from the de Rham differential $d$ by
$\frak m_{1,0}(h) = (-1)^{n+\deg h +1}dh$.
See \cite[Remark 3.5.8]{fooobook}.
\end{defn}
By (\ref{pushout}) we easily find that if $h \in
\Omega(L(u))$ is $T^n$-invariant  then
$\delta_{\beta,\ell}^{\frak b,b}(h)
=  \rho^{b}(\partial\beta)\delta_{\beta,\ell}^{\frak b}(h)$ is also $T^n$-invariant .
Therefore
$$
\delta^{\frak b,b}(H(L(u);\Lambda_0))
\subset H(L(u);\Lambda_0).
$$
We define
\begin{equation}\label{canboundary}
\delta_{\text{\rm can}}^{\frak b,b}: H(L(u);\Lambda_0)
\to H(L(u);\Lambda_0)
\end{equation}
as the restriction of $\delta^{\frak b,b}$.
\begin{rem}
In fact, we have $\delta_{\text{\rm can}}^{\frak b,b} = \frak m_1^{\frak b,b}$
where the right hand side is in Definition \ref{defmrho}.
We will prove it in Section \ref{sec:deltaisthesame}.
\end{rem}
Let
$$
\Pi: \Omega(L(u)) \to H(L(u);\R)  \hookrightarrow \Omega(L(u))
$$
be the harmonic projection.

Now we have the following:
\begin{thm}\label{squizthm}
We have
\begin{equation}\label{d2iszerofirst}
\delta^{\frak b,b} \circ \delta^{\frak b,b} = 0.
\end{equation}
Moreover $\Pi$ is a chain homotopy equivalence.
\end{thm}
\begin{proof}
The proof is based on the study of boundary of
$\mathcal M_{k+1;\ell}^{\text{\rm main}}(\beta;\text{\bf p})$.
Using Lemma \ref{boundaryno2}, Condition \ref{perturbforq} and the definition, we can easily show the
following formula.
\begin{equation}
\frak m_{1,0}\circ \delta_{\beta;\ell}^{\frak b}
+\delta_{\beta;\ell}^{\frak b}\circ \frak m_{1,0}
= -\sum_{\beta_1+\beta_2=\beta}\sum_{\ell_1,\ell_2=0}^{\infty}
\frac{\ell!}{\ell_1!\ell_2!}\delta_{\beta_1;\ell_1}^{\frak b}
\circ \delta_{\beta_2;\ell_2}^{\frak b}.
\end{equation}
Therefore
$$\aligned
(\delta^{\frak b,b})^2
= \frak m_{1,0}^2 &+ \sum_{\beta}\sum_{\ell=0}^{\infty}
\frac{\rho^{b}(\partial\beta)}{\ell!}(\delta_{\beta;\ell}^{\frak b}
\circ \frak m_{1,0} + \frak m_{1,0}\circ \delta_{\beta;\ell}^{\frak b})\\
&+ \sum_{\beta_1+\beta_2=\beta}\sum_{\ell_1,\ell_2=0}^{\infty}
\frac{\rho^{b}(\partial\beta_1)\rho^{b}(\partial\beta_2)}{\ell_1!\ell_2!}
\delta_{\beta_1;\ell_1}^{\frak b}
\circ \delta_{\beta_2;\ell_2}^{\frak b}
= 0.
\endaligned$$
We next show that $\Pi$ is a chain map. We observe that
\begin{equation}\label{Piintegral}
\Pi(h) =
\int_{g \in T^n} g^*h dg,
\end{equation}
where $g : L \to L$ for $g\in T^n$ is obtained by $T^n$ action and
$dg$ is the normalized Haar measure. The formula (\ref{Piintegral})
and $T^n$ invariance of $\mathcal
M_{2;\ell}^{\text{main}}(\beta;\text{\bf f}_2)$ imply that
$$
\Pi \circ \delta^{\frak b,b} =
\delta_{\rm can}^{\frak b,b} \circ \Pi.
$$
Therefore $\Pi$ is a chain map. It is then easy to show that it is a chain homotopy equivalence.
\end{proof}
\par
\section{Well-definedness of Kodaira-Spencer map}
\label{sec:welldef}

In Section \ref{sec:statements} we used a choice of a section
$H(X;\Lambda_0) \to \mathcal A(\Lambda_0)$ of $\pi$ to give the
definition of the Kodaira-Spencer map
$$
{\frak{ks}}_{\frak b}: H(X;\Lambda_0) \to \text{\rm Jac}(\frak{PO}_{\frak b}).
$$
In this section, we prove that the definition of ${\frak{ks}}_{\frak
b}$ does not depend on the choice of the section.
\par
To make this statement precise, we need to set-up some notations.
Let $w_0,\dots,w_B$ be a coordinates of $\frak b$ with respect to
the basis $\{\frak f_i\}_{i =0}^B$, i.e., $\frak b = \sum_i w_i\frak
f_i$. We consider the potential function $
\frak{PO}(w_0,\dots,w_B;y_1,\dots,y_n) $ for $\frak b \in \mathcal
A(\Lambda_0)$ and define a map
$ \widetilde{\frak{ks}}_{\frak b}:
T_{\frak b} \mathcal A(\Lambda_0) \to
\Lambda\langle\!\langle y,y^{-1}\rangle\!\rangle^{\overset{\circ}P}_0$ by
$$
\widetilde{\frak{ks}}_{\frak b} \left(
\frac{\partial}{\partial w_i}\right)
=
\frac{\partial\frak{PO}_{\frak b}}{\partial w_i},
$$
and
$$
\overline{\frak w}_i\frac{\partial}{\partial \overline{\frak w}_i} = \frac{\partial}{\partial w_i}
$$
when $\deg w_i =2$.
Let $\pi: \mathcal A(\Lambda_0) \to H(X;\Lambda_0)$ be the
canonical map given in (\ref{Htohomology}). 
\begin{thm}\label{indepenceKS}
The following diagram commutes:
\begin{equation}
\begin{CD}
\mathcal A(\Lambda_0) @ > {\widetilde{\frak{ks}}_{\frak b}} >>
\Lambda\langle\!\langle y,y^{-1}\rangle\!\rangle_0^{\overset{\circ}P}  \\
@ V{\pi}VV @ VVV\\
H(X;\Lambda_0) @ > {{\frak{ks}}_{\frak b}} >> \text{\rm Jac}(\frak{PO}_{\frak
b}).
\end{CD}
\end{equation}
In particular, the definition of $\frak{ks}_{\frak b}$ is independent
of the section $H(X;\Lambda_0) \to \mathcal A(\Lambda_0)$ used in
Section $\ref{sec:statements}$.
\end{thm}
\par
The rest of this section is occupied with the proof of this
theorem.
\par
Let $Q = \sum_iw_i(Q) \text{\bf f}_i$ with $[Q] = 0$ in
$H(X;\Lambda_0)$. We will prove
\begin{equation}\label{conclusionsec4}
\widetilde{\frak{ks}}_{\frak b}(Q)
\in \left(y_i \frac{\partial\frak{PO}_{\frak b}}{\partial y_i}
: i=1,\dots,n\right).
\end{equation}
\begin{rem}
Note (\ref{conclusionsec4}) is slightly stronger than Theorem \ref{indepenceKS}
since we do not take closure in the right hand side.
\end{rem}
We may assume
$w_0(Q) = 0$. (In fact, $\mathcal A^0 \to H^0(L(u);\Z)$ is an
isomorphism.) We consider the case of $\beta$, $\frak b$ and $\text{\bf f}_i$ for which the moduli spaces
\begin{equation}\label{moduli71}
\mathcal M^{\text{\rm main}}_{1;\ell+1}(L(u),\beta;\frak b^{\otimes
\ell} \otimes \text{\bf f}_i)
\end{equation}
have dimension $n$. Then the evaluation map
$$
\text{\rm ev}_0: \mathcal M^{\text{\rm main}}_{1;\ell+1}(L(u),\beta;\frak
b^{\otimes \ell} \otimes \text{\bf f}_i) \to L(u)
$$
defines a homology class
\begin{equation}\label{deg0formninaru}
\text{\rm ev}_{0 *}(\mathcal M^{\text{\rm main}}_{1;\ell+1}(L(u),\beta;\frak b^{\otimes \ell} \otimes \text{\bf f}_i))
\in H_n(L(u);\Lambda_0) \cong H^0(L(u);\Lambda_0) = \Lambda_0.
\end{equation}
\begin{rem}
Actually $\frak b = \sum \frak b_i T^{\lambda_i}$ is not a space but contains $T$.
So (\ref{moduli71}) is not a space but a kind of formal power series
with space as coefficients. Namely
$$ \aligned
&\mathcal M^{\text{\rm main}}_{1;\ell+1}(L(u),\beta;\frak b^{\otimes
\ell} \otimes \text{\bf f}_i) \\
&= \sum_{i_1=1}^{\infty}\cdots\sum_{i_{\ell}=1}^{\infty}
\mathcal M^{\text{\rm main}}_{1;\ell+1}(L(u),\beta;\frak b_{i_1}\otimes
\dots \frak b_{i_{\ell}}\otimes \text{\bf f}_i) T^{\lambda_{i_1}+
\cdots + \lambda_{i_{\ell}}}.
\endaligned$$
Then the left hand side of (\ref{deg0formninaru}) has an obvious
sense. By an abuse of notation we use the expression such as
(\ref{moduli71}) frequently for the rest of this paper.
\end{rem}
\begin{lem}\label{KScalcu} Let $b = \sum_{i=1}^n x_i {\bf e}_i$ and
$y_i = e^{x_i}$. Then we have
$$\aligned
\widetilde{\frak{ks}}_{\frak b}(Q)(y) =
\sum_{i=1}^B\sum_{\ell=0}^{\infty}\sum_{\beta\in H_2(X;L(u);\Z)}
&\frac{w_i(Q)}{\ell!} T^{\omega\cap \beta/2\pi}
\rho^{b}(\partial \beta) \\
&\text{\rm ev}_{0 *}(\mathcal M^{\text{\rm main}}_{1;\ell+1}(L(u),\beta;\frak
b^{\otimes \ell} \otimes \text{\bf f}_i)),
\endaligned$$
where $Q = \sum w_i(Q) \text{\bf f}_i$.
Note the right hand side is an element of $\Lambda_0$ by the
isomorphism $(\ref{deg0formninaru})$.
\end{lem}
\begin{proof}
This follows by differentiating
\begin{equation}\label{POformula}
\frak{PO}(\frak b;y) = \sum_{\ell=0}^{\infty}\sum_{\beta\in
H_2(X;L(u);\Z)} \frac{1}{\ell!} T^{\omega\cap \beta/2\pi}
\rho^{b}(\partial \beta) \text{\rm ev}_{0 *}(\mathcal M^{\text{\rm
main}}_{1;\ell}(L(u),\beta;\frak b^{\otimes \ell}))
\end{equation}
with respect to $w_i$'s.
(\ref{POformula}) follows from
(\ref{PPformula}).
\end{proof}
For the convergence of our series, we
slightly rewrite Lemma \ref{KScalcu}.
(See the proof of  \cite[Proposition 11.4]{fooo09}.)
We put $\frak b = \sum_{a=1}^B
\frak b_a \text{\bf f}_a$. (We need not consider $\text{\bf f}_0$
because it does not affect the Jacobian ring.) We put
\begin{equation}
\frak b_{\rm high} = \sum_{a=m+1}^B \frak b_a \text{\bf f}_a,
\quad
\frak b_2 = \sum_{a=1}^m \frak b_a \text{\bf f}_a.
\end{equation}
Namely they are the degree $\ge 4$ part and the degree $2$ part, respectively.
\begin{lem}\label{deg2separate}
We have
$$\aligned
\widetilde{\frak{ks}}_{\frak b}(Q)(y) =
\sum_{i=1}^{B}\sum_{\ell=0}^{\infty}\sum_{\beta\in
H_2(X;L(u);\Z)} &\frac{w_i(Q)}{\ell!} T^{\omega\cap \beta/2\pi}\exp
(\frak b_2 \cap \beta)
\rho^{b}(\partial \beta) \\
&\text{\rm ev}_{0 *}(\mathcal M^{\text{\rm main}}_{1;\ell+1}(L(u),\beta;\frak b_{\rm high}^{\otimes \ell} \otimes \text{\bf f}_i)).
\endaligned$$
\end{lem}
\par
We are now ready to give the proof of Theorem \ref{indepenceKS}. We
split the proof into two parts. One is consideration of the case $Q \in
\bigoplus_{k\ge 4}\mathcal A^k(\Lambda_0)$ and the other is of the case
$Q \in \mathcal A^2(\Lambda_0)$.
\par
We start with the first case. We use our assumption and take a
singular chain $R$ of $X$ such that $\partial R = Q$. (Here we fix a
triangulation of $Q$ and regard it as a smooth singular chain.) We
emphasize it is {\it impossible} to choose such an $R$ that is
$T^n$-invariant.
\par
We now consider
$$
\mathcal M^{\text{\rm main}}_{k+1;\ell+1}(L(u),\beta;\text{\bf p} \otimes R)
= \mathcal M^{\text{\rm main}}_{k+1;\ell+1}(L(u),\beta) {}_{\text{\rm ev}^{\text{int}}} \times
(\text{\bf p}(1) \times \dots \times \text{\bf p}(\ell) \times R).
$$
Here $\text{\bf p}(1), \dots, \text{\bf p}(\ell)$ are
$T^n$-invariant  cycles. We remark that in \cite{fooo09} we never
considered chains on $X$ which are not $T^n$-invariant. Here we
however need to use the chains that are not $T^n$-invariant.
Because of this, there are several points we need to be careful in
the subsequent discussions henceforth.
\begin{lem}\label{splitKuranishilem}
There exists a Kuranishi structure with boundary on the moduli space $\mathcal
M^{\text{\rm main}}_{k+1;\ell+1}(L(u),\beta;\text{\bf p} \otimes
R)$. Its boundary is a disjoint union of the following three types
of fiber product. Here $\beta_1+\beta_2 = \beta$, $k_1+k_2=k+1$,
$i=1,\dots,k_2$ and $(\text{\bf p}_1,\text{\bf p}_2) = \text{\rm
Split}((\mathbb L_1,\mathbb L_2),\text{\bf p})$ for some $(\mathbb L_1,\mathbb L_2) \in
\text{\rm Shuff}(\ell)$.
\begin{equation}\label{splitKuranishi1}
\mathcal M^{\text{\rm main}}_{k_1+1;\vert\text{\bf p}_1\vert+1}(L(u),\beta_1;\text{\bf p}_1 \otimes R)
{}_{\text{\rm ev}_0}\times_{\text{\rm ev}_i} \mathcal M^{\text{\rm main}}_{k_2+1;\vert\text{\bf p}_2\vert}
(L(u),\beta_2;\text{\bf p}_2),
\end{equation}
\begin{equation}\label{splitKuranishi2}
\mathcal M^{\text{\rm main}}_{k_1+1;\vert\text{\bf p}_1\vert}(L(u),\beta_1;\text{\bf p}_1)
{}_{\text{\rm ev}_0}\times_{\text{\rm ev}_i} \mathcal M^{\text{\rm main}}_{k_2+1;\vert\text{\bf p}_2\vert+1}(L(u),\beta_2;\text{\bf p}_2 \otimes R),
\end{equation}
\begin{equation}\label{splitKuranishi3}
\mathcal M^{\text{\rm main}}_{k+1;\ell+1}(L(u),\beta;\text{\bf p} \otimes Q).
\end{equation}
Here and hereafter we write $\vert \text{\bf p}\vert = \ell$ if
$\text{\bf p} \in E_{\ell}\mathcal A(\Lambda_0)$.
\end{lem}
This follows from  \cite[Lemmata 6.4 and 6.5]{fooo09}.
\par
We next define a multisection (perturbation). We remark that the
multisection $\frak q$ on $\mathcal M^{\text{\rm
main}}_{k+1;\ell}(L(u),\beta;\text{\bf p})$ was already fixed in
Lemma \ref{lem65toricII} (=  \cite[Lemma 6.5]{fooo09}), which we use here. The multisection we have
chosen there is $T^n$-equivariant which we use to prove that $\text{\rm ev}_0:
\mathcal M^{\text{\rm main}}_{k+1;\ell}(L(u),\beta;\text{\bf
p})^{\frak q} \to L(u)$ is a submersion. 
\begin{rem}
Note $\mathcal
M^{\text{\rm main}}_{k+1;\ell}(L(u),\beta;\text{\bf p})^{\frak q} $
is a zero set of {\it multi-valued} perturbation. Although the zero
set does not define a smooth manifold, we can define the notion of
submersion from such a space to a manifold as follows. 
\par
Let $(V,E,\Gamma,\psi,s)$ be a chart of a good coordinate system 
of our Kuranishi structure. 
Namley $U=V/\Gamma$ is an orbifold, $E$ is an orbibundle on it, $s$ is a section of 
$E$ and $\psi : s^{-1}(0)/\Gamma \to 
\mathcal M^{\text{\rm main}}_{k+1;\ell}(L(u),\beta;\text{\bf
p})^{\frak q}$  is a homeomorphism onto an open set.
Our multisection $s_{\epsilon}$ on this chart in given as follows.
\par
For each $[x] \in U$ there exists an orbifold chart $U_x = V_x/\Gamma_x$
where $V_x$ is a smooth manifold on which a finite group $\Gamma_x$ 
acts effectively. The restriction of $E$ to $U_x$ is regarded as 
$(E_x \times V_x)/\Gamma_x$. The restriction of our multisection 
$s_{\epsilon}$ to $U_x$ is given by 
$$
s_{\epsilon} = (s_{\epsilon,1},\dots,s_{\epsilon,m}) :
V_x \to E_x^m
$$
that induces a $\Gamma_x$ equivariant map $V_x \to S^m(E_x)$,
where $S^m$ denotes the symmetric power.
We call $s_{\epsilon,i}$ the branch. Each branch is of $C^{\infty}$ class 
and is transversal to $0$.
\par
Let $y \in s_{\epsilon,i}^{-1}(0)$. The tangent space 
$T_ys_{\epsilon,i}^{-1}(0)$ is defined.
The map $\text{\rm ev}_0$ induces a smooth map 
$\tilde{\text{\rm ev}}_0 : V_x \to L(u)$.
We say $\text{\rm ev}_0$ is a submersion if 
$$
d_y\tilde{\text{\rm ev}}_0 : T_ys_{\epsilon,i}^{-1}(0) \to T_{\text{\rm ev}_0(y)}L(u)
$$
is a submersion for all $y$. 
(See also \cite{FO}, \cite{fooo06}.)
 \end{rem}
For the case of our moduli space $\mathcal M^{\text{\rm
main}}_{k+1;\ell}(L(u),\beta;\text{\bf p}\otimes R)$, we do not have
$T^n$-action on it. So we can no longer use this technique.
Therefore the evaluation map $\text{\rm ev}_0: \mathcal M^{\text{\rm
main}}_{k+1;\ell}(L(u),\beta;\text{\bf p}\otimes R) \to L(u)$ is not
necessarily submersive. (It may occur for example the case where
its dimension is strictly smaller than $n$.) On the other hand, we
are using de Rham cohomology here. So we need to define the
integration along the fiber of our evaluation map and we need our
evaluation map to be submersive. We can do this by using the notion
of continuous family of multisections. This notion was explained in
detail in \cite[ Section 12]{fooo09}.
\begin{lem}\label{perturbwellks}
There exists a continuous family of multisections on $\mathcal
M^{\text{\rm main}}_{k+1}(L(u),\beta;\text{\bf p}\otimes R)$ with
the following properties:
\begin{enumerate}
\item It is compatible with the decomposition into $(\ref{splitKuranishi1})$,
$(\ref{splitKuranishi2})$, $(\ref{splitKuranishi3})$ of its boundary.
\item
The restriction of the evaluation map $\text{\rm ev}_0$ to the zero set
of the family of multisections is a submersion to $L(u)$.
\item
It is compatible with the forgetful map
$$
\frak{forget}:
\mathcal M^{\text{\rm main}}_{k+1;\ell+1}(L(u),\beta;\text{\bf p}\otimes R)
\to \mathcal M^{\text{\rm main}}_{1;\ell+1}(L(u),\beta;\text{\bf p}\otimes R)
$$
of boundary marked points.
\end{enumerate}
\end{lem}
\begin{proof}
The proof is by an induction over the symplectic area $\beta\cap
\omega$ and the number of boundary marked points, using  
\cite[Lemma 12.19]{fooo09}.
\end{proof}
\begin{lem}\label{unitvanish}
We have an equality
\begin{equation}\label{ev0*}
\text{\rm ev}_{0*}\left(
\mathcal M^{\text{\rm main}}_{1;\vert\text{\bf p}_1\vert}(L(u),\beta_1;\text{\bf p}_1)
{}_{\text{\rm ev}_0}\times_{\text{\rm ev}_1} \mathcal M^{\text{\rm main}}_{2;\vert\text{\bf p}_2\vert+1}
(L(u),\beta_2;\text{\bf p}_2 \otimes R)\right) = 0
\end{equation}
as differential forms on $L(u)$.
\end{lem}
\begin{proof} By definition, it is easy to see that the differential form \eqref{ev0*} is
the integration along the fiber of the differential form
$$
\text{\rm ev}_{0*}\left(
\mathcal M^{\text{\rm main}}_{1;\vert\text{\bf p}_1\vert}(L(u),\beta_1;\text{\bf p}_1)\right)
$$
by the correspondence
$$
(\text{\rm ev}_0,\text{\rm ev}_1):\mathcal M^{\text{\rm main}}_{2;\vert\text{\bf p}_2\vert+1}(L(u),\beta_2;\text{\bf p}_2 \otimes R)
\to L \times L.
$$
On the other hand, it follows from  \cite[Corollary 6.6]{fooo09} that the dimension of $\mathcal M^{\text{\rm
main}}_{1;\vert\text{\bf p}_1\vert}(L(u),\beta_1;\text{\bf p}_1)$ is
$n$ or higher. Therefore integration along the fiber gives rise to
a function, that is a degree $0$ form on $L(u)$. (See (\ref{deg0formninaru}). Actually it is a constant because of $T^n$
invariance.)
\par
By Lemma \ref{perturbwellks}.3, the evaluation map $\text{\rm ev}_1$ factors through
$$
\mathcal M^{\text{\rm main}}_{2;\vert\text{\bf p}_2\vert+1}(L(u),\beta_2;\text{\bf p}_2 \otimes R)
\overset{\frak{forget}}{\longrightarrow}  \mathcal M^{\text{\rm main}}_{1;\vert\text{\bf p}_2\vert+1}(L(u),
\beta_2;\text{\bf p}_2 \otimes R) \overset{\text{\rm ev}_0}{\longrightarrow} L(u)
$$
and hence the fiber of $\text{\rm ev}_1:\mathcal M^{\text{\rm main}}_{2;\vert\text{\bf p}_2\vert+1}(L(u),\beta_2;\text{\bf p}_2 \otimes R)  \to L$ has
positive dimension and hence \eqref{ev0*} follows.
\end{proof}
\begin{rem}
We may rephrase the argument used in the above proof as follows:
Since $b$ is a weak bounding cochain,
$\frak m_{0,\beta}^{\frak b,b}(1) = c \text{\rm PD}[L]$ ($c\in \Lambda_0$) by definition.
On the other hand $ \text{\rm PD}[L]$ is a strict unit of our filtered $A_{\infty}$ algebra
in the de Rham model. Therefore
$\frak m_{1,\beta}^{\frak b,b}( \text{\rm PD}[L]) = 0$.
\par
We remark that if we had used singular cohomology instead of de Rham cohomology,
then $\text{\rm PD}[L]$ would become only a homotopy unit rather than a strict unit which would
require more technical and complicated arguments.
This is one of reasons why we use the de Rham cohomology model here.
\end{rem}
\begin{defn} We denote
\begin{equation}\label{frakXdef}
\aligned
\tilde{\frak X}^u(\frak b,b) =
\sum_{\ell=0}^{\infty}\sum_{\beta\in H_2(X;L(u);\Z)} &\frac{1}{\ell!} T^{\omega\cap \beta/2\pi}
\rho^{b}(\partial \beta) \exp(\frak b_2\cap \beta) \\
&\text{\rm ev}_{0 *}
(\mathcal M^{\text{\rm main}}_{1;\ell+1}(L(u),\beta;\frak b_{\rm high}^{\otimes \ell} \otimes R)).
\endaligned
\end{equation}
This is an element of $\Omega(L(u)) \widehat{\otimes} \Lambda_0$ that is
formally obtained by replacing $Q$ by $R$ in the formula in Lemma
\ref{deg2separate}. We also define
$$
{\frak X}^u(\frak b,b) = \Pi(\tilde{\frak X}^u(\frak b,b))
\in H(L(u);\Lambda_0).
$$
\end{defn}
We note that here we fix $u$. We fix $\frak b$ and $u$ and regard
$\tilde{\frak X}^u(\frak b,b)$, ${\frak X}^u(\frak b,b)$ as functions of $b$.
Let $2^{\{1,\dots,n\}}$ be the set of all
subsets of $\{1,\dots,n\}$. For $I=\{i_1,\dots,i_k\} \in
2^{\{1,\dots,n\}}$,
$i_j<i_{j+1}$, we put $\text{\bf e}_I = \text{\bf e}_{i_1} \wedge \dots
\wedge \text{\bf e}_{i_k}$. It forms a basis of $H(L(u);\Z)$.
We put
$b = \sum_{I\in 2^{\{1,\dots,n\}}} x_I(u)\text{\bf e}_I$, $y_i(u) = e^{x_i(u)}$ and
define
$$
\tilde{\frak X}_{\frak b}^u(y_1(u),\dots,y_n(u))
= \tilde{\frak X}^u(\frak b,b),
\quad
{\frak X}_{\frak b}^u(y_1(u),\dots,y_n(u))
= {\frak X}^u(\frak b,b).
$$
We fix a diffeomorphism $\psi_u: T^n \cong L(u)$ which induces an identification
$H^*(T^n) \cong H(L(u))$ and $\Omega(T^n) \cong \Omega (L(u))$.
We also recall that $y_i$, $y(u)_i$ are the coordinates of $H^1(T^n;\C)$, $H^1(L(u);\C)$
respectively and satisfies
$$
y_i(u) = T^{-u_i} y_i.
$$
(See  \cite[(3.10)]{fooo09} for the precise meaning of this relation.)
\par
We regard each of $\tilde{\frak X}_{\frak b}^u$,
${\frak X}_{\frak b}^u$ as a formal sum of the form
$$
\sum a_i T^{\lambda_i} y^{v_i}
$$
where $a_i \in \Omega(T^n)$ or $a_i \in H(T^n;\C)$ respectively.
The next proposition claims that they converge in
$d_{\overset{\circ}P}$ topology.
\begin{prop}\label{fraXconv}
\begin{enumerate}
\item We have
$$
\tilde{\frak X}_{\frak b}^u \in \Omega(T^n) \widehat{\otimes}
\Lambda\langle\!\langle y,y^{-1}\rangle\!\rangle_0^{\overset{\circ}P}
$$
and
$$
{\frak X}_{\frak b}^u \in H(T^n;\Lambda\langle\!\langle y,y^{-1}\rangle\!\rangle_0^{\overset{\circ}P} ).
$$
\item In case $\frak b \in \mathcal A(\Lambda_+)$ we can
take $\Lambda\langle\!\langle y,y^{-1}\rangle\!\rangle_0^{P}$ in place of 
$\Lambda\langle\!\langle y,y^{-1}\rangle\!\rangle_0^{\overset{\circ}P}$.
\end{enumerate}
\end{prop}
\begin{proof}
We remark that Gromov compactness immediately
implies that the series (\ref{frakXdef}) converges in $\frak v_T^u$-norm.
It implies
$$
\tilde{\frak X}_{\frak b}^u \in \Lambda_0\langle\!\langle y(u),y(u)^{-1}\rangle\!\rangle .
$$
To prove Proposition \ref{fraXconv}.1 (resp. Proposition \ref{fraXconv}.2) it suffices to show that
(\ref{frakXdef}) converges in $\frak v_T^{u'}$-norm for {\it any} $u' \in \text{\rm Int}P$ (resp. $u' \in P$).
We remark that $\tilde{\frak X}_{\frak b}^u$ {\it depends} on $u$.
(In a similar situation in  \cite[Section 7]{fooo09}
we proved that the potential function (with bulk) is
independent of $u$. In that situation, our moduli space is
$T^n$-invariant. However in the current situation,
the moduli space $\mathcal M^{\text{\rm main}}_{1;\ell+1}(L(u),\beta;\frak b^{\otimes \ell} \otimes R)$ is {\it not} $T^n$-invariant. This is because the chain $R$ is not $T^n$-invariant.
Nevertheless, we will still be able to prove Proposition \ref{fraXconv}.)
\par
We remark that there exists a biholomorphic action of $(\C \setminus \{0\})^n$
on $X$ that extends the given $T^n$ action.
Moreover for each $u,u' \in \text{Int} P$, there exists 
$g \in (\C \setminus \{0\})^n$ such that $g L(u) = L(u')$.
Therefore
$$
\mathcal M^{\text{\rm main}}_{1,\ell+1}(L(u'),\beta;g\frak b_{\text{\rm high}}^{\otimes \ell}\otimes gR)
\cong
\mathcal M^{\text{\rm main}}_{1,\ell+1}(L(u),\beta;\frak b_{\text{\rm high}}^{\otimes \ell}\otimes R).
$$
We apply Gromov compactness to the left hand side to show that 
(\ref{frakXdef}) converges in 
$\frak v_T^{u'}$-norm for any $u' \in \text{\rm Int}P$.
The proof of Proposition \ref{fraXconv}.1 is complete.
\begin{rem}
We thank the referee who pointed out that using this $(\C \setminus \{0\})^n$
we can simplify the proof of Proposition \ref{fraXconv}.1.
\end{rem}
We next prove Proposition \ref{fraXconv}.2.
\par
Suppose $\text{\rm ev}_{0 *}
(\mathcal M^{\text{\rm main}}_{1;\ell+1}(L(u),\beta;\frak b_{\rm high}^{\otimes \ell} \otimes R)) \ne 0$.
Then by  \cite[Theorem 11.1]{fooo08} we can write
\begin{equation}\label{whatisbeta}
\beta = \sum_{i=1}^m k_i(\beta)\beta_i + \alpha(\beta)
\end{equation}
so that $k_i(\beta) \ge 0$, $\sum k_i(\beta) > 0$, and that $\alpha(\beta)$
is realized by a sum of holomorphic spheres.
\par
Since $\partial \beta = \sum_i k_i(\beta) \partial \beta_i$, it follows
from \eqref{zjdef} that
$$
y(u)_1^{\partial \beta \cap \text{\bf e}_1}\cdots
y(u)_n^{\partial \beta \cap \text{\bf e}_n}
= T^{c(\beta)} z_1^{k_1(\beta)}\cdots z_m^{k_m(\beta)}
$$
for
\begin{equation}\label{cbetadecide}
c(\beta) = - \sum_{i=1}^m k_i(\beta)\ell_i(u).
\end{equation}
Combining (\ref{whatisbeta}) and \eqref{cbetadecide}, we obtain
\begin{equation}\label{betacapomega}
\beta \cap \omega/2\pi = \alpha(\beta) \cap \omega/2\pi - c(\beta).
\end{equation}
We write the contribution from $(\beta,\ell)$ in (\ref{frakXdef}) as:
\begin{equation}\label{tildefrakX}
\widetilde{\frak X}_{\frak b,\beta,\ell}^u =  \frac{1}{\ell!} T^{\alpha(\beta) \cap \omega/2\pi}
 z_1^{k_1(\beta)}\cdots z_m^{k_m(\beta)} a(\beta,\ell)
\end{equation}
where
\begin{equation}
a(\beta,\ell) =
\exp(\frak b_0\cap \beta) \text{\rm ev}_{0 *}
(\mathcal M^{\text{\rm main}}_{1;\ell+1}(L(u),\beta;\frak b_{\rm high}^{\otimes \ell} \otimes R)) \in \Lambda_0 \widehat{\otimes} \Omega(T^n).
\end{equation}

\begin{lem}\label{unboundedX1}
Suppose that there are infinitely many different
$(\beta_c,\ell_c)$'s $(c=1,2,\dots)$ with $a(\beta_c,\ell_c)\ne 0$.
If $\frak b \in \mathcal A(\Lambda_+)$, then $\frak v_T^P(\widetilde{\frak X}_{\frak b,\beta,\ell}^u)$
are unbounded.
\end{lem}
\begin{proof} Suppose to the contrary that $\frak v_T^P(\widetilde{\frak X}_{\frak b,\beta,\ell}^u)$
are bounded. Then since we have
$\frak v_T^P(\widetilde{\frak X}_{\frak b,\beta,\ell}^u) \ge \ell_c \frak v_T(\frak b_{\rm high})$,
it follows that
 $\ell_c$ are bounded. By the hypothesis that the numbers of pairs $(\beta_c,\ell_c)$
are assumed to be infinite, there must be infinitely many
different $\beta_c$'s. We have the (virtual) dimension
\begin{equation}\label{virdim}
n + \mu(\beta_c) - (\text{deg }\frak b_{\rm high} -2)\ell_c - (\text{deg }R -2) -2
\end{equation}
for $\mathcal M^{\text{\rm main}}_{1;\ell_c+1}(L(u),\beta_c;\frak b_{\rm high}^{\otimes \ell} \otimes R)$.
Therefore $a(\beta_c,\ell_c) \neq 0$ only when
\begin{equation}\label{0virdimn}
0 \leq n + \mu(\beta_c) - (\text{deg }\frak b_{\rm high} -2)\ell_c - (\text{deg }R -2) -2 \leq n.
\end{equation}
Therefore since $\ell_c$ is bounded, the Maslov indices of $\beta_c$'s must
be bounded.

From \eqref{whatisbeta} and $\mu(\beta_i) = 2$, we have
$$
\mu(\beta_c) = \sum_{i=1}^m 2 k_i(\beta_c) + 2c_1(\alpha(\beta_c)).
$$
We split our consideration into two cases:
\begin{enumerate}
\item If $\sum k_i(\beta_c)$ are unbounded, then the Chern numbers of $\alpha(\beta_c)$
are unbounded.
\item
If $\sum k_i(\beta_c)$ are bounded, then since there are infinitely many $\beta_c$'s
there must be infinitely many different $\alpha(\beta_c)$'s.
\end{enumerate}
We thus find that there are always infinitely many different $\alpha(\beta_c)$'s.
Therefore we can conclude from Gromov compactness that the values $\omega \cap \alpha(\beta_c)$
are unbounded. But by taking the $\frak v_T^P$-valuation of \eqref{tildefrakX}, we have
$$
\frak v_T^P(\widetilde{\frak X}_{\frak b,\beta,\ell}^u)
\ge \omega \cap \alpha(\beta_c)/2\pi,
$$
which gives rise to a contradiction. This finishes the proof.
\end{proof}
Lemma \ref{unboundedX1} implies that 
(\ref{frakXdef}) converges in 
$\frak v_T^{u'}$-norm for any $u' \in P$ if 
$\frak b \in \mathcal A(\Lambda_+)$.
The proof of Proposition \ref{fraXconv}.2 is complete.
\end{proof}
\begin{prop}\label{bdryandks}
We have
\begin{equation}
\widetilde{\frak{ks}}_{\frak b}(Q)\cdot \text{\rm PD}[L(u)]
= \delta_{\text{\rm can}}^{\frak b,b}({\frak X}_{\frak b}^u).
\end{equation}
\end{prop}
\begin{proof}
By Theorem \ref{squizthm} we have:
\begin{equation}\label{orio7.132tochu}
\delta_{\text{\rm can}}^{\frak b,b}({\frak X}_{\frak b}^u)
= \delta_{\text{\rm can}}^{\frak b,b}(\Pi(\widetilde{\frak X}_{\frak b}^u)) \\
= \Pi(\delta^{\frak b,b}(\widetilde{\frak X}_{\frak b}^u)).
\end{equation}
It follows from Lemmata \ref{splitKuranishilem}, \ref{perturbwellks} and the
definition that the chain $\frak m_{1,\beta_0}(\widetilde{\frak X}_{\frak b}^u))$
is a sum of the terms corresponding to
(\ref{splitKuranishi1})-(\ref{splitKuranishi3}).
(Here $\text{\bf p} = \frak b\otimes \dots \otimes \frak b$.)
\par
By Lemma \ref{unitvanish} the term corresponding to (\ref{splitKuranishi2}) vanishes.
\par
By definition, the term corresponding to (\ref{splitKuranishi1}) is
$
\sum_{\beta \ne \beta_0} \frak m_{1,\beta}(\widetilde{\frak X}_{\frak b}^u).
$
The term corresponding to (\ref{splitKuranishi3}) is $\widetilde{\frak{ks}}_{\frak b}(Q)$
by Lemma \ref{KScalcu}. Therefore we have
$$
\widetilde{\frak{ks}}_{\frak b}(Q)\cdot \text{\rm PD}(L(u)) = \delta^{\frak b,b}(\widetilde{\frak X}_{\frak b}^u).
$$
The proposition now follows from (\ref{orio7.132tochu}).
\end{proof}
\begin{prop}\label{generatedede}
The image of the boundary operator $ \delta_{\text{\rm can}}^{\frak b,b}$ 
is contained in the Jacobian ideal.
More precisely, we have the following:
\begin{enumerate}
\item
$$
\text{\rm Im} \delta_{\text{\rm can}}^{\frak b,b}
\subset 
\left( y_i\frac{\partial \frak{PO}_{\frak b}}
{\partial y_i}: i=1,\dots,n\right)H(T^n;\Lambda
\langle\!\langle y,y^{-1}\rangle\!\rangle_0^{\overset{\circ}P} ).
$$
\item
If $\frak b \in \mathcal A(\Lambda_+)$, then 
we have
$$
\text{\rm Im} \delta_{\text{\rm can}}^{\frak b,b}
\subset 
\left( y_i\frac{\partial \frak{PO}_{\frak b}}
{\partial y_i}: i=1,\dots,n\right)H(T^n;\Lambda\langle\!\langle y,y^{-1}\rangle\!\rangle_0^{P}).
$$
\end{enumerate}
\end{prop}
\begin{defn}\label{defg0}
We denote by $G_0$ the discrete submonoid generated by
$$
\{\alpha\cap \omega/2\pi \mid \mathcal M(\alpha) \ne \emptyset\}.
$$
For $\frak b \in \mathcal A(\Lambda_0)$, we denote by $G(\frak b)$
the smallest discrete monoid containing $G_0$ such that $\frak b$ is
$G(\frak b)$-gapped in the sense of Definition \ref{gapping}.
\end{defn}
\begin{proof}[Proof of Proposition \ref{generatedede}]
$G(\frak b)$ defined above is a discrete submonoid such that
$\frak b$ is $G(\frak b)$-gapped and $G(\frak b) \supset G_0$.
We put $G(\frak b) = \{\lambda_0,\lambda_1,\dots\}$ so that
$0 = \lambda_0 < \lambda_1 < \cdots$.
\par
We first prove Statement 1.
Let $\frak I_k$ be the ideal of
$\Lambda
\langle\!\langle y,y^{-1}\rangle\!\rangle_0^{\overset{\circ}P} $ generated by elements
$$
T^{\lambda_l} z_1^{j_1}\cdots z_m^{j_m}
$$
such that $\sum_{i=1}^m j_i + l \ge k$.
We remark that the topology defined by this filtration 
is the same as the topology defined by the metric $d_{\overset{\circ}P}$.
(This is a consequence of Lemma \ref{surjhomfromz}.)
\par
We recall from Definition \ref{defmrho} (see Section \ref{sec:deltaisthesame}) that:
$$
\frak m_2^{\frak b,b}(h_1,h_2)
= \sum_{\ell,\beta}\frac{T^{\beta\cap\omega/2\pi}}{\ell!}
\exp(\beta \cap \frak b_2)
\rho^b(\partial \beta)
\frak q_{\beta;2,\ell}(\frak b^{\otimes \ell}_{\rm high} ;
h_1\otimes h_2).
$$
It defines a ring structure
\begin{equation}\label{cupfrakb}
h_1 \cup^{\frak b,b} h_2 = (-1)^{\deg h_1(\deg h_2 + 1)}
\frak m_2^{\frak b,b}(h_1,h_2)
\end{equation}
on $H(T^n;\Lambda\langle\!\langle y,y^{-1}\rangle\!\rangle_0^{\overset{\circ}P})$. Then
$\cup^{\frak b,b}$ is associative and $G(\frak b)$-gapped.
Consider the generators $\text{\bf e}_i$ ($i=1,\dots,n$) of
$H(T^n;\Z)$.
\begin{lem}\label{densesubroing}
The subring of $(H(T^n;\Lambda\langle\!\langle y,y^{-1}\rangle\!\rangle_0^{\overset{\circ}P}),\cup^{\frak b,b})$
generated by $\{\text{\bf e}_i \mid i=1,\dots,n\}$ is dense.
\end{lem}
Here we use the metric induced by the metric $d_{\overset{\circ}P}$.
\begin{proof}
Let $\Lambda\langle\!\langle y,y^{-1}\rangle\!\rangle_0^{\overset{\circ}P} \langle n\rangle$
be the free associative algebra which is generated by $n$ elements over
the ring $\Lambda\langle\!\langle y,y^{-1}\rangle\!\rangle_0^{\overset{\circ}P}$.
\par
We can define the notion of degree of an element of $\Lambda\langle\!\langle y,y^{-1}\rangle\!\rangle_0^{\overset{\circ}P} \langle n\rangle$
in an obvious way.
Let
$\Lambda\langle\!\langle y,y^{-1}\rangle\!\rangle_0^{\overset{\circ}P} \langle n\rangle_{\le n'}$ 
be the 
$\Lambda\langle\!\langle y,y^{-1}\rangle\!\rangle_0^{\overset{\circ}P}$-submodule of 
$\Lambda\langle\!\langle y,y^{-1}\rangle\!\rangle_0^{\overset{\circ}P} \langle n\rangle$ 
consisting of the elements of degree $\le n'$.

For each $\mathcal P \in \Lambda\langle\!\langle y,y^{-1}\rangle\!\rangle_0^{\overset{\circ}P} \langle n\rangle$
we can define $\mathcal P(\text{\bf e}_1,\dots,\text{\bf e}_n)
\in H(T^n;\Lambda\langle\!\langle y,y^{-1}\rangle\!\rangle_0^{\overset{\circ}P})$
using the product $\cup^{\frak b,b}$ in an obvious way.
\par
The lemma will be an immediate consequence of the following:
\begin{sublem}\label{proddense}
For each given $k$ and $\frak P \in H(T^n;\Lambda\langle\!\langle y,y^{-1}\rangle\!\rangle_0^{\overset{\circ}P})$
there exists $\mathcal P \in \Lambda\langle\!\langle y,y^{-1}\rangle\!\rangle_0^{\overset{\circ}P} \langle n\rangle_{\le n}$ such that
$$
\frak P - \mathcal P(\text{\bf e}_1,\dots,\text{\bf e}_n)
\in \frak I_k H(T^n;\Lambda\langle\!\langle y,y^{-1}\rangle\!\rangle_0^{\overset{\circ}P}).
$$
\end{sublem}
Here and hereafter we write
$$
H(T^n;\frak I_k \Lambda\langle\!\langle y,y^{-1}\rangle\!\rangle_0^{\overset{\circ}P})
=
\frak I_k H(T^n;\Lambda\langle\!\langle y,y^{-1}\rangle\!\rangle_0^{\overset{\circ}P}).
$$
\begin{proof}
We prove this sublemma by an induction over $k$. Since
$$
h_1 \cup^{\frak b,b} h_2 - h_1 \wedge h_2
\in  \frak I_1 H(T^n;\Lambda\langle\!\langle y,y^{-1}\rangle\!\rangle_0^{\overset{\circ}P}),
$$
the case $k=1$ is a consequence of the fact that the usual
cohomology ring $H(T^n;\Z)$ is generated by $\text{\bf e}_i$.
(We also use the fact that any element of $H(T^n;\Z)$ 
is a $\Z$ linear combination of the monomials of $\text{\bf e}_i$'s of degree $\le n$.)
\par
Suppose that Sublemma \ref{proddense} holds for $k-1$ and let
$$
\aligned
\frak P - \mathcal P(\text{\bf e}_1,\dots,\text{\bf e}_n)
& \equiv \sum_{ l,j_1,\dots,j_m}
C(l;j_1,\dots,j_m) T^{\lambda_l} z_1^{j_1}\cdots z_m^{j_m} \\
& \quad \mod \frak I_{k} H(T^n;\Lambda\langle\!\langle y,y^{-1}\rangle\!\rangle_0^{\overset{\circ}P}),
\endaligned
$$
where the summation is taken over $l,j_1,\dots,j_m$,
such that $\sum_{i=1}^m j_i + l= k-1$ and
$C(l;j_1,\dots,j_m)
\in H(T^n;\Lambda\langle\!\langle y,y^{-1}\rangle\!\rangle_0^{\overset{\circ}P})$.
Applying the case $k=1$, we can find
$\mathcal P(l;j_1,\dots,j_m) \in \Lambda\langle\!\langle y,y^{-1}\rangle\!\rangle_0^{\overset{\circ}P} \langle n\rangle_{\le n}$ such that
$$
C(l;j_1,\dots,j_m)  - \mathcal P(l;j_1,\dots,j_m)
(\text{\bf e}_1,\dots,\text{\bf e}_n)
\in \frak I_1 H(T^n;\Lambda\langle\!\langle y,y^{-1}\rangle\!\rangle_0^{\overset{\circ}P}).
$$
If we put
$$
\mathcal P'
= \mathcal P +
\sum_{l,j_1,\dots,j_m}
T^{\lambda_l} z_1^{j_1}\cdots z_m^{j_m}\mathcal P(l;j_1,\dots,j_m),
$$
we have
$$
\frak P - \mathcal P'(\text{\bf e}_1,\dots,\text{\bf e}_n)
\in \frak I_{k} H(T^n;\Lambda\langle\!\langle y,y^{-1}\rangle\!\rangle_0^{\overset{\circ}P})
$$
by definition. Sublemma \ref{proddense} is proved.
\end{proof}
Hence the proof of Lemma \ref{densesubroing} is completed.
\end{proof}
\par
We next remark the identity
\begin{equation}\label{derivation}
\delta_{\text{\rm can}}^{\frak b,b}(h_1 \cup^{\frak b,b} h_2)
= \pm \delta_{\text{\rm can}}^{\frak b,b}(h_1) \cup^{\frak b,b} h_2
\pm h_1 \cup^{\frak b,b} \delta_{\text{\rm can}}^{\frak b,b}(h_2).
\end{equation}
We note that (\ref{derivation}) is a consequence of the filtered $A_{\infty}$
relation of $\frak m^{\frak b,b}$ and the fact $\frak m_0(1)$ and so
$\frak m_0^{\frak b,b}(1)$ is proportional to the (exact) unit.
\par
We also remark the equality
\begin{equation}\label{unitiscycle}
\delta_{\rm can}^{\frak b,b}(\text{\bf e}_0) = 0
\end{equation}
where $\text{\bf e}_0$ is the unit, which is the Poincar\'e dual to the 
fundamental cycle. 
(In fact (\ref{unitiscycle}) holds for any canonical model\index{canonical model} of 
homotopically unital filtered $A_{\infty}$ algebra.
\cite[Theorem $5.4.2'$]{fooobook2}.)
\par
We also have:
\begin{equation}\label{deltaanddiff}
\delta_{\rm can}^{\frak b,b}(\text{\bf e}_i) = y_i\frac{\partial\frak{PO}_{\frak b}}
{\partial y_i}\text{\bf e}_0.
\end{equation}
In fact, by definition, we have:
$$
\frak{PO}_{\frak b}(b)\text{\bf e}_0 = \sum_{k=0}^{\infty}\frak m_{k}^{\frak b}(b,\dots,b).
$$
Therefore using $b = \sum_{i=1}^n x_i\text{\bf e}_i$, 
$y_i = e^{x_i}$, ($y_i\frac{\partial}{\partial y_i} = \frac{\partial}{\partial x_i}$) we have
$$
 y_i\frac{\partial\frak{PO}_{\frak b}}
{\partial y_i}\text{\bf e}_0
= \sum_{k_1=0}^{\infty}\sum_{k_2=0}^{\infty} \frak m_{k_1+k_2+1}^{\frak b}(b^{\otimes k_1},\text{\bf e}_i,b^{\otimes k_1})
= \delta_{\rm can}^{\frak b,b}(\text{\bf e}_i).
$$
\begin{lem}
For any
$\mathcal P \in  \Lambda\langle\!\langle y,y^{-1}\rangle\!\rangle_0^{\overset{\circ}P} \langle n\rangle$
we have:
\begin{equation}\label{iretemod0}
\delta_{\text{\rm can}}^{\frak b,b}(\mathcal P(\text{\bf e}_1,\dots,\text{\bf e}_n))
\in \left( y_i\frac{\frak{PO}_{\frak b}}
{\partial y_i}: i=1,\dots,n\right)H(T^n;\Lambda\langle\!\langle y,y^{-1}\rangle\!\rangle_0^{\overset{\circ}P}).
\end{equation}
\par
If we assume 
$$
\mathcal P \in \frak I_{k}\Lambda\langle\!\langle y,y^{-1}\rangle\!\rangle_0^{\overset{\circ}P} \langle n\rangle_{\le n}
$$
in addition, 
then there exists
$$
\mathcal R_{i;\mathcal P} \in \frak I_{k}\Lambda\langle\!\langle y,y^{-1}\rangle\!\rangle_0^{\overset{\circ}P} \langle n\rangle_{\le n}
$$
such that
\begin{equation}\label{iretemod0degincl}
\delta_{\text{\rm can}}^{\frak b,b}(\mathcal P(\text{\bf e}_1,\dots,\text{\bf e}_n))
= \sum_{i=1}^n y_i\frac{\frak{PO}_{\frak b}}
{\partial y_i} \mathcal R_{i;\mathcal P}(\text{\bf e}_1,\dots,\text{\bf e}_n) .
\end{equation}
\end{lem}
\begin{proof}
Using (\ref{derivation}), (\ref{unitiscycle}) and (\ref{deltaanddiff}),
we can prove the lemma by induction over the degree of  $\mathcal P$.
\end{proof}
Now we are in the position to complete the proof of Proposition \ref{generatedede}.
Let $\frak P \in H(T^n;\Lambda\langle\!\langle y,y^{-1}\rangle\!\rangle_0^{\overset{\circ}P})$.
By Lemma \ref{densesubroing} and Sublemma \ref{proddense} we can write it as a series 
$$
\frak P = \sum_{j=1}^{\infty} \mathcal P_j(\text{\bf e}_1,\dots,\text{\bf e}_n)
$$
where
$$
 \mathcal P_j \in \frak I_{k_j}\Lambda\langle\!\langle y,y^{-1}\rangle\!\rangle_0^{\overset{\circ}P} \langle n\rangle_{\le n}.
$$
By (\ref{iretemod0degincl}) we have $\mathcal R_{j,i} 
\in \frak I_{k_j}\Lambda\langle\!\langle y,y^{-1}\rangle\!\rangle_0^{\overset{\circ}P} \langle n\rangle_{\le n}$ such that
$$
\delta_{\text{\rm can}}^{\frak b,b}(\mathcal P_j(\text{\bf e}_1,\dots,\text{\bf e}_n))
=
\sum_{i=1}^ny_j\frac{\frak{PO}_{\frak b}}
{\partial y_i} \mathcal R_{j,i} (\text{\bf e}_1,\dots,\text{\bf e}_n)
$$
Note $\Lambda\langle\!\langle y,y^{-1}\rangle\!\rangle_0^{\overset{\circ}P} \langle n\rangle_{\le n}$
as a $\Lambda\langle\!\langle y,y^{-1}\rangle\!\rangle_0^{\overset{\circ}P}$ module is 
a direct sum of finitely many copies of $\Lambda\langle\!\langle y,y^{-1}\rangle\!\rangle_0^{\overset{\circ}P}$.
We define a topology on $\Lambda\langle\!\langle y,y^{-1}\rangle\!\rangle_0^{\overset{\circ}P} \langle n\rangle_{\le n}$
using this isomorphism.
Then, since $\lim_{j\to \infty}k_j = \infty$ the series
$
\sum_{j=1}^{\infty} \mathcal R_{j,i}
$
converges in $\Lambda\langle\!\langle y,y^{-1}\rangle\!\rangle_0^{\overset{\circ}P} \langle n\rangle_{\le n}$.
We put
$$
\mathcal R_i =  \sum_{j=1}^{\infty} \mathcal R_{j,i} \in \Lambda\langle\!\langle y,y^{-1}\rangle\!\rangle_0^{\overset{\circ}P} \langle n\rangle_{\le n}.
$$
Then we have
$$
\delta_{\text{\rm can}}^{\frak b,b}(\frak P)
=
\sum_{i=1}^{n} y_i\frac{\frak{PO}_{\frak b}}
{\partial y_i}\mathcal R_{i}(\text{\bf e}_1,\dots,\text{\bf e}_n).
$$
We have thus proved Proposition \ref{generatedede}.1.
\par
We next prove Proposition \ref{generatedede}.2.
\begin{lem}\label{densesubroing2}
If $\frak b \in \mathcal A(\Lambda_+)$ then
the subring of $(H(T^n;\Lambda\langle\!\langle y,y^{-1}\rangle\!\rangle_0^{P}),\cup^{\frak b,b})$
generated by $\{\text{\bf e}_i \mid i=1,\dots,n\}$ is dense.
\end{lem}
We use the metric $d_P$ to define the topology on $H(T^n;\Lambda\langle\!\langle y,y^{-1}\rangle\!\rangle_0^{P})$.
\begin{proof}
Let $\frak J_k$ be the ideal $T^{\lambda_k}\Lambda\langle\!\langle y,y^{-1}\rangle\!\rangle_0^{P}$ 
of $\Lambda\langle\!\langle y,y^{-1}\rangle\!\rangle_0^{P}$.
\begin{sublem}\label{proddense2}
For each given $k$ and $\frak P \in H(T^n;\Lambda\langle\!\langle y,y^{-1}\rangle\!\rangle_0^{P})$
there exists $\mathcal P \in \Lambda\langle\!\langle y,y^{-1}\rangle\!\rangle_0^{P} \langle n\rangle_{\le n}$ such that
$$
\frak P - \mathcal P(\text{\bf e}_1,\dots,\text{\bf e}_n)
\in \frak J_k H(T^n;\Lambda\langle\!\langle y,y^{-1}\rangle\!\rangle_0^{P}).
$$
\end{sublem}
\begin{proof}
We first prove the case $k=1$.
Let $\lambda$ be the symplectic area of the non-constant pseudoholomorphic sphere of 
the smallest area in $X$.
Note that  $\lambda_1 \le 2\pi \lambda$.  
\par
We consider the operator 
\begin{equation}\label{m2terms}
\aligned
&\frak m_2^{\beta,b}(h_1,h_2) \\
&= \sum_{\ell=0}^{\infty}\sum_{k_1,k_2,k_3=0}^{\infty}
\frac{T^{\beta\cap \omega/2\pi}}{\ell !}
\rho^{b}(\partial\beta)
\frak q_{\ell;k_1+k_2+k_3+2;\beta}(\frak b^{\ell};b^{k_1},h_1,b^{k_2},h_2,b^{k_3})
\endaligned
\end{equation}
for $\beta \ne 0$. Suppose 
the term for $\ell = \ell_0$ in the right hand side of  (\ref{m2terms}) is nonzero.
Then the moduli space $\mathcal M_{1;\ell_0}(L(u),\beta;\frak b^{\otimes \ell_0})$ is nonempty.
\par
If $\ell_0 >0$ in addition, then the corresponding term in the right hand side of  (\ref{m2terms}) is 
zero modulo $\lambda_1$.
(This is because $\frak b \equiv 0 \mod T^{\lambda_1}$.)
\par
We next consider the case when $\ell_0 =0$.
If Maslov index of $\beta$ is strictly smaller than $2$, then \cite[Theorem 11.1]{fooo08} implies that
$$
\beta \cap \omega \ge \lambda.
$$ 
(In fact an element of $\mathcal M_{1;0}(L(u),\beta)$ must have a nontrivial sphere bubble in this case.)
\par
In case the Maslov index of $\beta$ is not smaller than $2$,  
we have
$$
\frak m_{k_1+k_2+k_3+2;\beta}(\frak b^{\ell};b^{k_1},h_1,b^{k_2},h_2,b^{k_3}) \in \bigoplus_{d < \deg h_1 + \deg h_2} H^d(L(u);\Q)
$$
for $h_1, h_2 \in H(L(u);\Q)$.
In sum we have
$$
\aligned
h_1 \cup^{\frak b,b} h_2 - h_1 \wedge h_2
\in  &\frak J_1 H(T^n;\Lambda\langle\!\langle y,y^{-1}\rangle\!\rangle_0^P)\\
&\oplus 
\bigoplus_{d < \deg h_1 + \deg h_2} H^d(L(u);\Lambda\langle\!\langle y,y^{-1}\rangle\!\rangle_0^P).
\endaligned
$$
Using it we can prove the case $k=1$, by a downward induction of the cohomology degree of the element $\frak P$.
\par
Using the case $k=1$, we can prove  the case when $k$ is general in the same way as Sublemma \ref{proddense}.
\end{proof}
Lemma \ref{densesubroing2} follows from Sublemma \ref{proddense2}.
\end{proof}
Once Lemma \ref{densesubroing2} and Sublemma \ref{proddense2} are proved, the rest of the 
proof of  Proposition \ref{generatedede}.2 is similar to the  Proposition \ref{generatedede}.1 and is omitted.
\end{proof}
Therefore Propositions \ref{generatedede} and \ref{bdryandks} imply
Theorem \ref{indepenceKS} in the case when $\deg Q > 2$.

We finally consider the case when degree of $Q$ is $2$.
We put $Q = \sum_{i=1}^m Q_i\text{\bf f}_i$ and
$Q^* = \sum_{i=1}^m Q_iD_i \in H_{2n-2}(D;R)$.
\begin{lem}\label{deg2casecalculate}
We have
$$\aligned
\widetilde{\frak{ks}}_{\frak b}(Q)(y)
= \sum_{\ell=0}^{\infty}\sum_{\beta \in H_2(X;L(u);\Z)} &\frac{1}{\ell!}
(Q^* \cap \beta) \rho^{b}(\partial \beta)
\exp (\frak b_2 \cap \beta) T^{\omega\cap \beta/2\pi}
\\
&\text{\rm ev}_{0 *}(\mathcal M^{\text{\rm main}}_{1;\ell}(L(u),\beta;\frak b_{\rm high}^{\otimes \ell})).
\endaligned
$$
Here $Q^* \cap \beta$ is defined by the intersection pairing
$$
\cap: H_{2n-2}(D;\C) \otimes H_2(X;L(u);\C) \to \C.
$$
This is well defined since $D \cap L(u) = \emptyset$.
\end{lem}
\begin{proof}
This is a consequence of Lemma \ref{deg2separate}.
\end{proof}
Since $[Q^*] = 0$ in $H_{2n-2}(X;\C)$ by the hypothesis, there exists $R$ such that
$\partial R = Q^*$. We have
$[R] \in H_{2n-1}(X, D;\C)$,
using the exactness and compatiblity of the intersection pairing of
\begin{equation}\label{XLuhomoexact}
\aligned
& H_2(X;\C) \to H_2(X,L(u);\C)\to H_1(L(u);\C) \\
& H_{2n-2}(X;\C) \leftarrow H_{2n-2}(D;\C) \leftarrow H_{2n-1}(X,D;\C), 
\endaligned
\end{equation}
we can find $d_i$ such that
\begin{equation}\label{cuptopocalcu}
Q^* \cap \beta = \sum_{i=1}^n d_i(\partial \beta \cap \text{\bf e}_i).
\end{equation}
By differentiating (\ref{POformula}), we have
\begin{equation}
\aligned
y_i \frac{\partial\frak{PO}_{\frak b}}{\partial y_i}
= \sum_{\ell=0}^{\infty}\sum_{\beta\in H_2(X;L(u);\Z)} &\frac{T^{\omega\cap \beta/2\pi}}{\ell!}
(\partial \beta \cap \text{\bf e}_i)\rho^{b}(\partial \beta) \exp (\frak b_2 \cap \beta)
\\
&\text{\rm ev}_{0 *}(\mathcal M^{\text{\rm main}}_{1;\ell}(L(u),\beta;\frak b_{\rm high}^{\otimes \ell})).
\endaligned
\end{equation}
Therefore
$$
\widetilde{\frak{ks}}_{\frak b}(Q) =  \sum_{i=1}^n d_i y_i \frac{\partial\frak{PO}_{\frak b}}{\partial y_i} .
$$
The proof of Theorem \ref{indepenceKS} is now complete.
\qed
\section{Well-definedness of potential function}
\label{sec:uptovariables}
\par
Let $\frak b, \frak b' \in \mathcal A(\Lambda_0)$, which
represent the same cohomology class in $H(X;\Lambda_0)$.
In this section we prove that $\mathfrak{PO}_{\frak b}$
coincides with  $\mathfrak{PO}_{\frak b'}$ after an
appropriate coordinate change of variables.
Some kinds of such statement follow from the
well-definedness of potential function up to filtered
$A_{\infty}$ homotopy equivalence, which we established in
\cite{fooo06}, for arbitrary symplectic manifold and
its relatively spin Lagrangian submanifold.
(See also \cite{Aur07} for some discussion related to it.)
On the other hand, for the purpose of this paper, we need to
specify the type of coordinate change we use. For example,
we want them to converge everywhere in the moment polytope.
For this purpose we use the special feature of our situation.
\par
We start with the definition of the class of coordinate changes.

\begin{defn}\label{coordinatechange}
We consider $n$ elements $y_i' \in \Lambda\langle\!\langle y,y^{-1}\rangle\!\rangle^{\overset{\circ}P}$ $(i=1,\dots,n)$.
See Definition \ref{definitonPcirczero} for this notation.
\begin{enumerate}
\item We say that $y' = (y'_1,\dots,y'_n)$ is a {\it coordinate change
converging on} $\text{\rm Int} P$
(or a {\it coordinate change on $\text{\rm Int}P$}) if
\begin{equation}\label{leadingcondition}
y'_i \equiv c_i y_i  \mod y_i\Lambda\langle\!\langle y,y^{-1}\rangle\!\rangle_+^{\overset{\circ}P}
\end{equation}
$c_i \in \C \setminus \{0\}$.\index{coordinate change of $y_i$}
\item We say that the coordinate change is {\it strict} if $c_i =1$ for all $i$.
\item We say that the {\it coordinate change
converges on $P$} if
$y_i' \in \Lambda\langle\!\langle y,y^{-1}\rangle\!\rangle^{P}$ $(i=1,\dots,n)$
in addition.
Its strictness is defined in the same way.
We also say that $y'$ is
a {\it coordinate change on $P$}.
\item
Coordinate change $y'$ on $P$ is said to be
{\it $G$-gapped} if $(y'_i - c_iy_i)/y_i$ is $G$-gapped.
(See Definition \ref{def24}.)
\item
{\it Coordinate change converging at $u$} (or {\it coordinate change at $u$}) is
defined as an element
$$
y(u)'_i \in \Lambda_0\langle\!\langle y(u),y(u)^{-1}\rangle\!\rangle
$$
such that
$$
y'_i(u) \equiv c_i y_i(u) \mod y_i(u)\Lambda\langle\!\langle y(u),y(u)^{-1}\rangle\!\rangle_+^u.
$$
Its strictness is defined in the same way.
\end{enumerate}
\end{defn}
Let $y' = (y'_1,\dots,y'_n)$ be a coordinate change
converging on $\text{\rm Int} P$.
If a series
$$
\mathcal P = \sum_i a_iT^{\lambda_i}y^{v_i}
$$
converges in
$\Lambda\langle\!\langle y,y^{-1}\rangle\!\rangle_0^{\overset{\circ}P}$,
then
\begin{equation}\label{composition}
\sum_i a_iT^{\lambda_i}(y')^{v_i}
\end{equation}
also converges.  Here
$v_i = (v_i^{(1)},\dots,v_i^{(n)}) \in \Z^n$ and
$$
(y')^{v_i} = (y'_1)^{v_i^{(1)}}\cdots (y'_n)^{v_i^{(n)}}.
$$
We remark that
$$
(y_i')^{-1}
= c_i^{-1}y_i^{-1}\left(1+\frac{y'_i - c_iy_i}{c_iy_i}\right)^{-1}
= c_i^{-1}y_i^{-1}\sum_{k=0}^{\infty}
(-1)^k\left(\frac{y'_i - c_iy_i}{c_iy_i}\right)^k
$$
converges in $\frak v_T^{\overset{\circ}P}$ norm by
(\ref{leadingcondition}).
Therefore (\ref{composition}) makes sense.
\begin{defn}
We denote (\ref{composition}) by
$\mathcal P(y')$ and call it the {\it coordinate change} of
$\mathcal P$ by $y'$.
There are the versions corresponding to
Definition \ref{coordinatechange}.3 or 4, which
are defined in the same way.
We can define a composition of coordinate
change by applying (\ref{composition}) in each
factor. We write the composition
of $y'$ and $y''$ by $y'\circ y''$.
\end{defn}

The proof of the following lemma is easy and so omitted.

\begin{lem}\label{cctrivial}
Coordinate change induces a continuous endomorphism
on $\Lambda\langle\!\langle y,y^{-1}\rangle\!\rangle_0^{\overset{\circ}P}$,
$\Lambda\langle\!\langle y,y^{-1}\rangle\!\rangle_0^{P}$,
and $\Lambda\langle\!\langle y,y^{-1}\rangle\!\rangle_0^u$ respectively.
Therefore composition thereof induces a coordinate change again.
Composition of coordinate change is associative.
Composition of $G$-gapped coordinate changes are $G$-gapped.
\end{lem}

\begin{conv}
Let $Z_1, \dots, Z_m$ be formal variables.
Let $z_1, \dots, z_m$ be the elements in $\Lambda_0\langle\!\langle y,y^{-1}\rangle\!\rangle $
defined by \eqref{zjdef}. For each $R \in \Lambda_0[[Z_1,\dots,Z_m]]$, we put
$$
R(z) = R(z_1,\dots,z_m).
$$
This is well-defined as an element of $\Lambda\langle\!\langle y,y^{-1}\rangle\!\rangle_0^{\overset{\circ }P}$
by Lemma \ref{surjhomfromz}.
\end{conv}

\begin{lem}\label{ccinverse}
Coordinate change has an inverse. Namely
for $y'$ there exists $y''$ such that
$y' \circ y''(y) = y'' \circ y'(y) = y$.
\end{lem}
\begin{proof}
It suffices to prove the lemma for the case where $y'$ is strict.
We first consider the case $y'$ is a coordinate
change on $P$. We assume that $y'$ is $G$-gapped
and put $G=\{\lambda_0,\lambda_1,\dots\}$.
($\lambda_0 =0 < \lambda_1 < \cdots$.)
For each $k=1,2,\dots$, we are going to find a $G$-gapped coordinate change on $P$,
$y''_{(k)} = (y''_{(k),1},\dots,y''_{(k),n})$ such that
\begin{enumerate}
\item
$y''_{(k)}$ is $G$-gapped and strict.
\item
$
((y''_{(k)}\circ y')_i-y_i)/y_i \in T^{\lambda_k}\Lambda\langle\!\langle y,y^{-1}\rangle\!\rangle_0^{P}.
$
\end{enumerate}
We prove this by induction over $k$. We assume
$$
((y''_{(k)}\circ y')_i -y_i)/y_i \equiv a_i(z)
T^{\lambda_k} \mod T^{\lambda_{k+1}}\Lambda\langle\!\langle y,y^{-1}\rangle\!\rangle_0^{P}
$$
where $a_i \in R[Z_1,\dots,Z_m]$. We put
$$
y''_{(k+1),i} = (1 - T^{\lambda_k}a_i)y^{\prime\prime}_{(k),i}.
$$
Since $y'_i \equiv y''_i \equiv y_i \mod y_i\Lambda
\langle\!\langle y,y^{-1}\rangle\!\rangle_+^{P}$, it is easy to see that
$y''_{(k+1),i}$ satisfies (1),(2).
\par
Now $\lim_{k\to\infty} y''_{(k+1)}$ converges and gives a left
inverse to $y'$.
Existence of the right inverse can be proved in the
same way. Then the right and the left inverse coincide
by a standard argument of group theory.
\par
The case of coordinate change at $u$ is the same.
\par
For the case of coordinate change on $\text{\rm Int} P$
we proceed as follows.
We note that if $y'$ is a coordinate change
on $P$, it induces an isomorphism
on $\Lambda\langle\!\langle y,y^{-1}\rangle\!\rangle_0^{P}$.
Now let $y'$ be a coordinate change on $\text{\rm Int} P$.
Then for each sufficiently small $\e > 0$, the coordinate change $y'$ induces an isomorphism
on $\Lambda\langle\!\langle y,y^{-1}\rangle\!\rangle_0^{P_{\e}}$.
Therefore by taking the projective limit of these isomorphisms
as $\e \to 0$, which obviously exists, the coordinate change $y'$ induces an isomorphism on
$\Lambda\langle\!\langle y,y^{-1}\rangle\!\rangle_0^{\overset{\circ}P}$ by Lemma \ref{Poisprojlim}.
Now it is easy to find its inverse.
\end{proof}
\begin{rem}
The set of all coordinate changes forms a solvable group.
It is related to the group of self homotopy equivalences of
filtered $A_{\infty}$ algebra and also relevant to the
nilpotent group discussed in \cite{KS2} and \cite{grospand}.
We remark however in \cite{KS2} and \cite{grospand} the coordinate changes that
preserve the volume element $dy_1\cdots dy_n/(y_1\cdots y_n)$ are studied.
We do not prove here that the coordinate change we produce
(for example in Theorem \ref{cchangethem} below) has this property.
This point is related to the Poincar\'e duality and cyclic symmetry.
So if we use the cyclically symmetric version of the filtered $A_{\infty}$
structure that we produce in Section \ref{sec:cyclic}, we may obtain a coordinate change which preserves the volume element.
\end{rem}
For $\mathcal P \in \Lambda\langle\!\langle y,y^{-1}\rangle\!\rangle_0^{\overset{\circ}P}$,
we define its Jacobian ring by
$$
\text{\rm Jac}(\mathcal P)
=
\frac{\Lambda\langle\!\langle y,y^{-1}\rangle\!\rangle_0^{\overset{\circ}P} }
{\text{Clos}_{d_{\overset{\circ}P}}\left(y_i\frac{\partial \mathcal P}{\partial y_i}: i=1,\dots,n\right)}.
$$
It follows from Lemmata \ref{cctrivial}, \ref{ccinverse}
that a coordinate change induces an isomorphism
of Jacobian rings.
\par
The following is the main result of this section.
\begin{thm}\label{cchangethem}
Let $\frak b,\frak b' \in \mathcal A(\Lambda_0)$.
We assume that $[\frak b] = [\frak b']
\in H(X;\Lambda_0)$.
\par
Then there exists a coordinate change $y'$ on $\text{\rm Int} P$, such that
\begin{equation}\label{POcoordinatechage}
\mathfrak{PO}_{\frak b}(y')
= \mathfrak{PO}_{\frak b'}(y).
\end{equation}
\par
If $\frak b - \frak b' \in \mathcal A(\Lambda_+)$,
then $y'$ can be taken to be strict.
\par
If both $\frak b,\frak b' \in \mathcal A(\Lambda_+)$,
then $y'$ can be taken to be
a strict coordinate change on $P$.
\end{thm}
\begin{exm}\label{coorinateexam}
Let $X=\C P^2$ with the moment polytope $P$ as in Section \ref{sec:example}.
Put
$$
P_i = \{(u_1,u_2) \in P \mid u_i = 0\}, \quad \text{\bf f}_i = {\rm PD} (\pi^{-1}(P_i))
$$
and let $c\in \R\setminus\{0\}$. We consider $\frak b= c\text{\bf f}_1$ and  $\frak b'= c\text{\bf f}_2$.
Then we have
$$
\mathfrak{PO}_{\frak b} = e^{c}z_1 + z_2 + z_3,
\qquad \mathfrak{PO}_{\frak b'} = z_1 + e^{c}z_2 + z_3.
$$
We put $y'_1 = e^{-c} y_1$, $y'_2 = e^{c} y_2$.
It is easy to see (\ref{POcoordinatechage}) holds.
\end{exm}
\begin{proof}
We first consider the case $\frak b - \frak b' \in \bigoplus_{k\ge 4}\mathcal A^k(\Lambda_0)$.
Moreover we also consider $\frak b, \frak b' \in \bigoplus_{k\ge 4}\mathcal A^k(\Lambda_+)$
since the same argument applies to the case with $\frak b, \frak b' \in \mathcal A(\Lambda_0)$
with $\Lambda\langle\!\langle y, y^{-1}\rangle\!\rangle_0^{P}$ being replaced by $\Lambda\langle\!\langle y,y^{-1}\rangle\!\rangle_0^{\overset{\circ}P}$.

We take a new formal variable $s$ and put
$$
\frak b(s) = s \frak b + (1-s) \frak b'.
$$
Let $R$ be a smooth singular chain
$
\partial R = \frak b - \frak b'.
$
For given $u \in \text{Int}P$, we define
\begin{equation}\label{frakXSdef}
\aligned
\tilde{\frak X}^u_S(s,b) =
\sum_{\ell=0}^{\infty}\sum_{\beta \in H_2(X;L(u);\Z)} &\frac{ T^{\omega\cap \beta/2\pi}}{\ell!}
\rho^{b}(\partial \beta) \exp(\frak b_2 \cap \beta) \\
&\text{\rm ev}_{0 *}
(\mathcal M^{\text{\rm main}}_{1;\ell+1}(L(u),\beta;\frak b(s)_{\rm high}^{\otimes \ell} \otimes R)),
\endaligned\end{equation}
\begin{equation}\label{frakPXSdef}
{\frak X}^u_S(s,b) = \Pi(\tilde{\frak X}^u_S(s,b)),
\end{equation}
where $\Pi$ is as in (\ref{harmonicproj}).
Here $\frak b_{\rm high}$ is degree $\ge 4$ part and $\frak b_2$ is the
degree $2$ part. We remark that $\frak b_2 = \frak b'_2 = \frak b(s)_2$ by
the hypothesis.
\begin{defn}
We denote by  $\Lambda\langle\!\langle s, y,y^{-1}\rangle\!\rangle^P_0$
the set of all formal power series $\sum_{k=0}^{\infty} P_k s^k$
where 
$P_k \in \Lambda\langle\!\langle  y,y^{-1}\rangle\!\rangle^P_0$
with 
$$
\lim_{k\to \infty}\frak v_T^{u}(P_k) = 0
$$
for any $u \in P$.
\par
We say an element $\sum_{k=0}^{\infty} P_k s^k 
\in \Lambda\langle\!\langle s, y,y^{-1}\rangle\!\rangle^P_0$ is 
an element of
$\Lambda\langle\!\langle s, y,y^{-1}\rangle\!\rangle^P_+$
if $P_k \in \Lambda\langle\!\langle  y,y^{-1}\rangle\!\rangle^P_+$
for all $k$.
\end{defn}
\par
In the same way as the proof of
Proposition \ref{fraXconv}, we can prove that
(\ref{frakXSdef}) converges in $\frak v_T^P$ norm.
Hence ${\frak X}^u_S(s,b)$ can be regarded as an element
of $H(L(u);\Lambda\langle\!\langle s,y,y^{-1}\rangle\!\rangle^P_+)$.
In the same way as the proof of \eqref{conclusionsec4} we have
\begin{equation}\label{ksintheimage}
\widetilde{\mathfrak{ks}}_{\frak b(s)}(\frak b - \frak b')
=
\delta^{\frak b(s),b}({\frak X}^u_S(s,b)).
\end{equation}
We write
$$
\frak X_{S}^u(s,b) = \sum_{k=0}^{\infty} T^{\lambda_k} P_k s^k
$$
where $P_k \in  \Lambda\langle\!\langle  y,y^{-1}\rangle\!\rangle^P_0$, 
$\lambda_k \to +\infty$.
We then apply Proposition \ref{generatedede} to each  $P_k$ and
find
$$
\frak X_{S,i}^u(s,b) \in \Lambda \langle\!\langle s,y,y^{-1}\rangle\!\rangle_+^P
$$
such that
$$
\delta^{\frak b,b}(\frak X_{S}^u(s,b))
= \sum_i \frak X_{S,i}^u(s,b) y_i \frac{\partial\frak{PO}_{\frak b(s)}}{\partial y_i}.
$$
\par
We note that
\begin{equation}\label{poderived}
\widetilde{\mathfrak{ks}}_{\frak b(s)}(\frak b - \frak b')
= \frac{\partial \frak {PO}_{\frak b(s)}}{\partial s}.
\end{equation}
Let
$$
y_i(s) \in \Lambda\langle\!\langle s,y,y^{-1}\rangle\!\rangle^{P}
$$
satisfy the same properties as the coordinate change on $P$.
Using \eqref{poderived}, we calculate
$$
\frac{\partial}{\partial s}\left(
\frak {PO}_{\frak b(s)}(y(s))\right)
= \widetilde{\mathfrak{ks}}_{\frak b(s)}(\frak b - \frak b')
\circ y(s)
+ \sum_i
\left(\frac{\partial\frak {PO}_{\frak b(s)}}{\partial y_i}
\frac{\partial y_i(s)}{\partial s}\right)\circ y(s).
$$
The right hand side is zero if
\begin{equation}\label{ODEy'}
y_i(s)  {\frak X}^u_{S,i} + \frac{\partial y_i(s)}{\partial s}
= 0.
\end{equation}
We will assume
$y_i(s) \equiv y_i \mod y_i\Lambda\langle\!\langle s,y,y^{-1}\rangle\!\rangle_+^{P}$.
\par
Motivated by \eqref{leadingcondition} and \eqref{ODEy'}, we set
$$
Y_i(s) = \log (y_i(s)/y_i) \in \Lambda\langle\!\langle y,y^{-1}\rangle\!\rangle_+^P.
$$
Then (\ref{ODEy'}) is equivalent to
\begin{equation}\label{ODEY}
\frac{d}{ds}Y_i(s) +  {\frak X}^u_{S,i}  = 0.
\end{equation}
This is an ordinary differential equation
and ${\frak X}^u_{S,i} \in \Lambda\langle\!\langle s,y,y^{-1}\rangle\!\rangle_+^{P}$. Therefore it has a unique solution with
$Y_i(0) = 0$, which we again denote by $Y_i$. We now define
$$
y_i(s) = \exp(Y_i(s)) y_i.
$$
Then we have
$$
\frac{\partial}{\partial s}\left(
\frak {PO}_{\frak b(s)}(y(s))\right)
= 0.
$$
Therefore
$$
\frak {PO}_{\frak b'}(y(1))
= \frak {PO}_{\frak b(1)}(y(1))
= \frak {PO}_{\frak b(0)}(y(0))
=  \frak {PO}_{\frak b}(y).
$$
Hence $y' = y(1)$ is the coordinate change we are looking for.
\par
We finally consider the case $\frak b' - \frak b \in \mathcal A^2(\Lambda_0)$.
We use the notation of the proof of the corresponding case of
Theorem \ref{indepenceKS}. We put $Q = \frak b' - \frak b$ and put
$\partial R = Q^*$.  
In the same way as (\ref{cuptopocalcu}) we can
find $d_i$ such that
\begin{equation}\label{cuptopocalcu2}
Q^* \cap \beta = \sum_{i=1}^n d_i(\partial \beta \cap \text{\bf e}_i).
\end{equation}
We then put
\begin{equation}\label{coordinatechange2}
y'_i = e^{d_i} y_i.
\end{equation}
Now we have (see (\ref{PPformula}))
$$
\aligned
\frak{PO}_{\frak b}(y')
&= \sum_{\ell=0}^{\infty}\sum_{\beta\in H_2(X;L(u);\Z)} \frac{T^{\omega\cap \beta/2\pi}}{\ell!}
\exp(\frak b_2 \cap \beta)
y'(u)_1^{\partial \beta \cap \text{\bf e}_1}\cdots
y'(u)_n^{\partial \beta \cap \text{\bf e}_n}\\
&\qquad\quad\qquad \text{\rm ev}_{0 *}(\mathcal M^{\text{\rm main}}_{1;\ell}(L(u),\beta;\frak b_{\rm high}^{\otimes \ell})) \\
&=\sum_{\ell=0}^{\infty}\sum_{\beta\in H_2(X;L(u);\Z)} \frac{T^{\omega\cap \beta/2\pi}}{\ell!}
\exp(\frak b'_2 \cap \beta)
y(u)_1^{\partial \beta \cap \text{\bf e}_1}\cdots
y(u)_n^{\partial \beta \cap \text{\bf e}_n}\\
&\qquad\quad\qquad \text{\rm ev}_{0 *}(\mathcal M^{\text{\rm main}}_{1;\ell}(L(u),\beta;\frak b_{\rm high}^{\otimes \ell})) \\
&= \frak{PO}_{\frak b'}(y).
\endaligned$$
The proof of Theorem \ref{cchangethem} is complete.
\end{proof}
\par
\section{Kodaira-Spencer map is a ring homomorphism}
\label{sec:ringhomo}
\par
In this section we prove the following theorem. 
\begin{thm}\label{multiplicative}
The map
$
\frak{ks}_{\frak b}: (H(X;\Lambda_0),\cup^{\frak b})
\to \text{\rm Jac}(\frak{PO}_{\frak b})
$
is a ring homomorphism.
\end{thm}
\begin{rem}
Various statements related to this theorem have
been discussed in the previous literature sometimes in a greater
generality. We explain some of them in Section \ref{sec:QMHOch}.
\end{rem}
\par
We fix a basis $\overline{\text{\bf f}}_i$,
$i = 0,\dots,B'$ of $H(X;\C)$ as in Section \ref{sec:statements}.
We define $c^c_{a(1) a(2)} \in \Lambda_0$ by
\begin{equation}\label{defcabc}
\overline{\text{\bf f}}_{a(1)} \cup^{\frak b} {\overline{\text{\bf f}}}_{a(2)}
= \sum_{c=0}^{B'} c^c_{a(1) a(2)} \overline{\text{\bf f}}_c.
\end{equation}
We will prove
\begin{equation}\label{eq:ks-ring}
\aligned
\widetilde{\frak{ks}}_{\frak b}(\overline{\text{\bf f}}_{a(1)})
\widetilde{\frak{ks}}_{\frak b}({\overline{\text{\bf f}}}_{a(2)})
-  \sum_{c=0}^{B'} c^c_{a(1) a(2)}&\widetilde{\frak{ks}}_{\frak b}(\overline{\text{\bf f}}_c)
\equiv 0 
\\
&\mod
\text{\rm Clos}_{d_{\overset{\circ}P}}\left( y_i\frac{\partial \frak{\frak PO}_{\frak b}}
{\partial y_i}
: i = 1,\dots,n\right).
\endaligned
\end{equation}
Here
$\widetilde{\frak{ks}}_{\frak b}(\overline{\text{\bf f}}_{a(1)})
\widetilde{\frak{ks}}_{\frak b}({\overline{\text{\bf f}}}_{a(2)})$
denotes multiplication as the elements of $\Lambda
\langle\!\langle y,y^{-1}\rangle\!\rangle_0^{\overset{\circ}P}$.
\par
We decompose $\frak b = \frak b_2 + \frak b_{\rm high}$
where $\frak b_2$ is degree $2$ and $\frak b_{\rm high}$ is
degree $\ge 4$. (We do not need to discuss degree $0$ component since
it does not affect Jacobian ring.)
\par
Let $\mathcal M^{\text{main}}_{1;\ell+2}(\beta)$
be the moduli space of stable maps from genus zero
semi-stable curve with $1$ boundary and
$\ell+2$ interior marked points and of homology class $\beta$.
We denote the fibre product
$$
\mathcal M^{\text{main}}_{1;\ell+2}(\beta) \times_{X^{\ell+2}}
(\overline{\text{\bf f}}_{a(1)} \times {\overline{\text{\bf f}}}_{a(2)} \times
\frak b_{\rm high}^{\ell})
$$
by $\mathcal M^{\text{main}}_{1;\ell+2}(\beta;\overline{\text{\bf f}}_{a(1)} \otimes
{\overline{\text{\bf f}}}_{a(2)} \otimes \frak b_{\rm high}^{\otimes \ell})$.
We recall from  \cite[Lemma 6.5]{fooo09} that we have already made a choice of
a multisection thereon that is $T^n$-equivariant and transversal.
We call this Kuranishi structure and multisection the $\frak q$-multisection.
(They are used to define the operator $\frak q$.)
\par
We consider $\mathcal M^{\text{main}}_{1;\ell+2}(\beta;\overline{\text{\bf f}}_{a(1)} \otimes
{\overline{\text{\bf f}}}_{a(2)} \otimes \frak b_{\rm high}^{\otimes \ell})$ and
focus on its components of (virtual) dimension $n$: Recall that
$\frak b_{\rm high}$ lies in $H(X;\Lambda_0)$ and is a formal power series
whose summands may have different degrees. For the simplicity of
exposition, we will just denote union of these components by
$\mathcal M^{\text{main}}_{1;\ell+2}(\beta;\overline{\text{\bf f}}_{a(1)} \otimes
{\overline{\text{\bf f}}}_{a(2)} \otimes \frak b_{\rm high}^{\otimes \ell})$
as long as there is no danger of confusion.
\par
We consider the forgetful map
\begin{equation}\label{forgetmap}
\frak{forget}:
\mathcal M^{\text{main}}_{1;\ell+2}(\beta;\overline{\text{\bf f}}_{a(1)} \otimes
{\overline{\text{\bf f}}}_{a(2)} \otimes \frak b_{\rm high}^{\otimes \ell})
\to \mathcal M^{\text{\rm main}}_{1;2}.
\end{equation}
This is a composition of the forgetful map
$\mathcal M^{\text{main}}_{1;\ell+2}(\beta;\overline{\text{\bf f}}_{a(1)} \otimes
{\overline{\text{\bf f}}}_{a(2)} \otimes \frak b_{\rm high}^{\otimes \ell})
\to\mathcal M^{\text{\rm main}}_{1;\ell + 2}$
which forgets maps and then shrinks the resulting unstable components if any,
followed by the forgetful map
$\mathcal M^{\text{\rm main}}_{1;\ell + 2} \to \mathcal M^{\text{\rm main}}_{1;2}$,
forgetting 3rd,\dots,$(\ell+1)$-th marked points.
\par
We remark that $\mathcal M^{\text{\rm main}}_{1;2}$ is homeomorphic
to the disk $D^2$. In fact, $\mathcal M^{\text{\rm main}}_{1;1}$ is a
one point and the fiber of the forgetful map $\mathcal M^{\text{\rm
main}}_{1;2} \to \mathcal M^{\text{\rm main}}_{1;1}$ is a disk,
which is parameterized by the position of the last interior marked
point we forget.
\begin{lem}\label{stratifyD12}
$\mathcal M^{\text{\rm main}}_{1;2}$ has a stratification
such that the interior of each stratum is described as follows.
\begin{enumerate}
\item $\text{\rm Int}D^2 \setminus \{0\}$.
\item Two copies of $(-1,1)$. We call them  $(-1,1)_1$,
$(-1,1)_2$.
\item Three points. $[\Sigma_0]$, $[\Sigma_{12}]$, $[\Sigma_{21}]$.
\end{enumerate}
The point $0 \in D^2$ becomes $[\Sigma_0]$
in the closure of the stratum $(1)$.
The two boundary points $\pm 1$ in the closure of $(-1,1)_i$
$(i=1,2)$ become $[\Sigma_{12}]$ and $[\Sigma_{21}]$.
Thus the union of the four strata $(-1,1)_1$, $(-1,1)_2$,
$[\Sigma_{12}]$, $[\Sigma_{21}]$  is a circle, which becomes $\partial D^2$
in the closure of the stratum $(1)$.
\par
Each of the above strata corresponds to the combinatorial type
of the elements of $\mathcal M^{\text{\rm main}}_{1;2}$
\end{lem}
\begin{proof}
In addition to the boundary, $\mathcal M^{\text{\rm main}}_{1;2}$
has one singular point. That is the point where two interior marked points encounter.
We denote this point by $[\Sigma_0]$. (Here $\Sigma_0$ is
the corresponding stable bordered curve.)
\par
The boundary of $\mathcal M^{\text{\rm main}}_{1;2}$
is identified with two copies of $\mathcal M^{\text{\rm main}}_{2;1}$.
They correspond to the configuration $\Sigma
= \Sigma' \cup \Sigma''$ where $[\Sigma'] \in \mathcal M^{\text{\rm main}}_{1;1}$ and $[\Sigma''] \in \mathcal M^{\text{\rm main}}_{2;1}$.
The two components $\Sigma'$ and $\Sigma''$ are
glued at their boundary marked points.
There are two components depending which interior marked
points of  $\mathcal M^{\text{\rm main}}_{1;2}$ becomes
the interior marked point of $\Sigma'$.
\par
$\mathcal M^{\text{\rm main}}_{1;1}$ is one point.
We see that $\mathcal M^{\text{\rm main}}_{2;1}$
is an interval as follows:
We can normalize an element of it $(z_0,z_1^{\text{\rm int}},z_2^{\text{\rm int}})$
so the boundary marked point is $+1=z_0$ and
the real part of the interior marked points $z_i^{\text{\rm int}}$
is zero. Moreover we may choose $z_1^{\text{\rm int}} + z_2^{\text{\rm int}}=0$.
Thus $\mathcal M^{\text{\rm main}}_{2;1}$
is identified with
$$\{(z_1^{\text{\rm int}},z_2^{\text{\rm int}}) \mid
z_1^{\text{\rm int}} = \sqrt{-1}a,~ z_2^{\text{\rm int}} = -\sqrt{-1}a,~ a\in [-1,1]\} \cong [-1,1].$$
\par
We write as $(-1,1)_i$ the interior of the component
such that $i$-th
interior marked
point becomes the interior marked point of $\mathcal M^{\text{\rm main}}_{1;1}$.
\par
By the above identification, $+1$ of $\partial[-1,1]_2$
corresponds to $-1$ of $\partial[-1,1]_1$.
We write this point $\Sigma_{12}$.
The other point is $\Sigma_{21}$.
\end{proof}
\par
\hskip0.2cm
\epsfbox{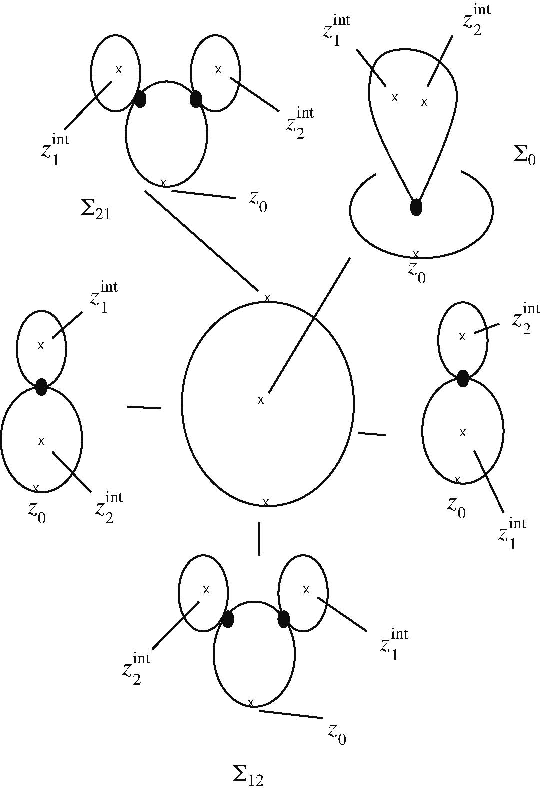}
\par
\centerline{\bf Figure 2.6.1}
\par
The domain of (\ref{forgetmap}) also  has a stratification according to the
combinatorial type of the source.
The next lemma is obvious from construction.
\begin{lem}\label{smoothmaptoDM} The map
$(\ref{forgetmap})$ is continuous and is a stratawise smooth submersion.
\end{lem}
\begin{rem}
When we use the smooth structure we put on
$\mathcal M^{\text{main}}_{1;\ell+2}(\beta;\overline{\text{\bf f}}_{a(1)} \otimes
{\overline{\text{\bf f}}}_{a(2)} \otimes \frak b_{\rm high}^{\otimes \ell})$
in  \cite[Section A1.4]{fooobook2}, then $(\ref{forgetmap})$ will not be
smooth in the normal direction to the strata.
See \cite[pages 777-778]{fooobook2}.
\end{rem}
We consider the preimages
$$
\aligned
&\frak{forget}^{-1}([\Sigma_0]) \subset \mathcal M^{\text{main}}_{1;\ell+2}(\beta;\overline{\text{\bf f}}_{a(1)} \otimes
{\overline{\text{\bf f}}}_{a(2)} \otimes \frak b_{\rm high}^{\otimes \ell}),\\
&\frak{forget}^{-1}([\Sigma_{12}])
\subset \mathcal M^{\text{main}}_{1;\ell+2}(\beta;\overline{\text{\bf f}}_{a(1)} \otimes
{\overline{\text{\bf f}}}_{a(2)} \otimes \frak b_{\rm high}^{\otimes \ell}).
\endaligned
$$
The proof of Theorem \ref{multiplicative} consists of the following 3 steps.
\begin{enumerate}
\item[Step I:]
In this step we show that the weighted sum of 
${\rm ev}_{0 *}(\frak{forget}^{-1}([\Sigma_0]))$ 
for various 
\begin{equation}\label{variousev0}
{\rm ev}_0 : \mathcal M^{\text{main}}_{1;\ell+2}(\beta;\overline{\text{\bf f}}_{a(1)} \otimes
{\overline{\text{\bf f}}}_{a(2)} \otimes \frak b_{\rm high}^{\otimes \ell})
\to L(u)
\end{equation}
with appropriate weight coincides with 
$
\sum_{c=0}^{B'} c^c_{a(1) a(2)}\widetilde{\frak{ks}}_{\frak b}(\overline{\text{\bf f}}_c)
$
modulo an element of the Jacobian ideal $\text{\rm Clos}_{d_{\overset{\circ}P}}\left( y_i\frac{\partial \frak{\frak PO}_{\frak b}}
{\partial y_i}
: i = 1,\dots,n\right)$.
\item[Step II:]
In this step we show that an appropriate weighted sum of ${\rm ev}_{0 *}(\frak{forget}^{-1}([\Sigma_0]))$ 
for various ${\rm ev}_0$ as in 
(\ref{variousev0})
coincides with an appropriate weighted sum of ${\rm ev}_{0 *}(\frak{forget}^{-1}([\Sigma_{12}]))$
for various 
${\rm ev}_0$ as in 
(\ref{variousev0})
modulo an element of the Jacobian ideal $\text{\rm Clos}_{d_{\overset{\circ}P}}\left( y_i\frac{\partial \frak{\frak PO}_{\frak b}}
{\partial y_i}
: i = 1,\dots,n\right)$.
\item[Step III:]
In this step we show that an appropriate weighted sum of ${\rm ev}_{0 *}(\frak{forget}^{-1}([\Sigma_{12}]))$ 
for various ${\rm ev}_0$ as in 
(\ref{variousev0})
coincides with $\widetilde{\frak{ks}}_{\frak b}(\overline{\text{\bf f}}_{a(1)})
\widetilde{\frak{ks}}_{\frak b}({\overline{\text{\bf f}}}_{a(2)})$ 
modulo  an element of the Jacobian ideal  $\text{\rm Clos}_{d_{\overset{\circ}P}}\left( y_i\frac{\partial \frak{\frak PO}_{\frak b}}
{\partial y_i}
: i = 1,\dots,n\right)$.
\end{enumerate}
\par
\hskip2.2cm
\epsfbox{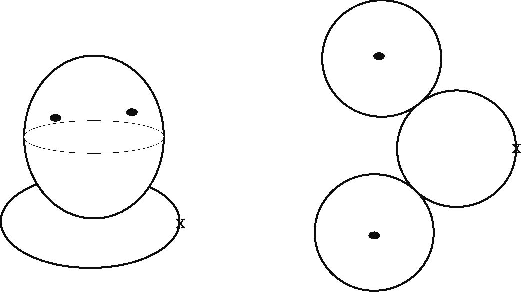}
\par
\centerline{\bf Figure 2.6.2}
\par
Note the geometric intuition behind the proof of Steps I and III is clear from 
the Figure 2.6.2 above.
The main part of its proof consists of describing the way how we perturb the 
moduli spaces involved so that the operators obtained after perturbation 
indeed satisfies the equalities expected from the geometric intuition.
\par
Intuitively Step II is a consequence of the standard cobordism argument.
Namely it follows from the fact that $\frak{forget}^{-1}([\Sigma_0])$ 
is cobordant to $\frak{forget}^{-1}([\Sigma_{12}])$.
There are two points we need to take care of to work out this cobordism argument:
\par
One is describing the way how we perturb the moduli spaces so that 
the perturbation is consistent with the cobordism argument.
\par
The other is to go around the problem that the map $\frak{forget}$ is only 
stratawise smooth. 
In fact the two points $[\Sigma_0]$, $[\Sigma_{12}]$ consist of single 
strata each. So the smoothness of $\frak{forget}$ breaks down there.
\par
The rest of this section is occupied with a detailed discussion on those points.
\par
Our first task is to provide appropriate multisections on certain neighborhoods
of $\frak{forget}^{-1}([\Sigma_0])$ and $\frak{forget}^{-1}([\Sigma_{12}])$ respectively. For this
purpose, we now describe the stratification of our space in these neighborhoods in detail.

For $\alpha\in H_2(X;\Z)$, let $\mathcal M_{\ell}(\alpha)$
be the moduli space of stable maps from genus zero
closed Riemann surface with $\ell$ marked points and
of homology class $\alpha$.
Let
$$
\text{\rm ev} = (\text{\rm ev}_1,\dots,\text{\rm ev}_{\ell}):
\mathcal M_{\ell}(\alpha) \to X^{\ell}
$$
be the evaluation map. We define
\begin{equation}
\aligned
&\mathcal M_{\ell_1+3}(\alpha; \overline{\text{\bf f}}_{a(1)}
\otimes  {\overline{\text{\bf f}}}_{a(2)}\otimes \frak b_{\rm high}^{\otimes\ell_1})\\
&=
\mathcal M_{\ell_1+3}(\alpha) 
{}_{(\text{\rm ev}_1,\dots,\text{\rm ev}_{\ell_1+2})}\times_{X^{\ell_1+2}}
(\overline{\text{\bf f}}_{a(1)}
\times  {\overline{\text{\bf f}}}_{a(2)}\times \frak b_{\rm high}^{\ell_1}).
\endaligned\end{equation}
$\text{\rm ev}_{\ell_1+3}$ defines an evaluation map
$$
\text{\rm ev}_{\ell_1+3}:
\mathcal M_{\ell_1+3}(\alpha; \overline{\text{\bf f}}_{a(1)}
\otimes  {\overline{\text{\bf f}}}_{a(2)}\otimes \frak b_{\rm high}^{\otimes\ell_1})
\to X.
$$
We next consider $\mathcal M_{k+1;\ell_2}^{\text{\rm main}}(\beta)$,
that is the moduli space of bordered genus zero stable
map with $k+1$ boundary and $\ell_2$ interior marked points and of
homology class $\beta \in H_2(X;L(u);\Z)$.
We consider
$$
\mathcal M_{k+1;\ell_2+1}^{\text{\rm main}}(\beta;\frak b_{\rm high}^{\otimes\ell_2})
=
\mathcal M_{k+1;\ell_2+1}^{\text{\rm main}}(\beta) 
{}_{(\text{\rm ev}^{\text{\rm int}}_1,\dots,\text{\rm ev}^{\text{\rm int}}_{\ell_2})}\times_{X^{\ell_2}} \frak b_{\rm high}^{\ell_2}
$$
and the evaluation map
$$
\text{\rm ev} = (\text{\rm ev}^{\text{\rm int}}_{\ell_2+1},\text{\rm ev}_0,\text{\rm ev}_1,\dots,\text{\rm ev}_k)
: \mathcal M_{k+1;\ell_2+1}^{\text{\rm main}}(\beta;\frak b_{\rm high}^{\otimes\ell_2})
\to X \times L(u)^{k+1}.
$$
\begin{defn}\label{defn96}
Let $K$ be a smooth singular chain of $X^2$.
We define
$$
\mathcal M_{k+1}(\alpha;\beta;\overline{\text{\bf f}}_{a(1)},
{\overline{\text{\bf f}}}_{a(2)};\ell_1,\ell_2;K)
$$
to be the fiber product (over $X^2$) of
\begin{equation}\label{katawarefiber}
\mathcal M_{\ell_1+3}(\alpha; \overline{\text{\bf f}}_{a(1)}
\otimes  {\overline{\text{\bf f}}}_{a(2)}\otimes \frak b_{\rm high}^{\otimes\ell_1})
\times
\mathcal M_{k+1;\ell_2+1}^{\text{\rm main}}(\beta;\frak b_{\rm high}^{\otimes\ell_2})
\end{equation}
with
$K$
by the map $({\text{\rm ev}_{\ell_1+3},\text{\rm ev}^{\text{\rm int}}_{\ell_2+1}})$.
Using the evaluation maps at the boundary marked points of
the second factor of (\ref{katawarefiber}), we have an
evaluation map
$$\text{\rm ev} = (\text{\rm ev}_0,\dots,\text{\rm ev}_{k}):
\mathcal M_{k+1}(\alpha;\beta;\overline{\text{\bf f}}_{a(1)};
{\overline{\text{\bf f}}}_{a(2)};\ell_1,\ell_2;K)  \to L(u)^{k+1}.
$$
\end{defn}
If $K$ is a manifold with corner, then
$\mathcal M_{k+1}(\alpha;\beta;\overline{\text{\bf f}}_{a(1)},{\overline{\text{\bf f}}}_{a(2)};\ell_1,\ell_2;K)$ has a
Kuranishi structure. 
Note the fiber product we take to define 
$
\mathcal M_{k+1}(\alpha;\beta;\overline{\text{\bf f}}_{a(1)},
{\overline{\text{\bf f}}}_{a(2)};\ell_1,\ell_2;K)
$
is a particular case of fiber product between two spaces with 
Kuranishi structures, since a manifold with corner can 
be regarded as a space with Kuranishi structure (with trivial obstruction bundle).
Such fiber product is defined for example in \cite[Definition 4.9 (2)]{tech2}.
In general, when $K$ has a triangulation $K = \cup K_i$,
we can equip
$\mathcal M_{k+1}(\alpha;\beta;\overline{\text{\bf f}}_{a(1)},{\overline{\text{\bf f}}}_{a(2)};\ell_1,\ell_2;K_i)$
with  Kuranishi structures for each of the singular simplices $K_i$
so that they coincide on the overlapped parts.
\begin{rem}
Kuranishi space with corners are naturally stratified by the canonical stratification of manifold with corners.
In case we take the fiber product of $X$ with a simplex  $\Delta$, and the fiber product with $\Delta_i$ (its faces), 
the Kuranishi space that is a fiber product of $X$ with  $\Delta_i$ is a closure of a stratum of  
the fiber product of $X$ with $\Delta$.
We can formulate compatibility of the structures by using this fact.
\end{rem}
By an abuse of notation, we simply say that
$\mathcal M_{k+1}(\alpha;\beta;\overline{\text{\bf f}}_{a(1)},{\overline{\text{\bf f}}}_{a(2)};\ell_1,\ell_2;K)$ has
a Kuranishi structure. 
\begin{lem}\label{boudaryMippai}
As a space with Kuranishi structure, the boundary of the moduli space
$\mathcal M_{k+1}(\alpha;\beta;\overline{\text{\bf f}}_{a(1)},{\overline{\text{\bf f}}}_{a(2)};\ell_1,\ell_2;K)$
is the union of the spaces of one of the following three types:
\begin{enumerate}
\item
$\mathcal M_{k+1}(\alpha;\beta;\overline{\text{\bf f}}_{a(1)},{\overline{\text{\bf f}}}_{a(2)};\ell_1,\ell_2;\partial K)$.
\item
$
\mathcal M_{k'+1}(\alpha;\beta';\overline{\text{\bf f}}_{a(1)},{\overline{\text{\bf f}}}_{a(2)};\ell_1,\ell'_2;K)
_{\text{\rm ev}_0}\times_{\text{\rm ev}_i}
\mathcal M_{k''+1;\ell''_{2}+1}^{\text{\rm main}}
(\beta'';\frak b_{\rm high}^{\otimes\ell''_2}),
$
\par
where $k'+k'' = k+1$, $\ell'_2+\ell''_2 = \ell_2$,
$\beta'\#\beta'' = \beta$,
$i \in \{1,\dots,k''\}$.
We put $\frak q$-Kuranishi structure and multisection on the second factor.
\item
$
\mathcal M_{k'+1;\ell'_{2}+1}^{\text{\rm main}}
(\beta';\frak b_{\rm high}^{\otimes\ell'_2})
_{\text{\rm ev}_0}\times_{\text{\rm ev}_i}
\mathcal M_{k''+1}(\alpha;\beta'';\overline{\text{\bf f}}_{a(1)},{\overline{\text{\bf f}}}_{a(2)};\ell_1,\ell''_2;K),
$
\par
where $k'+k'' = k+1$, $\ell'_2+\ell''_2 = \ell_2$,
$\beta'\#\beta'' = \beta$,
$i \in \{1,\dots,k''\}$.
We put $\frak q$-Kuranishi structure and multisection on the first factor.
\end{enumerate}
\end{lem}
\begin{proof}
Our space $\mathcal M_{k+1}(\alpha;\beta;\ell_1,\ell_2;K)$
is defined by taking fiber product. Therefore
its boundary is described by taking the boundary
of one of the factors.
The moduli space of stable maps from a {\it closed}
curve has no boundary. The boundary of
$\mathcal M_{k+1;\ell_2}^{\text{\rm main}}(\beta)$ is union of
$\mathcal M_{k'+1;\ell'_2}^{\text{\rm main}}(\beta')
_{\text{\rm ev}_0}\times_{\text{\rm ev}_i} \mathcal M_{k''+1;\ell''_2}^{\text{\rm main}}(\beta'')$
where $k'+k'' = k+1$, $\ell'_2+\ell''_2 = \ell_2$,
$\beta'\#\beta'' = \beta$,
$i \in \{1,\dots,k''\}$. Hence the lemma.
\end{proof}
We now consider the case $K = \Delta$ which is the diagonal
$$
\Delta = \{(x,x) \mid x\in X\} \subset X^2,
$$
and the moduli space $\mathcal M_{k+1}(\alpha;\beta;\overline{\text{\bf f}}_{a(1)},{\overline{\text{\bf f}}}_{a(2)};\ell_1,\ell_2;\Delta)$.
\begin{lem}\label{deltaidentify}
There exists a surjective map
\begin{equation}\label{gluefromunion}
\text{\rm Glue}:
\bigcup_{\alpha\# \beta'=\beta}
\bigcup_{\ell_1+\ell_2=\ell}\mathcal M_{k+1}(\alpha;\beta';\overline{\text{\bf f}}_{a(1)},
{\overline{\text{\bf f}}}_{a(2)};\ell_1,\ell_2;\Delta)
\to \frak{forget}^{-1}([\Sigma_0])
\end{equation}
which defines an isomorphism outside codimension $2$ strata as a space with
Kuranishi structure. Here we have
$$
\frak{forget}^{-1}([\Sigma_0])\subset \mathcal M^{\text{\rm main}}_{k+1;\ell_1+\ell_2+2}(\beta'\#\alpha;\overline{\text{\bf f}}_{a(1)} \otimes
{\overline{\text{\bf f}}}_{a(2)} \otimes \frak b_{\rm high}^{\otimes \ell}).
$$
\end{lem}
\begin{proof}
Let
$$
\aligned
(w_1;\Sigma^{\text{\rm int}};\vec z^{\rm int}) &\in
\mathcal M_{\ell_1+3}(\alpha; \overline{\text{\bf f}}_{a(1)}
\otimes  {\overline{\text{\bf f}}}_{a(2)}\otimes \frak b_{\rm high}^{\otimes\ell_1}),
\\
(w_2;\Sigma^{\text{\rm bdy}};\vec z',\vec z^{\prime \rm int}) &\in
\mathcal M_{k+1;\ell_2+1}^{\text{\rm main}}(\beta';\frak b_{\rm high}^{\otimes\ell_2}).
\endaligned$$
We assume that $\text{\rm ev}_{\ell_1+3}(w_1;\Sigma^{\text{\rm int}};\vec z^{\rm int})
= \text{\rm ev}^{\text{\rm int}}_{\ell_2+1}(w_2;\Sigma^{\text{\rm bdy}};\vec z',\vec z^{\prime,\rm int})$.
We then glue $\Sigma^{\text{\rm int}}$ and $\Sigma^{\text{\rm bdy}}$
at $z^{\text{\rm int}}_{\ell_1+3} \in \Sigma^{\text{\rm int}}$ and
$z^{\prime \text{\rm int}}_{\ell_2+1} \in \Sigma^{\text{\rm bdy}}$
to obtain $\Sigma$.
\par
We put $w = w_1$ on $\Sigma^{\text{\rm int}}$ and
$w = w_2$ on $\Sigma^{\text{\rm bdy}}$. We then obtain
$w: \Sigma \to X$.
The marked points on $\Sigma$ are determined from those on
$\Sigma^{\text{\rm int}}$ and $\Sigma^{\text{\rm bdy}}$
in an obvious way. We have thus defined
\begin{equation}\label{componentglue}
\text{\rm Glue}:
\mathcal M_{k+1}(\alpha;\beta';\overline{\text{\bf f}}_{a(1)},
{\overline{\text{\bf f}}}_{a(2)};\ell_1,\ell_2;\Delta)
\to \frak{forget}^{-1}([\Sigma_0]).
\end{equation}
On such a component, the map $\text{\rm Glue}$ is an isomorphism
onto its image, as a space with Kuranishi structure.
In fact, the element $(w;\Sigma;\vec z,\vec z^{\rm int})$ in
$\frak{forget}^{-1}([\Sigma_0])$ is in the image of  {\it unique}
point by the map (\ref{gluefromunion}) if the
following two conditions are satisfied.
\begin{conds}\label{1tai1cond}
\begin{enumerate}
\item Either the $1$st or the $2$nd marked point of $\Sigma^{\text{\rm int}}$
lies in the same irreducible component as the $(\ell_1+3)$-th interior marked point.
\item
The $(\ell_2+1)$-th interior marked point of  $\Sigma^{\text{\rm bdy}}$
is on a disk component, (that is not on the sphere bubble).
\end{enumerate}
\end{conds}
Clearly the set of points which do not satisfy Condition
\ref{1tai1cond} is of codimension 2 or greater.
Hence the lemma.
\end{proof}
\begin{exm}\label{doublepointexam}
An example for which Condition
\ref{1tai1cond} is not satisfied can be described as follows:
(See Figure 2.6.3.) 
Consider $\Sigma$ which has one
disk component $D^2$ and two sphere components $S^2_1$,
$S^2_2$. We assume $D^2 \cap S^2_1 =\{p_1\} = $ one point,
$D^2 \cap S^2_2 = \emptyset$, $S^2_1 \cap S^2_2 =\{p_2\} = $ one point.
Moreover we assume the first and 2nd interior marked points
$z_1^{\text{\rm int}}$, $z_2^{\text{\rm int}}$ are on
$S^2_2$. Let $w: (\Sigma,\partial\Sigma) \to (X,L(u))$
be holomorphic. Then we have
$$
(w;\Sigma;z_0,z_1^{\text{\rm int}},z_2^{\text{\rm int}})
\in \frak{forget}^{-1}([\Sigma_0])
\subset \mathcal M_{1;2}^{\text{\rm main}}(\beta).
$$
We put $\beta_{(0)} = w_*([D^2])$, $\alpha_1 = w_*[S^2_1]$,
$\alpha_2 = w_*[S^2_2]$. ($\beta = \beta_{(0)} \# \alpha_1 \#
\alpha_2$.) Then we have
$$
((w;S^2_1\#S^2_2;z_1^{\text{\rm int}},z_2^{\text{\rm int}},p_1),
(w;D^2;z_0,p_1))
\in \mathcal M_{1}(\alpha_1\#\alpha_2;\beta_{(0)};\overline{\text{\bf f}}_{a(1)},
{\overline{\text{\bf f}}}_{a(2)};\ell_1,\ell_2;\Delta)
$$
On the other hand,
$$
((w;S^2_2;z_1^{\text{\rm int}},z_2^{\text{\rm int}},p_2),
(w;D^2\#S^2_1;z_0,p_2))
\in \mathcal M_{1}(\alpha_2;\alpha_1\#\beta_{(0)};\overline{\text{\bf f}}_{a(1)},
{\overline{\text{\bf f}}}_{a(2)};\ell_1,\ell_2;\Delta).
$$
Both of them go to
$(w;\Sigma;z_0,z_1^{\text{\rm int}},z_2^{\text{\rm int}})
\in \frak{forget}^{-1}([\Sigma_0])$
by the map $\text{\rm Glue}$.
\end{exm}
\par
\hskip4.2cm
\epsfbox{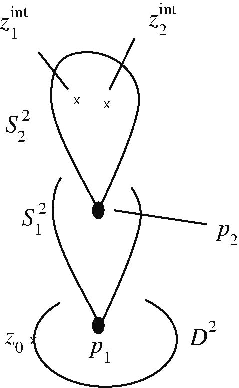}
\par
\centerline{\bf Figure 2.6.3}
\par
Let us describe $\frak{forget}^{-1}([\Sigma_0])$ in more detail.
We consider a point $(w;\Sigma)$ on it.
(For simplicity we omit marked points from our notation.)
We define the decomposition
\begin{equation}\label{componentdecompose}
\Sigma = \Sigma^{\text{\rm int}}
\cup \Sigma^{\text{\rm glue}}
\cup \Sigma^{\text{\rm bdy}}.
\end{equation}
Let $\Sigma^{\text{\rm int},0}$ be the smallest
union of components of $\Sigma$ that is connected
and contains the first and second interior marked points.
Let $\Sigma^{\text{\rm bdy},0}$ be the smallest
union of components of $\Sigma$ that is connected
and contains all the disk components.
Let $\Sigma^{\text{\rm glue},0}$ be the union of sphere
components in $\Sigma \setminus \Sigma^{\text{\rm int},0
} \setminus \Sigma^{\text{\rm bdy},0}$ such that the following holds:
$$
\text{All the path joining $\Sigma^{\text{\rm int},0}$
with $\Sigma^{\text{\rm bdy},0}$ should intersect this
sphere component.}
$$
The union $\Sigma^{\text{\rm int},0}
\cup \Sigma^{\text{\rm glue},0} \cup \Sigma^{\text{\rm bdy},0}$
is connected.
The complement of this union consists of disjoint union of
trees of sphere components. We include those trees to
one of $\Sigma^{\text{\rm int},0}$,
$\Sigma^{\text{\rm glue},0}$, $\Sigma^{\text{\rm bdy},0}$
according to which it is rooted.
We have thus obtained $\Sigma^{\text{\rm int}}$,
$\Sigma^{\text{\rm glue}}$, $\Sigma^{\text{\rm bdy}}$
and the decomposition (\ref{componentdecompose}).
\par
\begin{exm}
In the case of $\Sigma$ in  Example \ref{doublepointexam}, we have
$\Sigma^{\text{\rm int}} = S^2_2$,
$\Sigma^{\text{\rm glue}} = S^2_1$,
$\Sigma^{\text{\rm bdy}} = D^2$.
\end{exm}
\begin{lem}\label{normalcrossing}
$(w;\Sigma)$ satisfies Condition \ref{1tai1cond} if and only
if $\Sigma^{\text{\rm glue}}$ is empty.
\par
Suppose $\Sigma^{\text{\rm glue},0}$ is non-empty and
contains $N$ sphere components. Then the
cardinality of $\text{\rm Glue}^{-1}(w;\Sigma)$ is $N+1$.
\par
Moreover the map  $\text{\rm Glue}$ in the neighborhood of
$\text{\bf x} \in \frak{forget}^{-1}([\Sigma_0])$ is described as follows:
There is a Kuranishi neighborhood $U$ of $\text{\bf x}$
in $\mathcal M^{\text{\rm main}}_{k+1;\ell_1+\ell_2+2}(\beta;\overline{\text{\bf f}}_{a(1)} \otimes
{\overline{\text{\bf f}}}_{a(2)} \otimes \frak b_{\rm high}^{\otimes \ell})$
such that
$$
U = \overline U \times \C^{N+1}
$$
where $\text{\bf x} = (o,0)$.
\par
We remark that $\text{\bf x}$ is an element of the intersection
$$
\bigcap_{i=0}^N
\text{\rm Glue} \left(\mathcal M_{k+1}(\alpha_i;\beta_i;\overline{\text{\bf f}}_{a(1)},
{\overline{\text{\bf f}}}_{a(2)};\ell_{i,1},\ell_{i,2};\Delta)\right)
$$
if $\text{\bf x} = \text{\rm Glue}(\text{\bf x}_i)$ with
$\text{\bf x}_i \in \mathcal M_{k+1}(\alpha_i;\beta_i;\overline{\text{\bf f}}_{a(1)},
{\overline{\text{\bf f}}}_{a(2)};\ell_{i,1},\ell_{i,2};\Delta)$.
A Kuranishi neighborhood $U_i$ of $\text{\bf x}_i$ is identified with
$$
U_i = \overline U \times (\C^{i}\times \{0\} \times \C^{N-i}),
$$
where $U_i$ is embedded in $U$ by the obvious inclusion.
\par
The space
$\frak{forget}^{-1}([\Sigma_0])$ is described in our
Kuranishi neighborhood by
$$
\left(\left(\bigcup_{i=0}^{N} U_i\right)  \cap s^{-1}_{\text{\bf x}}(0) \right)/
\Gamma_{\text{\bf x}}.
$$
Here $s_{\text{\bf x}}$ is the Kuranishi map and $\Gamma_{\text{\bf x}}$
is the finite group of symmetry of $(w;\Sigma)$.
\end{lem}
\begin{rem}
We can summarize Lemma \ref{normalcrossing} as follows:
Our space
$\frak{forget}^{-1}([\Sigma_0])$ is a `normal crossing divisor' and
$N+1$ components intersect transversally at $\text{\bf x}$.
(We put the phrase `normal crossing divisor' in quote since
it comes with additional data of a Kuranishi map and an obstruction bundle.)
\end{rem}
\begin{proof}
Let $S^2_i$ ($i=1,\dots,N$) be the component of $\Sigma^{\text{\rm glue},0}$.
We include each of the tree of components in
$\Sigma^{\text{\rm glue}} \setminus \Sigma^{\text{\rm glue},0}$
to the component $S^2_i$ where this tree is rooted.
We then have a decomposition
$$
\Sigma^{\text{\rm glue}}
= \bigcup_{i=1}^N \Sigma^{\text{\rm glue}}_i.
$$
We may enumerate them so that
$$\aligned
&\Sigma^{\text{\rm glue}}_1 \cap \Sigma^{\text{\rm bdy}}
= \{z_0\}, \\
&\Sigma^{\text{\rm glue}}_i \cap \Sigma^{\text{\rm glue}}_{i+1}
= \{z_i\}, \quad i = 1,\dots,N-1\\
&\Sigma^{\text{\rm glue}}_N \cap \Sigma^{\text{\rm int}}
= \{z_N\}.
\endaligned$$
For $i = 0,\dots,N$, we put
$$
\beta_{(i)} = w_*([\Sigma^{\text{\rm bdy}} ]) + \sum_{j=1}^{i}
w_*([\Sigma^{\text{\rm glue}}_j]),
$$
$$
\alpha_i = \sum_{j=i+1}^{N}
w_*([\Sigma^{\text{\rm glue}}_j]) + w_*([\Sigma^{\text{\rm int}}]).
$$
If we cut $\Sigma$ at $z_i$, then we obtain an element of
$\mathcal M_{k+1}(\alpha_i;\beta_i;\overline{\text{\bf f}}_{a(1)},
{\overline{\text{\bf f}}}_{a(2)};\ell_{i,1},\ell_{i,2};\Delta)$.
It is clear that those elements are all the elements of
$\text{\rm Glue}^{-1}(w;\Sigma)$.
\par
The factor $\overline U$ is the product of the Kuranishi
neighborhoods of $(w;\Sigma^{\text{\rm bdy}})$,
$(w;\Sigma^{\text{\rm glue}}_i)$,
$(w;\Sigma^{\text{\rm int}})$ in the moduli spaces in which they are
contained.
Other than $\overline U$, the Kuranishi neighborhood of
$\text{\bf x}$ in
$\mathcal M^{\text{main}}_{k+1;\ell_1+\ell_2+2}(\beta;\overline{\text{\bf f}}_{a(1)} \otimes
{\overline{\text{\bf f}}}_{a(2)} \otimes \frak b_{\rm high}^{\otimes \ell})$
has a parameter. That is the parameter to glue
$N+2$ components at the singular points $z_0,\dots,z_N$.
They give the factor $\C^{N+1}$.
\par
If we resolve all the singularities $z_j$ ($j=0,\dots,N$) other than
$z_i$, we still obtain an element of $\frak{forget}^{-1}([\Sigma_0])$.
Moreover it is in the image of
$\mathcal M_{k+1}(\alpha_i;\beta_i;\overline{\text{\bf f}}_{a(1)},
{\overline{\text{\bf f}}}_{a(2)};\ell_{i,1},\ell_{i,2};\Delta)$.
Therefore the lemma follows.
\par
(We remark that our Kuranishi structure is constructed
inductively over the strata  in the way we have described above.
See \cite{FO},\cite{fooo06}.)
\end{proof}
\begin{rem}\label{rem:1014}
It is difficult to find a $T^n$-equivariant multisection that is  
compatible with
the fiber product description in Lemmata
\ref{boudaryMippai}, \ref{deltaidentify}, and \ref{normalcrossing}.
See \cite[Remark 11.4]{fooo08} .
It seems that it is possible to take the Kuranishi structure that is $T^n$-equivariant and 
is compatible with
the fiber product description. However we do not prove it in Section \ref{sec:cyclicKura}.
The proof seems to be slightly more cumbersome than the one in Section \ref{sec:cyclicKura}.
So we do not use it.
\end{rem}
\par
To go around the problem mentioned in Remark \ref{rem:1014}, we use 
a Kuranishi structure and multisection on 
a neighborhood of 
$\frak{forget}^{-1}([\Sigma_0])$ in $\mathcal M^{\text{\rm main}}_{k+1;\ell_1+\ell_2+2}(\beta'\#\alpha;\overline{\text{\bf f}}_{a(1)} \otimes
{\overline{\text{\bf f}}}_{a(2)} \otimes \frak b_{\rm high}^{\otimes \ell})$
that is different from the $\frak q$-multisection
and is {\it not} $T^n$-equivariant. 
On the other hand, since we are working on de Rham theory, we need to take
integration along the fiber by the map $\text{\rm ev}_0$. 
Since our multisection is not $T^n$ equivariant, we can not use the trick that 
$T^n$-equivariant map to a $T^n$ orbit is automatically a submersion.
So to achieve submersivity of the evaluation map, we use the technique of continuous family of
multisections. 
(We discussed the technique of continuous family of multisections and 
its application to define smooth correspondence in de Rham theory  in detail in  
\cite[Section 12]{fooo09}.)
We denote this family of multisections by $\frak s$.
We describe the property of $\frak s$ in the next lemma.
\par
Recall that we use the 1st and 2nd interior marked points 
to take fiber product with 
$\overline{\text{\bf f}}_{a(1)}$, 
${\overline{\text{\bf f}}}_{a(2)}$.
We consider the $\mathcal M^{\text{main}}_{1;\ell+2}(\beta)$ 
factor $(\Sigma,u)$ of this fiber product.
We say that a disk component thereof is {\it bubble disk component} if 
it does not contain neither the first nor the second interior marked point.
(The `disk component' above includes the tree of sphere 
components rooted on it. See the beginning of Section \ref{sec:equikuracot}.)
\begin{lem}\label{perturbnbdsigmasec9}
There exist a Kuranishi structure and a system of  
continuous families of  multisections $\frak s$,
on a neighborhood of $\frak{forget}^{-1}([\Sigma_0])$, 
with the following properties:
\begin{enumerate}
\item $\frak s$ is transversal to $0$.
\item
$\frak s$  is compatible at the boundaries described in
Lemma \ref{boudaryMippai}.
(Note $\partial \Delta = \emptyset$ so  Lemma \ref{boudaryMippai}.1
does not occur.)
\item
On fiber product factor corresponding to the bubble disk component,
$\frak s$ coincides 
with $\frak q$-multisection.
\item $\frak s$ is compatible with the forgetful map of the 1st, \dots, $k$-th boundary marked points.
\item $\frak s$ is invariant under the arbitrary permutations of the interior marked points.
\item $\frak s$ is transversal to the map $(\ref{componentglue})$ in the sense we describe below:
\item $\text{\rm ev}_0$ is a strata-wise submersion on the intersection of $\frak{forget}^{-1}([\Sigma_0])$ 
and the zero set of $\frak s$.
\end{enumerate}
\end{lem}
Let us state Item 6 precisely.
Let $U'$ be a Kuranishi neighborhood of a point $\text{\bf x}$ in
$\mathcal M_{k+1}(\alpha_i;\beta_i;\overline{\text{\bf f}}_{a(1)},
{\overline{\text{\bf f}}}_{a(2)};\ell_{i,1},\ell_{i,2};\Delta)$.
Let $\text{\rm Glue}
(\text{\bf x}) = \text{\bf y} \in \mathcal M^{\text{main}}_{k+1;\ell_1+\ell_2+2}
(\beta;\overline{\text{\bf f}}_{a(1)} \otimes
{\overline{\text{\bf f}}}_{a(2)} \otimes \frak b_{\rm high}^{\otimes \ell})$. Let $U$ be a Kuranishi neighborhood of $\text{\bf y}$
and $E$  an obstruction bundle.
Let $\frak s_i$ be a branch of $\frak s_w$ that is a member of our family of multisections $\frak s = \{\frak s_w\}$. We require
that $(\frak s_i)^{-1}(0)$ is transversal to
the map $: U' \to U$, which is induced by
$\text{\rm Glue}$.
\begin{rem}
The compatibility with $(\ref{componentglue})$ we claimed in Item 6 above 
implies certain compatibility with the sphere bubble there.
We do not require any other compatibility with the sphere bubble. 
See  Section \ref{sec:equikuracot}.
\end{rem}
\par
The proof of  Lemma \ref{perturbnbdsigmasec9}
is given  in Section \ref{sec:equikuracot}.
\begin{rem}\label{interiorinconsistent}
The family of multisections $\frak s$ on a neighborhood of $\frak{forget}^{-1}([\Sigma_0])$
which we have just described induces a continuous family of multisections on
$\mathcal M_{k+1}(\alpha;\beta;\overline{\text{\bf f}}_{a(1)},
{\overline{\text{\bf f}}}_{a(2)};\ell_{1},\ell_{2};\Delta)$
by the pull-back under the map $\text{\rm Glue}$.
On the other hand, by definition,
$\mathcal M_{k+1}(\alpha;\beta;\overline{\text{\bf f}}_{a(1)},
{\overline{\text{\bf f}}}_{a(2)};\ell_{1},\ell_{2};\Delta)$ is the fiber product of two spaces
\begin{equation}\label{98}
\mathcal M_{\ell_1+3}(\alpha; \overline{\text{\bf f}}_{a(1)}
\otimes  {\overline{\text{\bf f}}}_{a(2)}\otimes \frak b_{\rm high}^{\otimes\ell_1}),
\quad
\mathcal M_{k+1;\ell_2+1}^{\text{\rm main}}(\beta;\frak b_{\rm high}^{\otimes\ell_2})
\end{equation}
over the diagonal. We already have multisections that is the
$\frak q$-multisections
(see Definition \ref{defn66}.) on the second component in
(\ref{98}). However, the multisection induced by $\text{\rm Glue}$
does not necessarily coincide with the fiber product multisection.
In fact, it seems difficult to construct a system of multisections
($\frak q$-multisections) so that the fiber product multisection of (\ref{98})
along the diagonal is transversal. This is related to
 \cite[Remark 11.4]{fooo08}.
\end{rem}
We next construct a multisection of a neighborhood of
$\frak{forget}^{-1}([\Sigma_{12}])$. 

We first define a map:
\begin{equation}\label{gluediscmap}
\aligned
\text{\rm Glue}:
&\left(\mathcal M_{k_1+1;\ell_1+1}(\beta_{(1)};\overline{\text{\bf f}}_{a(1)}
\otimes \frak b_{\rm high}^{\otimes\ell_1}) \times
\mathcal M_{k_2+1;\ell_2+1}(\beta_{(2)};
{\overline{\text{\bf f}}}_{a(2)}
\otimes \frak b_{\rm high}^{\otimes\ell_2})\right)
\\
&\qquad{}_{(\text{\rm ev}_0,\text{\rm ev}_0)}\times_{(\text{\rm ev}_i,\text{\rm ev}_j)}\mathcal M_{k_3+3;\ell_3}(\beta_{(0)};\frak b_{\rm high}^{\otimes\ell_3})
\to \frak{forget}^{-1}([\Sigma_{12}]).
\endaligned
\end{equation}

Let $w_1 : \Sigma_{z_1^{\text
{\rm int}}} \to X$ together with marked points on 
$\Sigma_{z_1^{\text
{\rm int}}}$ define an element of 
$\mathcal M_{k_1+1;\ell_1+1}(\beta_{(1)};\overline{\text{\bf f}}_{a(1)}
\otimes \frak b_{\rm high}^{\otimes\ell_1})$
and let 
$w_2 : \Sigma_{z_2^{\text
{\rm int}}} \to X$ together with marked points on 
$\Sigma_{z_2^{\text
{\rm int}}}$ define an element of 
$\mathcal M_{k_2+1;\ell_2+1}(\beta_{(2)};\overline{\text{\bf f}}_{a(1)}
\otimes \frak b_{\rm high}^{\otimes\ell_2})$.
In particular. we have
\begin{equation}\label{homologyclasseach}
w_{1*}([\Sigma_{z_1^{\text
{\rm int}}}]) = \beta_{(1)},
\quad
w_{2*}([\Sigma_{z_2^{\text
{\rm int}}}]) = \beta_{(2)}.
\end{equation}
Here $\Sigma_{z_1^{\text{\rm int}}}$ and 
$\Sigma_{z_2^{\text
{\rm int}}}$
come with interior marked points $z_1^{\text
{\rm int}}$ and $z_2^{\text
{\rm int}}$ respecively.
They are regarded as the first and the second interior marked points 
and we take a fiber products with 
$\overline{\text{\bf f}}_{a(1)}$, 
${\overline{\text{\bf f}}}_{a(2)}$, at those marked points, respectively.
(We take fiber product with $\frak b_{\rm high}$ at the other 
interior marked points.)

Let $w_0 : \Sigma_{z_0} \to X$ together  
with marked points on 
$\Sigma_{z_0}$ define an element of 
$\mathcal M_{k_3+3;\ell_3}(\beta_{(0)};\frak b_{\rm high}^{\otimes\ell_3})$.
We have
$$
w_{0*}([\Sigma_{z_0}]) = \beta_{(0)}.
$$
Here $\Sigma_{z_0}$ has $\ell_3$ interior marked points at which 
we take fiber product with $\frak b_{\rm high}$.
$z_0$ is the $0$-th boundary marked point of 
$\Sigma_{z_0}$.

We assume that the evaluation map 
at the $i$-th (resp. $j$-th) boundary marked point of $w_0$ 
coincides with the 
evaluation map at the 0-th boundary marked point of $w_1$
(resp. $w_2$).

We put 
$\Sigma = \Sigma_{z_0}\cup \Sigma_{z_1^{\text
{\rm int}}}\cup \Sigma_{z_2^{\text
{\rm int}}}$
where we glue $i$-th (resp. $j$-th) boundary marked point of 
$\Sigma_{z_0}$ with the $0$-th boundary marked point of 
$\Sigma_{z_1^{\text
{\rm int}}}$
(resp. $\Sigma_{z_2^{\text
{\rm int}}}$).
The maps 
$w_1,w_2,w_0$ induce 
$w : \Sigma \to X$. 
We thus obtain
$$
[w;\Sigma] \in \frak{forget}^{-1}([\Sigma_{12}])
 \subset \mathcal M^{\text{main}}_{k+1;\ell+2}(\beta;\overline{\text{\bf f}}_{a(1)} \otimes
{\overline{\text{\bf f}}}_{a(2)} \otimes \frak b_{\rm high}^{\otimes \ell}).
$$
\begin{rem}\label{nonuniquedec}
The decomposition
$\Sigma = \Sigma_{z_0}\cup \Sigma_{z_1^{\text
{\rm int}}}\cup \Sigma_{z_2^{\text
{\rm int}}}$ is {\it not} unique.
\end{rem}

\begin{defn}
We denote by $\frak X(\beta_{(1)},\beta_{(2)},\beta_{(0)})$ the
union of the images of (\ref{gluediscmap}) for various $i,j$.
We thus have a decomposition
\begin{equation}\label{homologyclasseachdecomp}
\frak{forget}^{-1}([\Sigma_{12}])
=
\bigcup_{\beta_{(1)}+\beta_{(2)}+\beta_{(0)}=\beta}\frak X(\beta_{(1)},\beta_{(2)},\beta_{(0)}).
\end{equation}
\end{defn}

Since the decomposition $\Sigma = \Sigma_{z_0}\cup \Sigma_{z_1^{\text
{\rm int}}}\cup \Sigma_{z_2^{\text
{\rm int}}}$  is not unique, the images of the maps $\text{\rm Glue}$
could overlap. (See Figure 2.6.4.)
\par
\hskip0.5cm
\epsfbox{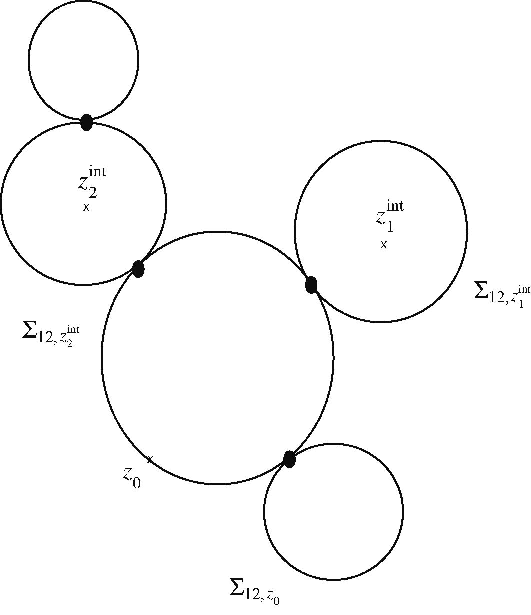}
\par
\centerline{\bf Figure 2.6.4}
\par
\begin{exm}\label{Examle16}
Let $\Sigma$ be the union of $4$ disks $D_0,\dots,D_3$.
We assume $D_1 \cap D_3 = \{p_{13}\}$, $D_0 \cap D_3= \{p_{03}\}$, $D_0 \cap D_2= \{p_{02}\}$
and other intersections are empty.
We take $z_1^{\text{\rm int}} \in D_1$, $z_2^{\text{\rm int}} \in D_2$
and $z_0 \in D_0$.
Let $w: \Sigma \to X$ be a holomorphic map.
Denote $\beta_{(i)} = w_*([D_i])$.
Then
\begin{equation}\label{associative1}
\aligned
&(((w;D_1\#D_3;z_1^{\text{\rm int}};p_{03}),
(w;D_2;z_2^{\text{\rm int}};p_{02})),
(w;D_0;z_0,p_{03},p_{02}))
\\
&\in (\mathcal M_{1;1}(\beta_{(1)}+\beta_{(3)};\overline{\text{\bf f}}_{a(1)}) \times
\mathcal M_{1;1}(\beta_{(2)};{\overline{\text{\bf f}}}_{a(2)})) {}_{(\text{\rm ev}_0,\text{\rm ev}_0)}\times_{(\text{\rm ev}_1,\text{\rm ev}_2)}
\mathcal M_{3;0}(\beta_{(0)}).
\endaligned\end{equation}
On the other hand,
\begin{equation}\label{associative2}\aligned
&(((w;D_1;z_1^{\text{\rm int}};p_{13})),
(w;D_2;z_2^{\text{\rm int}};p_{02})),
(w;D_3\#D_0;z_0,p_{13},p_{02}))
\\
&\in (\mathcal M_{1;1}(\beta_{(1)};\overline{\text{\bf f}}_{a(1)}) \times
\mathcal M_{1;1}(\beta_{(2)};{\overline{\text{\bf f}}}_{a(2)})) {}_{(\text{\rm ev}_0,\text{\rm ev}_0)}
\times_{(\text{\rm ev}_1,\text{\rm ev}_2)}
\mathcal M_{3;0}(\beta_{(0)}+\beta_{(3)}).
\endaligned\end{equation}
They both map to the same element under the map $\text{\rm Glue}$.
(See Figure 2.6.5.)
\par
\hskip2.4cm
\epsfbox{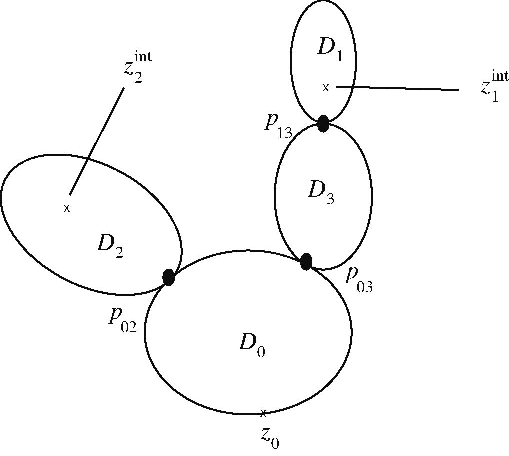}
\par
\centerline{\bf Figure 2.6.5}
\par
Here, however, the situation is different from the one for the case of
$\frak{forget}^{-1}([\Sigma_{0}])$. Namely,
we can choose our multisection on a neighborhood of
$\frak{forget}^{-1}([\Sigma_{12}])$ so that
(\ref{gluediscmap}) preserves the multisection.
(Compare Remark \ref{interiorinconsistent}.)
In fact, the left hand side of (\ref{associative1}) is contained in
\begin{equation}\label{M11ibanme}
\aligned
(\mathcal M_{1;1}(\beta_{(1)};\overline{\text{\bf f}}_{a(1)}) {}_{\text{\rm ev}_0}\times_{\text{\rm ev}_1} \mathcal M_{2;0}(\beta_{(3)})) 
&\times
\mathcal M_{1;1}(\beta_{(2)};{\overline{\text{\bf f}}}_{a(2)}))\\ 
&{}_{(\text{\rm ev}_0,\text{\rm ev}_0)}\times_{\text{(\rm ev}_1,\text{\rm ev}_2)}
\mathcal M_{3;0}(\beta_{(0)}),
\endaligned
\end{equation}
and  the left hand side of (\ref{associative2}) is contained in
\begin{equation}\label{M112banme}
\aligned
&(\mathcal M_{1;1}(\beta_{(1)};\overline{\text{\bf f}}_{a(1)}) ) \times
\mathcal M_{1;1}(\beta_{(2)};{\overline{\text{\bf f}}}_{a(2)})) \\
&{}_{(\text{\rm ev}_0,\text{\rm ev}_0)}\times_{(\text{\rm ev}_1,\text{\rm ev}_2)}
(\mathcal M_{2;0}(\beta_{(3)}){}_{\text{\rm ev}_1}\times_{\text{\rm ev}_1} \mathcal M_{3;0}(\beta_{(0)})).
\endaligned
\end{equation}
These two spaces are isomorphic together with their Kuranishi structures.
(This is because of the associativity of the fiber product.)
Moreover the $\frak q$-multisections on them coincide with each other.
We note that  $\frak q$-multisections are constructed inductively over
the number of {\it disk} components (and of symplectic area).
Such induction was possible because $\text{\rm ev}_0$ (the evaluation map at the
boundary marked point) automatically becomes a submersion by $T^n$ equivariance.
Since the fiber products appearing in both (\ref{M11ibanme}) and
(\ref{M112banme}) correspond to gluing at the boundary singular
point, they are compatible with the $\frak q$-multisections.
\end{exm}
We have the following lemma which corresponds to Lemma
\ref{normalcrossing}.
\begin{lem}\label{cornergluelocal}
Each $\frak X(\beta_{(1)},\beta_{(2)},\beta_{(0)})$
has Kuranishi structure with corners.
Let $\text{\bf x} \in \frak{forget}^{-1}([\Sigma_{12}])$.
Then there exist positive integers $N_1$, $N_2$ such that the following
holds.
\begin{enumerate}
\item
A Kuranishi neighborhood $U$ of $\text{\bf x}$ in
$\mathcal M^{\text{\rm main}}_{k+1;\ell+2}(\beta;\overline{\text{\bf f}}_{a(1)} \otimes  {\overline{\text{\bf f}}}_{a(2)} \otimes \frak b^{\otimes \ell})$
is $\overline U \times [0,1)^{N_1+1} \times [0,1)^{N_2+1}$.
\item
The space $\frak{forget}^{-1}([\Sigma_{12}])$ in this
Kuranishi neighborhood is
$$
\left(\left(\overline U \times \partial([0,1)^{N_1}) \times \partial([0,1)^{N_2})\right)
\cap s^{-1}(0)\right) / \Gamma_{\text{\bf x}}.
$$
\item
$\text{\bf x}$ is contained in the union of $(N_1+1) \times (N_2+1)$ different
strata, that is, in $\frak X(\beta_{(m_1,m_2,1)},\beta_{(m_1,m_2,2)},\beta_{(m_1,m_2,0)})$.
($m_1 = 1,\dots,N_1$, $m_2 = 1,\dots,N_2$.)
\item
$\text{\bf x}$ lies on the codimension $N_1+N_2$ corner on each
stratum in Item 3, that is  $\frak X(\beta_{(m_1,m_2,1)},\beta_{(m_1,m_2,2)},\beta_{(m_1,m_2,0)})$.
\item
The Kuranishi neighborhood of $\text{\bf x}$
in $\frak X(\beta_{(m_1,m_2,1)},\beta_{(m_1,m_2,2)},\beta_{(m_1,m_2,0)})$ can be identified
with
$$
\overline U \times ([0,1)^{m_1} \times \{0\} \times
[0,1)^{N_1-m_1}) \times ([0,1)^{m_2} \times \{0\} \times
[0,1)^{N_2-m_2}).
$$
\item The obstruction bundle of $\frak X(\beta_{(m_1,m_2,1)},\beta_{(m_1,m_2,2)},\beta_{(m_1,m_2,0)})$
is obtained by restricting the obstruction bundle of
$\mathcal M^{\text{\rm main}}_{k+1;\ell+2}(\beta;\overline{\text{\bf f}}_{a(1)} \otimes  {\overline{\text{\bf f}}}_{a(2)} \otimes \frak b_{\rm high}^{\otimes \ell})$.
\end{enumerate}
\end{lem}
The proof of this lemma is similar to that of Lemma
\ref{normalcrossing} and omitted.
\begin{lem}\label{perturbconsistence}
The multisections on $\frak X(\beta_{(m_1,m_2,1)},\beta_{(m_1,m_2,2)},\beta_{(m_1,m_2,0)})$
induced by $(\ref{gluediscmap})$ from the $\frak q$-multisections 
(for various $m_1, m_2$) 
coincide each other on the overlapped part.
They also coincide with the $\frak q$-multisection of
$\mathcal M^{\text{\rm main}}_{k+1;\ell+2}(\beta;\overline{\text{\bf f}}_{a(1)}
\otimes  {\overline{\text{\bf f}}}_{a(2)} \otimes \frak b_{\rm high}^{\otimes \ell})$.
Hence it was already extended to the neighborhood of
$\frak{forget}^{-1}([\Sigma_{12}])$.
\end{lem}
\begin{proof}
The  $\frak q$-multisection and the Kuranishi structure (for which $\frak q$-multisection is defined) 
is obtained as the fiber product of those given on each of the disk component.
(See the paragraph right after Example \ref{Examle16}.)
This is a consequence of our construction of $\frak q$-multisections.
\end{proof}
\begin{rem}
We note that $\frak q$-multisection is not necessarily a fiber product 
of those of $S^2$ and $D^2$ components. 
\end{rem}
Thus we have defined (family of) multisections on 
a neighborhood of $\frak{forget}^{-1}([\Sigma_{12}])$ and of 
$\frak{forget}^{-1}([\Sigma_{0}])$.
We can
then extend them to the whole 
$\mathcal M^{\text{main}}_{k+1;\ell_1+\ell_2+2}
(\beta;\overline{\text{\bf f}}_{a(1)} \otimes
{\overline{\text{\bf f}}}_{a(2)} \otimes \frak b_{\rm high}^{\otimes \ell})$ 
so that Items 1,3,4,5,7 of Lemma \ref{perturbnbdsigmasec9} are 
satisfied.
We also denote this system of multisections by $\frak s$.
\par
From now on we consider the moduli space with $k=0$, $k_1,k_2, k_3=0$.
Namely we consider the case when there is only one boundary marked point
except the situation of (\ref{gluediscmap}).
In the situation of (\ref{gluediscmap}) we consider the case $k_3=0$ in the factor 
$\mathcal M_{k_3+3;\ell_3}(\beta_{(0)};\frak b_{\rm high}^{\otimes\ell_3})$. 
So we consider the case when we have three boundary marked points.
Studying those cases is enough to prove Theorem \ref{multiplicative}.
\par
We use a family of multisections  $\frak s$ to
define integration along the fiber of the map
$$
\text{\rm ev}_0: \mathcal M^{\text{main}}_{1;\ell+2}
(\beta;\overline{\text{\bf f}}_{a(1)} \otimes
{\overline{\text{\bf f}}}_{a(2)} \otimes \frak b_{\rm high}^{\otimes \ell})
\to L(u)
$$
restricted to
$$
\frak s^{-1}(0) \cap \frak{forget}^{-1}([\Sigma_0]),
\quad
\frak s^{-1}(0) \cap \frak{forget}^{-1}([\Sigma_{12}]).
$$
Hereafter we omit $\frak s$ from the notation, in case no confusion can occur. Namely we will just write as
$\frak{forget}^{-1}([\Sigma_0])$ etc.
in place of $\frak s^{-1}(0) \cap \frak{forget}^{-1}([\Sigma_0])$
etc.
\begin{prop}\label{0is12} We have
$$
(\text{\rm ev}_0)_* \left(\frak{forget}^{-1}([\Sigma_0])\right)
= (\text{\rm ev}_0)_* \left(\frak{forget}^{-1}([\Sigma_{12}])\right).
$$
\end{prop}
Here we recall that in this section we consider the case where 
the virtual dimension of $\frak{forget}^{-1}([\Sigma_0])$
is $n$ and so
$$
(\text{\rm ev}_0)_* \left(\frak{forget}^{-1}([\Sigma_0])\right)
\in H_n(L(u);\Lambda_0) \cong H^0(L(u);\Lambda_0) = \Lambda_0.
$$
Similarly $(\text{\rm ev}_0)_* \left(\frak{forget}^{-1}([\Sigma_{12}])\right)
\in \Lambda_0$.
\par
Proposition \ref{0is12} is Step II we mentioned at the beginning of the 
proof of Theorem \ref{multiplicative}.

\begin{proof}
We first take an element $[\Sigma_{2}] \in \mathcal M_{1;2}^{\text{\rm main}}$
which is in the stratum 1 in Lemma \ref{stratifyD12} and is sufficiently close to $[\Sigma_{12}]$.
\begin{lem}\label{genericpointinverse}
$$
(\text{\rm ev}_0)_* \left(\frak{forget}^{-1}([\Sigma_2])\right)
= (\text{\rm ev}_0)_* \left(\frak{forget}^{-1}([\Sigma_{12}])\right).
$$
\end{lem}
\begin{proof}
Let $\frak z \in L(u)$ be a generic point.
We note that
$(\text{\rm ev}_0)_* \left(\frak{forget}^{-1}([\Sigma_{12}])\right)$ may
be regarded as a $T^n$-invariant  differential $0$ form, that is
nothing but a constant function. This implies that
the order counted with sign of the moduli space
$(\text{\rm ev}_0)^{-1}(\frak z) \cap \frak{forget}^{-1}([\Sigma_{12}])$
is the right hand side of Lemma \ref{0is12}. Namely:
$$
\# \left( (\text{\rm ev}_0)^{-1}(\frak z) \cap \frak{forget}^{-1}([\Sigma_{12}])  \right)
= (\text{\rm ev}_0)_* \left(\frak{forget}^{-1}([\Sigma_{12}])\right).
$$
In the same way we have
$$
\# \left( (\text{\rm ev}_0)^{-1}(\frak z) \cap \frak{forget}^{-1}([\Sigma_{2}])  \right)
= (\text{\rm ev}_0)_* \left(\frak{forget}^{-1}([\Sigma_{2}])\right).
$$
By choosing $\frak z$ generic, we may assume that
any element of $ (\text{\rm ev}_0)^{-1}(\frak z) \cap \frak{forget}^{-1}([\Sigma_{12}]  )$ satisfies
analogous condition as Condition \ref{1tai1cond}.
We may also assume that any element of $
 (\text{\rm ev}_0)^{-1}(\frak z) \cap\frak{forget}^{-1}([\Sigma_{12}])$ is a map from a union of three disks. Then,
by a standard gluing and cobordism argument, we have
$$
\# \left( (\text{\rm ev}_0)^{-1}(\frak z) \cap \frak{forget}^{-1}([\Sigma_{12}])  \right)
=
\# \left( (\text{\rm ev}_0)^{-1}(\frak z) \cap \frak{forget}^{-1}([\Sigma_{2}] ) \right).
$$
Lemma \ref{genericpointinverse} follows.
\end{proof}
We next take a sequence of elements  $[\Sigma_{3,i}] \in \mathcal M_{1;2}^{\text{\rm main}}$
which converges to $[\Sigma_{0}]$ and
is in the stratum (1) in Lemma \ref{stratifyD12}.
\begin{lem}\label{sublemma:Sigma3limit}
\begin{equation}\label{3moonaji}
\lim_{i\to \infty}(\text{\rm ev}_0)_* \left(\frak{forget}^{-1}([\Sigma_{3,i}])\right)
= (\text{\rm ev}_0)_* \left(\frak{forget}^{-1}([\Sigma_{0}])\right).
\end{equation}
\end{lem}
The proof is more involved than the proof of Lemma \ref{genericpointinverse}
since our Kuranishi structure and family of multisections are not $T^n$-equivariant
in a neighborhood of $\frak{forget}^{-1}([\Sigma_{0}])$.
We also note that the forgetful map is neither smooth nor a submersion at the fiber of 
$[\Sigma_0]$.
The proof will be given in Section \ref{sec:cornerestimate}.
\begin{lem}
\begin{equation}\label{4moonaji}
(\text{\rm ev}_0)_* \left(\frak{forget}^{-1}([\Sigma_{3,i}])\right)
= (\text{\rm ev}_0)_* \left(\frak{forget}^{-1}([\Sigma_{2}])\right).
\end{equation}
\end{lem}
\begin{proof}
We take a generic smooth path  $\gamma: [0,1] \to \mathcal M^{\text{\rm main}}_{1;2}$ joining $[\Sigma_{2}]$ and
$[\Sigma_{3,i}]$ and contained in the stratum 1 in Lemma \ref{stratifyD12}. By Lemma \ref{smoothmaptoDM}
$$
[0,1] {}_{\gamma}\times_{\frak{forget}}
\mathcal M^{\text{main}}_{1;\ell+2}
(\beta;\overline{\text{\bf f}}_{a(1)}\otimes
{\overline{\text{\bf f}}}_{a(2)}\otimes \frak b_{\rm high}^{\otimes \ell})
$$
has a Kuranishi structure and we may choose $\frak s$
and $\gamma$ so that we can use Stokes' theorem
(\cite[Lemma 12.18]{fooo09}) on
$
[0,1] {}_{\gamma}\times_{\frak{forget}} \frak s^{-1}(0)
$.
It implies (\ref{4moonaji}).
\end{proof}
Lemmata \ref{genericpointinverse}, \ref{sublemma:Sigma3limit}, \ref{4moonaji}  immediately imply Proposition \ref{0is12}.
\end{proof}
\begin{lem}\label{atsigma12}
If $a(1),a(2)\ne 0$, we have:
$$\aligned
&(\text{\rm ev}_0)_* \left(\frak{forget}^{-1}([\Sigma_{12}])\right) \\
&= \sum_{\beta_{(1)}+\beta_{(2)}=\beta}
\sum_{\ell_1+\ell_2=\ell}\frac{\ell!}{\ell_1!\ell_2!}
(\text{\rm ev}_0)_* \left(\mathcal M_{1,\ell_1+1}^{\text{\rm main}}
(\beta_{(1)};\overline{\text{\bf f}}_{a(1)}\otimes \frak b_{\rm high}^{\otimes \ell_1})
\right)\\
&\qquad\qquad\qquad\qquad
\times (\text{\rm ev}_0)_* \left(\mathcal M_{1,\ell_2+1}^{\text{\rm main}}
(\beta_{(2)};{\overline{\text{\bf f}}}_{a(2)}\otimes \frak b_{\rm high}^{\otimes \ell_2})\right).
\endaligned$$
Here summation in the right hand side is taken over those
$\mathcal M_{1,\ell_1+1}^{\text{\rm main}}
(\beta_{(1)};\overline{\text{\bf f}}_1\otimes \frak b_{\rm high}^{\otimes \ell_1})$
which are $n$-dimensional.
(We recall $\overline{\text{\bf f}}_0 = \text{\rm PD}[X]$.)
\par
If $a(1)=0$ and $a(2)\ne 0$, then
$$
(\text{\rm ev}_0)_* \left(\frak{forget}^{-1}([\Sigma_{12}])\right)
= (\text{\rm ev}_0)_* \left(\mathcal M_{1,\ell+1}^{\text{\rm main}}
(\beta);{\overline{\text{\bf f}}}_{a(2)}\otimes \frak b_{\rm high}^{\otimes \ell})\right).
$$
The case when $a(1)\ne 0$ and $a(2)= 0$ is similar.
\par
If $a(2)=a(1)=0$, then
$$
(\text{\rm ev}_0)_* \left(\frak{forget}^{-1}([\Sigma_{12}])\right) = 1.
$$
\end{lem}
\begin{rem}
We note that we do not break $T^n$ symmetry in a neighborhood of 
$\frak{forget}^{-1}([\Sigma_{12}])$.
\end{rem}
\begin{proof}
Using the decomposition (\ref{homologyclasseachdecomp}), we study
the contributions of various components of
$\frak X(\beta_{(1)},\beta_{(2)},\beta_{(0)})$ to the left hand side
of Lemma \ref{atsigma12}.
We first consider the case $\beta_{(0)}=0$.
Using Lemma \ref{perturbconsistence} it is easy to see that
$$\aligned
& (\text{\rm ev}_0)_*(\frak X(\beta_{(1)},\beta_{(2)},0))\\
&=
\sum_{\ell_1+\ell_2=\ell}\frac{\ell!}{\ell_1!\ell_2!}
 (\text{\rm ev}_0)_* \left(\mathcal M_{1,\ell_1+1}^{\text{\rm main}}
(\beta_{(1)};\overline{\text{\bf f}}_{a(1)}\otimes \frak b_{\rm high}^{\otimes \ell_1})
\right)
\times (\text{\rm ev}_0)_* \left(\mathcal M_{1,\ell_2+1}^{\text{\rm main}}
(\beta_{(2)};{\overline{\text{\bf f}}}_{a(2)}\otimes \frak b_{\rm high}^{\otimes \ell_2})\right),
\endaligned$$
if $a(1),a(2)\ne 0$.
(Note $\frac{\ell!}{\ell_1!\ell_2!}$ counts the way how to divide $\ell$ interior 
marked points to the two factors.)
\par
If $a(1)=0$, then only the case $\beta_{(1)} = 0 = \ell_1$ can contribute and
$$
(\text{\rm ev}_0)_* (\mathcal M_{1,1}^{\text{\rm main}}
(\beta_{(1)};\overline{\text{\bf f}}_0\otimes \frak b_{\rm high}^{\otimes 0})) = 1.
$$
In case $a(1)=a(2)=0$, only the case $\beta_{(1)} = \beta_{(2)} =0 = \ell_1= \ell_2$
contributes and it gives $1$.
\par
We next study the case $\beta_{(0)} \ne 0$.
We observe that the differential form
\begin{equation}\label{pushfirstcomp}
(\text{\rm ev}_0)_* \left(\mathcal M_{1,\ell_1+1}^{\text{\rm main}}
(\beta_{(1)};\overline{\text{\bf f}}_1\otimes \frak b_{\rm high}^{\otimes \ell_1})
\right)
\end{equation}
is degree $0$ and constant, (unless it is zero).
This is a consequence of the fact that
$$
\text{\rm ev}_0: \mathcal M_{1,\ell_1+1}^{\text{\rm main}}
(\beta_{(1)};\overline{\text{\bf f}}_{a(1)}\otimes \frak b_{\rm high}^{\otimes \ell_1})
\to L(u)
$$
is a submersion. (We again remark that we did not break the
$T^n$ symmetry in a neighborhood of $[\Sigma_{12}]$.)
We write (\ref{pushfirstcomp}) as $c_{1,\ell_1,\beta_{(1)}} 1_{L(u)}$. We define
$c_{2,\ell_2,\beta_{(2)}}$ in a similar way.
Now by  \cite[Lemma 12.20]{fooo09} we have
$$\aligned
& (\text{\rm ev}_0)_*(\frak X(\beta_{(1)},\beta_{(2)},\beta_{(0)}))\\
&=
\sum_{\ell_0+\ell_1+\ell_2=\ell}
c_{1,\ell_1,\beta_{(1)}}c_{2,\ell_2,\beta_{(2)}}
(\text{\rm ev}_0)_* (\text{\rm ev}_1,\text{\rm ev}_2)^*(1_{L(u)} \times  1_{L(u)}).
\endaligned
$$
Here
\begin{equation}\label{ev1ev2}
(\text{\rm ev}_1,\text{\rm ev}_2):
\mathcal M_{3,\ell_0}^{\text{\rm main}}
(\beta_{(0)}; \frak b_{\rm high}^{\otimes \ell_0})
\to L(u)^2.
\end{equation}
Since $1_{L(u)}$ is the strict unit, this is zero
unless $\beta_{(0)} = 0$.
(In other words, (\ref{ev1ev2}) factors through
$
\mathcal M_{3,\ell_0}^{\text{\rm main}}
(\beta_{(0)}; \frak b_{\rm high}^{\otimes \ell_0})
\to
\mathcal M_{1,\ell_0}^{\text{\rm main}}
(\beta_{(0)}; \frak b_{\rm high}^{\otimes \ell_0})
$. So this is zero unless $\beta_{(0)} = 0$.)
The proof of Lemma \ref{atsigma12}
is complete.
\end{proof}
It is easy to see that an appropriate weighted sum of the right hand side of Lemma \ref{atsigma12}
becomes the first product term in the left hand side of
equality  (\ref{eq:ks-ring}).
(See (\ref{PPformula}).)
Namely we have
\begin{equation}\label{for922}
\aligned
&\widetilde{\frak{ks}}_{\frak b}(\overline{\text{\bf f}}_{a(1)})\widetilde{\frak{ks}}_{\frak b}({\overline{\text{\bf f}}}_{a(2)})
\\=
&\sum_{\ell_1,\ell_2}\sum_{\beta_{(1)},\beta_{(2)}}\\
& \rho_b(\partial\beta_{(1)})
\frac{\exp(\frak b_2\cap \beta_{(1)})T^{\omega\cap \beta_{(1)}/2\pi}}{\ell_1!}
(\text{\rm ev}_0)_* \left(\mathcal M_{1,\ell_1+1}^{\text{\rm main}}
(\beta_{(1)};\overline{\text{\bf f}}_{a(1)}\otimes \frak b_{\rm high}^{\otimes \ell_1})
\right)\\
& \rho_b(\partial\beta_{(2)})
\frac{\exp(\frak b_2\cap \beta_{(2)})T^{\omega\cap \beta_{(2)}/2\pi}}{\ell_2!}
(\text{\rm ev}_0)_* \left(\mathcal M_{1,\ell_2+1}^{\text{\rm main}}
(\beta_{(2)};{\overline{\text{\bf f}}}_{a(2)}\otimes \frak b_{\rm high}^{\otimes \ell_2})\right).
\endaligned
\end{equation}
Here we put $b= \sum_{i=1}^n x_i\text{\bf e}_i$ and $y_i = T^{u_i}e^{x_i}$.
Then $\rho_b(\partial\beta_{(1)}), \rho_b(\partial\beta_{(2)})$ are monomials 
of $y_i$'s and the right hand side of (\ref{for922}) 
is a formal power series of $y_i$'s.
\par\smallskip
In the rest of this section we will check that the weighted sum of 
the left hand side of Proposition \ref{0is12} becomes
the second term of (\ref{eq:ks-ring})
modulo an element of the defining ideal of the Jacobian ring. 
The proof is a combination of the proofs of
associativity of the quantum cup product and of commutativity of
the diagram in Theorem \ref{indepenceKS}.

We already defined the
moduli space
$\mathcal M_{1}(\alpha;\beta;\overline{\text{\bf f}}_{a(1)},
{\overline{\text{\bf f}}}_{a(2)};\ell_1,\ell_2;K)$ for a smooth singular chain
$K \subset X^2$.
So far we have considered the case $K = \Delta$, which is $T^n$-invariant.
In the next step, we need to use $K$ which is not $T^n$-invariant.
It prevents us from taking a $T^n$-equivariant multisection on
$\mathcal M_{1}(\alpha;\beta;\overline{\text{\bf f}}_{a(1)},
{\overline{\text{\bf f}}}_{a(2)};\ell_1,\ell_2;K)$ any more.
This is another reason why we need to break the $T^n$ symmetry and  use the technique of continuous family of
multisections. 
\begin{conds}\label{compatiblefamilies}
We consider systems of continuous families of multisections
$\frak s$ on various
$\mathcal M_{1}(\alpha;\beta;\overline{\text{\bf f}}_{a(1)},{\overline{\text{\bf f}}}_{a(2)};\ell_1,\ell_2;K)$
and
$\mathcal M_{1}(\alpha;\beta;\overline{\text{\bf f}}_{a(1)},{\overline{\text{\bf f}}}_{a(2)};\ell_1,\ell_2;\partial K)$
with the following properties.
\begin{enumerate}
\item
The evaluation maps 
$$\text{\rm ev}_0: \mathcal M_{1}(\alpha;\beta;\overline{\text{\bf f}}_{a(1)},{\overline{\text{\bf f}}}_{a(2)};\ell_1,\ell_2;K)^{\frak s}
\to L(u)
$$
and
$$\text{\rm ev}_0: \mathcal M_{1}(\alpha;\beta;\overline{\text{\bf f}}_{a(1)},{\overline{\text{\bf f}}}_{a(2)};\ell_1,\ell_2;\partial K)^{\frak s}
\to L(u)$$ are submersions.
Here $\mathcal M_{1}(\alpha;\beta;\overline{\text{\bf f}}_{a(1)},{\overline{\text{\bf f}}}_{a(2)};\ell_1,\ell_2; K)^{\frak s}$
denotes the zero set of the multisection $\frak s$ of the
Kuranishi structure of the moduli space
$$\mathcal M_{1}(\alpha;\beta;\overline{\text{\bf f}}_{a(1)},{\overline{\text{\bf f}}}_{a(2)};\ell_1,\ell_2;K).$$
\item
They are compatible at the boundaries described in Lemma
\ref{boudaryMippai}.
\item
They are compatible with the forgetful maps of the
boundary marked points.
\item 
They coincide with the $\frak q$-multisection on the bubble disk components.
\end{enumerate}
\end{conds}
Since $\partial\partial K = 0$, we do not
have the stratum 1 of Lemma \ref{boudaryMippai} for the case of $\partial K$.
\begin{lem}\label{prodcontfamiext}
If we have a system of continuous families of multisections satisfying
Condition \ref{compatiblefamilies} for $(\partial K,0)$,
then we can find a system of continuous families of multisections
for $(K,\partial K)$ satisfying Condition \ref{compatiblefamilies} which extends
the given one on $\partial K$.
\end{lem}
\begin{proof}
By an induction over $\beta\cap \omega$,  $\ell_1$, $\ell_2$,
and $k$, the lemma follows easily from  \cite[Lemma 12.14]{fooo09}.
\end{proof}
Now we go back to the proof of Theorem \ref{multiplicative}.
Let
$\Delta = \{(x,x) \in X^2 \mid x \in X\}$ be the diagonal.
Since $\overline{\text{\bf f}}_i$ ($i=0,\dots,B'$)
form a basis of $H(X;\C)$, there exists a smooth
singular chain $R$ on $X^2$ such that
\begin{equation}\label{diagonaldecomp}
\partial R = \Delta - \sum_{i,j=0}^{B'}
(-1)^{\text{\rm deg} \overline{\text{\bf f}}_i\text{\rm deg} \overline{\text{\bf f}}_j}g^{ij}
(\overline{\text{\bf f}}_i \times \overline{\text{\bf f}}_j),
\end{equation}
where
$
g_{ij} = \langle \overline{\text{\bf f}}_i, \overline{\text{\bf f}}_j\rangle_{\text{\rm PD}_X}$
with $0 \leq i,  j \leq B'$,
and $(g^{ij})$ is the inverse of the matrix $(g_{ij})$.
(See Lemma \ref{LemmaE} for the sign in (\ref{diagonaldecomp}).)
\begin{rem}
We note that $R$ can not be taken to be $T^n$-invariant.
\end{rem}
\par
We consider
$\mathcal M_{1}(\alpha;\beta;\overline{\text{\bf f}}_{a(1)},
{\overline{\text{\bf f}}}_{a(2)};\ell_1,\ell_2;\partial R)$.
We use a family of multisections on
$\mathcal M_{1}(\alpha;\beta;\overline{\text{\bf f}}_{a(1)},
{\overline{\text{\bf f}}}_{a(2)};\ell_1,\ell_2;\Delta)$
such that $\text{\rm ev}_0$ is a submersion.
This family of multisections is induced from  $\frak s$.
\par
By Definition \ref{defn96} we have
\begin{equation}\label{abkakeruc}
\aligned
&\mathcal M_{1}(\alpha;\beta;\overline{\text{\bf f}}_{a(1)},
{\overline{\text{\bf f}}}_{a(2)};\ell_1,\ell_2;\overline{\text{\bf f}}_i \times \overline{\text{\bf f}}_j) \\
&= \mathcal M_{\ell_1+3}(\alpha;\overline{\text{\bf f}}_{a(1)}
\otimes {\overline{\text{\bf f}}}_{a(2)} \otimes \overline{\text{\bf f}}_i
\otimes \frak b_{\rm high}^{\otimes\ell_1})
\times \mathcal M_{1;\ell_2+1}(\beta;\overline{\text{\bf f}}_j \otimes \frak b_{\rm high}^{\otimes\ell_2}).
\endaligned
\end{equation}
Here the product is the {\it direct} product.
On the second factor of (\ref{abkakeruc}), the
(single) multisection
(that is the $\frak q$-multisection)
was already chosen in  \cite[Lemma 6.5]{fooo09} so that
$\text{\rm ev}_0$ is a submersion. We take a (single) multisection on
the first factor, such that it becomes transversal to $0$.
\par
Therefore we apply Lemma \ref{prodcontfamiext} to obtain a system of
continuous families of multisections $\frak s$ on
$\mathcal M_{1}(\alpha;\beta;\overline{\text{\bf f}}_{a(1)},
{\overline{\text{\bf f}}}_{a(2)};\ell_1,\ell_2;R)$ compatible
with other choices we made and such that $\text{\rm ev}_0$ is
a submersion on its zero set.
\par
Therefore we can take an integration along the fiber for
$\mathcal M_{1}(\alpha;\beta;\overline{\text{\bf f}}_{a(1)},
{\overline{\text{\bf f}}}_{a(2)};\ell_1,\ell_2;R)$.
(We omit $\frak s$ from notation here.) (See  \cite[Definition 12.11]{fooo09}.)
\begin{defn}
We put
\begin{equation}\label{fraXproddef}
\aligned
\widetilde{\frak X}(\overline{\text{\bf f}}_{a(1)},{\overline{\text{\bf f}}}_{a(2)};R;\frak b)(b)
= \sum_{\ell_1,\ell_2}\sum_{\alpha,\beta}&
\frac{T^{(\alpha\#\beta)\cap \omega/2\pi}}{(\ell_1+\ell_2)!}
\rho^b(\partial \beta)\exp(\frak b_2\cap (\alpha\#\beta)) \\
&(\text{\rm ev}_{0})_*(
\mathcal M_{1}(\alpha;\beta;\overline{\text{\bf f}}_{a(1)},
{\overline{\text{\bf f}}}_{a(2)};\ell_1,\ell_2;R)),
\endaligned
\end{equation}
and
\begin{equation}\label{tildeXprod}
{\frak X}(\overline{\text{\bf f}}_{a(1)},{\overline{\text{\bf f}}}_{a(2)};R;\frak b)(b)
= \Pi( \widetilde{\frak X}
(\overline{\text{\bf f}}_{a(1)},{\overline{\text{\bf f}}}_{a(2)};R;\frak b)).
\end{equation}
Here $\Pi$ is the harmonic projection in (\ref{harmonicproj}).
\end{defn}
We put $b= \sum_{i=1}^n x_i\text{\bf e}_i$ and $y_i = T^{u_i}e^{x_i}$. 
We then regard $\widetilde{\frak X}(\overline{\text{\bf f}}_{a(1)},{\overline{\text{\bf f}}}_{a(2)};R;\frak b)$,
${\frak X}(\overline{\text{\bf f}}_{a(1)},{\overline{\text{\bf f}}}_{a(2)};R;\frak b)$
as differential form valued formal functions
of $y_i$. More precisely, we have:
\begin{lem}\label{convergenceforR}
The formal power serieses
$\widetilde{\frak X}(\overline{\text{\bf f}}_{a(1)},{\overline{\text{\bf f}}}_{a(2)};R;\frak b)$
and ${\frak X}(\overline{\text{\bf f}}_{a(1)},{\overline{\text{\bf f}}}_{a(2)};R;\frak b)$ are contained in
$\Lambda\langle\!\langle y,y^{-1}\rangle\!\rangle_{0}^{\overset{\circ}P} \widehat{\otimes}\Omega(L(u))$, $H(L(u);\Lambda\langle\!\langle y,y^{-1}\rangle\!\rangle_{0}^{\overset{\circ}P} )$, respectively.
They are contained in $\Lambda\langle\!\langle y,y^{-1}\rangle\!\rangle_{0}^{P}
\widehat{\otimes}\Omega(L(u))$,
$H(L(u);\Lambda\langle\!\langle y,y^{-1}\rangle\!\rangle_{0}^{P} )$ if $\frak b \in \mathcal A(\Lambda_+)$.
\end{lem}
The proof is the same as the proof of Proposition \ref{fraXconv} and is omitted.
We put
$$
\frak{forget}^{-1}([\Sigma_0])_{\beta,\ell}
=
 \mathcal M^{\text{main}}_{1;\ell+2}
(\beta;\overline{\text{\bf f}}_{a(1)} \otimes
{\overline{\text{\bf f}}}_{a(2)} \otimes \frak b_{\rm high}^{\otimes \ell})
\cap
\frak s^{-1}(0) \cap \frak{forget}^{-1}([\Sigma_0]).
$$
\begin{lem}\label{prodmainlemma}
\begin{equation}\label{multiplicativeatosukosi}
\aligned
&\sum_{\beta,\ell}
\frac{T^{\beta\cap \omega/2\pi}\exp(\frak b_2\cap \alpha)}{\ell!}\rho_b(\partial\beta)
\Pi (\text{\rm ev}_0)_* \left(\frak{forget}^{-1}([\Sigma_0])_{\beta,\ell}\right) \\
&\qquad-  \sum_c c^c_{a(1) a(2)}\widetilde{\frak{ks}}_{\frak b}(\overline{\text{\bf f}}_c)
\\
&= - \delta
_{\text{\rm can}}
^{\frak b,b}({\frak X}(\overline{\text{\bf f}}_{a(1)},{\overline{\text{\bf f}}}_{a(2)};R;\frak b)).
\endaligned
\end{equation}
\end{lem}
\begin{proof} Using Lemma \ref{boudaryMippai},
we study the boundary of $\mathcal M_{1}(\alpha;\beta;\overline{\text{\bf f}}_{a(1)},
{\overline{\text{\bf f}}}_{a(2)};\ell_1,\ell_2;R)$,
when its virtual dimension is $n$.
We then use Stokes' theorem (\cite[Lemma 12.13]{fooo09}) to show
that the sum vanishes. Then we put a weighted sum as in \eqref{fraXproddef}, which will
give rise to the equality (\ref{multiplicativeatosukosi}).
The details are in order.
\par
There are three types of components, 1, 2 and 3 in Lemma \ref{boudaryMippai}.
The case 1 that is the space
$\mathcal M_{1}(\alpha;\beta;\overline{\text{\bf f}}_{a(1)},
{\overline{\text{\bf f}}}_{a(2)};\ell_1,\ell_2;\partial R)$  becomes
the left hand side of Lemma \ref{prodmainlemma} after
taking integration along the fiber and the weighted sum as in
(\ref{fraXproddef}). This follows from (\ref{abkakeruc}),
Lemmata \ref{deltaidentify}, \ref{0is12} and the formula
$$
\overline{\text{\bf f}}_{a(1)} \cup^{\frak b}{\overline{\text{\bf f}}}_{a(2)}
=
\sum_{\alpha,\ell,i,j}g^{ij}\frac{T^{\alpha\cap \omega/2\pi}\exp(\frak b_2\cap \alpha)}{\ell!}
\#\left(\mathcal M_{\ell+3}(\alpha;\overline{\text{\bf f}}_{a(1)}
\otimes {\overline{\text{\bf f}}}_{a(2)} \otimes \overline{\text{\bf f}}_i
\otimes \frak b_{\rm high}^{\otimes\ell})\right)  \overline{\text{\bf f}}_j.
$$
\par
We next study the other cases, 2,3 of Lemma \ref{boudaryMippai}.
\par
For the case $K = R$, the boundary of type 2 becomes
\begin{equation}
\aligned
&\frak Y(\alpha,\beta',\beta'',\ell_1,\ell'_2,\ell''_2)
\\
&:=
\mathcal M_{1}(\alpha;\beta';\overline{\text{\bf f}}_{a(1)},
{\overline{\text{\bf f}}}_{a(2)};\ell_1,\ell'_2;R)
{}_{\text{\rm ev}_0}\times_{\text{\rm ev}_1}
\mathcal M_{2}(\beta'';\frak b_{\rm high}^{\otimes \ell''_2}).
\endaligned
\end{equation}
We put a (single) multisection on
$\mathcal M_{2}(\beta'';\frak b_{\rm high}^{\otimes \ell''_2})$ which is
$T^n$-equivariant.
(Namely the $\frak q$-multisection.)
By definition, Theorem \ref{squizthm} and
Formula (\ref{grouplikeelement}), we have
$$\aligned
&\sum_{\alpha,\beta',\beta''}\sum_{\ell_1,\ell'_2,\ell''_2}
\frac{T^{(\alpha\#\beta'\#\beta'')\cap \omega/2\pi}}{\ell_1!\ell'_2!\ell''_2!}
\rho^b(\partial (\beta'\#\beta''))\exp(\frak b_2\cap (\alpha\#\beta'\#\beta''))
\\
&\qquad\qquad\qquad\qquad\qquad\times
(\text{\rm ev}_0)_*(\frak Y(\alpha,\beta',\beta'',\ell_1,\ell'_2,\ell''_2))
\\
&= \delta^{\frak b,b}(\widetilde{\frak X}(\overline{\text{\bf f}}_{a(1)},
{\overline{\text{\bf f}}}_{a(2)};R;\frak b))
\endaligned
$$
\par
We finally study the boundary of type 3. In the present case this
becomes
$$
\mathcal M_{1;\ell'_{2}+1}^{\text{\rm main}}
(\beta';\frak b_{\rm high}^{\otimes\ell'_2})
_{\text{\rm ev}_0}\times_{\text{\rm ev}_i}
\mathcal M_{2}(\alpha;\beta'';\overline{\text{\bf f}}_{a(1)},{\overline{\text{\bf f}}}_{a(2)};\ell_1,\ell''_2;R).
$$
We remark that $\mathcal M_{1;\ell'_{2}+1}^{\text{\rm main}}
(\beta';\frak b_{\rm high}^{\otimes\ell'_2})$ is $T^n$-invariant and we have
chosen a $T^n$-equivariant multisection
($\frak q$-multisection) on it.
Therefore
$$
(\text{\rm ev}_0)_*
(\mathcal M_{1;\ell'_{2}+1}^{\text{\rm main}}
(\beta';\frak b_{\rm high}^{\otimes\ell'_2}))
$$
is a $T^n$-invariant differential $0$-form on $L(u)$,
which is nothing but a constant function.
Using Condition \ref{compatiblefamilies}.3, this
implies that the contribution of the boundary of type 3
vanishes.
\par
Combining the above argument, we obtain
\begin{equation} \label{multiplicativeatosukosi20}
\aligned
&\sum_{\beta,\ell}
\frac{T^{\beta\cap \omega/2\pi}\exp(\frak b_2\cap \alpha)}{\ell!}
\rho_b(\partial \beta)
\\
&
\qquad (\text{\rm ev}_0)_* \left(\frak{forget}^{-1}([\Sigma_0])_{\beta,\ell}\right)
-  \sum_c c^c_{a(1) a(2)}\widetilde{\frak{ks}}_{\frak b}(\overline{\text{\bf f}}_c)
\\
&=
- \delta^{\frak b,b}(\widetilde{\frak X}(\overline{\text{\bf f}}_{a(1)},
{\overline{\text{\bf f}}}_{a(2)};R;\frak b)).
\endaligned
\end{equation}
Note that the averaging of differential forms on $L(u)$ under the $T^n$-action coincides
with the harmonic projection $\Pi$.
In the construction of $\delta^{\frak b,b}$, we use ${\frak q}$-perturbation, which is $T^n$-equivariant. Thus 
$\delta^{\frak b,b}$ is $T^n$-equivariant.
Since the left hand side is $T^n$-invariant, we obtain \eqref{multiplicativeatosukosi}
after taking the average of both sides of \eqref{multiplicativeatosukosi20} under the $T^n$-action.
The proof of Lemma \ref{prodmainlemma} is now complete.
\end{proof}

The proof of Theorem \ref{multiplicative} is complete.
\qed
\par
\section{Surjectivity of Kodaira-Spencer map}
\label{sec:surf}

In this section, we will prove the following theorem.

\begin{thm}\label{surj}
The Kodaira-Spencer map
$\frak{ks}_{\frak b}: H(X;\Lambda_0) \to \text{\rm Jac}(\frak{PO}_{\frak b})$
is surjective.
\end{thm}
We begin with the corresponding result for the $\C$-reduction of
$\frak{ks}_{\frak b}$, which is Proposition \ref{specialfiber} below.
\par
We have the isomorphism
$$
H(X;\C)\cong {H(X;\Lambda_0)}/{H(X;\Lambda_+)}.
$$
We denote
$$
\overline{\text{\rm Jac}}(\frak{PO}_{\frak b}) = \frac{{\text{\rm Jac}}(\frak{PO}_{\frak b})}{
\Lambda_+{\text{\rm Jac}}(\frak{PO}_{\frak b})}.
$$
Then $\frak{ks}_{\frak b}$ induces a map
\begin{equation}\label{redksmap}
\overline{\frak{ks}}_{\frak b}: H(X;\C) \to \overline{\text{\rm Jac}}(\frak{PO}_{\frak b}),
\end{equation}
which is a ring homomorphism by Theorem \ref{multiplicative}.
(Here $H(X;\C)$ is a ring by the classical cup product without quantum corrections.)
\begin{prop}\label{specialfiber}
\begin{enumerate}
\item The map $\overline{\frak{ks}}_{\frak b}$ is surjective.
\item If $\frak b \in \mathcal A(\Lambda_+)$, then it is an isomorphism.
\end{enumerate}
\end{prop}
\begin{rem}
Actually Proposition \ref{specialfiber}.2 holds without the assumption $\frak b \in \mathcal A(\Lambda_+)$.
The proof of this will be completed in Section \ref{sec:injgen}.
\end{rem}
We defined $G_0$ and $G(\frak b)$ in Definition \ref{defg0}.
By definition, $G_0 \subset G(\frak b)$.
We put $G(\frak b) = \{\lambda_0,\lambda_2,\dots\}$, $0=\lambda_0 < \lambda_1 < \cdots$.
\begin{proof}[Proof of Proposition \ref{specialfiber}]
By \cite[Lemma 3.10]{fooo09}, the elements $z_j$ ($j=1,\dots,m$)
$\mod \Lambda_+$ generate $\overline{\text{\rm Jac}}(\frak{PO}_{\frak b})$.
\par
Let $\frak b = \sum_{a=1}^B \frak b_a \text{\bf f}_a$,
$\frak b_a \equiv \frak b_{a,0} \mod \Lambda_+$ with $\frak b_{a,0} \in \C$.
\begin{lem}\label{leadingterm}
There exists $\lambda_1 > 0$ such that
$$
\frak{PO}_{\frak b} \equiv \sum_{j=1}^m e^{\frak b_{j,0}} z_j
+
\sum_{j, a^j_1,\dots,a^j_{l_j}}\frac{1}{l_j!} \frak b_{a_1^j,0}\cdots \frak b_{a_{l_j}^j,0}z^{\beta_j} (\text{\rm ev}_0)_*
(\mathcal M_{1;l_j}(\beta_j;\text{\bf f}_j))
\mod T^{\lambda_1}.
$$
Notations appearing in the second term is as follows.
$$
\text{\bf f}_j = \text{\bf f}_{a_1^j} \otimes \dots \otimes  \text{\bf f}_{a_{l_j}^j},
$$
and
$$
\text{\bf f}_{a_i^j} = \text{\bf f}_{a_i^j(1)} \cup \dots \cup \text{\bf f}_{a_i^j(d_{i,j})}
$$
with $a_i^j(l) \in \{1,\dots,m\}$, $2d_{i,j} = \deg \text{\bf f}_{a_i^j}$. We put
$d_j = d_{1,j}+\dots+d_{l_j,j}$ and
$$
\beta_j = \beta_{c_j(1)} + \dots + \beta_{c_j(d_{j})}
\quad
z^{\beta_j} = z_{c_j(1)}\cdots z_{c_j(d_{j})}.
$$
Also
$$
a^j_1 = \cdots = a^j_{l_j} \in \{1,\ldots,m\}
$$
never occurs. (This case corresponds to the first term.)
\end{lem}
\begin{proof}
Let $\lambda_1$ be the smallest nonzero element of $G(\frak b)$.
We consider $\beta\in H_2(X;L(u);\Z)$, $\beta\ne 0$ with $\mathcal M_{1;\ell}(\beta) \ne \emptyset$.
By  \cite[Theorem 11.1]{fooo08}, we find
$$
\beta = \sum_{j=1}^m k_j \beta_j + \alpha
$$
$k_j \ge 0$, $\sum k_j >0$ and $\alpha$ is realized by
a sum of holomorphic spheres.
The term induced by  $\mathcal M_{1;\ell}(\beta)$ is
$$
c T^{\beta \cap \omega/2\pi} y(u)_1^{k'_1}\cdots y(u)_n^{k'_n}
= c T^{\alpha \cap \omega/2\pi}z_{1}^{k_1}\cdots z_m^{k_m},
$$
where $c \in \C$.
Therefore it suffices to consider the case $\alpha = 0$.
(In fact otherwise $\alpha \cap \omega/2\pi \ge \lambda_1$.)
In this case the Maslov index of $\beta$ is $2\sum_{i=1}^m k_i$.
By dimension counting we have $d_j = \sum_{i=1}^m k_i$, $j=1,\cdots ,m$.
Therefore it suffices to consider the case $\beta = \beta_{a(1)} + \dots + \beta_{a(d)}$.
In case $\beta = \beta_j$, we have
$$
(\text{\rm ev}_0)_*
(\mathcal M_{1;1}(\beta_j;\text{\bf f}_j)) = 1
$$
by  \cite[Theorem 11.1]{fooo08}.
The lemma follows from (\ref{PPformula}).
\end{proof}
We now prove Proposition \ref{specialfiber}.2.
In case $\frak b \in\mathcal A(\Lambda_+)$, Lemma \ref{leadingterm} implies
$$
\frak{PO}_{\frak b} \equiv z_1 +\dots + z_m \mod \Lambda\langle\!\langle y,y^{-1}\rangle\!\rangle_+^P.
$$
Therefore we have
\begin{equation}
y_i \frac{\partial \frak{PO}_{\frak b}}{\partial y_i}
\equiv
\sum_{j=1}^m v_{j,i}z_j
\mod \Lambda\langle\!\langle y,y^{-1}\rangle\!\rangle_+^P ,
\end{equation}
where $v_{j,i} \in \Z$ is defined by $d\ell_j=(v_{j,1},\dots,v_{j,n})$.
(Here $\ell_ j : P \to \R$ is the affine function as in 
(\ref{defpolytope}).  Note we use the symbol $l_j$ 
to denote the number of interior marked points in Lemma \ref{leadingterm}.
We used the symbol $\ell$ for the number of interior marked points 
in several other places.
We changed the notation here to avoid confusion with the affine function 
$\ell_j$.)
\par
On the other hand, it is easy to see that
$z_j$ satisfies the Stanley-Reisner relation
modulo $\Lambda\langle\!\langle y,y^{-1}\rangle\!\rangle_+^P$.
Furthermore
\begin{equation}\nonumber
\overline{\frak{ks}}_{\frak b}(\text{\bf f}_j) \equiv z_j
\mod \Lambda\langle\!\langle y,y^{-1}\rangle\!\rangle_+^P .
\end{equation}
Now it is a consequence of classical result of Stanley
that
$\overline{\frak{ks}}_{\frak b}$ is an isomorphism.
(See \cite[p106 Proposition]{fulton}, for example.)
\par
Next we prove Proposition \ref{specialfiber}.1.
Actually, we will prove a slightly sharper result.
We first recall the following:
\begin{conv}
For $R \in \Lambda_0[[Z_1,\dots,Z_m]]$, we substitute $z_j$ into $Z_j$ and denote
$$
R(z) = R(z_1,\dots,z_m).
$$
This is well-defined by Lemma \ref{surjhomfromz}.
\end{conv}

\begin{prop}\label{redsurjstrong}
Let $R \in \Lambda\langle\!\langle y,y^{-1}\rangle\!\rangle_0^{\overset{\circ}P}$.
Then there exist $\rho \in H(X;\Lambda_0)$ and $U \in
\Lambda\langle\!\langle y,y^{-1}\rangle\!\rangle_0^{\overset{\circ}P}$ such that
$$
R - T^{\lambda_1}U = {\frak ks}_{\frak b}(\rho)
$$
in $\text{\rm Jac}(\frak{PO}_{\frak b})$.
\end{prop}

\begin{proof}
By Lemma \ref{leadingterm} we have
\begin{equation}\label{imksformal}
{\frak{ks}}_{\frak b}(\text{\bf f}_j) \equiv c_jz_j
+ P_j(z)
\mod T^{\lambda_1}\Lambda\langle\!\langle y,y^{-1}\rangle\!\rangle_0^{\overset{\circ}P}
\end{equation}
for $j = 1,\dots,m$, $c_j \in \C \setminus\{0\}$ and
$P_j$ are formal power series of $Z_j$ with
$\Lambda_0$ coefficients, such that each term of it has degree $\ge 2$.
(In fact, the second term of Lemma \ref{leadingterm} has degree $\ge 2$ in $z_i$.)
\par
Using the fact that ${\frak{ks}}_{\frak b}$ is a ring homomorphism
and that $\text{deg} \text{\bf f}_j=2$ and $\prod \text{\bf f}_j^{r_j} = 0$
(where $\prod$ is the ordinary cup product)
whenever $\sum r_j \geq \dim X/2+1$, we can inductively apply (\ref{imksformal}) to find
$\rho \in H(X;\Lambda_0)$, $R'\in \Lambda_0[[Z_1,\dots,Z_m]]$ such that
\begin{equation}\label{eq:R-ksbrho}
R(z) - {\frak{ks}}_{\frak b}(\rho) = R'(z)
+ T^{\lambda_1} U_1
\end{equation}
in ${\text{\rm Jac}}(\frak{PO}_{\frak b})$, $U_1 \in \Lambda
\langle\!\langle y,y^{-1}\rangle\!\rangle_0^{\overset{\circ}P}$ and the degree of each term of $R'$
is greater than $\dim X/2$.

\begin{lem}\label{hdegcans}
If the degree of each monomial in $R' \in
\Lambda_0[[Z_1,\dots,Z_m]]$ is greater than $\dim X/2$, then there
exists $U_2 \in \Lambda\langle\!\langle y,y^{-1}\rangle\!\rangle_0^{\overset{\circ}P}$ such that
the following equality holds in $\text{\rm Jac}(\frak{PO}_{\frak b})$.
$$
R'(z) = T^{\lambda_1}U_2.
$$
\end{lem}
\begin{proof}
We put
$$
c_N = \left[\frac{2N}{\dim X+1}\right]
$$
where $[a]$ is the greatest integer less than equal to $a$. We note that
$c_N \to \infty$ as $N \to \infty$.
\begin{sublem}\label{sublemmamendou}
Let $N$ be a positive integer with $N > \dim X/2$. If the degree of
each term of $R(N) \in \Lambda_0[[Z_1,\dots,Z_m]]$ is greater than
$N$, then there exist $R(N+1) \in \Lambda_0[[Z_1,\dots,Z_m]]$ and
$U(N) \in \Lambda\langle\!\langle y,y^{-1}\rangle\!\rangle_0^{\overset{\circ}P}$ such that the
degree of each term of $R(N+1)$ is greater than $N+1$ and that the
equality
$$
R(N+1)(z) - R(N)(z) = T^{\lambda_1c_N} U(N)
$$
holds in $\text{\rm Jac}(\frak{PO}_{\frak b})$.
\end{sublem}
\begin{proof}
We may assume that  $R(N)$ is a monomial $R_0$ of degree $N$. We can then write
$$
R_0 = R_1 R_2 \cdots R_{c_N}
$$
so that each $R_i$ is a monomial of degree $> \dim X/2$.

Then applying (\ref{imksformal}) for each $i = 1, \ldots, c_N$, we obtain a polynomial ${\mathcal R}_i(z)$
such that each monomial appearing in ${\mathcal R}_i(z)$ is of degree
greater than that of $R_i(z)$ and satisfies
$$
\frak{ks}_{\frak b}(R_i(\text{\bf f}_1/c_1,\dots,\text{\bf f}_m/c_m))
- R_i(z) = {\mathcal R}_i(z) + T^{\lambda_1}V_i
$$
in $\text{\rm Jac}(\frak{PO}_{\frak b})$.
Since the ordinary cup product of
$1+(\dim X/2)$ elements vanishes on $H(X;\C)$, there exists a polynomial $Q_i(z)$
such that
$$
\frak{ks}_{\frak b}(R_i(\text{\bf f}_1/c_1,\dots,\text{\bf f}_m/c_m)) = T^{\lambda_1} Q_i(z).
$$
We now set
$$
R(N+1)(z)=-\prod_{i=1}^{c_N} (R_i(z) + {\mathcal R}_i(z)) +  R(N)(z).
$$
Then we find that the degree of each monomial in $R(N+1)(z)$ is bigger than
that of $R(N)(z)$ and that
$$
R(N+1)(z)-R(N)(z) = -\prod_{i=1}^{c_N} (R_i(z) + {\mathcal R}_i(z))
= -T^{c_N \lambda_1} \prod_{i=1}^{c_N}(V_i(z)+Q_i(z)).
$$
By setting $U(N) = -\prod_{i=1}^{c_N}(V_i(z)+Q_i(z))$, we have finished the proof.
\end{proof}

We apply Sublemma \ref{sublemmamendou} inductively by putting $R' = R(N_0)$.
($N_0 = \dim X/2 + 1$.)
We have
$$
R(N)(z) = R(N_0)(z) - \sum_{k=N_0}^{N-1} T^{\lambda_1c_k} U(k).
$$
Furthermore
$$
\lim_{N \to \infty}\sum_{k=N_0}^{N} T^{\lambda_1c_k} U(k)
$$
converges. Moreover $R(N)(z)$ converges to $0$ since
$R(N)$ converges to zero in the formal power series ring.
Therefore $R'(z) = R(N_0)(z)$ can be divided by $T^{\lambda_1}$ in
$\text{\rm Jac}(\frak{PO}_{\frak b})$. This finishes the proof of Lemma \ref{hdegcans}.
\end{proof}
We have finished the proof of Proposition \ref{redsurjstrong}
by applying Lemma \ref{hdegcans} to \eqref{eq:R-ksbrho}.
\end{proof}
\par
Proposition \ref{specialfiber}.1 is an immediate consequence of
Proposition \ref{redsurjstrong}.
\end{proof}
\par
Now we are ready to complete the proof of Theorem \ref{surj}.
\begin{proof}
Let
$R \in \Lambda\langle\!\langle y,y^{-1}\rangle\!\rangle_0^{P} $.
We use Proposition \ref{redsurjstrong} to find
$\rho_1 \in H(X;\Lambda_0)$ and $R_1$ such that
$
R - \frak{ks}_{\frak b}(\rho_1) = T^{\lambda_1}R_1.
$
We find inductively $\rho_k \in H(X;\Lambda_0)$ and $R_k$
such that
$
R_{k-1} - \frak{ks}_{\frak b}(\rho_k) = T^{\lambda_1}R_k.
$
Clearly
$
\rho =
\sum_{k=1}^{\infty} T^{\lambda_1(k-1)}\rho_k
$
converges and
$
R = \frak{ks}_{\frak b}(\rho).
$
\end{proof}
\section{Versality of the potential function with bulk}
\label{sec:versality}
Let $z_j$ $(j=1,\dots ,m)$ be as in Definition \ref{zidef}.\index{coordinate change of $y_i$}
\begin{thm}\label{versality}
Let $\frak P \in \Lambda\langle\!\langle y,y^{-1}\rangle\!\rangle_0^{\overset{\circ}P}$. We
assume that $\frak P$ is $G$-gapped for a discrete monoid $G$
containing $G_0$.
\begin{enumerate}
\item If $\frak P \equiv z_1+\dots +z_m \mod
\Lambda\langle\!\langle y,y^{-1}\rangle\!\rangle_+^{\overset{\circ}P}$, then there exists
a $G$-gapped element $\frak b \in \mathcal A(\Lambda_+)$ and a
$G$-gapped strict coordinate change $y'$ \emph{on $\text{\rm Int} P$}
such that
$$
\frak P(y') = \frak{PO}_{\frak b}(y).
$$
\item
If $\frak P \equiv c_1z_1+\dots +c_mz_m \mod
\Lambda\langle\!\langle y,y^{-1}\rangle\!\rangle_{+}^{\overset{\circ}P}$ with $c_i \in \C\setminus \{0\}$,
then there exist a $G$-gapped element
$\frak b \in \mathcal A(\Lambda_0)$ and a $G$-gapped coordinate change $y'$ on
$\text{\rm Int} P$ such that
$$
\frak P(y') = \frak{PO}_{\frak b}(y).
$$
\end{enumerate}
If $\frak P \in \Lambda\langle\!\langle y,y^{-1}\rangle\!\rangle_0^{P}$, then $y'$ can be taken to
be a strict $G$-gapped coordinate change \emph{on $P$}.
\end{thm}
(Note that $\Lambda\langle\!\langle y,y^{-1}\rangle\!\rangle_{+}^{\overset{\circ}P}$ is
the ideal of $\Lambda\langle\!\langle y,y^{-1}\rangle\!\rangle_{0}^{\overset{\circ}P}$
which is defined in Definition \ref{def24}.3.)
\begin{proof}
We first prove Theorem \ref{versality}.1. We focus on the case of $\frak P \in
\Lambda\langle\!\langle y,y^{-1}\rangle\!\rangle_0^{\overset{\circ}P}$ since the case of $\frak P
\in \Lambda\langle\!\langle y,y^{-1}\rangle\!\rangle_0^{P}$ is similar. We enumerate the elements
of $G$ and write $G = \{\lambda_0, \lambda_1, \dots\}$ so that
$\lambda_0  = 0$, $\lambda_i < \lambda_{i+1}$. We note that $\lambda_i \to \infty$ as
$i \to \infty$.
\par
We will construct a cycle $\frak b(k) \in \CA(\Lambda_+)$ and a
coordinate change $y'(k)$ inductively over $k$ with the following
properties:
\begin{enumerate}
\item They are $G$-gapped.
\item $y'(k)$ is a $G$-gapped and strict coordinate change on $\text{\rm Int}P$.
\item $y'(k)$ satisfies
$$
\frak P(y'(k)) - \frak{PO}_{\frak b(k)}(y) \equiv 0
\mod T^{\lambda_{k+1}}
\Lambda\langle\!\langle y,y^{-1}\rangle\!\rangle_0^{\overset{\circ}P}.
$$
\item $(y'(k+1)_i - y'(k)_i)/y_i \equiv 0 \mod T^{\lambda_{k+1}}
\Lambda\langle\!\langle y,y^{-1}\rangle\!\rangle_0^{\overset{\circ} P}$.
\end{enumerate}
\par
When $k = 0$, we set $\frak b(0) = 0$ and $y'(0) = y$. Then the hypothesis
$\frak P \equiv z_1+ \dots + z_m \mod \Lambda\langle\!\langle y,y^{-1}\rangle\!\rangle^{\overset{\circ}P}_+$
implies the above properties for $k = 0$.
\par
Now suppose we have constructed $\frak b(k)$ and $y'(k)$ so that
$$
\frak P(y'(k)) - \frak{PO}_{\frak b(k)}(y)
\equiv T^{\lambda_{k+1}}\mathcal P
\mod T^{\lambda_{k+2}} \Lambda\langle\!\langle y,y^{-1}\rangle\!\rangle_0^{\overset{\circ}P}
$$
holds.
By Proposition \ref{specialfiber} we may take $\frak z_{k+1} \in H(X;\C)$,
$W(k,i) \in \C[[Z_1,\dots,Z_m]]$ ($i=1,\dots,n$) so that they satisfy
$$
\mathcal P - \frak{ks}_{\frak b(k)}(\frak z_{k+1}) \equiv \sum_{i=1}^n
W(k,i)(z) y_i\frac{\partial\frak{PO}_{\frak b(k)}} {\partial y_i}
\mod T^{\lambda_{1}} \Lambda\langle\!\langle y,y^{-1}\rangle\!\rangle_0^{\overset{\circ}P}.
$$
(Here we use the fact that everything is $G$-gapped.) We define
$$
\frak b(k+1) = \frak b(k) + T^{\lambda_{k+1}}\frak z_{k+1},
$$
and
$$
y'(k+1)_i = y'(k)_i -  T^{\lambda_{k+1}}W(k,i)(z) y_i, \quad i =
1,\dots, n.
$$
It is easy to check Conditions 1 - 4 above for these choices. By
construction, the limits $\frak b = \lim_{k\to \infty} \frak b(k)$,
$y' = \lim_{k\to \infty} y'(k)$ exist and satisfy the properties
required in the statement of Theorem \ref{versality}.1.
\par
We next prove Theorem \ref{versality}.2. After applying an explicit coordinate change 
$\frak b(0) = \sum (\log c_i) \text{\bf f}_i$, we have
$$
\frak P(y'(0)) = \frak{PO}_{\frak b(0)} (y)\mod T^{\lambda_{1}}
\Lambda\langle\!\langle y,y^{-1}\rangle\!\rangle_0^{\overset{\circ}P} .
$$
Once we have arranged this, we can adapt the same inductive argument as for
the case Theorem \ref{versality}.1.
\end{proof}
\begin{rem}
We do not attempt to remove the gapped-ness from the hypotheses of
Theorem \ref{versality} since non-gapped elements of
$\Lambda\langle\!\langle y,y^{-1}\rangle\!\rangle_{0}^{\overset{\circ}P}$ do not seem  to appear
in geometry.
\end{rem}

\section{Algebraization of Jacobian ring}
\label{sec:algebra}
In this section we prove that there exists a coordinate change
transforming the potential function to a Laurent \emph{polynomial}.
\begin{defn}
We denote \index{$\Lambda[y,y^{-1}]_0^P$}
$$
\Lambda[y,y^{-1}]_0^P
= \Lambda[y,y^{-1}] \cap \Lambda\langle\!\langle y,y^{-1}\rangle\!\rangle_0^P
= \{w \in \Lambda[y,y^{-1}] \mid \frak v_T^P(w) \ge 0\}.
$$
We note that $\Lambda[y,y^{-1}]_0^P$ coincides with
$$
\Lambda[y,y^{-1}] \cap \Lambda\langle\!\langle y,y^{-1}\rangle\!\rangle_0^{\overset{\circ}P} .
$$
In other words, $\Lambda[y,y^{-1}]_0^P$ is the image of
$\Lambda_0[Z_1,\dots,Z_m]$ under the homomorphism in Lemma
\ref{surjhomfromz}. 
We also emphasize that
$
\Lambda[y,y^{-1}]_0^P \ne \Lambda_0[y,y^{-1}].
$ 
For $\epsilon >0$ we put
$$
P_{\epsilon}
= \{ u\in \R^n \mid \ell_j(u) \ge \epsilon,  ~ j=1,\dots,m \}.
$$
See Section \ref{sec:valuation}.
(Note
$
P
= \{ u\in \R^n \mid \ell_j(u) \ge 0, ~ j=1,\dots,m \}.
$)
\end{defn}
\begin{thm}\label{algebraization}
\begin{enumerate}
\item If $\frak b \in \mathcal A(\Lambda_0)$, then
for each $\epsilon>0$, there exists a strict
coordinate change $y'$ on $P_{\epsilon}$ such that\index{coordinate change of $y_i$}
\begin{equation}\label{algebraika}
\frak{PO}_{\frak b}(y') \in \Lambda[y,y^{-1}]^{P_{\epsilon}}_0.
\end{equation}
\item If $\frak b \in \mathcal A(\Lambda_+)$, then $y'$ can
be taken to be a strict coordinate change on $P$ and
$$
\frak{PO}_{\frak b}(y') \in \Lambda[y,y^{-1}]^{P}_0.
$$
If $\frak b$ is
$G$-gapped, so is $y'$.
\end{enumerate}
\end{thm}
\begin{proof}
We first prove Theorem \ref{algebraization}.2.
Let $\frak b \in \mathcal A(\Lambda_+)$ be $G$-gapped for
a discrete submonoid $G \supseteq G_0$. We put $G =
\{\lambda_0,\lambda_1,\dots\}$ as before.
We first consider a polynomial
$$
\frak{P}_0 = z_1 + \dots + z_m.
$$
This satisfies
$\frak{PO}_{\frak b}\equiv \frak{P}_0
\mod T^{\lambda_{1}}\Lambda\langle\!\langle y,y^{-1}\rangle\!\rangle_0^{P} $ for
$\frak b \in \mathcal A(\Lambda_+)$. 
We put
$$
\overline{\text{\rm Jac}}(\frak P_0)
= \frac{\text{\rm Jac}(\frak P_0)}{\Lambda_+ \text{\rm Jac}(\frak P_0)}.
$$
Since the property
$\frak{PO}_{\frak b}\equiv \frak{P}_0
\mod T^{\lambda_{1}}\Lambda\langle\!\langle y,y^{-1}\rangle\!\rangle_0^{P}$ implies the equality of ideals
in $\Lambda\langle\!\langle y,y^{-1}\rangle\!\rangle_0^{P}$
$$
\Lambda\langle\!\langle y,y^{-1}\rangle\!\rangle_+^{P}
+
\left( y_i \frac{\partial
\frak{PO}_{\frak b}}{\partial y_i} ; i=1,\dots,n
\right)
=
\Lambda\langle\!\langle y,y^{-1}\rangle\!\rangle_+^{P}
+
\left( y_i
\frac{\partial \frak{P}_{0}}{\partial y_i} ; i=1,\dots,n
\right),
$$
we find
$$
\overline{\text{\rm Jac}}(\frak {PO}_{\frak b})
\cong\overline{\text{\rm Jac}}(\frak P_0).
$$
(However, there is no canonical map between
$\text{\rm Jac}(\frak {PO}_{\frak b})
$ and $\text{\rm Jac}(\frak P_0)$.)
Through this isomorphism,
we have a ring isomorphism which we also denote by
$\overline{\frak{ks}}$:
$$
\overline{\frak{ks}}~:~ H(X;\C) \to
\overline{\text{\rm Jac}}(\frak {PO}_{\frak b})
\cong
\overline{\text{\rm Jac}}
(\frak P_0)
$$
which sends
$\overline{\text{\bf f}}_i$ to $\overline{z}_i$.
Then we have
\begin{lem}\label{polygenerate}
The polynomials of degree $\le \dim X/2$ of
$\overline{z}_1,\dots,\overline{z}_m$
generate $\overline{\text{\rm Jac}}(\frak P_0)$.
\end{lem}

We now inductively construct a coordinate change $y'(k)$ and a polynomial $B_k
\in \Lambda_0[Z_1,\dots,Z_m]$ for $k = 0, \dots $ with the
following properties:
\begin{enumerate}
\item $B_k$ is a polynomial of degree $\le \dim X/2$.
It is $G$-gapped.
\item $y'(k)$ is a $G$-gapped strict coordinate change on $P$.
\item
$$
\frak{PO}_{\frak b}(y'(k)) - (\frak P_0 + B_k(z)) \equiv 0
\mod T^{\lambda_{k+1}}\Lambda\langle\!\langle y,y^{-1}\rangle\!\rangle_0^{P} .
$$
\item
$$
\frac{y'(k+1)_i - y'(k)_i}{y_i} \equiv 0 \mod T^{\lambda_{k+1}}\Lambda\langle\!\langle y,y^{-1}\rangle\!\rangle_0^{P} .
$$
\item
$$
B_{k+1} - B_k \equiv 0
\mod T^{\lambda_{k+1}}\Lambda_0[Z_1,\dots,Z_m].
$$
\end{enumerate}

Theorem \ref{versality} enables us to choose $B_0=0$. Now suppose we
have constructed $B_\ell$ and $y'(\ell)$ for $\ell =k$. Then
$$
\frak{PO}_{\frak b}(y'(k)) - (\frak P_0 + B_k(z))
\equiv T^{\lambda_{k+1}} \mathcal P_k(z)
\mod T^{\lambda_{k+2}}\Lambda\langle\!\langle y,y^{-1}\rangle\!\rangle_0^{P}
$$
for some $\mathcal P_k \in \C[Z_1,\dots,Z_m]$.
(Here Lemmata \ref{PPnotPO} and \ref{surjhomfromz} show
that we can take $\mathcal P_k$ as a polynomial.)
By Lemma
\ref{polygenerate} there exist polynomials $\Delta B_{k+1} \in
\C[Z_1,\dots,Z_m]$ of degree $\le \dim X/2$ and $W_{k+1,i} \in
\C[Z_1,\dots,Z_m]$ such that
$$
\mathcal P_k
\equiv \Delta B_{k+1}(z) + \sum_{i=1}^n  W_{k+1,i}(z)
y_i\frac{\partial \frak P_0}{\partial y_i}
\mod T^{\lambda_1}\Lambda\langle\!\langle y,y^{-1}\rangle\!\rangle_0^{P} .
$$
We put
$$
y'(k+1)_i = y'(k)_i - T^{\lambda_{k+1}}W_{k+1,i}(z)y_i,
\qquad
B_{k+1} = B_k + T^{\lambda_{k+1}}\Delta B_{k+1}.
$$
It is easy to see that they have the required properties and hence we
finish the construction of $B_k$ and
$y'(k)$ for $k \in \Z_{\geq 0}$.

Now we consider the limits $B = \lim_{k\to\infty} B_{k}(z)$ and
$y'_{\infty} = \lim_{k\to\infty} y'(k)$ which converge
by
construction.
We also find that
$y'_{\infty} \in
\Lambda\langle\!\langle y,y^{-1}\rangle\!\rangle^{P}$.
\par
If we put $\frak B = B(z) \in \Lambda[y,y^{-1}]_0^{P}$, the equality
$\frak{PO}_{\frak b}(y'_{\infty}) = \frak{P}_{0} + \frak{B}$ follows. We have thus
proved Thorem \ref{algebraization}.2.
\par
We next prove Thorem \ref{algebraization}.1.
By Lemma \ref{leadingterm}, we may put
\begin{equation}\label{PObLT0}
\frak{PO}_{\frak b}  \equiv c_1 z_1 + \cdots + c_mz_m  + R(z_1,\dots,z_m) \mod \Lambda[y,y^{-1}]_+^P
\end{equation}
where $R \in \Lambda_0[[Z_1,\dots,Z_m]]$ each summand of which has degree $\ge 2$,
and $c_i \in \C\setminus \{0\}$.
We put
$$
z_i^{\epsilon} = T^{-\epsilon}z_i.
$$
Note that\index{$\frak v_T^{P_{\epsilon}}$}
$$
\frak v_T^{P_{\epsilon}}(z_i^{\epsilon})
=\inf _{u\in P_{\e}}\frak v_T^{u}(z_i^{\epsilon})
=0.
$$
In fact, $z_i^{\epsilon} $ plays a role of $z_i$
for our polytope $P_{\epsilon}$.
There exists a surjective homomorphism
$$
\Lambda_0\langle\!\langle Z_1,\dots,Z_m \rangle\!\rangle
\to \Lambda\langle\!\langle  y,y^{-1}\rangle\!\rangle_{0}^{P_{\epsilon}}
$$
which sends $Z_i$ to $z_i^{\epsilon}$.
We consider
$
T^{-\epsilon}\frak{PO}_{\frak b} \in  \Lambda\langle\!\langle  y,y^{-1}\rangle\!\rangle_{0}^{P_{\epsilon}}.
$
From (\ref{PObLT0}) it is easy to see that
\begin{equation}\label{Tepsilonb}
T^{-\epsilon}\frak{PO}_{\frak b}
\equiv c_1z^{\epsilon}_1 + \dots + c_mz^{\epsilon}_m
\mod  \Lambda\langle\!\langle  y,y^{-1}\rangle\!\rangle_{+}^{P_{\epsilon}}.
\end{equation}
Then we can repeat the proof of Thorem \ref{algebraization}.2 and obtained the conclusion.
\end{proof}

\begin{rem}
Theorem \ref{algebraization} should be related to the classical results
in the deformation theory of isolated singularity. See \cite{artin},
\cite{Elkik}.
\end{rem}

We next discuss the relationships between several variations of the
Jacobian rings.
\begin{conds}\label{alegbraikacond}
Let $\frak P \in \Lambda[y,y^{-1}]^P_0$. We assume that $ \frak P -
\sum_{i=1}^{m} c_iz_i $ is contained in the ideal
$\Lambda[y,y^{-1}]_+^P$.
\end{conds}
\begin{rem}
For $\frak b \in \mathcal A(\Lambda_+)$,
$\frak {PO}_{\frak b}$ satisfies this condition.
For $\frak b \in \mathcal A(\Lambda_0)$,
Formula (\ref{Tepsilonb}) implies that the condition is
satisfied for $T^{-\epsilon}\frak{PO}_{\frak b} $ if we replace
$P$ and $z_i$ by $P_{\epsilon}$ and $z_i^{\epsilon}$ respectively.
\end{rem}
\begin{defn}
Suppose that $\frak P$ satisfies Condition $\ref{alegbraikacond}$.
We put
$$
\text{\rm PJaco}(\frak P) = \frac{\Lambda[y,y^{-1}]^P_0}{\left(
y_i\frac{\partial \frak P}{\partial y_i}\right)},
\qquad
\text{\rm PJac}_{\frak v^P_T}(\frak P) =
\frac{\Lambda[y,y^{-1}]^P_0}{\text{Clos}_{\frak v_T^{P}}\left(
y_i\frac{\partial \frak P}{\partial y_i}\right)},
$$
$$
\text{\rm Jaco}_{\frak v^P_T}(\frak P) = \frac{\Lambda\langle\!\langle y,y^{-1}\rangle\!\rangle_0^{P} }{\left(
y_i\frac{\partial \frak P}{\partial y_i}\right)},
\qquad
\text{\rm Jac}_{\frak v^P_T}(\frak P) = \frac{\Lambda\langle\!\langle y,y^{-1}\rangle\!\rangle_0^{P} }{\text{Clos}_{\frak v_T^{P}}\left(
y_i\frac{\partial \frak P}{\partial y_i}\right)},
$$
$$
\text{\rm Jaco}(\frak P) = \frac{\Lambda\langle\!\langle y,y^{-1}\rangle\!\rangle_0^{\overset{\circ}P} }{\left(
y_i\frac{\partial \frak P}{\partial y_i}\right)},
\qquad
\text{\rm PJac}(\frak P) = \frac{\Lambda[y,y^{-1}]^P_0}{\text{Clos}_{d_{\overset{\circ}P}}\left(
y_i\frac{\partial \frak P}{\partial y_i}\right)}.
$$
Here we recall that $\text{Clos}_{\frak v_T^{P}}$ means the closure
with respect to the norm $\frak v_T^{P}$ and
$\text{Clos}_{d_{\overset{\circ}P}}$ means the closure with respect
to the metric $d_{\overset{\circ}P}$. (See Definition \ref{toplamdaP0}.)
We also note that we already defined
$$
\text{\rm Jac}(\frak P) = \frac{\Lambda\langle
\!\langle y,y^{-1}\rangle\!\rangle_0^{\overset{\circ}P} }{\text{Clos}_{d_{\overset{\circ}P}}\left(
y_i\frac{\partial \frak P}{\partial y_i}\right)}
$$
before in Definition \ref{def28}.
\end{defn}
In the above notations of various Jacobian rings, the letter `P' in front of `Jac' stands for `polynomial',
the letter `o' after `Jac' stands for `taking the quotient \emph{not} by the closure' and the subscript
$\frak v^P_T$ stands for `using the topology with respect to the norm $\frak v^P_T$'.
\par
There are natural homomorphisms described in the following diagram.
\begin{equation}
\xymatrix{ & \text{\rm PJaco}(\frak P)\ar[d]\ar[ddr] & & \\
& \text{\rm PJac}(\frak P) \ar[rr] &  &\text{\rm Jac}(\frak P) \\
& &\text{\rm Jaco}_{\frak v^P_T}(\frak P) \ar[d] & \\
\text{\rm PJac}_{\frak v^P_T}(\frak P)\ar[uur] \ar[rr] & &
\text{\rm Jac}_{\frak v^P_T}(\frak P)\ar[uur]}
\end{equation}
\par
\begin{lem}\label{polyapprox}
Suppose $\frak P$ satisfies Condition \ref{alegbraikacond}.
Then we have the following:
\begin{enumerate}
\item The homomorphism $\text{\rm PJac}_{\frak v^P_T}(\frak P) \to \text{\rm Jac}_{\frak v^P_T}(\frak P)$
is an isomorphism.
\item The homomorphism $\text{\rm PJac}(\frak P) \to \text{\rm Jac}(\frak P)$
is an isomorphism.
\item
The homomorphism $\text{\rm Jaco}_{\frak v^P_T}(\frak P) \to \text{\rm Jac}_{\frak v^P_T}(\frak P)$
is an isomorphism if $\text{\rm Jaco}_{\frak v^P_T}(\frak P)$ is torsion free.
\item
The homomorphism $\text{\rm Jac}_{\frak v^P_T}(\frak P) \to \text{\rm Jac}(\frak P)$
is an isomorphism if $\text{\rm Jac}_{\frak v^P_T}(\frak P)$ is torsion free.
\end{enumerate}
\end{lem}
\begin{rem}
\begin{enumerate}
\item
Theorem \ref{Mirmain} implies that  $\text{\rm Jac}(\frak P)$
is always torsion free. However, while we are proving
Theorem \ref{Mirmain}, we need to check the assumption to use Lemma \ref{polyapprox}.4.
\item Finally (at the end of Section \ref{sec:injgen}) we have
$$
\text{\rm PJac}(\frak P) = \text{\rm Jac}(\frak P) = \text{\rm PJac}_{\frak v^P_T}(\frak P) = \text{\rm Jac}_{\frak v^P_T}(\frak P) =
\text{\rm Jaco}_{\frak v^P_T}(\frak P).
$$
On the other hand, $\text{\rm PJaco}(\frak P) \ne \text{\rm PJac}(\frak P)$ in general.
(See Example 8.2 \cite{fooo08}.)  We will not further study the
relationship between $\text{\rm Jaco}(\frak P)$ and $\text{\rm Jac}(\frak P)$, since such
a study will not be needed in this paper.
\end{enumerate}
\end{rem}
\begin{proof}[Proof of Lemma \ref{polyapprox}]
We choose $G \supset G_0$ such that
$\frak P$ is $G$-gapped and set $G = \{\lambda_0,\lambda_1,
\dots\}$.
\par
We first prove Lemma \ref{polyapprox}.2. We remark that Lemma \ref{polygenerate} holds
for $\overline{\text{\rm Jac}}(\frak P)$ also. (We use Theorem \ref{versality}
to prove it.) Moreover Sublemma \ref{sublemmamendou} also holds for
${\text{\rm Jac}}(\frak P)$.
\par
Injectivity of the map in Lemma \ref{polyapprox}.2 is clear and so we now prove its
surjectivity. Let $\frak B \in {\text{\rm Jac}}(\frak P)$. By an induction over
$k$, we will find $B_k \in \Lambda_0[Z_1, \dots,Z_n]$ such that
\begin{enumerate}
\item $\deg B_k \le \dim X/2$.
\item $\frak B - B_k(z) \in T^{\lambda_k}\text{\rm Jac}(\frak P)$.
\item $B_{k+1} - B_{k} \in T^{\lambda_k}\Lambda_0[[Z_1,
\dots,Z_n]]$.
\end{enumerate}
Suppose we have constructed $B_k$ and put
$$
\frak B - B_k(z) = T^{\lambda_k} F_k(z)
$$
with $F_k \in \Lambda_0[[Z_1,
\dots,Z_n]]$. Let $F'_k$ be the sum of the terms of
$F_k$ of degree $\le \dim X/2$.
Then
$$
F_k(z) - F'_{k}(z) \in T^{\lambda_1}\text{\rm Jac}(\frak P)
$$
by Sublemma \ref{sublemmamendou}. It is easy to see that
$$
B_{k+1} = B_k + T^{\lambda_k}F'_k
$$
has the required properties.
\par
Since $\deg B_k \leq \dim X/2$ is bounded and $\lambda_k \to \infty$, $B_k$ converges to a
\emph{polynomial} in $\Lambda_0[Z_1,\dots,Z_m]$ which we denote by
$B$. Since $B_k(z)$ converges to $\frak B$ in $\text{\rm Jac}(\frak P)$ by
Condition 2, we have proved $B(z) = \frak B$. By varying $\frak B$
in $\text{\rm Jac}(\frak P)$, the surjectivity of ${\rm PJac}(\frak P) \to \text{\rm Jac}(\frak
P)$ follows.
\par
The proof of Lemma \ref{polyapprox}.1 is easier and left for the reader.
\par
We next prove Lemma \ref{polyapprox}.3. This time the surjectivity is obvious and so
it suffices to show the injectivity.
Let $x \in \Lambda\langle\!\langle y,y^{-1}\rangle\!\rangle_0^P$ be an element which projects down to zero
in $\text{\rm Jac}_{\frak v^P_T}(\frak P)$, i.e., $x \in \text{Clos}_{\frak v_T^{P}}\left(y_i\frac{\del\frak P}{\del y_i}\right)$.
By definition of $\text{Clos}_{\frak v_T^{P}}\left(y_i\frac{\del\frak P}{\del y_i}\right)$,
there exists a sequence $x_{\ell} \in \left(
y_i\frac{\partial \frak P}{\partial y_i}\right)$ such that
$x - x_{\ell}$ converges to zero in $\frak v_T^P$ topology.
Take a monoid $G \supseteq G_0$ such that $\frak P$ and $x$ are $G$-gapped and put
$G = \{ \lambda_0,\lambda_1,\dots\}$.
\par
We will inductively construct $\frak Z_{k,i} \in \Lambda\langle\!\langle y,y^{-1}\rangle\!\rangle_0^{P}$
which is $G$-gapped and satisfies
$$
\aligned
x - \sum_{i=1}^n \frak Z_{k,i} y_i\frac{\partial \frak P}{\partial y_i}
&\equiv 0 \mod T^{\lambda_k}\Lambda\langle\!\langle y,y^{-1}\rangle\!\rangle_0^{P},
\\
\frak Z_{k+1,i} - \frak Z_{k,i} &\equiv 0 \mod T^{\lambda_k}\Lambda\langle\!\langle y,y^{-1}\rangle\!\rangle_0^{P}.
\endaligned$$
Suppose we have $\frak Z_{k,i}$. We put
\begin{equation}\nonumber
x -  \sum_{i=1}^n \frak Z_{k,i} y_i\frac{\partial \frak P}{\partial y_i}
= T^{\lambda_k}\frak Z(k)
\end{equation}
for some $\frak Z(k) \in \Lambda\langle\!\langle y,y^{-1}\rangle\!\rangle_0^{P}$.
We take a sufficiently large ${\ell}$ so that
$$
\frak v^P_T(x_{\ell} - x) > \lambda_{k+1}.
$$
Then there exists $\frak U_{k+1} \in \Lambda\langle\!\langle y,y^{-1}\rangle\!\rangle_0^P$ such that
$$
(x_{\ell} - x) = T^{\lambda_{k+1}}\frak U_{k+1}.
$$
Then
$$
x_{\ell} - \sum_{i=1}^n \frak Z_{k,i} y_i\frac{\partial \frak P}{\partial y_i}
= T^{\lambda_k}(\frak Z(k) + T^{\lambda_{k+1} - \lambda_k}\frak U_{k+1}).
$$
Since $x_\ell \in \left(y_i\frac{\del\frak P}{\del y_i}\right)$, this implies that
$
T^{\lambda_k}(\frak Z(k) + T^{\lambda_{k+1}-\lambda_k}\frak U_{k+1})
$
lies in $\left(y_i\frac{\del\frak P}{\del y_i}\right)$
and hence is zero in $\text{\rm Jaco}_{\frak v^P_T}(\frak P)$.
By the hypothesis that $\text{\rm Jaco}_{\frak v^P_T}(\frak P)$ is torsion free, we find that
$\frak Z(k) + T^{\lambda_{k+1}-\lambda_k}\frak U_{k+1}$ is zero in
$\text{\rm Jaco}_{\frak v^P_T}(\frak P)$.
Therefore there exists $Z(k,i) \in \Lambda_0\langle\!\langle Z_1,\dots,Z_m\rangle\!\rangle$ such that
$$
\frak Z(k) + T^{\lambda_{k+1}-\lambda_{k}}\frak U_{k+1}\equiv \sum_{i=1}^n Z(k,i)(z) y_i\frac{\partial \frak P}{\partial y_i}
\mod T^{\lambda_{k+1} -\lambda_k}\Lambda\langle\!\langle y,y^{-1}\rangle\!\rangle_+^{P} .
$$
Hence $\frak Z_{k+1,i}:= \frak Z_{k,i} + T^{\lambda_k}Z(k,i)(z)$,
has the required properties for $k+1$.
\par
If we put $\lim_{k\to\infty} \frak Z_{k,i} = \frak Z_i$, then we have
$$
x = \sum_{i=1}^n \frak Z_{i} y_i\frac{\partial \frak P}{\partial y_i}
\in \left(y_i\frac{\del\frak P}{\del y_i}\right)
$$
and hence $x = 0$ in $\text{\rm Jaco}_{\frak v^P_T}(\frak P)$ as required.
\par
Finally we prove Lemma \ref{polyapprox}.4. Surjectivity of the map immediately follows from
Lemma \ref{polyapprox}.2 and so the injectivity is the only non-obvious part of the proof. By the
torsion-freeness assumption of $\text{\rm Jac}_{\frak v^P_T}(\frak P)$,
it suffices to show the injectivity after taking
the tensor product $\otimes_{\Lambda_0}\Lambda$.
Then by the versality given in Theorem \ref{versality}, it is
enough to consider the case of the potential function $\frak P = \frak{PO}_{\frak b}$.
For this case, we use the decomposition given Proposition \ref{Morsesplit} for $\frak{PO}_{\frak b}$:
it is easy to see that $\text{\rm Jac}_{\frak v^P_T}(\frak{PO}_{\frak b})$ has the same
kind of decomposition such as the one for $\text{\rm Jac}(\frak{PO}_{\frak b})$. Recall
if $\frak y \in \text{\rm Crit}(\frak{PO}_{\frak b})$ and if we set
$u:= \frak v_T(\frak y)$ $u$ lies in $\text{\rm Int} P$.
Therefore
\begin{equation}\label{locandval}
\text{\rm Jac}_{\frak v^P_T}(\frak{PO}_{\frak b};\frak y)
\cong
\text{\rm Jac}(\frak{PO}_{\frak b};\frak y).
\end{equation}
Once we have this, we can obtain the factorization
$$
\frac{\Lambda\langle\!\langle y,y^{-1}\rangle\!\rangle_0^{P} }{\text{Clos}_{\frak v_T^P}\left(
y_i\frac{\partial \frak{PO}_{\frak b}}{\partial y_i}\right)}
\otimes_{\Lambda_0} \Lambda
\cong \prod_{\frak y \in \text{\rm Crit}(\frak{PO}_{\frak b})} \text{\rm Jac}(\frak{PO}_{\frak b};\frak y)
$$
as for $\text{\rm Jac}(\frak{PO}_{\frak b})$ in Proposition \ref{Morsesplit} in the same way.
(Note we take the closure with respect to the norm $\frak v_T^P$ in the left
hand side of the factorization.)
Therefore the injectivity in Lemma \ref{polyapprox}.4 follows because the last product is
isomorphic to $\text{\rm Jac}(\frak{PO}_{\frak b})\otimes_{\Lambda_0} \Lambda$ itself.
Now the proof of Lemma \ref{polyapprox} is complete.
\end{proof}
\par
\section{Seidel homomorphism and a result by McDuff-Tolman}
\label{sec:seidel}
In this section, we use the results by Seidel \cite{seidel:auto} and McDuff-Tolman \cite{mc-tol} to prove Theorem \ref{qhalgebraization}.
To describe the results of \cite{mc-tol}, \cite{seidel:auto}, we need some notations.
\par
Let $\text{\bf f}_i$ $(i=1,\dots,m)$ be the divisors
which generate $\mathcal A^2$. We put
$$
D_i = \text{\rm PD}(\text{\bf f}_i) \in H^2(X;\Z).
$$
Let $\mathcal P \subset \{1,\dots,m\}$ be a primitive
collection. (See for example  \cite[Definition 2.4]{fooo08}.)
We have
$
\bigcap_{i \in \mathcal P} D_i = \emptyset
$
and we have $\mathcal P' \subset \{1,\dots,m\} \setminus \mathcal P$
and $k_{i'}$ ($i' \in \mathcal P'$) such that
$$
\sum_{i\in \mathcal P}d\ell_i = \sum_{i'\in \mathcal P'}k_{i'}d\ell_{i'},
\qquad
\bigcap_{i' \in \mathcal P'}  D_{i'} \ne \emptyset.
$$
($k_i \in \Z_{\ge 0}$.) We define
$$
\omega(\mathcal P) = \sum_{i\in \mathcal P} \ell_i(u)
- \sum_{i'\in \mathcal P'} k_{i'}\ell_{i'}(u).
$$
Note this is independent of the choice of $u$ by  \cite[Lemma 6.2]{fooo08}.
We take $m$ formal variables $Z_i$.
\begin{defn}
The {\it quantum Stanley-Reisner relation associated with $\mathcal P$} is
\begin{equation}\label{QSRrel}
\prod_{i \in \mathcal P} Z_i - T^{\omega(\mathcal P)}
\prod_{i'\in \mathcal P'}Z_{i'}^{k_{i'}} = 0.
\end{equation}
We denote the left hand side of (\ref{QSRrel}) by $\frak I(\mathcal P)$.
The {\it quantum Stanley-Reisner ideal}\index{quantum Stanley-Reisner ideal} is the ideal of
$\Lambda_{0}^P[Z_1,\dots,Z_m]$ generated by $\frak I(\mathcal P)$
for various $\mathcal P$. We denote it by
$SR_{\omega}$.
\end{defn}
\begin{thm}[McDuff-Tolman \cite{mc-tol}]
\label{mcduftolman}
There exists an element $z'_i$ of $H(X;\Lambda_0)$ with
the following properties:
\begin{enumerate}
\item $z'_i - D_i \equiv 0 \mod H(X;\Lambda_+)$.
\item
$z'_i$ satisfies the quantum Stanley-Reisner relation in the
(small) quantum cohomology ring.
\end{enumerate}
\end{thm}
This theorem is proved in  \cite[Section 5.1]{mc-tol}. There our $z'_i$ is
written as $y_i$.
\begin{rem}
We take the quantum cup product without bulk deformation
in Theorem \ref{mcduftolman}.2 above. We can prove a similar result for the quantum cup
product with bulk by using Seidel homomorphism with bulk. We do not prove this here
since we do not use it. 
See \cite{fooospectr}.
\end{rem}
Let $d\ell_j = (v_{j,1},\dots,v_{j,n}) \in \Z^n$ for
$j=1,\dots ,m$. We put
\begin{equation}\label{132}
\frak P^{(0)}_i = \sum_{j=1}^m v_{j,i}Z_j \in \Lambda_0[Z_1,\dots,Z_m],
\qquad i=1,\ldots,n.
\end{equation}
\begin{thm}\label{qhalgebraization}
There exist $\frak P_i \in \Lambda_0[Z_1,\dots,Z_m]$
($i=1,\ldots,n$) such that
\begin{enumerate}\label{Pjzero}
\item
$\frak P_i(z'_1,\dots,z'_m) = 0$ in $H(X;\Lambda_0)$.
Here we use the quantum cup product with $\frak b= \text{\bf 0}$ in the
left hand side.
\item The ring homomorphism
$$
\psi_1: \frac{\Lambda_0[Z_1,\dots,Z_m]}{\text{\rm Clos}_{\frak v_T}
\left(SR_{\omega}\cup\{\frak P_i \mid i=1,\dots,n\}\right)}
\to QH(X;\Lambda_0)
$$
which sends $Z_j$ to $z'_j$ is an isomorphism. 
See Remark \ref{rem:135} for the definition of 
${\text{\rm Clos}}_{\frak v_T}$.
\item
$\frak P_i \equiv \frak P^{(0)}_i \mod \Lambda_+[Z_1,\dots,Z_m]$.
\end{enumerate}
\end{thm}
\begin{rem}\label{rem:135} 
\begin{enumerate}
\item
We recall from Section \ref{sec:statements} that 
we extend the $\frak v_T$ norm on $\Lambda_0$ to 
$\Lambda_0[Z_1, \dots ,Z_m]$ by setting $\frak v_T(Z_i)=0$ 
$(i=1,\dots ,m)$ and denote by 
$\Lambda_0\langle\!\langle Z_1,\dots,Z_m\rangle\!\rangle$
the completion of 
$\Lambda_0[Z_1, \dots ,Z_m]$ with respect to this $\frak v_T$ norm. 
In the theorem above, 
${\text{\rm Clos}}_{\frak v_T}$\index{${\text{\rm Clos}}_{\frak v_T}$} stands for the closure 
with respect to this $\frak v_T$ norm in 
$\Lambda_0\langle\!\langle Z_1,\dots,Z_m\rangle\!\rangle$. 
\item
There is a similar discussion in  \cite[Section 5.1]{mc-tol}. We give the detail
of its proof below for completeness' sake since in \cite{mc-tol} the
process of taking the closure of the {\it ideal} is not explicitly discussed. 
\end{enumerate}
\end{rem}
\begin{proof}
We begin with the following proposition.
\begin{prop}\label{psi1}
There exists an isomorphism
$$
\psi_2: \frac{\Lambda_0\langle\!\langle Z_1,\dots,Z_m\rangle\!\rangle }{\text{\rm Clos}_{\frak v_T}
\left(SR_{\omega}\right)}
\to \Lambda\langle\!\langle y,y^{-1}\rangle\!\rangle_0^P
$$
such that $\psi_2(Z_j) = z_j$.
\end{prop}
\begin{proof}
By Lemma \ref{surjhomfromz} we have a surjective and continuous homomorphism
$$
\tilde{\psi}_{2}: \Lambda_0\langle\!\langle Z_1,\dots,Z_m\rangle\!\rangle
\to \Lambda\langle\!\langle y,y^{-1}\rangle\!\rangle_0^P
$$
such that $\tilde{\psi}_2(Z_j) = z_j$.
Since $z_j$ satisfies the quantum Stanley-Reisner relation,
$\tilde{\psi}_2$ induces a continuous and surjective homomorphism
$\psi_2$. It suffices to show the injectivity of $\psi_2$.
Let $Z^w = Z_1^{w_1}\cdots Z_m^{w_m}$ be a
monomial with $w = (w_1,\dots,w_m) \in \Z_{\ge 0}^m$. We say $Z^w$ 
is \emph{in minimal expression} if $\frak v_T^{P}(z^w) = 0$.
\begin{lem}\label{minimal}
For any  $Z^w$ there exists a
minimal expression $Z^{w'}$ and $\lambda \in \R$ such that
\begin{enumerate}
\item  $Z^w \equiv T^{\lambda}Z^{w'} \mod SR_{\omega}$.
\item  $\lambda \ge 0$ and the equality holds only when $w' = w$.
\end{enumerate}
\end{lem}
\begin{proof}
See the proof of  \cite[Proposition 6.9]{fooo08}.
\end{proof}
Now let $P \in \Lambda_0\langle\!\langle Z_1,\dots,Z_m\rangle\!\rangle $.
Then by the definition of $\Lambda_0\langle\!\langle Z_1,\dots,Z_m\rangle\!\rangle $
which is the completion with respect to the norm $\frak v_T$,
we can write $P$ as the sum $P = \sum_{j=1}^{\infty} a_j Z^{w_{(j)}}$
with $w_{(j)} = (w_{(j),1},\dots,w_{(j),m})$
and $\lim \frak v_T(a_j) = \infty$.
\begin{lem} Suppose $P = \sum_{j=1}^{\infty} a_j Z^{w_{(j)}}$
is in the kernel of $\psi_2$. Choose $w'_{(j)}$ and $\lambda'_j$
as in Lemma \ref{minimal} so that
$
P \equiv \sum_j a_j T^{\lambda'_j} Z^{w'_{(j)}} 
\mod {\text{\rm Clos}}_{\frak v_T} SR_{\omega}
$
and that $Z^{w'_{(j)}}$ are in minimal expression. Then we have
$$
\sum_j a_j T^{\lambda'_j} Z^{w'_{(j)}} = 0.
$$
\end{lem}
\begin{proof} Suppose to the contrary that we
have $\sum_j a_j T^{\lambda'_j} Z^{w'_{(j)}} \neq 0$.
\par
We can rewrite the sum $\sum_j a_j T^{\lambda'_j} Z^{w'_{(j)}}$ as
$
\sum_j b_j T^{\lambda'_j} Z^{w''_{(j)}}
$
so that $w''_{(j)} \ne w''_{(j')}$ if $j\ne j'$. We put
$
c = \frak v_T\left(\sum_j b_j T^{\lambda'_j} Z^{w''_{(j)}}\right).
$
(Note that we use the norm $\frak v_T$ in $\Lambda_0\langle\!\langle Z_1,\dots,Z_m\rangle\!\rangle$ here.)
Take
$
J = \{j \mid \frak v_T(b_j) +  \lambda'_j= c\}.
$
\par
We decompose
$$
b_j = T^{\frak v_T(b_j)} b_{j,0} + b_{j,+}
$$
where $b_{j,0} \in \C$ and $b_{j,+} \equiv 0 \mod T^{\frak v_T(b_j)}\Lambda_+$.
Then
$$
\sum_j b_j T^{\lambda'_j} Z^{w''_{(j)}}
\equiv
\sum_{j\in J} b_{j,0} T^{c} Z^{w''_{(j)}} \mod T^{c}\Lambda_+.
$$
Lemma \ref{minimal}.2 implies that the right hand side 
is nonzero.
Therefore
$J$ is nonempty.
\par
For $j \in J$, we write
$
Z^{w''_{(j)}} = \prod_{i=1}^{m} Z_i^{k_{i,j}}.
$
We put $I(j) = \{i_{1,j},\dots,i_{l(j),j}\} = \{i\mid k_{i,j} \ne 0\}$.
By definition of the minimal expression, we find that
$
P_{I(j)} = \bigcap_{a=1}^{l(j)} P_{i_{a,j}}
$
is nonempty. Let $u_j$ be an interior point of
$P_{I(j)}$. We recall $\frak v_T^{u_j}(z^{w''_{(j)}}) = 0$ since
$u_j \in P_{I(j)} \subset \partial P$
and $\frak v_T^{u}(z^{w''_{(j)}})>0$ for an $u \notin P_{I(j)}$.
\par
We choose $j_{\rm min}$ so that $P_{I(j_{\rm min})}$ is minimal in that
$P_{I(j_{\rm min})} \supset P_{I(j')}$ does not hold for any 
$j'$ with $I(j')\ne I(j_{\rm min})$.
Then, we have
$\frak v_T^{u_{j'}}(z^{w''_{(j')}}) > 0$ for any $j'$ with $I(j') \ne I(j_{\rm min})$.
In fact there is $i \in \{ i_{1,j'}, \dots, i_{\ell(j'),j'}\}$ such that $k_{i,j'} > 0$ and $k_{i,j_{\rm min}} = 0$.
Since $\frak v_T^{u_{j_{\rm min}}}(z_i) > 0$ for such $i$, it follows that
$\frak v_T^{u_{j_{\rm min}}} (z^{w''_{(j')}}) > 0$.
\par
Therefore 
$i_{u_{j_{\rm min}}} \psi_2(\sum_j b_j T^{\lambda'_j} Z^{w''_{(j)}}) \equiv 
\sum_{j;I(j)=I(j_{\rm min})} b_jT^{\lambda'_j}z^{w''_{(j)}}$
modulo $T^c\Lambda\langle\!\langle y,y^{-1}\rangle\!\rangle_0^{u_{j_{\rm min}}}$. Here
$$
i_{u_{j_{\rm min}}}: \Lambda\langle\!\langle y,y^{-1}\rangle\!\rangle_0^{P}
\to \Lambda\langle\!\langle y,y^{-1}\rangle\!\rangle_0^{u_{j_{\rm min}}} 
$$
is the obvious inclusion.
This contradicts to the hypothesis $P \in \operatorname{Ker} \psi_2$.
\end{proof}
The proof of Proposition \ref{psi1} is now complete.
\end{proof}
Construction of $\frak P_i$ required in Theorem \ref{Pjzero}
will be finished by the following lemma. 
\begin{lem}\label{Pjkinductconstr}
Let $G \supset G_0$ such that $z'_j$ are $G$ gapped.
Then
there exist $\frak P_i^{(k)} \in \Lambda_0[Z_1,\dots,Z_m]$ with the following
properties:
\begin{enumerate}
\item $\deg \frak P_i^{(k)} \le \dim X/2$ and $\frak P_i^{(k)}$ is $G$-gapped.
\item $\frak P_i^{(k+1)}  \equiv \frak P_i^{(k)} \mod T^{\lambda_{k+1}}\Lambda_0[Z_1,\dots,Z_m]$.
\item $\frak P_i^{(0)}$ is as in $(\ref{132})$.
\item
$
\frak P_i^{(k)}(z'_1,\dots,z'_m) \equiv 0 \mod T^{\lambda_{k+1}}H(X;\Lambda_0).
$
\end{enumerate}
\end{lem}
\begin{proof}
Let $\frak P_i^{(0)}$ be as in $(\ref{132})$. Theorem \ref{mcduftolman} implies
$$
\frak P_i^{(0)}(z'_1,\dots,z'_m) =
\sum_{j} v_{j,i} D_j\equiv 0 \mod  T^{\lambda_1}H(X;\Lambda_0).
$$
Hence we can start the induction.
We define $c_i \in H(X;\C)$ by the formula
$$
\frak P_i^{(k)}(z'_1,\dots,z'_m) \equiv T^{\lambda_{k+1}}c_{i} \mod T^{\lambda_{k+2}}H(X;\Lambda_0).
$$
Since $D_i$ generate $H(X;\C)$ as a ring,
we have $P_{i,k} \in \C[Z_1,\dots,Z_m]$ such that
$$
c_i = P_{i,k}(D_1,\dots,D_m).
$$
We put $\frak P_i^{(k+1)} = \frak P_i^{(k)} - T^{\lambda_{k+1}}P_{i,k}$.
It is easy to check that this has the required properties.
\end{proof}
We put $\frak P_i = \lim_{k\to\infty} \frak P_i^{(k)}$.
Properties 1 and 3 are satisfied by Lemma \ref{Pjkinductconstr}.
The proof that $\psi_1$ is an isomorphism is the same as the proof of
Lemma \ref{polyapprox}.2 and 3. (We remark that $H(X;\Lambda_0)$ is a
free $\Lambda_0$ module. Therefore the assumption of Lemma \ref{polyapprox}.3
is satisfied in our case.)
Thus the next lemma completes the proof of Theorem \ref{qhalgebraization}.
\end{proof}
\begin{lem}\label{13saigo}
We have an isomorphism 
$$
\frac{\Lambda_0[Z_1,\dots,Z_m]}{\text{\rm Clos}_{\frak v_T}
\left(SR_{\omega}\cup\{\frak P_i \mid i=1,\dots,n\}\right)}
\cong \frac{\Lambda_0\langle\!\langle Z_1,\dots,Z_m\rangle\!\rangle }{\text{\rm Clos}_{\frak v_T}
\left(SR_{\omega}\cup\{\frak P_i \mid i=1,\dots,n\}\right)}.
$$
\end{lem}
\begin{proof}
We note that
$\psi_1(Z_i) = z'_i$.
Since $z'_i \equiv \text{\rm PD}([D_i]) \mod T^{\lambda_1}$
and quantum cohomology is the same as singular cohomology ring
modulo $ T^{\lambda_1}$, we have
\begin{equation}\label{Zhigherpowervanish}
\psi_1(Z_1^{k_1}\cdots Z_m^{k_m})
\equiv 0
\mod T^{\lambda_1(k_1+\cdots+k_m-n)/n}H(X;\Lambda_0).
\end{equation}
Using (\ref{Zhigherpowervanish}) in place of Sublemma \ref{sublemmamendou}
we can repeat the proof of Lemma \ref{polyapprox}.2
to prove Lemma \ref{13saigo}.
\end{proof}
\par
\section{Injectivity of Kodaira-Spencer map}
\label{sec:injgen}
\subsection{The case $\frak b \in \mathcal A(\Lambda_+)$}
\label{subsec:caseofLambdazero0}
In this subsection, we prove the injectivity for the case
$\frak b \in \mathcal A(\Lambda_+)$. We discuss the case
$\frak b \in \mathcal A(\Lambda_0)$ in the next subsection.
\begin{thm}\label{injectivity}
If $\frak b \in \mathcal A(\Lambda_+)$, then
$$
\frak{ks}_{\frak b}: H(X;\Lambda_0) \to \text{\rm Jac}(\mathfrak{PO}_{\frak b})
$$
is injective.
\end{thm}
The main part of the proof is the following proposition.
We put
$$\aligned
\text{\rm Jac}_{\frak v^P_T}(\mathfrak{PO}_{\frak b})_{\Lambda} &=
\text{\rm Jac}_{\frak v^P_T}(\mathfrak{PO}_{\frak b}) \otimes_{\Lambda_0} \Lambda,
\\
\text{\rm Jac}(\mathfrak{PO}_{\frak b})_{\Lambda} &=
\text{\rm Jac}(\mathfrak{PO}_{\frak b}) \otimes_{\Lambda_0} \Lambda.
\endaligned$$
\begin{prop}\label{bettiest}
$$
\dim_{\C} H(X;\C) \le \dim_{\Lambda}\text{\rm Jac}_{\frak v^P_T}(\mathfrak{PO}_{\frak b})_{\Lambda}.
$$
\end{prop}
\begin{proof}[Proposition \ref{bettiest} $\Rightarrow$ Theorem \ref{injectivity}]
Let $\overline{\text{\bf f}}_i$ ($i = 1,\dots, B'$) be the basis of $H(X;\Lambda_0)$.
By Theorem \ref{surj}, $\frak{ks}_{\frak b}(\overline{\text{\bf f}}_i)$
generates $\text{\rm Jac}_{\frak v^P_T}(\mathfrak{PO}_{\frak b})$.
Therefore by Proposition \ref{bettiest}
they are linearly independent in $\text{\rm Jac}_{\frak v^P_T}(\mathfrak{PO}_{\frak b})_{\Lambda}$.
Thus
$$
\frak{ks}^P_{\frak b}: H(X;\Lambda_0) \to \text{\rm Jac}_{\frak v^P_T}(\mathfrak{PO}_{\frak b})
$$
is injective. We can also prove that $\frak{ks}^P_{\frak b}$
is surjective in the same way as
Theorem \ref{surj}.
It follows that $\text{\rm Jac}_{\frak v^P_T}(\mathfrak{PO}_{\frak b})$ is torsion free. Therefore
$\text{\rm Jac}_{\frak v^P_T}(\mathfrak{PO}_{\frak b}) \cong \text{\rm Jac}(\mathfrak{PO}_{\frak b})$ by
Lemma \ref{polyapprox}.4. Hence Theorem \ref{injectivity}.
\end{proof}
The proof of Proposition \ref{bettiest} occupies the rest of this subsection.
Let $y'$ be as in Theorem \ref{algebraization}.2.
We put
$$
\mathcal P = \mathfrak{PO}_{\frak b}(y')  \in \Lambda[y,y^{-1}]_0^P.
$$
Since $y'$ is strict and $\frak b \in \mathcal A(\Lambda_+)$,
we have
$$
\mathcal P \equiv z_1 + \dots + z_m \mod \Lambda[y,y^{-1}]_+^P.
$$
Let $\frak P_i \in \Lambda_0[Z_1,\dots,Z_m]$ be as in Theorem \ref{qhalgebraization}.
We have
$
\frak P_i(z) \in \Lambda[y,y^{-1}]_0^P.
$
\begin{defn}
\begin{equation}
\frak B_i= s\frak P_i(z) + (1-s) y_i \frac{\partial \mathcal P}{\partial y_i}
\in \Lambda[s,y,y^{-1}]_0^P.
\end{equation}
\begin{equation}
\mathcal R = \frac{\Lambda[s,y,y^{-1}]}{(\frak B_i:i=1,\dots,n)},
\qquad
\underline{\mathcal R} = \frac{\Lambda[s,y,y^{-1}]}{\text{\rm Clos}_{\frak v_T^P}(\frak B_i:i=1,\dots,n)}.
\end{equation}
\end{defn}
There are morphisms of schemes over $\Lambda$:
$$
Spec(\underline{\mathcal R}) \longrightarrow Spec({\mathcal R}) \overset{\pi}\longrightarrow Spec(\Lambda[s])
= \Lambda.
$$
We denote $I = (\frak B_1,\dots,\frak B_n )$
and by $\overline{I}$ the closure of $I$ with respect to ${\frak v}_T^P$.
We take irredundant primary decompositions of the ideals $I$ and $\overline{I}$:
\begin{equation}
I=\bigcap_{i \in \frak I_+} {\frak q_i}, \qquad \overline{I}=\bigcap_{j \in  \frak J} {\frak Q_j}.   \label{primary}
\end{equation}
Here an irredundant primary decomposition represents the given ideal
as the intersection of primary ideals such that any $\frak q_i$ (resp. $\frak Q_j$) cannot be
eliminated and the set of associated primes $\{ \sqrt{q_i} \}$ (resp. $\{ \sqrt{\frak Q_j} \}$) consists
of mutually distinct prime ideals.
Note that, for an irredundant primary decomposition,  primary ideals belonging
to  minimal associated primes are uniquely determined.
We decompose $Spec({\mathcal R})$
(resp. $Spec(\underline{\mathcal R}$) ) into irreducible components.
Namely
$$
\vert Spec({\mathcal R})\vert = \bigcup_{i \in \frak I_+}
\left\vert Spec\left(\frac{\Lambda[s,y,y^{-1}]}{\frak q_i}\right)\right\vert,
\qquad
\vert Spec(\underline{\mathcal R})\vert = \bigcup_{j \in  \frak J}
\left\vert Spec\left(\frac{\Lambda[s,y,y^{-1}]} {\frak Q_j}\right)\right\vert.
$$
Here and hereafter in this section, we put
$\vert Spec A\vert$ for the set of $\Lambda$-valued points of the $\Lambda$-scheme $Spec A$.
We put
\begin{equation}\label{defJ}
\frak I = \left\{ i \in \frak I_+ ~\left\vert~
\left\vert Spec\left(\frac{\Lambda[s,y,y^{-1}]}{\frak q_i}\right) \right\vert
\cap \left\vert Spec(\underline{\mathcal R}) \right\vert \ne \emptyset \right\} \right.
\end{equation}
and
$$
\frak X_i = Spec\left(\frac{\Lambda[s,y,y^{-1}]}{\frak q_i}\right),
\qquad
\frak Y_j =Spec\left(\frac{\Lambda[s,y,y^{-1}]} {\frak Q_j}\right).
$$
Let $(\frak X_i)_{\rm red}$ (resp. $(\frak Y_j)_{\rm red}$)
be the associated reduced schemes to  $\frak X_i$ (resp. $\frak Y_j$).
Namely
$$
(\frak X_i)_{\rm red} = Spec\left(\frac{\Lambda[s,y,y^{-1}]}{\sqrt{\frak q_i}}\right),
\qquad
(\frak Y_j)_{\rm red} =Spec\left(\frac{\Lambda[s,y,y^{-1}]} {\sqrt{\frak Q_j}}\right).
$$
We also put
\begin{equation}\label{I0decompose}
I_0 = \bigcap_{i\in \frak I} \frak q_i
\end{equation}
and
$$
\frak X = Spec\left( \frac{\Lambda[s,y,y^{-1}]}{I_0}\right),  \qquad
\frak Y = Spec\left( \frac{\Lambda[s,y,y^{-1}]}{\bigcap_{j \in \frak J} \frak Q_j}\right).
$$
Then we have
$$
\vert \frak X \vert = \bigcup_{i\in \frak I}
\left\vert \frak X_i \right\vert,
\qquad
\vert \frak Y \vert = \bigcup_{j\in \frak J}
\left\vert \frak Y_j \right\vert.
$$
\begin{lem}\label{underRcomplete}
We have
$$
\underline{\mathcal R} \cong \frac{\Lambda\langle\!\langle y,y^{-1}\rangle\!\rangle^P [s]}{\text{\rm Clos}_{\frak v_T^P}(\frak B_i:i=1,\dots,n)}.
$$
\end{lem}
\begin{proof}
We use the following:
\begin{sublem}\label{SRlinearisenough}
Let $z_i$ be as in $(\ref{zjdef})$. Then we have
$$
z_1^{k_1}\cdots z_m^{k_m} \equiv 0 \mod T^{\lambda_1(k_1+\cdots+k_m-n)/n}\frac{\Lambda\langle\!\langle y,y^{-1}\rangle\!\rangle_0^P [s]}{\text{\rm Clos}_{\frak v_T^P}(\frak B_i:i=1,\dots,n)}.
$$
\end{sublem}
\begin{proof}
We note that $z_i$ satisfies the Stanley-Reisner relation and linear relation
(that is
$
\sum_{j=1}^m v_{ij}z_j =0
$,
$i=1,\ldots,n$, see (\ref{132}))
modulo $T^{\lambda_1}H(X;\Lambda_0)$.
The cohomology ring $H(X;\Q)$ is a quotient of the polynomial ring over $z_i$ by the
Stanley-Reisner relation and linear relation.
 (See \cite[p106 Proposition]{fulton}, for example.) On $H(X;\Q)$ the product of $(n+1)$ of $z_i$'s are zero.
 Therefore we have
$$
z_1^{k_1}\cdots z_m^{k_m} \equiv 0 \mod T^{\lambda_1(k_1+\dots+k_m-n)/n}H(X;\Lambda_0)
$$
if $k_1+\cdots+k_m > n$. The sublemma follows easily.
\end{proof}
Once we have Sublemma \ref{SRlinearisenough},
the rest of the proof is the same as the proof of Lemma \ref{13saigo}.
\end{proof}
\begin{lem}\label{inmopoly}
If $(\frak y_1,\dots,\frak y_n,\frak s) \in \vert Spec(\underline{\mathcal R})\vert$,
then
\begin{equation}\label{inmompoly}
(\frak v_T(\frak y_1),\dots,\frak v_T(\frak y_n)) \in P.
\end{equation}
If $(\frak y_1,\dots,\frak y_n,\frak s)\in \vert Spec({\mathcal R})\vert$
and $(\ref{inmompoly})$ is satisfied, then
$(\frak y_1,\dots,\frak y_n,\frak s)\in \vert Spec(\underline{\mathcal R})\vert$.
\end{lem}
\begin{proof}
The first half is a consequence of Lemmata \ref{pointP} (B) and \ref{underRcomplete}.
The second half follows from the fact that
$y_i \mapsto \frak y_i$ induces a continuous map from
$(\Lambda[y,y^{-1}],\frak v_T^P)$ to $\Lambda$ if $(\frak v_T(\frak y_1),\dots,\frak v_T(\frak y_n)) \in P$.
\end{proof}
\begin{cor}\label{I0Icoro}
We have
$$
\overline I_0 = \overline I,
$$
where $\overline I_0$ is the closure of $I_0$.
\end{cor}
\begin{proof}
Let $i \notin \frak I$. We will prove $\overline{\frak q}_i
= \Lambda[s,y,y^{-1}]$. 
In fact, since $i \notin \frak I$, we have
$$
\left\vert Spec \left( \frac{\Lambda[s,y,y^{-1}]}{\overline{\frak q}_i}\right) \right\vert
\cap \vert Spec(\underline{\mathcal R})\vert
\subseteq
\left\vert Spec \left( \frac{\Lambda[s,y,y^{-1}]}{\frak q_i}\right) \right\vert
\cap \vert Spec(\underline{\mathcal R})\vert =\emptyset.
$$
\par
On the other hand, since  $\overline{\frak q}_i \supseteq \overline I$, it follows from
Lemma \ref{inmopoly} that
$$
\left\vert Spec \left( \frac{\Lambda[s,y,y^{-1}]}{\overline{\frak q}_i}\right) \right\vert
\subseteq \vert Spec(\underline{\mathcal R})\vert.
$$
Therefore
$$
\left\vert Spec \left( \frac{\Lambda[s,y,y^{-1}]}{\overline{\frak q}_i}\right) \right\vert  = \emptyset.
$$
Namely $\overline{\frak q}_i
= \Lambda[s,y,y^{-1}]$ as claimed. 
Therefore $1 \in \overline{\frak q}_i $ for $i \notin \frak I$.
The corollary follows easily.
\end{proof}
\begin{prop}\label{componentmajiwari}
For each $i \in \frak I$ there
exists $j\in \frak J$ such that $\vert \frak X_i \vert \subseteq \vert\frak Y_j  \vert$.  \end{prop}
\begin{proof}
Let $\frak Z$ be one of the
irreducible components of $\frak X_i \cap Spec(\underline{\mathcal R})$.
We will prove the coincidence of geometric points in $\frak X_i$ and $\frak Z$, 
by contradiction. 
The proposition then will follow.
\begin{lem}\label{formalization} Let $\frak Z$, $\frak X_i$ be as above,
and $\vert\frak Z  \vert \ne \vert \frak X_i \vert$.
Then, there exists a commutative diagram of $\Lambda$-schemes:
\begin{equation}\label{nonliftdiagram1}
\begin{CD}
(\frak X_i)_{\rm red} @<<< \frak Z_{\rm red} \\
@ AA{\frak f}A @ AA{\overline{\frak f}}A\\
Spec(B) @<<<  Spec(\Lambda)
\end{CD}
\end{equation}
with the following properties.
Let $\tilde{\frak y} \in \vert Spec(B) \vert$ be the geometric point corresponding to the
lower horizontal arrow.
\begin{enumerate}
\item $Spec(B)$ is smooth at $\tilde{\frak y}$.
\item $\frak f$ is a composition of an embedding and normalization.
\item  $\vert Spec(B)\vert \cap \vert Spec(\underline{\mathcal R})\vert$ is zero dimensional.
In particular, $\frak f$ cannot be lifted to a morphism $Spec(B) \to
\frak Z_{\rm red}$.
\end{enumerate}
\end{lem}
\begin{proof}
We may take a 1-dimensional irreducible and reduced scheme
$\frak C =Spec(A) \subset (\frak
X_i)_{\rm red}$, which is
a subscheme of $Spec(\mathcal R)$
with $\dim (\frak C \cap \frak Z_{\rm red}) = 0$.
Let $\frak y$ be the intersection point.
Then we normalize $\frak C$ at $\frak y$.
We pick a point $\tilde{\frak y}$ in the inverse image
of $\frak y$ under the normalization map.
Its appropriate open subscheme is the required
$Spec(B)$.
\end{proof}
Let $\frak m_{\tilde{\frak y}}$ be a maximal ideal of $B$ corresponding to $\tilde{\frak y}$.
The localization $B_{\frak m_{\tilde{\frak y}}}$ is a one dimensional regular local ring.
Let $S \in B$ be a generator of the ideal
$\frak m_{\tilde{\frak y}}$ in
$B_{\frak m_{\tilde{\frak y}}}$.
(Note there is no relation between $S$ and $s$.)
The inclusion $\Lambda[S] \to B$ can be regarded as a map
$Spec(B) \to Spec(\Lambda[S])$. Its differential at
$\tilde{\frak y}$ is nonzero.
We next want to apply `inverse mapping theorem' for this map at $\tilde{\frak y}$.
Of course, in algebraic geometry we do not have one.
In other words,
the inverse does not exist on a Zariski open subset.
However, if we consider the infinitesimal neighborhood of $S=0$ in
$Spec(\Lambda[S])$, we have an inverse.
The detail will be as follows.
\par
We regard $y_i, s \in B$ as $\Lambda$-valued functions on
$\vert Spec (B)\vert$.
We write $s = y_0$ and
define $\frak y_i = y_i(\tilde{\frak y})$
($i=0,\dots ,n$).
Here $\frak y_i\ne 0$ for $i\ne 0$.
\par
Now we take generators $z_1, \dots , z_h$ of $B$
as a $\Lambda$-algebra
such that $z_i (\tilde{\frak y})=0$ ($i=1,\dots ,h$).
Since $z_i \in \frak m_{\tilde{\frak y}} \subset B_{\frak m_{\tilde{\frak y}}}$, we have
\begin{equation}\label{21118}
z_i = S^k \frac{P_i(z_1,\dots,z_h)}{Q_i(z_1,\dots,z_h)}
\end{equation}
with $P_i(0,\dots,0) \ne 0$, $Q_i(0,\dots,0) \ne 0$, $k\in \Z_{>0}$.
We put
\begin{equation}\label{defFi}
\aligned
F_i(Z_1,\dots,Z_h,S)
&= Z_i - S^k\frac{P_i(Z_1,\dots,Z_h)}{Q_i(Z_1,\dots,Z_h)} \\
&= Z_i - S^kU_i(Z_1,\dots,Z_h).
\endaligned
\end{equation}
\begin{rem}
Note $P_i, Q_i \in \Lambda[[Z_1,\dots,Z_h]]$ 
are formal power series. The equality (\ref{defFi}) 
is one in the ring $\Lambda[[Z_1,\dots,Z_h,S]]$.
We substitute $Z_i$ by $z_i$ to obtain an element of $B$. 
So (\ref{21118}) is an equality in the ring $B_{{\frak m}_{\tilde{\frak y}}}$.
\end{rem}
In the next lemma we use the norm $\frak v_T$ on 
$\Lambda[[Z_1,\dots,Z_h]]$ such that $\frak v_T(Z_i) = 0$.
\begin{lem}
Under the situation above,
we have $K, C>0$ such that if
\begin{equation}\label{adiccond}
\frak v_T(a_i) \ge K \quad (i=1,\dots ,h), \quad a_i \in \Lambda_0,
\end{equation}
then $Q_i(a_1,\dots,a_h) \ne 0$ and
the following inequalities hold
\begin{equation}\label{adicconcl1}
\frak v_T\left( \frac{\partial U_i}{\partial Z_j} (a_1,\dots,a_h)\right) > -C,
\end{equation}
\begin{equation}\label{adicconcl2}
\frak v_T\left( \frac{\partial^2 U_i}{\partial Z_j\partial Z_k} (a_1,\dots,a_h)\right) > -C.
\end{equation}
\end{lem}
The proof is easy and is left to the reader.
We note that if both of $(a_1,\dots ,a_h)$ and
$(a'_1,\dots ,a'_h)$ satisfy
(\ref{adiccond}), then (\ref{adicconcl2}) yields
\begin{equation}\label{adicconcl3}
\aligned
&\frak v_T\left( \frac{\partial U_i}{\partial Z_j} (a_1,\dots,a_h)
-
\frac{\partial U_i}{\partial Z_j} (a'_1,\dots,a'_h) \right)\\
&> \inf_{i=1,\dots,h}\frak v_T(a_i-a'_i)- C.
\endaligned
\end{equation}
\par
We choose $K>2C$.
Now we consider $\frak S$ with $\frak v_T(\frak S) \ge K$.
Then we have
\begin{equation}\label{adicconcl4}
\frac{\partial F_i}{\partial Z_j} (0,\dots,0,\frak S)
\equiv
\begin{cases}
1 \mod T^{K-C}  &\text{if $i=j$}\\
0 \mod T^{K-C}  &\text{if $i\ne j$}.
\end{cases}
\end{equation}
This is a consequence of (\ref{adicconcl1}).
We also have
$$
F_i(0,\dots,0,\frak S) \equiv 0 \mod T^{K-C}.
$$
Now we iteratively define
$a_i^{\ell} \in \Lambda_0 ~(i=1,\dots ,h)$
by induction on $\ell =0,1,\dots$ so that
$$
\left\{
\aligned
(\Delta_{\ell+1,i})_{i=1}^h
&= -\left(\left[
\frac{\partial F_i}{\partial Z_j}(a^{\ell}_1,\dots,a^{\ell}_h,\frak S)
\right]_{i=1,~j=1}^{i=h,~j=h}\right)^{-1}
(F_{i}(a^{\ell}_1,\dots,a^{\ell}_h,\frak S))_{i=1}^h \\
a_i^{\ell} &= \Delta_{1,i} + \dots +\Delta_{\ell,i},
\quad a_i^0 = 0.
\endaligned
\right.
$$
We can prove by induction on $\ell$ that
$$
\aligned
F_i(a^{\ell}_1,\dots,a^{\ell}_h,\frak S)
& \equiv 0 \mod T^{(\ell+1)(K-C)} \\
\Delta_{\ell+1,i} & \equiv 0 \mod T^{(\ell+1)(K-C)}.
\endaligned
$$
\par
We define
\begin{equation}\label{defai}
a_i(\frak S) = \lim_{\ell\to \infty}  a^{\ell}_i.
\end{equation}
Then they satisfy the equations
$$
F_i(a_1(\frak S),\dots,a_h(\frak S),\frak S) = 0, \quad i=1,\dots ,h.
$$
On the other hand, we can write
\begin{equation}\label{defRi}
y_i = R_i (z_1, \dots ,z_h), \quad i=0,\dots ,n
\end{equation}
for some $R_i \in \Lambda[Z_1,\dots ,Z_h]$. 
We define
$$
 y_i(\frak S) = R_i(a_1(\frak S),\dots,
a_h(\frak S)), \quad 
y(\frak S) = 
(y_1(\frak S), \dots , y_n(\frak S)).
$$
\begin{lem} \label{inBlem}
For any $\frak S$ with
sufficiently large $\frak v_T(\frak S)$,
$y(\frak S)$
is a geometric point of $Spec(B)$.
\end{lem}
We will prove it later.
\par
We consider $\frak S$ with sufficiently large
$\frak v_T(\frak S)$ as in Lemma \ref{inBlem}.
Note $\frak v_T(y_i(\frak S)) = \frak v_T(y_i(0)) = \frak v_T(\frak y_i)$ holds
if $\frak v_T(\frak S)$ is sufficiently large.
\par
On the other hand, since $\frak y$ is in $\vert Spec(\underline{\mathcal R})\vert$,
Lemma \ref{inmompoly} implies
$$
({\frak v}_T(\frak y_1),\dots,{\frak v}_T(\frak y_n)) \in P.
$$
Therefore $y(\frak S) $ is also a geometric
point of $\vert Spec(\underline{\mathcal R})\vert$.
Moreover, from $U_i(0,\dots,0) \ne 0$, we have $a_i(\frak S) \ne 0$ for at least one $i$.
Since
$Spec (B)$ is a normalization of $\frak C =Spec (A)$, the ring
$B$ is a finite $A$-module and so the mapping
$$
(a_1,\dots a_h) \mapsto (R_0(a_1, \dots ,a_h), \dots , R_n
(a_1,\dots ,a_h))
$$
is finite to one.
Thus not only $a(\frak S)=(a_1(\frak S), \dots , a_h(\frak S))$ but also $y(\frak S)$ are
non-constant.
\par
Hence $Spec(B)$ has infinitely many geometric points in
$\vert Spec(\underline{\mathcal R})\vert$.
This is a contradiction.
It implies that the assumption 
$\vert \frak Z \vert \ne \vert \frak X_i\vert$ in Lemma \ref{formalization}  is not correct.
Namely $\vert \frak Z \vert = \vert \frak X_i\vert$.
The proof of Proposition \ref{componentmajiwari} is completed modulo Lemma \ref{inBlem}.
\end{proof}
\begin{proof}[Proof of Lemma \ref{inBlem}]
We may identify $B$ with 
$$
\frac{\Lambda[Z_1,\dots,Z_h,S]}{\frak p}
$$
where $Z_i$ corresponds to $z_i$ in $B$.
We denote by $B_{\frak m_{\tilde{\frak y}}}$
the localization of $B$ by the maximal ideal
$\frak m_{\tilde{\frak y}}$ corresponding to the point
$\tilde{\frak y}$. 
We next take the completion
$B_{\frak m_{\tilde{\frak y}}}^{\wedge}$
of $B_{\frak m_{\tilde{\frak y}}}$
by the ideal
$\frak m_{\tilde{\frak y}}$.
It is a complete discrete valuation ring over the field
$\Lambda$
(with respect to the $\frak m_{\tilde{\frak y}}$ topology),
hence isomorphic to the formal power series ring $\Lambda [[t]]$.
(See for example, \cite[Proposition 10.16]{Eisen}, \cite[Theorem 37 in (17.H)]{Ma70}.)
(Here we used the fact that $\Lambda$ is algebraically closed in order to
show
$\Lambda \cong B/{\frak m_{\tilde{\frak y}}}
\cong
B_{{\frak m}_{\tilde{\frak y}}}^{\wedge} /
\frak m_{\tilde{\frak y}}
B_{{\frak m}_{\tilde{\frak y}}}^{\wedge}
$.
(See \cite[Lemma A.1]{fooo08}.))
Note that we may take $S \in \Lambda[[t]]$ as the generator
of $\frak m_{\tilde{\frak y}}$ we took before.
Thus we may take $t=S$.
\par
We have the following diagram.
\begin{equation}
\begin{CD}
B \cong
\Lambda[Z_1,\dots,Z_h,S]/\frak p @= \Lambda[Z_1,\dots,Z_h,S]/\frak p \\
@VVV @VVV \\
B_{{\frak m}_{\tilde{\frak y}}}  @>>>
B/\frak m_{\tilde{\frak y}}
\\
@VVV  @VV{\widetilde{\frak f}}V \\
B_{{\frak m}_{\tilde{\frak y}}}^{\wedge} \cong \Lambda[[S]] @>>>
B_{{\frak m}_{\tilde{\frak y}}}^{\wedge} /
\frak m_{\tilde{\frak y}} B_{{\frak m}_{\tilde{\frak y}}}^{\wedge}
 \cong \Lambda
\end{CD}
\end{equation}
\par
The element $z_i \in B$  induces
 $[z_i] \in
B_{{\frak m}_{\tilde{\frak y}}}^{\wedge}
\cong \Lambda[[S]]$.
We denote this element by
$\frak f_i(S)$.
For any $W \in \frak p \subset \Lambda[Z_1,\dots,Z_h,S]$
we have
\begin{equation}\label{mitashu}
W(\frak f_1(S),\dots,\frak f_h(S),S)=0
\in \Lambda[[S]].
\end{equation}
We remark that
$$
Q_i(Z_1,\dots,Z_h)F_i(Z_1,\dots,Z_h,S) \in \frak p.
$$
Therefore by (\ref{mitashu}) we have
$$
Q_i(\frak f_1(T^{K}S),\dots,\frak f_h(T^{K}S))
F_i(\frak f_1(T^{K}S),\dots,\frak f_h(T^{K}S), T^K S)=0.
$$
Since
$Q_i(\frak f_1(T^{K}S),\dots,\frak f_h(T^{K}S))\ne 0$ for sufficiently large $K$, we obtain
$$
F_i(\frak f_1(T^{K}S),\dots,\frak f_h(T^{K}S), T^K S)=0.
$$
In other words, by (\ref{defFi}) we have
\begin{equation}\label{inductiveU}
\aligned
\frak f_i(T^{K}S) &= T^{kK}S^k\frac{P_i(\frak f_1(T^{K}S),\dots,\frak f_h(T^{K}S))}
{Q_i(\frak f_1(T^{K}S),\dots,\frak f_h(T^{K}S))}\\
&= T^{kK}S^kU_i(\frak f_1(T^{K}S),\dots,\frak f_h(T^{K}S)).
\endaligned
\end{equation}
We observe that $Q_i$ is nonzero at $Z_1=\dots = Z_h =0$, therefore
we may regard
$$
U_i \in \Lambda[[Z_1,\dots,Z_h]].
$$
Note $\frak f_i(S) \in S\Lambda[[S]]$.
Therefore the substitution
$$
U_i(\frak f_1(T^{K}S),\dots,\frak f_h(T^{K}S))
$$
makes sense using $S$-adic topology. (Namely as a formal power series of $S$.)
Now we observe that
the equation
\begin{equation}\label{againUandS}
\frak f_i(T^{K}S) = T^{kK}S^kU_i(\frak f_1(T^{K}S),\dots,\frak f_h(T^{K}S))
\end{equation}
with $\frak f_i(0) = 0$, {\it uniquely} characterizes
$\frak f_i(S) \in \Lambda[[S]]$.
\par
In fact,
we define $X_{i,m}$ $(i=1,\dots ,h)$ inductively by
$X_{i,0} = 0$ and
\begin{equation}\label{Xinduction}
X_{i,m+1} = T^{kK}S^kU_i(X_{1,m},\dots,X_{h,m}).
\end{equation}
Then $\{ X_{i,m}\}_{m=0,1,\dots}$
converges to $\frak f_i(T^KS)$
in $S$-adic topology and is necessarily the unique solution of
(\ref{againUandS}).
\par
Now we consider the same inductive construction of the
solution of (\ref{againUandS}), but this time with $T$-adic topology.
Namely we use $\frak v_T$ with $\frak v_T(S) = 0$.
Then for a sufficiently large $K$,
we have a solution of (\ref{Xinduction})
in the way we described during the iterative argument of the proof of Proposition \ref{componentmajiwari}.
It gives
$a_i(T^KS)$.
Namely,
if we substitute $S$ with element $c$ of $\Lambda_0$,
we obtain
$a_i(cT^K)$ which we defined before in
(\ref{defai}) with $\frak S=cT^K$ satisfying
$\frak v_T(\frak S) \ge K$.
\par
Now by the uniqueness of the solution of the recursive equation
(\ref{Xinduction}) we have
$$
a_i(T^KS) = \frak f_i(T^{K}S).
$$
In particular,
$$
W(a_1(T^KS),\dots,a_h(T^KS),T^kS)=0
$$
for all $W\in \frak p$. This proves the lemma.
\end{proof}
\begin{rem}
Suppose that $\frak f_i(S) \in T^{-C} \Lambda_0[[S]]$ for some $C$.
(We do not know whether this is true or not.)
In that case we can discuss as follows.
We put
$\frak g_i(S) = R_i (\frak f_1 (S),\dots , \frak f_h(S))$,
where $R_i \in \Lambda[Z_1, \dots ,Z_h]$ as in (\ref{defRi}).
We then can write
$$
\frak g_i(S) = \sum_{j=0}^{\infty} T^{\lambda_{i,j}}\frak g_{i,j}(S).
$$
Here $\frak g_{i,j}(S) \in \C[[S]] \setminus\{0\}$ and $\lambda_{i,j} \in \R$.
We may assume $\lambda_{i,j} < \lambda_{i,j+1}$.
\par
We put
$$
j_{i}^0 = \inf\{ j \mid \frak g_{i,j}(0) \ne 0\}.
$$
It is easy to see that $\frak y_i = \frak g_i(0)$.
Therefore we have
\begin{equation}\label{inIntP}
{\frak v}_T(\frak y)=(\lambda_{1,j_1^0},\dots,\lambda_{n,j_n^0}) = u \in P.
\end{equation}
Now we write
$$
{\frak g}_{i,j}(S)=d_{i,j}S^{\rho_{i,j}} + {\text{higher order terms}}
$$
with $d_{i,j}\in \C\setminus 0$ and $\rho_{i,j}>0$. For $\kappa > 0$ we
define
$$
\lambda_{i}(\kappa) = \inf\{ \lambda_{i,j} + \rho_{i,j}\kappa \mid j=0,\dots, j_i^0\}.
$$
We take a complex number $c \in \C$ which satisfies
$$
\sum_{j =0,\dots, j_i^0; \atop \lambda_{i,j}+\rho_{i,j}\kappa=\lambda_i(\kappa)} d_{i,j}c^{\rho_{i,j}} \ne 0.
$$
We note that we can choose $c$ independent of $\kappa$,
because it is enough to take $c$ so that
$
\sum_{j \in {\mathcal J}}
d_{i,j}c^{\rho_{i,j}} \ne 0
$
for all ${\mathcal J} \subseteq \{0,\dots, j_i^0 \}$.
Then we have
$$
\aligned
{\vec {\frak v}}_T(\frak g(cT^{\kappa})) & =
{\vec {\frak v}}_T({\frak g}_1(cT^{\kappa}), \dots , {\frak g}_n(cT^{\kappa})) \\
& = (\lambda_1(\kappa), \cdots ,\lambda_n(\kappa)),
\endaligned
$$
which we denote by $\lambda(\kappa)$.
The function
$\lambda(\kappa)$ is piecewise affine. In particular, it is continuous.
(\ref{inIntP}) implies that
$\lambda(\kappa) \in P$ for sufficiently large $\kappa$.
Indeed, $\lambda(\kappa)$ is constant
for sufficiently large $\kappa$.
It actually implies $\lambda(\kappa) \in P$ for all $\kappa > 0$.
In fact if $\lambda(\kappa) \notin P$ for some $\kappa$ then
there exists $\kappa_0$ such that $\lambda(c\kappa_0) \in \partial P$.
This implies that $(\frak y'_1,\dots,\frak y'_n) = \frak g(cT^{\kappa_0})$ is a
solution of
$
\frak B_i = 0
$
and its valuation lies on the boundary of moment polytope.
This is impossible since the leading order equation of $\frak B_i$ is the same as
the leading order equation of $y_i\frac{\partial \frak{PO}_{\frak b}}{\partial y_i}$.
\par
In a similar situation, \cite[Section 5]{fooo09}  provides an example where
$\kappa \mapsto (\lambda_1(\kappa),\dots,\lambda_n(\kappa))$ is a
non-constant affine map. We can also find an example
where it is piecewise affine in a similar way.
\end{rem}
\par
The following is a consequence of
Proposition \ref{componentmajiwari}.
\begin{cor}\label{fiberdim0}
If $i \in \frak I$, then for any $\frak s \in \Lambda$ we have:
\begin{equation}\label{fiberdim0form1}
{\rm rank}_{\Lambda} \frac{\Lambda [s,y,y^{-1}]}{\frak q_i  (s-\frak s)} < \infty.
\end{equation}
We also have
\begin{equation}\label{fiberdim0form2}
\dim \frak X \le 1.
\end{equation}
\end{cor}
\begin{proof}
Sublemma \ref{SRlinearisenough}
implies that the fiber $\pi^{-1}(\frak s) \cap \frak Y$
is $0$ dimensional. Then Proposition
\ref{componentmajiwari} implies that $\pi^{-1}(\frak s) \cap \frak X$ is $0$ dimensional.
This implies (\ref{fiberdim0form1}).  (\ref{fiberdim0form2}) is its immediate consequence.
\end{proof}
\begin{prop}\label{compatnessfiber0}
$\pi: \frak X \to Spec(\Lambda[s])$ is projective.
\end{prop}
\begin{proof}
(The argument below is of flavor similar to the valuative criterion of properness
(See  \cite[Theorem 4.7 Chapter II]{Hart}).)
\par
We denote by $\P^1_{y_i}$ the projective line $Proj(\Lambda[y_0,y_i])$.
We have $\prod_{i=1}^n  \P^1_{y_i} \times Spec(\Lambda[s]) \supset \frak X$.
(Here the product is a fiber product over $Spec(\Lambda)$.)
Let $\overline{\frak X}$ be the Zariski closure of $\frak X$ in
$\prod_{i=1}^n  \P^1_{y_i} \times Spec(\Lambda[s])$.
It suffices to show
$\overline{\frak X} = \frak X$.
We will prove this by contradiction.
\begin{lem}
\begin{equation}\label{differencefinite}
\# (\vert\overline{\frak X}\vert \setminus
\vert \frak X \vert ) < \infty.
\end{equation}
\end{lem}
\begin{proof}
Let $\frak X_a$ be an irreducible component of $\frak X$.
Then (\ref{fiberdim0form2}) implies that $\frak X_a$ is a one dimensional irreducible variety 
and is Zariski closed in the affine space $\vert Spec(\Lambda[s,y,y^{-1}]) \vert$.
We denote by $\overline{\frak X}_a$ the Zariski closure of $\frak X_a$ in 
$\prod_{i=1}^n \P_{y_i}^1 \times \vert Spec(\Lambda[s])\vert$.
Clearly $\overline{\frak X}_a$ is  one dimensional irreducible algebraic variety and 
$\overline{\frak X}_a \setminus \vert Spec(\Lambda[s,y,y^{-1}])\vert$ is a
proper subvariety. Therefore 
$\overline{\frak X}_a \setminus \vert Spec(\Lambda[s,y,y^{-1}]) \vert$ is $0$ dimensional 
and hence the set of its geometric points is of finite order.
Note $\vert \overline{\frak X}_a\vert \cap \vert Spec(\Lambda[s,y,y^{-1}])\vert
= \vert \frak X_a\vert$. The lemma follows.
\end{proof}
Let $(\frak y,\frak s) \in \overline{\frak X} \setminus \frak X$.
We may replace the coordinate $y_i$ by
$y_i^{-1}$ and may assume that
$\frak y \in Spec(\Lambda[y_1,\dots,y_n])$.
(The map $y_i \mapsto y_i^{-1}$ is an automorphism of $\P^1_{y_i}$.)
\par
We put $\frak y = (\frak y_1,\ldots,\frak y_n)$.
Since $(\frak y, \frak s) \notin \frak X$, there exists $i$ such that
$\frak y_i = 0$.
We take one dimensional scheme $\frak C \subset \overline{\frak X}$
such that $(\frak y, \frak s) \in \frak C$ and $\frak C \setminus \{(\frak y, \frak s)\}
\subset \frak X$.
Replacing $\frak C$ by its normalization $Spec(B)$ we may assume that
$(\frak y, \frak s)$ is a smooth point of $Spec(B)$.
(In the proof of Proposition \ref{componentmajiwari}
we denote it by $\tilde{\frak y}$.)
Let $\frak m_{(\frak y, \frak s)}$ be the maximal ideal of $B$ corresponding to
$(\frak y, \frak s)$ and let $S$ be a generator of 
$\frak m_{(\frak y, \frak s)}B_{\frak m_{(\frak y, \frak s)}}$.
\par
In the same way as the proof of Proposition
\ref{componentmajiwari} we can find $K$ such that if $\kappa >K$ then
there exists a geometric point 
$\tilde{\frak y}(T^{\kappa}) = (\frak s(T^{\kappa}),\frak y_1(T^{\kappa}),\dots,\frak y_n(T^{\kappa}))$ of 
$Spec (B)$ such that
$\tilde{\frak y} - \tilde{\frak y}(T^{\kappa}) \equiv 0 \mod T^{\kappa}$.
If $\frak y_i \ne 0$, then $\frak v_T(\frak y_i(T^{\kappa}) ) = \frak v_T(\frak y_i)$ for sufficiently
large $\kappa$.
If $\frak y_i = 0$, then $\frak v_T(\frak y_i(T^{\kappa}) ) \ge \kappa$.
Since $\frak y_i = 0$ for at least one $i$, it implies that
$\vec{\frak v}_T(\frak y(T^{\kappa})) \notin P$ for sufficiently large $\kappa$.
Therefore $Spec (B)$ has infinitely many geometric points outside $\frak X$.
This contradicts to $Spec (B) \subset \overline{\frak X}$ and (\ref{differencefinite}).
\end{proof}
\begin{prop}\label{dimX}
We have
\begin{enumerate}
\item $\dim \frak X = 1$.
\item $\pi: \frak X \to Spec(\Lambda[s])$ is flat.
\item Let $\frak m_{(\frak y,\frak s)}$ be the
maximal ideal corresponding to $(\frak y,\frak s)\in \vert\frak X\vert$. Then  $\frak B_1,\dots,\frak B_n$  is a regular sequence
in $\Lambda[s,y,y^{-1}]_{\frak m_{(\frak y,\frak s)}}$.
\end{enumerate}
\end{prop}
\begin{proof}
Proof of 3:
Let $\frak x = (\frak y,\frak s) \in \vert\frak X\vert$ and
$\frak p_{\frak x}$ or $\frak m_{\frak x}$ the ideals
of $\mathcal R$ or $\Lambda[s,y,y^{-1}]$ corresponding to $\frak x$, respectively.
By Corollary \ref{fiberdim0}, we have
$$
0 = \dim \frac{\mathcal R_{\frak p_{\frak x}}}{(s-\frak s)} =
\dim \frac{\Lambda[s,y,y^{-1}]_{\frak m_{\frak x}}}{(\frak B_1,\dots,\frak B_n,s-\frak s)}.
$$
It follows that $(\frak B_1,\dots,\frak B_n,s-\frak s)$ is a system of parameter of
a Cohen-Macaulay ring $\Lambda[y,y^{-1},s]_{\frak m_{\frak x}}$.
(See  \cite[12.J]{Ma70}.) Therefore it is a regular sequence.
(See  \cite[Theorem 31]{Ma70}.)
Statement 3 follows.
\par
Proof of 2:
Note that the flatness of $\pi$ is equivalent to the flatness for all localizations at all maximal ideals of
$Spec(\Lambda[s])$.
Since the fibers of $\pi: \frak X \to Spec(\Lambda[s])$ are zero dimensional,
it suffices to show that
$s-\frak s$ is not a zero
divisor in $\mathcal R_{\frak p_{\frak x}}$ for any $\frak x = (\frak y,\frak s) \in \vert\frak X\vert$.
This follows from the fact that $\frak B_1,\dots,\frak B_n,s-\frak s$ is a regular sequence
in $\Lambda[s,y,y^{-1}]_{\frak m_{\frak x}}$.
\par
Proof of 1:
By Proposition \ref{dimX}.3, $\frak B_1,\dots,\frak B_n$
is a regular seqence in $\Lambda[s,y,y^{-1}]_{\frak m_{\frak x}}$,
where $\frak m_{\frak x}$ is any maximal ideal 
of $\frak X$.
Therefore 
$\dim \frak X = \dim \Lambda[s,y,y^{-1}] - n = 1$.
(See  \cite[Theorem 29]{Ma70}.)
\end{proof}
\begin{prop}\label{flattsuika}
For any $i\in \frak I$ we have the following.
\begin{enumerate}
\item $\dim \frak X_i = 1$.
\item
$
\pi : \frak X_i \to Spec(\Lambda[s])
$
is flat.
\end{enumerate}
\end{prop}
\begin{proof}
Proof of 1:
By Proposition \ref{dimX}.1, $\dim \frak X_i $ is either $1$ or $0$.
Suppose $\dim \frak X_i = 0$.
Let $\frak x \in \vert \frak X_i\vert$.
By Proposition \ref{dimX}.3,
$\frak B_1,\dots,\frak B_n,s-\frak s$ is a regular sequence
in $\Lambda[y,y^{-1},s]_{\frak m_{\frak x}}$.
By Lemma \ref{primalyremainlocal} below, $\frak q_i \Lambda[s,y,y^{-1}]_{\frak m_{\frak x}}$
appears in the irredundant primaly decomposition
of $(\frak B_1,\dots,\frak B_n) \Lambda[s,y,y^{-1}]_{\frak m_{\frak x}}$.
Then by Theorem \ref{unmix} below,
the height of $\frak q_i$ in $\Lambda[s,y,y^{-1}]_{\frak m_{\frak x}}$ is $n$.
Namely
$$
\dim \frac{\Lambda[s,y,y^{-1}]_{\frak m_{\frak x}}}{\frak q_i\Lambda[s,y,y^{-1}]_{\frak m_{\frak x}}}
= 1.
$$
This contradicts to $\dim \frak X_i = 0$.
Thus we obtain 1.
\par
To prove 2 we prepare some elementary lemmata and a basic fact.
\begin{lem}\label{primalyremainlocal}
Let $R$ be an integral domain and $\frak q$, $\frak Q$  primaly ideals of $R$.
Suppose that $\frak p$ is a prime ideal containing both $\frak q$ and $\frak Q$.
If
$
\frak q R_{\frak p} \subseteq \frak Q R_{\frak p},
$
then $\frak q \subseteq \frak Q$.
\end{lem}
\begin{proof}
Let $x \in \frak q$. Then by assumption we have
$a \in R$, $b \in R \setminus \frak p$ and $c \in \frak Q$ such that
$
x = c a/b
$
in $R_{\frak p}$.
It follows that $xb = ca$ in $R$. (We use the fact that $R$ is an integral domain here.)
Therefore $xb \in \frak Q$. Since $b \in R \setminus \frak p$ and $\frak p \supseteq \frak Q$
is a prime ideal, it follows that $x \in \frak Q$. We thus prove $\frak q \subseteq \frak Q$.
\end{proof}
\begin{lem}\label{torsion}
Let $R$ be a commutative ring with unit, $A$ a commutative $R$-algebra and $J_0$ an ideal in $A$.
Suppose that $J_0=\bigcap_{i\in \frak I} {\frak q}_i$ is an irredundant primary decomposition
without embedded primes.
Then the following two conditions are equivalent:
\begin{enumerate}
\item  $A/J_0$ is torsion-free as an $R$-module.
\item For any $i$, $A/{\frak q}_i$ is torsion-free as an $R$-module.
\end{enumerate}
\end{lem}
\begin{proof}
Suppose that $A/{\frak q}_i$ has a torsion, i.e., there exist $a \in R \setminus \{0\}$
and $f \in A \setminus {\frak q}_i$ such that $f \cdot a \in {\frak q}_i$.
Since the decomposition $J_0=\bigcap_i {\frak q}_i$ has no embedded primes,
we can pick $g_j \in {\frak q}_j \setminus \sqrt{{\frak q}_i}$ for each $j \neq i$.
Set $h=f \prod_{j\neq i} g_j$.  Then we find that $h \notin {\frak q}_i$, hence $h \notin J_0$,
and $ah=(af) \prod_{j \notin i}g_j \in J_0$.  Namely, $h$ is a torsion element in $A/J_0$.
The converse implication is clear.
\end{proof}
\begin{thm}[{\rm Macauley's unmixedness theorem. cf.  \cite{Eisen} Corollary 18.14}]\label{unmix}
Let $R$ be a Cohen-Macauley ring and $I$ an ideal generated
by an $R$-regular sequence $(x_1,\ldots,x_k)$, then $I$ has no embedded primes and
the height of $I$ is $k$.
\end{thm}
\par
\noindent
{\it Proof of Proposition \ref{flattsuika}.2:}
Let $\frak x = (\frak y,\frak s) \in \vert \frak X_i\vert$.
By definition (\ref{defJ})
we have
$$
\aligned
\mathcal R_{\frak m_{\frak x}} & = \frac{\Lambda[s,y,y^{-1}]_{\frak m_{\frak x}}}
{(\frak B_1,\dots,\frak B_n)_{{\frak m}_{\frak x}}}
=
\frac{\Lambda[s,y,y^{-1}]_{\frak m_{\frak x}}}{\bigcap_{i \in {\frak I}_+} (\frak q_i)_{{\frak m}_{\frak x}}} \\
& =
\frac{\Lambda[s,y,y^{-1}]_{\frak m_{\frak x}}}{\bigcap_{i \in {\frak I}} (\frak q_i)_{{\frak m}_{\frak x}}}
=\frac{\Lambda[s,y,y^{-1}]_{\frak m_{\frak x}}}{\bigcap_{\frak q_i \subset {\frak m}_{\frak x}} (\frak q_i)_{{\frak m}_{\frak x}}},
\endaligned$$
because for $i \in {\frak I}_+ \setminus {\frak I}$
we have $(\frak q_i)_{{\frak m}_{\frak x}}=(1)$. 
By Proposition \ref{dimX}.3, $\frak B_1,\dots,\frak B_n$
is a regular sequence in $\Lambda[s,y,y^{-1}]_{\frak m_{\frak x}}$.
Therefore, Proposition \ref{dimX}.2, Theorem \ref{unmix} and Lemma \ref{torsion} above
imply the flatness of $
\pi : \frak X_i \to Spec(\Lambda[s])
$.
\end{proof}
During the proof of Proposition \ref{flattsuika}, we have showed that
irredundant primary decompositions (\ref{I0decompose}) has no embedded primes.
\par
\begin{lem}\label{XYgepoint}
In the situation of Proposition \ref{componentmajiwari}
we have
$\vert \frak X_i\vert =  \vert \frak Y_j \vert$.
\end{lem}
\begin{proof}
By Corollary \ref{I0Icoro} we have
\begin{equation}\label{uderRtoR}
\vert Spec(\underline{\mathcal R}) \vert \subseteq \vert \frak X\vert.
\end{equation}
Therefore $\dim \frak Y_j  \le 1$.
On the other hand, $\dim \frak X_i = 1$ by Proposition \ref{flattsuika}.
Therefore the lemma follows from Proposition \ref{componentmajiwari}.
\end{proof}
For each $i\in \frak I$ we take $ j(i)$ such that
$\frak X_i$ and $\frak Y_{j(i)}$ are as in Proposition \ref{componentmajiwari}.
\begin{prop}\label{coincidenceasschemes}
We have $\frak X_i = \frak Y_{j(i)}$, i.e., $\frak q_i = \frak Q_{j(i)}$.
\end{prop}
\begin{proof}
By Corollary \ref{I0Icoro} we have $I_0 \subseteq \overline I$. It follows
from Lemma \ref{XYgepoint} that
$$
\vert \frak X \vert = \vert Spec(\underline{\mathcal R}) \vert.
$$
Therefore by Lemma \ref{XYgepoint} again we have
$\dim \frak Y_j = 0$, i.e.,the associated prime of  $\frak Q_j$ is an
embedded
prime, if $j \notin \{ j(i) \mid i \in \frak I\}$.

We put
$$
B_1=  \bigcup_{j \notin \{ j(i) \mid i \in \frak I\}} \pi(\vert \frak Y_j
\vert).
$$
$B_1 \subset \Lambda$ is a finite set. We put
$$
B_2 = \bigcup_{i_1,\ne i_2; i_1,i_2 \in \frak I}  \pi(\vert \frak X_{i_1}
\vert \cap \vert \frak X_{i_j} \vert).
$$
$B_2 \subset \Lambda$ is also a finite set.  Let $B = B_1 \cup B_2$.
We take $c \in \Lambda \setminus \pi(B)$.
Let $\frak p$ be the prime ideal $\sqrt{{\frak q}_i}=\sqrt{{\frak
Q}_{j(i)}}$.
Pick a maximal ideal ${\mathfrak m}$ containing $\frak p$ and $s-c$.
In other words, the maximal ideal $\mathfrak m={\mathfrak m}_{\mathfrak
x}$ defines a point
$${\mathfrak x}=(c,{\mathfrak y})
\in
\vert {\frak X}_i \vert \cap \pi^{-1}(c) = \vert {\mathfrak Y}_{j(i)}
\vert \cap \pi^{-1}(c).
$$
We denote by ${\mathfrak m}_{\mathfrak y}$ the maximal ideal in
$\Lambda[y,y^{-1}]$, which defines  $\frak y \in (\Lambda^{\times})^n$.
\par
Now we compare the rings $\{{\mathcal R}/(s-c){\mathcal R}\}_{\frak m_{\frak y}}$ and
$\{\underline{\mathcal R}/(s-c)\underline{\mathcal R}\}_{\frak m_{\frak y}}$.
By Lemma \ref{inmompoly} and Proposition \ref{componentmajiwari}, we have
$\vec{\frak v}_T(\frak y) \in P$.
In the same way as in  \cite[Section 7]{fooo08}, these rings are
characterized as
the common generalized eigenspaces of the multiplication by $y_i$ with
eigenvalues $\frak y_i$.
Then by the same way as in the proof of  \cite[Proposition 8.6]{fooo08} we
can show that the natural projection
${\mathcal R} \to
\underline{\mathcal R}$
induces an isomorphism
$$
\left({\mathcal R}/{(s-c) {\mathcal R}}\right)_{\frak m_{\frak y}} \cong
\left({\underline{\mathcal R}}/{(s-c)\underline{\mathcal R}}\right)_{\frak
m_{\frak y}}.
$$
We first consider
the kernel of the projection
${\mathcal R}_{{\mathfrak m}_{\mathfrak x}} \to
\underline{\mathcal R}_{{\mathfrak m}_{\mathfrak x}}$
which is $\overline{I}_{{\mathfrak m}_{\mathfrak x}}/I_{{\mathfrak m}_{\mathfrak
x}}$.
We have the following:

\begin{lem}\label{qQ_m}
We have
$$
I_{{\mathfrak m}_{\mathfrak x}} = ({\mathfrak q}_i)_{{\mathfrak
m}_{\mathfrak x}}, \qquad
\overline{I}_{{\mathfrak m}_{\mathfrak x}}=({\mathfrak
Q}_{j(i)})_{{\mathfrak m}_{\mathfrak x}}.
$$
In particular,
$$
\overline{I}_{{\mathfrak m}_{\mathfrak x}}/I_{{\mathfrak m}_{\mathfrak x}}
\cong
({\mathfrak Q}_{j(i)})_{{\mathfrak m}_{\mathfrak x}}/({\mathfrak
q}_i)_{{\mathfrak m}_{\mathfrak x}}.
$$
\end{lem}

\begin{proof} Recall from (\ref{primary}) that
$
I=\bigcap_{i \in \frak I_+} {\frak q_i}$ and
$\overline{I}=\bigcap_{j \in  \frak J} {\frak Q_j}$.
We note that ${\frak q}_k$ for
$k \in {\frak I}_+ \setminus {\frak I}$ cannot be 
contained in ${\mathfrak m}_{\mathfrak x}$,
because $\vert {\mathfrak X}_k \vert \cap \vert Spec (\underline{\mathcal
R})\vert = \emptyset$.
Therefore we have
$$
(I_0)_{{\mathfrak m}_{\mathfrak x}}=I_{{\mathfrak m}_{\mathfrak x}} =
\bigcap_{\mathfrak q_j  \subset {\mathfrak m}_{\mathfrak x}}
(q_j)_{{\mathfrak m}_{\mathfrak x}}.
$$
We also have
$$
\overline{I}_{{\mathfrak m}_{\mathfrak x}}=
\bigcap_{\mathfrak Q_k \subset
{\mathfrak m}_{\mathfrak x}} (\mathfrak Q_k)_{{\mathfrak m}_{\mathfrak
x}}.
$$

By the choice of $c$, we find that
${\mathfrak X}_i$, resp. ${\mathfrak Y}_{j(i)}$,
is the only irreducible component of $\vert {\mathfrak X} \vert$, resp.
$\vert {\mathfrak Y} \vert $
containing the $\Lambda$-point ${\mathfrak x}$, and
${\mathfrak q}_i$, resp. ${\mathfrak Q}_{j(i)}$,
is the unique primary ideal in the irredundant primary decompositions
(\ref{primary})  of $I$, resp.
$\overline{I}$, which is contained in ${\mathfrak m}_{\mathfrak x}$.
Hence we obtain
$$I_{{\mathfrak m}_{\mathfrak x}} = (\mathfrak q_{i})_{{\mathfrak
m}_{\mathfrak x}} \subset \overline{I}_{{\mathfrak m}_{\mathfrak
x}}=(\mathfrak Q_{j(i)})_{{\mathfrak m}_{\mathfrak x}}.$$
\end{proof}

The rings ${\mathcal R}/(s-c){\mathcal R}$ and
$\underline{\mathcal R}/(s-c) \underline{\mathcal R}$ are
$\Lambda[s,y,y^{-1}]$-algebras.  By substituting $c$ to $s$, we can also
regard them as $\Lambda[y,y^{-1}]$-algebras.  Then we have
$$
({\mathcal R}/(s-c){\mathcal R})_{{\mathfrak m}_{\mathfrak y}}
\cong {\mathcal R}_{{\mathfrak m}_{\mathfrak x}}/ (s-c) {\mathcal
R}_{{\mathfrak m}_{\mathfrak x}}
$$
and
$$
(\underline{\mathcal R}/(s-c)\underline{\mathcal R})_{{\mathfrak
m}_{\mathfrak y}}
\cong \underline{\mathcal R}_{{\mathfrak m}_{\mathfrak x}}/ (s-c)
\underline{\mathcal  R}_{{\mathfrak m}_{\mathfrak x}}.
$$

By utilizing Lemma \ref{qQ_m}, we have
the following commutative diagram such that the horizontal and vertical
lines are exact sequences.

\begin{equation}
\begin{CD}
 @.  0 @>>> 0 @>>> 0 @. \\
@. @VVV @VVV @VVV @.  \\
0 @>>>  {\text{\rm Ker }} q @>>> ({\mathfrak Q}_{j(i)}/{\mathfrak q}_i)_{{\mathfrak
m}_{\mathfrak x}} @>>> 0 @>>> 0 \\
@VVV @VVV @VVV @VVV @VVV \\
0 @ >>> (s-c) {\mathcal R}_{{\mathfrak m}_{\mathfrak x}} @>>>
{\mathcal R}_{{\mathfrak m}_{\mathfrak x}} @>>>
({\mathcal R}/(s-c) {\mathcal R})_{{\mathfrak m}_{\mathfrak x}} @>>> 0  \\
@VVV @ VqVV @VVV  @VVV @VVV \\
0 @ >>> (s-c) \underline{\mathcal R}_{{\mathfrak m}_{\mathfrak x}} @>>>
\underline{\mathcal R}_{{\mathfrak m}_{\mathfrak x}} @>>>
(\underline{\mathcal R}/(s-c) \underline{\mathcal R})_{{\mathfrak
m}_{\mathfrak x}} @>>> 0 \\
@. @VVV @VVV @VVV @.\\
@.  0 @>>> 0 @>>> 0 @.
\end{CD} \nonumber
\end{equation}
Note that
$$ {\text{\rm Ker }} q = (s-c) {\mathcal R}_{{\mathfrak m}_{\mathfrak x}} \cap \left(
{\mathfrak Q}_{j(i)}/{\mathfrak q}_i \right)_{{\mathfrak m}_{\mathfrak x}}
\cong \left({\mathfrak q}_i + (s-c) \Lambda[s,y,y^{-1}] \cap {\mathfrak
Q}_{j(i)} \right)_{{\mathfrak m}_{\mathfrak x}}/({\mathfrak
q}_i)_{{\mathfrak m}_{\mathfrak x}}.$$
Since ${\mathfrak Q}_{j(i)}$ is a primary ideal and $s-c \notin
\sqrt{{\mathfrak Q}_{j(i)}}$, we have
$$(s-c) \Lambda[s,y,y^{-1}] \cap {\mathfrak Q}_{j(i)} =(s-c) {\mathfrak
Q}_{j(i)}.$$
Therefore we have
$$
{\text{\rm Ker }} q \cong (s-c) ({\mathfrak Q}_{j(i)}/{\mathfrak q}_i)_{{\mathfrak
m}_{\mathfrak x}}.
$$
It follows from the above diagram that
$$(s-c) ({\mathfrak Q}_{j(i)}/{\mathfrak q}_i)_{{\mathfrak m}_{\mathfrak
x}} \cong
({\mathfrak Q}_{j(i)}/{\mathfrak q}_i)_{{\mathfrak m}_{\mathfrak x}}.
$$

Now we set $M={\mathfrak Q}_{j(i)}/{\mathfrak q}_i$.
What we have showed is
$$
M_{{\mathfrak m}_{\mathfrak x}}/
(s-c)M_{{\mathfrak m}_{\mathfrak x}} =
\left( M/(s-c)M\right) _{{\mathfrak m}_{\mathfrak x}}
= 0.
$$
Then we obtain
\begin{equation}\label{Mmzero}
M_{{\mathfrak m}_{\mathfrak y}}/
(s-c)M_{{\mathfrak m}_{\mathfrak y}}=
\left( M/(s-c)M\right) _{{\mathfrak m}_{\mathfrak y}} \cong
\left( M/(s-c)M\right) _{{\mathfrak m}_{\mathfrak x}}
= 0.
\end{equation}
When
we localize $M/(s-c)M$ at the maximal ideal
${\mathfrak m}_{\mathfrak x}$
in $\Lambda[s,y,y^{-1}]$,
we regard it as a $\Lambda[s,y,y^{-1}]$-module
such that $s$ acts on it by the multiplication by $c$.
When we localize $M/(s-c)M$ at the maximal ideal
${\mathfrak m}_{\mathfrak y}$
in $\Lambda[y,y^{-1}]$,
we regard it as a
$\Lambda[y,y^{-1}]$-module.
Therefore we obtain the isomorphism of the second equality
in (\ref{Mmzero}).
\par
To apply Nakayama's lemma, we show the following lemma whose
proof is given at the end of the proof of Porposition
\ref{coincidenceasschemes}.

\begin{lem}\label{Mmfg}
$M_{{\mathfrak m}_{\mathfrak y}}$ is a finitiely generated module over
$\Lambda[s]$.
\end{lem}
We localize $\Lambda[s]$ and $M_{{\mathfrak m}_{\mathfrak y}}$ at the
maximal ideal $(s-c) \subset \Lambda[s]$ to obtain
$$(M_{{\mathfrak m}_{\mathfrak y}})_{(s-c)}=(s-c)
(M_{{\mathfrak m}_{\mathfrak y}})_{(s-c)}$$
by (\ref{Mmzero}). Then
Nakayama's lemma implies that
$(M_{{\mathfrak m}_{\mathfrak y}})_{(s-c)}=0$.

We can regard $(M_{{\mathfrak m}_{\mathfrak y}})_{(s-c)}$ as
a $\Lambda[s,y,y^{-1}]$-module and localize it at
the maximal ideal ${\mathfrak m}_{\mathfrak x}$ in $\Lambda[s,y,y^{-1}]$. Then
by noticing the fact that
${\mathfrak m}_{\mathfrak y}= {\mathfrak m}_{\mathfrak x} \cap
\Lambda[y,y^{-1}]$ and
$(s-c) \Lambda[s] = {\mathfrak m}_{\mathfrak x} \cap \Lambda[s]$,
we find that
$$M_{{\mathfrak m}_{\mathfrak x}} =
\left( (M_{{\mathfrak m}_{\mathfrak y}})_{(s-c)} \right)_{{\mathfrak
m}_{\mathfrak x}}.
$$
Therefore we have
$M_{{\mathfrak m}_{\mathfrak x}}=0$, i.e.,
$({\mathfrak q}_i)_{{\mathfrak m}_{\mathfrak x}}
=({\mathfrak Q}_{j(i)})_{{\mathfrak m}_{\mathfrak x}}$.
Then  Lemma \ref{primalyremainlocal} implies that
${\mathfrak q}_i = {\mathfrak Q}_{j(i)}$.
Hence we obtain $\frak X_i=\frak Y_{j(i)}$ as schemes.
This finishes the proof of Proposition \ref{coincidenceasschemes} modulo Lemma \ref{Mmfg}.
\end{proof}
\begin{proof}[Proof of Lemma \ref{Mmfg}]
By Corollary \ref{fiberdim0},
${\rm dim}_{\Lambda} \frac{\Lambda [s,y,y^{-1}]}{\frak q_i + (s-\frak s)}$
is finite.
Thus it is an Artinian ring. 
Hence the Krull dimension 
$
\left( \frac{\Lambda [s,y,y^{-1}]}{\frak q_i + (s-\frak s)
}\right)_{{\mathfrak m}_{\mathfrak y}}
$ 
is $0$.
Since $\{{\frak B}_1, \dots, {\frak B}_n, s-c\}$ is a regular sequence,
the Krull dimension of
$\left(\frac{\Lambda [s,y,y^{-1}]}{\frak q_i}
\right)_{{\mathfrak m}_{\mathfrak y}}$
is $1$.
Since $s$ is transcendental over $\Lambda$, we find that
$\left(\frac{\Lambda [s,y,y^{-1}]}{\frak q_i}
\right)_{{\mathfrak m}_{\mathfrak y}}$
is a finitely generated $\Lambda[s]$-module.
From the exact sequence
$$
0 \to  \left({\mathfrak Q}_{j(i)}/{\mathfrak q}_i
\right)_{{\mathfrak m}_{\mathfrak
y}} \to \left(\frac{\Lambda [s,y,y^{-1}]}{\frak q_i}
\right)_{{\mathfrak m}_{\mathfrak
y}} \to
\left(\frac{\Lambda [s,y,y^{-1}]}{\frak Q_j(i)}
\right)_{{\mathfrak m}_{\mathfrak y}}
\to 0,
$$
we find that
$M_{{\mathfrak m}_{\mathfrak y}}=
\left({\mathfrak Q}_{j(i)}/
{\mathfrak q}_i
\right)_{{\mathfrak m}_{\mathfrak y}}$
is a finitely generated
$\Lambda[s]$-module.
\end{proof}
\begin{cor}\label{I0andJ}
$I_0 = \overline I_0= \overline I$. $\frak X = Spec(\underline{\mathcal R})$.
\end{cor}
\begin{proof}
By Proposition \ref{coincidenceasschemes}
$
I_0 = \bigcap_{j(i) : i \in \frak I} \frak Q_{j(i)}.
$
On the other hand we have
$
 \bigcap_{j(i) : i \in \frak I} \frak Q_{j(i)}\supseteq  \overline I  =\overline I_0
\supseteq I_0.
$
Hence
$
\overline I  =\bigcap_{j(i) : i \in \frak I} \frak Q_{j(i)} = I_0,
$
as required.
\end{proof}
We are now in the position to complete the proof of Proposition \ref{bettiest}.
By Proposition \ref{compatnessfiber0},
the coordinate ring $\underline{\mathcal R}$ of $\frak X$ is a finitely generated $\Lambda[s]$ module.
On the other hand, since $\Lambda[s]$ is a principal ideal domain, Proposition \ref{dimX} implies that
$\underline{\mathcal R}$ is a finitely generated
free module over $\Lambda[s]$.
In particular, we find that $\dim_{\Lambda} \underline{\mathcal R}/(s-0)=\dim_{\Lambda}\underline{\mathcal R}/(s-1)$.
Therefore
$$\aligned
\dim_{\Lambda} H(X;\Lambda)
&= \dim_{\Lambda} \underline{\mathcal R}/(s-1)
= \dim \underline{\mathcal R}/(s-0)
\\
&= \dim_{\Lambda} \text{\rm PJac}_{\frak v^P_T}(\frak{PO}_{\frak b})_{\Lambda}
= \dim_{\Lambda} \text{\rm Jac}_{\frak v^P_T}(\frak{PO}_{\frak b})_{\Lambda}.
\endaligned$$
Here we use Lemma \ref{polyapprox} to prove the last
equality.
\qed\par
In geometric terms, we may also put it as follows.
We have (\cite[Theorems 9.9 III and 12.8 III]{Hart})
$$\aligned
\dim_{\Lambda} H(X;\Lambda)
&= \dim_{\Lambda} H^0(\pi^{-1}(1);\frak O_{\pi^{-1}(1)})
= \dim H^0_{\Lambda}(\pi^{-1}(0);\frak O_{\pi^{-1}(0)})
\\
&= \dim_{\Lambda} \text{\rm PJac}_{\frak v^P_T}(\frak{PO}_{\frak b})_{\Lambda}
= \dim_{\Lambda} \text{\rm Jac}_{\frak v^P_T}(\frak{PO}_{\frak b})_{\Lambda}.
\endaligned$$
\par

\subsection{The case $\frak b \in \mathcal A(\Lambda_0)$}
\label{subsec:caseofLambdaz}
In this subsection, we
generalize Theorem \ref{injectivity} to the case $\frak b \in \mathcal A(\Lambda_0)$
and complete the proof of Theorem \ref{Mirmain}.1.
We choose any $\frak P'\in \Lambda[y,y^{-1}]_0^P$ such that
$\frak P'$ is transformed to an element $\frak{PO}_{\frak b'}$
for some $\frak b' \in \mathcal A(\Lambda_+)$.
(Theorem \ref{algebraization}.2.)
We have
$$
\dim_{\Lambda} \text{\rm Jac}(\frak P')_{\Lambda} = \sum_{k}
\dim_{\Lambda} H^k(X;\Lambda)
$$
by Theorem \ref{injectivity}. Hence the next proposition
completes the proof of Theorem \ref{Mirmain}.1.
\qed
\begin{prop}\label{rankequalgeneral}
$$
\dim_{\Lambda} \text{\rm Jac}(\frak{PO}_{\frak b})_{\Lambda}
=\dim_{\Lambda} \text{\rm Jac}(\frak P')_{\Lambda}.
$$
\end{prop}
\begin{rem}\label{PeANDpotential}
During the proof of Proposition \ref{rankequalgeneral} 
we need restrict the potential function 
$\frak{PO}_{\frak b}$ 
from the interior of the polytope $\overset{\circ}P$
to the polytope 
$P_{\epsilon}=\{ u \in \R^n \mid \ell_j(u) \ge \epsilon, ~j=1,\dots , m\}$ for $\epsilon >0$.
We can do it as follows. 

Consider the following diagram:

\begin{equation}
\begin{CD}
\Lambda_0 [[Z_1, \dots, Z_m]] @>{\phi}>> \Lambda
\langle\!\langle y,y^{-1} \rangle\!\rangle_0^{\overset{\circ}P} \\
@V{\widetilde{j}}VV @V{j}VV  \\
\Lambda_0 \langle\!\langle Z_1^{(\epsilon)}, \dots, Z_m^{(\epsilon)}
\rangle\!\rangle @>{\phi^{(\epsilon)}}>>
\Lambda \langle\!\langle y,y^{-1} \rangle\!\rangle_0^{P_{\epsilon}}
\end{CD} \nonumber
\end{equation}
where $\phi$ and $\phi^{(\epsilon)}$ are the homomorphisms we have used:
$$\phi(Z_j)=T^{-\lambda_j} y_1^{v_{j,1}} \cdots y_n^{v_{j,n}}, \quad
\phi^{(\epsilon)}(Z_j^{(\epsilon)})=T^{-(\lambda_j+\epsilon)}
y_1^{v_{j,1}} \cdots y_n^{v_{j,n}},$$
$j$ is the inclusion, i.e., $j(y_i)=y_i$, and
$\widetilde{j}(Z_i)= T^{\epsilon} Z_i^{(\epsilon)}$.
It is straightforward to check the commutativity of the above diagram.
We write $z_j=\phi(Z_j)$ and
$z_j^{(\epsilon)}=\phi^{(\epsilon)}(Z_j^{(\epsilon)})$.
Then if ${\mathcal P} \in \Lambda
\langle\!\langle y,y^{-1} \rangle\!\rangle_0^{\overset{\circ}P}$ has no constant terms, 
we have
$$
j({\mathcal P}(z_1, \dots, z_m))={\mathcal P}(T^{\epsilon}
z_1^{(\epsilon)}, \dots, T^{\epsilon} z_m^{(\epsilon)})
$$
which belongs to $T^{\epsilon} \Lambda \langle\!\langle
y,y^{-1} \rangle\!\rangle_0^{P_{\epsilon}}$.
\end{rem}

\begin{proof}[Proof of Proposition \ref{rankequalgeneral}]
By Lemma \ref{leadingterm}, we may put
\begin{equation}\label{PObLT}
\frak{PO}_{\frak b}  \equiv c_1 z_1 + \cdots + c_mz_m  + R(z_1,\dots,z_m) \mod \Lambda[y,y^{-1}]_+^P
\end{equation}
where $R \in \Lambda_0[[Z_1,\dots,Z_m]]$ each summand of which has degree $\ge 2$,
and $c_i \in \C\setminus \{0\}$. 
We firstly take $\e >0$ so that
there is no critical point $\frak y$ of
$\frak {PO}_{\frak b}$ such that
$\vec{\frak v}_T(\frak y) \in \text{\rm Int} P \setminus P_{\epsilon}$.
By Theorem \ref{algebraization}
we have a strict coordinate change $y'(y)$ which converge on $P_{\epsilon}$
such that
$$
\frak{PO}_{\frak b} (y') = \frak P(y) \in \Lambda[y,y^{-1}]_0^{P_{\epsilon}}.
$$
(Namely it is a Laurent polynomial.)
Then
by Proposition \ref{Morsesplit},  we have
$$
\text{\rm Jac}(\frak{PO}_{\frak b} )_{\Lambda}
\cong
\frac{\Lambda\langle\!\langle y,y^{-1}\rangle\!\rangle^{P_{\epsilon}}}
{{\rm Clos}_{\frak v_T^{P_{\epsilon}}}\left(
y_i \frac{\partial \frak{PO}_{\frak b} }{\partial y_i}
: i=1,\dots,n
\right)}
\cong
\frac{\Lambda\langle\!\langle y,y^{-1}\rangle\!\rangle^{P_{\epsilon}}}
{{\rm Clos}_{\frak v_T^{P_{\epsilon}}}\left(
y_i \frac{\partial \frak P}{\partial y_i}
: i=1,\dots,n
\right)}.
$$
(Note we already proved Lemma \ref{Jacfinitedimension} when we proved the surjectivity of
Kodaira-Spencer map. Therefore the proof of Proposition \ref{Morsesplit} is already completed by now.)
Moreover
$$
\text{\rm Jac}(\frak P')_{\Lambda}
\cong
\frac{\Lambda\langle\!\langle y,y^{-1}\rangle\!\rangle^{P_{\epsilon}}}
{{\rm Clos}_{\frak v_T^{P_{\epsilon}}}\left(
y_i \frac{\partial \frak P'}{\partial y_i}
: i=1,\dots,n
\right)}.
$$
Thus it suffices to show that
$$
\dim_{\Lambda} \frac{\Lambda\langle\!\langle y,y^{-1}\rangle\!\rangle^{P_{\epsilon}}}
{{\rm Clos}_{\frak v_T^{P_{\epsilon}}}\left(
y_i \frac{\partial \frak P}{\partial y_i}
: i=1,\dots,n
\right)}
 = \dim_{\Lambda} \frac{\Lambda\langle\!\langle y,y^{-1}\rangle\!\rangle^{P_{\epsilon}}}
{{\rm Clos}_{\frak v_T^{P_{\epsilon}}}\left(
y_i \frac{\partial \frak P'}{\partial y_i}
: i=1,\dots,n
\right)}.
$$
We define
$$
\frak P_{\text{\rm para}} = \alpha s\frak P + (1-s)\frak P'
\in \Lambda[s,y,y^{-1}]_0^{P_{\epsilon}}.
$$
(Here $\alpha \in \C \setminus \{1/c_1,\dots,1/c_m,0\}$ is a constant.)
In fact, we have
$$
T^{-\epsilon}\frak P_{\text{\rm para}}
\in \Lambda[s,y,y^{-1}]_0^{P_{\epsilon}}.
$$
From now on we choose and take $\epsilon>0$  small so that the following holds:
For any $s \in \C$ and $i$ ($i=1,\dots,m$) with
$
\alpha sc_i + (1-s) \ne 0
$, 
there is no solution
of leading term equation of
$\frak P_{\text{\rm para}}$ on
$\text{\rm Int} P \setminus P_{\epsilon}$.
(See \cite[Definitin 4.3]{fooo09}  for the definition of leading term equation.)
We also assume that the leading term equation of
$\alpha \frak P - \frak P'$ also has no  solution on $\text{\rm Int} P \setminus P_{\epsilon}$.
(We use the assumption  $\alpha \in \C \setminus 
\{1/c_1,\dots,1/c_m,0\}$ here.)
\par
We put
$$
\frak  R =
\frac{\Lambda[s,y,y^{-1}]}{\left(
y_i \frac{\partial \frak P_{\text{\rm para}}}{\partial y_i}: i=1,\dots,n\right)},
\qquad
\underline{\frak  R} =
\frac{\Lambda[s,y,y^{-1}]}{\text{\rm Clos}_{\frak {v}_T^{P_{\e}}}\left(
y_i \frac{\partial \frak P_{\text{\rm para}}}{\partial y_i}: i=1,\dots,n\right)}.
$$
We have natural morphisms
$$
Spec(\underline{\frak R}) \longrightarrow Spec({\frak R}) \overset{\pi}\longrightarrow Spec(\Lambda[s])
= \Lambda.
$$
\par
If we put
$$
z_i^{\epsilon} = T^{-\epsilon}z_i,
$$
we derive from (\ref{PObLT})  that
\begin{equation}\label{PleadingtermR}
T^{-\epsilon}\frak P_{\text{\rm para}}
\equiv \sum_{i=1}^m (\alpha sc_i + (1-s)) z^{\epsilon}_i  \mod \Lambda[y,y^{-1}]_+^{P_{\epsilon}}.
\end{equation}
\par
Denote
$$
\aligned
B &= \{b_1,\ldots,b_A\} = \{ s \mid \exists i \,\,s\alpha c_i + (1-s)= 0\} \subset \C, \\
B_+ &= \{ b \in \Lambda \mid \exists i \,\,b \equiv b_i, \mod \Lambda_+, i=1,\dots,m\},  \quad
B_+^c = \Lambda \setminus B_+.
\endaligned$$
\par
Let
$$
I = \left(
y_i \frac{\partial \frak P_{\text{\rm para}}}{\partial y_i}: i=1,\dots,n\right)
$$
and let $\overline I$ be its closure with respect to the $\frak v_T^{P_{\epsilon}}$-norm.
We take their irredundant primary decompositions:
\begin{equation}\label{primary2}
I=\bigcap_{i \in \frak I_+} {\frak q_i}, \qquad \overline{I}=\bigcap_{j \in  \frak J_+} {\frak Q_j}.
\end{equation}
We define
\begin{equation}\label{defJ2}
\aligned
\frak I &= \left\{ i \in \frak I_+ ~\left\vert~
\left\vert Spec\left(\frac{\Lambda[s,y,y^{-1}]}{\frak q_i}\right) \right\vert
\cap \left\vert Spec(\underline{\frak R}) \right\vert
\cap \pi^{-1}(B_+^c) \ne \emptyset \right\} \right.,
\\
\frak J &= \left\{ j \in \frak J_+ ~\left\vert~
\left\vert Spec\left(\frac{\Lambda[s,y,y^{-1}]}{\frak Q_j}\right) \right\vert
\cap \pi^{-1}(B_+^c) \ne \emptyset \right\} \right..
\endaligned\end{equation}
We also put
$$
\frak X_i = Spec\left(\frac{\Lambda[s,y,y^{-1}]}{\frak q_i}\right),
\qquad
\frak Y_j =Spec\left(\frac{\Lambda[s,y,y^{-1}]} {\frak Q_j}\right).
$$
\begin{equation}\label{I0decompose2}
I_0 = \bigcap_{i\in \frak I} \frak q_i,
\qquad
J_0 = \bigcap_{j\in \frak J} \frak Q_j,
\end{equation}
and
$$
\frak X = Spec\left( \frac{\Lambda[s,y,y^{-1}]}{I_0}\right),  \qquad
\frak Y = Spec\left( \frac{\Lambda[s,y,y^{-1}]}{J_0 }\right).
$$
$\frak X_i$ and $\frak Y_j$ are irreducible components of
$\frak X$ and $\frak Y$, respectively.
\begin{lem}\label{inmopoly2}
If $(\frak y_1,\dots,\frak y_n,\frak s) \in \vert \frak Y\vert$
and $\frak s \in B_+^c$,
then
\begin{equation}\label{inmompoly2}
(\frak v_T(\frak y_1),\dots,\frak v_T(\frak y_n)) \in P_{\e}.
\end{equation}
\par
If $(\frak y_1,\dots,\frak y_n,\frak s)\in \vert \frak X \vert$,
$\frak s  \in B_+^c$
and $(\ref{inmompoly2})$ is satisfied, then
$(\frak y_1,\dots,\frak y_n,\frak s)\in \vert \frak Y\vert$.
\end{lem}
\begin{proof}
Let $\frak s \in B_+^c$. Then we can show
\begin{equation}
\frac{\underline{\frak R}}{(s-\frak s)}
\cong \frac{\Lambda\langle\!\langle y,y^{-1}\rangle\!\rangle^{P_{\epsilon}}}
{\text{\rm Clos}_{\frak v_T^{P_{\epsilon}}}\left(y_i\frac{\partial \left(\frak P_{\text{\rm para}}\vert_{s=\frak s}\right)}{\partial y_i}:i=1,\dots,\right)}
\end{equation}
in the same way as the proof of Lemma \ref{underRcomplete}. 
In fact, $\frak P_{\rm para}\vert_{s=\frak s} \equiv \sum(\alpha\frak s c_i+(1-\frak s)) z_i \mod T^{\epsilon} \Lambda[y,y^{-1}]^{P_{\epsilon}}_+$
and $\alpha\frak s c_i+(1-\frak s) \notin \Lambda_+$ if $\frak s \in B_+^c$.
(In case $\frak s \in \Lambda \setminus \Lambda_0$ we have
$\frak P_{\rm para}\vert_{s=\frak s} \equiv T^{-\rho}\sum(\alpha c_i  - 1) z_i \mod T^{-\rho}\Lambda 
[y,y^{-1}]_+^{P_{\epsilon}}$ 
and $\alpha \ne 1/c_i$.)
\par
The rest of the proof of the first half is the same as the proof of Lemma \ref{inmompoly}
using Lemma \ref{polyapprox}.2 in place of Lemma \ref{underRcomplete}.
(See (\ref{PleadingtermR}).)
\par
Let $(\frak y_1,\dots,\frak y_n,\frak s)\in \vert \frak X \vert$,
$\frak s  \in B_+^c$.
Then $y_i \mapsto \frak y_i$ and $s \mapsto \frak s$ define a
continuous homomorphism from $(\Lambda[y,y^{-1},s],\frak v_T^{P_{\epsilon}})$.
Therefore
$(\frak y_1,\dots,\frak y_n,\frak s)\in \vert Spec(\underline{\frak R}) \vert$.
Let $Spec(\Lambda[s,y,y^{-1}]/\frak Q_j)$ be an irreducible component
containing $(\frak y_1,\dots,\frak y_n,\frak s)$.
Since $\frak s \in B_+^c$, $j\in \frak J$.
Therefore $(\frak y_1,\dots,\frak y_n,\frak s)\in \vert \frak Y \vert$
as required.
\end{proof}
We define a finite set $\mathcal B_1 \subset B_+$ as follows.
Let $i\in {\frak I}_{+} \setminus \frak I$. Then,
$$
\vert Spec(\Lambda[s,y,y^{-1}]/\frak q_i)\vert \cap \pi^{-1}(B^c_+) \cap (\vec{\frak v}_T)^{-1}(P_{\epsilon}) =
\emptyset.
$$
Therefore, by Lemma \ref{inmopoly2},
$$
\vert Spec(\Lambda[s,y,y^{-1}]/\overline{\frak q}_i)\vert \subset \pi^{-1}(B_+).
$$
Since
\begin{equation}\label{piimageclosure}
\pi\left(
\vert Spec(\Lambda[s,y,y^{-1}]/\overline{\frak q}_i)\vert
\right)
\end{equation}
is a constructible set contained in $B_+$ by Chevalley's theorem
(See  \cite[page 42]{Ma70}), it is easy to see that (\ref{piimageclosure})
is a finite set.
\par
In the same way we find that if $j\in \frak J_+ \setminus \frak J$
then
$
\pi\left(
\vert Spec(\Lambda[s,y,y^{-1}]/\frak Q_j)\vert
\right)
$
is a finite subset of $B_+$. 
\par
We define
\begin{equation}
\mathcal B_1 = \bigcup_{i\in \frak I_+\setminus \frak I}\pi\left(
\vert Spec(\Lambda[s,y,y^{-1}]/\overline{\frak q}_i)\vert
\right)
\cup
\bigcup_{j\in \frak J_+ \setminus \frak J}\pi\left(
\vert Spec(\Lambda[s,y,y^{-1}]/\frak Q_j)\vert
\right).
\end{equation}
Hereafter we denote by $A_{\mathcal B_1}$ the
localization of $A$ by the multiplicative set generated by
$\{s-\frak s \mid \frak s \in\mathcal B_1\}$.
\begin{cor}\label{I0Icoro2}
We have
$$
(\overline I_0)_{\mathcal B_1} = (\overline I)_{\mathcal B_1} = (J_0)_{\mathcal B_1},
$$
where $\overline I_0$ is the closure of $I_0$.
\end{cor}
\begin{proof}
The proof is similar to the proof of Corollary \ref{I0Icoro}.
If $j \in \frak J_+ \setminus \frak J$, then
$1 \in (\frak Q_j)_{\mathcal B_1}$.
(This is  because $\vert Spec(\Lambda[s,y,y^{-1}]
/\frak Q_j)\vert \cap \pi^{-1}(\Lambda \setminus \mathcal B_1) =
\emptyset$.)
This implies
$(\overline I)_{\mathcal B_1} = (J_0)_{\mathcal B_1}$.
If $i\in \frak I_+\setminus \frak I$, then
$1 \in (\overline{\frak q}_i)_{\mathcal B_1}$.
This implies $(\overline I_0)_{\mathcal B_1} = (\overline I)_{\mathcal B_1}$.
\end{proof}
\begin{prop}\label{componentmajiwari2}
For each $i \in \frak I$ there
exists $j\in \frak J$ such that $\vert \frak X_i \vert \subseteq \vert\frak Y_j  \vert$.  \end{prop}
\begin{proof}
Let $\frak Z$ be one of the
irreducible components of  $\frak X_i \cap \frak Y$.
We may assume
$$
\vert\frak Z\vert \cap (\vec{\frak v}_T)^{-1}(P_{\epsilon}) \cap \pi^{-1}(B_+^c) \ne \emptyset.
$$
In fact,
\begin{equation}\label{kouhoy}
\vert\frak X_i\vert \cap (\vec{\frak v}_T)^{-1}(P_{\epsilon}) \cap \pi^{-1}(B_+^c) \ne \emptyset,
\end{equation}
by definition. Any point on (\ref{kouhoy}) is also contained in $\frak Y$ by Lemma \ref{inmopoly2}.
Thus we take $\frak Z$ which contains it.
\par
We will prove the coincidence of geometric points in $\frak X_i$ and $\frak Z$.
To prove this we prove an analogue of Lemma \ref{formalization}
with an additional condition that
the point corresponding to $\overline{\frak f}$ is
contained in $\pi^{-1}(B_+^c)$.
We can prove it in the following way.
Assume $\vert\frak Z  \vert \ne \vert \frak X_i \vert$.
We consider the set of points
$\frak y \in \vert \frak Z\vert$
such that there exists a 1-dimensional irreducible and reduced scheme $\frak C \subset (\frak
X_i)_{\rm red}$, which is
a subscheme of $Spec(\frak R)$,
with $\dim (\frak C \cap \frak Z_{\rm red}) = 0$ and $\frak y \in \vert\frak C \cap \frak Z_{\rm red}\vert$.
We can easily see that this set is Zariski open in
$\vert \frak Z\vert$.
Therefore there exists such $\frak y$ in  $\pi^{-1}(B_+^c)$.
The rest of the proof is the same as the proof of Lemma \ref{formalization}.
\par
Now we observe that for $\frak s = \pi(\frak y)$ where $\frak y$ corresponds to
$\overline{\frak f}$ we have
$$
T^{-\epsilon} \frak P_{\frak s} \equiv \sum_{j=1}^m (\alpha \frak sc_i + (1-\frak s))z_i^{\epsilon}
\mod \Lambda[y,y^{-1}]_+^{P_{\epsilon}}
$$
and $\alpha \frak sc_i + (1-\frak s) \notin \Lambda_+$.
We can use this fact and proceed in the same way as in the proof of Proposition \ref{componentmajiwari}
to prove Proposition \ref{componentmajiwari2}.
\end{proof}
\begin{lem}\label{fiberdim02}
If $i \in \frak I_0$, then for any $\frak s \in B_+^c$ we have:
$$
{\rm rank}_{\Lambda} \frac{\Lambda [s,y,y^{-1}]}{\frak q_i  (s-\frak s)} < \infty.
$$
\end{lem}
\begin{proof}
By Proposition \ref{componentmajiwari2} it suffices to show that
$\frak Y  \cap \pi^{-1}(\frak s)$ is zero dimensional.
For this purpose it suffices to show that
$$
\dim_{\Lambda} {\rm Jac}(T^{-\epsilon}\frak P_{\frak s}) < \infty.
$$
This follows from  (\ref{PleadingtermR}), Theorem \ref{surj} and Theorem \ref{versality}.2.
\end{proof}
\begin{cor}\label{dim11}
If $i \in \frak I_0$, then $\dim \frak X_i \le 1$.
\end{cor}
This is immediate from Lemma \ref{fiberdim02}.
\begin{prop}\label{compactnessfiber2}
There exists a finite subset $\mathcal B_2  \supset \mathcal B_1$ of
$B_+$ such that
$$
\frak X \cap \pi^{-1}(\Lambda \setminus \mathcal B_2) \to \Lambda \setminus \mathcal B_2
$$
is projective.
\end{prop}
\begin{proof}
We consider the scheme $\prod_{i=1}^n \P^1_{y_i} \times Spec(\Lambda[s])$.
Let $\overline{\frak X}$ be the closure of $\frak X$ in
$\prod_{i=1}^n \P^1_{y_i} \times Spec(\Lambda[s])$.
We put
$$
D = \left.\left\{ (y_1,\dots:y_n,s) \in
\left\vert \prod_{i=1}^n \P^1_{y_i} \times Spec(\Lambda[s])\right\vert
~\right\vert~ \exists i ~ y_i \in \{0,\infty\}\right\}.
$$
By Corollary \ref{dim11} the intersection
$\vert\overline{\frak X}\vert \cap D$ is a finite set.
We put
$$
\mathcal B _2= \pi (\vert\overline{\frak X}\vert \cap D)   \cup \mathcal B_1 \subset \Lambda.
$$
We can show that $\mathcal B_2 \subset  B_+$ in the same way as the
proof of Proposition \ref{compatnessfiber0}.
\end{proof}
\begin{prop}\label{dimXflat}
There exists a finite subset $\mathcal B_3 \supset \mathcal B_2$ of
$B_+$ such that the following holds.
We denote $\mathcal B_3^c = \Lambda \setminus \mathcal B_3$.
\begin{enumerate}
\item $\dim \frak X = 1$.
\item $\pi: \frak X  \cap \pi^{-1}(\mathcal B_3^c)\to \mathcal B_3^c$ is flat.
\item Let $\frak m_{(\frak y,\frak s)}$ be the
maximal ideal corresponding to $(\frak y,\frak s)\in \vert\frak X\vert \cap \pi^{-1}(\mathcal B_3^c)$.
Then  $y_1\frac{\partial \frak P_{\rm para}}{\partial y_1},\dots,
y_n\frac{\partial \frak P_{\rm para}}{\partial y_n}$  is a regular sequence
in $\Lambda[s,y,y^{-1}]_{\frak m_{(\frak y,\frak s)}}$.
\end{enumerate}
\end{prop}
\begin{proof}
We put
$$
\vert\frak X' \vert
= \bigcup_{i \in \frak I_+\setminus \frak I}
\left\vert
Spec\left(
\frac{\Lambda[s,y,y^{-1}]}{\frak q_i}
\right)\right\vert.
$$
Since $\dim \frak X \le 1$,
the intersection $\vert\frak X' \vert \cap \vert\frak X \vert$ is a finite set.
We put
$$
\mathcal B_3
= \mathcal B_2 \cup \pi(\vert\frak X' \vert \cap \vert\frak X \vert).
$$
Then we can prove Proposition \ref{dimXflat}.3 in the same way as the proof of Proposition \ref{dimX}.3 as follows.
Let $ (\frak y,\frak s) \in \vert\frak X\vert \cap \pi^{-1}(\mathcal B_3^c)$,
and $\frak m_{ (\frak y,\frak s)}$ the corresponding ideal of $\Lambda[s,y,y^{-1}]$.
We find that
$I\Lambda[s,y,y^{-1}]_{\frak m_{\frak x}} = I_0\Lambda[s,y,y^{-1}]_{\frak m_{\frak x}}$.
Therefore Lemma \ref{fiberdim02} implies that
 $y_1\frac{\partial \frak P_{\rm para}}{\partial y_1},\dots,
y_n\frac{\partial \frak P_{\rm para}}{\partial y_n}, s-\frak s$
is a regular sequence of $\frak m_{ (\frak y,\frak s)}$.
Proposition \ref{dimXflat}.3 follows.
\par
We can then prove Proposition \ref{dimXflat}.2 and 1 in the same way as  Proposition \ref{dimX}.
\end{proof}
\begin{prop}\label{flattsuika2}
For any $i\in \frak I$ we have the following.
\begin{enumerate}
\item $\dim \frak X_i = 1$.
\item $\pi : \frak X_i  \cap \pi^{-1}(\mathcal B_3^c) \to \mathcal B_3^c
$
is flat.
\end{enumerate}
\end{prop}
\begin{proof}
Proposition \ref{flattsuika2} follows from Proposition \ref{dimXflat},
Lemmata \ref{primalyremainlocal}, \ref{torsion}
and Theorem \ref{unmix}
in the same way as the proof of Proposition \ref{flattsuika}.
\end{proof}
\begin{lem}\label{XYgepoint2}
In the situation of Proposition \ref{componentmajiwari2}, we have
$\vert \frak X_i\vert =  \vert \frak Y_j \vert$.
\end{lem}
The proof is by the same argument as in Lemma \ref{XYgepoint}.
For each $i\in \frak I$ we take $ j(i)$ such that
$\frak X_i$ and $\frak Y_{j(i)}$ are as in Proposition \ref{componentmajiwari2}.
\begin{prop}\label{coincidenceasschemes2}
We have $\frak X_i = \frak Y_{j(i)}$, i.e., $\frak q_i = \frak Q_{j(i)}$.
\end{prop}
\begin{proof}
The proof is similar to the proof of Proposition \ref{coincidenceasschemes}.
\par
We put
$$
B_1=  \bigcup_{j  \in \frak J_+ \setminus \{ j(i) \mid i \in \frak I\}} \pi(\vert \frak Y_j
\vert)
$$
and
$$
B_2 = \bigcup_{i_1,\ne i_2; i_1,i_2 \in \frak I}  \pi(\vert \frak X_{i_1}
\vert \cap \vert \frak X_{i_j} \vert).
$$
They are finite sets. We take
$c \in \Lambda \setminus \pi(B_1\cup B_2) \setminus B_+$.
\par
We take $\frak x = (c,\frak y) \in \vert \frak X_i\vert \cap \pi^{-1}(c)$.
Let $\frak m_{\frak y}$ be the maximal ideal of $\Lambda[y,y^{-1}]$ corresponding to $\frak y$.
By the same argument as in Proposition \ref{coincidenceasschemes} we can
show
$$
(\frak R/(s-c))_{\frak m_y} \cong (\underline{\frak R}/(s-c))_{\frak m_y}.
$$
We can now apply the proof of Proposition \ref{coincidenceasschemes} to obtain
\begin{equation}\label{localizationcoince}
\frak R_{\frak m_{\frak x}} = \underline{\frak R}_{\frak m_{\frak x}}.
\end{equation}
Here $\frak m_{\frak x}$ is the maximal ideal of $\Lambda[s,y,y^{-1}]$ corresponding to $\frak x$.
\par
In fact, we note that (\ref{localizationcoince}) is an equality of the ideals which are
localized at $\frak x$.
Also the argument of this equality in the proof of  Proposition \ref{coincidenceasschemes}
works only on those local rings and their quotients.
Therefore the proof goes without change, using Propositions \ref{dimXflat}
and \ref{flattsuika2}.
\par
Once (\ref{localizationcoince}) is proved, we have
$\frak q_i = \frak Q_{j(i)}$ by Lemma \ref{primalyremainlocal}.
Proposition \ref{coincidenceasschemes2} follows.
\end{proof}
\begin{cor}\label{I0andJ2}
$\frak X = \frak Y$.
\end{cor}
The proof is the same as one of Corollary \ref{I0andJ} using Proposition \ref{coincidenceasschemes2}.
Now we obtain
$$\aligned
\dim_{\Lambda} \text{\rm PJac}(\frak{P}')_{\Lambda}
&= \dim_{\Lambda} \text{\rm Jac}(\frak{P}')_{\Lambda}
= \dim_{\Lambda} H^0(\pi^{-1}(1);\frak O_{\pi^{-1}(1)})
\\
&= \dim H^0_{\Lambda}(\pi^{-1}(0);\frak O_{\pi^{-1}(0)})
= \dim_{\Lambda} \text{\rm PJac}(\frak{P})_{\Lambda}
= \dim_{\Lambda} \text{\rm Jac}(\frak{P})_{\Lambda}.
\endaligned$$
Hence the proof of Proposition \ref{rankequalgeneral}.
\end{proof}
\begin{rem}
Comparing with the argument for the case of
$\frak b \in \mathcal {A}(\Lambda_+)$,
we replace $P$ by $P_{\epsilon}$ in the above argument.
In fact, in our situation, there may be a solution of $y_k\frac{\partial \frak P_{\rm para}}{\partial y_k} = 0$
whose valuation lies
on $\partial P$.
Moreover
$
\frak{PO}_{\frak b} \equiv z_1 + \cdots + z_m \mod
T^{\lambda_1}\Lambda_{0}\langle\!\langle s,y,y(u)^{-1}\rangle\!\rangle
$
may not be true.
(See  Example \ref{boundarycrit}.)
By this reason we prove algebraization of $\frak{PO}_{\frak b}$
by a coordinate transformation converging on $P_{\e}$ only.
We can easily prove that there exists a neighborhood $U$
of $\partial  P$ such that there exists no solution of $y_k\frac{\partial \frak P_{\rm para}}{\partial y_k}= 0$
whose valuation lies on $U \cap \text{\rm Int}P$.
In fact, if $u \in U\cap \text{\rm Int}P$, then
$\frak P_{\rm para} \equiv z_1 + \cdots + z_m \mod T^{\lambda_1}\Lambda_{0}[s,y(u),y(u)^{-1}]$.
We used this fact to find $P_{\epsilon}$.
\end{rem}

\begin{exm}\label{boundarycrit}
Let $X = \C P^2$ and $D = \{(0,0)\} \subset \C^2 \subset  \C P^2$.
We put $\frak b = c\text{\rm PD}[D]$, $c \in \C\setminus \{0\}$.  Then we have
$$
\frak{PO}_{\frak b} \equiv y_1 + y_2 + cy_1y_2 + \sum_{k_1,k_2\ge 0 \atop k_1+k_2 \ge 3}c_{k_1k_2}y^{k_1}_1y^{k_2}_2 
\mod T\Lambda\langle\!\langle y,y^{-1}\rangle\!\rangle^{\overset{\circ}P}_0
$$
where $c_{k_1,k_2} \in \R$ and $c_{k_1,k_2} = c_{k_2,k_1}$.
We prove that one of the following holds:
\begin{enumerate}
\item
$\frak{PO}_{\frak b} \notin \Lambda\langle\!\langle y,y^{-1}\rangle\!\rangle^{P}_0$.
\item
$$
\text{\rm rank}_{\Lambda} \frac{\Lambda\langle\!\langle y,y^{-1}\rangle\!\rangle^P }{\text{\rm Clos}_{\frak v_T^P}\left(y_1\frac{\partial \frak{PO}_{\frak b}}{\partial y_1},
y_2\frac{\partial \frak{PO}_{\frak b}}{\partial y_2} \right)}
\ne  3 = \text{\rm rank}_{\Q} H(\C P^2;\Q).
$$
\end{enumerate}
In fact suppose $\frak{PO}_{\frak b} \in \Lambda\langle\!\langle y,y^{-1}\rangle\!\rangle^{P}_0$.
Then $\sum_{k_1,k_2\ge 0 \atop k_1+k_2 \ge 3}c_{k_1k_2}y^{k_1}_1y^{k_2}_2$ is necessary a polynomial.
The equation
$$
0 = 1 + cy  + \sum_{k_1,k_2\ge 0 \atop k_1+k_2 \ge 3}k_1c_{k_1k_2}y^{k_1+k_2-1}
$$
has a root $\frak y \in \C \setminus \{0\}$. There exists 
$$
y_1 = y_2 = y(T) \equiv \frak y \mod T
$$
that is a cricical point of $\frak{PO}_{\frak b}$. Therefore $y_i \mapsto y(T)$ defines a homomorphism
$$
\frac{\Lambda\langle\!\langle y,y^{-1}\rangle\!\rangle^P}{\text{\rm Clos}_{\frak v_T^P}\left(y_1\frac{\partial \frak{PO}_{\frak b}}{\partial y_1},
y_2\frac{\partial \frak{PO}_{\frak b}}{\partial y_2} \right)}
\to \Lambda.
$$
It follows that
$$
\text{\rm rank}_{\Lambda} \frac{\Lambda\langle\!\langle y,y^{-1}\rangle\!\rangle^P }{\text{\rm Clos}_{\frak v_T^P}\left(y_1\frac{\partial \frak{PO}_{\frak b}}{\partial y_1},
y_2\frac{\partial \frak{PO}_{\frak b}}{\partial y_2} \right)}
\ge 4 > 3 = \text{\rm rank}_{\Q} H(\C P^2;\Q).
$$
On the other hand,
$$
\text{\rm rank}_{\Lambda} \frac{\Lambda\langle\!\langle y,y^{-1}\rangle\!\rangle^{\overset{\circ}{P}} }{\text{\rm Clos}_{d_{\overset{\circ}{P}}}\left(y_1\frac{\partial \frak{PO}_{\frak b}}{\partial y_1},
y_2\frac{\partial \frak{PO}_{\frak b}}{\partial y_2} \right)}
= 3.
$$
In fact, $\sum_{k=1}^{\infty} y_1^k$ converges in $d^{\overset{\circ}{P}}$-topology. Therefore
$y_i \mapsto  y(T)$ does not define a homomorphism
$$
\frac{\Lambda\langle\!\langle y,y^{-1}\rangle\!\rangle^{\overset{\circ}{P}} }{\text{\rm Clos}_{d_{\overset{\circ}{P}}}\left(y_1\frac{\partial \frak{PO}_{\frak b}}{\partial y_1},
y_2\frac{\partial \frak{PO}_{\frak b}}{\partial y_2} \right)}
\to \Lambda.
$$
\end{exm}
\begin{rem}
The authors thank H. Iritani, a discussion with him clarifies an error in a previous version 
of this paper related to this example.
\end{rem}
\par

\section{The Chern class $c_1$ and critical values of $\frak{PO}_{\frak b}$}
\label{sec:c1}
In this section we prove Theorem \ref{c1iscrit}.
We regard $\frak{PO}_{\frak b}$ as an element in $\Lambda_0\langle\!\langle y,y^{-1}
\rangle\!\rangle$.
\begin{prop}\label{wherec1goes}
If we assume $\frak b \in H^2(X;\Lambda_0)$, then we have
$$
\frak{ks}_{\frak b}(c_1(X)) \equiv \frak{PO}_{\frak b}
 \mod
\left( y_i\frac{\partial \frak{\frak PO}_{\frak b}}
{\partial y_i}
: i = 1,\dots,n\right).
$$
\end{prop}
\begin{proof} We decompose $b = b_0 + x \in H^1(L(u);\Lambda_0)$ where $\frak v_T(b_0) = 0$ and $\frak v_T(x) > 0$.
Then, by definition, we have
\begin{eqnarray}\label{POformula2}
 & & \frak{PO}_{\frak b}(y_1,\dots,y_n) \text{\rm PD}[L(u)] \nonumber\\
&=&\sum_{k=0}^{\infty}\sum_{\beta} T^{\omega\cap \beta/2\pi} \frak m_{k,\beta}^{\frak b,\rho}(b,\dots,b)
\nonumber\\
& = & \sum_{k=0}^{\infty} \sum_{\beta} \exp(\frak b \cap \beta) T^{\omega\cap \beta/2\pi}
\rho^{b_0}(\partial \beta) \frak m_{k,\beta}(x,\dots,x).
\end{eqnarray}
Here we use the assumption $\frak b \in H^2(X;\Lambda_0)$.
It is well known (see \cite[Lemma in p. 109]{fulton}) that the first
Chern class $c_1(X)$ can be expressed as
$$
c_1(X) = \sum_{i=1}^m \text{\rm PD}[D_{\text{\bf p}_i}].
$$
By a dimension counting argument using the choice $\frak b \in H^2(X;\Lambda_0)$, we derive that only
the classes $\beta \in H_2(X,L(u);\Z)$ in (\ref{POformula2}) satisfying
$\mu(\beta) =2$ nontrivially contribute in the sum of the right hand side.
Therefore, by the Maslov index formula and the classification theorem
of holomorphic disks stated in \cite[p.781]{cho-oh}, we obtain
$c_1(X) \cap \beta = 1$.

Regard $c_1(X)$ as a constant vector field on $H(X;\Lambda_0)$. Then,
by the definition of $\frak{ks}_{\frak b}$ given in Definition \ref{defn:KSmap},
we have
\begin{equation}\label{eq:ksc1X}
\frak{ks}_{\frak b}(c_1(X)) = c_1(X)[\frak{PO}_{\frak b}]
\mod
\left( y_i\frac{\partial \frak{\frak PO}_{\frak b}}
{\partial y_i}
: i = 1,\dots,n\right)
\end{equation}
where $c_1(X)[\frak{PO}_{\frak b}]$ denotes the directional derivative of
the potential function $\frak{PO}_{\frak b}$ with respect to the vector field $c_1(X)$.
Noting that only the term $\exp(\frak b \cap \beta)$ depends on the variable $\frak b$ and
$$
c_1(X)[\exp(\frak b \cap \beta)] = \sum_{i=1}^m \text{\rm PD}[D_{\text{\bf p}_i}][\exp(\frak b \cap \beta)]
= \left(\sum_{i=1}^m \text{\rm PD}[D_{\text{\bf p}_i}][\frak b \cap \beta] \right)\exp(\frak b \cap \beta).
$$
A straightforward computation, or the formula
$$
(\vec v[x_1], \dots, \vec v[x_n]) = (v_1, \dots, v_n) \quad \mbox{for }  \vec v = (v_1, \dots, v_n),
$$
gives rise to
$$
\text{\rm PD}[D_{\text{\bf p}_i}] [\frak b \cap \beta] = [D_{\text{\bf p}_i}] \cap \beta
$$
and hence
$$
c_1(X)[\exp(\frak b \cap \beta)] = (c_1(X) \cap \beta) \exp(\frak b \cap \beta) = \exp(\frak b \cap \beta).
$$
Substituting this into \eqref{eq:ksc1X} and \eqref{POformula2}, the proposition follows.
\end{proof}
Theorem \ref{c1iscrit} follows from Proposition \ref{wherec1goes} and
Theorem \ref{mtintro}.1 by standard commutative algebra.
\qed
\par
\section{Hirzebruch surface $F_2$: an example}
\label{sec:hirze}
In this section, we study an example using Theorem \ref{c1iscrit}.
Consider the Hirzebruch surface $F_2(\alpha)$ whose moment polytope is given by
\begin{equation}
P_{\alpha} =
\{
(u_1,u_2)  \mid
0 < u_1,u_2, 
u_1 + 2u_1 < 2,
 u_2 < 1-\alpha
\}.
\end{equation}
We put
$$
\aligned
&\partial_1P = \{(u_1,u_2) \in P \mid u_1 = 0\},
\quad
\partial_2P = \{(u_1,u_2) \in P \mid u_2 = 0\}, \\
&\partial_3P = \{(u_1,u_2) \in P \mid u_1+2u_2 = 2\},\quad
\partial_4P = \{(u_1,u_2) \in P \mid u_2 =1-\alpha\},
\endaligned
$$
and $D_i = \pi^{-1}(\partial_iP)$, $i=1,2,3,4$.
We put
\begin{equation}
\frak b = \sum_{i=1}^4 w_i \text{\rm PD}([D_i]) \in H^2(F_2(\alpha);\Lambda_0).
\end{equation}
\begin{thm}\label{F2them}
The potential function of $F_2(\alpha)$ is given by
\begin{equation}\label{F2poten}
\aligned
\frak{PO}_{\frak b}^u
=
T^{u_1}e^{w_1}y_1+ T^{u_2}e^{w_2}y_2
&+ T^{2-u_1-2u_2}e^{w_3}y_1^{-1}y_2^{-2} \\
&+ T^{1-\alpha-u_2}e^{w_4}(1+T^{2\alpha})y_2^{-1}.
\endaligned
\end{equation}
\end{thm}
\begin{warn}\label{warn:y}
In this section, we write the variables $y_i(u)$ which depend on $u \in {\rm Int} P$, as $y_i$ to simplify the notations. This is the same as ones used in \cite{fooo10}. 
So $y_i$ in this section is $y_i(u)$, not $y_i$ in \eqref{yy(u)identify}.  
\end{warn}
\begin{rem} 
\begin{enumerate}
\item
For $\frak b =0$, the leading order potential function,
which is defined and denoted by $\frak{PO}_0$ in \cite[(4.9)]{fooo08}  is
given by
$$
T^{u_1}y_1+ T^{u_2}y_2
+ T^{2-u_1-2u_2}y_1^{-1}y_2^{-2}
+ T^{1-\alpha-u_2}y_2^{-1}.
$$
We note that this is different from the full potential function $\frak{PO}^u_{\frak b}$
with $\frak b = \text{\bf 0}$
for the case of $F_2(\alpha)$.
\item Theorem \ref{F2them} was previously proved by D. Auroux \cite[Proposition 3.1]{Aur09} 
by a different method.
Note \cite{Aur09} only states the case $\frak b = 0$. However the case $\frak b \ne 0$ 
can also follow from the classification of the pseudoholomorphic disk given in 
 \cite{Aur09}.
\end{enumerate}
\end{rem}
Theorem \ref{F2them} is \cite[Theorem 3.1]{fooo10}. We reproduce the proof below
for reader's convenience.
\begin{proof}
We first prove the following lemma.
\begin{lem}
There exists $c \in \Lambda_+$ such that
\begin{equation}
\aligned
\frak{PO}_{\frak b}^u
=
T^{u_1}e^{w_1}y_1+ T^{u_2}e^{w_2}y_2
&+ T^{2-u_1-2u_2}e^{w_3}y_1^{-1}y_2^{-2} \\
&+ T^{1-\alpha-u_2}e^{w_4}(1+c)y_2^{-1}.
\endaligned
\end{equation}
\end{lem}
\begin{proof}
Since $\deg \frak b = 2$, only the moduli spaces
$\mathcal M_{k+1}(\beta)$ with $\mu(\beta) = 2$ contribute to
$\frak{PO}_{\frak b}^u$. (See  \cite[(6.9)]{fooo09}.)
We remark that
$
c_1(D_i) > 0
$ for $i = 1,2,3$ and $c_1(D_4) = 0$.
Suppose $\mathcal M_{k+1}(\beta) \ne \emptyset$.
Then by \cite[Theorem 11.1]{fooo08}  there exist $k_i \ge 0$ and
$\ell_j \ge 0$ with $\sum k_i > 0$ such that
$$
\beta = \sum_{i=1}^4 k_i \beta_i + \sum_{j=1}^4 \ell_j D_j.
$$
Here $\beta_i \in H_2(X,L(u);\Z)$ with $\beta_i \cap D_j = 1$ ($i=j$),  $\beta_i \cap D_j =0$ ($i\ne j$).
We remark $\mu(\beta_i) =2$.
Therefore, using $\mu(\beta) = 2$, we find that $\sum k_i = 1$ and $\ell_j = 0$ for $j\ne 4$.
Since $\beta$ must be represented by a connected genus zero bordered curve,
it is easy to see that $k_4 =1$ if $\ell_4 \ne 0$.
Thus only the classes $\beta_i$ ($i=1,2,3,4$) and $\beta_4 + \ell D_4$ contribute.
The lemma follows easily.
\end{proof}
Thus it suffices to show $c=T^{2\alpha}$.
(We remark that $2\alpha = D_4 \cap \omega/2\pi$. See \cite[Theorem 8.1]{cho-oh}  or rather its proof.)
\par
Let $S^2(A)$ be a symplectic $2$ sphere with area $A$. 
The following lemma is proved in \cite[Proposition 5.1]{fooo10}.
\begin{lem}
$F_2(\alpha)$ is symplectomorphic to $S^2(1-\alpha) \times S^2(1+\alpha)$.
\end{lem}
The moment polytope of $S^2(1-\alpha) \times S^2(1+\alpha)$ is
$$
P'(\alpha) = \{(u_1,u_2) \mid 0 < u_1 < 1-\alpha,  0 < u_2 < 1+\alpha\}.
$$
There is a unique balanced fiber $u = ((1-\alpha)/2,(1+\alpha)/2)$, whose potential function is
\begin{equation}\label{S2pot}
 T^{(1-\alpha)/2}(y_1+y_1^{-1}) +  T^{(1+\alpha)/2}(y_2+y_2^{-1}).
\end{equation}
It has 4 critical points $y_1,y_2 = \pm 1$ with the corresponding critical values given by
\begin{equation}\label{crival1}
\pm 2T^{(1-\alpha)/2}(1 \pm T^{\alpha})
\end{equation}
respectively.
\par
On the other hand, the balanced fiber of $F_2(\alpha)$ is located at
$((1+\alpha)/2,(1-\alpha)/2)$, which has the potential function
\begin{equation}
\frak{PO}^u
=
T^{(1-\alpha)/2}(y_2 + (1+c)y_2^{-1})  +
T^{(1+\alpha)/2}(y_1 +y_1^{-1}y_2^{-2}).
\end{equation}
The critical point equation for $\frak{PO}^u$ is
\begin{eqnarray}
0 &=& 1 - y_1^{-2}y_2^{-2}.  \label{crit1}\\
0 &=& 1 -2 T^{\alpha}y_1^{-1}y_2^{-3} - (1+c)y_2^{-2}.\label{crit2}
\end{eqnarray}
The equation (\ref{crit1}) implies $y_1y_2 = \pm 1$. We consider the two
cases separately.
\par
\emph{Case 1}:  $y_1y_2 = -1$. Equation (\ref{crit2}) implies $y_2^2 = 1 + c - 2T^{\alpha}$.
Then the corresponding critical values are
\begin{equation}\label{crival2}
\pm 2 T^{(1-\alpha)/2} \sqrt{1+c - 2T^{\alpha}}.
\end{equation}
\par
\emph{Case 2}:   $y_1y_2 = 1$.  (\ref{crit2}) implies $y_2^2 = 1 + c + 2T^{\alpha}$.
Then the critical values are
\begin{equation}\label{crival3}
\pm 2 T^{(1-\alpha)/2} \sqrt{1+c +2T^{\alpha}}.
\end{equation}
Thus we must have
$$
1 \pm T^{\alpha} = \sqrt{1+c \pm2T^{\alpha}}
$$
from which $c=T^{2\alpha}$ immediately follows.
The proof of Theorem \ref{F2them} is complete.
\end{proof}
We have thus calculated the potential function of $F_2(\alpha)$ using
Theorem \ref{c1iscrit}. (We use Theorem \ref{mtintro}.1 but not 2
in the proof of Theorem \ref{c1iscrit}.)
Below we check Theorem \ref{mtintro}.2 for the case of $F_2(\alpha)$.
For simplicity we consider the case $\frak b = \text{\bf 0}$.
We remark that $F_2(\alpha)$ and $S^2(1-\alpha) \times S^2(1+\alpha)$
are both nef.
Therefore it follows from Theorem \ref{cliffordZ}.3 that
the residue pairing is given by the reciprocal of the determinant of the
Hessian matrix.
\par
We first calculate the Hessian determinant of the potential function (\ref{S2pot})
of $S^2(1-\alpha) \times S^2(1+\alpha)$ at the critical points $y_1 = \pm 1$, $y_2 = \pm 1$.
(We take $x_i = \log y_i$ as the coordinate.)
It is easy to see that the Hessian matrix is diagonal with entries
$$
\left( y_1\frac{\partial}{\partial y_1}\right)^2(\frak{PO}) = T^{(1-\alpha)/2}(y_1+y_1^{-1}),
\left( y_2\frac{\partial}{\partial y_2}\right)^2(\frak{PO}) = T^{(1+\alpha)/2}(y_2+y_2^{-1}).
$$
Therefore its determinant is:
\begin{equation}\label{4T}
\pm 4 T.
\end{equation}
\par
We next consider the case of $F_2(\alpha)$.
The Hessian matrix of (\ref{F2poten}) is
$$
\left[
\begin{matrix}
T^{(1+\alpha)/2}(y_1+y_1^{-1}y_2^{-2})   &   2T^{(1+\alpha)/2} y_1^{-1}y_2^{-2} \\
2T^{(1+\alpha)/2} y_1^{-1}y_2^{-2}            &
T^{(1-\alpha)/2}(y_2+(1+T^{2\alpha})y_2^{-1}) + 4T^{(1+\alpha)/2}y_1^{-1}y_2^{-2}
\end{matrix}
\right]
$$
For the case $y_1y_2=-1$, $y_2 = \pm(1 - T^{\alpha})$, the determinant of this matrix is
$$
-\frac{2T^{(1+\alpha)/2}}{y_2} \cdot \frac{2T^{(1-\alpha)/2}(1-3T^{\alpha}+T^{2\alpha})}{y_2}
- \frac{4T^{1+\alpha}}{y_2^2} = \mp 4T.
$$
For the case $y_1y_2=1$, $y_2 = \pm(1 + T^{\alpha})$, the determinant of this matrix is
$$
\frac{2T^{(1+\alpha)/2}}{y_2} \cdot \frac{2T^{(1-\alpha)/2}(1+3T^{\alpha}+T^{2\alpha})}{y_2}
- \frac{4T^{1+\alpha}}{y_2^2} = \pm 4T.
$$
Thus they coincide with (\ref{4T}).
\par
In other words, the Poincar\'e duality pairing between the units of the factors of
the decomposition of $QH(F_2(\alpha);\Lambda) \cong \Lambda^4$ gives rise to
$1/4T, 1/4T, - 1/4T$ and $ -1/4T$ respectively.
\par
We have thus illustrated Theorem \ref{mtintro}.2 in this case.
\par

\chapter{Coincidence of pairings}
\section{Operator $\frak p$ and Poincar\'e duality}
\label{sec:operatorp}
Sections \ref{sec:operatorp} - \ref{sec:PDRes} are devoted to the
proof of  Theorem \ref{Mirmain}.2.
In Sections \ref{sec:operatorp}, \ref{sec:operatorptoric}, we
construct a homomorphism\index{$i_{\ast,\text{\rm qm},(\frak b,b,u)}$}
\begin{equation}\label{qsharpmap}
i_{\ast,\text{\rm qm},(\frak b,b,u)}: HF((L(u),\frak b,b);(L(u),\frak b,b);\Lambda_0)
\to H(X;\Lambda_0).
\end{equation}
We now explain how the map $\frak{ks}_{\frak b}$ is related
to a homomorphism\index{$i^{\ast}_{\text{\rm qm},(\frak b,b,u)}$}
\begin{equation}\label{eq:iqmbu}
i^{\ast}_{\text{\rm qm},(\frak b,b,u)}: H(X;\Lambda_0)
\to HF((L(u),\frak b,b);(L(u),\frak b,b);\Lambda_0).
\end{equation}
In case $\frak{PO}_{\frak b}$ is a
Morse function, we have a splitting
$$
\text{\rm Jac}(\frak{PO}_{\frak b})
\otimes_{\Lambda_0} \Lambda \cong \prod_{\frak y \in\text{\rm Crit}(\frak{PO}_{\frak b})}
\Lambda 1_{\frak y}.
$$
See Proposition \ref{Morsesplit} and Definition \ref{1ydef}.
We note that the proof of Proposition \ref{Morsesplit}
is completed at this stage since the proof of
Theorem \ref{Mirmain}.1 is completed in
Section \ref{sec:injgen}.
Each critical point $\frak y$ corresponds to
a pair $(u,b)$ such that
$HF((L(u),\frak b,b);(L(u),\frak b,b);\Lambda_0)
\cong H(T^n;\Lambda_0)$.
(\cite[Theorem 3.12]{fooo09}.) So we identify $\frak y$ to such
$(u,b)$ and regard $\text{\rm Crit}(\frak{PO}_{\frak b})$ as
the set of such $(u,b)$'s.
We write $1_{(u,b)}$ in place of $1_{\frak y}$.
We now decompose $\frak{ks}_{\frak b}(Q) \in \text{\rm Jac}(\frak{PO}_{\frak b})
\otimes_{\Lambda_0} \Lambda$ into
\begin{equation}
\frak{ks}_{\frak b}(Q) =
\sum_{(u,b) \in \text{\rm Crit}(\frak{PO}_{\frak b})} \frak{ks}_{\frak b;(u,b)}(Q)
1_{(u,b)}.
\end{equation}
\begin{lem}\label{lem171}
\begin{equation}\label{isharpandks}
i^{\ast}_{\text{\rm qm},(\frak b,b,u)}(Q)
=  \frak{ks}_{\frak b;(u,b)}(Q)\cdot \text{\rm PD}([L(u)]).
\end{equation}
\end{lem}
This lemma is a direct consequence of the
definition. 
See \cite[Theorem 20.23]{fooospectr} for the detail of the proof. 
We are going to construct the adjoint operator $i_{\ast,\text{\rm qm},(\frak b,b,u)}$
of $i^\ast_{\text{\rm qm}, (\frak b,b,u)}$, i.e., the operator satisfying
\begin{equation}\label{dualize}
\langle i^{\ast}_{\text{\rm qm},(\frak b,b,u)}(Q),P
\rangle_{\text{PD}_{L(u)}}
=
\langle Q,i_{\ast,\text{\rm qm},(\frak b,b,u)}(P)
\rangle_{\text{PD}_{X}}.
\end{equation}
Here $\langle \cdot,\cdot
\rangle_{\text{PD}_{L(u)}}$ and
$\langle \cdot,\cdot
\rangle_{\text{PD}_{X}}$ are Poincar\'e duality pairings on
$L(u)$ and $X$ respectively, and
$Q \in H(X;\Lambda_0)$ and
$P \in H(L(u);\Lambda_0) \cong HF((L(u),\frak b,b);(L(u),\frak b,b);\Lambda_0)$.
\par
For the proof of Theorem \ref{Mirmain}.2 we need
to realize the operator
$i_{\ast,\text{\rm qm},(\frak b,b,u)}$ geometrically.
We use the operator $\frak p$, which was introduced
in  \cite[Section 3.8]{fooobook}, \cite[Section 7.4]{fooobook2},
for this purpose.
We review it in this section.
Let $C$ be a graded free $\Lambda_{0,\text{\rm nov}}$ module.
Here $\Lambda_{0,\text{\rm nov}}$\index{$\Lambda_{0,\text{\rm nov}}$} is the
universal Novikov ring:
$$
\Lambda_{0,\text{\rm nov}} =
\left.\left\{\sum_{i=0}^{\infty}
a_iT^{\lambda_i}e^{n_i} \in \Lambda_{\text{\rm nov}}
~\right\vert~
\lambda_i \ge 0, \lim_{i\to\infty}\lambda_i = \infty,
a_i \in R\right\}.
$$
($\deg e =2$.)
We define its degree shift $C[1]$ by $C[1]^k = C^{k+1}$. We
define $B_k(C[1])$, $\widehat B(C[1])$ as in Section \ref{sec:frakqreview}.
We also consider a map $\Delta^{k-1}: BC \to (BC)^{\otimes k}$
as in Section \ref{sec:frakqreview} and use the expression:
$$
\Delta^{k-1}(\text{\bf x}) = \sum_c \text{\bf x}^{k;1}_c \otimes
\dots \otimes \text{\bf x}^{k;k}_c.
$$
\par
We define
$
\text{cyc}: B_k(C[1]) \to B_k(C[1])
$
by\index{$\text{cyc}$}
\begin{equation}\label{defcyc}
\text{cyc}(x_1\otimes\dots\otimes x_k) = (-1)^{\deg'x_k\times
(\sum_{i=1}^{k-1}\deg'x_i)}
x_k\otimes x_1\otimes\dots\otimes x_{k-1}.
\end{equation}
Here and hereafter we put
$
\deg' x = \deg x - 1.
$
We also put
$\deg' (x_1\otimes \dots \otimes x_k) =
\sum_{i=1}^k \deg'x_i$.
\begin{defn}\label{def:cyclicB}
$B^{\text{cyc}}_k(C[1])$ is the quotient of $B_k(C[1])$ by the
submodule generated by $\text{cyc}(\text{\bf x}) - \text{\bf x}$
for $\text{\bf x} \in B_k(C[1])$. We put:
\index{$\widehat B^{\text{cyc}}(C[1])$}
$$
\widehat B^{\text{cyc}}(C[1]) =\widehat{\bigoplus_{k=0}^{\infty}}
B^{\text{cyc}}_k(C[1])
$$
the completed direct sum of them.
\end{defn}
For a filtered $A_{\infty}$ algebra
$(C,\{\frak m_k\}_{k=0}^{\infty})$ we define
$
\widehat{d}: \widehat B(C[1])
\to \widehat B(C[1])
$
by
\begin{equation}\label{def:dhat}
\widehat{d}(\text{\bf x})
= \sum_c (-1)^{\deg' \text{\bf x}^{3;1}}\text{\bf x}^{3;1}_c
\otimes \frak m_{\vert\text{\bf x}^{3;2}_c\vert}(\text{\bf x}^{3;2}_c) \otimes \text{\bf x}^{3;3}_c,
\end{equation}
where we define the nonnegative integer $\vert\text{\bf x}\vert$ to be the integer $k$ such that  $\text{\bf x} \in B_{k}(C[1])$.
\par
The equation $\widehat d \circ \widehat d = 0$ holds by definition of
filtered $A_{\infty}$ algebra.
\par
We also define\footnote{
The authors thank to Dr.Ganatra who explained that one needs to
subtract the term $ \frak m_0(1) \otimes \text{\bf x}$ 
in the differential
$\delta^{\rm cyc}$ during the conference on Cyclic symmetry in symplectic geometry at
American Institute of Math. 2009 October.}
\begin{subequations}\label{hatd}
\begin{equation}\label{delhatd}
\aligned
\delta^{\rm cyc}(\text{\bf x})
= &\sum_c (-1)^{\deg' \text{\bf x}^{3;1}}\text{\bf x}^{3;1}_c
\otimes \frak m_{\vert\text{\bf x}^{3;2}_c\vert}(\text{\bf x}^{3;2}_c) \otimes \text{\bf x}^{3;3}_c \\
& + \sum_{c: \text{\bf x}_c^{3;1} \ne 1, \text{\bf x}_c^{3;3} \ne 1} (-1)^{(\deg' \text{\bf x}^{3;1}_c+\deg' \text{\bf x}^{3;2}_c)
\deg' \text{\bf x}^{3;3}_c} \\
&\qquad\qquad\qquad \frak m_{\vert\text{\bf x}^{3;1}_c\vert+\vert\text{\bf x}^{3;3}_c\vert}(\text{\bf x}^{3;3}_c \otimes  \text{\bf x}^{3;1}_c)
\otimes\text{\bf x}^{3;2}_c 
\\
&- \frak m_0(1) \otimes  \text{\bf x}.
\endaligned
\end{equation}
In other words
\begin{equation}
\aligned
& \delta^{\rm cyc}(x_1 \otimes \dots \otimes x_k) \\
& =
\sum_{1\le i \le j+1 \le k+1}(-1)^* x_1 \otimes \dots \otimes
\frak m_{j-i+1}(x_i,\dots,x_j) \otimes \dots \otimes x_k \\
&\quad + \sum_{1\le i < j \le k}
(-1)^{**}
\frak m_{i+k-j+1}(x_j,\dots,x_k,x_1,\dots,x_i)
\otimes x_{i+1} \otimes \dots \otimes x_{j-1} \\
\endaligned
\end{equation}
\end{subequations}
where
$$
* =  \deg'x_1+\dots + \deg'x_{i-1},
\quad ** = (\deg'x_j+\dots + \deg'x_k)(\deg'x_1+\dots + \deg'x_{j-1}).
$$
\begin{lem} $\delta^{\rm cyc}\circ\delta^{\rm cyc}=0$.
\end{lem}
The proof is a straightforward calculation which we omit.
\par
We can prove easily that $\delta^{\rm cyc}$ induces an operator on
$\widehat B^{\text{\rm cyc}}(C[1])$, which we denote by the same symbol. Therefore
$(\widehat B^{\text{\rm cyc}}(C[1]),\delta^{\rm cyc})$ is a chain complex.
\begin{rem}
We do not need any kind of cyclic symmetry of $\frak m$
in order to prove the well-definedness of $\delta^{\rm cyc}$ on $\widehat B^{\text{\rm cyc}}(C[1])$.
\end{rem}
\begin{rem}
In \cite[Lemma \& Definition 3.8.4]{fooobook}, the differential on $\widehat B^{\text{\rm cyc}}(C[1])$
is defined by restricting $\widehat d$.
In this article we defined $\widehat B^{\text{\rm cyc}}(C[1])$ as a quotient of $\widehat B(C[1])$.
So we need to define differential in a different way.
(See Remark \ref{quotientsub}.)
However they are isomorphic as we show in 
Remark \ref{remarkoninvandsub}.
\end{rem}
We can generalize the construction of $\widehat d$ as follows. Let $H$ be a graded free $\Lambda_{0,\text{\rm nov}}$ module.
We define  $E_k(H[2])$, $\widehat E(H[2])$ as in Section \ref{sec:frakqreview}.
\par
Using a sequence of operators
\begin{equation}\label{qeqation}
\frak q_{\ell;k}:
E_{\ell} (H[2]) \otimes B_k(H^*(L;R)[1]) \to  H^*(L;R)[1]
\end{equation}
$k,\ell = 0,1,2,\dots$ introduced in
Section \ref{sec:frakqreview},
we obtain $\widehat{\frak q}: E (H[2]) \otimes \widehat B(C[1]) \to \widehat B(C[1])$ by
\begin{equation}
{\widehat{\frak q}}_{\ell}(\text{\bf y};\text{\bf x})
= \sum_c (-1)^{\deg' \text{\bf x}^{3;1}_c(\deg'\text{\bf y}+1)}
\text{\bf x}^{3;1}_c
\otimes \frak q_{\ell;\vert\text{\bf x}^{3;2}_c\vert}(
\text{\bf y};\text{\bf x}^{3;2}_c)
\otimes \text{\bf x}^{3;3}_c.
\end{equation}
We also define $\widehat{\widehat{\frak q}}$ by
\begin{equation}
\aligned
\widehat{\widehat{\frak q}}
_{\ell}(\text{\bf y};\text{\bf x})
= &\sum_c (-1)^{\deg' \text{\bf x}_c^{3;1}(\deg'\text{\bf y}+1)}
\text{\bf x}^{3;1}_c
\otimes \frak q_{\ell;\vert\text{\bf x}^{3;2}_c\vert}(
\text{\bf y};\text{\bf x}^{3;2}_c)
\otimes \text{\bf x}^{3;3}_c \\
& +  \sum_{c: \text{\bf x}_c^{3;1} \ne 1, \text{\bf x}_c^{3;3} \ne 1} (-1)^{(\deg' \text{\bf x}_c^{3;1}+\deg' \text{\bf x}_c^{3;2})
\deg' \text{\bf x}_c^{3;3}} \\
&\qquad\qquad \frak q_{\vert\text{\bf x}^{3;1}_c\vert+\vert\text{\bf x}^{3;3}_c\vert}(
\text{\bf y};\text{\bf x}^{3;3}_c \otimes  \text{\bf x}^{3;1}_c)
\otimes\text{\bf x}^{3;2}_c \\
&- \frak q_0(\text{\bf y};1) \otimes  \text{\bf x}.
\endaligned
\end{equation}
It induces
$\widehat{\frak q}^{\text{\rm cyc}}: E (H[2]) \otimes\widehat B^{\text{\rm cyc}}(C[1]) \to \widehat B^{\text{\rm cyc}}(C[1])$.
\par
Now we consider $C = H(L;\Lambda_{0,\text{\rm nov}})$,
$H = H(X;\Lambda_{0,\text{\rm nov}})$. Here $L$
is a relatively spin Lagrangian submanifold of a
symplectic manifold $X$.
(We do {\it not} assume that $X$ is a toric manifold or $L$ is a torus in this section.)
In this case, the operator $\frak q$ is defined in 
\cite[Sections 3.8]{fooobook}, \cite[Section 7.4]{fooobook2}. It was explained in  \cite[Section 2]{fooo09}.
In fact, we have
\begin{equation}\label{POPDgeneral}
\frak{PO}_{\frak b}(b) {\rm PD}([L])= \sum_{\ell=0}^{\infty}\sum_{k=0}^{\infty}
\frac{1}{\ell !}
\frak q_{\ell;k}
(\frak b^{\otimes\ell} ; b^{\otimes k})
\end{equation}
and
$$
i^{\ast}_{\text{\rm qm},(\frak b,b)}(Q) =
(-1)^{\deg Q}\sum_{\ell=0}^{\infty}\sum_{k=0}^{\infty}
\frac{1}{\ell !}\frak q_{\ell+1;k} (
[Q\otimes \frak b^{\otimes\ell}] ; b^{\otimes k}).
$$
Here $\frak b \in H(X;\Lambda_{{\rm nov},0})$ and $b$ is a weak 
bounding cochain of $\frak b$-deformed filtered $A_{\infty}$ algebra associated with $L$.
(The latter is defined by the condition that the right hand side of (\ref{POPDgeneral}) 
is proportional to the unit.)
\par
We define the operator
$
GW_{0,\ell+1}: E_{\ell+1}(H(X;\Lambda_0)) \to H(X;\Lambda_0)
$
by
$$
\langle
\frak a_0,GW_{\ell+1}(\frak a_1 \otimes \dots \otimes \frak a_{\ell+1})\rangle_{\text{\rm PD}_X}
= GW_{\ell+2}(\frak a_0\otimes \frak a_1 \otimes \dots \otimes \frak a_{\ell+1})
$$
where the right hand side is explained in Section \ref{sec:statements}.
Let
$
[x_1\otimes \dots \otimes x_k]$ be the
equivalence class of $x_1\otimes\dots\otimes x_k
\in B_k(C[1])$
in $B^{\text{\rm cyc}}_k(C[1])$.
Let $(S(X;\Lambda_0),\delta)$
be an appropriate cochain complex whose cohomology
group is the cohomology group of $X$. (In \cite{fooobook} we used
a countably generated subcomplex of singular chain complex.
In the next section we will use de Rham complex.)
\index{operator $\frak p$}
\begin{thm}[\cite{fooobook}, Theorem 3.8.9]\label{pmain}
There exists a family of the operators
$$
\frak p_{\ell;k}:
E_{\ell} (H(X;\Lambda_{0,\text{\rm nov}})[2])
\otimes B^{\text{\rm cyc}}_k(H(L;\Lambda_{0,\text{\rm nov}})[1])
\to S(X;\Lambda_{0,\text{\rm nov}})
$$
for $k,\ell = 0,1,\dots$ with the following properties.
\begin{enumerate}
\item Let $k>0$ and $\text{\bf y} \in \widehat E(H(X;\Lambda_{0,\text{\rm nov}})[2])$
and $\text{\bf x} \in \widehat B^{\text{\rm cyc}}_k(H(L;\Lambda_{0,\text{\rm nov}})[1])$. Then, we have
$$
\delta(\frak{p}( \text{\bf y}; \text{\bf x})) + \sum_{c} \frak{p}( \text{\bf y}^{(2;1)}_{c} ;
(\widehat{\frak q}^{\text{\rm cyc}} (\text{\bf y}^{(2;2)}_{c} ; \text{\bf x})
))= 0.
$$
\item
Let $1 \in \widehat B^{\text{\rm cyc}}_0(H(L;\Lambda_{0,\text{\rm nov}})[1])$ and
$\text{\bf y} \in \widehat E(H(X;\Lambda_{0,\text{\rm nov}})[2])$.
Then we have
$$\aligned
\delta(\frak{p}( \text{\bf y} ; 1)) &+
\sum_{c}  (-1)^{\deg \text{\bf
y}^{(2;1)}_{c}}\frak{p}( \text{\bf y}^{(2;1)}_{c} ;
(
\frak q (\text{\bf y}^{(2;2)}_{c} ; 1)
))\\
&+ \widetilde{GW}_{0,\ell+1}(X)(\text{\bf y}; \text{\rm PD}[L]) = 0.
\endaligned$$
Here $\widetilde{GW}_{0,\ell+1}$ is a homomorphism to $S(X;\Lambda_{0,\text{\rm nov}})$
which realizes ${GW}_{0,\ell+1}$ in the cohomology level.
\item
Let $\text{\bf e}_L = \text{\rm PD}[L]$ be the unit, that is the
Poincar\'e dual to the fundamental class of $L$. Then we have
$$
\frak{p}_{\ell;k}
(\text{\bf y};[\text{\bf e}_L\otimes x_1\otimes \dots\otimes x_{k-1}])
= 0.
$$
\item Let $\text{\bf e}_X = \text{\rm PD}[X]$ be the
Poincar\'e dual to the fundamental class of $X$. Then we have the following:
\begin{itemize}
\item If $k>2$, or $\ell \ge 1$, then
$$
\frak{p}_{\ell+1;k}
([\text{\bf e}_X\otimes y_1 \otimes \dots \otimes y_{\ell}] ; \text{\bf x}
) = 0.
$$
\item For $\ell = 1$, $k = 0$,
$$
\frak{p}_{1;0}(\text{\bf e}_X ; 1) = i_!(\text{\bf e}_L).
$$
\item For $\ell = 1$, $k = 1$,
$$
\frak{p}_{1;1}(\text{\bf e}_X ; x) = 0.
$$
\end{itemize}
\item
For $\ell \ge 1$ we have
$
\frak{p}_{k;\ell} \equiv 0 \mod \Lambda_{0,{\rm nov}}^+.$
(Here  $\Lambda_{0,{\rm nov}}^+$ is the set of all elements of $\Lambda_{0,{\rm nov}}$ such that all the
exponents of $T$ is positive.)
\end{enumerate}
\end{thm}
\begin{rem}
\begin{enumerate}
\item
\cite[Theorems 3.8.9 and Proposition 3.8.78]{fooobook} are the chain level
versions of Theorem \ref{pmain}. It is reduced to the cohomology version
in  \cite[Theorem 7.4.192]{fooobook2}, from which Theorem \ref{pmain} immediately follows.
\item In  \cite[(3.8.10.6)]{fooobook} it is claimed that, in case $\ell = 0$ and
$k=2$,
$
\frak{p}_{2;0}(1 ; [\text{\bf e}_L\otimes x])
= x.
$
This is actually an error as is also mentioned in \cite[Remark 15.7]{fooooverZ}.
Namely the right hand side should be zero. 
Then it will be a special case of Item 3 above.
(We thank a referee who pointed out an error of this formula.)
\item In \cite[Theorem 3.8.9]{fooobook}, there is a
factor $1/\ell!$ in front of the term $\widetilde{GW}_{0,\ell+1}(X)(\text{\bf y}; \text{\rm PD}[L])$.
We do not have one here since we take the quotient of $BC$ by the action of symmetric group for the definition of $EC$, instead of taking the invariant subset of $BC$.
See Remark \ref{quotientsub}.
\item
We mention some of the related works in Section \ref{sec:QMHOch}.
\end{enumerate}
\end{rem}
Let $x \in C[1]$ be an element of an even (shifted) degree.
We consider the element $e^x = \sum_{k=0}^{\infty} x^{\otimes k}
\in B(C[1])$. For $w \in C$,
we consider
$$
e^x w e^x \in B(C[1]), \quad [we^x] \in B^{\text{\rm cyc}}C.
$$
Let $\frak b \in H^2(X;\Lambda_{0,\text{\rm nov}}^+)$ and $b \in H^1(L;\Lambda_{0,\text{\rm nov}}^+)$.
We recall $e^{\frak b} = \sum \frak b^{\otimes\ell}/\ell!$.
We assume that $(\frak b,b)$ is a weak bounding cochain. Namely we assume
$
\frak q(e^{\frak b} ; e^b) = \frak{PO}_{\frak b}(b) \text{\bf e}_L,
$
where $ \frak{PO}_{\frak b}(b) \in \Lambda_{0,\text{\rm nov}}^+$.
For $x \in H(L;\Lambda_{0,\text{\rm nov}})$
we define an operator
$
\delta^{\frak b,b}_{\text{\rm can}}: H(L;\Lambda_{0,\text{\rm nov}}) \to H(L;\Lambda_{0,\text{\rm nov}})
$
of degree $1$ by
$
\delta_{\text{\rm can}}^{\frak b,b}(w) = \frak q(e^{\frak b} ; [w e^b]).
$
The basic properties of $\frak q$ (see  \cite[Theorem 3.8.32]{fooobook} or
\cite[Theorem 2.1]{fooo09}) imply $\delta_{\text{\rm can}}^{\frak b,b} \circ
\delta_{\text{\rm can}}^{\frak b,b} = 0$.
\begin{defn}\label{defi*}
We define a homomorphism
$$
i_{\#,\text{\rm qm},(\frak b,b)}
: H(L;\Lambda_{0,\text{\rm nov}}) \to S(X;\Lambda_{0,\text{\rm nov}})
$$ of degree $0$ by
\begin{equation}\label{fromptoi}
i_{\#,\text{\rm qm},(\frak b,b)}(x)
= \frak p(e^{\frak b} ; [x e^b]).
\end{equation}
\end{defn}
\begin{lem}\label{i*chainmap}
$i_{\#,\text{\rm qm},(\frak b,b)}$ is a chain map.
Namely we have
$$
\delta\circ i_{\#,\text{\rm qm},(\frak b,b)}+
i_{\#,\text{\rm qm},(\frak b,b)} \circ \delta_{\text{\rm can}}^{\frak b,b} = 0.
$$
\end{lem}
Note that the boundary operator on $H(X;\Lambda_{0,\text{\rm nov}})$ is $0$.
\begin{proof}
We observe that
$$
\Delta(e^{\frak b}) = e^{\frak b} \otimes e^{\frak b}.
$$
It follows that
$$\aligned
&\widehat{\frak q}^{\text{cyc}} (e^{\frak b}; [x e^b]) \\
&=  [\frak q(e^{\frak b}; e^b x e^b) e^b]
+ (-1)^{\deg'x} [xe^b{\frak q}^{\text{cyc}} (e^{\frak b} ; e^b ) e^b]\\
&= [\delta_{\text{\rm can}}^{\frak b,b}(x)e^b] +
(-1)^{\deg' x}\frak{PO}(\frak b,b)[xe^b \text{\bf e}_L e^b].
\endaligned$$
Therefore Theorem \ref{pmain} (1) and (3) imply
$$\aligned
-(\delta \circ  \frak{p})(e^{\frak b} ; [xe^b])
&= \frak{p}(e^{\frak b} ;
[\delta_{\text{\rm can}}^{\frak b,b}(x)e^b]) = (i_{\#,\text{\rm qm},(\frak b,b)} \circ \delta_{\text{\rm can}}^{\frak b,b})(x),
\endaligned$$
as required.
\end{proof}
We recall\index{$HF((L(u),\frak b,b);(L(u),\frak b,b);\Lambda_0)$}
$$
HF((L(u),\frak b,b);(L(u),\frak b,b);\Lambda_0)
= \frac{\text{\rm Ker} \delta_{\text{\rm can}}^{\frak b,b}}{\text{\rm Im} \delta_{\text{\rm can}}^{\frak b,b}}.
$$
Therefore $i_{\#,\text{\rm qm},(\frak b,b)}$ induces a map (\ref{qsharpmap}).
We will discuss (\ref{dualize}) in Section \ref{sec:operatorptoric}.
\begin{rem}
In this section, where we review general story from \cite{fooobook} \cite{fooobook2},
we use the universal Novikov ring $\Lambda_{0,\text{\rm nov}}$ which contains two
formal parameters $T$ and $e$. In other sections of this current paper we use $\Lambda_{0}$ which
contains only $T$. The parameter $e$ is used so that all the operators have the
well defined degree (for example $\frak m_k$ has degree $1$).
Stories over $\Lambda_{0,\text{\rm nov}}$ coefficients and over $\Lambda_0$ are basically the same
if we take enough care of the degree.
Since the ring $\Lambda_0$ behaves
better than $\Lambda_{0,\text{\rm nov}}$ from the point of view of commutative algebra, we use $\Lambda_0$ in our toric case.
\end{rem}
\par
\section{Cyclic symmetry in the toric case}
\label{sec:cyclic}
The discussion on the operator $\frak p$
in Section \ref{sec:operatorp} shows that it is essential to use a  system of perturbations of the moduli spaces
$\mathcal M_{k+1;\ell}^{\text{\rm main}}(\beta)$ that is invariant under
the cyclic permutation of the boundary marked points.
The cyclic symmetry is also important when we study the moduli spaces
of holomorphic annuli in Section \ref{sec:annuli}. 
\par
Such a perturbation is constructed in \cite{fooo091}.
In this section we adapt this construction to the present situation and explain how we can apply a similar scheme to construct the cyclically
symmetric analog $\frak q^{\frak c}$ to the operator $\frak q$. We also study
the relationship between the operators $\frak q$ and $\frak q^{\frak c}$.
\par
A review of the perturbation to construct the operator $\frak q$ is
in Section \ref{sec:frakqreview}.
We next discuss the cyclically symmetric perturbation of
$\mathcal M_{k+1;\ell}^{\text{\rm main}}(\beta;\text{\bf p})$.
Let
\begin{equation}
\text{\rm cyc} : \mathcal M_{k+1;\ell}^{\text{\rm main}}(\beta;\text{\bf p}) \to \mathcal M_{k+1;\ell}^{\text{\rm main}}
(\beta;\text{\bf p})
\end{equation}
be the map induced by a cyclic permutation of marked points that sends $(0,1, \dots, k)$ to $(1,\dots, k, 0)$.
It generates the action of $\Z_{k+1}$. We can construct a
Kuranishi structure of $\mathcal M_{k+1;\ell}^{\text{\rm main}}(\beta;\text{\bf p})$ to which the
$\Z_{k+1}$ action extends.
\begin{prop}\label{existcycpert}
There exists a system of continuous families of multisections satisfying
the properties $1$-$5$ stated in Condition \ref{perturbforq}. In addition, it also satisfies:
\par
$6.$ \quad  It is invariant under the above $\Z_{k+1}$-action.
\end{prop}
\begin{rem}
The notion of continuous family of multisections and
its application to the smooth correspondence in
de Rham theory is discussed in \cite[Section 12]{fooo09}.
In our case where there is a $T^n$ action, we can
perform the construction there so that the family of multisections
is $T^n$-equivariant. We need to include the case where the parameter space
of the family of multisections admits a $T^n$ action.
We will describe it in Sections \ref{sec:equimulticot}.
\par
Actually we can work out the whole argument without requiring $T^n$-
equivariance of the family of multisections in Proposition \ref{existcycpert}. If we perform the construction in a
$T^n$-equivariant way, construction of the canonical model\index{canonical model}
becomes much simpler. This is the reason why we use the $T^n$-equivariance
in Proposition \ref{existcycpert}.
On the contrary, in the construction of the operator $\frak q$ where we use
a single multisection, $T^n$ equivariance is used to ensure
surjectivity of the evaluation map.
This surjectivity is used in \cite[Proposition 6.6]{fooo09}  and to
construct the embedding $H^1(L(u);\Lambda_0)
\subset \mathcal M_{\text{\rm weak}}(L(u);\Lambda_0)$. In this regard,
the latter usage is more serious than that of the case of
Proposition \ref{existcycpert}.
\end{rem}
Except the statement on $T^n$ equivariance,
Proposition \ref{existcycpert} follows from
\cite[Corollary 3.1]{fooo091}. We can perform the construction
in a $T^n$-equivariant way.  The proof is given in detail in Sections \ref{sec:cyclicKura} - \ref{sec:equimulticot}.
\begin{rem}
Here is an important remark to be made.
In \cite{fooo091} one of the key idea to construct a cyclically symmetric  (family of) perturbations
is to include $\mathcal M_{0;\ell}(\beta;\text{\bf p})$ also.
In \cite{fooo08} the construction of the
$\frak q$-perturbation on
$\mathcal M_{k+1;\ell}(\beta;\text{\bf p})$ uses the fact that
the $T^n$ action there is free. This does not hold for the case $k+1=0$.
Therefore the argument of \cite{fooo08} is not enough to construct cyclically
symmetric perturbation.
It seems that it is necessary to use a {\it continuous family} of multisections
for the proof of Proposition \ref{existcycpert}.
\end{rem}
\begin{defn}
We call the system of continuous family of multisections (perturbations) in Proposition \ref{existcycpert}
the {\it $\frak c$-perturbation}\index{multisection (perturbation)!$\frak c$-multisection} and write  the corresponding
perturbed moduli space as $\mathcal M_{0;\ell}(\beta;\text{\bf p})^{\frak c}$.
\end{defn}
As we show later in this section, we can use the $\frak c$-perturbation in place
of the $\frak q$-perturbation to define an operator that
is not only cyclically symmetric but shares most of the properties
of $\frak q$.
\par
However there is one property of the $\frak q$-perturbation which is not shared
by the $\frak c$-perturbation. The next lemma is related to this point.
Let $\text{\bf p} = (\text{\bf p}(1),\dots,\text{\bf p}(\ell))$ be $T^n$-invariant cycles as in Section
\ref{sec:frakqreview}.
\begin{conds}\label{nefcond}
$X$ is nef and $\deg \text{\bf p}(i) = 2$ for all $i$.
(Note $\dim \text{\bf p}(i) = 2n - \deg  \text{\bf p}(i)$.)
\end{conds}
\begin{lem}\label{dimandempty}
Assume Condition \ref{nefcond}.
We may take a $\frak c$-perturbation so that the following holds.
$\mathcal M_{k+1;\ell}(\beta;\text{\bf p})^{\frak c}$
is empty for $k \ge 0$ if one of the following is satisfied:
\begin{enumerate}
\item $\mu(\beta) - \sum_i(2n- \dim \text{\bf p}(i) - 2) < 0$.
\item $\mu(\beta) - \sum_i (2n- \dim \text{\bf p}(i) - 2) = 0$ and $\beta \ne 0$.
\end{enumerate}
\end{lem}
\begin{proof}
Since $\deg {\bf p}(i)=2$ for all $i$, the condition 1 (resp. 2) 
is equivalent to 
$\mu(\beta)<0$ (resp. $\mu(\beta) =0$ for $\beta \ne 0$).
\par
We assume that $\mathcal M_{1;0}(\beta)$ is nonempty.
By \cite[Theorem 11.1]{fooo08}  we have
$$
\beta = k_{i_1}\beta_{i_1} + \dots + k_{i_1}\beta_{i_l} + \alpha, \quad k_{i_j} \in \Z_{\ge 0},
$$
where $i_1,\dots,i_l \in \{1,\dots,m\}$ with $\mu(\beta_{i_j})=2$,
$\alpha \in H_2(X;\Z)$ and $\alpha$ is realized by
a sum of holomorphic spheres.
Since $X$ is nef, we have $c_1(\alpha) \ge 0$.
This leads a contradiction to the condition 1 or 2.
Thus $\mathcal M_{1;\ell}(\beta;\text{\bf p}) = \emptyset$.
\par
Here we note that we do not need to perturb the moduli space. 
Hence it is also empty after
perturbation.
\end{proof}
We recall that the virtual dimension of 
$\mathcal M_{k+1;\ell}(\beta;\text{\bf p})^{\frak c}$ is 
given by 
$$
\dim \mathcal M_{k+1;\ell}(\beta;\text{\bf p})^{\frak c}
= 
n+\mu(\beta)+ k-2 -\sum_i(2n- \dim \text{\bf p}(i) - 2).
$$
So in the situation of Lemma \ref{dimandempty} we have
$$
\dim \mathcal M_{1;\ell}(\beta;\text{\bf p})^{\frak c}
\le n-2.
$$
\begin{rem}
Lemma \ref{dimandempty} is similar to \cite[Corollary 6.6]{fooo09}.
However in \cite[Corollary 6.6]{fooo09}  we do not assume 
Condition $\ref{nefcond}$. In the situation of \cite{fooo09} 
we can use a {\it single} multisection. Then 
the moduli space with negative (virtual) dimension can be perturbed and becomes an empty set. 
Hence we do not need Condition \ref{nefcond}. 
But if we  
use a {\it continuous family} of multisections to achieve 
the cyclic symmetry, we need Condition 
\ref{nefcond}.
\end{rem}
\begin{defn}\index{operator $\frak q^{\frak c}$}
We use the perturbed moduli space $\mathcal M_{k+1;\ell}(\beta;\text{\bf p})^{\frak c}$
to define
$$
\frak q^{\frak c}_{\ell;k+1;\beta} :
E_{\ell}(\mathcal A[2]) \otimes B_k(H(L(u);\R)[1]) \to H(L(u);\R)[1]
$$
by
\begin{equation}\label{qcdef}
\frak q^{\frak c}_{\ell;k;\beta}([\text{\bf p}] ; h_1\otimes\dots\otimes h_k)
= \text{\rm ev}_{0 *}(\text{\rm ev}^*(h_1\times\dots\times h_k);
\mathcal M_{k+1;\ell}^{\text{\rm main}}(\beta;\text{\bf p})^{\frak c})
\end{equation}
where
$$
(\text{\rm ev}_0,\text{\rm ev}): \mathcal M_{k+1;\ell}^{\text{\rm main}}(\beta;\text{\bf p})^{\frak c}
\to L(u) \times L(u)^k
$$
and $h_i$ are $T^n$-invariant forms.
(Here we use Convention \ref{ev0kakikataconv}.)
The right hand side of (\ref{qcdef}) is independent of the choice of the
representative $\text{\bf p}$ but depends only on the
equivalence class $[\text{\bf p}]$. This is a consequence of
Condition \ref{perturbforq}.5.
\par
Next we define
$$
\frak q_{\ell;k;\beta}^{\rho,\frak c}([\text{\bf p}] ; h_1\otimes\dots\otimes h_k) = 
\rho(\partial\beta)
\frak q^{\frak c}_{\ell;k;\beta}([\text{\bf p}] ; h_1\otimes\dots\otimes h_k) 
$$
and 
$\frak q_{\ell;k}^{\rho,\frak c} 
= \sum_{\beta}{T^{\beta \cap \omega/2\pi}}\frak q_{\ell;k;\beta}^{\rho,\frak c}.
$
\end{defn}
\begin{thm}\label{qcprop}
The operators $\frak
q_{\ell;k;\beta}^{\rho,\frak c}$ have the following properties:
\begin{enumerate}
\item
{\rm (\ref{qmaineq0})} holds. (We replace $\frak q$ by $\frak q^{\rho,\frak c}$.)
\item {\rm (\ref{unital0})} and {\rm(\ref{unital20})} hold.
(We replace $\frak q$ by $\frak q^{\rho,\frak c}$.)
\item The operator $\frak q_{\ell;k;\beta}^{\rho,\frak c}$ is cyclic, i.e., satisfies
\begin{equation}\label{formulacyclic}
\aligned
&\langle\frak q_{\ell;k;\beta}^{\rho,\frak c}(\text{\bf y};h_1,\dots,h_k),h_0\rangle_{\rm{cyc}}
\\
&=
(-1)^{\deg' h_0(\deg' h_1 + \dots + \deg' h_k)}
\langle\frak q_{\ell;k;\beta}^{\rho,\frak c}(\text{\bf y};h_0,h_1,\dots,h_{k-1}),h_k\rangle_{\rm{cyc}}.
\endaligned
\end{equation}
Here $\langle\cdot , \cdot \rangle_{\rm{cyc}}$ is defined by \eqref{PDandCYC}, see also
Subsection \ref{subsec:signpreliminary}.
\end{enumerate}
\end{thm}
\begin{proof}
The statements 1 and 2  follow from Condition $\ref{perturbforq}$.1-5 in the same way as the case of
$\frak q$. The statement 3 is a consequence of Proposition \ref{existcycpert}.6.
\end{proof}
\begin{rem}
\begin{enumerate}
\item
As we mentioned in Remark \ref{rem235}, we need to fix an energy level $E_0$ and restrict the
construction of the family of multisections for the
moduli space of maps of homology class $\beta$ with $\beta \cap \omega \le E_0$,
by the reason explained in \cite[Subsection 7.2.3]{fooobook2}.
(After that, we can use homological algebra to extend construction of
the operators to all $\beta$'s.)
\par
This process is discussed in \cite[Section 14]{fooo091}, following
\cite[Section 7.2]{fooobook2}.
In the current
$T^n$-equivariant case, we can perform the same construction
in a $T^n$-equivariant way.
We will not repeat this kind of remarks in other similar situations in the rest of this paper.
\item
The sign of (\ref{formulacyclic}) looks different from the one in \cite{fooobook}
at glance but is mathematically consistent.
We here use the sign convention of \cite{fooo091}.
See \cite[Lemma 6.1]{fooo091}  about the relationship between
two conventions.
\end{enumerate}
\end{rem}
\begin{defn}\index{$\frak b = \frak b_0 + \frak b_2 + \frak b_{\text{\rm high}}$}
Let $\frak b \in H^{even}(X;\Lambda_0)$ and $\rho : H_1(L(u);\Z) \to \C^*$ a
local system. We write
$$\frak b = \frak b_0 + \frak b_2 + \frak b_{\text{\rm high}}$$ where
$\frak b_0 \in H^{0}(X;\Lambda_0)$, $\frak b_2 \in H^{2}(X;\Lambda_0)$,
$\frak b_{\text{\rm high}} \in H^{2m}(X;\Lambda_0)$ $(m > 1)$.
We define
$
\frak m_{k;\beta}^{\frak c,\frak b,\rho} :
B_k(H(L;\R)[1]) \to H(L;\R)[1]
$
by
\begin{equation}\label{mcyclicdef}
\aligned
&\frak m_{k}^{\frak c,\frak b,\rho}(h_1,\dots,h_k) \\
&= \sum_{\beta}\sum_{\ell=0}^{\infty}
\frac{T^{\beta \cap \omega/2\pi}}{\ell!}
\rho(\partial \beta)
\exp(\frak b_2\cap \beta)
\frak q_{\ell;k;\beta}^{\frak c}(\frak b_{\text{\rm high}}^{\ell};h_1,\dots,h_k).
\endaligned\end{equation}
Let $b \in H^{odd}(L(u);\Lambda_0)$. We define\index{$b = b_0 + b_+$}
$$b = b_0 + b_+$$ where
$b_0 \in H^1(L(u);\C)$ and $b_+ \in H^1(L(u);\Lambda_+) \oplus \bigoplus_{m>0}H^{2m+1}(L(u);\Lambda_0)$.
We define $\rho :  H_1(L(u);\Z) \to \C^*$ by $\gamma \mapsto e^{\gamma \cap b_0}$.
We define
\begin{equation}\label{mcyclicdef2}
\frak m_{k}^{\frak c,\frak b,b}(h_1,\dots,h_k) =
\sum_{l_0,\dots,l_{k}}\frak m_{k+ \sum l_i;\beta}^{\frak c,\frak b,\rho}(b_+^{l_0},h_1,b_+^{l_1},\dots,b_+^{l_{k-1}},
h_k,b_+^{l_k}).
\end{equation}
\end{defn}
In the rest of this paper
we consider $T^n$-invariant differential forms unless otherwise stated
and hence we work with the canonical model.
\begin{rem}
Convergence of the series in the right hand side of (\ref{mcyclicdef2}) is proved in the same way as in the
proof of Lemma \ref{unboundedX1}.
\par
We handle degree 2 classes $\frak b_2$ in a different way from higher degree
classes $\frak b_{\text{\rm high}}$. We do so for the proof of convergence of (\ref{mcyclicdef2})
and similar other series.
\end{rem}
\begin{prop}
$(H(L(u);\Lambda_0),\{\frak m_{k}^{\frak c,\frak b,\rho}\}_{k=0}^{\infty})$ and
$(H(L(u);\Lambda_0),\{\frak m_{k}^{\frak c,\frak b,b}\}_{k=0}^{\infty})$ are
unital and cyclically symmetric filtered $A_{\infty}$ algebras.
\end{prop}
\begin{proof}
This follows from Theorem \ref{qcprop}.
\end{proof}
\begin{lem}\label{weakMCsol}
If $X$ is nef and $\frak b \in H^2(X;\Lambda_0)$, then
$$
H^1(L(u);\Lambda_+) \subset \widehat{\mathcal M}_{\text{\rm weak}}(
(H(L(u);\Lambda_0),\{\frak m_{k}^{\frak c,\frak b,\rho}\}_{k=0}^{\infty});\Lambda_+).
$$
\end{lem}
\begin{proof}
This follows from Lemma \ref{dimandempty} in the same way as  \cite[Proposition 3.1]{fooo09}
follows from  \cite[Corollary 6.6]{fooo09}.
\end{proof}
We note that without the assumption that $X$ is nef and $\frak b \in H^2(X;\Lambda_0)$
it is not yet clear at this stage whether
$(H(L(u);\Lambda_0),\{\frak m_{k}^{\frak c,\frak b,\rho}\}_{k=0}^{\infty})$ is
weakly unobstructed. We will prove this below by comparing $\frak q^{\frak c}$ with $\frak q$.
We use the bifurcation method for this purpose.
(See \cite[Subsection 7.2.14]{fooobook2}  about the cobordism method and
the bifurcation method.)
\par
Let $H = H(L(u);\Lambda_0)$ and $\overline H =  H(L(u);\R)$.
We consider the set of formal sums
\begin{equation}\label{elementdt}
a(t) + dt \wedge b(t)
\end{equation}
where $a(t) \in C^{\infty}([0,1],\overline H^k)$,
$b(t) \in C^{\infty}([0,1],\overline H^{k-1})$.
We denote the set of such expressions (\ref{elementdt}) by
$
C^{\infty}([0,1] \times \overline H)^k.
$
We will define a filtered $A_{\infty}$ structure on 
$
C^{\infty}([0,1] \times H)
$ (in Definition \ref{def:paraAinftyoperators}) 
as follows:
\par
We assume that, for each $t\in [0,1]$, we have operations:
\begin{equation}\label{param}
\frak q^{t}_{\ell;k,\beta} : E_{\ell}(\mathcal A[2]) \otimes B_k(\overline H[1]) \to \overline H[1]
\end{equation}
of degree $-\mu(\beta)+1$ and
\begin{equation}\label{parac}
\frak{Q}^t_{\ell;k,\beta} : E_{\ell}(\mathcal A[2]) \otimes B_k(\overline H[1]) \to \overline H[1]
\end{equation}
of degree $-\mu(\beta)$.
We put 
$$
\aligned
\frak q^t_{k,\beta} & =\sum_{\ell =0}^{\infty}
\frak q^t_{\ell, k,\beta}, \quad 
\frak q^t_{\beta}=\sum_{k =0}^{\infty}
\frak q^t_{k,\beta}, \\
\frak Q^t_{k,\beta} & =\sum_{\ell =0}^{\infty}
\frak Q^t_{\ell, k,\beta}, \quad 
\frak Q^t_{\beta}=\sum_{k =0}^{\infty}
\frak Q^t_{k,\beta}. \\
\endaligned
$$
\begin{defn}\label{smoothmtandct}
We say $\frak q^t_{\ell;k,\beta}$ is {\it smooth} if for
each given $y_1,\dots, y_{\ell}, x_1,\ldots,x_k$
$$
t\mapsto \frak q^t_{\ell;k,\beta}(y_1,\dots, y_{\ell};x_1,\ldots,x_k)
$$
is an element of $C^{\infty}([0,1],\overline H)$. 
The smoothness of $\frak Q^t_{\ell;k,\beta}$ is defined in the same way.
\end{defn}
\begin{defn}\label{pisotopydef}
We say $(H,\{\frak q^t_{\ell;k,\beta}\},\{\frak Q^t_{\ell;k,\beta}\})$
is a {\it pseudo-isotopy} of unital $G$-gapped filtered $A_{\infty}$ algebras 
\index{pseudo-isotopy of $A_{\infty}$ algebras}
with bulk if the following holds:
\begin{enumerate}
\item  $\frak q^t_{\ell;k,\beta}$ and $\frak Q^t_{\ell;k,\beta}$ are smooth.
\item For each (but fixed) $t$, the pair  $(H,\{\frak q^t_{k,\beta}\})$ satisfies
the conclusion of Proposition $\ref{qcprop}$.
\item  Let $\text{\bf x} \in B_k(\overline H[1])$, $\text{\bf y} \in E_{\ell}(\mathcal A[2])$.
We put
$$
\Delta^2\text{\bf x} = \sum_{c_1 \in C_1} \text{\bf x}_{c_1;1}
\otimes \text{\bf x}_{c_1;2}\otimes \text{\bf x}_{c_1;3},
\quad
\Delta_{\rm shuff}\text{\bf y} = \sum_{c_2 \in C_2} \text{\bf y}_{c_2;1}
\otimes \text{\bf y}_{c_2;2}.
$$
(Here $\Delta_{\rm shuff}$ is the shuffle coproduct on $E(\mathcal A[2])$. See Remark \ref{quotientsub}.)
Then we have
\begin{equation}\label{isotopymaineq0}
\aligned
&\frac{d}{dt} \frak q_{\beta}^t(\text{\bf y};\text{\bf x}) \\
&+ \sum_{c_1 \in C_1}\sum_{c_2 \in C_2}\sum_{\beta_1+\beta_2=\beta}
(-1)^{*}\frak Q^t_{\beta_1}(\text{\bf y}_{c_2;1};\text{\bf x}_{c_1;1},
\frak q_{\beta_2}^t(\text{\bf y}_{c_2;2};\text{\bf x}_{c_1;2}),\text{\bf x}_{c_1;3}) \\
&- \sum_{c_1 \in C_1}\sum_{c_2 \in C_2}\sum_{\beta_1+\beta_2=\beta}
(-1)^{**}\frak q^t_{\beta_1}(\text{\bf y}_{c_2;1};\text{\bf x}_{c_1;1},
\frak Q_{\beta_2}^t(\text{\bf y}_{c_2;2};\text{\bf x}_{c_1;2}),\text{\bf x}_{c_1;3})\\
&=0.
\endaligned
\end{equation}
Here the sign is given by $* = \deg' \text{\bf y}_{c_2;1}+ \deg' \text{\bf x}_{c_1;1} +
\deg' \text{\bf x}_{c_1;1}\deg' \text{\bf y}_{c_2;2}$ and
$** = \deg' \text{\bf x}_{c_1;1}\deg' \text{\bf y}_{c_2;2}$ .
\item
$\frak q_{\ell;k,\beta_0}^t$ is independent of $t$, 
and $\frak Q_{k,\beta_0}^t = 0$.
(Here $\beta_0 = 0$ in $H_2(X;L(u))$.)
\item
$\frak q_{\ell;k,\beta}^0$ satisfies (\ref{qmaineq0}).
\end{enumerate}
\end{defn}
\begin{prop}\label{qisoconstr}
There exists a pseudo-isotopy $(H,\{\frak q^t_{k,\beta}\},\{\frak Q^t_{k,\beta}\})$
of unital $G$-gapped filtered $A_{\infty}$ algebras with bulk
such that:
\begin{enumerate}
\item $\frak q_{\ell;k,\beta}^0 = \frak q_{\ell;k,\beta}^{\rho}$.
\item $\frak q_{\ell;k,\beta}^1 = \frak q_{\ell;k,\beta}^{\rho,\frak c}$.
\end{enumerate}
\end{prop}
\begin{proof}
The proof is based on the following lemma.
\begin{lem}\label{multipar}
There exists a system of continuous family of multisections $\frak{par}$
on $[0,1]\times \mathcal M_{\ell;k+1;\beta}^{\text{\rm main}}(\beta;\text{\bf p})$
with the following properties. We use $s$ to parametrize $[0,1]$.
\begin{enumerate}
\item The $[0,1]$ parametrized version of Condition \ref{perturbforq} holds.
\item At $s=0$ our family of multisections $\frak{par}$ coincides with $\frak q$
perturbation.
\item At $s=1$ our family of multisections $\frak{par}$ coincides with $\frak c$
perturbation.
\end{enumerate}
\end{lem}
\begin{proof}
This is a straightforward analogue of the construction of $\frak q$ perturbation.
\end{proof}
We use the perturbed moduli space $([0,1]\times \mathcal M_{k+1;\ell}(\beta;\text{\bf p}))^{\frak{par}}$
to define
$$
{\frak{qQ}}_{\ell;k;\beta} :
E_{\ell}(\mathcal A[2]) \otimes B_k(H(L(u);\R)[1]) \to C^{\infty}([0,1],\Lambda^{*}([0,1]) \otimes H(L(u);\R))
$$
(where $\Lambda^{*}([0,1])$ is the vector bundle of differential forms on $[0,1]$)
by
\begin{equation}\label{Qqcdef}
\aligned
& \frak{qQ}_{\ell;k;\beta}([\text{\bf p}] ; h_1\otimes\dots\otimes h_k) \\
&= (\text{\rm ev}_s,\text{\rm ev}_0)_*(\text{\rm ev}^*(h_1\times\dots\times h_k);\mathcal M_{k+1;\ell}^{\text{\rm main}}(\beta;\text{\bf p})^{\frak{par}}).
\endaligned
\end{equation}
Here
$$
(\text{\rm ev}_0,\text{\rm ev}): ([0,1]\times \mathcal M_{k+1;\ell}^{\text{\rm main}}(\beta;\text{\bf p}))^{\frak {par}}
\to L(u) \times L(u)^k
$$
are the evaluation maps at boundary marked points and
$$
\text{\rm ev}_s : ([0,1]\times \mathcal M_{k+1;\ell}^{\text{\rm main}}(\beta;\text{\bf p}))^{\frak {par}}
\to [0,1]
$$
is the projection.
\par
We next define
$$\aligned
&\frak{qQ}_{\ell;k}^{\rho}([\text{\bf p}] ; h_1\otimes\dots\otimes h_k) \\
&= \sum_{\beta\in H_2(X;L(u);\Z)}
{T^{\beta \cap \omega/2\pi}}
\rho(\partial\beta)
\frak {qQ}_{\ell;k;\beta}([\text{\bf p}]; h_1\otimes\dots\otimes h_k).
\endaligned$$
We restrict it to $T^n$-invariant forms and obtain $\frak{qQ}_{\ell;k}^{\rho,\text{\rm can}}$.
By decomposing it into the part which contains $ds$ and into the part which does not
contain $ds$, we obtain $\frak Q^t$ and $\frak q^t$ respectively which have the required properties. 
This finishes the proof of Proposition \ref{qisoconstr}.
\end{proof}
\begin{lem}\label{nefcondconclupara}
Suppose that Condition \ref{nefcond} holds.
Then
$$
\frak q^t_{\ell;k;\beta}([\text{\bf p}];h_1,\dots,h_k) = 0, \quad
\frak Q^t_{\ell;k;\beta}([\text{\bf p}];h_1,\dots,h_k) = 0
$$
if one of the following conditions is satisfied.
\begin{enumerate}
\item $\mu(\beta) - \sum_i(2n- \dim \text{\bf p}(i) - 2) < 0$.
\item $\mu(\beta) - \sum_i(2n- \dim \text{\bf p}(i) - 2) = 0$ and $\beta \ne 0$.
\end{enumerate}
\end{lem}
The proof is similar to that of Lemma \ref{dimandempty} and so omitted.
\begin{defn}\label{def:paraAinftyoperators}
Let $\frak b \in H^{even}(X;\Lambda_0)$ and $\rho : H_1(L(u);\Z) \to \C^*$. We write
$\frak b = \frak b_0 + \frak b_2 + \frak b_{\text{\rm high}}$ where
$\frak b_0 \in H^{0}(X;\Lambda_0)$, $\frak b_2 \in H^{2}(X;\Lambda_0)$,
$\frak b_{\text{\rm high}} \in H^{2m}(X;\Lambda_0)$ $(m > 1)$.
We define
\begin{equation}\label{oneparamcyclicdef}
\aligned
&\frak m_{k}^{t,\frak b,\rho}(h_1,\dots,h_k) \\
&= \sum_{\beta}\sum_{\ell=0}^{\infty}
\frac{T^{\beta \cap \omega/2\pi}}{\ell!}
\rho(\partial \beta)
\exp(\frak b_2\cap \beta)
\frak q_{\ell;k;\beta}^{t}(\frak b_{\text{\rm high}}^{\ell};h_1,\dots,h_k), \\
&\frak c_{k,\frak b}^{t,\frak b,\rho}(h_1,\dots,h_k) \\
&= \sum_{\beta}\sum_{\ell=0}^{\infty}
\frac{T^{\beta \cap \omega/2\pi}}{\ell!}
\rho(\partial \beta)
\exp(\frak b_2\cap \beta)
\frak Q_{\ell;k;\beta}^{t}(\frak b_{\text{\rm high}}^{\ell};h_1,\dots,h_k).
\endaligned\end{equation}
Let $b \in H^{odd}(L(u);\Lambda_0)$. We define $b = b_0 + b_+$ where
$b_0 \in H^1(L(u);\C)$ and $b_+ \in H^1(L(u);\Lambda_+) \oplus \bigoplus_{m>0}H^{2m+1}(L(u);\Lambda_0)$.
We define a local system $\rho :  H_1(L(u);\Z) \to \C^*$ by
$$
\rho ~:~\gamma \mapsto e^{\gamma \cap b_0}
$$
and define
\begin{equation}\label{oneparamcyclicdef2}
\aligned
\frak m_{k}^{t,\frak b, b}(h_1,\dots,h_k) &=
\sum_{l_0,\dots,l_{k}}\frak m_{k+ \sum l_i;\beta}^{t,\frak b,\rho}(b_+^{l_0},h_{1},b_+^{l_1},\dots,b_+^{l_{k-1}},
h_k,b_+^{l_k}), \\
\frak c_{k}^{t,\frak b, b}(h_1,\dots,h_k) &=
\sum_{l_0,\dots,l_{k}}\frak Q_{k+ \sum l_i;\beta}^{t,\frak b,\rho}(b_+^{l_0},h_{1}1,b_+^{l_1},\dots,b_+^{l_{k-1}},
h_k,b_+^{l_k}).
\endaligned
\end{equation}
\end{defn}
We consider $x_i(t) + dt \wedge y_i(t) = \text{\bf x}_i \in C^{\infty}([0,1],\overline C)$.
We define
$$
\widehat{\frak m}^{t,\frak b,\rho}_{k,\beta}(\text{\bf x}_1,\ldots,\text{\bf x}_k)
= x(t) + dt \wedge y(t),
$$
where
\begin{subequations}\label{combineainf}
\begin{equation}
x(t) = {\frak m}^{t,\frak b,\rho}_{k,\beta}(x_1(t),\ldots,x_k(t))
\end{equation}
\begin{equation}\label{combineainfmain}
\aligned
y(t)
=
& \frak c^{t,\frak b,\rho}_{k,\beta}
(x_1(t),\ldots,x_k(t)) \\
&-\sum_{i=1}^k (-1)^{*_i} \frak m^{t,\rho}_{k,\beta}
(x_1(t),\ldots,x_{i-1}(t),y_i(t),x_{i+1}(t),\ldots,x_k(t))
\endaligned\end{equation}
if $(k,\beta) \ne (1,\beta_0)$
($\beta_0= 0 \in H_2(X;L(u))$) and
\begin{equation}
y(t) = \frac{d}{dt} x_1(t) + \frak m_{1,\beta_0}^{t,\frak b,\rho} (y_1(t))
\end{equation}
if $(k,\beta) = (1,\beta_0)$. Here $*_i$ in (\ref{combineainfmain}) is
$*_i = \deg' x_1 +\cdots+\deg'x_{i-1}$.
\end{subequations}
We define $\widehat{\frak m}^{t,\frak b,b}_{k,\beta}$ in a similar way.
\par
We put
$$
\aligned
&(\text{\rm Eval}_{t=0})_1(x(t) + dt \wedge y(t)) = x(0),
\quad (\text{\rm Eval}_{t=1})_1(x(t) + dt \wedge y(t)) = x(1),
\\
&(\text{\rm Eval}_{t=0})_k = (\text{\rm Eval}_{t=1})_k = 0,
\quad \text{for $k\ne 1$.}
\endaligned
$$
\begin{lem}\label{hmoequi}
$(C^{\infty}([0,1] \times  H),\{\widehat{\frak m}^{t,\frak b,\rho}_{k,\beta}\})$
is a $G$ gapped unital filtered $A_{\infty}$ algebra.
\par
$\text{\rm Eval}_{t=0}$ is a strict unital filtered $A_{\infty}$ homomorphism from
$(C^{\infty}([0,1] \times  H),\{\widehat{\frak m}^{t,\frak b,\rho}_{k,\beta}\})$
to $(H,\{{\frak m}^{\frak b,\rho}_{k,\beta}\})$.
\par
$\text{\rm Eval}_{t=1}$ is a strict unital filtered $A_{\infty}$ homomorphism from
$(C^{\infty}([0,1] \times  H),\{\widehat{\frak m}^{t,\frak b,\rho}_{k,\beta}\})$
to $(H,\{{\frak m}^{\frak c,\frak b,\rho}_{k,\beta}\})$.
\end{lem}
The proof is a straightforward calculation from Proposition \ref{qisoconstr}.
\begin{cor}\label{isowithc}
$(H,\{{\frak m}^{\frak c,\frak b,\rho}_{k,\beta}\})$
is unital isomorphic to
$(H,\{{\frak m}^{\frak b,\rho}_{k,\beta}\})$.
\end{cor}
\begin{rem}
Here an `isomorphism' between two filtered $A_{\infty}$ algebras 
means a fitered $A_{\infty}$ homomorphism which has an inverse.
Note such a homomorphism is in general nonlinear.
By this reason it is called `quasi-isomorphism' sometimes.
We use the word `isomorphism' here for the consistency with \cite{fooobook}.
Note it is an isomorphism in the usual sense of category theory 
in the category whose object is a filtered $A_{\infty}$ 
algebra and whose morphism is a filtered $A_{\infty}$
homomorphism.
\end{rem}
Since $\text{\rm Eval}$ induces an isomorphism on $\frak m_{1,0}$-cohomology, this follows from Lemma \ref{hmoequi} and
\cite[Theorem 4.2.45]{fooobook}.
(Note in the our case $\frak m_{1,0} = 0$ so homotopy equivalence
has an inverse. See \cite[Proposition 5.4.5]{fooobook}.)
See also \cite[Proposition 4.1]{fooo091}. The isomorphism constructed there
from the pseudo-isotopy in Lemma \ref{hmoequi}
turns out to be the isomorphism obtained by the cobordism method as in
\cite[Subsection 4.6.1]{fooobook}  using the time-ordered-product moduli space.
\par
Hereafter we fix an isomorphism
\begin{equation}
\frak F : (H,\{{\frak m}^{\frak b,\rho}_{k,\beta}\})
\to (H,\{{\frak m}^{\frak c,\frak b,\rho}_{k,\beta}\}).
\end{equation}
We note that
\begin{equation}
\frak F \equiv \text{\rm identity} \mod \Lambda_+.
\end{equation}
Corollary \ref{isowithc} implies
\begin{equation}\label{weakchaninbyF}
\frak F_*(H^1(L(u);\Lambda_+)) \subset \widehat{\mathcal M}_{\text{\rm weak}}
((H,\{{\frak m}^{\frak c,\frak b,\rho}_{k,\beta}\});\Lambda_+)
\end{equation}
and twisting the boundary map by non-unitary local systems $\rho$
enables us to extend the domain of $\frak F$ to $H^1(L(u);\Lambda_0)$.
\begin{rem}
Since the isomorphism $\frak F$ may increase the degree, the image
$\frak F_*(H^1(L(u);\Lambda_+))$ may not necessarily be contained
in $H^1(L(u);\Lambda_+)$ again. The authors do not know whether we can take $\frak F_*$
so that its image is contained in $H^1(L(u);\Lambda_+)$.
As a matter of fact, it is in general very hard to calculate the map
$\frak m^{\frak c,\frak b,\rho}_{k,\beta}$.
\end{rem}

However, there are two particular cases for which $\frak F_*(H^1(L(u);\Lambda_+))$ is contained
in the degree one cohomology $H^1(L(u);\Lambda_0)$. We discuss these two cases in detail now.

\begin{lem}\label{H1isinMweak2}
If dimension of $L(u)$ is $2$, then
$$
H^1(L(u);\Lambda_+) \subset \widehat{\mathcal M}_{\text{\rm weak}}((H,\{{\frak m}^{\frak c,\frak b,\rho}_{k,\beta}\});\Lambda_+).
$$
\end{lem}
\begin{proof}
We note that $\frak F_*$ preserves the parity of the degree of the elements.
 For the case of dimension $2$, the elements of odd degree are necessarily of degree $1$. The lemma then follows
from (\ref{weakchaninbyF}).
\end{proof}
\begin{lem}\label{1630}
If $X$ is nef and $\frak b \in H^2(X;\Lambda_0)$, then
$$
C^{\infty}([0,1] \times H^1(L(u);\Lambda_+)) \subset
\widehat{\mathcal M}_{\text{\rm weak}}((C^{\infty}([0,1] \times  H),\{\widehat{\frak m}^{t,\frak b,\rho}_{k,\beta}\});\Lambda_+).
$$
\end{lem}
\begin{proof}
This follows from Lemma \ref{nefcondconclupara} in the same way as  \cite[Proposition 3.1]{fooo09}
follows from  \cite[Corollary 6.6]{fooo09}.
\end{proof}
\begin{cor}\label{H1identity}
If $X$ is nef and $\frak b \in H^2(X;\Lambda_0)$, then $\frak F_*$ is an identity map on $H^1(L(u);\Lambda_+)$.
In particular,
$$
H^1(L(u);\Lambda_+) \subset \widehat{\mathcal M}_{\text{\rm weak}}((H,\{{\frak m}^{\frak c,\frak b,\rho}_{k,\beta}\});\Lambda_+).
$$
\end{cor}
This is immediate from Lemma \ref{1630}.
\begin{defn}\label{1827}
Let $b = b_0 + b_+ \in H^1(L(u);\Lambda_0)$, where $b_0 \in H^1(L(u);\C)$ and
$b_+ \in H^1(L(u);\Lambda_+)$. The term $b_0$ induces the representation $\rho = \rho_{b_0}$ as before.
We put $b_+^{\frak c} = \frak F_*(b_+)$ and  $b^{\frak c} = b_0 + b_+^{\frak c}$.
\end{defn}
We next consider the homomorphism
$$
 i^{\ast}_{\text{\rm qm},(\frak b,b,u)} : H(X;\Lambda_0) \to H(L(u);\Lambda_0)
$$
in (\ref{qsharpmap}). It is defined as:
\begin{equation}\label{1822}
\aligned
& i^{\ast}_{\text{\rm qm},(\frak b,b,u)}(Q)  \\
 &= \sum_{\beta}\sum_{\ell}\sum_k
\frac{T^{\beta\cap \omega/2\pi}}{\ell!}\rho(\partial\beta)e^{\frak b_2\cap \beta}\frak q_{\ell+1;k;\beta}(Q
 \otimes \frak b_{\text{\rm high}}^{\ell};
 b_+,\dots,b_+).
\endaligned
\end{equation}
We replace $\frak q$ by $\frak q^{\frak c}$ above and $b_+$ by $b_+^{\frak c}$ and obtain:
$$
 i^{\frak c \ast}_{\text{\rm qm},(\frak b,b^{\frak c},u)} : H(X;\Lambda_0) \to H(L(u);\Lambda_0).
$$
\begin{rem}
In case $\deg b_+ = 1$, we can rewrite (\ref{1822}) to
\begin{equation}\label{1823}
\sum_{\beta}\sum_{\ell}
\frac{T^{\beta\cap \omega/2\pi}}{\ell!}\rho^b(\partial\beta)e^{\frak b_2\cap \beta}\frak q_{\ell+1;0;\beta}(Q
 \otimes \frak b_{\text{\rm high}}^{\ell};1).
\end{equation}
See Section \ref{sec:deltaisthesame}.
In the current circumstance, $b^{\frak c}_+$ may not have degree one. Because of this
 we can not replace (\ref{1822}) by (\ref{1823}) in our definition of
$ i^{\frak c \ast}_{\text{\rm qm},(\frak b,b^{\frak c},u)}$.
\end{rem}
\begin{lem}\label{Isharprelation}
If $\frak F_*(b_+) = b_+^{\frak c}$, then
\begin{equation}
\frak F_*\circ  i^{\ast}_{\text{\rm qm},(\frak b,b,u)} =  i^{\frak c \ast}_{\text{\rm qm},(\frak b,b^{\frak c},u)}.
\end{equation}
\end{lem}
\begin{proof}
We take  $\tilde b_+$ which goes to $b_+$ and $b_+^{\frak c}$ by $\text{\rm Eval}_{t=0}$
and $\text{\rm Eval}_{t=1}$ respectively.
We use $
\frak{qQ}_{\ell;k;\beta} :
E_{\ell}(\mathcal A[2]) \otimes B_k(H(L(u);\R)[1]) \to C^{\infty}([0,1],H(L(u);\Lambda_0))[1]
$
to define
$
 i^{\ast \frak{para}}_{\text{\rm qm},(\frak b,b,u)} : H(X;\Lambda_0) \to C^{\infty}([0,1] \times H(L(u);\Lambda_0))
$
in the same way.  Then we have a commutative diagram:
\begin{equation}
\begin{CD}
H(X;\Lambda_0)
@>{ i^{ \frak c\ast}_{\text{\rm qm},(\frak b,b^{\frak c},u)}}>> H(L(u);\Lambda_0) \\
\Vert & &@AA{\text{\rm Eval}_{t=1}}A \\
H(X;\Lambda_0)  @>{ i^{\ast \frak{para}}_{\text{\rm qm},(\frak b,\tilde b,u)}}>>  C^{\infty}([0,1] \times H(L(u);\Lambda_0)) \\
\Vert & &@VV{\text{\rm Eval}_{t=0}}V \\
H(X;\Lambda_0)
@>{ i^{\ast}_{\text{\rm qm},(\frak b,b,u)}}>> H(L(u);\Lambda_0)
\end{CD}
\end{equation}
By taking the cohomology of $C^{\infty}([0,1] \times H(L(u);\Lambda_0))$,
the right hand side becomes an isomorphism which is nothing but $\frak F_*$.
The lemma follows.
\end{proof}
\begin{rem}\label{rem:wallcross}
We describe the morphism $\frak F$ for the case
in which $n=2$ but $X$ is not nef, and explain how $\frak q^{\frak c}$ could be different
from $\frak q$. Namely we compare $\frak{PO}$ and $\frak{PO}^{\frak c} = \frak{PO}
\circ \frak F^{-1}_*$.
Note that 
\begin{equation}\label{POcdesribe}
\frak{PO}^{\frak c}(0;b)
= \sum_{\beta}\sum_{k\ge 0}
\rho(\partial\beta) e^{\frak b_2 \cap \beta}
\frak q_{0;k+1;\beta}^{\frak c}(1;b_+,\dots,b_+)T^{\beta \cap \omega/2\pi}.
\end{equation}
Here $b = b_0 + b_+ \in H^1(L(u);\Lambda_0)$,
$b_0 \in H^1(L(u);\C), b_+ \in H^1(L(u);\Lambda_+)$,
and $\rho$ is induced by $b_0$ as before.
\par
We consider the moduli space
$[0,1]  \times \mathcal M_{k+1;0}^{\text{\rm main}}(\beta)$ and its
family of multisections $\frak{par}$ as in Lemma \ref{multipar}.
\par
We recall that in (\ref{POcdesribe})  we take the sum over $\beta$
with $\mu(\beta) = 2$.
On the other hand,
\begin{equation}
\dim( [0,1] \times \mathcal M_{1:0}^{\text{\rm main}}(\beta))
= n+ \mu(\beta)+1-3+1 = \mu(\beta) + 1.
\end{equation}
(Here we use $n=2$.)
Therefore, in case $\mu(\beta) = 0,2$, this moduli space can be nonempty
and
$$\aligned
&\text{\rm ev}_{0 *}
\left(
([0,1] \times \mathcal M_{1:0}^{\text{\rm main}}(\beta))^{\frak{par}}
\right)
\\
&\in C^{\infty}([0,1],H^{2-\mu(\beta)}(L(u);\R)) \oplus dt \wedge C^{\infty}([0,1],H^{1-\mu(\beta)}(L(u);\R)).
\endaligned$$
We put the left hand side as
$$
c(\beta)(t) + dt \wedge e(\beta)(t).
$$
\par
Now let us consider an element $b_+ \in H^1(L(u);\Lambda_+)$ of form
$$
b_+ = \sum_{\beta: \mu(\beta) = 2} T^{\beta \cap \omega/2\pi} b_{\beta}, \quad b_{\beta} \in H^1(L(u);\C).
$$
We solve the ordinary differential equation
$$
\frac{db_{\beta}(t)}{dt} + \sum_{\beta_1,\beta_2 \atop \beta_1+\beta_2=\beta,
\mu(\beta_1) = 0}
\exp (b_{\beta_2}(t)\cap \partial b_1) e(\beta_1)(t) = 0
$$
with initial value $b_{\beta}(0) = b_{\beta}$. 
Then we can prove that
$$
\frak F_*(b_+) = \sum_{\beta: \mu(\beta) = 2} T^{\beta \cap \omega/2\pi} b_{\beta}(1).
$$
(There is a related argument by Auroux \cite{Aur07}.)
\end{rem}
\section{Operator $\frak p$ in the toric case}
\label{sec:operatorptoric}
We now start the definition of $\frak p$.
Consider the moduli space
$\mathcal M_{k;\ell+1}^{\text{\rm main}}(\beta)$ and the evaluation map
$$
\text{\rm ev} = (\text{\rm ev}^{\text{\rm int}},\text{\rm ev},\text{\rm ev}_0^{\rm int}):
\mathcal M_{k;\ell+1}^{\text{\rm main}}(\beta)
\to X^{\ell} \times L(u)^k \times X.
$$
We note that the $0$-th marked point is an {\it interior} marked
point in this case. So the target of $\text{\rm ev}^{\rm int}_0$ is $X$.
Let $\text{\bf p} = (\text{\bf p}(1),\dots,\text{\bf p}(\ell))$ be
$T^n$ invariant cycles as in Section \ref{sec:frakqreview}.
We put
$$
\mathcal M_{k;\ell+1}^{\text{\rm main}}(\beta;\text{\bf p})
= \mathcal M_{k;\ell+1}^{\text{\rm main}}(\beta)
{}_{\text{\rm ev}^{\text{\rm int}}}\times
(\text{\bf p}(1)\times \dots \times \text{\bf p}(\ell)).
$$
\begin{lem}\label{pmodulikura}
The moduli spaces
$\mathcal M_{k;\ell+1}^{\text{\rm main}}(\beta;\text{\bf p})$
have $T^n$-equivariant Kuranishi structures with boundaries and corners, 
which we call $\frak p$-Kuranishi structures.
There exists a system of families of multisections, 
which we call $\frak p$-multisections,\index{multisection (perturbation)!$\frak p$-multisection} 
on them such that:
\begin{enumerate}
\item It is $T^n$-equivariant and disk-component-wise.
(We define the notion,  `disk-component-wise' in Definition \ref{cpwkura}.)
\item It is invariant under the cyclic permutation of boundary marked points.
\item It is compatible with the forgetful map of boundary marked points.
\item It is transversal to zero.
\item $\text{\rm ev}^+_0 : \mathcal M_{k;\ell+1}^{\text{\rm main}}(\beta;
\text{\bf p})^{\frak p} \to X$
is a submersion.
\item
If $k> 0$, then its boundary
is described by the following fiber product:
$$
\mathcal M_{k_1+1;\ell_1}^{\text{\rm main}}(\beta_1;\text{\bf p}_1)^{\frak c}
{}_{\text{\rm ev}_0}\times_{\text{\rm ev}_i}
\mathcal M_{k_2;\ell_2+1}^{\text{\rm main}}(\beta_2;\text{\bf p}_2)^{\frak p}.
$$
Here $\beta_1+\beta_2 = \beta$, $k_1+k_2 = k+1$, $(\text{\bf p}_1,\text{\bf p}_2) =
\text{\rm Split}(\text{\bf p},(\mathbb L_1,\mathbb L_2))$,
$\vert\text{\bf p}_i\vert = \ell_i$, $i = 1,\dots,k_2$.
\item
If $k=0$ and $\partial \beta=0 \in H_1(L(u);\Z)$, there is another type of boundary component
$$
\mathcal M_{0;\ell+2}(X;\widetilde{\beta};\text{\bf p})^{\frak p}
{}_{\text{\rm ev}_{\ell+1}} \times_X L(u).
$$
Here $\widetilde{\beta} \in H_2(X;\Z)$ with $i_{\ast}(\widetilde{\beta})=\beta$ 
under the natural map $i_* : H_2(X;\Z) \to H_2(X,L(u);\Z)$ and 
$\mathcal M_{0;\ell+2}(X;\widetilde{\beta})$ is 
the moduli space of genus zero stable maps (without boundary) to $X$ with $\ell+2$ marked points representing 
the class $\widetilde{\beta}$ and 
$\mathcal M_{0;\ell+2}(X;\widetilde{\beta};\text{\bf p})
= \mathcal M_{0;\ell+2}(X;\widetilde{\beta}) {}_{{\rm ev}}\times_{X^{\ell}}(\text{\bf p}(1) \times \dots \times 
\text{\bf p}(\ell))$.
(See  \cite[Proposition 3.8.27]{fooobook} for the boundary of type $7$.)
\end{enumerate}
\end{lem}
The lemma can be proved in the same way as  \cite[Theorem 3.1 and Corollary 3.1]{fooo091} 
except the  statement on $T^n$ equivariance.
The detail of the proof is given in Sections \ref{sec:cyclicKura} and \ref{sec:equimulticot}.
We note that on the `bubble' components in Item 6 we always
use the $\frak c$-perturbation.
\par
We note that the boundary of type 7 corresponds to the term 
$\widetilde{GW}_{0,\ell+1}(X)(\text{\bf y}; \text{\rm PD}[L])$
in Theorem \ref{pmain}.2.
In our case of toric fiber, this term is zero in the homology 
level since $[L(u)] = 0$ in $H(X;\Z)$.
The term $\widetilde{GW}_{0,\ell+1}(X)(\text{\bf y}; \text{\rm PD}[L])$ 
appears only in Theorem \ref{pmain}.2, which is the case $k=0$, and not in Theorem \ref{pmain}.1, 
where $k>0$.
We note that the case $k=0$ is never used in this paper.
In fact, we use $\frak p$ to obtain the map $i_{\#,\text{\rm qm},(\frak b,b)}$ as in Definition \ref{defi*}, 
where the case $k=0$ is not used.
\par
To study the Poincar\'e duality pairing (\ref{dualize}), it is convenient to
use de Rham complex of $X$ as the target of the operator $\frak p$.
This leads us to use a continuous family of multisections.
We note that we have already chosen multisections on
$\mathcal M_{k+1;\ell}^{\text{\rm main}}(\beta;\text{\bf p})^{\frak c}$.
\index{operator $\frak p$}
\begin{defn}\label{def:p}
For $k > 0$, we define
$$
\frak p_{\ell;k;\beta}:
E_{\ell}(\mathcal A(\Lambda_0)[2])
\otimes B_k^{\text{\rm cyc}}(H(L(u);\Lambda_0)[1]) \to \Omega(X)
$$
by
\begin{equation}\label{deffrapbeta}
\frak p_{\ell;k;\beta}([\text{\bf p}];
[h_1\otimes \dots \otimes h_k])
=
\text{\rm ev}^+_{0*} (\text{\rm ev}^* (h_1\times \dots \times h_k);\mathcal M_{k;\ell+1}^{\text{\rm main}}(\beta;\text{\bf p})^{\frak p}).
\end{equation}
(We use Convention \ref{ev0kakikataconv}.) 
The right hand side is independent of the choice of representative $h_1\otimes \dots \otimes h_k$
but depends only on its equivalence class in $B^{\text{\rm cyc}}_k(H(L(u);\Lambda_0)[1])$
by Lemma \ref{pmodulikura}.2.
(Hereafter we do not repeat similar remark.)
We put
$$
\frak p_{\ell;k}^{(u,b)}
= \sum_{\beta\in H_2(X;L(u);\Z)}
{T^{\beta \cap \omega/2\pi}}
y(u)_1^{\partial \beta \cap \text{\bf e}_1}\cdots
y(u)_n^{\partial \beta \cap \text{\bf e}_n}
\frak p_{\ell;k;\beta}.
$$
This defines a map
$
E_{\ell}\mathcal A(\Lambda_0)[2]
\otimes B^{\text{\rm cyc}}_k(H(L(u);\Lambda_0)[1]) \to \Omega(X) \widehat{\otimes} \Lambda\langle\!\langle y(u),y(u)^{-1}
\rangle\!\rangle
$.
\par
When $\text{\bf p}$ contains $\text{\bf f}_0 = \text{\rm PD}[X] =
\text{\bf e}_{X}$, we use Theorem \ref{pmain}.4 as
a part of the
definition of $\frak p$.
\end{defn}
We fix $u\in \overset{\circ}P$. (In Sections \ref{sec:operatorptoric}-\ref{sec:PDRes}
it is unnecessary to move $u$ and regard the operations as the ones over
$\Lambda\langle\!\langle y,y^{-1}\rangle\!\rangle_0^{\overset{\circ}P}$.)
\begin{thm} The map
$\frak p$ defined above satisfies Theorem \ref{pmain}.1-5,
where we take $\frak q^{\frak c}$ instead of $\frak q$.
\end{thm}
\begin{proof}
Below we prove Theorem \ref{pmain} Items 1-5.
\par
Proof of Items 1 \& 2: These follow from Lemma \ref{pmodulikura} and Definition \ref{def:p}.
\par
Proof of Item 3: Note $\text{\bf e}_L = 1 \in H^0(L(u);\Lambda_0)$. On the
other hand, Lemma \ref{pmodulikura}.3 implies that
the following diagram commutes:
\begin{equation}
\begin{CD}
\mathcal M_{k;\ell+1}^{\text{\rm main}}(\beta;\text{\bf p})^{\frak p}
@>{\text{\rm ev}_0^+}>> X \\
@ VV{\frak{forget}}V& \hskip-0.25cm \Vert\\
\mathcal M_{k-1;\ell+1}^{\text{\rm main}}(\beta;\text{\bf p})^{\frak p}  @>{\text{\rm ev}_0^+}>>  X
\end{CD}
\end{equation}
Item 3 follows easily.
\par
Item 4 is a consequence of the definition.
\par
Item 5 follows from the fact that
the moduli space for $\beta=0$ contributes only to $\frak p_{2,0}$
which gives the usual wedge product.
\end{proof}
Now we use $\frak p^{(u,b)}$ to define a map 
$$
i_{\#,\text{\rm qm},(\frak b,b,u)}:
H(L(u);\Lambda_0) \to \Omega(X).
$$
Write $\frak b = \frak b_2 + \frak b_{\rm high}$ where $\frak b_2$ has degree
$2$ and $\frak b_{\rm high}$ has degree higher than $2$.

\begin{defn} We define\index{$i_{\#,\text{\rm qm},(\frak b,b,u)}$}
$$
i_{\#,\text{\rm qm},(\frak b,b,u)}(h)
= \sum_{\beta}\sum_{\ell}\frac{T^{\beta\cap \omega/2\pi}}{\ell!}
\exp{(\beta\cap \frak b_2)}
\rho^{b}(\partial \beta)
\frak p_{\ell;1;\beta}(\frak b_{\rm high}^{\otimes\ell}; h).
$$
Here $\rho^b$ is the representation of $\pi_1(L(u))$ obtained from $b$ as before.
\end{defn}
\begin{rem}
If $b = b_{0} + b_+$ where $b_0$ is of degree $1$ and the degree of
$b_+$ is $>1$, then we have
$$
i_{\#,\text{\rm qm},(\frak b,b,u)}(h)
= \sum_{\beta}\sum_{\ell,k} \frac{T^{\beta\cap \omega/2\pi}}{\ell!}
\exp{(\beta\cap \frak b_2)}\exp(b_0\cap \partial\beta)
\frak p_{\ell;k+1;\beta}
(\frak b_{\rm high}^{\otimes\ell}
; [h \otimes b_+^{\otimes k}]).
$$
This expression is closer to (\ref{fromptoi}) than the above definition. See Section  \ref{sec:deltaisthesame}.
\end{rem}
We consider the boundary operator defined from $\frak q_{\ell;k;\beta}^{\frak c}$ 
in the same way as (\ref{canboundary}) on $H(L(u);\Lambda_0)$.
We consider the usual derivation $d$ on $\Omega(X)$. 
\begin{lem}
$i_{\#,\text{\rm qm},(\frak b,b,u)}$ is a chain map.
\end{lem}
The proof is the same as the proof of Lemma \ref{i*chainmap}.
Thus $i_{\#,\text{\rm qm},(\frak b,b,u)}$ induces 
a map
\begin{equation}
i_{\ast,\text{\rm qm},(\frak b,b,u)}: HF((L(u);(\frak b,b)),(L(u);(\frak b,b));\Lambda_0)
\to H(X;\Lambda_0).
\end{equation}

We now consider the pairing $\langle \cdot,\cdot \rangle_{\text{\rm PD}_{L(u)}}$ by
\begin{equation}\label{196}
\langle h', h \rangle_{\text{\rm PD}_{L(u)}} = \int_{L(u)}h'\wedge h
\end{equation}
and $\langle \cdot,\cdot \rangle_{\text{\rm PD}_{X}}$ by
\begin{equation}\label{197}
\langle Q,v \rangle_{\text{\rm PD}_{X}} = \int_Q v = \int_X \text{\rm PD}_X(Q) \wedge v.
\end{equation}
We also recall the map
$$
i^{\frak c \ast}_{\text{\rm qm},(\frak b,b,u)}: H(X;\Lambda_0) \to H(L(u);\Lambda_0)
$$
introduced in Section \ref{sec:cyclic}.

\begin{lem}\label{unitprop} For any $Q\in \mathcal A(X)$,
the element $i^{\frak c \ast}_{\text{\rm qm},(\frak b,b,u)}(Q)$
is contained in $\Lambda_0\cdot\{e_L\}$.
\end{lem}
\begin{proof} By definition, $i^{\ast}_{\text{\rm qm},(\frak b,b,u)}(Q)$ is proportional to the
unit which is a degree-zero form. (See Lemma \ref{lem171}.) Since $\frak F$ is unital
filtered $A_{\infty}$ homomorphism, $\frak F_*$ sends the unit to the unit.
Then Lemma \ref{unitprop} follows from Lemma \ref{Isharprelation}.
\end{proof}

The following theorem describes the relationship between the
maps $i^{\frak c \ast}_{\text{\rm qm},(\frak b,b,u)}$ and
$i_{\ast,\text{\rm qm},(\frak b,b,u)}$.

\begin{thm}\label{pqpoincare} Let $Q \in \CA(X)$ and $h \in H^*(L(u);\Lambda_0)$. Then we have
\begin{equation}\label{popandqformula}
\langle i^{\frak c \ast}_{\text{\rm qm},(\frak b,b,u)}(Q),h
\rangle_{\text{\rm PD}_{L(u)}}
=
\langle Q,i_{\ast,\text{\rm qm},(\frak b,b,u)}(h)
\rangle_{\text{\rm PD}_{X}}.
\end{equation}
\end{thm}

\begin{cor}\label{degreemustn}
If the degree $n$ component of $h$ is zero, then the right hand side of $(\ref{popandqformula})$ is zero.
\end{cor}
\begin{proof}
This follows from Theorem \ref{pqpoincare} and Lemma \ref{unitprop}.
\end{proof}

\begin{proof}[Proof of Theorem \ref{pqpoincare}]
Let $\text{\bf p} \in E_{\ell}(\mathcal A[2])$. We put
$$
\Delta \text{\bf p} = \sum_c \text{\bf p}_{1;c} \otimes \text{\bf p}_{2;c},
$$
where $ \text{\bf p}_{i;c} \in E_{\ell_{c;i}}(\mathcal A[2])$. The following is the
key geometric result on which the proof of \eqref{popandqformula} is based.
\begin{lem}\label{bdrydescription}
$[0,1] \times \mathcal M_{k+1;\ell+1}^{\text{\rm main}}(\beta;\text{\bf p})$
has a Kuranishi structure and its boundary is the union of
following 4 types of fiber products:
\begin{enumerate}
\item
$\{0\} \times \mathcal M_{k+1;\ell+1}^{\text{\rm main}}(\beta;\text{\bf p})$.
\item
$\{1\} \times \mathcal M_{k+1;\ell+1}^{\text{\rm main}}(\beta;\text{\bf p})$.
\item
$\mathcal M_{k_1+1;\ell_{c;1}}(\beta_1;\text{\bf p}_{1;c})
{}_{\text{\rm ev}_0}\times_{\text{\rm ev}_j}
\left([0,1] \times  \mathcal M_{k_2+2;\ell_{c;2}+1}^{\text{\rm main}}(\beta_2;\text{\bf p}_{2;c})\right)$.
Here $k_1+k_2=k$, $\beta_1+\beta_2=\beta$, $j=1,\dots,k_2+1$.
\item
$\left([0,1] \times  \mathcal M_{k_1+1;\ell_{c;1}+1}(\beta_1;\text{\bf p}_{1;c}) \right)
{}_{\text{\rm ev}_0}\times_{\text{\rm ev}_j}
\mathcal M_{k_2+2;\ell_{c;1}}^{\text{\rm main}}(\beta_2;\text{\bf p}_{2;c}
)$. Here  $k_1+k_2=k$, $\beta_1+\beta_2=\beta$, $j=1,\dots,k_2+1$.
\end{enumerate}
And there exists a system $\frak{parap}$ of a continuous family of multisections on
$[0,1] \times \mathcal M_{k+1;\ell}^{\text{\rm main}}
(\beta;\text{\bf p})$ with the following properties:
\begin{enumerate}
\item[(a)]
It is transversal to zero and $T^n$-equivariant.
\item[(b)]
It is invariant under the cyclic permutation of the boundary marked
points and under the permutations of the interior marked points.
\item[(c)]
It is compatible with the forgetful map of the boundary marked points.
\item[(d)]
On the boundary component of type 1 above
it coincides with $\frak c$ perturbation, and 
on the boundary component of type
2 above it coincides with $\frak p$ perturbation.
\item[(e)]
On the boundary component of type 3 above it is compatible.
Here we put the $\frak c$ perturbation for the first factor.
\item[(f)]
On the boundary component of type 4 above it is compatible.
Here we put the $\frak c$ perturbation for the
second factor.
\end{enumerate}
\end{lem}
\begin{proof}
We recall that the $0$-th marked point is the interior marked point when we define
$\frak p$ and $\frak{parap}$.
So when we consider the codimension one boundary of
$\mathcal M_{k,\ell+1}^{\text{\rm main}}(\beta;\text{\bf p})$ the
$\frak p$ or $\frak{parap}$-perturbation is applied to the
factor containing the $0$-th marked point. Once this point is understood,
the proof is the same as other similar statements which we discuss in Sections \ref{sec:cyclicKura}, \ref{sec:compwiseplusTnequiv}.
(We remark that we use the cyclic symmetry of  $\frak{parap}$ (Lemma \ref{bdrydescription} (b)) 
together with cyclic symmetry of $\frak p$, here.)
\end{proof}
Going back to the proof of the theorem, we apply Stokes' formula and obtain
\begin{equation}\label{stokespara}
\aligned
&\int_{\partial([0,1] \times \mathcal M_{k+1;\ell+1}^{\text{\rm main}}(\beta;\text{\bf p}))^{\frak{parap}}}
(\text{\rm ev},\text{\rm ev}_0)^*((b_+^{\frak c})^{\otimes k} \times h)  \\
&= \int_{([0,1] \times \mathcal M_{k+1;\ell+1}^{\text{\rm main}}(\beta;\text{\bf p}))^{\frak{parap}}}
(\text{\rm ev},\text{\rm ev}_0)^*(d((b_+^{\frak c})^{\otimes k} \times h))
= 0.
\endaligned\end{equation}
\begin{rem}
In the proofs of the identity (\ref{stokespara}) and of similar others, we use the integration of a differential
form (or integration along the fiber of a map) on the zero set of a {\it continuous family of}
multisections. The dimension of this zero set is given by the sum of
the (virtual) dimension of the space with Kuranishi structure and
that of the parameter space $W$. (Note $W$ is defined only locally.)
We choose a ($T^n$-invariant) smooth form $\omega$ with compact support on $W$ that
has degree $\dim W$ and total mass $1$, and do integration of a
differential form after we take a wedge product with this form.
So strictly speaking, the integrand of the left hand side of (\ref{stokespara}) should
actually be $(\text{\rm ev},\text{\rm ev}_0)^*((b_+^{\frak c})^{\otimes k} \times h) \wedge \omega$, for example.
By abuse of notations, however, we write this integration just as (\ref{stokespara}) here and henceforth.
See \cite[Section 12]{fooo09}  for the precise definition.
\end{rem}
We decompose the left hand side of (\ref{stokespara}) according to the
decomposition of the boundary to Lemma \ref{bdrydescription}.1-4.
Let us write the term corresponding to Lemma \ref{bdrydescription}.$m$ as
$$
C(\beta;\text{\bf p};k;\ell;m), \quad m=1,2,3,4.
$$
We put
\begin{equation}\label{betasplit}
\frak b_{\rm high}^{\otimes\ell} = \sum_{I \in \frak I_{\ell}} c_{I}
\text{\bf f}_{i_1}
\otimes \dots \otimes \text{\bf f}_{i_{\ell}}
= \sum_{I \in \frak I_{\ell}} c_{I} \text{\bf f}_{I}.
\end{equation}
Here $\frak I_\ell$ is the set of multi-indices $I = (i_1,\dots,i_{\ell})$ and
$c_{I} \in \Lambda_0$.
\par
We now calculate
\begin{equation}\label{termforj}
\sum_{\ell,k,\beta} \sum_{I \in \frak I_{\ell}}
\frac{c_I}{\ell!}\rho(\partial\beta)e^{\frak b_2\cap \omega}
C(\beta;\text{\bf f}_{I} \otimes Q;k;\ell;m)
\end{equation}
for $m=1,2,3,4$.
\par
 The cases $m=1$ and $m=2$ give rise to the left and the right hand sides of Theorem \ref{pqpoincare} by Lemma \ref{bdrydescription}
(d) respectively.
\par
Let us consider the case $m=3$. Decompose $I = I_1 \sqcup I_2$
for $I_j \in \frak I_{\ell_j}$, $j = 1,  2$ with $\ell_1 + \ell_2 = \ell$ and
$
k = k_1 + k_2,  \beta = \beta_1 + \beta_2,  \text{\bf p} = \text{\bf p}_1 \otimes \text{\bf p}_2.
$
Then, thanks to $\Delta_{\rm shuff}(e^{\frak b})= e^{\frak b} \otimes e^{\frak b}$ (see Remark \ref{quotientsub})
we have
\begin{equation}\label{coalgecoeff}
\Delta_{\rm shuff}(e^{\frak b}) 
= \sum_{\ell_1,\ell_2}\sum_{I_1 \in \frak I_{\ell_1}}\sum_{I_2 \in \frak I_{\ell_2}} c_{I_1} \text{\bf f}_{I_1} \otimes c_{I_2}\text{\bf f}_{I_2}.
\end{equation}
We put
\begin{equation}\label{Edefine}
E =
\sum_{\beta_1,\ell_1,k_1}\sum_{I_1\in \frak I_{\ell_1}}
\frac{c_{I_1}}{\ell_1!}\text{\rm ev}_{0*}
\left(
(\text{\rm ev})^*((b_+^{\frak c})^{\otimes k_1});\mathcal M_{k_1+1;\ell_1+1}^{\text{\rm main}}(\beta;\text{\bf f}_{I_1} \times Q)^{\frak c}
\right).
\end{equation}
This is an element of $H(L(u);\Lambda_0)$.
Using (\ref{coalgecoeff}), the case $m=3$ of
(\ref{termforj}) becomes:
\begin{equation}\label{contri3}
\aligned
&\sum_{\beta_2,\ell_2,k_2}\sum_{j=1}^{k_2}\sum_{I_2\in \frak I_{\ell_2}}
\frac{c_{I_2}}{\ell_2!}\rho(\partial\beta)e^{\frak b_2\cap \omega} \\
&\int_{([0,1] \times \mathcal M_{k_2+1;\ell_2+1}^{\text{\rm main}}(\beta_2;\text{\bf p}_2))^{\frak{parap}}}
(\text{\rm ev},\text{\rm ev}_0)^*((b_+^{\frak c})^{j} \times E \times (b_+^{\frak c})^{k_2-j}
\times h).
\endaligned
\end{equation}
In the same way as the proof of Lemma \ref{unitprop}, we can prove that
$E$ is proportional to the unit.
Therefore Lemma \ref{bdrydescription} (c) implies that (\ref{contri3}) is $0$.
Namely (\ref{termforj}) is zero for the case $m=3$.
\par
We finally consider the case $m=4$.
We define
\begin{equation}
\aligned
F =
&\sum_{\beta_1,\ell_1,k_1}\sum_{I_1\in \frak I_{\ell_1}}
\frac{c_{I_1}}{\ell_1!}\rho(\partial\beta)e^{\frak b_2\cap \omega}
\\
&
\text{\rm ev}_{0*}\left(
\text{\rm ev}^*(b_+^{\frak c})^{k_1};[0,1]\times
\mathcal M_{k_1+1;\ell_1+1}^{\text{\rm main}}(\beta;\text{\bf f}_{I_1}\times Q)^{\frak{parap}}
\right).
\endaligned
\end{equation}
Then, by using (\ref{coalgecoeff}),
the case $m=4$ of (\ref{termforj}) is
\begin{equation}
\aligned
\sum_{\beta_2,\ell_2,k_2}&\sum_{I_2\in \frak I_{\ell_2}}\sum_{j=1}^{k_2}
\frac{c_{I_2}}{\ell_2!}\rho(\partial\beta)e^{\frak b_2\cap \omega}\\
&\text{\rm ev}_{0*}\left(
(\text{\rm ev})^*((b_+^{\frak c})^{j}\times F \times (b_+^{\frak c})^{(k_2-j)}
\times h);
\mathcal M_{k_2+1;\ell_2+1}^{\text{\rm main}}(\beta;\text{\bf f}_{I_2})^{\frak c}
\right).
\endaligned\end{equation}
By definition this is
$$
\langle \delta^{\frak b,b^{\frak c}}_{\frak c}(F),h\rangle_{\text{\rm PD}_{L(u)}}
$$
and is
zero.
This is because $(u,b^{\frak c})$ is a critical
point of $\frak{PO}_{\frak b}^{\frak c}$ and so $\delta^{\frak b,b^{\frak c}}_{\frak c}=0$ and Corollary \ref{isowithc}.
The proof of Theorem \ref{pqpoincare} is complete.
\end{proof}
\begin{rem}\label{cycandPD}
We used the cyclic symmetry in the proof of Theorem \ref{pqpoincare}
above.
In other words the above proof does not work to prove the same
conclusion (\ref{popandqformula}) when we replace
the left hand side by
$\langle i^{\ast}_{\text{\rm qm},(\frak b,b,u)}(Q),h
\rangle_{\text{\rm PD}_{L(u)}}$. 
The reason is as follows: 
To construct the $\frak p$-perturbation
that is compatible with the $\frak q$-perturbation,
we need to take the $\frak q$-perturbation for the first factor in
Lemma \ref{pmodulikura}.6.
Moreover we have to take the $0$-th marked point
of this factor to the point which intersects with the
second factor. (Namely the point which becomes the singular
point after gluing.)
\par
On the other hand, for Formula (\ref{popandqformula})
to hold, we need to take the
$0$-th marked point so that
we pull back the form $h$ by ${\rm ev}$ at the $0$-th marked point.
\par
We cannot realize both of them at the same time in the situation of
Lemma \ref{bdrydescription}.4.
\end{rem}
\par
\section{Moduli space of holomorphic annuli}
\label{sec:annuli}
Let $(u,b) \in \text{\rm Crit}(\frak{PO}_{\frak b})$. We take a basis $\text{\bf e}_I$ with
$I \in 2^{\{1,\dots,n\}}$ of $H(L(u);\Q)$. Let
$$
g_{IJ} = \langle \text{\bf e}_I, \text{\bf e}_J \rangle_{\text{\rm PD}_{L(u)}},
\quad
[g^{IJ}]
= [g_{IJ}]^{-1}.
$$
Here $[g^{IJ}]$, $[g_{IJ}]$ are $2^n \times 2^n$ matrices.
We recall that the operator
$$
\frak m_2^{\frak c,\frak b,b}:
HF((L(u),\frak b,b),(L(u),\frak b,b);\Lambda_0)^{\otimes 2}
\to
HF((L(u),\frak b,b),(L(u),\frak b,b);\Lambda_0)
$$
is defined by
\begin{equation}
\aligned
\frak m_2^{\frak c,\frak b,b}(h_1,h_2) &=
\sum_{k_1,k_2,k_3}\sum_{\ell,\beta} {T^{\beta\cap \omega/2\pi}}
\exp(\beta\cap \frak b_2)\rho(\partial \beta)\\
&\frak q_{\ell;2+k_1+k_2+k_3;\beta}^{\frak c}(\frak b_{\rm high}^{\otimes\ell}
; (b_+^{\frak c})^{k_1} \times h_1
\times (b_+^{\frak c})^{k_2}\times h_2\times (b_+^{\frak c})^{k_3}).
\endaligned
\end{equation}
(See (\ref{mcyclicdef}), (\ref{mcyclicdef2}).)
In the next theorem, we involve the harmonic forms
on $L(u)$ of all degrees beside those of top degree $n$, although  in the end
only the top degree forms turn out to give non-trivial contributions in the formula below.
\begin{thm}\label{annulusmain}
Let $\frak v, \frak w \in H(L(u);\Lambda_0)$. Then we have
\begin{equation}\label{21mainformula}
\aligned
&\langle i_{\ast,\text{\rm qm},(\frak b,b^{\frak c})}(\frak v),
i_{\ast,\text{\rm qm},(\frak b,b^{\frak c})}(\frak w)
\rangle_{\text{\rm PD}_{X}} \\
&=
\sum_{I,J \in 2^{\{1,\dots,n\}}} (-1)^{*}g^{IJ}
\langle \frak m_2^{\frak c,\frak b,b^{\frak c}}(\text{\bf e}_I,\frak v),
\frak m_2^{\frak c,\frak b,b^{\frak c}}(\text{\bf e}_J,\frak w)
\rangle_{\text{\rm PD}_{L(u)}}.
\endaligned\end{equation}
\end{thm}
Here $* = \frac{n(n-1)}{2}$ in case the degrees of $\frak v, \frak w$ are $n$.
We prove it in
Subsection \ref{subsec:signproof}.
(See (\ref{eqpreconclusion}).)
For the sign $*$ in the general case, see Proposition \ref{propconclusion}.
\begin{rem}
We do not assume nondegeneracy of the critical point $(u,b)$
of $\frak{PO}_{\frak b}$
in Theorem \ref{annulusmain}.
Theorem \ref{clifford} implies that,  in the case
of degenerate $(u,b)$ and of Condition \ref{fano2jigen}.1 or 2 being satisfied,
the right hand side of (\ref{21mainformula}) is $0$ and so the left hand side is
also $0$.
\end{rem}
The right hand side is closely related to the invariant
$Z(\frak b,b)$ in Definition \ref{res2}.
\par
The proof of Theorem \ref{annulusmain} is based on the study of
the moduli space of holomorphic maps from the annuli. We begin with defining the
corresponding moduli space.
\par
Consider a bordered semi-stable curve $\Sigma$ of genus $0$ with two
boundary components $\partial_1 \Sigma$,
$\partial_2 \Sigma$ and $k_1+1$, $k_2+1$ boundary marked points
$\vec z_1 = (z^1_{0},z_{1,1},\dots,z_{1,k_1})$,
$\vec z_2 = (z^2_{0},z_{2,1},\dots,z_{2,k_2})$ on each
of the components $\partial_1 \Sigma$,
$\partial_2 \Sigma$ of the boundary, and $\ell$ interior marked points
$\vec z^{+} = (z^+_1,
\dots,z^+_{\ell})$.
We denote it by $(\Sigma;\vec z_1,\vec z_2,\vec z^{+})$
or sometimes by $\Sigma$ for short.
We say that $(\Sigma;\vec z_1,\vec z_2,\vec z^{+})$
is in the {\it main component} if $\vec z_i$ respects the cyclic order of the
boundary.
(When we regard $\Sigma$ as an annulus in $\C$,
we take the counter clockwise cyclic order for the outer circle and the
clockwise cyclic order for the inner circle.)
Let $(\Sigma;\vec z_1,\vec z_2,\vec z^{+})$  be as above and
consider a holomorphic map
$w: (\Sigma,\partial \Sigma) \to (X,L(u))$.
We say $(w;\Sigma;\vec z_1,\vec z_2,\vec z^{+})$
is a {\it stable map} if the group of automorphisms is a finite group.
\begin{defn}
We put $\beta \in H_2(X,L(u);\Lambda_0)$.
We denote by
\index{$\mathcal M_{(k_1+1,k_2+1);\ell}^{\text{ann;main}}
(\beta)$}
$$
\mathcal M_{(k_1+1,k_2+1);\ell}^{\text{ann;main}}
(\beta)
$$
the set of all isomorphism classes of stable maps $(w;\Sigma;\vec z_1,
\vec z_2,\vec z^{+})$
such that $w_*[\Sigma] = \beta$ and
$(\Sigma;\vec z_1,\vec z_2,\vec z^{+})$
is in the main component, which we described above.
\end{defn}
In the same way as  \cite[Section 3]{fooo00} and  \cite[Section 7.1]{fooobook2}
we can prove that
$\mathcal M_{(k_1+1,k_2+1);\ell}^{\text{ann;main}}$ has
a Kuranishi structure with boundary and corners.
There are obvious evaluation maps
$$\aligned
\text{\rm ev} = & (\text{\rm ev}^1,\text{\rm ev}^2,\text{\rm ev}^{\text{\rm int}},\text{\rm ev}^1_0,\text{\rm ev}^2_0) \\
&: \mathcal M_{(k_1+1,k_2+1);\ell}^{\text{ann;main}}(\beta)
\to L(u)^{k_1} \times L(u)^{k_2} \times X^{\ell} \times L(u)
\times L(u).
\endaligned$$
Here $\text{\rm ev}^1_0$, $\text{\rm ev}^2_0$ are the evaluation maps
at $z^1_0$, $z^2_0$, and
$\text{\rm ev}^1$, $\text{\rm ev}^2$ are evaluation maps at
$(z_{1,1},\dots,z_{1,k_1})$ and $(z_{2,1},\dots,z_{2,k_2})$,
respectively.
We consider the case $k_1\ge 0$, $k_2 \ge 0$ only, in this
paper.
\par
Let $\mathcal M_{(1,1);0}^{\text{ann;main}}$ be the moduli space
$\mathcal M_{(0+1,0+1);0}^{\text{ann;main}}(0)$
which corresponds to the case $X=L(u) = $one point.
There exists an obvious forgetful map
\begin{equation}\label{forgetfromann}
\frak{forget}: \mathcal M_{(k_1+1,k_2+1);\ell}^{\text{ann;main}}(\beta)
\to \mathcal M_{(1,1);0}^{\text{ann;main}}.
\end{equation}
(Here we forget $(z_{1,1},\dots,z_{1,k_1})$ and $(z_{2,1},\dots,z_{2,k_2})$
among the boundary marked points.)
\begin{conv}
Hereafter we write
$$
\frak{forget}^{-1}(K) \cap
\mathcal M_{(k_1+1,k_2+1);\ell}^{\text{ann;main}}(\beta)
$$
for the inverse image of $K \subset \mathcal M_{(1,1);0}^{\text{ann;main}}$
by (\ref{forgetfromann}).
\end{conv}
\par
In a way similar to Lemma \ref{stratifyD12} we can
describe $\mathcal M^{\text{\rm main}}_{(1,1);0}$ as follows.
\begin{lem}\label{stratifyAn11}
$\mathcal M^{\text{\rm ann;main}}_{(1,1);0}$ has a stratification
such that the interior of each stratum are as follows.
\begin{enumerate}
\item $\text{\rm Int}D^2 \setminus \{0\}$.
\item An arc $(-1,1)$.
\item Two points. $[\Sigma_1]$, $[\Sigma_{2}]$.
\end{enumerate}
\par
The point $0 \in D^2$ becomes $[\Sigma_1]$
in the closure of the stratum 1.
The two boundary points $\pm 1$ in the closure of $(-1,1)$
become $[\Sigma_{2}]$.
Thus the union of the two strata $(-1,1)$, $[\Sigma_{2}]$ becomes a circle,
which turns out to be $\partial D^2$
in the closure of the stratum 1.
\par
Each of the above strata corresponds to the combinatorial type
of the elements of $\mathcal M^{\text{\rm main}}_{1;2}$.
\end{lem}
\begin{proof}
We first define two points $[\Sigma_1]$,
$[\Sigma_2]$ of
$\mathcal M_{(1,1);0}^{\text{ann;main}}$.
\par
Let $D^2_1$ and $D^2_2$ be two copies of disk $D^2 = \{z \in \C \mid
\vert z\vert \le 1\}$. We identify $0 \in D^2_1$ and $0 \in D^2_2$
to obtain $\Sigma_1$. We put
$\partial_1\Sigma_1 = \partial D^2_1$, $\partial_2\Sigma_1 = \partial D^2_2$,
$z^1_0 = 1 \in \partial D^2_1$, $z^2_0 = 1 \in \partial D^2_2$.
Thus we have $(\Sigma_1;z^1_1,z^1_2) \in \mathcal M_{(1,1);0}^{\text{ann;main}}$.
\par
We next define $\Sigma_2$. We consider two copies
$D^2_1$, $D^2_2$ of disk.
We put $\frak z^i_k = e^{2k\pi\sqrt{-1}/3} \in \partial D^2_i$ for $i=1,2$,
$k=0,1,2$.
\par
We identify $\frak z^1_1$ with $\frak z^2_1$ and
$\frak z^1_2$ with $\frak z^2_2$ respectively and obtain
$\Sigma_2$.
We put $\frak z^i_0 = z^1_i$ and $i=1,2$. We call
$\partial_i \Sigma_2$ the component which contains $z_i^1$.
We have thus defined
$(\Sigma_2;z^1_1,z^1_2) \in \mathcal M_{(1,1);0}^{\text{\rm ann;main}}$.
\par
We next describe the stratum 2.
We consider the moduli space $\mathcal M_{4;0}^{\text{\rm main}}$ of
disks with $4$ boundary marked points without
interior marked points. It is an arc $[0,1]$.
Let $(D^2;\frak z_1,\frak z_2,\frak z_3,\frak z_4)$ be
a point in the interior of $\mathcal M_{4;0}^{\text{\rm main}}$.
We identify $\frak z_2$ with $\frak z_4$ to obtain $\Sigma$ and
put $z_1 = \frak z_1$, $z_2 =\frak z_4$.
Then $(\Sigma;z_1,z_2) \in \mathcal M_{(1,1);0}^{\text{ann;main}}$.
\par
It is easy to see that as $(D^2;\frak z_1,\frak z_2,\frak z_3,\frak z_4)$ approaches
to the boundary of $\mathcal M_{4;0}^{\text{\rm main}}$
$(\Sigma;z_1,z_2)$
converges to $[\Sigma_2]$ in
$\mathcal M_{(1,1);0}^{\text{ann;main}}$.
Thus the union of the strata 2,3 we have described above becomes a circle.
It is easy to see that it gives the stratification as in Lemma \ref{stratifyAn11}.
\end{proof}
\begin{lem}\label{bdryMannu}
\index{Kuranishi structure!Kuranishi structure on moduli of holomorphic annuli}
The space $\mathcal M_{(k_1+1,k_2+1);\ell}^{\text{\rm ann;main}}(\beta)$
has a Kuranishi structure with the following properties.
\begin{enumerate}
\item[(1)] It is transversal to $0$.
\item[(2)] It is  invariant under the permutation of $\ell$ interior marked points.
\item[(3)] It is compatible with the forgetful map of $1$st - $k_1$-th 
and $1$-st - $k_2$-th marked points on the two boundary components, 
in the sense of \cite[Section 3]{fooo091}.
\end{enumerate}
Moreover the boundary of $\mathcal M_{(k_1+1,k_2+1);\ell}^{\text{\rm ann;main}}(\beta)$
is the union of the following 5 types of fiber products:
\begin{enumerate}
\item
$$
\mathcal M_{k'_1+1;\ell'}^{\rm main}(\beta')^{\frak c} {}_{\text{\rm ev}_0}\times_{\text{\rm ev}^1_i}
\mathcal M_{(k''_1+1,k_2+1);\ell''}^{\text{\rm ann;main}}(\beta'').
$$
Here $k'_1+k''_1 = k_1+1$, $\ell'+\ell'' = \ell$,
$\beta'+\beta'' = \beta$, $i = 1,\dots,k''_1$.
\item
$$
\mathcal M_{k'_2+1;\ell'}^{\rm main}(\beta')^{\frak c} {}_{\text{\rm ev}_0}\times_{\text{\rm ev}^2_i}
\mathcal M_{(k_1+1,k''_2+1);\ell''}^{\text{\rm ann;main}}(\beta'').
$$
Here $k'_2+k''_2 = k_2+1$, $\ell'+\ell'' = \ell$,
$\beta'+\beta'' = \beta$, $i = 1,\dots,k''_2$.
\item
$$
\mathcal M_{(k'_1+1,k_2+1);\ell'}^{\text{\rm ann;main}}(\beta')
 {}_{\text{\rm ev}^1_0}\times_{\text{\rm ev}_i} 
\mathcal M_{k''_1+1;\ell''}^{\rm main}(\beta'')
^{\frak c}.
$$
Here $k'_1+k''_1 = k_1+1$, $\ell'+\ell'' = \ell$,
$\beta'+\beta'' = \beta$, $i = 1,\dots,k''_1$.
\item
$$
\mathcal M_{(k_1+1,k'_2+1);\ell'}^{\text{\rm ann;main}}(\beta')
 {}_{\text{\rm ev}^2_0}\times_{\text{\rm ev}_i} 
 \mathcal M_{k''_2+1;\ell''}^{\rm main}(\beta'')
^{\frak c}.
$$
Here $k'_2+k''_2 = k_1+1$, $\ell'+\ell'' = \ell$,
$\beta'+\beta'' = \beta$, $i = 1,\dots,k''_2$.
\item
$\frak{forget}^{-1}(\partial \mathcal M_{(1,1);0}^{\text{\rm ann;main}})
\cap \mathcal M_{(k_1+1,k_2+1);\ell}^{\text{\rm ann;main}}(\beta)$.
\end{enumerate}
\end{lem}
We note that we put $\frak c$ in the first factors of 1, 2 and
third factors of 3, 4. This is to clarify we always use the $\frak c$-perturbation  for such factor.
This applies to many other similar cases. Namely we
use the $\frak c$-perturbation for the component of 
disk bubble.
\par
The proof of Lemma \ref{bdryMannu} is given in Subsection \ref{subsec:kuraanul}. See also the end of Subsection \ref{subsec:kuraanul3}.
\begin{rem}
We do not claim our perturbation in Lemma \ref{bdryMannu} is cyclically symmetric.
\end{rem}
For $\text{\bf p}$ as in Section \ref{sec:frakqreview},
we define\index{$\mathcal M_{(k_1+1,k_2+1);\ell}^{\text{ann;main}}
(\beta;\text{\bf p})$}
\begin{equation}
\mathcal M_{(k_1+1,k_2+1);\ell}^{\text{ann;main}}
(\beta;\text{\bf p})
= \mathcal M_{(k_1+1,k_2+1);\ell}^{\text{ann;main}}
(\beta) {}_{\text{\rm ev}^{\text{\rm int}}} \times_{X^{\ell}}
(\text{\bf p}(1) \times \dots \times \text{\bf p}(\ell)).
\end{equation}
(See Subsection \ref{subsec:fiberproduct} for the way how we take the fiber product.)
We can describe its boundary by using Lemma \ref{bdryMannu},
in a way similar to the cases we have already discussed in several other occasions.

The proof of the next lemma is easy and so omitted.
\begin{lem}\label{forgetsmooth2}
The forgetful map $(\ref{forgetfromann})$ is continuous and is a smooth submersion
on each stratum.
\end{lem}
The proof of Theorem \ref{annulusmain} consists of a series of propositions.
We write
\begin{equation}
[\frak v;k] = \frak v \otimes\underbrace{b_+^{\frak c} \otimes\dots \otimes b_+^{\frak c}}_{k}
\in B_{k+1}(H(L(u);\Lambda_0)[1])
\end{equation}
for the simplicity of notations. Hereafter we apply the following:
\begin{conv}\label{convint}
When we integrate the pull back of $[\frak v;k_1]$ or $[\frak w;k_2]$
we always pull back the forms $\frak v$, $\frak w$ by the evaluation maps 
at the 0th marked points and $\frak b_+^{\frak c}$ by the 
evaluation map at 1st - $k$th marked points.
\end{conv}
\begin{prop}\label{doublehatandpdX}
$$\aligned
&\sum_{\beta}\sum_{\ell}\sum_{k_1,k_2}\frac{T^{\beta\cap \omega/2\pi}}{\ell!}
\exp(\beta \cap \frak b_2)\rho(\partial \beta)
\\
&\qquad\qquad\qquad
\int_{\frak{forget}^{-1}([\Sigma_1])
\cap \mathcal M_{(k_1+1,k_2+1);\ell}^{\text{\rm ann;main}}
(\beta;\frak b_{\rm high}^{\otimes\ell})}(\text{\rm ev}^1,\text{\rm ev}^2)^*([\frak v;k_1] \times [\frak w;k_2]) \\
&=
\sum_{\beta_1, \beta_2}\sum_{\ell_1, \ell_2
\atop : \ell_1+\ell_2=\ell}\sum_{k_1,k_2}
\frac{T^{(\beta_1+\beta_2)\cap \omega}}{\ell_1!\ell_2!}
\exp((\beta_1+\beta_2) \cap \frak b_2)\rho(\partial (\beta_1+\beta_2)) \\
&\qquad\qquad\qquad
\langle
\frak p_{\ell_1,\beta_1}(\frak b_{\rm high}^{\otimes\ell_1}
; [\frak v;k_1]),
\frak p_{\ell_2,\beta_2}(\frak b_{\rm high}^{\otimes\ell_2}
; [\frak w;k_2])
\rangle_{\text{\rm PD}_{X}}.
\endaligned$$
\end{prop}
\begin{prop}\label{doublehatistwodisk}
$$\aligned
&\sum_{\beta}\sum_{\ell}\sum_{k_1,k_2}\frac{T^{\beta\cap \omega/2\pi}}{\ell!}
\exp(\beta\cap \frak b_2)\rho(\partial \beta)\\
&\qquad\qquad\qquad
\int_{\frak{forget}^{-1}([\Sigma_1])
\cap \mathcal M_{(k_1+1,k_2+1);\ell}^{\text{\rm ann;main}}
(\beta;\frak b_{\rm high}^{\otimes\ell})}(\text{\rm ev}^1,\text{\rm ev}^2)^*([\frak v;k_1] \times [\frak w;k_2]) \\
&= \sum_{\beta}\sum_{\ell}\sum_{k_1,k_2}\frac{T^{\beta\cap \omega/2\pi}}{\ell !}
\exp(\beta\cap \frak b_2)\rho(\partial \beta)
\\
&\qquad\qquad\qquad
\int_{\frak{forget}^{-1}([\Sigma_2])
\cap \mathcal M_{(k_1+1,k_2+1);\ell}^{\text{\rm ann;main}}
(\beta;\frak b_{\rm high}^{\otimes\ell})}(\text{\rm ev}^1,\text{\rm ev}^2)^*([\frak v;k_1] \times [\frak w;k_2]).
\endaligned$$
\end{prop}
\begin{prop}\label{twodiskandm2}
$$
\aligned
&\sum_{\beta}\sum_{\ell}\sum_{k_1,k_2}\frac{T^{\beta\cap \omega/2\pi}}{\ell!}
\exp(\beta\cap \frak b_2)\rho(\partial \beta)\\
&\qquad\qquad\qquad
\int_{\frak{forget}^{-1}([\Sigma_2])
\cap \mathcal M_{(k_1+1,k_2+1);\ell}^{\text{\rm ann;main}}
(\beta;\frak b_{\rm high}^{\otimes\ell})}(\text{\rm ev}^1,\text{\rm ev}^2)^*([\frak v;k_1] \times [\frak w;k_2]) \\
&=
\sum_{\beta_1, \beta_2}\sum_{\ell_1, \ell_2
\atop : \ell_1+\ell_2=\ell}\sum_{I,J \in 2^{\{1,\dots,n\}}}
g^{IJ}
\frac{T^{(\beta_1+\beta_2)\cap \omega}}{\ell_1!\ell_2!}
\exp((\beta_1+\beta_2)\cap \frak b_2)\rho(\partial (\beta_1+\beta_2))
\\
&\qquad\qquad\qquad\qquad
\langle \frak q^{\frak c,\rho}_{\ell_1;2;\beta_1}(\frak b_{\rm high}^{\otimes\ell_1}
;\frak v \otimes \text{\bf e}_I),
\frak q^{\frak c,\rho}_{\ell_2;2;\beta_2}(\frak b_{\rm high}^{\otimes\ell_2}
; \text{\bf e}_J\otimes \frak w)
\rangle_{\text{\rm PD}_{L}}.
\endaligned
$$
\end{prop}
Here in the right hand side of Proposition \ref{twodiskandm2}, the operator 
${\frak q}^{\frak c,\rho}_
{\ell;2;\beta}$ is defined by
\begin{equation}
\aligned
&\frak q^{\frak c,\rho}_{\ell;2;\beta}([\text{\bf p}];x_1,x_2) \\
&=
\sum_{k_1,k_2,k_3}
\frak q^{\frak c,\rho}_{\ell;2+k_1+k_2+k_3;\beta}([\text{\bf p}];
(b_+^{\frak c})^{k_1} \otimes x_1\otimes (b_+^{\frak c})^{k_2} \otimes x_2
 \otimes (b_+^{\frak c})^{k_3}).
\endaligned
\end{equation}
\par
It is straightforward to see that Propositions \ref{doublehatandpdX},
\ref{doublehatistwodisk} and \ref{twodiskandm2} imply Theorem \ref{annulusmain}.
The rest of this section is occupied with the proofs of Propositions \ref{doublehatandpdX},
\ref{doublehatistwodisk} and \ref{twodiskandm2}.
Actually the geometric idea behind the validity of these
propositions is clear. The main task is to construct
appropriate perturbations (a family of multisections) of the
moduli spaces involved so that they are consistent with
the perturbations we fixed in the earlier stage to define various other
operators considered. We explain the details of this construction
in this section.
In the course of doing this, we also determine our choice of perturbations
made for $\frak{forget}^{-1}([\Sigma_1])
\cap \mathcal M_{(1,1);\ell}^{\text{\rm ann;main}}
(\beta;\frak b_{\rm high}^{\otimes\ell})$
and for $\frak{forget}^{-1}([\Sigma_2])
\cap \mathcal M_{(1,1);\ell}^{\text{\rm ann;main}}
(\beta;\frak b_{\rm high}^{\otimes\ell})$.
\par
\begin{proof}[Proof of Proposition \ref{doublehatandpdX}]
The arguments used in this proof is similar to the arguments used
around Lemmata \ref{deltaidentify}, \ref{prodmainlemma} etc.. 
We start with the following:\index{$\mathcal M(k_1,k_2;\ell_1,\ell_2;\beta_1,\beta_2;\text{\bf p}_1,\text{\bf p}_2)$}
\begin{defn}\label{twofiberproductofp}
$$\aligned
\mathcal M(k_1,k_2;\ell_1,\ell_2;\beta_1,\beta_2;\text{\bf p}_1,\text{\bf p}_2)
= \mathcal M^{\text{\rm main}}_{k_1;\ell_1+1}(\beta_1;\text{\bf p}_1)
{}_{\text{\rm ev}^{\text{\rm int}}_{\ell_1+1}} \times_{\text{\rm ev}^{\text{\rm int}}_{\ell_2+1}}
\mathcal M^{\text{\rm main}}_{k_2;\ell_2+1}(\beta_2;\text{\bf p}_2).
\endaligned$$
Here the fiber product is taken over $X$.
Since evaluation maps are weakly submersive, this space has a
fiber product Kuranishi structure.
(See \cite[Section A.1.2]{fooobook2}.)
The Kuranishi structure we put in the right hand side is $\frak p$-Kuranishi structure.
\end{defn}
Let
$\text{Split}((\mathbb L_1,\mathbb L_2),\text{\bf p})
= (\text{\bf
p}_1,\text{\bf p}_2)$.
(See Section \ref{sec:frakqreview} for the notation.)
We can glue elements of
$\mathcal M(k_1,k_2;\ell_1,\ell_2;\beta_1,\beta_2;\text{\bf p}_1,
\text{\bf p}_2)$
at  $\ell_1+1$-th and $\ell_2+1$-th marked points. We then
get an element of
$$
\frak{forget}^{-1}([\Sigma_1]) \cap \mathcal M_{(k_1,k_2);\ell}^{\text{\rm ann;main}}
(\beta_1+\beta_2;\text{\bf p})
$$
where $\ell = \ell_1+\ell_2$.
We have thus obtained a map
\begin{equation}\label{gluemapannul}
\aligned
\text{\rm Glue}
~: ~&\mathcal M(k_1,k_2;\ell_1,\ell_2;\beta_1,\beta_2;\text{\bf p}_1,\text{\bf p}_2)\\
&
\to \frak{forget}^{-1}([\Sigma_1]) \cap \mathcal M_{(k_1,k_2);\ell}^{\text{\rm ann;main}}
(\beta_1+\beta_2;\text{\bf p}).
\endaligned\end{equation}
Let $\beta$, $k_1, k_2$, $\ell$ be given.
We consider the disjoint union of the moduli spaces
$\mathcal M(k_1,k_2;\ell_1,\ell_2;\beta_1,\beta_2;\text{\bf p}_1,\text{\bf p}_2)$
where
$(\text{\bf
p}_1,\text{\bf p}_2)
= \text{Split}((\mathbb L_1,\mathbb L_2),\text{\bf p})$ with
$(\mathbb L_1,\mathbb L_2) \in \text{\rm Shuff}(\ell)$
and $\beta_1 +\beta_2=\beta$.
\begin{lem}\label{glueisolemann}
The maps $(\ref{gluemapannul})$ define a map from the above (disjoint) union of the
moduli spaces
$\mathcal M(k_1,k_2;\ell_1,\ell_2;\beta_1,\beta_2;\text{\bf p}_1,\text{\bf p}_2)$
to $\frak{forget}^{-1}([\Sigma_1]) \cap
\mathcal M_{(k_1,k_2);\ell}^{\text{\rm ann;main}}(\beta;\text{\bf p})$.
It is surjective and is an isomorphism outside the codimension $2$ strata.
\end{lem}
The proof is the same as the proof of Lemma \ref{deltaidentify} and omitted.
The reason why the map in Lemma \ref{glueisolemann} is not an
isomorphism is the same as in the case of Lemma \ref{smoothmaptoDM}.
Namely it should be regarded as a `normalization' of the `normal
crossing divisor' $\frak{forget}^{-1}([\Sigma_1])$.
\par
A problem similar to one mentioned in Remark \ref{rem:1014}
occurs here. This is the reason why we do not require our multiseciton to be $T^{n}$-equivariant
in the next lemma.
\par
\begin{lem}\label{perturbnbdsigma1}
\index{Kuranishi structure!Kuranishi structure on moduli of holomorphic annuli}
There exists a Kuranishi structure and a system of continuous families of
multisections defined on a neighborhood of $\frak{forget}^{-1}([\Sigma_1])$
with the following properties:
\begin{enumerate}
\item It is transversal to $0$.
\item
It is compatible at the boundaries described in
Lemma \ref{bdryMannu}.
\item It is compatible with the forgetful map of the boundary marked points.
\item It is invariant under the cyclic permutations of boundary marked points and
arbitrary permutations of the interior marked points.
\item
It is compatible with the fiber product description of $\frak{forget}^{-1}([\Sigma_1])$.
\item It is transversal to the map $(\ref{gluemapannul})$ in the sense we describe below.
\end{enumerate}
\end{lem}
We now make the statements of Condition 5 precise.
Let $U$ be a Kuranishi neighborhood of a point $\frak x$ in
$\mathcal M(k_1,k_2;\ell_1,\ell_2;\beta_1,\beta_2;\text{\bf p}_1,\text{\bf p}_2)$.
Denote
$$
\text{\rm Glue}_{k_1,k_2;\ell_1,\ell_2;\beta_1,\beta_2;\text{\bf p}_1,\text{\bf p}_2}
(\frak x) = \frak y \in \mathcal M_{(k_1,k_2);\ell}^{\text{\rm ann;main}}
(\beta_1+\beta_2;\text{\bf p}).$$
Let $U'$ be a Kuranishi neighborhood of $\frak y$
and $E'$ an obstruction bundle.
Let $\frak s_i(x,w)$ be a branch of a member of our continuous family of multisections.
Here $w \in W$ is an element of the parameter space
of our continuous family of multisections and
$x \in U'$ is an element of the Kuranishi neighborhood. We require
that $\frak s_i^{-1}(0)$ is transversal to
the map $: U\times W \to U' \times W$ induced by
$\text{\rm Glue}_{k_1,k_2;\ell_1,\ell_2;\beta_1,\beta_2;\text{\bf p}_1,\text{\bf p}_2}$.
\begin{proof}
The proof can be done by an induction over $\beta\cap \omega$ and the number of
boundary marked. See Subsection \ref{subsec:kuraanul3} and 
also Section \ref{sec:equikuracot}.
\end{proof}
The choice of a system of Kuranishi structures and continuous families of multisections
satisfying Lemma \ref{perturbnbdsigma1} induces a system of
Kuranishi structures and continuous families of multisections on
$\mathcal M(k_1,k_2;\ell_1,\ell_2;\beta_1,\beta_2;\text{\bf p}_1,\text{\bf p}_2)$.
We call this a \emph{pull-back multisection} and denote it by $\frak{back}$.
We denote its zero set by
$$
\mathcal M(k_1,k_2;\ell_1,\ell_2;\beta_1,\beta_2;\text{\bf p}_1,\text{\bf p}_2)^{\frak {back}}.
$$
\par
Lemma \ref{glueisolemann} immediately implies the
equality
\begin{equation}\label{forgetandfiber}
\aligned
&\int_{\frak{forget}^{-1}([\Sigma_1])
\cap \mathcal M_{(k_1+1,k_2+1);\ell}^{\text{\rm ann;main}}
(\beta;\text{\bf p})}(\text{\rm ev}^1,\text{\rm ev}^2)^*([\frak v;k_1] \times [\frak w;k_2]) \\
&=
\sum_{\beta_1+\beta_2=\beta}\sum_{(\mathbb L_1,\mathbb L_2) \in \text{\rm Shuff}(\ell)}
\\
&\qquad\qquad\qquad
\int_{\mathcal M(k_1+1,k_2+1;\ell_1,\ell_2;\beta_1,\beta_2;\text{\bf p}_1,
\text{\bf p}_2)^{\frak {back}}}(\text{\rm ev}^1,\text{\rm ev}^2)^*([\frak v;k_1] \times [\frak w;k_2]),
\endaligned
\end{equation}
where $(\text{\bf
p}_1,\text{\bf p}_2)
= \text{Split}((\mathbb L_1,\mathbb L_2),\text{\bf p})$.
\par
We next study the right hand side of
Proposition \ref{doublehatandpdX}.
Let $\text{\bf
p}_i \in E_{\ell_i}(\mathcal A[2])$. We assume that
they do not contain the degree zero element $\text{\bf f}_0$.
We fixed a continuous family of multisections on
$\mathcal M_{k_i;\ell_i+1}(\beta_i;\text{\bf p}_i)$ when we defined
an operator $\frak p$. It is chosen so that $\text{\rm ev}^{\rm int}_{\ell_i+1}$
 becomes a submersion. Therefore it induces a continuous family of multisections on
$\mathcal M(k_1,k_2;\ell_1,\ell_2;\beta_1,\beta_2;\text{\bf p}_1,\text{\bf p}_2)$.
We call this multisection the
{\it $\frak p$-multisection}\index{multisection (perturbation)!$\frak p$-multisection} and denote its zero set by
$\mathcal M(k_1,k_2;\ell_1,\ell_2;\beta_1,\beta_2;
\text{\bf p}_1,\text{\bf p}_2)^{\frak p}$.
Then the next formula is immediate from definition
\begin{equation}\label{pandfiberprod}
\aligned
&\int_{\mathcal M(k_1+1,k_2+1;\ell_1,\ell_2;\beta_1,\beta_2;\text{\bf p}_1,\text{\bf p}_2)^{\frak p}}
(\text{\rm ev}^1,\text{\rm ev}^2)^*([\frak v;k_1] \times [\frak w;k_2]) \\
&\qquad\qquad
=
\langle\frak p_{\ell_1,\beta_1}(\text{\bf p}_1
\otimes [\frak v;k_1]),
\frak p_{\ell_2,\beta_2}(\text{\bf p}_2
\otimes[\frak w;k_2])
\rangle_{\text{\rm PD}_{X}}.
\endaligned
\end{equation}
\begin{lem}\label{existsmultion01times}
There exists a system of continuous families of multisections $\frak s$ on the space
$[0,1] \times \mathcal M(k_1,k_2;\ell_1,\ell_2;\beta_1,\beta_2;\text{\bf p}_1,\text{\bf p}_2)$
with the following properties.
\begin{enumerate}
\item It is transversal to $0$.
\item It is invariant under the cyclic permutations of
$k_1$ (resp. $k_2$) boundary marked points
on each of the boundary components.
It is invariant under
arbitrary permutations of the interior marked points.
\item It is compatible with the forgetful map of the boundary marked points.
\item At $\{1\} \times \mathcal M(k_1,k_2;\ell_1,\ell_2;\beta_1,\beta_2;\text{\bf p}_1,\text{\bf p}_2)$ it becomes 
the fiber product of the $\frak p$-Kuranishi structures as in Definition \ref{twofiberproductofp}.
\item
At $\{0\} \times \mathcal M(k_1,k_2;\ell_1,\ell_2;\beta_1,\beta_2;\text{\bf p}_1,\text{\bf p}_2)$ it becomes $\frak{back}$.
\item It is compatible at the boundary components
\begin{enumerate}
\item
$$
\mathcal M^{\text{\rm main}}_{k'_1+1,\ell'_1}(\beta'_1;\text{\bf p}'_1)
{}_{\text{\rm ev}_0}^{\frak c} \times_{\text{\rm ev}^1_i}
([0,1]\times \mathcal M(k''_1,k_2;
\ell''_1,\ell_2;\beta''_1,\beta_2;\text{\bf p}''_1,\text{\bf p}_2))
$$
where $\ell'_1 + \ell''_1  = \ell_1$, $k'_1 + k''_1  = k_1+1$,
$\beta'_1 + \beta''_1  = \beta_1$, $(\text{\bf p}'_1,\text{\bf p}''_1) =
\text{\rm Split}(\text{\bf p}_1,(\mathbb L_1,\mathbb L_2))$,
$i = 1,\dots,k''_1$. We put $\frak c$ perturbation on the first factor,
and
\item $$
([0,1]\times \mathcal M(k''_1,k_2;
\ell''_1,\ell_2;\beta''_1,\beta_2;\text{\bf p}''_1,\text{\bf p}_2))
{}_{\text{\rm ev}^1_i} \times_{\text{\rm ev}_0}
\mathcal M^{\text{\rm main}}_{k'_1+1,\ell'_1}(\beta'_1;\text{\bf p}'_1)^{\frak c}
$$
where $\ell'_1 + \ell''_1  = \ell_1$, $k'_1 + k''_1  = k_1+1$,
$\beta'_1 + \beta''_1  = \beta_1$, $(\text{\bf p}'_1,\text{\bf p}''_1) =
\text{\rm Split}(\text{\bf p}_1,(\mathbb L_1,\mathbb L_2))$,
$i = 1,\dots,k''_1$. We put $\frak c$ perturbation on the second factor.
\end{enumerate}
\item It is compatible at the boundary components
\begin{enumerate}
\item
$$
\mathcal M_{k'_2+1,\ell'_2}^{\text{\rm main}}(\beta'_2;\text{\bf p}'_2)^{\frak c}{}
{}_{\text{\rm ev}_0} \times_{\text{\rm ev}^{2}_i}
([0,1]\times \mathcal M(k_1,k''_2;
\ell_1,\ell''_2;\beta_1,\beta''_2;\text{\bf p}_1,\text{\bf p}''_2))
$$
where $\ell'_2 + \ell''_2  = \ell_2$, $k'_2 + k''_2  = k_2+1$,
$\beta'_2 + \beta''_2  = \beta_2$, $(\text{\bf p}'_2,\text{\bf p}''_2) =
\text{\rm Split}(\text{\bf p}_2,(\mathbb L_1,\mathbb L_2))$,
$i = 1,\dots,k''_2$. We put $\frak c$ perturbation on the first factor, and
\item
$$
([0,1]\times \mathcal M(k_1,k''_2;
\ell_1,\ell''_2;\beta_1,\beta''_2;\text{\bf p}_1,\text{\bf p}''_2))
{}_{\text{\rm ev}^{2}_i} \times_{\text{\rm ev}_0}
\mathcal M_{k'_2+1,\ell'_2}^{\text{\rm main}}(\beta'_2;\text{\bf p}'_2)^{\frak c}
$$
where $\ell'_2 + \ell''_2  = \ell_2$, $k'_2 + k''_2  = k_2+1$,
$\beta'_2 + \beta''_2  = \beta_2$, $(\text{\bf p}'_2,\text{\bf p}''_2) =
\text{\rm Split}(\text{\bf p}_2,(\mathbb L_1,\mathbb L_2))$,
$i = 1,\dots,k''_2$. We put $\frak c$ perturbation on the second factor.
\end{enumerate}
\end{enumerate}
\end{lem}
\begin{proof}
The proof can be done in the same way as we discuss in several other places.
See Section \ref{sec:equikuracot}, Subsections \ref{subsec:lemma426}--\ref{subsec:kuraanul3}.
\end{proof}
We use Lemma \ref{existsmultion01times} to prove the next lemma.
\begin{lem}\label{backisp}
$$
\aligned
&\sum_{\beta_1,\beta_2}\sum_{\ell_1+\ell_2=\ell}\sum_{k_1,k_2}
\frac{T^{(\beta_1+\beta_2)\cap \omega/2\pi}}{\ell_1!\ell_2!}
\exp((\beta_1+\beta_2)\cap \frak b_2)\rho(\partial (\beta_1+\beta_2))
\\
&\qquad \qquad
\int_{(\mathcal M(k_1,k_2;\ell_1,\ell_2;\beta_1,\beta_2;\frak b_{\rm high}^{\otimes\ell_1},\frak b_{\rm high}^{\otimes\ell_2})^{\frak p}}
(\text{\rm ev}^1,\text{\rm ev}^2)^*([\frak v;k_1]\times [\frak w;k_2]) \\
&=
\sum_{\beta_1,\beta_2}\sum_{\ell_1+\ell_2=\ell}\sum_{k_1,k_2}
\frac{T^{(\beta_1+\beta_2)\cap \omega/2\pi}}{\ell_1!\ell_2!}
\exp((\beta_1+\beta_2)\cap \frak b_2)\rho(\partial (\beta_1+\beta_2))
\\
&\qquad \qquad
\int_{(\mathcal M(k_1,k_2;\ell_1,\ell_2;\beta_1,\beta_2;\frak b_{\rm high}^{\otimes\ell_1},\frak b_{\rm high}^{\otimes\ell_2})^{\frak {back}}}
(\text{\rm ev}^1,\text{\rm ev}^2)^*([\frak v;k_1]\times [\frak w;k_2]).
\endaligned
$$
\end{lem}
\begin{proof}
We put
\begin{equation}\label{betasplit2}
\frak b_{\rm high}^{\otimes\ell} = \sum_{I \in \frak I_{\ell}} c_I \text{\bf f}_{i_1}
\otimes \dots \otimes \text{\bf f}_{i_{\ell}}
= \sum_{I \in \frak I_{\ell}} c_{I,\ell} \text{\bf f}_{I,\ell}.
\end{equation}
We apply Stokes' theorem (\cite[Lemma 12.18]{fooo09}) on
\begin{equation}\label{zeroplusookii}
([0,1] \times 
{\mathcal M}(k_1+1,k_2+1;\ell_1,\ell_2;\beta_1,\beta_2;
\text{\bf f}_{I_1,{\ell_1}},\text{\bf f}_{I_2,{\ell_2}}))^{\frak s}
\end{equation}
to the pull-back $(\text{\rm ev}^1,\text{\rm ev}^2)^*([\frak v;k_1]\times [\frak w;k_2])$,
which is a closed form.
\par
Let us study the integration of $(\text{\rm ev}^1,\text{\rm ev}^2)^*([\frak v;k_1]\times [\frak w;k_2])$
on each component of the boundary of (\ref{zeroplusookii}).
\par
The case of Lemma \ref{existsmultion01times}.4 and 5 give the left and right hand
sides of the formula in Lemma \ref{backisp}, respectively after taking weighted sum with
weight
\begin{equation}\label{thisisweight}
\frac{T^{\beta\cap \omega/2\pi}}{\ell_1!\ell_2!}
\exp(\beta\cap \frak b_2)\rho(\partial \beta).
\end{equation}
\par
Let us study Lemma \ref{existsmultion01times}.6.
We write
$$
[\frak v;k_1]
=[\frak v\otimes (b_+^{\frak c})^{k_1}] \in B_{k_1+1}^{\rm cyc}(H(L(u));\Lambda_0)[1]).
$$
(Note we did not break cyclic symmetry in Lemmata \ref{perturbnbdsigma1},
\ref{existsmultion01times}.)
We integrate it on the boundary component described by Lemma \ref{existsmultion01times}.6. We divide it into two cases.
\par
(Case 1)
 Lemma \ref{existsmultion01times}.6
(a)
This is the case when the boundary marked point corresponding to $\frak v$
is on the second factor.
\par
(Case 2)   Lemma \ref{existsmultion01times}.6
(b).
This is the case when the boundary marked point corresponding to $\frak v$
is on the first factor.
\par
We can show that contribution of both cases vanishes in the same way
as in the proof of Theorem \ref{pqpoincare}.
Namely we can study (Case 1) in the same way as Lemma \ref{bdrydescription}.3
and (Case 2) in the same way as Lemma \ref{bdrydescription}.4.
We repeat the detail below for completeness' sake.
\par
We first study Case 1.  We define $E$ by (\ref{Edefine}).
Then the contribution of Case 1 is
$$\aligned
&\sum_{k_1,k_2,
\atop \ell_1,\ell_2}\sum_{I_1\in\frak I_{\ell_1} \atop I_2\in\frak I_{\ell_2}}\sum_{i=1}^{k_1+1}
\sum_{\beta_1,\beta_2}
 T^{(\beta_1+\beta_2)\cap \omega/2\pi} \exp ((\beta_1+\beta_2)\cap \frak b_2)
 \frac{c_{I_1}c_{I_2}}{\ell_1!\ell_2!} \rho(\partial(\beta_1+\beta_2)) \\
&\int_{([0,1] \times \mathcal M(k_1,k_2;\ell_1;\ell_2;\beta_1;\beta_2;
\text{\bf p}_1,\text{\bf p}_2))^{\frak s}}
(\text{\rm ev}_1,\text{\rm ev}_2)^*([\frak v  \times (b_+^{\frak c})^{i-1}\times E \times (b_+^{\frak c})^{k_1+1-i}],[\frak w;k_2]).
\endaligned$$
In the same way as the proof of Lemma \ref{unitprop}, we can prove that
$E$ is proportional to the unit. Therefore this contribution
vanishes by Lemma \ref{existsmultion01times}.3.
\par
We next study Case 2.
We define $R$ by
\begin{equation}\aligned
R = &
\sum_{k_1,k_2,
\atop \ell_1,\ell_2}\sum_{I_1\in\frak I_{\ell_1} \atop I_2\in\frak I_{\ell_2}}
\sum_{\beta_1,\beta_2}
T^{\beta\cap\omega/2\pi}\exp(\beta\cap \frak b_2) \frac{c_I}{\ell!}
\rho(\partial\beta) \\
&(\text{\rm ev}_0)_*
\big(
((\text{\rm ev}_1,\text{\rm ev}_2,\dots,\text{\rm ev}_{k_1}),\text{\rm ev}_2)^*((b^{\frak c}_+)^{k_1}
,[\frak w;k_2]); \\
&\qquad\quad([0,1] \times \mathcal M(k_1+1,k_2;\ell_1;\ell_2;\beta_1;\beta_2;
\text{\bf p}_1,\text{\bf p}_2))^{\frak s}
\big).
\endaligned
\end{equation}
By the definition of $R$, the contribution of Case 2 is $\langle\delta^{\frak b,b^{\frak c}}(R),\frak v\rangle = 0$.  
By the cyclic symmetry, we have 
$$
\langle\delta^{\frak b,b^{\frak c}}(R),\frak v\rangle = \pm \langle R, \delta^{\frak b,b^{\frak c}} (\frak v) \rangle= 0.
$$
Note that $R$ may not be $T^n$ invariant but $\frak v$ is $T^n$ 
invariant.
We can show that the contribution of Lemma \ref{existsmultion01times}.7
vanishes in the same way.
The proof of Lemma \ref{backisp} is now complete.
\end{proof}
Proposition \ref{doublehatandpdX} follows from Formulae (\ref{forgetandfiber}),
(\ref{pandfiberprod}) and Lemma \ref{backisp}.
\end{proof}
Before proceeding to the proofs of Propositions \ref{doublehatistwodisk} and
\ref{twodiskandm2}, we fix the choice of perturbation on
$\frak{forget}^{-1}([\Sigma_2])
\cap \mathcal M_{(k_1+1,k_2+1);\ell}^{\text{\rm ann;main}}
(\beta;\frak b_{\rm high}^{\otimes\ell})$.
\par
Let $(\text{\bf p}_1,\text{\bf p}_2)  =
\text{\rm Split}(\text{\bf p},(\mathbb L_1,\mathbb L_2))$
and
$
\ell_j = \vert \text{\bf p}_j \vert$,
$j =1,2.$
\par
We consider $\mathcal M^{\text{\rm main}}_{k_j+2;\ell_j}(\beta_j;
\text{\bf p}_j)$.
Let $a_j \in \{1,\dots,k_{j}+1\}$. We consider
\begin{equation}\label{ddoublefiber1}
(\text{\rm ev}_0,\text{\rm ev}_{a_1}):
\mathcal M^{\text{\rm main}}_{k_1+2;\ell_1}(\beta_1;
\text{\bf p}_1)
\to L(u) \times L(u)
\end{equation}
and
\begin{equation}\label{ddoublefiber2}
(\text{\rm ev}_{a_2},\text{\rm ev}_0):
\mathcal M^{\text{\rm main}}_{k_2+2;\ell_2}(\beta_2;
\text{\bf p}_2)
\to L(u) \times L(u).
\end{equation}
It is easy to see that $\frak{forget}^{-1}([\Sigma_2])
\cap \mathcal M_{(k_1+1,k_2+1);\ell}^{\text{\rm ann;main}}
(\beta;\frak b_{\rm high}^{\otimes\ell})$ 
is the fiber product of the
two maps $(\ref{ddoublefiber1})$ and
$(\ref{ddoublefiber2})$.
However, these two maps are not
transversal.

\begin{exm}
We consider the case $k_1=k_2=1$, $\ell_1=\ell_2 = 0$, 
and $\beta_1 = \beta_2 = 0 = \beta_0$.
Note $\mathcal M^{\text{\rm main}}_{3;0}(\beta_0)
= L(u)$.
The obstruction bundle of $\frak c$-perturbation is trivial.
The maps $(\ref{ddoublefiber1})$ and
$(\ref{ddoublefiber2})$ 
are embedding $x \mapsto (x,x)$.
So the fiber product is not transversal on the Kuranishi neighborhood of $\frak c$-perturbation.
\par
We note that
$\frak{forget}^{-1}([\Sigma_2])
\cap \mathcal M_{(1,1);0}^{\text{\rm ann;main}}
(\beta_0) = L(u)$, but this moduli space is obstructed.
The virtual dimension of this moduli space is $0$, hence its virtual fundamental class is a rational number.  
By an argument similar to the one in the proof of Theorem \ref{annulusmain} given in this section, this number is equal to 
the intersection number of $[L(u)]$ with itself in $X$
and is zero.
\end{exm}
Therefore we need to take a different perturbation from 
$\frak c$-perturbation.
For $k_1,k_2,k_3 \ge 0$ we denote 
$$
\mathcal M^{\text{\rm main}}_{k_1,k_2,k_3;\ell}(\beta;
\text{\bf p})
:= 
\mathcal M^{\text{\rm main}}_{k_1+k_2+k_3+3;\ell}(\beta;
\text{\bf p}).
$$
We distinguish $0$-th, $(k_1+1)$-th and $(k_1+k_2+2)$-th (boudnary) marked 
points from other marked points. 
We call those three marked points 
the {\it unforgetable marked points}\index{unforgetable marked points}  
and the other marked points
{\it forgetable marked points}\index{forgetable marked points}.
We write the evaluation map as
\begin{equation}\label{form316}
\aligned
{\rm ev}= &({\rm ev}^0_0,{\rm ev}^0_1,{\rm ev}^0_2;{\rm ev}_1,{\rm ev}_2,{\rm ev}_3) : \\
&\mathcal M^{\text{\rm main}}_{k_1,k_2,k_3;\ell}(\beta;
\text{\bf p})
\to L(u) ^3 \times L(u)^{k_1} \times L(u)^{k_2} \times L(u)^{k_3},
\endaligned
\end{equation}
where ${\rm ev}^0_0,{\rm ev}^0_1,{\rm ev}^0_3$ are the evaluation maps at unforgetable marked points.
We write ${\rm ev}_i = ({\rm ev}_{i,1},\dots,{\rm ev}_{i,k_i})$.
\par
Using the terminology `unforgetable', we rewrite Convention
\ref{convint} as follows:
\begin{conv}\label{convint2}
When we integrate the pull back of $[\frak v;k_1]$ or $[\frak w;k_2]$
we alway pull back the forms $\frak v$, $\frak w$ by the evaluation maps 
at the unforgetable marked points and $\frak b_+^{\frak c}$ by the 
evaluation map at the forgetable marked points.
\end{conv}
 
\begin{lem}\label{existstkura}
There exists a system of Kuranishi structures and multisections on  
$\mathcal M^{\text{\rm main}}_{k_1,k_2,k_3;\ell}(\beta;
\text{\bf p})$ with the following properties.
We call this Kuranishi structure and multisection the 
{\it $\frak t$-Kuranishi structure}\index{Kuranishi structure!$\frak t$-Kuranishi structure} and
{\it $\frak t$-multisection}\index{multisection (perturbation)!$\frak t$-multisection}, respectively.
\begin{enumerate}
\item
They are transversal to $0$ and are $T^n$ equivariant.
\item
They are compatible with the forgetful map of the forgetable boundary marked points.
\item
They are invariant under the permutation of the interior marked points together with the 
permutation of ${\bf p}_i$.
\item
The evaluation map
$$
 ({\rm ev}^0_0,{\rm ev}^0_1,{\rm ev}^0_2) :\mathcal M^{\text{\rm main}}_{k_1,k_2,k_3;\ell}(\beta;
\text{\bf p})^{\frak t}
\to L(u)^3
$$
at the unforgetable marked points is a submersion on the zero set of 
$\frak t$-multisection.
\item
The boundary of $\mathcal M^{\text{\rm main}}_{k_1,k_2,k_3;\ell}(\beta;
\text{\bf p})^{\frak t}$ is described as the union of the following two types 
of fiber products.
\begin{enumerate}
\item
\begin{equation}
\mathcal M^{\text{\rm main}}_{k'_i+k'_{i+1}+2;\ell'}(\beta_1;
\text{\bf p}_1)^{\frak c}
\,{}_{{\rm ev}_0}\times_{{\rm ev}^0_i}
\mathcal M^{\text{\rm main}}_{k''_1,k''_2,k''_3;\ell''}(\beta_2;
\text{\bf p}_2)^{\frak t}.
\end{equation}
Here $i$ is one of $1$,$2$,$3$. ($i+1$ means $1$ if $i=3$.)
Moreover $k'_i + k''_{i} = k_i$, $k'_{i+1} + k''_{i+1} = k_{i+1}$
and $k''_j = k_j$ if $j\ne i,i+1$.
$\beta = \beta_1+\beta_2$, $\ell = \ell_1 + \ell_2$ and $\text{\bf p} = (\text{\bf p},\text{\bf p}_2)$.
\par
The $i$-th unforgetable marked point of $\mathcal M^{\text{\rm main}}_{k_1,k_2,k_3;\ell}(\beta;
\text{\bf p})$ becomes $k'_i+1$-th marked point of the first factor.
The $j$-th unforgetable marked point of $\mathcal M^{\text{\rm main}}_{k_1,k_2,k_3;\ell}(\beta;
\text{\bf p})$ 
becomes $j$-th unforgetable marked point of the second factor.
\par
We note that we take $\frak c$-multisection for the first factor 
and $\frak t$-multisection for the second factor.
(See Figure 3.4.1.)
\item
\begin{equation}
\mathcal M^{\text{\rm main}}_{k'_i+1;\ell'}(\beta_1;
\text{\bf p}_1)^{\frak c}
\,{}_{{\rm ev}_0}\times_{{\rm ev}_{i,j}}
\mathcal M^{\text{\rm main}}_{k''_1,k''_2,k''_3;\ell''}(\beta_2;
\text{\bf p}_2)^{\frak t}.
\end{equation}
Here $i$ is one of $1$,$2$,$3$. ($i+1$ means $1$ if $i=3$.)
Moreover $k'_i + k''_{i} -1 = k_i$, 
and $k''_j = k_j$ if $j\ne i$. $j \in \{1,\dots,k''_i\}$, $k''_i \ge 1$. 
${\rm ev}_i = ({\rm ev}_{i,1},\dots,{\rm ev}_{i,j},\dots,{\rm ev}_{i,k''_i})$ where ${\rm ev}_i$ is as in 
(\ref{form316}).
\par
$\beta = \beta_1+\beta_2$, $\ell = \ell_1 + \ell_2$ and $\text{\bf p} = (\text{\bf p},\text{\bf p}_2)$.
\par
The $i$-th unforgetable marked point of $\mathcal M^{\text{\rm main}}_{k_1,k_2,k_3;\ell}(\beta;
\text{\bf p})$ becomes the $i$-th unforgetable marked point of the second factor.
\par
We note that we take $\frak c$-multisection for the first factor 
and $\frak t$-multisection for the second factor.
(See Figure 3.4.1.)
\end{enumerate}
\end{enumerate}
\end{lem}
\par
\epsfbox{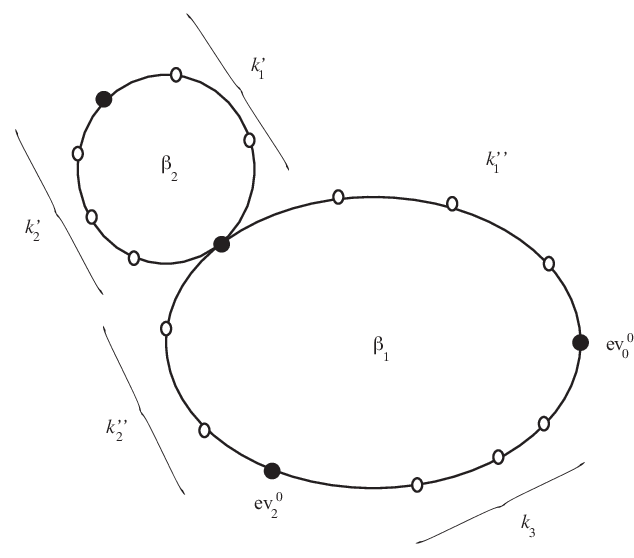}
\par
\centerline{\bf Figure 3.4.1}
\par
The proof is given in Subsection \ref{subsec:tkuranishi}.
Hereafter we write $\vec k = (k_1,k_2,k_3)$ etc..
We consider 
\begin{equation}\label{ddoublefiber12}
(\text{\rm ev}^0_1,\text{\rm ev}^0_{2}):
\mathcal M^{\text{\rm main}}_{\vec k^{(1)};\ell_1}(\beta_1;
\text{\bf p}_1)^{\frak t}
\to L(u) \times L(u)
\end{equation}
and
\begin{equation}\label{ddoublefiber22}
(\text{\rm ev}^0_1,\text{\rm ev}^0_{2}):
\mathcal M^{\text{\rm main}}_{\vec k^{(2)};\ell_2}(\beta_2;
\text{\bf p}_2)^{\frak t}
\to L(u) \times L(u)
\end{equation}
\begin{defn}\label{gluedmoduliqann}
We denote by
$
\mathcal N(\vec k^{(1)},\vec k^{(2)};\ell_{1},\ell_{2};
\beta_1,\beta_2;\text{\bf p}_1;\text{\bf p}_2)
$\index{$\mathcal N(\vec k^{(1)},\vec k^{(2)};\ell_{1},\ell_{2};
\beta_1,\beta_2;\text{\bf p}_1;\text{\bf p}_2)$}
the fiber product of the maps $(\ref{ddoublefiber12})$ and
$(\ref{ddoublefiber22})$.
\par
There is an obvious map (see Figure 3.4.2)
\begin{equation}\label{Ngluingmap}
\aligned
&\text{Glue}: \mathcal N(\vec k^{(1)},\vec k^{(2)};\ell_{1},\ell_{2};
\beta_1,\beta_2;\text{\bf p}_1;\text{\bf p}_2)
\\
&\to
\mathcal M_{(k^{(1)}_2+ k^{(2)}_1 + k^{(2)}_3+1,k^{(2)}_2+ k^{(1)}_1 + k^{(1)}_3+1);\ell_1+\ell_2}^{\text{ann;main}}(\beta_1+\beta_2;
\text{\bf p})\cap \frak{forget}^{-1}([\Sigma_2]).
\endaligned\end{equation}
\end{defn}
\par
\epsfbox{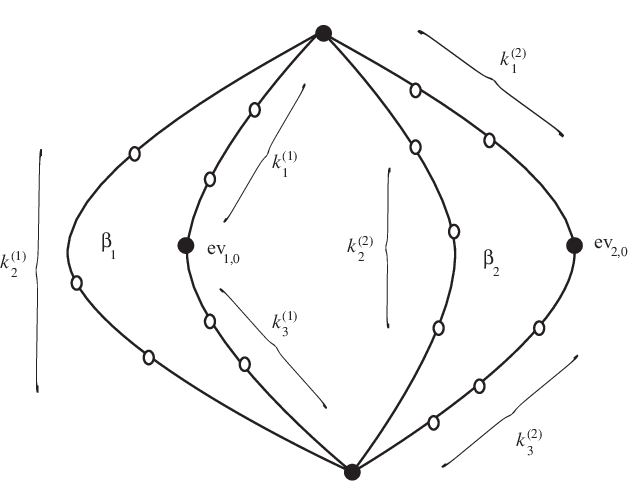}
\par
\centerline{\bf Figure 3.4.2}
\par
We denote the composition 
$$
{\rm ev}_1 \circ \text{Glue}:
\mathcal N(\vec k^{(1)},\vec k^{(2)};\ell_{1},\ell_{2};
\beta_1,\beta_2;\text{\bf p}_1;\text{\bf p}_2)\to L(u)^{k_1+1}
$$
by ${\rm ev}_1$ also. The notation ${\rm ev}_2$ is used in the same way.
\begin{lem}\label{doublegluelem}
The map $\text{\rm Glue}$ defines an isomorphism onto the image.
The moduli space
$\mathcal M_{(k_1+1,k_2+1);\ell}^{\text{\rm ann;main}}(\beta;
\text{\bf p}) \cap \frak{forget}^{-1}([\Sigma_2])$ is decomposed into the
images of $\text{\rm Glue}$'s from difference sources.
If $\text{\bf x} \in \mathcal M_{(k_1+1,k_2+1);\ell}^{\text{\rm ann;main}}(\beta;
\text{\bf p}) \cap \frak{forget}^{-1}([\Sigma_2])$ is an
image of a codimension $N$ corner point of
$\mathcal N(\vec k^{(1)},\vec k^{(2)};\ell_{1},\ell_{2};
\beta_1,\beta_2;\text{\bf p}_1;\text{\bf p}_2)$,
then it is an image of $N$ points from the union
of the sources of $\text{\rm Glue}$.
\par
The $\frak t$-Kuranishi structures of the sources
are compatible at the overlapped part of the target.
\par
Furthermore, the $\frak t$-multisections  are also compatible.
\end{lem}
\begin{proof}
In fact, the proof is the same as those of Lemmata \ref{cornergluelocal},
\ref{perturbconsistence}.
We however emphasis that cyclic symmetry of $\frak c$-perturbation is crucial for this lemma to hold.
To elaborate this point, let us consider the subset consisting of elements
$((D^2,u_1),(D^2,u_2),(D^2,u_3))$ of
$\mathcal M_3^{\text{\rm main}}(\beta_1) \times \mathcal M_2^{\text{\rm main}}(\beta_2) \times
\mathcal M_3^{\text{\rm main}}(\beta_3)$ such that
\begin{equation}
\text{\rm ev}_1(u_1) = \text{\rm ev}_1(u_2), \quad
\text{\rm ev}_2(u_1) = \text{\rm ev}_2(u_3), \quad
\text{\rm ev}_0(u_2) = \text{\rm ev}_1(u_3).
\end{equation}
\par
\hskip0.3cm
\epsfbox{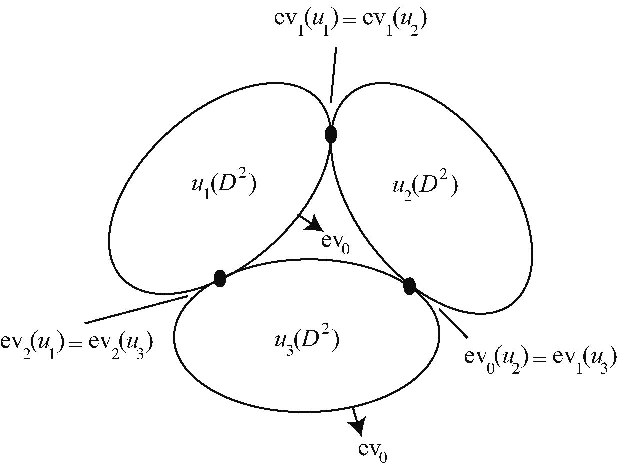}
\par
\centerline{\bf Figure 3.4.3}
\par
See Figure 3.4.3. This set is a subset of
$\mathcal M^{\text{ann;main}}_{(1,1);0}(\beta_1+\beta_2+\beta_3)
\cap \frak{forget}^{-1}(\Sigma_2)
$ and lies in its codimension $2$ strata.
\par
We consider two (among three) ways to resolve singularities.
Namely at $\text{\rm ev}_{1}(u_1) = \text{\rm ev}_{1}(u_2)$ or at $\text{\rm ev}_0(u_2) = \text{\rm ev}_1(u_3)$.
\par
We note that for $((D^2,u_1)$, $(D^2,u_3))$ the choice of the $0$-th marked
point (among the three unforgetable boundary marked points) is canonical.
Namely there is a marked point which does not meet other irreducible
component. We take it as the $0$-th one.
\par
On the other hand, there is no good way to decide which (among two)
boundary marked point is the $0$-th one for $(D^2,u_2)$. 
In fact, if we insist on the consistency with the process to resolve $\text{\rm ev}_1(u_1) = \text{\rm ev}_1(u_2)$,
the marked point of $(D^2,u_2)$ that intersects with $(D^2,u_1)$
must be the $0$-th one.
If we insist on the consistency with the process to resolve $\text{\rm ev}_0(u_2) = \text{\rm ev}_1(u_3)$,
the marked point of $(D^2,u_2)$ that intersects with $(D^2,u_3)$
must be the $0$-th one.
\par
Thus for the consistency of both, we need to take our perturbation of
$\mathcal M_2^{\text{\rm main}}(\beta_2)$ to be invariant under cyclic
symmetry.
Actually this is one of the reasons we take the $\frak c$-perturbation 
on $\mathcal M_2^{\text{\rm main}}(\beta_2)$.
Note we use the $\frak t$-perturbation on $\mathcal M_3^{\text{\rm main}}(\beta_1)$
and $\mathcal M_3^{\text{\rm main}}(\beta_3)$.
(See however \cite{abouz}, where M.  Abouzaid uses less symmetry than ours.
We use cyclic symmetry also in Sections  \ref{sec:operatorptoric}
and \ref{sec:clifford}. See Remarks \ref{cycandPD} and \ref{rem224}.)
\par
Once this point is understood, the proof is the same as 
those of Lemmata \ref{cornergluelocal},
\ref{perturbconsistence}.
\end{proof}
By Lemma \ref{doublegluelem} we can define
a system of continuous families of multisections $\frak s$
on $\mathcal M_{(k_1,k_2);\ell}^{\text{ann;main}}(\beta;
\text{\bf p}) \cap \frak{forget}^{-1}([\Sigma_2])$  so that
it coincides with the fiber product of $\frak t$-multisection on
each factor of
$\mathcal N(\vec k^{(1)},\vec k^{(2)};\ell_{1},\ell_{2};
\beta_1,\beta_2;\text{\bf p}_1;\text{\bf p}_2)$.
\par
We have thus fixed our family of multisections for
$\mathcal M_{(k_1,k_2);\ell}^{\text{ann;main}}(\beta;
\text{\bf p}) \cap \frak{forget}^{-1}([\Sigma_2])$.
It has the following properties.
\begin{conds}\label{cond21}
\begin{enumerate}
\item It is transversal to $0$ and is $T^n$-equivariant.
\item It is invariant under 
arbitrary permutations of the interior marked points.
\item It is compatible with the forgetful map of the boundary marked points.
\item It is transversal to the map $(\ref{Ngluingmap})$ in a sense similar to Lemma \ref{perturbnbdsigmasec9}.5.
\item It is compatible at the boundary components.
When there is a disk bubble(s), the multisection is a fiber product.
Here we take $\frak c$-perturbation on the disk bubble.
\item On $\mathcal M_{(k_1+1,k_2+1);\ell}^{\text{ann;main}}(\beta;
\text{\bf p}) \cap \frak{forget}^{-1}([\Sigma_2])$ it coincides with the fiber product of $\frak t$-multisections.
\end{enumerate}
\par
The boundary components of various irreducible components of
$\mathcal M_{(k_1,k_2);\ell}^{\text{ann;main}}(\beta;
\text{\bf p}) \cap \frak{forget}^{-1}([\Sigma_2])$ are divided into 6 cases as in Figure 3.4.4 below.
\end{conds}
\begin{rem}
We take $\frak c$-perturbation on the disk bubble in Condition \ref{cond21}. 
So because of the cyclic symmetry of $\frak c$-perturbation it does not matter which marked points 
we regard as the $0$-th one. 
The same remark applies to the disk components appearing in Condition \ref{cond21}.6.
\end{rem}
\begin{rem}
Note that we break cyclic symmetry of the boundary marked points of the annulus.
In fact, among $k_1 + 1$ and $k_2 +1$ boundary marked points of the two boundary 
components, $k_1$ and $k_2$ are forgetable and $1$ and $1$ are unforgetable.
Those two kinds of the boundary marked points are not symmetric.
\end{rem}
\hskip-0.5cm
\epsfbox{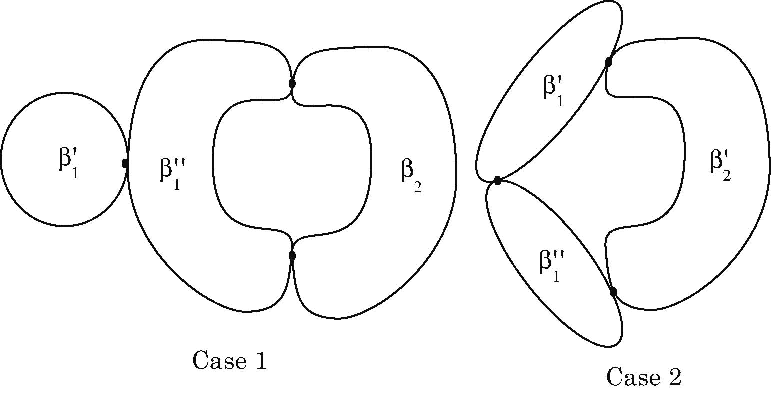}
\par
\hskip0.7cm
\epsfbox{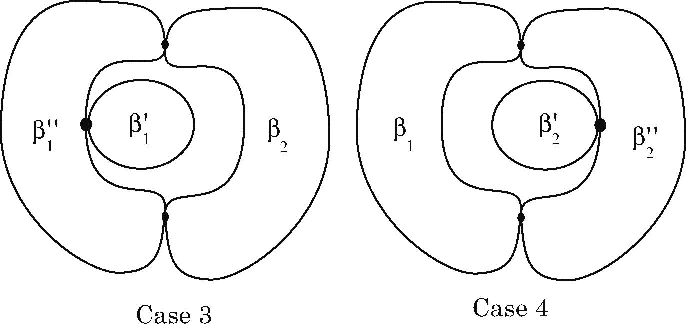}
\par
\hskip-0.5cm
\epsfbox{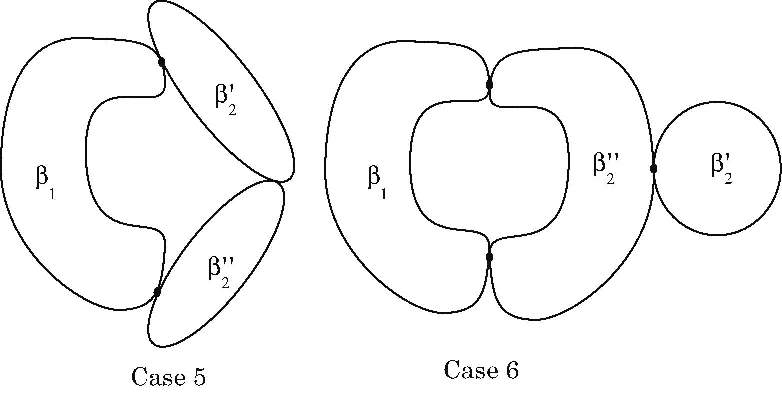}
\par
\centerline{\bf Figure 3.4.4}
\par
\begin{proof}[Proof of Proposition \ref{doublehatistwodisk}]
We extend the above multisection to a neighborhood of
$\mathcal M_{(k_1,k_2);\ell}^{\text{ann;main}}(\beta;
\text{\bf p}) \cap \frak{forget}^{-1}([\Sigma_2])$ in
$\mathcal M_{(k_1,k_2);\ell}^{\text{ann;main}}(\beta;
\text{\bf p})$. Then
for each $[\Sigma'] \in \mathcal M_{(1,1);0}^{\text{ann;main}}$ sufficiently close to $[\Sigma_2]$
we obtain a multisection on $\mathcal M_{(k_1,k_2);\ell}^{\text{ann;main}}(\beta;
\text{\bf p}) \cap \frak{forget}^{-1}([\Sigma'])$ which
satisfies the same conditions as Condition \ref{cond21} above.
We choose a sequence $[\Sigma'_i] \in \mathcal M_{(1,1);0}^{\text{ann;main}}$
which is in stratum 1 of Lemma \ref{stratifyAn11}
and
$$
\lim_{i\to\infty}[\Sigma'_i] = [\Sigma_2].
$$
\begin{lem}\label{limitinterior}
$$\aligned
&\lim_{i\to\infty} \int_{\mathcal M_{(k_1+1,k_2+1);\ell}^{\text{\rm ann;main}}(\beta;
\text{\bf p}) \cap \frak{forget}^{-1}([\Sigma'_i])}(\text{\rm ev}^1,\text{\rm ev}^2)^*(
[\frak v;k_1]
,[\frak w;k_2]) \\
&=
\int_{\mathcal M_{(k_1+1,k_2+1);\ell}^{\text{\rm ann;main}}(\beta;
\text{\bf p}) \cap \frak{forget}^{-1}([\Sigma_2])}(\text{\rm ev}^1,\text{\rm ev}^2)^*(
[\frak v;k_1]
,[\frak w;k_2]).
\endaligned$$
\end{lem}
Note we apply Convention \ref{convint} here.
\begin{proof}
This lemma is slightly nontrivial since $\frak{forget}$ is not smooth at the fiber of
$[\Sigma_2]$.
Let us decompose
$\mathcal M_{(k_1+1,k_2+1);\ell}^{\text{ann;main}}(\beta;\text{\bf p})
\cap \frak{forget}^{-1}([\Sigma_2]))$
to the union of
$\mathcal N(\vec k^{(1)},\vec k^{(2)};\ell_{1},\ell_{2};
\beta_1,\beta_2;\text{\bf p}_1;\text{\bf p}_2)$ by Lemma \ref{doublegluelem}.
We take an increasing sequence of compact subsets
$$
\mathcal N_K(\vec k^{(1)},\vec k^{(2)};\ell_{1},\ell_{2};
\beta_1,\beta_2;\text{\bf p}_1;\text{\bf p}_2), \quad K = 1,2, \ldots, \to \infty
$$
of the union
$\text{Int}\,\,
\mathcal N(\vec k^{(1)},\vec k^{(2)};\ell_{1},\ell_{2};
\beta_1,\beta_2;\text{\bf p}_1;\text{\bf p}_2)$
of the interior of its irreducible components such that
$$\aligned
\bigcup_{K}
\mathcal N_K(\vec k^{(1)},\vec k^{(2)};\ell_{1},\ell_{2};
\beta_1,\beta_2;\text{\bf p}_1;\text{\bf p}_2) 
= \text{Int}\,\,
\mathcal N(\vec k^{(1)},\vec k^{(2)};\ell_{1},\ell_{2};
\beta_1,\beta_2;\text{\bf p}_1;\text{\bf p}_2).
\endaligned$$
We consider the irreducible components of dimension
$\deg [\frak v;k_1] + \deg [w;k_2]$ only.
For each $K$ we can construct a map
$$\aligned
\Phi_{K,i}: \mathcal N(\vec k^{(1)},\vec k^{(2)};\ell_{1},\ell_{2};
\beta_1,\beta_2;\text{\bf p}_1;\text{\bf p}_2) 
\to \mathcal M_{(k_1+1,k_2+1);\ell}^{\text{\rm ann;main}}(\beta;
\text{\bf p}) \cap \frak{forget}^{-1}([\Sigma'_i]),
\endaligned$$
where $k_1 = k^{(1)}_2+ k^{(2)}_1 + k^{(2)}_3$,$k_2 =  k^{(2)}_2+ k^{(1)}_1 + k^{(1)}_3$,
such that
$$
\lim_{i\to\infty} \Phi_{K,i}^*(\text{\rm ev}^1,\text{\rm ev}^2)^*([\frak v;k_1]
, [\frak w;k_2])
= (\text{\rm ev}^1,\text{\rm ev}^2)^*([\frak v;k_1]
, [\frak w;k_2]).
$$
By transversality and definition, the contribution to the
integration of the right hand side of Lemma \ref{limitinterior}
of the codimension one stratum is $0$.
Therefore the integration of
$(\text{\rm ev}^1,\text{\rm ev}^2)^*([\frak v;k_1]
, [\frak w;k_2])$
on the union of the images of $\Phi_{K,i}$
can be chosen arbitrary close to the right hand side
if we take $K, i$ large.
\par
On the other hand, the complement of the union
of the images of $\Phi_{K,i}$ are obtained by
gluing at least three pieces of disks.
(Those configurations are described as in Figure 3.4.5.)
So there are at least 3 parameters $T_1,T_2,T_3$ for gluing.
When we restrict these 3 parameters by requiring that the
resulting source becomes the given $\Sigma'_i$, we
still have one more parameter $T$.
We can estimate the differential of $(\text{\rm ev}^1,\text{\rm ev}^2)$
with respect to $T$ in the same way as in  \cite[Lemma A1.58]{fooo06}
and find that it is of order $e^{-cT}$.
Therefore the integration of $(\text{\rm ev}^1,\text{\rm ev}^2)^*([\frak v;k_1]
, [\frak w;k_2])$
on the complement of the union of images of $\Phi_{K,i}$
goes to zero as $K,i \to \infty$.
Hence the lemma.
The argument here is similar to the proof of Lemma  \ref{sublemma:Sigma3limit},
which we discuss in detail in Section \ref{sec:cornerestimate}.
\par
\hskip1.8cm
\epsfbox{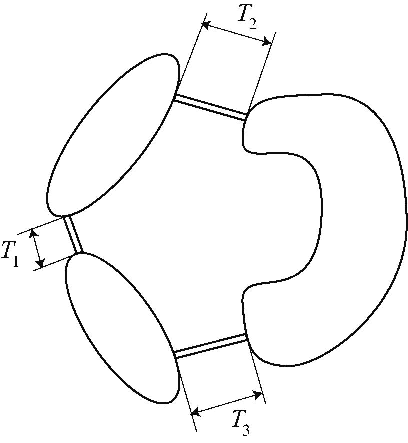}
\par
\centerline{\bf Figure 3.4.5}
\par
\end{proof}
We next choose a sequence $[\Sigma''_i]\in \mathcal M_{(1,1);0}^{\text{ann;main}}$
which is in stratum 1 in Lemma \ref{stratifyAn11}
and
$$
\lim_{i\to\infty}[\Sigma''_i] = [\Sigma_1].
$$
We take a perturbation of
$\mathcal M_{(1,1);\ell}^{\text{\rm ann;main}}(\beta;
\text{\bf p}) \cap \frak{forget}^{-1}([\Sigma''_i])$
by restricting the perturbation given in Lemma \ref{perturbnbdsigma1}.
\begin{lem}\label{limitinterior2}
$$\aligned
&\lim_{i\to\infty} \int_{\mathcal M_{(k_1+1,k_2+1);\ell}^{\text{\rm ann;main}}(\beta;
\text{\bf p}) \cap \frak{forget}^{-1}([\Sigma''_i])}(\text{\rm ev}^1,\text{\rm ev}^2)^*([\frak v;k_1]
, [\frak w;k_2]) \\
&=
\int_{\mathcal M_{(k_1+1,k_2+1);\ell}^{\text{\rm ann;main}}(\beta;
\text{\bf p}) \cap \frak{forget}^{-1}([\Sigma_1])}(\text{\rm ev}^1,\text{\rm ev}^2)^*([\frak v;k_1]
, [\frak w;k_2]).
\endaligned$$
\end{lem}
The proof is the same as the proof of Lemma \ref{limitinterior}.
\begin{lem}\label{joininterior}
For any $i$ we have
$$
\aligned
\sum_{\beta}\sum_{\ell}\sum_{k_1,k_2}&\frac{T^{\beta\cap \omega/2\pi}}{\ell!}
\exp(\beta\cap \frak b_2)\rho(\partial \beta) \\
&\int_{\frak{forget}^{-1}([\Sigma'_i])
\cap \mathcal M_{(k_1+1,k_2+1);\ell}^{\text{\rm ann;main}}
(\beta;\frak b_{\rm high}^{\otimes\ell})}(\text{\rm ev}^1,\text{\rm ev}^2)^*([\frak v;k_1]
, [\frak w;k_2])\\
= \sum_{\beta}\sum_{\ell}\sum_{k_1,k_2}&\frac{T^{\beta\cap \omega/2\pi}}{\ell!}
\exp(\beta\cap \frak b_2)\rho(\partial \beta) \\
&\int_{\frak{forget}^{-1}([\Sigma''_i])
\cap \mathcal M_{(k_1+1,k_2+1);\ell}^{\text{\rm ann;main}}
(\beta;\frak b_{\rm high}^{\otimes\ell})}(\text{\rm ev}^1,\text{\rm ev}^2)^*([\frak v;k_1]
, [\frak w;k_2]).
\endaligned
$$
\end{lem}
\begin{proof}
We take a path $\gamma:  [0,1] \to \mathcal M_{(1,1);0}^{\text{\rm ann;main}}$
joining $[\Sigma'_i]$ with $[\Sigma''_i]$ in the stratum 1 in Lemma \ref{stratifyAn11}.
We take the fiber product
$$
[0,1] {}_{\gamma}\times_{\frak{forget}}
\mathcal M_{(k_1,k_2);\ell}^{\text{\rm ann;main}}(\beta;
\text{\bf p})
$$
and denote it by
$\mathcal M_{(k_1,k_2);\ell}^{\text{\rm ann;main}}(\beta;
\text{\bf p};\gamma)$. This space has a Kuranishi structure and
its boundary is described as follows:
\begin{enumerate}
\item
$\frak{forget}^{-1}([\Sigma'_i])
\cap \mathcal M_{(k_1,k_2);\ell}^{\text{\rm ann;main}}
(\beta;\text{\bf p})$.
\item
$\frak{forget}^{-1}([\Sigma''_i])
\cap \mathcal M_{(k_1,k_2);\ell}^{\text{\rm ann;main}}
(\beta;\text{\bf p})$.
\item
$$
\mathcal M_{k'_1+1;\ell'}^{\rm main}(\beta';\text{\bf p}_1)^{\frak c} {}_{\text{\rm ev}_0}\times_{\text{\rm ev}^1_i}
\mathcal M_{(k''_1+1,k_2+1);\ell''}^{\text{\rm ann;main}}(\beta'';\text{\bf p}_2;\gamma).
$$
Here $k'_1+k''_1 = k_1+1$, $\ell'+\ell'' = \ell$,
$\beta'+\beta'' = \beta$, $i = 1,\dots,k''_1$.
\item
$$
\mathcal M_{k'_2+1;\ell'}^{\rm main}(\beta';\text{\bf p}_1)^{\frak c} {}_{\text{\rm ev}_0}\times_{\text{\rm ev}^2_i}
\mathcal M_{(k_1+1,k''_2+1);\ell''}^{\text{\rm ann;main}}(\beta'';\text{\bf p}_2;\gamma).
$$
Here $k'_2+k''_2 = k_2+1$, $\ell'+\ell'' = \ell$,
$\beta'+\beta'' = \beta$, $i = 1,\dots,k''_2$.
\item
$$
\mathcal M_{(k'_1+1,k_2+1);\ell'}^{\text{\rm ann;main}}(\beta';\text{\bf p}_1)
 {}_{\text{\rm ev}^1_0}\times_{\text{\rm ev}_i} 
\mathcal M_{k''_1+1;\ell''}^{\rm main}(\beta';\text{\bf p}_2;\gamma')
^{\frak c}.
$$
Here $k'_1+k''_1 = k_1+1$, $\ell'+\ell'' = \ell$,
$\beta'+\beta'' = \beta$, $i = 1,\dots,k''_1$.
\item
$$
\mathcal M_{(k_1+1,k'_2+1);\ell'}^{\text{\rm ann;main}}(\beta';\text{\bf p}_1)
 {}_{\text{\rm ev}^2_0}\times_{\text{\rm ev}_i} 
\mathcal M_{k''_2+1;\ell''}^{\rm main}(\beta'';\text{\bf p}_2;\gamma)
^{\frak c}.
$$
Here $k'_2+k''_2 = k_1+1$, $\ell'+\ell'' = \ell$,
$\beta'+\beta'' = \beta$, $i = 1,\dots,k''_2$.
\end{enumerate}
We define a multisection $\frak s$ on
$\mathcal M_{(k_1,k_2);\ell}^{\text{\rm ann;main}}(\beta;
\text{\bf p};\gamma)$ so that it is compatible at the above boundaries.
We also assume $\frak s$ is transversal to $0$.
Moreover it is symmetric with respect to the 
permutation of interior marked points.
\begin{rem}
We do not claim our Kuranishi structure and multisection $\frak s$ is 
$T^n$ equivariant. We do not claim $\frak s$ is invariant uder the cyclic permutation of the 
boundary marked points.
Actually in most of the part of our moduli space, we can keep both of the symmetries.
However in a neighborhood of $\frak{forget}^{-1}([\Sigma_1])$ it is hard to keep 
$T^n$-equivariance, and in a neighborhood of  $\frak{forget}^{-1}([\Sigma_2])$ it is hard to keep 
cyclic symmetry.
\end{rem}
\par
Now we put
\begin{equation}\label{betasplit3}
\frak b_{\rm high}^{\otimes\ell} = \sum_{I \in \frak I_{\ell}} c_I \text{\bf f}_{i_1}
\otimes \dots \otimes \text{\bf f}_{i_{\ell}}
= \sum_{I \in \frak I_{\ell}} c_{I,\ell} \text{\bf f}_{I,\ell}
\nonumber\end{equation}
and apply Stokes's theorem
for $(\text{\rm ev}^1,\text{\rm ev}^2)^*([\frak v;k_1]
, [\frak w;k_2])$
on
$\mathcal M_{(k_1,k_2);\ell}^{\text{\rm ann;main}}(\beta;
\text{\bf f}_{I,\ell};\gamma)^{\frak s}$.
The rest of the proof is the same as the proof of
Lemma \ref{backisp}.
\end{proof}
Proposition \ref{doublehatistwodisk} follows from
Lemmata \ref{limitinterior}, \ref{limitinterior2}, and
\ref{joininterior}.
\end{proof}
\begin{proof}[Proof of  Proposition \ref{twodiskandm2}]
We begin with the following lemma.
Let $g \in T^n$. We let it act on $L(u) \times L(u)$ by
$g(x,y) = (gx,y)$. Let $dg$ be the Haar measure of $T^n$ and let
$T_{\Delta}$ be an distributional $n$ form on $L(u) \times L(u)$
which is Poincar\'e dual to the diagonal.
We regard $\text{\bf e}_I$ (the generator of $H(L(u);\C)$) as a
harmonic form on $L(u)$. Let $g^{IJ}$ be as in Theorem \ref{annulusmain}.
\begin{lem}\label{diagoavelage}
$$
\int_{g \in T^n} g_*T_{\Delta} dg =
\sum_{I,J} (-1)^{\vert I \vert \vert J\vert}g^{IJ} \text{\bf e}_I \times \text{\bf e}_J.
$$
\end{lem}
See Lemma \ref{LemmaE} for the sign.
(Precisely speaking, in Lemma \ref{LemmaE} we work on chains instead of cochains.)
Here the equality is as distributional $n$ forms on $L(u) \times L(u)$.
\begin{proof}
The left hand side is a $T^n \times T^n$ invariant $n$ form. Therefore
it is smooth and harmonic. It is also cohomologous to the Poincar\'e dual to the
diagonal. The lemma follows.
\end{proof}
Let $g_1,g_2 \in T^n$. We define
\begin{equation}\label{ddoublefiber1prime}
(g_1\text{\rm ev}^0_1,g_2\text{\rm ev}^0_2):
\mathcal M^{\text{\rm main}}_{\vec k^{(1)};\ell_1}(\beta_1;
\text{\bf p}_1)^{\frak t}
\to L(u) \times L(u)
\end{equation}
by
$$
(g_1\text{\rm ev}^0_1,g_2\text{\rm ev}^0_2)(\text{\bf x}) = (g_1 \text{\rm ev}^0_1(\text{\bf x}),
g_2 \text{\rm ev}^0_2(\text{\bf x})).
$$
\begin{lem}\label{transversedoublefiber2}
The two maps $(\ref{ddoublefiber1prime})$ and
$(\ref{ddoublefiber22})$ are transversal.
\end{lem}
This is a consequence of Lemma \ref{existstkura}.
We denote the fiber product of $(\ref{ddoublefiber1prime})$ and
$(\ref{ddoublefiber22})$ by
$
\mathcal N(\vec k^{(1)},\vec k^{(2)};\ell_{1},\ell_{2};
\beta_1,\beta_2;\text{\bf p}_1;\text{\bf p}_2;g_1,g_2)
$. 
Note that 
we use $\frak t$ perturbation in Lemmas \ref{diagomodosi}, \ref{avelageddiagm}, \ref{mainlemmaGavelage} below.
\begin{lem}\label{diagomodosi}
$$
\aligned
&\frac{1}{\ell!}\int_{\mathcal M_{(k_1,k_2);\ell}^{\text{\rm ann;main}}(\beta;
\frak b_{\rm high}^{\otimes\ell}) \cap \frak{forget}^{-1}([\Sigma_2])}(\text{\rm ev}^1,\text{\rm ev}^2)^*([\frak v;k_1], [\frak w;k_2])\\
&=
\sum_{\beta_1+\beta_2=\beta \atop \ell_1+\ell_2=\ell}
\sum_{k_1 = k^{(1)}_2+ k^{(2)}_1 + k^{(2)}_3
\atop k_2 =  k^{(2)}_2+ k^{(1)}_1 + k^{(1)}_3}\frac{1}{\ell_1!\ell_2!}\int_{\mathcal N(\vec k^{(1)},\vec k^{(2)};\ell_{1},\ell_{2};
\beta_1,\beta_2;\frak b_{\rm high}^{\otimes\ell_1};\frak b_{\rm high}^{\otimes\ell_2};e,e)}(\text{\rm ev}^1,\text{\rm ev}^2)^*([\frak v;k_1], [\frak w;k_2]).
\endaligned
$$
Here $e$ denotes the unit of the group $T^n$.
\end{lem}
\begin{proof}
This is immediate from Lemma \ref{doublegluelem} and the definition.
\end{proof}
\begin{lem}\label{avelageddiagm}
\begin{equation}
\aligned
&\sum_{\beta_1+\beta_2=\beta \atop \ell_1+\ell_2=\ell}
\sum_{k_1 = k^{(1)}_2+ k^{(2)}_1 + k^{(2)}_3
\atop k_2 =  k^{(2)}_2+ k^{(1)}_1 + k^{(1)}_3}\frac{T^{(\beta_1+\beta_2)\cap \omega}}{\ell_1!\ell_2!}
\exp((\beta_1+\beta_2)\cap \frak b_2)\rho(\partial \beta)\\
&\int_{(g_1,g_2) \in (T^n)^2}dg_1dg_2\int_{\mathcal N(\vec k^{(1)},\vec k^{(2)};\ell_{1},\ell_{2};
\beta_1,\beta_2;\frak b_{\rm high}^{\otimes\ell_1};\frak b_{\rm high}^{\otimes\ell_2};g_1,g_2)}(\text{\rm ev}^1,\text{\rm ev}^2)^*([\frak v;k_1], [\frak w;k_2]) \\
&=
\sum_{\beta_1,\beta_2}\sum_{\ell_1+\ell_2}\sum_{I,J \in 2^{\{1,\dots,n\}}}
g^{IJ}
\frac{T^{(\beta_1+\beta_2)\cap \omega}}{\ell_1!\ell_2!}
\exp((\beta_1+\beta_2)\cap \frak b_2)\rho(\partial (\beta_1+\beta_2))
\\
&\qquad\qquad\qquad\qquad
\langle \frak q^{\frak c,b_+^{\frak c}}_{2,\ell_1;\beta_1}(\frak b_{\rm high}^{\otimes\ell_1}
; (\frak v \otimes \text{\bf e}_I)),
\frak q^{\frak c,b_+^{\frak c}}_{2,\ell_2;\beta_2}(\frak b_{\rm high}^{\otimes\ell_2}
; ( \text{\bf e}_J\otimes \frak w))
\rangle_{\text{\rm PD}_{L}}.
\endaligned
\nonumber\end{equation}
\end{lem}
We will prove Lemma \ref{avelageddiagm} at the end of this section.
\begin{lem}\label{mainlemmaGavelage}
$$
\aligned
\sum_{\beta_1+\beta_2=\beta \atop \ell_1+\ell_2=\ell}
&\sum_{k_1 = k^{(1)}_2+ k^{(2)}_1 + k^{(2)}_3
\atop k_2 =  k^{(2)}_2+ k^{(1)}_1 + k^{(1)}_3}\frac{T^{(\beta_1+\beta_2)\cap \omega}}{\ell_1!\ell_2!}
\exp((\beta_1+\beta_2)\cap \frak b_2)\rho(\partial \beta)\\
&\int_{\mathcal N(\vec k^{(1)},\vec k^{(2)};\ell_{1},\ell_{2};
\beta_1,\beta_2;\frak b_{\rm high}^{\otimes\ell_1};\frak b_{\rm high}^{\otimes\ell_2};g_1,g_2)}(\text{\rm ev}^1,\text{\rm ev}^2)^*([\frak v;k_1], [\frak w;k_2]) \endaligned
$$
is independent of $g_1, g_2 \in T^n$.
\end{lem}
\begin{proof}
We consider the case $(g_1,g_2) =(g,e)$ for simplicity.
Let $g,g' \in T^n$ and
$\gamma: [0,1] \to T^n$ be a path joining them.
We put
$$\aligned
&\mathcal N(\vec k^{(1)},\vec k^{(2)};\ell_{1},\ell_{2};
\beta_1,\beta_2;\frak b_{\rm high}^{\otimes\ell_1};\frak b_{\rm high}^{\otimes\ell_2};\gamma,e) \\
&=
\bigcup_{t \in [0,1]} \{t\}
\times \mathcal N(\vec k^{(1)},\vec k^{(2)};\ell_{1},\ell_{2};
\beta_1,\beta_2;\frak b_{\rm high}^{\otimes\ell_1};\frak b_{\rm high}^{\otimes\ell_2};\gamma(t),e).
\endaligned
$$
Note by Lemma \ref{transversedoublefiber2} the $\frak t$-perturbation
is transversal for {\it any}
$$
\mathcal N(\vec k^{(1)},\vec k^{(2)};\ell_{1},\ell_{2};
\beta_1,\beta_2;\frak b_{\rm high}^{\otimes\ell_1};\frak b_{\rm high}^{\otimes\ell_2};\gamma(t),e).
$$
So we have already chosen the
perturbation of the moduli space
$$
\mathcal N(\vec k^{(1)},\vec k^{(2)};\ell_{1},\ell_{2};
\beta_1,\beta_2;\frak b_{\rm high}^{\otimes\ell_1};\frak b_{\rm high}^{\otimes\ell_2};\gamma,e).
$$
The boundary of
$
\mathcal N(\vec k^{(1)},\vec k^{(2)};\ell_{1},\ell_{2};
\beta_1,\beta_2;\frak b_{\rm high}^{\otimes\ell_1};\frak b_{\rm high}^{\otimes\ell_2};\gamma,e)
$
is
$$
\mathcal N(\vec k^{(1)},\vec k^{(2)};\ell_{1},\ell_{2};
\beta_1,\beta_2;\frak b_{\rm high}^{\otimes\ell_1};\frak b_{\rm high}^{\otimes\ell_2};g,e)
$$
and
$$
\mathcal N(\vec k^{(1)},\vec k^{(2)};\ell_{1},\ell_{2};
\beta_1,\beta_2;\frak b_{\rm high}^{\otimes\ell_1};\frak b_{\rm high}^{\otimes\ell_2};g',e),
$$
and those described by the following Figure 3.4.6.
\par
\epsfbox{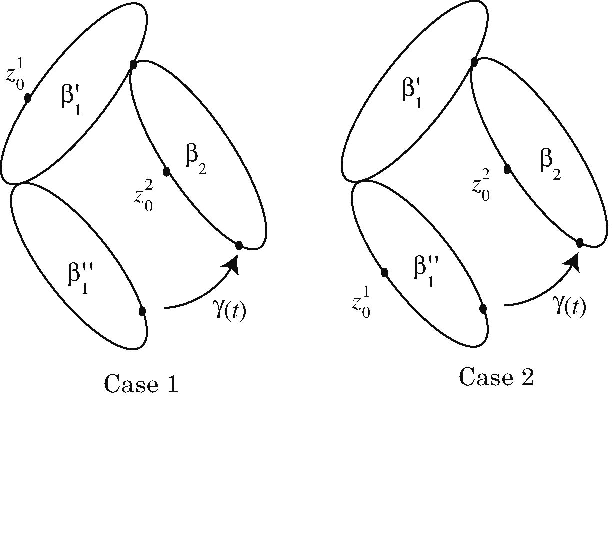}
\par
\hskip-1.4cm
\epsfbox{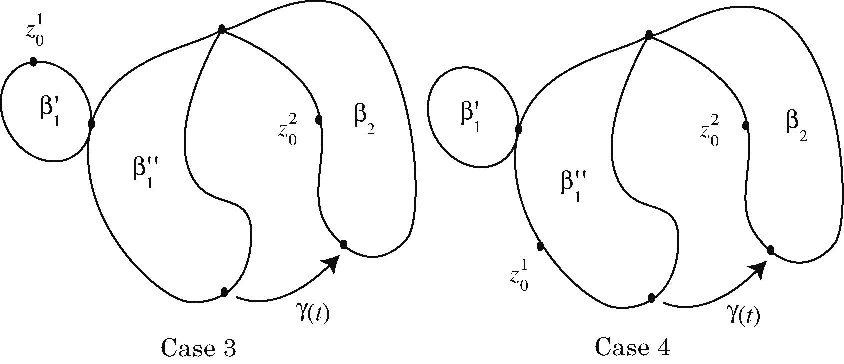}
\par
\hskip-1.0cm
\epsfbox{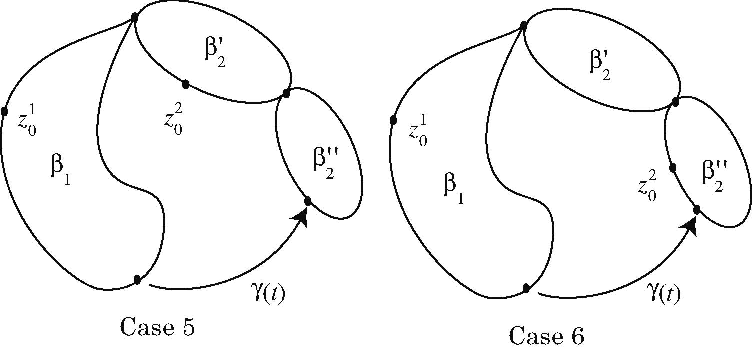}
\par
\epsfbox{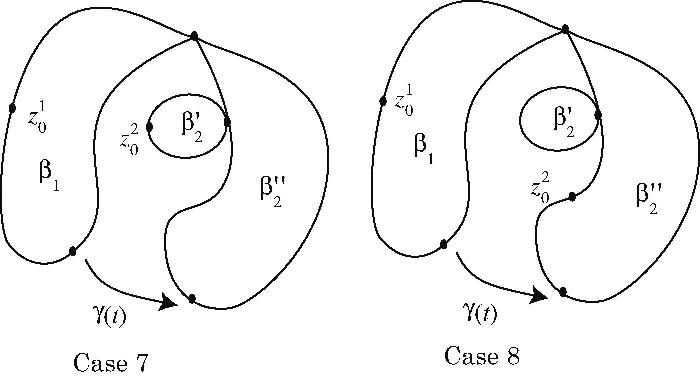}
\par
\centerline{\bf Figure 3.4.6}
\par
We apply Stokes' theorem to $(\text{\rm ev}^1,\text{\rm ev}^2)^*([\frak v;k_1], [\frak w;k_2])$.
It suffices to show that the contribution from the boundaries
of Figure 3.4.6 becomes zero after taking the weighted sum.
\par
The boundary Case 1 cancels with those of Case 5, after taking sum
over $\beta_1,\beta_2$.
\par
Similarly the boundary Case 2 cancels with those of Case 6, after taking sum
over $\beta_1,\beta_2$. (We use the $T^n$-equivariance of the
moduli space and $\frak c$-perturbation here.)
\par
The boundaries of Cases 4 and 8 are zero. This is by the same reason as
Lemma \ref{bdrydescription}.3.
(Namely unitality of the element $E$ there and compatibility with
the forgetful map of our perturbations.)
(Note that there are some other cases where bubbles with no boundary marked points
occur at the different positions. They all can be handled in the same way.)
\par
Finally the contribution from the boundary
components of Cases 3 and 7 becomes zero after
taking weighted sum. This is by the same reason as 
Case 4 of Lemma \ref{bdrydescription}.
(The vanishing of the boundary operator.)
Lemma \ref{mainlemmaGavelage} is proved.
\end{proof}
Lemmata \ref{diagomodosi}, \ref{avelageddiagm}, \ref{mainlemmaGavelage}
imply Proposition \ref{twodiskandm2}.
\end{proof}
We finally prove Lemma \ref{avelageddiagm}.
\begin{proof}[Proof of Lemma \ref{avelageddiagm}]
The idea of the proof is to interpolate $\frak c$ and $\frak t$ perturbation and 
use Stokes' theorem.
We first consider the map
\begin{equation}\label{ddoublefiber1prime22}
\frak E :
(T^n)^2 \times \mathcal M^{\text{\rm main}}_{\vec k^{(1)};\ell_1}(\beta_1;
\text{\bf p}_1)
\to L(u) \times L(u)
\end{equation}
defined by
$$
\frak E((g_1,g_2),{\bf x}) = (g_1{\rm ev}^0_1,g_2{\rm ev}^0_2).
$$

Since the map $\frak E$ in \eqref{ddoublefiber1prime22} is a submersion, we have the following 

\begin{lem}\label{ctOKOK}
The map (\ref{ddoublefiber1prime22}) is transversal to (\ref{ddoublefiber22}) for both
$\frak c$-perturbation and $\frak t$-perturbation.
\end{lem}

\par
We note that the fiber product of (\ref{ddoublefiber1prime22}) and (\ref{ddoublefiber22})
is identified with the union
\begin{equation}\label{Tparafiber}
\bigcup_{(g_1,g_2)\in (T^n)^2}
\{(g_1,g_2)\} \times \mathcal N(\vec k^{(1)},\vec k^{(2)};\ell_{1},\ell_{2};
\beta_1,\beta_2;\frak b_{\rm high}^{\otimes\ell_1};\frak b_{\rm high}^{\otimes\ell_2};g_1,g_2)
\end{equation}
We denote the space (\ref{Tparafiber}) by
$\mathcal N(\vec k^{(1)},\vec k^{(2)};\ell_{1},\ell_{2};
\beta_1,\beta_2;\frak b_{\rm high}^{\otimes\ell_1};\frak b_{\rm high}^{\otimes\ell_2};T^n,T^n)$.
By Lemma \ref{transversedoublefiber2}
the $\frak t$ (resp. $\frak c$) perturbation induces a perturbation of this space. 
We denote the corresponding space by 
$\mathcal N(\vec k^{(1)},\vec k^{(2)};\ell_{1},\ell_{2};
\beta_1,\beta_2;\frak b_{\rm high}^{\otimes\ell_1};\frak b_{\rm high}^{\otimes\ell_2};T^n,T^n)^{\frak t}$ 
(resp. $\mathcal N(\vec k^{(1)},\vec k^{(2)};\ell_{1},\ell_{2};
\beta_1,\beta_2;\frak b_{\rm high}^{\otimes\ell_1};\frak b_{\rm high}^{\otimes\ell_2};T^n,T^n)^{\frak c}$), 
when we specify the perturbation. 
We denote the projection from 
$\mathcal N(\vec k^{(1)},\vec k^{(2)};\ell_{1},\ell_{2};
\beta_1,\beta_2;\frak b_{\rm high}^{\otimes\ell_1};\frak b_{\rm high}^{\otimes\ell_2};T^n,T^n)$
to $T^{2n}$ by $\pi_{T^{2n}}$.
\begin{lem}\label{lhs3428}
The left hand side of Lemma \ref{avelageddiagm} is
\begin{equation}\label{3426aa}
\aligned
&\sum_{\beta_1+\beta_2=\beta \atop \ell_1+\ell_2=\ell}
\sum_{k_1 = k^{(1)}_2+ k^{(2)}_1 + k^{(2)}_3
\atop k_2 =  k^{(2)}_2+ k^{(1)}_1 + k^{(1)}_3}\frac{T^{(\beta_1+\beta_2)\cap \omega}}{\ell_1!\ell_2!}
\exp((\beta_1+\beta_2)\cap \frak b_2)\rho(\partial \beta)
\\
&\int_{\mathcal N(\vec k^{(1)},\vec k^{(2)};\ell_{1},\ell_{2};
\beta_1,\beta_2;\frak b_{\rm high}^{\otimes\ell_1};\frak b_{\rm high}^{\otimes\ell_2};T^n,T^n)^{\frak t}}
\pi_{T^{2n}}^*(dg_1dg_2) \wedge (\text{\rm ev}^1,\text{\rm ev}^2)^*([\frak v;k_1], [\frak w;k_2]). 
\endaligned
\end{equation}
Here $dg_i$ is the normalized volume form of $T^n$.
\end{lem}
This lemma is obvious from definition.
(Note we use $\frak t$-perturbation here.)
\begin{lem}\label{rhs3428}
The right hand side of Lemma \ref{avelageddiagm} is
\begin{equation}\label{rhs3428form}
\aligned
&\sum_{\beta_1+\beta_2=\beta \atop \ell_1+\ell_2=\ell}
\sum_{k_1 = k^{(1)}_2+ k^{(2)}_1 + k^{(2)}_3
\atop k_2 =  k^{(2)}_2+ k^{(1)}_1 + k^{(1)}_3}\frac{T^{(\beta_1+\beta_2)\cap \omega}}{\ell_1!\ell_2!}
\exp((\beta_1+\beta_2)\cap \frak b_2)\rho(\partial \beta)
\\
&\int_{\mathcal N(\vec k^{(1)},\vec k^{(2)};\ell_{1},\ell_{2};
\beta_1,\beta_2;\frak b_{\rm high}^{\otimes\ell_1};\frak b_{\rm high}^{\otimes\ell_2};T^n,T^n)^{\frak c}}
\pi_{T^{2n}}^*(dg_1dg_2) \wedge (\text{\rm ev}^1,\text{\rm ev}^2)^*([\frak v;k_1], [\frak w;k_2]). 
\endaligned
\end{equation}
\end{lem}
\begin{proof}
Note we use $\frak c$-perturbation in (\ref{rhs3428form}).
The lemma then follows from Lemma \ref{diagoavelage}.
\end{proof}
So to complete the proof of Lemma \ref{avelageddiagm} it suffices to show 
that (\ref{3426aa}) is equal to (\ref{rhs3428form}).
\begin{lem}\label{existstkura22}
There exist systems of Kuranishi structures and multisections on  
$[0,1]\times \mathcal M^{\text{\rm main}}_{k_1,k_2,k_3;\ell}(\beta;
\text{\bf p})$ with the following properties.
We call this multisection $\frak{para}$-multisection.
\begin{enumerate}
\item
They are transversal to $0$ and are $T^n$ equivariant.
\item
They are compatible with the forgetful map of the forgetable boundary marked points.
\item
They are invariant under the permutation of the interior marked points together with the 
permutation of ${\bf p}_i$.
\item
The evaluation map (with $[0,1]$ component)
$$
(\pi_{[0,1]}\times  {\rm ev}^0_i) :  ([0,1]\times \mathcal M^{\text{\rm main}}_{k_1,k_2,k_3;\ell}(\beta;
\text{\bf p}))^{\frak{para}}
\to [0,1] \times L(u)
$$
($i=0,1,2$)
at one of the unforgetable marked points is a submersion on the zero set of 
this multisection.
\item
The boundary of $([0,1]\times \mathcal M^{\text{\rm main}}_{k_1,k_2,k_3;\ell}(\beta;
\text{\bf p}))^{\frak{para}}$ is described as the union of the following three types 
of fiber products.
\begin{enumerate}
\item
\begin{equation}
\mathcal M^{\text{\rm main}}_{k'_i+k'_{i+1}+2;\ell'}(\beta_1;
\text{\bf p}_1)^{\frak c}
\,{}_{{\rm ev}_0}\times_{{\rm ev}^0_i}
([0,1]\times \mathcal M^{\text{\rm main}}_{k''_1,k''_2,k''_3;\ell''}(\beta_2;
\text{\bf p}_2))^{\frak{para}}.
\end{equation}
Here $i$ is one of $1$,$2$,$3$. ($i+1$ means $1$ if $i=3$.)
Moreover $k'_i + k''_{i} = k_i$, $k'_{i+1} + k''_{i+1} = k_{i+1}$
and $k''_j = k_j$ if $j\ne i,i+1$.
\par
$\beta = \beta_1+\beta_2$, $\ell = \ell_1 + \ell_2$ and $\text{\bf p} = (\text{\bf p},\text{\bf p}_2)$.
\par
The $i$-th unforgetable marked point of $[0,1]\times \mathcal M^{\text{\rm main}}_{k_1,k_2,k_3;\ell}(\beta;
\text{\bf p})$ becomes $k'_i+1$-th marked point of the first factor 
$\mathcal M^{\text{\rm main}}_{k'_i+k'_{i+1}+2;\ell'}(\beta_1;
\text{\bf p}_1)^{\frak c}$.
The $j$-th unforgetable marked point of $[0,1]\times \mathcal M^{\text{\rm main}}_{k_1,k_2,k_3;\ell}(\beta;
\text{\bf p})$ 
becomes $j$-th unforgetable marked point of the second factor.
\par
We remark that we take $\frak c$-multisection for the first factor 
and $\frak{para}$-multisection for the second factor.
\item
\begin{equation}
\mathcal M^{\text{\rm main}}_{k'_i+1;\ell'}(\beta_1;
\text{\bf p}_1)^{\frak c}
\,{}_{{\rm ev}_0}\times_{{\rm ev}_{i,j}}
([0,1]\times \mathcal M^{\text{\rm main}}_{k''_1,k''_2,k''_3;\ell''}(\beta_2;
\text{\bf p}_2))^{\frak{para}}.
\end{equation}
Here $i$ is one of $1$,$2$,$3$. ($i+1$ means $1$ if $i=3$.)
Moreover $k'_i + k''_{i} -1 = k_i$, 
and $k''_j = k_j$ if $j\ne i$. $j \in \{1,\dots,k''_i\}$, $k''_i \ge 1$. 
${\rm ev}_i = ({\rm ev}_{i,1},\dots,{\rm ev}_{i,j},\dots,{\rm ev}_{i,k''_i})$ where $ev_i$ is as in 
(\ref{form316}).
\par
$\beta = \beta_1+\beta_2$, $\ell = \ell_1 + \ell_2$ and $\text{\bf p} = (\text{\bf p},\text{\bf p}_2)$.
\par
The $i$-th unforgetable marked point of $\mathcal M^{\text{\rm main}}_{k_1,k_2,k_3;\ell}(\beta;
\text{\bf p})$ becomes the $i$-th unforgetable marked point of the second factor.
\par
We remark that we take $\frak c$-multisection for the first factor 
and $\frak t$-multisection for the second factor.
\item
$$
(\{0,1\}\times \mathcal M^{\text{\rm main}}_{k_1,k_2,k_3;\ell}(\beta;
\text{\bf p}))^{\frak{para}}.
$$
We put $\frak c$-perturbation for
$
(\{0\}\times \mathcal M^{\text{\rm main}}_{k_1,k_2,k_3;\ell}(\beta;
\text{\bf p}))^{\frak{para}}
$
and $\frak t$-perturbation for
$
(\{1\}\times \mathcal M^{\text{\rm main}}_{k_1,k_2,k_3;\ell}(\beta;
\text{\bf p}))^{\frak{para}}.
$
\end{enumerate}
\end{enumerate}
\end{lem}
The proof is similar to the proof of Lemma \ref{bdryMannu}
which is given in
 Subsection \ref{subsec:tkuranishi}.
\end{proof}
Now we use $\frak{para}$-perturbation in a similar way as Lemma  \ref{ctOKOK} 
to obtain a Kuranishi structure on 
$[0,1] \times \mathcal N(\vec k^{(1)},\vec k^{(2)};\ell_{1},\ell_{2};
\beta_1,\beta_2;\frak b_{\rm high}^{\otimes\ell_1};\frak b_{\rm high}^{\otimes\ell_2};T^n,T^n)$
such that it becomes $\frak c$ (resp. $\frak t$) perturbation 
on 
$\{0\} \times \mathcal N(\vec k^{(1)},\vec k^{(2)};\ell_{1},\ell_{2};
\beta_1,\beta_2;\frak b_{\rm high}^{\otimes\ell_1};\frak b_{\rm high}^{\otimes\ell_2};T^n,T^n)^{\frak c}$
(resp. 
$\{1\} \times \mathcal N(\vec k^{(1)},\vec k^{(2)};\ell_{1},\ell_{2};
\beta_1,\beta_2;\frak b_{\rm high}^{\otimes\ell_1};\frak b_{\rm high}^{\otimes\ell_2};T^n,T^n)^{\frak t}$.)
Using it we apply Stokes' theorem in the same way as the proof of Lemma \ref{backisp}, 
to show that (\ref{3426aa}) is equal to (\ref{rhs3428form}).
It implies Lemma \ref{avelageddiagm}.
\qed
\par
The proof of Theorem \ref{annulusmain} is now complete.
\qed
\begin{rem}
We actually did not use $T^n$-equivariance of $\frak t$-perturbation.
But we did use $T^n$-equivariance of $\frak c$-perturbation.
\end{rem}
\begin{rem}
Theorem \ref{annulusmain} is stated and proved only for the case of
fibers of toric manifolds. However various parts of the proof
(especially the idea to use the moduli space of pseudo-holomorphic
annuli and cyclic symmetry) can be applied to more general situation.
\par
We choose to prove this statement only for the toric case in this paper
since in toric case we can take short cut in several places, while
the proof for the general case is rather involved.
For example, in various places we use $T^n$-equivariant perturbations in the present paper
and so do not need to use general differential forms on $L(u) = T^n$
other than the $T^n$-invariant ones. In the general situation we need to
use a canonical model and so the moduli spaces involved become more complicated to define.
Also we use the technique of averaging differential forms
over the $T^n$ action in our proof of Proposition \ref{twodiskandm2}.
Since this technique is not available for other cases, we need to
use more sophisticated homological algebra instead in the non-toric cases.
\par
Another reason why we restrict ourselves to the toric case here is that
at the time of writing this article we do not have any other particular examples
to which we may apply the ideas used in this paper. However we have little doubt that many
non-toric examples, to which many ideas presented here can be applied,
will be found in the near future.
\end{rem}
\par
\section{Poincar\'e duality is residue pairing}
\label{sec:PDRes}
We take
$$
\text{\rm vol}_{L(u)} \in H^n(L(u);\Lambda_0)
\cong HF((L(u),\frak b,b^{\frak c}),(L(u),\frak b,b^{\frak c});\Lambda_0),
$$
which is the volume form of $L(u)$. We normalize it so that $\int_{L(u)}\text{\rm vol}_{L(u)} = 1$.
Note we identify $L(u)$ with a standard torus $T^n$ by using the fact that $L(u)$ is a $T^n$ orbit.
So using the standard orientation of $T^n$ we orient $L(u)$. 
(See \cite[p127]{fooo08}.)
\begin{rem}
It may be confusing that $\text{\rm vol}_{L(u)} \ne 1$ ($={\bf e}_L$).
We note that
$1$ is the Poincar\'e dual to the fundamental {\it homology}
class $[L(u)]$. The degree of $\text{\rm vol}_{L(u)}$ is $n$ and the
degree of $1$ is $0$.
\end{rem}
The results of Sections
\ref{sec:operatorptoric},
\ref{sec:annuli} imply:
\begin{prop}\label{mainprop} Let $(u,b) \in \text{\rm Crit}(\frak{PO}_{\frak b})$.
Then we have\index{trace}\index{Frobenius algebra!trace}\index{$Z(C)$}
$$
\langle i_{\ast,\text{\rm qm},(\frak b,b,u)}(\text{\rm vol}_{L(u)}),
i_{\ast,\text{\rm qm},(\frak b,b,u)}(\text{\rm vol}_{L(u)})
\rangle_{\text{\rm PD}_{X}}
= Z(\frak b,b).
$$
\end{prop}
\begin{rem}
We still do not assume that $(u,b)$ is nondegenerate.
When $(u,b)$ is a degenerate critical point, Proposition \ref{mainprop}
implies that the left hand side is $0$.
\end{rem}
\begin{proof}
This is immediate from  Theorem \ref{annulusmain} and
Corollary \ref{degreemustn}.
See Subsection \ref{subsec:signproof} for the proof of the sign
in the above formula.
\end{proof}
Now we are ready to complete the proof of Theorem \ref{Mirmain}.2.
We assume that $\frak{PO}_{\frak b}$ is a Morse function.
We first consider two elements $(u,b), (u',b') \in \text{\rm Crit}(\frak{PO}_{\frak b})$,
$(u,b) \ne (u',b')$ corresponding to the unit elements of some factors
$\text{\rm Jac}(\frak{PO}_{\frak b};\eta)$ in the factorization of
$\text{\rm Jac}(\frak{PO}_{\frak b})$ given in Proposition \ref{Morsesplit}.
We denote the corresponding unit elements of $\text{\rm Jac}(\frak{PO}_{\frak b})$
by $1_{(u,b)},  1_{(u',b')}$ respectively. Then
$$
\langle 1_{(u,b)}, 1_{(u',b')} \rangle_{\text{\rm PD}_{X}}
= \langle 1_{(u,b)}, 1_{(u',b')} \cup^{\frak b} 1 \rangle_{\text{\rm PD}_{X}}
= \langle 1_{(u,b)}\cup^{\frak b}  1_{(u',b')} , 1 \rangle_{\text{\rm PD}_{X}}
= 0.
$$
Here $\cup^{\frak b}$ is the quantum cup product deformed by $\frak b$.
(See Section \ref{sec:statements}.)
(The second equality is nothing but the well known identity for the 
Frobenius pairing in Gromov-Witten theory. See
for example \cite{Manin:qhm}.)
So it suffices to study the case $(u,b) = (u',b')$.
We identify $1_{(u,b)} \in \text{\rm Jac}(\frak{PO}_{\frak b})\otimes_{\Lambda_0}
\Lambda$ with an
element of $H(X;\Lambda)$ by the isomorphism given in Theorem \ref{Mirmain}.1.
Since $1_{(u,b)}$ is idempotent and
$\frak {ks}_{\frak b}$ is a ring isomorphism it follows from Lemma \ref{lem171} that
$$
i^{\ast}_{\text{\rm qm},(u,b)}(1_{(u,b)}) = 1 \in H^0(L(u);\Lambda).
$$
Therefore Theorem \ref{pqpoincare} and the equality 
$\langle 1,\text{vol}_{L(u)}\rangle_{\text{\rm PD}_{L(u)}} = 1$ imply
\begin{equation}\label{kocchikite1}
\langle 1_{(u,b)},i_{\ast,\text{\rm qm},(\frak b,b,u)}(\text{vol}_{L(u)})\rangle_{\text{\rm PD}_{X}}
= 1.
\end{equation}
Since
$$
i^{\ast}_{\text{\rm qm},(u,b)}(1_{(u',b')}) = 0
$$
for $(u',b') \ne (u,b)$, it follows from (\ref{kocchikite1}) that
\begin{equation}\label{kocchikite2}
i_{\ast,\text{\rm qm},(\frak b,b,u)}(\text{vol}_{L(u)})
= \frac{1}{\langle 1_{(u,b)},1_{(u,b)}\rangle_{\text{\rm PD}_{X}}}1_{(u,b)}.
\end{equation}
Therefore Proposition \ref{mainprop} implies
$$
\langle 1_{(u,b)},1_{(u,b)}\rangle_{\text{\rm PD}_{X}}
Z(\frak b,b)
= 1.
$$
The proof of Theorem \ref{Mirmain}.2 is now complete.
\qed
\par
\section{Clifford algebra and Hessian matrix}
\label{sec:clifford}
In Sections \ref{sec:clifford} and \ref{sec:ResHess},
we will prove Theorem \ref{cliffordZ}. 
In \cite{Cho05II}, Cho proved that Floer cohomology
of the toric fiber is isomorphic to the Clifford algebra\index{Clifford algebra}  
as a ring in the toric Fano case (with $b =0,  \frak b=0$). In this section, we follow Cho's argument
to prove the following Theorem \ref{clifford}, from which Theorem \ref{cliffordZ}.2 and 3 follow.
\par
Let $X_1,\dots,X_n$ be formal variables and $d_i \in \Lambda$
($i=1,\dots,n$). We consider relations
\begin{equation}\label{clifrel}
\left\{
\aligned
X_iX_j+X_jX_i&=0,\qquad \qquad i \ne j\\
X_i X_i&=d_i 1.
\endaligned
\right.
\end{equation}
We take a free (noncommutative) $\Lambda$ algebra generated by
$X_1,\dots,X_n$ and divide it by the two-sided ideal generated
by (\ref{clifrel}). We denote it by
$\text{\rm Cliff}_\Lambda(n;\vec d)$, where we set $\vec d = (d_1,\dots,d_n)$.
\par
Let $I = (i_1,\dots,i_k)$, $1 \le i_1<\dots<i_k \le n$.
We write the set of such $I$'s by $2^{\{1,\dots,n\}}$.
We put
$$
X_I = X_{i_1} X_{i_2} \cdots X_{i_{k-1}} X_{i_k} \in
\text{\rm Cliff}_\Lambda(n;\vec d).
$$
It is well known and can be easily checked that
$\{X_I \mid I \in 2^{\{1,\dots,n\}}\}$
forms a basis of $\text{\rm Cliff}_\Lambda(n;\vec d)$ as a
$\Lambda$ vector space. 
For $I \in 2^{\{1,\dots,n\}}$, we denote its
complement by $I^c$.
\par
Now we consider $(u,b) \in \text{\rm Crit}(\frak{PO}_{\frak b})$.
We take $(u,b^{\frak c})$ which is the
weak bounding cochain of the cyclic filtered $A_{\infty}$
algebra $(H(L(u);\Lambda_0),\{\frak m_k^{\frak c,\frak b}\})$
corresponding to $b$.
When Condition $\ref{fano2jigen}$.2 is satisfied, we have
$b=b^{\frak c}$.
If   Condition $\ref{fano2jigen}$.1 or 2
is satisfied, we have $b^{\frak c} \in H^1(L(u);\Lambda_0)$.
\par
We put $y(u)_i = e^{x_i}$ and regard $\frak{PO}_{\frak b}$ as
a function of $x_i$. Denote
$b = \sum_{i=1}^n x_i(b) \text{\bf e}_i$.
\begin{defn}
If Condition $\ref{fano2jigen}$.2 is satisfied,
we put $\frak{PO}_{\frak b}^{\frak c} = \frak{PO}_{\frak b}$.
If Condition $\ref{fano2jigen}$.1 is satisfied,
we put
\begin{equation}\label{FcandF}
\frak{PO}_{\frak b}^{\frak c} = \frak{PO}_{\frak b}\circ \frak F_*^{-1}
: H^1(L(u);\Lambda_0) \to \Lambda_0.
\end{equation}
\end{defn}
We recall that Corollary \ref{H1identity} shows
$\frak F_* =$ identity on $H^1(L(u);\Lambda_0)$
when Condition $\ref{fano2jigen}$.2 is satisfied.
\par
We now take an $n\times n$ matrix $A = [a_{ij}]_{i,j=1}^{i,j=n}$ ($a_{ij}
\in \Lambda$) such that
\begin{enumerate}
\item
\begin{equation}\label{def:di}
{}^tA \left[ \frac{\partial^2\frak{PO}_{\frak b}}{\partial x_i\partial x_j}
(x(b))\right]_{i,j=1}^{i,j=n} A = 2
\left[
\begin{matrix}
d_1  & 0   & \dots & 0 \\
0    & d_2 & \dots & 0 \\
\vdots & \vdots & \ddots & \vdots \\
0    & 0 & \dots & d_n
\end{matrix}
\right].
\end{equation}
\item $\det A = 1$.
\end{enumerate}
We put $\vec d=(d_1, \dots ,d_n)$. We note that 
$\vec d$ depends on the choice of $A$.
(If the potential function is a real valued function, we can take $A$ as an orthogonal matrix. But in general, 
$A$ is not an orthogonal matrix.) 
But the product $2^n d_1\cdots d_n$, the determinant of the Hessian matrix, is independent of the choice of $A$. 
\par
In Sections \ref{sec:PDRes} and \ref{sec:clifford},  we identify
$$
H(L(u),\Lambda)
\cong
HF((L(u);(\frak b,b));(L(u);(\frak b,b));\Lambda).
$$
Then by definition, the product $\cup^{\frak c,Q}$ is given by
$$
h_1 \cup^{\frak c,Q} h_2 = (-1)^{\deg h_1(\deg h_2+1)}\frak m_2^{\frak c,\frak b,b^{\frak c}}(h_1,h_2)
$$
where $\frak m_2^{\frak c,\frak b,b^{\frak c}}$ is as defined in \eqref{mcyclicdef2}.
Then the (filtered) $A_{\infty}$ relation implies that
$(H(L(u);\Lambda),\cup^{\frak c,Q})$ defines a ring.

\begin{thm}\label{clifford}
Suppose that one of the two conditions in Condition \ref{fano2jigen}
is satisfied. Then there exist a basis $\{\text{\bf e}'_i
\mid i=1,\dots,n\} \subset H^1(L(u);\Lambda)$
and a map
$$
\Phi: \text{\rm Cliff}_\Lambda(n;\vec d) \to (H(L(u);\Lambda),\cup^{\frak c,Q})
$$
with the following properties:
\begin{enumerate}
\item  $\Phi$ is a $\Lambda$ algebra isomorphism.
\item  $\Phi(X_i) = \text{\bf e}'_i$.
\item We put $\Phi(X_I) = \text{\bf e}'_I$. Then
$$
\langle \text{\bf e}'_I, \text{\bf e}'_J \rangle_{\text{\rm PD}_{L(u)}}
=
\begin{cases}
(-1)^{*(I)}& J = I^c,\\
0& \text{\rm otherwise}.
\end{cases}
$$
Here
\begin{equation}\label{starI}
*(I) = \# \{ (i,j) \in I \times I^c \mid j < i\}.
\end{equation}
\item If $I \ne \{1,\dots,n\}$, then
$$
\text{\bf e}'_I \in \bigoplus_{d<n} H^d(L(u);\Lambda).
$$
\item
$$
\text{\bf e}'_{\{1,\dots,n\}} - \text{\rm vol}_{L(u)}
\in \bigoplus_{d<n} H^d(L(u);\Lambda).
$$
\end{enumerate}
\end{thm}
\begin{proof}
We first consider the case $\dim L(u) =2$.
(That is the case where Condition $\ref{fano2jigen}$.1 holds.)
By Lemma \ref{H1isinMweak2}, we have
$b^{\frak c} \in H^1(L(u);\Lambda_0)$.
\par
Let $x \in H^1(L(u);\Lambda_0)$. We have
$$
\langle x\cup^{\frak c,Q} x, 1\rangle_{\text{\rm PD}_{L(u)}}
= \langle x , x\cup^{\frak c,Q} 1\rangle_{\text{\rm PD}_{L(u)}}
=  \langle x , x\rangle_{\text{\rm PD}_{L(u)}}
= 0.
$$
Therefore using $\dim L(u) =2$ again there exists
$\frak Q : H^1(L(u);\Lambda_0) \to \Lambda_0$ such that
$$
x\cup^{\frak c,Q} x = \frak Q(x) 1.
$$
It is easy to see that $\frak Q$ is a quadratic form.
In fact, we have
\begin{equation}
\frak Q(x) = \frac{1}{2}\frac{d^2}{dt^2} \frak{PO}_{\frak b}^{\frak c}(b^{\frak c} + tx).
\end{equation}
Therefore there exists a basis
$\text{\bf e}'_1$, $\text{\bf e}'_2$ of $H^1(L(u);\Lambda_0)$ such that
$$
\frak Q(x_1\text{\bf e}'_1+x_2\text{\bf e}'_2) = d_1x_1^2 + d_2x_2^2,
$$
where $d_1, d_2 \in \Lambda_0$ are 
as in \eqref{def:di}
for the Hessian of $\frak{PO}_{\frak b}^{\frak c}$
at $b^{\frak c}$.
By replacing $\text{\bf e}'_1$ by $c\text{\bf e}'_1$ for some $c \in \Lambda_0
\setminus \Lambda_+$ we may assume
$
\text{\bf e}'_1\text{\bf e}'_2 - \text{\rm vol}_{L(u)}
\in \bigoplus_{d<n} H^d(L(u);\Lambda).
$
It is easy to see that $\text{\bf e}'_1$, $\text{\bf e}'_2$ have required properties.
\par
We next consider the case Condition $\ref{fano2jigen}$.2
is satisfied.
We start with proving two lemmata. Consider the classical cup product
$h_1\wedge h_2$.
\begin{lem}\label{gradedprod}
Assume  Condition \ref{fano2jigen}.2. If $h_i \in
H^{d_i}(L(u);\Lambda)$, then
$$
h_1 \cup^{\frak c,Q} h_2 - h_1 \wedge h_2
\in \bigoplus_{d<d_1+d_2} H^d(L(u);\Lambda).
$$
\end{lem}
\begin{proof}
By Lemma \ref{nefcondconclupara}, we have
$\dim \mathcal M_{3;\ell}^{\text{\rm main};\text{\bf p}}(\beta;\text{\bf p}) \ge n+2$
if degree of $\text{\bf p}(i) = 2$. Therefore
$$\aligned
\deg \frak q_{\ell;2;\beta}(\frak b^{\otimes\ell},h_1,h_2)
&= \deg \text{\rm ev}_{0*}((\text{\rm ev}_1,\text{\rm ev}_2)^* (h_1\times h_2);
\mathcal M_{3;\ell}^{\text{\rm main}}(\beta;\frak b^{\otimes\ell}   )) \\
&= \deg h_1 + \deg h_2 + n - \dim \mathcal M_{3;\ell}^{\text{\rm
main};}(\beta;\frak b^{\otimes\ell}) \le d_1 + d_2 - 2.
\endaligned$$
The lemma follows.
\end{proof}
\begin{lem}\label{calcum2}
Assume  Condition \ref{fano2jigen}.2. Then we have
$$
\text{\bf e}_i \cup^{\frak c,Q} \text{\bf e}_j
+ \text{\bf e}_j \cup^{\frak c,Q} \text{\bf e}_i
= \left(\frac{\partial^2\frak{PO}_{\frak b}}{\partial x_i\partial x_j}
(x(b))\right) 1.
$$
\end{lem}
\begin{proof}
By definition
$$
\frak{PO}_{\frak b}^u(x)
= \sum_{\beta,\text{\bf p}}{c_{\text{\bf p}}}
T^{\beta\cap\omega/2\pi}
c(\beta;\text{\bf p})
\prod_{i=1}^n\exp((\text{\bf e}_i \cap \partial \beta)x_i).
$$
Here we put
$$
\sum_{\ell=0}^{\infty} \frac{1}{\ell!}\frak b^{\otimes\ell}
= \sum_{\text{\bf p}} c_{\text{\bf p}}\text{\bf p}
$$
and
$$
c(\beta;\text{\bf p}) = \text{\rm ev}_{0*}(\mathcal M_{1;\vert\text{\bf p}\vert}
(\beta;\text{\bf p})).
$$
The summation is taken
for all $\beta$, $\text{\bf p}$ with
$\dim \mathcal M_{1;\vert\text{\bf p}\vert}
(\beta;\text{\bf p}) = n$.
(See  \cite[(6.9), (7.2)]{fooo09}.) It is then easy to see that
\begin{equation}\label{secondderivative}
\frac{\partial^2\frak{PO}_{\frak b}}{\partial x_i\partial x_j}
= \sum_{\beta,\text{\bf p}}{c_{\text{\bf p}}}
T^{\beta\cap\omega/2\pi}
c(\beta;\text{\bf p})
(\text{\bf e}_i \cap \partial \beta)
(\text{\bf e}_j \cap \partial \beta)
\prod_{a=1}^n\exp((\text{\bf e}_a \cap \partial \beta)x_a).
\end{equation}
Let $h_i$ be the harmonic one form whose cohomology
class is $\text{\bf e}_i$.
We can easily prove that
\begin{equation}\label{correspondense3formula}
\aligned
&\text{\rm ev}_{0*}((\text{\rm ev}_1,\text{\rm ev}_2)^*(h_i \times h_j)
+ (\text{\rm ev}_1,\text{\rm ev}_2)^*(h_j \times h_i); \mathcal M_{3;\vert\text{\bf p}\vert}^{\text{\rm main};\text{\bf p}}(\beta;\text{\bf p})) \\
&= c(\beta;\text{\bf p}) (\text{\bf e}_i \cap \partial \beta)
(\text{\bf e}_j \cap \partial \beta).
\endaligned
\end{equation}
(See the proof of  \cite[Lemma 10.8]{fooo08} and Section \ref{sec:deltaisthesame} of the present paper
for the relevant argument.)
We observe 
\begin{equation}\label{q2formula}
\aligned
&\frak q^{\frak c}_{\ell;2}([\text{\bf p}] ;h_i\otimes h_j)\\
&= \Big(\sum_{\beta,\text{\bf p}}{c_{\text{\bf p}}}
T^{\beta\cap\omega/2\pi}  \left(\prod_{a=1}^n\exp((\text{\bf e}_a \cap \partial \beta)x_i)\right)\\
& \qquad
\times \text{\rm ev}_{0*}((\text{\rm ev}_1,\text{\rm ev}_2)^*(h_i \times h_j)
; \mathcal M_{3;\vert\text{\bf p}\vert}^{\text{\rm main};\text{\bf p}}(\beta;\text{\bf p}))^{\frak c}\Big) 1.
\endaligned\end{equation}
By taking Lemma \ref{dimandempty} into account, 
Lemma \ref{calcum2} follows immediately from
(\ref{secondderivative}), (\ref{correspondense3formula}),
(\ref{q2formula}).
\end{proof}
By the cyclicity of the $\frak c$-perturbation we have cyclic property
\begin{equation}\label{cupcycl}
\langle h_1 \cup^{\frak c,Q} h_2, h_3 \rangle_{\text{\rm PD}_{L(u)}}
= \langle h_1, h_2\cup^{\frak c,Q} h_3 \rangle_{\text{\rm PD}_{L(u)}}
\end{equation}
for $h_1,h_2,h_3 \in H(L(u);\Lambda_0)$.
\par
We recall that we take an $n\times n$ matrix $A = [a_{ij}]_{i,j=1}^{i,j=n}$ ($a_{ij}
\in \Lambda$) satisfying \eqref{def:di} and 
$\det A = 1$. 
Using the matrix $A$, we put
\begin{equation}\label{def:eiprime}
x'_i = \sum_j a_{ji} x_i, \qquad \text{\bf e}'_i = \sum_{ij} a_{ij}\text{\bf e}_i.
\end{equation}
Lemma \ref{calcum2} implies that
$$
\text{\bf e}'_i \cup^{\frak c,Q} \text{\bf e}'_j
+ \text{\bf e}'_j \cup^{\frak c,Q} \text{\bf e}'_i = 2d_i \delta_{ij}.
$$
Therefore there exists a ring homomorphism
$
\Phi: \text{\rm Cliff}_\Lambda(n;\vec d) \to (H(L(u);\Lambda),\cup^Q)
$
satisfying Theorem \ref{clifford}.2.
\par
We define a filtration on $\text{\rm Cliff}_\Lambda(n;\vec d)$ by
$$
F^l\text{\rm Cliff}_\Lambda(n;\vec d) = \bigoplus_{\# I \le l} \Lambda X_I.
$$
The associated graded algebra $gr(\text{\rm Cliff}_\Lambda(n;\vec d))$ is the exterior algebra.
We define a filtration on $H(L(u);\Lambda)$ by
$$
F^lH(L(u);\Lambda) = \bigoplus_{k \le l} H^k(L(u);\Lambda).
$$
Lemma \ref{gradedprod} implies that the quantum product $\cup^Q$ respects
this filtration and
$gr(H(L(u);\Lambda))$ is isomorphic to the exterior
algebra as a ring.
\par
Therefore $\Phi$ induces an isomorphism
$gr(\text{\rm Cliff}_\Lambda(n;\vec d)) \to gr(H(L(u);\Lambda))$.
It implies that $\Phi$ itself is an isomorphism.
Hence follows Statement 1.
\par
Statement 4 follows immediately from Lemma  \ref{gradedprod} and the definition.
By Lemma \ref{gradedprod} we have
$$
\aligned
\text{\bf e}'_{\{1,\dots, n\}}
&= \text{\bf e}'_{1} \cup^{\frak c,Q} \dots \cup^{\frak c,Q}\text{\bf e}'_{n}
\equiv \text{\bf e}'_{1} \wedge \dots \wedge \text{\bf e}'_{n} \\
&\equiv \det A \text{\bf e}_{1} \wedge \dots \wedge \text{\bf e}_{n}
= \text{\rm vol}_{L(u)}
\qquad \mod \bigoplus_{d<n} H^d(L(u);\Lambda).
\endaligned$$
Statement 5 follows.
\par
For $I = \{i_1,\dots,i_k\} \in 2^{\{1,\dots,n\}}$ we put
$d^I = d_{i_1}\cdots d_{i_k}$.
By the Clifford relation (\ref{clifrel}), we have
\begin{equation}\label{prodinclif}
\text{\bf e}'_{I} \cup^{\frak c,Q} \text{\bf e}'_{J}
= (-1)^* d^{I\cap J} \text{\bf e}'_{I \ominus J}
\end{equation}
where
$I \ominus J = (I\cup J) \setminus (I\cap J)$.
(We do not calculate the sign $*$ in (\ref{prodinclif}) since we do not need it.)
(\ref{prodinclif}) and Statement 4 imply that
$\text{\bf e}'_{I} \cup^{\frak c,Q} \text{\bf e}'_{J}$ has no degree $n$
component if $J\ne I^c$.
Therefore using (\ref{cupcycl}), we have:
$$
\langle \text{\bf e}'_I, \text{\bf e}'_J \rangle_{\text{\rm PD}_{L(u)}}
= \langle \text{\bf e}'_I, \text{\bf e}'_J \cup^{\frak c,Q} 1\rangle_{\text{\rm PD}_{L(u)}}
=
\langle \text{\bf e}'_I \cup^{\frak c,Q} \text{\bf e}'_J, 1 \rangle_{\text{\rm PD}_{L(u)}}
= \int_{L(u)} \text{\bf e}'_I \cup^{\frak c,Q} \text{\bf e}'_J = 0.
$$
\begin{rem}\label{rem224}
This is another point for which we use the cyclic symmetry in an essential way.
\end{rem}
(\ref{prodinclif}) implies
$$\aligned
\text{\bf e}'_{I} \wedge \text{\bf e}'_{I^c}
\equiv \text{\bf e}'_{I} \cup^{\frak c,Q} \text{\bf e}'_{I^c}
&\equiv (-1)^{*(I)} \text{\bf e}'_{\{1,\dots,n\}} \\
&\equiv (-1)^{*(I)} \text{vol}_{L(u)}
\mod \bigoplus_{d<n} H^d(L(u);\Lambda).
\endaligned$$
Therefore we have
$$
\langle \text{\bf e}'_I, \text{\bf e}'_{I^c} \rangle_{\text{\rm PD}_{L(u)}}
= \int_{L(u)} \text{\bf e}'_{I} \wedge \text{\bf e}'_{I^c} = (-1)^{*(I)}.
$$
Here $*(I)$ is as in (\ref{starI}).
We have thus proved Statement 2.
The proof of Theorem \ref{clifford} is now complete.
\end{proof}
\par
\section{Residue pairing and Hessian determinant}
\label{sec:ResHess}
In this section we complete the proof of Theorem \ref{cliffordZ}.
We first prove Theorem \ref{cliffordZ}.2 and 3.
Using Theorem \ref{clifford}.3 and (\ref{prodinclif}) we calculate
$$
Z((H(L(u);\Lambda),\langle \cdot, \cdot\rangle_{\text{\rm
PD}_{L(u)}}, \cup^{\frak c,\frak b,b},PD[L(u)]))
$$
defined as
(\ref{defnformulaZ}). 
\index{Frobenius algebra!trace}
We consider the nonzero terms in
(\ref{defnformulaZ}). Since $g^{I_3\emptyset}g^{J_3\emptyset} \ne 0$, it follows
from Theorem \ref{clifford}.3 that $I_3 = J_3 = \{1,\dots,n\}$.
Then $\langle \text{\bf e}_{I_1} \cup^{\frak c,Q}  \text{\bf
e}_{I_2},
 \text{\bf e}_{I_3}\rangle \ne 0$  implies
$I_1=I_2$, which we put $I$.
\par
Similarly we have $J_1 = J_2$.
$g^{I_1J_1} \ne 0$ implies $J_1 = J_2 = I^c$.
Thus we have

\begin{equation}\label{Zcalcuform}
\aligned
&Z((H(L(u);\Lambda),\langle \cdot, \cdot\rangle_{\text{\rm PD}_{L(u)}},
\cup^{\frak c,\frak b,b},PD[L(u)])) \\
&= \sum_{I \in 2^{\{1,\dots,n\}}}(-1)^{*_{I,1}}
\langle \text{\bf e}_{I} \cup^{\frak c,Q}  \text{\bf e}_{I},
 \text{\bf e}_{1,\dots,n}\rangle_{PD_{L(u)}}
\langle \text{\bf e}_{I^c} \cup^{\frak c,Q}  \text{\bf e}_{I^c},
 \text{\bf e}_{1,\dots,n}\rangle_{PD_{L(u)}} \\
&= \sum_{I \in 2^{\{1,\dots,n\}}}(-1)^{*_{I,2}}
d_1\dots d_n.
\endaligned
\end{equation}
If $\#I = p$, then
$$
*_{I,1} = p(n-p) + n(n-1)/2.
$$
(This sign comes from the sign $*$ in  (\ref{defnformulaZ}).)
We then have
$$
*_{I,2} = p(p-1)/2 + (n-p)(n-p-1)/2 + *_{I,1} \equiv 0 \mod 2.
$$
Therefore (\ref{Zcalcuform}) is equal to
\begin{equation}\label{Zcalcuform2}
2^n d_1\cdots d_n =
\det\left[ \frac{\partial^2\frak{PO}^{\frak c}_{\frak b}}{\partial x_i\partial x_j}
(x(b))\right]_{i,j=1}^{i,j=n}.
\end{equation}
Thus Theorem \ref{cliffordZ}.3 holds. Theorem \ref{cliffordZ}.2 also follows using
equality $\frak{PO}_{\frak b}^{\frak c} = \frak{PO}_{\frak b}$
((\ref{FcandF}) and Corollary \ref{H1identity}).
\par
We next prove Theorem \ref{cliffordZ}.1.
Let $X$ be a compact toric manifold and $\frak b \in H(X;\Lambda_0)$.
Let $(u,b) \in \frak M(X,\frak b)$.
\par
Here we use $\frak q$-perturbation
(and not $\frak c$-perturbation) to define $\frak M(X,\frak b)$.
We study the operation $\frak m_2^{\frak b,b}$
(and not $\frak m_2^{\frak c,\frak b,b}$) which is
defined by $\frak q$-perturbation at the first step of the proof.
So, at the first step of the proof, we consider $b \in H^1(L(u);\Lambda_0)$.
\par
We assume that $b$ is a nondegenerate critical point of
$\frak{PO}_{\frak b}$.
We define $\cup^Q$ by
$$
h_1 \cup^Q h_2 =  (-1)^{\deg h_1(\deg h_2+1)}\frak m_2^{\frak b,b}(h_1,h_2).
$$
\par
We first prove the following variant of Theorem \ref{clifford}.
\begin{prop}\label{clifford2}
Let $X$ be a compact toric manifold. 
(It does not necessarily satisfy Condition \ref{fano2jigen}.) 
There exist a basis $\{\text{\bf e}'_i \mid i=1,\dots,n\} \subset H^1(L(u);\Lambda)$
 and
a map
$$
\Phi: \text{\rm Cliff}_\Lambda(n;\vec d) \to (H(L(u);\Lambda),\cup^Q)
$$
with the following properties:
\begin{enumerate}
\item  $\Phi$ is a $\Lambda$ algebra isomorphism.
\item  $\Phi(X_i) = \text{\bf e}'_i$.
\item $\{\text{\bf e}'_i \mid i=1,\dots,n\} \subset H^1(L(u);\Lambda_0)$.
Moreover $\{\text{\bf e}'_i \mod  H^1(L(u);\Lambda_+) \mid
i=1,\dots,n \}$
forms a basis of $H^1(L(u);\C)$.
\item
$$
\text{\bf e}'_{\{1,\dots,n\}} - \text{\rm vol}_{L(u)}
\in \bigoplus_{d<n} H^d(L(u);\Lambda_+).
$$
\end{enumerate}
\end{prop}
\par
\begin{proof}
We start with the following lemma.
\begin{lem}\label{linearalgH}
Let $H$ be a symmetric $n\times n$ matrix with $\Lambda_0$ entries.
Then there exist $\text{\bf e}'_i \in \Lambda_0^n$, $i=1,\dots,n$,
with the following properties:
\begin{enumerate}
\item $\{\text{\bf e}'_i \mid i=1,\dots,n\}$ forms a basis of $\Lambda^n$.
\item $\{\text{\bf e}'_i \mod  \Lambda_+\mid i=1,\dots,n\}$ forms a basis of $\C^n$.
\item
$
( \text{\bf e}'_i, H \text{\bf e}'_j )
= 0
$
for $i\ne j$.
Here $( \cdot,\cdot) : \Lambda^n \times \Lambda^n \to
\Lambda$ is the standard quadratic form.
\end{enumerate}
\end{lem}
\begin{proof}
We take $\lambda \ge 0$ the largest number such that
all the entries of $T^{-\lambda}H$ are in $\Lambda_0$.
We may assume $\lambda = 0$ without loss of generality.
\par
Let
$\overline H \equiv H \mod \Lambda_+$ be the
complex $n\times n$ matrix. It is nonzero and symmetric.
We can take a basis $\overline{\text{\bf f}}_i$ of $\C^n$ such that
$
( \overline{\text{\bf f}}_i, \overline H \overline{\text{\bf f}}_j )
= 0
$
for $i\ne j$. We may assume furthermore that
$
( \overline{\text{\bf f}}_i, \overline H \overline{\text{\bf f}}_i )
= 1
$
for $i=1,\dots,k$ and
$
( \overline{\text{\bf f}}_i, \overline H \overline{\text{\bf f}}_i )
= 0
$ for $i=k+1,\dots,n$. ($k>0$.)
It is easy to lift $\overline{\text{\bf f}}_i$ to
${\text{\bf f}}_i \in \Lambda_0^n$ such that
 $
( {\text{\bf f}}_i, H {\text{\bf f}}_j )
= 0
$ holds
for any $i=1,\dots,k$ and $j\ne i$. 
Then, with respect to the splitting
$\Lambda_0^n \cong \Lambda_0^k \oplus \Lambda_0^{n-k}$,
$H$ also splits so that its restriction to the first factor has the
required basis $\text{\bf f}_i$, $i=1,\dots,k$.
Thus the lemma is proved by induction on $n$.
\end{proof}
We apply Lemma \ref{linearalgH} to the Hessian matrix
$\text{\rm Hess}_b(\frak{PO}_{\frak b})$ and obtain
$\text{\bf e}'_i$. We may replace  $\text{\bf e}'_1$ by $c\text{\bf e}'_1$
for some $c \in \Lambda_0 \setminus \Lambda_+$ so that
we have the following in addition:
\begin{equation}\label{detbase1}
\det [\text{\bf e}'_1 \dots \text{\bf e}'_n] = 1.
\end{equation}
(Here we regard each of $\text{\bf e}'_i$ as $1\times n$ matrix.) 
Now Statement 3 is satisfied.
\par
We next observe that the operator $\frak m^{\frak b,b}_2$, that is
constructed by the $\frak q$-perturbation,
satisfies the conclusion of Lemma \ref{gradedprod},
even when Condition $\ref{fano2jigen}$.1 or 2
is {\it not} assumed.
In fact, we can use \cite[Corollary 6.6]{fooo09} 
in place of Lemma \ref{dimandempty}. 
It implies
\begin{equation}\label{dimpluseromod0}
\text{\bf e}'_I  -  \text{\bf e}'_{i_1} \wedge \dots \wedge \text{\bf e}'_{i_k}
\in \bigoplus_{\ell < k}H^{\ell}(L(u);\Lambda_+)
\end{equation}
for $I = (i_1,\dots,i_k)$.
Here we put
$
{\text{\bf e}}'_{I} = {\text{\bf e}}'_{i_1} \cup^{Q} \dots \cup^{Q}
{\text{\bf e}}'_{i_k}.
$
We use the fact that the contribution of pseudo-holomorphic disks
to $\frak m_2$
is of positive energy to show that the difference in (\ref{dimpluseromod0}) is
$0$ modulo $\Lambda_+$.
\par
Lemma \ref{linearalgH}.3 and
(\ref{dimpluseromod0})  imply
$$
\text{\bf e}'_i \cup^{Q} \text{\bf e}'_j
+ \text{\bf e}'_j \cup^{Q} \text{\bf e}'_i =
\begin{cases}
2d_{i}    &\text{if $i=j$} \\
0         &\text{if $i \ne j$},
\end{cases}
$$
for some $d_i \in \Lambda_+$.
Then if we define the map $\Phi$ by $\Phi(X_i) = \text{\bf e}'_i$, Statement
1 follows.
\par
Statement 4 follows from (\ref{detbase1}) and (\ref{dimpluseromod0}).
The proof of Proposition \ref{clifford2} is complete.
\end{proof}

We recall that $d_i$ are defined by \eqref{def:di} and 
the coordinate $\text{\bf e}'_i$ defined by 
\eqref{def:eiprime}.
Namely if
we write
$b \in H^1(L(u);\Lambda_0)$ as
$b = \sum x'_i\text{\bf e}'_i$, then
$$
\det \left[ \frac{\partial^2\frak{PO}^{\frak b}}{\partial x'_i\partial x'_j}
(b)\right]_{i,j=1}^{i,j=n} =
2^nd_1\cdots d_n.
$$
In fact, the matrix in the left hand side is the
diagonal matrix with diagonal entries $d_i$. Then using (\ref{detbase1}) we have
\begin{equation}\label{detxprime}
\det \left[ \frac{\partial^2\frak{PO}^{\frak b}}{\partial x_i\partial x_j}
(b)\right]_{i,j=1}^{i,j=n} =2^n d_1\cdots d_n.
\end{equation}
Here $x_i$ is a coordinate of $H^1(L(u);\Lambda)$ with respect to the
basis $\text{\bf e}_i$.
\begin{rem}
Since we do not have
cyclic symmetry, Proposition \ref{clifford2} does {\it not} imply
$
\langle \text{\bf e}'_I,\text{\bf e}'_J\rangle_{PD_{L(u)}} =0
$
for $J \ne I^c$.
We only have
$$
\langle \text{\bf e}'_I, \text{\bf e}'_J \rangle_{\text{\rm PD}_{L(u)}}
\equiv
\begin{cases}
(-1)^{*(I)} & \mod \Lambda_+  \quad \mbox{for $J = I^c$}, \\
0 & \mod \Lambda_+  \quad  \mbox{\rm otherwise}.
\end{cases}
$$
Here we denote
\begin{equation}\label{starI2}
*(I) = \# \{ (i,j) \in I \times I^c \mid j < i\}.
\end{equation}
\end{rem}
Now we move from $\frak m^{\frak b,b}_2$
to $\frak m_2^{\frak c,\frak b,b^{\frak c}}$,
where  $b^{\frak c}$ is obtained from $b$ as in
Definition \ref{1827}. Then
by Corollary \ref{isowithc} we have a ring isomorphism
$$
\frak F_1 : (H(L(u);\Lambda_0),\cup^{Q})
\to (H(L(u);\Lambda_0),\cup^{\frak c,Q}).
$$
Here we remark that $\cup^{Q}$ is induced from $\frak m_2^{\frak b,b}$
and $\cup^{\frak c,Q}$ is induced from $\frak m_2^{\frak c,\frak b,b^{\frak c}}$.
\par
We now define
\begin{equation}\label{eprimeprime}
\text{\bf e}''_i = \frak F_1(\text{\bf e}'_i)
\end{equation}
and
\begin{equation}\label{eprimeprimeI}
{\text{\bf e}}''_{I} = {\text{\bf e}}''_{i_1} \cup^{\frak c,Q} \dots \cup^{\frak c,Q}
{\text{\bf e}}''_{i_k}.
\end{equation}
Now we have
\begin{prop}\label{clifford3}
There exists
a map
$$
\Psi: \text{\rm Cliff}_\Lambda(n;\vec d) \to (H(L(u);\Lambda),\cup^{\frak c,Q})
$$
with the following properties:
\begin{enumerate}
\item  $\Psi$ is a $\Lambda$ algebra isomorphism.
\item  $\Psi(X_i) = \text{\bf e}''_i$.
\item
$$
\text{\bf e}''_{\{1,\dots, n\}} - \text{\rm vol}_{L(u)}
\in \bigoplus_{d\le n} H^d(L(u);\Lambda_+).
$$
\end{enumerate}
\end{prop}
This is immediate from Proposition \ref{clifford2} and $\frak F_1 - id \equiv 0 \mod \Lambda_+$.
\par
Now we are in the position to complete the proof of Theorem \ref{cliffordZ}.1.
We have
$\frak o\in \Lambda_+$ such that
$$
(1+\frak o) {\rm vol}_{L(u)} - \text{\bf e}''_{\{1,\dots, n\}} \in \bigoplus_{k<n}H^k(L(u);\Lambda).
$$
Using Corollary \ref{degreemustn}
we find that
the formula (\ref{21mainformula}) for
$\frak v = \frak w = (1+\frak o) {\rm vol}_{L(u)}$ is equal
to   (\ref{21mainformula}) for $\frak v = \frak w =\text{\bf e}''_{\{1,\dots, n\}}$.
Namely
\begin{equation}\label{2310}
\aligned
&(1+\frak o)^2\langle i_{\ast,\text{\rm qm},(\frak b,b^{\frak c})}( {\rm vol}_{L(u)}),
i_{\ast,\text{\rm qm},(\frak b,b^{\frak c})}( {\rm vol}_{L(u)})
\rangle_{\text{\rm PD}_{X}} \\
&=
\sum_{I,J \in 2^{\{1,\dots,n\}}} (-1)^{*_{1,I,J}}g^{IJ}
\langle \frak m_2^{\frak c,\frak b,b^{\frak c}}(\text{\bf e}''_I,
\text{\bf e}''_{\{1,\dots, n\}}),
\frak m_2^{\frak c,\frak b,b^{\frak c}}(\text{\bf e}''_J,
\text{\bf e}''_{\{1,\dots, n\}})
\rangle_{\text{\rm PD}_{L(u)}}.
\endaligned\end{equation}
By the cyclic symmetry and Proposition \ref{clifford3},
Formula (\ref{2310}) is equal to
$$
\aligned
&\sum_{I,J \in 2^{\{1,\dots,n\}}} (-1)^{*_{2,I,J}}g^{IJ}
\langle \frak m_2^{\frak c,\frak b,b^{\frak c}}(\text{\bf e}''_{\{1,\dots, n\}},
\frak m_2^{\frak c,\frak b,b^{\frak c}}(\text{\bf e}''_J,
\text{\bf e}''_{\{1,\dots, n\}}))
,\text{\bf e}''_I\rangle_{\text{\rm PD}_{L(u)}} \\
&=
\sum_{I,J \in 2^{\{1,\dots,n\}}} (-1)^{*_{3,I,J}}g^{IJ}
d_{J}
\langle \frak m_2^{\frak c,\frak b,b^{\frak c}}(\text{\bf e}''_{\{1,\dots, n\}},
\text{\bf e}''_{J^c})
,\text{\bf e}''_I\rangle_{\text{\rm PD}_{L(u)}} \\
&=
\sum_{I,J \in 2^{\{1,\dots,n\}}} (-1)^{*_{4,I,J}}g^{IJ}
d_1\cdots d_n
\langle
\text{\bf e}''_{J}
,\text{\bf e}''_I\rangle_{\text{\rm PD}_{L(u)}}
= 2^nd_1\cdots d_n.
\endaligned
$$
Here $d_I = d_{i_1}\dots d_{i_{k}}$
for $I = (i_1,\dots,i_k)$.
\par
The calculation of the sign is as follows.
We put $\vert I\vert = k$, $\vert J\vert = \ell$. Note we only need to calculate the sign in case
$\ell \equiv n-k \mod 2$.
Now
$$\aligned
*_{1,I,J} &= \frac{n(n-1)}{2} \\
*_{2,I,J} &= *_{1,I,J} + \ell(k+1).
\endaligned$$
Here we get the sign $(n+k)(\ell + n + 1)$ to go from $\langle \cdot,\cdot\rangle_{PD_{L(u)}}$
to $\langle \cdot,\cdot\rangle_{\rm cyc}$ and $(k+1)(n+1+n+\ell+1)$ for the cyclic symmetry
and $\ell(k+1)$ to back from $\langle\cdot,\cdot \rangle_{\rm cyc}$ to $\langle  \cdot,\cdot\rangle_{PD_{L(u)}}$ .
The sum of these three is $\ell(k+1)$.
\par
We next have
$$
*_{3,I,J} = *_{2,I,J} + \ell(n+1) + \frac{\ell(\ell-1)}{2} + *(J).
$$
Here $*(J)$ is as in (\ref{starI2}).
In fact, we get $\ell(n+1)$ to go from $\frak m_2$ to the product,
and
$$
 \text{\bf e}''_{J} \cup^{\frak c,Q} \text{\bf e}''_{\{1,\dots, n\}}
= (-1)^{\frac{\ell(\ell-1)}{2} + *(J)} \text{\bf e}''_{J^c}
$$
by Clifford relation.
Similarly we have
$$
*_{4,I,J} = *_{3,I,J} + n(k+1) + \frac{(n-\ell)(n-\ell-1)}{2} + *(J).
$$
We can easily find that $*_{4,I,J} \equiv 0 \mod 2$ by using $k+\ell
\equiv n \mod 2$.
\par
The above calculation and (\ref{detxprime})
imply Theorem \ref{cliffordZ}.1.
\qed
\begin{rem}
We remark that $\text{\bf e}''_i \in H^1(L(u);\Lambda)$ may {\it not} hold.
$\text{\bf e}''_i$, ($i=1,\dots,n$) form a basis of the tangent space
$T_{b^{\frak c}}(\frak F_*(H^1(L(u);\Lambda_0))$.
\end{rem}
\par
\section{Cyclic homology and variation of the invariant $Z$}
\label{sec:cyclic cohomology}
The purpose of this section is to put the invariant $Z(H,\{\frak m_{k}\},\langle \cdot \rangle,e)$
(Definition \ref{invariantZ})
of cyclic filtered $A_{\infty}$ algebra $(H,\{\frak m_{k}\},\langle \cdot \rangle,e)$
in the perspective of several other stories such as cyclic cohomology.
\index{$Z(C)$}\index{trace}\index{Frobenius algebra!trace}
\par
We begin with a brief review of cyclic (co)homology of $A_{\infty}$ algebra which is a minor modification
of the story of cyclic (co)homology of associative algebra.
(\cite{conne}, \cite{loday}.)
\par
Let $K$ be a field of characteristic $0$ and
$(C,\{\frak m_k\}_{k=0}^{\infty},e)$ an $A_{\infty}$ algebra with unit $e$.
Throughout this section we assume the following condition.
\begin{conds}\label{ainfcond}
\begin{enumerate}
\item $C$ is a finite dimensional $K$ vector space.
\item $\frak m_0(1) = ce$ for some $c\in K$.
\item $\frak m_1 = 0$.
\end{enumerate}
Hereafter in this section, we replace $\frak m_0$ by $0$.
Namely we put $\frak m_0(1) = 0$.
Thanks to Item 2 above,  we still have a filtered
$A_{\infty}$ algebra, after this replacement .
\end{conds}
We recall the map  
$\text{\rm cyc}$ on $B_k(C[1])=C[1]^{\otimes k}$ defined by (\ref{defcyc}) and 
$B_k^{\text{\rm cyc}}(C[1])$ which is the quotient space 
of $B_k(C[1])$ defined by Definition 
\ref{def:cyclicB}.
We put
\begin{equation}
B^{\text{\rm cyc}}(C[1]) = \bigoplus_{k=0}^{\infty}B_k^{\text{\rm cyc}}(C[1]).
\end{equation}
The $A_{\infty}$ operation $\frak m_k$ defines a differential
$\delta^{\rm cyc}$ on it by (\ref{hatd}).
\index{cyclic homology}\index{$HC_{\ast}(C)$}\index{cyclic cohomology}\index{$HC^{\ast}(C)$}
\begin{defn}
The {\it cyclic homology} $HC_*(C)$ is the homology group of the
chain complex $(B^{\text{\rm cyc}}(C[1]),\delta^{\rm cyc})$.
\par
The {\it cyclic cohomology} $HC^*(C)$ is the cohomology group of the cochain complex
$(Hom(B^{\text{\rm cyc}}(C[1]),K),(\delta^{\rm cyc})^*)$.
\end{defn}
We put
\begin{equation}\label{hnumberfiltration}
F_kB^{\text{\rm cyc}}(C[1])  =  \bigoplus_{\ell=0}^{k}B_{\ell}^{\text{\rm cyc}}(C[1])
\end{equation}
and
\begin{equation}\label{chnumberfiltration}
F^kHom(B^{\text{\rm cyc}}(C[1]),K)
= \prod_{\ell=k}^{\infty}Hom(B_{\ell}^{\text{\rm cyc}}(C[1]),K).
\end{equation}
By Condition \ref{ainfcond}, that is $\frak m_0 = 0$ and 
$\frak m_1=0$,  we have
$$\aligned
&\delta^{\rm cyc}\left( F_kB^{\text{\rm cyc}}(C[1])\right)
\subset F_{k-1}B^{\text{\rm cyc}}(C[1]), \\
&(\delta^{\rm cyc})^*\left(F^kHom(B^{\text{\rm cyc}}(C[1]),K) \right) \subset
F^{k+1}Hom(B^{\text{\rm cyc}}(C[1]),K) .
\endaligned$$
Therefore we have homomorphisms
$$
H(F_kB^{\text{\rm cyc}}(C[1]),\delta^{\rm cyc}) \to HC_*(C), \quad
H(F^kHom(B^{\text{\rm cyc}}(C[1]),K),(\delta^{\rm cyc})^*) \to HC^*(C).
$$
We denote their images by $F_kHC_*(C)$,
$F^kHC^*(C)$, respectively.
\par
For a unital $A_{\infty}$ algebra $(C,\{\frak m_k\}_{k=0}^{\infty},e)$  satisfying Condition \ref{ainfcond},
$(C,\pm \frak m_2,e)$ becomes an (associative) algebra with unit.
(Here $\pm$ is as in (\ref{deformcup}).)
We call it the {\it underlying algebra}.\index{Clifford algebra}
\begin{prop}[Compare \cite{kassel} Section 6]\label{clifcyccalc}
If the underlying algebra of $(C,\{\frak m_k\}_{k=0}^{\infty},e)$ is isomorphic to the
Clifford algebra associated to a non-degenerate quadratic form, then
\begin{equation}\label{FkFkHCcalcu}
\frac{F_kHC_*(C)}{F_{k-1}HC_*(C)} \cong
\begin{cases}
0  & \text{$k$ is even}, \\
K  & \text{$k$ is odd}.
\end{cases}
\end{equation}
Moreover all the nonzero elements in the cyclic homology $HC_*(C)$ have degree $*
\equiv n-1\mod 2$. (Here we use the shifted degree.)
\end{prop}
\begin{proof}
Following \cite{kassel}, we first consider Hochschild homology. We consider
$BC = \bigoplus_{k=0}^{\infty} B_k(C[1])$ and define
$
\delta^{\rm cyc} : B(C[1]) \to B(C[1])$ by
(\ref{hatd}).
\begin{defn}\label{Hoch}\index{Hochschild homology}\index{$HH_{\ast}(C,C)$}
We call the homology group $H(B(C[1]),\delta^{\rm cyc})$, the
{\it Hochschild homology} and denote it by $HH_*(C,C)$.
\end{defn}
The boundary operator $\delta^{\rm cyc}$ respects the
number filtration $F_kB(C[1])=\bigoplus_{\ell=0}^{k}B_{\ell}(C[1])$. 
Therefore we have a homomorphism
$$
H(F_kB(C[1]),\delta^{\rm cyc}) \to HH_*(C,C).
$$
We denote its image by $F_kHH_*(C,C)$ which defines a filtration on $HH_*(C,C)$. 
\begin{lem}\label{HHvanish}
If the underlying algebra of $C$ is isomorphic to the Clifford algebra associated to a
non-degenerate quadratic form,
then $HH_*(C,C) \cong K$. It is generated by the element
\begin{equation}\label{1eelement}
e_{\{1,\dots, n\}} \in C = B_1(C[1]).
\end{equation}
In particular, $HH_*(C,C) = F_1HH_*(C,C)$.
\end{lem}
\begin{proof}
Let $C'$ be the underlying Clifford algebra regarded as an
$A_{\infty}$ algebra. (Namely the operator $\frak m_2$ of $C'$ is the same as
that of $C$ and all the higher operators $\frak m_k$ of $C'$ are zero.)
We consider the spectral sequence $E_k(C)$ induced by the
number filtration  $F_kB(C[1])$.
It is easy to see that
$E_2(C) \cong E_2(C')$.
On the other hand, as in \cite[page 696]{kassel}, 
$E_2(C')\cong E_{\infty}(C') \cong K$ and it is generated by
the element (\ref{1eelement}).
Therefore by degree reason, the spectral sequence for $C$ 
also degenerates at $E_2$ term.
The lemma follows.
\end{proof}
Now we use the following variant of 
Connes' exact sequence \cite{conne}.\index{Conne's exact sequence}
This $A_{\infty}$ version is given, for example,
in \cite{ChoLee}, \cite{HL}.
\begin{thm}\label{conneseq}
There exists a long exact sequence
\begin{equation}\label{periodexa}
{\longrightarrow} HC_*(C) \overset{S}{\longrightarrow} HC_*(C)
\overset{B}{\longrightarrow} HH_*(C,C)
\overset{I}{\longrightarrow} HC_*(C)  \overset{S}{\longrightarrow}
\end{equation}
such that
\begin{subequations}
\begin{eqnarray}
S(F_kHC_*(C) ) &\subset& F_{k-2}HC_*(C),
\label{exactandfilt1}\\ B(F_kHC_*(C) ) &\subset& F_{k+1}HH_*(C,C),
\label{exactandfilt2}\\
I(F_kHH_*(C,C) ) &\subset& F_{k}HC_*(C).
\label{exactandfilt3}\end{eqnarray}
\end{subequations}
\end{thm}
\begin{proof} 
We consider diagram of the bar complexes:
\begin{equation}\label{diagram0}
\begin{CD}
B(C[1]) @ < {1-\text{\rm cyc}} <<  B(C[1])
@<{N}<<   B(C[1])  @< {1-\text{\rm cyc}} << \dots \\
@ VV{\delta^{\rm cyc}}V @ VV{\widehat d}V  @VV{\delta^{\rm cyc}}V \\
B(C[1]) @ < {1-\text{\rm cyc}} <<  B(C[1])
@<{N}<<   B(C[1])  @< {1-\text{\rm cyc}} << \dots 
\end{CD}
\end{equation}
See \cite[page 54]{loday}. 
Here $N = \sum_{i=0}^{k-1} \text{\rm cyc}^i$
on $B_k(C[1])$. 
Note $\delta^{\rm cyc}$ is defined in (\ref{hatd}).
We denote by $B(C[1])_{\ell}$, $\ell =1,2,\dots$ 
the $\ell$-th column in the diagram. 
We define a map
$$
B(C[1])_{2\ell} \longrightarrow B(C[1])_{2\ell-1}
$$
as $1-\text{\rm cyc}$ and
$$
B(C[1])_{2\ell+1} \longrightarrow B(C[1])_{2\ell}
$$
as $N$. Moreover we define a map
$$
B(C[1])_{2\ell} \longrightarrow B(C[1])_{2\ell}
$$
as $\widehat d$ (see \eqref{def:dhat}) and
$$
B(C[1])_{2\ell-1} \longrightarrow B(C[1])_{2\ell-1}
$$
as  $\delta^{\rm cyc}$. The four assignments combined define
\begin{equation}\label{def:totald}
\widehat{\widehat d} : \bigoplus_{\ell=1}^{\infty}B(C[1])_{\ell} \to \bigoplus_{\ell=1}^{\infty}B(C[1])_{\ell}.
\end{equation}
Let us first check the commutativity of the diagram 
\eqref{diagram0} for completeness sake.
\begin{lem}\label{commutdiagramlem}
\begin{subequations}
\begin{eqnarray}
\delta^{\text{\rm cyc}}\circ (1-\text{\rm cyc})
& = &
(1-\text{\rm cyc}) \circ \widehat d, \label{commdiagram248}\\
\widehat d\circ N
& = &
N \circ \delta^{\text{\rm cyc}}. \label{commdiagram2482}
\end{eqnarray}
\end{subequations}
\end{lem}
\begin{rem}\label{remarkoninvandsub}
The commutativity of the Diagram (\ref{diagram0}) shows the 
equivalence of the two ways to define cyclic bar complex.
In this paper we take the {\it quotient}  of $(B(C[1]),\delta^{\text{\rm cyc}})$
with respect to the cyclic symmetry.
In \cite{fooobook} we take the {\it subset} of $(B(C[1]),\widehat d)$
that is invariant of the cyclic symmetry.
The choice of this paper corresponds to the 
cokernel of $1 - \text{cyc}$ and 
the choice of \cite{fooobook} corresponds 
to the kernel of $N$.
They are isomorphic together with their boundary operators 
(that are $\widehat d$ and $\delta^{\text{\rm cyc}}$, respectively)
by the commutativity of Diagram (\ref{diagram0}).
\end{rem}
\begin{proof}[Proof of Lemma \ref{commutdiagramlem}] 
We prove \eqref{commdiagram248}.
We apply the left hand side to $x_0 \otimes \cdots \otimes x_k$ and obtain
$$
\delta^{\text{\rm cyc}}
(x_0 \otimes \cdots \otimes x_k + (-1)^{*_1}x_k \otimes x_0 \otimes \cdots \otimes x_{k
-1})
$$
that is equal to 
\begin{eqnarray}
&&\!\!\!\!\!\!\!\!\!\!\!\!\!\!\!\!
\sum_{0 \le i \le j \le k} 
(-1)^{*_2}\overset{\times}x_0 \otimes \cdots \frak m_{j-i+1}(x_i,\dots,x_j) 
\cdots \otimes \overset{\times}x_k
\label{lhsfirst}\\
&&\!\!\!\!\!\!\!\!\!\!\!\!\!\!\!\!+
\sum_{0 \le i < j \le k} (-1)^{*_3}
\frak m_{i+k-j+2}(x_j,\dots,x_k,x_0,\dots,x_i) 
\otimes
\cdots \otimes \overset{\times}x_{j-1}
\label{lhssecond}\\
&&\!\!\!\!\!\!\!\!\!\!\!\!\!\!\!\!+
\sum_{0 \le i \le j \le k-1} (-1)^{*_4}
x_k \otimes \cdots \frak m_{j-i+1}(x_i,\dots,x_j) \cdots \otimes \overset{\times}x_{k-1}
\label{lhsthird}\\
&&\!\!\!\!\!\!\!\!\!\!\!\!\!\!\!\!+
\sum_{0  \le i \le k-1} (-1)^{*_5}
\frak m_{i+2}(x_k,x_0,\dots,x_i) \otimes \cdots\otimes 
\overset{\times}x_{k-1} \label{lhsfourth}\\
&&\!\!\!\!\!\!\!\!\!\!\!\!\!\!\!\!+
\sum_{0  \le i \le j \le k-1} (-1)^{*_6}
\frak m_{k-j+i+1}(x_j,\dots,x_k,\overset{\times}x_0,\dots,x_{i-1}) 
\otimes
\cdots \otimes \overset{\times}x_{j-1}.
\label{lhsfifth}\end{eqnarray}
Here the symbol $\times$ which we  put to $\overset{\times}x_k$ in (\ref{lhsfirst}) indicates the fact that 
$x_k$ is absent there in case $j=k$. 
The meaning of the symbol $\overset{\times}x_0$
etc. are similar.
\par
Applying the right hand side of (\ref{commdiagram248}) to 
 $x_0 \otimes \cdots \otimes x_k$ we obtain
$$ 
\sum_{0 \le i \le j \le k} 
(1 - \text{\rm cyc})((-1)^{*_7}
\overset{\times}x_0 \otimes \cdots \frak m_{j-i+1}(x_i,\dots,x_j) \cdots \otimes \overset{\times}x_k)
$$
that is equal to 
 \begin{eqnarray}
&&\!\!\!\!\!\!\!\!\!\!\!\!\!\!\!\!
\sum_{0 \le i \le j \le k} 
(-1)^{*_8}
\overset{\times}x_0 \otimes \cdots \frak m_{j-i+1}(x_i,\dots,x_j) \cdots \otimes \overset{\times}x_k
\label{rhsfirst}\\
&&\!\!\!\!\!\!\!\!\!\!\!\!\!\!\!\! +
\sum_{0 \le i \le j \le k-1} 
(-1)^{*_9}
x_k \otimes x_0 \otimes \cdots  \frak m_{j-i+1}(x_i,\dots,x_j) \cdots \otimes \overset{\times}x_{k-1}
\label{rhssecond}\\
&&\!\!\!\!\!\!\!\!\!\!\!\!\!\!\!\! +
\sum_{0 \le j \le k} 
(-1)^{*_{10}}
\frak m_{k-j+1}(x_j,\dots,x_k) \otimes \overset{\times}x_0 
\otimes \cdots \otimes x_{j-1}.
\label{rhsthird}\end{eqnarray}
\par
We first observe that (\ref{lhsfirst}) = (\ref{rhsfirst}).
In fact, the sign is
$$
*_2 
= 
\deg' x_0 + \dots + \deg'x_{i-1} = *_8.
$$
We next observe that (\ref{lhsthird}) = (\ref{rhssecond}).
In fact, the sign is
$$
*_4 
= 1+
(\deg' x_k)(\deg' x_0 + \dots + \deg' x_{k-1})
+ \deg' x_k + \deg' x_0 + \dots + \deg'x_{i-1} = *_9.
$$
We next observe that the case $j = k$ of (\ref{lhssecond})
cancels with (\ref{lhsfourth}). In fact, the sign is 
$$
*_3 = 
(\deg'x_k)(\deg' x_0 + \dots + \deg'x_{k-1})
= 1+ *_5.
$$
Note the extra minus sign (that is +1) in the right hand side comes from the 
minus sign we put to $1-\text{\rm cyc}$.
\par
We also observe that the case $j \ne  k$ of (\ref{lhssecond})
cancels with the $i \ne  0$ case of (\ref{lhsfifth}). In fact, the sign is
$
*_3 = 
(\deg' x_j + \dots + \deg'x_k)(\deg' x_0 + \dots + \deg'x_{j-1})
$ that becomes 
$
1+  *_6
$ after replacing $i$ by  $i-1$.
\par
Finally we observe that $i = 0$ case of (\ref{lhsfifth})
is equal to (\ref{rhsthird}). The sign is 
$$
*_6 
= 1+ (\deg' x_j + \dots + \deg'x_k)(\deg' x_0 + \cdots \deg'x_{j-1})
= *_{10}.
$$
Thus we have checked \eqref{commdiagram248}.
\par
Next we prove \eqref{commdiagram2482}. 
We apply the left hand side of (\ref{commdiagram2482}) 
to $x_0 \otimes \cdots \otimes x_k$ and obtain  
$\sum_{i=0}^k (-1)^{*_1} \widehat d (x_i \otimes \cdots \otimes x_{i-1})$
that is equal to
\begin{eqnarray}
&&\!\!\!\!\!\!\!\!\!\!\!\!\!\!\!\!  \sum_{0\le i\le j\le \ell \le k}
(-1)^{*_2} \overset{\times}x_i \otimes \cdots  
\frak m_{\ell - j +1}(x_j,\dots,x_{\ell}) \cdots 
\overset{\times}x_0 \cdots \otimes x_{i-1} 
\label{lfs1}
\\
&&\!\!\!\!\!\!\!\!\!\!\!\!\!\!\!\! +  \sum_{0\le \ell < i \le j \le k}
(-1)^{*_3} \overset{\times}x_i \otimes \cdots 
\frak m_{\ell +k - j +2}(x_j,\dots,x_0,\dots,x_{\ell}) \cdots \otimes \overset{\times}x_{i-1} 
\label{lfs2}
\\
&&\!\!\!\!\!\!\!\!\!\!\!\!\!\!\!\!+ \sum_{0\le j\le \ell < i \le k}
(-1)^{*_4} \overset{\times}x_i \otimes \cdots  
\overset{\times}x_0 \cdots
\frak m_{k+\ell - j +2}(x_j,\dots,x_{\ell}) \cdots 
\otimes \overset{\times}x_{i-1}.
\label{lfs3}
\end{eqnarray}
We apply the right hand side of (\ref{commdiagram2482}) 
to $x_0 \otimes \cdots \otimes x_k$ and obtain 
\begin{eqnarray}
&&\sum_{0\le i\le j \le k}
(-1)^{*_5} N(\overset{\times}x_0 \otimes \cdots 
\frak m_{j-i+1}(x_i,\dots,x_j)\dots \otimes\overset{\times}x_k)
\nonumber\\
&& +
\sum_{0\le i\le j \le k}
(-1)^{*_6} N(
\frak m_{i+k-j+2}(x_j,\dots,x_0,\dots,x_i)\otimes\dots \otimes\overset{\times}x_{j-1}).
\nonumber
\end{eqnarray}
It is equal to 
\begin{eqnarray}
&&\!\!\!\!\!\!\!\!\!\!\!\!\!\!\!\!\sum_{0\le i\le j < \ell \le k}
(-1)^{*_7}
x_{\ell} \otimes \cdots \overset{\times}x_0 \cdots 
\frak m_{j-i+1}(x_i,\dots,x_j)\dots \otimes\overset{\times}x_{\ell-1}
\label{rfs1}\\
&&\!\!\!\!\!\!\!\!\!\!\!\!\!\!\!\! +
\sum_{0\le i\le j \le k}
(-1)^{*_8} 
\frak m_{j-i+1}(x_i,\dots,x_j)\otimes\dots \otimes\overset{\times}x_{j-1}
\label{rfs2}
\\
&&\!\!\!\!\!\!\!\!\!\!\!\!\!\!\!\! +
\sum_{0\le \ell < i\le j \le k}
(-1)^{*_9} x_{\ell}\otimes \cdots
\frak m_{j-i+1}(x_i,\dots,x_j)\dots \otimes\overset{\times}x_{\ell-1}
\label{rfs3}
\\
&&\!\!\!\!\!\!\!\!\!\!\!\!\!\!\!\! +
\sum_{0\le i\le j \le k}
(-1)^{*_{10}} 
\frak m_{i+k-j+2}(x_j,\dots,x_0,\dots,x_i)\otimes\dots \otimes\overset{\times}x_{j-1}
\label{rfs4}
\\
&&\!\!\!\!\!\!\!\!\!\!\!\!\!\!\!\! +
\sum_{0\le i <\ell < j \le k}
(-1)^{*_{11}} x_{\ell}\otimes \cdots
\frak m_{i+k-j+2}(x_j,\dots,x_0,\dots,x_i)\dots \otimes\overset{\times}x_{\ell-1}.
\label{rfs5}
\end{eqnarray}
It is easy to see that the (\ref{lfs1}) = (\ref{rfs2}) + (\ref{rfs3}), 
(\ref{lfs2}) = (\ref{rfs4}) + (\ref{rfs5}), 
(\ref{lfs3}) = (\ref{rfs1}).
We note that coincidence of the sign is obvious since we always apply the Koszul rule here.
\par
Therefore the proof of Lemma \ref{commutdiagramlem} is complete.
\end{proof}
The commutativity of the diagram (\ref{diagram0}) implies that 
the operator in \eqref{def:totald} satisfies 
$\widehat{\widehat d}\circ \widehat{\widehat d} = 0$. 
Using the exactness of the horizontal lines of 
\eqref{diagram0} (see \cite[Theorem 2.1.5]{loday}), 
we can prove that
$$
HC(C) \cong  H\left(\bigoplus_{\ell=1}^{\infty}B(C[1])_{\ell}, \widehat{\widehat d} \right).
$$
We also note that the $\widehat d$ cohomology always vanishes.
\par
On the other hand,
we have an exact sequence of chain complecies
$$
0 \to\bigoplus_{\ell=1}^{2}B(C[1])_{\ell} \to\bigoplus_{\ell=1}^{\infty}B(C[1])_{\ell}
\to \bigoplus_{\ell=3}^{\infty}B(C[1])_{\ell}  \to 0.
$$
The homology of the second complex is the cyclic homology.
The homology of the third complex is the cyclic homology with
the degree shifted by $2$.
The homology of the first complex is the homology of $\delta^{\rm cyc}$ complex,
that is the Hochschild homology.
We thus obtain the exact sequence 
(\ref{periodexa}).
\par
We finally check the filtration properties 
(\ref{exactandfilt1}), (\ref{exactandfilt2}), (\ref{exactandfilt3}). 
We note the following diagram concerning filtrations 
defined as in \eqref{hnumberfiltration}. 
\begin{equation}\label{diagram}
\begin{CD}
@ VVV @ VVV  @VVV \\
F_3B(C[1]) @ < {1-\text{\rm cyc}} <<  F_3B(C[1])
@<{N}<<   F_3B(C[1])  @< {1-\text{\rm cyc}} << \dots \\
@ VV{\delta^{\rm cyc}}V @ VV{\widehat d}V  @VV{\delta^{\rm cyc}}V\\
F_2B(C[1]) @ < {1-\text{\rm cyc}} <<  F_2B(C[1])
@<{N}<<   F_2B(C[1])  @< {1-\text{\rm cyc}} << \dots \\
@ VV{\delta^{\rm cyc}}V @ VV{\widehat d}V  @VV{\delta^{\rm cyc}}V \\
C[1] @ < {1-\text{\rm cyc}} <<  C[1]
@<{N}<<   C[1]  @< {1-\text{\rm cyc}} << \dots \end{CD}
\end{equation}
We put 
$$
F_n (B_*(C[1])_*) = \bigoplus_{k+\ell \le n+1}B_k(C[1])_{\ell}.
$$
Here $B_k(C[1])_{\ell}$ is $B_k(C[1])=C[1]^{\otimes k}$ in the $\ell$-th column of the diagram \eqref{diagram0}. 
It defines a filtration on the total complex 
$B_*(C[1])_* = \bigoplus_{k,\ell}B_k(C[1])_{\ell}$. 
We define a filtration on homology $H(B_*(C[1]_*))$ by 
$$
F_n(H(B_*(C[1])_*)) = \text{\rm Image}(H( F_n(B_*(C[1])_*)) \to H(B_*(C[1])_*)).
$$
We recall that we have a filtration 
on $B^{\text{\rm cyc}}(C[1])$ defined by  
$$
F_n(B^{\text{\rm cyc}}(C[1]))
= \bigoplus_{k\le n}B^{\text{\rm cyc}}_k(C[1])
$$
and a filtration on the cyclic homology defined by 
$$
F_nHC(C) = \text{\rm Image}(H( F_n(B^{\rm{cyc}}(C[1])),\delta^{\rm cyc}) \to HC(C)).
$$
The chain map
$$
\pi : B_*(C[1])_* \to B^{\text{\rm cyc}}(C[1])
$$
given by $\pi((x_{k,\ell})) = [x_{k,1}]$ 
induces an isomorphism on the homology, 
since the horizontal lines of (\ref{diagram}) 
is exact and the cokernel of the first homomorphism 
is $B^{\text{\rm cyc}}(C[1])$.
Therefore, 
to prove (\ref{exactandfilt1}), (\ref{exactandfilt2}), (\ref{exactandfilt3})
it suffices to show the following.
\begin{lem}
$$
\pi_*(F_n(H(B_*(C[1])_*)))
= F_nHC(C).
$$
\end{lem}
\begin{proof}
Since $\pi$ preserves the filtrations, we have
$$\pi_*(F_n(H(B_*(C[1])_*)))
\subseteq F_nHC(C).
$$
We prove the inclusion of the opposite direction.
Let $(\overline x_k)_{k=1,\dots,n} \in B^{\text{\rm cyc}}(C[1])$
be a cycle.
We take $(x_k)_{k=1,\dots,n} \in B(C[1])$ whose equivalence 
class is  $(\overline x_k)$.
We put $x_{k,1} = x_k$.
Since $(\overline x_k)_{k=1,\dots,n}$ is a cycle,
we have 
$$
\delta^{\text{cyc}}(x_{k,1})_{k=1,\dots,n} 
\in \text{\rm Image}(1 - \text{cyc}).
$$
We take 
$(x_{k,2})_{k=1,\dots,n-1}$ such that
$$
(1 - \text{cyc})((x_{k,2})_{k=1,\dots,n-1})
=
\delta^{\text{cyc}}(x_{k,1})_{k=1,\dots,n}. 
$$
By commutativity of (\ref{diagram})
and the exactness of its horizontal lines we have 
an element 
$(x_{k,3})_{k=1,\dots,n-2}$ such that
$$
N((x_{k,3})_{k=1,\dots,n-2})
=
\widehat d(x_{k,2})_{k=1,\dots,n-1}. 
$$
So we find by induction an element 
$(x_{k,\ell})_{k+\ell \le n+1}$ such that
it is a cycle of $B_*(C[1])_*$ 
and $\pi((x_{k,\ell})_{k+\ell \le n+1})= (\overline x_k)_{k=1,\dots,n}$.
This proves the lemma.
\end{proof}
Therefore the proof of Theorem \ref{conneseq} is complete.
\end{proof}
Now we are ready to complete the proof of Proposition
\ref{clifcyccalc}.
By Lemma \ref{HHvanish} we have
$F_1HH_*(C,C) \cong K$.
By (\ref{periodexa}) and (\ref{exactandfilt2}), 
$I : HH_*(C,C) \to HC_*(C)$ is injective on $F_1HH_*(C,C)$.
Therefore the generator of $F_1HH_*(C,C)$ induces a nontrivial element of
$HH_*(C,C)$.
Then using $HH_*(C,C) \cong F_1HH_*(C,C) \cong K$ and  (\ref{periodexa})
we can prove that
$S$ induces an isomorphism
$$F_kHC_*(C)/F_{k-1}HC_*(C) \cong F_{k-2}HC_*(C)/F_{k-3}HC_*(C)$$
for $k\ge 3$.  (\ref{FkFkHCcalcu}) follows easily.
\par
We note that
$$
\deg' e_{\{1,\dots, n\}} = \deg  e_{\{1,\dots, n\}} - 1 = n-1.
$$
Moreover $S$ does not change the parity of the degree. The statement on degree follows.
\end{proof}
\par
It is well known among the experts that cyclic cohomology
controls deformations of cyclic $A_{\infty}$ algebras.
Explanation of this is now in order.
(In a similar way, deformations of $A_{\infty}$ algebras
are controlled by Hochschild cohomology. See \cite[Subsection 7.4.6]{fooobook2}, for example.)
\par
Let $(C,\{\frak m_k\},e)$ be a unital $A_{\infty}$ algebra
satisfying Condition \ref{ainfcond}.
Let $\langle ~\cdot ~\rangle : C^{k} \otimes C^{n-k} \to K$
be a nondegenerate bilinear form such that
\begin{equation}\label{mcyccond}
\aligned
&\langle \frak m_k(x_1,\dots,x_k), x_0  \rangle \\
&= (-1)^{\deg'x_0(\deg'x_1 + \dots + \deg' x_k)}
\langle \frak m_k(x_0,x_1,\dots,x_{k-1}), x_k \rangle.
\endaligned\end{equation}
In other words, $(C,\{\frak m_k\},\langle ~\cdot ~ \rangle,e)$ is a {\it cyclic}
unital $A_{\infty}$ algebra.
\par
If $\frak m_k$ is defined only for $k\le N$ and if the $A_{\infty}$ relation
\begin{equation}\label{Ainfityrelation}
\sum_{k_1,k_2 \atop k_1+k_2=k+1} \sum_i
(-1)^{\deg' x_1 + \dots + \deg'x_{i-1}} \frak m_{k_1}(x_1,\dots,\frak m_{k_2}(x_i,\dots),\dots, x_k) = 0
\end{equation}
is satisfied only for $k\le N$, we call it
a {\it cyclic unital $A_N$ algebra}.
\par
We next consider the ring
\begin{equation}
K_m
=
\begin{cases}
\displaystyle \frac{K[\theta]}{(\theta^{m+1})} &\text{if $m < \infty$}, \\
K[[\theta]]  &\text{if $m = \infty$}.
\end{cases}
\end{equation}
\begin{defn}
An {\it $A_{N}$ deformation of order $m$} of  $(C,\{\frak m_k\},\langle \cdot \rangle)$
is operations $\frak m^{\theta}_k$ on $C \otimes K_m$ such that
$(\ref{mcyccond})$ and $(\ref{Ainfityrelation})$ hold for $k \le N$ and that
$\frak m^{\theta}_k \equiv \frak m_k \mod \theta$.
\end{defn}
\begin{rem}\label{rem:A3deformation}
\begin{enumerate}
\item
We do not deform $\langle \cdot \rangle$ or $e$ but only the operations 
$\frak m_k$.
\item
The associative law in algebras is the $A_3$ relation for the product $\frak m_2$. Hence an $A_3$ deformation gives a deformation of algebras. 
\end{enumerate}
\end{rem}
A {\it pseudo-isotopy of cyclic $A_{\infty}$ algebras} 
\index{pseudo-isotopy of cyclic $A_{\infty}$ algebras} 
is defined by a family of operators
$$
\frak m^t_{k} : B_k(C[1]) \to C[1],
\quad
\frak c^t_{k} : B_k(C[1]) \to C[1]
$$
parametrized by $t \in [0,1]$
such that:
\begin{enumerate}
\item They are smooth with respect to $t$ in the sense of Definition \ref{smoothmtandct}.
\item For each $t$, $(C,\{\frak m^t_k\},\langle \cdot \rangle)$
is a cyclic $A_{\infty}$ algebra.
\item For each $x_i \in \overline{C}[1]$ (where $\overline{C}[1]$ is the $\C \cong \Lambda_0 /\Lambda_+$ reduction of $C[1]$)
\begin{equation}\label{inhomomainformula}
\langle
\frak c^t_{k,\beta}(x_1,\ldots,x_k),x_0
\rangle
=
(-1)^*\langle
\frak c^t_{k,\beta}(x_0,x_1,\ldots,x_{k-1}),x_{k}
\rangle.
\end{equation}
Here $* = \deg' x_0(\deg' x_1+\dots + \deg' x_k)$.
\item  For each $x_i \in \overline{C}[1]$
\begin{equation}\label{isotopymaineq}
\aligned
&\frac{d}{dt} \frak m_{k,\beta}^t(x_1,\ldots,x_k) \\
&+ \sum_{k_1+k_2=k+1}\sum_{\beta_1+\beta_2=\beta}\sum_{i=1}^{k-k_2+1}
(-1)^{*}\frak c^t_{k_1,\beta_1}(x_1,\ldots, \frak m_{k_2,\beta_2}^t(x_i,\ldots),\ldots,x_k) \\
&- \sum_{k_1+k_2=k+1}\sum_{\beta_1+\beta_2=\beta}\sum_{i=1}^{k-k_2+1}
\frak m^t_{k_1,\beta_1}(x_1,\ldots, \frak c_{k_2,\beta_2}^t(x_i,\ldots),\ldots,x_k)\\
&=0.
\endaligned
\end{equation}
Here $* = \deg' x_1 + \dots + \deg'x_{i-1}$.
\end{enumerate}
We can define a pseudo-isotopy of cyclic $A_N$ structures and
also of an algebra over $K_m$ in the same way.
\begin{defn}
Two $A_N$ deformations of order $m$ are said to be
{\it pseudo-isotopic} if there exists a pseudo-isotopy of cyclic $A_N$ structures over $K_m$
such that its operators $\frak m^t_k$ and $\frak c^t_k$ satisfy:
$$
\frac{d \frak m^t_k}{dt} \equiv 0 \mod \theta,
\quad \frak c^t_k \equiv 0 \mod \theta.
$$
\end{defn}
We can prove that a pseudo-isotopy is an equivalence relation in the same way
as in \cite[Lemma 8.2]{fooo091}.
\begin{prop}\label{classfysphoe}
 The set of
pseudo-isotopy classes of $A_N$ deformations of order 1
of $(C,\{\frak m_k\},\langle \cdot \rangle)$ corresponds one to one
to the set
$$
\aligned
&HC^*(C)/F^{N+2}HC^*(C)
\quad \text{(if $N<\infty$)},
\\
&\projlim_{m}HC^*(C)/F^{m+2}HC^*(C)
\quad \text{(if $N=\infty$)},
\endaligned$$
with
$
* \equiv n + 1 \mod 2.
$
\end{prop}
\begin{proof}
Let $\frak m_k^{\theta} = \frak m_k + \theta \Delta\frak m_k$  be an $A_N$ deformation of
order one.
We put
\begin{equation}\label{def:plus}
\Delta^+\frak m_k(x_1,\dots,x_k,x_0)
= \langle \Delta\frak m_k(x_1,\dots,x_k),x_0 \rangle .
\end{equation}
Then
\begin{equation}\label{deformationelement}
(\Delta^+\frak m_0,\dots,\Delta^+\frak m_{N})
\in \frac{\prod_{k=0}^{\infty} Hom(B_k^{\text{\rm cyc}}(C[1]),K)}
{\prod_{k=N+2}^{\infty} Hom(B_k^{\text{\rm cyc}}(C[1]),K)}.
\end{equation}
The
$A_{N}$ relation among them implies that
the left hand side of (\ref{deformationelement}) is a $(\delta^{\rm cyc})^*$ cocycle.
\par
Suppose $\frak m_k + \theta \Delta\frak m_k$ is pseudo-isotopic to
$\frak m_k + \theta \Delta'\frak m_k$. Take a
pseudo-isotopy $\frak m_k^t$, $\frak c_k^t$ among them.
Since $\frak c_k^t \equiv 0 \mod \theta$ we have
$\frak c_k^t = \overline{\frak c}_k^t \theta$. We put
$$
\frak e_k = \int_{0}^1  \overline{\frak c}_k^t dt,
$$
and define $\frak e^+_k$ 
in the same way as in \eqref{def:plus}.
It follows easily from (\ref{isotopymaineq}) that
$$
(\Delta^+\frak m_0,\dots,\Delta^+\frak m_N)
- (\Delta^{\prime +}\frak m_0,\dots,\Delta^{\prime +}\frak m_N)
= (\delta^{\rm cyc})^*(\frak e^+_0,\dots,\frak e^+_N).
$$
\par
We have thus defined a homomorphism to $HC^*(C)/F^{N+2}HC^*(C)$ from
the set of pseudo-isotopy classes of $A_N$ deformations of $(C,\{\frak m_k\},\langle \cdot \rangle)$ of order 1.
It is straightforward to check that this map gives
the required isomorphism. 
The proof in the case $N=\infty$ is similar.
\par
We note that if $\Delta^+\frak m_k(x_1,\dots,x_k,x_0) \ne 0$ then
$$
\sum_{i=1}^k \deg' x_i +1 + \deg' x_0 \equiv n+2 \mod 2.
$$
Therefore $\deg \Delta^+\frak m_k \equiv n+1 \mod 2$.
\end{proof}
\begin{prop}\label{liftdef1}
Let $(C,\{\frak m_k\},\langle \cdot \rangle,e)$ be a
unital and cyclic $A_{\infty}$ algebra 
satisfying Condition \ref{ainfcond}. 
Suppose that the underling 
algebra is isomorphic to the Clifford algebra 
associated to a non-degenerate quadratic form. 
Let $\alpha \in HC^*(C)$.
If $[\alpha] \in HC^*(C)/F^4HC^*(C)$ lifts to an $A_{3}$ deformation
of order 1, then
$\alpha$ lifts to an $A_{\infty}$ deformation
of order $\infty$.
\end{prop}
\begin{proof}
The proof is by the obstruction theory, similar
to, for example, \cite[Lemma 7.2.74]{fooobook2}.
We put $\alpha^+ = (\alpha_k^+)_{k=1}^{\infty}$,
$\alpha_k^+ = \sum_{j=1}^{\infty} \theta^{j}\alpha_{k,j}^+$,
$\alpha_{k,j}^+\in Hom(B_k^{\text{\rm cyc}}C,K)$.
We are given $\alpha_{3,1}^+$ and will find other $\alpha_{k,j}^+$.
The degree of $\alpha_{k,\ell}^+$ is $n-1$ modulo $2$.
\par
We put 
$$
B^{\text{cyc} \ast}_k C = Hom(B_k^{\text{\rm cyc}}C[1],K), 
\quad 
B^{\text{cyc} \ast} C = \prod_{k=0}^{\infty} Hom(B_k^{\text{\rm cyc}}C[1],K).
$$
On 
$
B^{\text{cyc} \ast}C
$
we have an IBL (involutive bi-Lie) algebra structure
(see \cite{CFL}).
As a part of the structure, we have an operation  
$$
[\cdot ,\cdot ]: B_{k_1}^{\text{cyc} \ast}C \otimes 
B_{k_2}^{\text{cyc} \ast}C \to 
B_{k_1+k_2-2}^{\text{cyc} \ast}C
$$ defined by
$$
\aligned
& [\alpha_{k_1}^+, \alpha_{k_2}^+](x_1,\dots ,x_{k_1+k_2-2}) \\
& = 
\sum_{a,b}\sum_{c=1}^{k_1+k_2-2}
\pm g^{ab}\alpha_{k_1}^+(e_a,x_c,\dots ,x_{c+k_1-2})
\alpha_{k_2}^+(e_b,x_{c+k_1-1},\dots ,x_{c-1})
\endaligned$$
for $\alpha_i^+ \in B_{k_i}^{\text{cyc}\ast}C$. 
Here $\{e_a \}$ is a $K$-basis of $C$ and 
$(g^{ab})$ is the transpose of the inverse matrix of 
$(g_{ab})=(\langle e_a, e_b\rangle )$ and we use the convention 
$x_{c}=x_{c'}$ if $c\equiv c' \mod k_1+k_2-2$. 
The operation $[\cdot ,\cdot ]$ does not depend on the choice of the basis of $C$.
If $\alpha_{k}^+ =\frak m_{k-1}^+$, the $A_{\infty}$ formula together with the cyclic symmetry of $\langle \cdot ,\cdot \rangle$ in the sense of 
\eqref{mcyccond} 
yields the following equation
\begin{equation}\label{LieandcyclicAinfty}
(\delta^{\rm cyc})^{\ast} \alpha^+_k  \pm
\frac{1}{2}
\sum_{k_1+k_2=k+2}[\alpha_{k_1}^+, \alpha_{k_2}^+]=0.
\end{equation}
See  \cite[Proposition 11.2]{CFL}. 
\par
Under the general facts as above understood, we go back to the proof of Proposition \ref{liftdef1}.
For each step we consider
$
Hom(B_k^{\text{\rm cyc}}C,K)
$
part of
$$
\sum_{k_1+k_2=k+2}\sum_{j_1+j_2=j} [\alpha_{k_1,j_1}^+, \alpha_{k_2,j_2}^+].
$$
By \eqref{LieandcyclicAinfty} this 
is a $(\delta^{\rm cyc})^{\ast}$ cocycle of degree $\equiv n \mod 2$.
Namely the obstruction to find them
inductively lies in the group $F^{k}HC^{*}(C)/F^{k+1}HC^*(C)$ with $* \equiv n \mod 2$.
This group vanishes by Proposition \ref{clifcyccalc}. Hence we finish the proof.
\end{proof}
\begin{rem}
Two cyclic $A_{\infty}$ algebras which are pseudo-isotopic
are homotopy equivalent as cyclic $A_{\infty}$ algebras. (This is proved in \cite[Theorem 8.2]{fooo091}.)
However the homotopy equivalence may not be strict
in the sense of \cite[Definition 3.2.29 (2)]{fooobook}.
In fact, we have the following example.
\par
Let us consider $C = \Lambda e_0 \oplus \Lambda e_1$ with
$\deg e_i = i$. Here $e_0$ is a strict unit and $\frak m_2(e_1,e_1) = e_0$.
We put $\langle e_0, e_1\rangle =1$ and $\frak m_k = 0$ for $k\ge 3$.
\par
Let $e^{i_1,\dots,i_k} \in Hom(B_kC,K)$ ($i_1,\dots,i_k \in \{0,1\}$)
be the elements such that $e^{i_1,\dots,i_k}(e_{i_1}\otimes \cdots \otimes e_{i_k}) = 1$
and that  $e^{i_1,\dots,i_k}$ are zero on other basis.
We put
$$
e^{(i_1,\dots,i_k)} = \sum_{j=1}^k \pm e^{i_j,\dots,i_k,i_1,\dots,i_{j-1}}.
$$
Here $\pm$ is chosen so that this element is cyclically symmetric.
It is easy to see that
$
(\delta^{\rm cyc})^{\ast} e^{(11)} = 0.
$
However $e^{(11)}$ is zero in cyclic cohomology. Namely
$
(\delta^{\rm cyc})^{\ast} e^{(0)} = 2e^{(11)}.
$
\par
This element $ e^{(11)}$ corresponds to the deformation of the
boundary operator $\frak m_1$ to $\frak m_1^{\theta}$ where
$
\frak m_1^{\theta}(e_1) = \theta/2\cdot e_1.
$
Proposition \ref{classfysphoe} implies that $(C,\frak m)$ is pseudo-isotopic to  $(C,\frak m^{\theta})$.
(Here pseudo-isotopy includes nonzero $\frak c_0^t$.)
\par
The (nonstrict) isomorphism between $(C,\frak m)$ and $(C,\frak m^{\theta})$
is given by $\varphi$ such that
$
\varphi_0(1) = \theta e_0.
$
Note $(C,\frak m)$ is not strictly homotopic to  $(C,\frak m^{\theta})$, since 
$H(C,\frak m_1) \neq H(C,\frak m_1^{\theta})$.
The pair $(C,\frak m)$,  $(C,\frak m^{\theta})$ is an example 
of pseudo-isopotic but not strictly homotopic pair of cyclic filtered $A_{\infty}$ algebras.
\par
We also note that in the above example we have
$
(\delta^{\rm cyc})^{\ast} e^{(01)} = -e^{(001)} + 3e^{(111)},
$
and $e^{(111)}$ gives a nonzero element in cyclic cohomology.
This corresponds to the deformation to $\frak m_2^{\prime \theta}$ where
$$
\frak m_2^{\prime \theta}(e_1,e_1) = e_0 + \theta e_0.
$$
We can show the nontriviality of this deformation by using, for example, the trace $Z$.
\par
On the other hand, the element $e^{(1111)}$ is a $(\delta^{\rm cyc})^{\ast}$ coboundary and
is so zero in cyclic cohomology. (This follows from Proposition \ref{clifcyccalc}.
We can also check it directly.)
\end{rem}

We remark that the first nontrivial element of the cyclic cohomology
lies in $F^{1}HC^*(C)/F^{2}HC^*(C)$. It deforms the operator $\frak m_0(1)$, that is
the critical value.
(This comes from Hochschild cohomology and so is an invariant of
$A_{\infty}$ algebra, whether or not it carries a compatible inner product.)
\par
We next discuss the relationship between the trace $Z(C,\{\frak m_k\},\langle~ \cdot~ \rangle,e)$
and cyclic homology of $(C,\{\frak m_k\},\langle ~\cdot~ \rangle,e)$.
This relationship becomes most transparent in case the following condition is satisfied.
\begin{defn}
Let $(C,\{\frak m_k\},\langle ~\cdot~ \rangle,e)$ be a
unital and cyclic $A_{\infty}$ algebra satisfying
Condition $\ref{ainfcond}$. 
We say that $\text{\bf e}_1,\dots,\text{\bf e}_n$ is a 
{\it cyclic Clifford basis}\index{cyclic Clifford basis}
if the following conditions are satisfied.
\begin{enumerate}
\item
$$
\text{\bf e}_i \cup^Q\text{\bf e}_j + \text{\bf e}_j \cup^Q\text{\bf e}_i =
\begin{cases}
2d_i  1 & i=j \\
0    & i \ne j,
\end{cases}
$$
where $d_i \in K \setminus 0$.
\item
$\text{\bf e}_I$, $I \in 2^{\{1,\dots,n\}}$ is a basis of $C$. Moreover
$$
\langle \text{\bf e}_I, \text{\bf e}_J \rangle
=
\begin{cases}
(-1)^*  & J=I^c \\
0    & \text{otherwise}.
\end{cases}
$$
Here $\text{\bf e}_I = \text{\bf e}_{i_1}\cup^Q \dots \cup^Q\text{\bf e}_{i_k}$
if $I = \{i_1,\dots,i_k\}$.
$* = \#\{ (i,j) \in I \times J \mid j<i\}$.
\end{enumerate}
\end{defn}
The basis $\{\text{\bf e}'_i\}$ in Theorem \ref{clifford} is a cyclic Clifford basis.
\par
\begin{lem}\label{Zclifordbases}
If $\{\text{\bf e}_1,\dots,\text{\bf e}_n\}$ is a cyclic Clifford basis of
$(C,\{\frak m_k\},\langle ~\cdot~ \rangle,e)$, then we have
$$
Z(C,\{\frak m_k\},\langle ~\cdot~ \rangle,e) = 2^n d_1\cdots d_n.
$$
\end{lem}
The proof is the same as (\ref{Zcalcuform}), (\ref{Zcalcuform2}).
\par
Let $\{\frak m^{\theta}_k\}$ be an $A_3$ deformation of
$(C,\{\frak m_k\},\langle ~\cdot~ \rangle,e)$ of order $1$  and consider
$$
Z(C,\{\frak m^{\theta}_k\},\langle ~\cdot~ \rangle,e) \in K_1.
$$
We put
$$
Z(C,\{\frak m^{\theta}_k\},\langle ~\cdot~ \rangle,e)
= Z(C,\{\frak m_k\},\langle ~\cdot~ \rangle,e)
+ \theta\Delta Z
$$
where $\Delta Z \in K$. We have thus defined 
\begin{equation}\label{deltaZ}
\Delta Z \in
Hom\left(\frac{HC^*(C)}{F^5HC^*(C)},K\right)
\cong Hom\left(\frac{HC^*(C)}{F^4HC^*(C)},K\right).
\end{equation}
Here we note that $F^4HC^*(C)/F^5HC^*(C)=0$.
Actually, the cohomological version of 
Proposition \ref{clifcyccalc} shows that 
\begin{equation}\label{FkFkHCcalcu2}
\frac{F^kHC^*(C)}{F^{k+1}HC^*(C)} \cong
\begin{cases}
0  & \text{$k$ is even}, \\
K  & \text{$k$ is odd}.
\end{cases}
\end{equation}
\begin{rem}
Here we use the fact that any isomorphism between unital algebras
preserves unit, automatically.
\end{rem}
\begin{thm}\label{Zandcyclic}
Let $(C,\{\frak m^{\theta}_k\},\langle \cdot \rangle,e)$ be a
unital cyclic
$A_{\infty}$ algebra that admits a cyclic Clifford basis.
Then the homomorphism
$$\Delta Z : \frac{F^3HC^*(C)}{F^5HC^*(C)}
\cong \frac{F^3HC^*(C)}{F^4HC^*(C)} \to K$$
is an isomorphism.
\end{thm}
\begin{proof}
We consider the following deformation:
We define $\frak m_2^{\theta}$ such that
the induced product structure $\cup^{Q, \theta}$ is 
$$
\text{\bf e}_i \cup^{Q, \theta}\text{\bf e}_j +\text{\bf e}_j  \cup^{Q, \theta}\text{\bf e}_i =
\begin{cases}
(2d_1 + \theta)1  & i=j=1 \\
2d_i 1  & i=j\ne1 \\
0    & i \ne j.
\end{cases}
$$
This is an $A_3$ deformation. 
So it lifts to a deformation of $A_{\infty}$
structure by Proposition
\ref{liftdef1}.
Since we do not deform $\frak m_0$ and $\frak m_1$, 
this deformation gives an element of
$$
X \in
\frac{F^3HC^*(C)}{F^4HC^*(C)}.
$$
By Lemma \ref{Zclifordbases}
$$
Z(C,\{\frak m^{\theta}_k\},\langle \cdot \rangle,e)
= 2^n(d_1+\theta)d_2\cdots d_n.
$$
Therefore
$$
\Delta Z(X) = 2^nd_2\cdots d_n \ne 0.
$$
Namely $\Delta Z \ne 0$.
Since  $\frac{F^3HC^*(C)}{F^4HC^*(C)} \cong K$ in our case,
we have proved Theorem \ref{Zandcyclic}.
\end{proof}
Going back to our geometric situation, we have the following:
\begin{thm}
If $b$ is a nondegenerate critical point of $\frak{PO}_{\frak b}$
which corresponds to $b^{\frak c}$ by Definition \ref{1827},
then
$(H(L(u);\Lambda),\{\frak m^{\frak c,\frak b,b^{\frak c}}_k\},\langle \cdot \rangle,e)$
has a cyclic Clifford basis.
\end{thm}
\begin{proof}
We will prove that the basis $\text{\bf e}''_i$  in Proposition \ref{clifford3}
is a cyclic Clifford basis up to scalar multiplication.
The next proposition is the main part of its proof.
\begin{prop}\label{2417}
If $J \ne I^c$, then $\langle \text{\bf e}''_I,\text{\bf e}''_J \rangle_{PD_{L(u)}} = 0$.
\end{prop}
\begin{proof}
We prove the following lemma.
\begin{lem}\label{2418}
If $K \subset 2^{\{1,\dots,n\}}$, $K \ne \{1,\dots,n\}$, then
\begin{equation}\label{lemmaindd}
\sum_{I,J \in 2^{\{1,\dots,n\}}} (-1)^{*}g^{IJ}
\langle \frak m_2^{\frak c,\frak b,b^{\frak c}}(\text{\bf e}''_I,
\text{\bf e}''_K),
\frak m_2^{\frak c,\frak b,b^{\frak c}}(\text{\bf e}''_J,\text{\bf e}''_{\{1,\dots, n\}})
\rangle_{\text{\rm PD}_{L(u)}} = 0.
\end{equation}
Here $* = n(n-1)/2$.
\end{lem}
\begin{proof}
We calculate the left hand side of (\ref{lemmaindd}) using cyclic symmetry to obtain
\begin{equation}\label{lemmaindd2}
\sum_{I,J \in 2^{\{1,\dots,n\}}} (-1)^{*}g^{IJ}
\langle \frak m_2^{\frak c,\frak b,b^{\frak c}}(
\text{\bf e}''_K,\frak m_2^{\frak c,\frak b,b^{\frak c}}
(\text{\bf e}''_J,\text{\bf e}''_{\{1,\dots, n\}})),
\text{\bf e}''_I\rangle_{\text{\rm PD}_{L(u)}}.
\end{equation}
We do not calculate the sign here but only remark that,
for given $J,K,n$, it depends only on the parity of $\vert I\vert$.
\par
We put
$$
\frak m_2^{\frak c,\frak b,b^{\frak c}}(
\text{\bf e}''_K,\frak m_2^{\frak c,\frak b,b^{\frak c}}(\text{\bf e}''_J,
\text{\bf e}''_{\{1,\dots, n\}}))
= \sum_{S\in 2^{\{1,\dots,n\}}} (-1)^*h_S(J,K) \text{\bf e}''_S.
$$
Then (\ref{lemmaindd2}) is equal to
$$
\sum_{I,J} (-1)^* h_S(J,K) g^{IJ}g_{SI}
= \sum_{J}  (-1)^*h_J(J,K).
$$
Here we use the fact that the sum is taken for $I$ such that the parity of $\vert I\vert$ is
the same as $n-\vert J\vert$ \par
On the other hand,
$$
\frak m_2^{\frak c,\frak b,b^{\frak c}}(
\text{\bf e}''_K,\frak m_2^{\frak c,\frak b,b^{\frak c}}(\text{\bf e}''_J,
\text{\bf e}''_{\{1,\dots, n\}}))
= \pm c \frak m_2^{\frak c,\frak b,b^{\frak c}}(\text{\bf e}''_K,\text{\bf e}''_{J^c})
= \pm  c'\text{\bf e}''_{K \ominus J^c}
$$
for $c,c' \in \Lambda$.
So $h_J(J,K) = 0$ for $J \ne K \ominus J^c$.
$J = K \ominus J^c$ does not hold for
$K \ne \{1,\dots,n\}$.
The proof of Lemma \ref{2418} is complete.
\end{proof}
\begin{cor}
$
\text{\bf e}''_K \in \bigoplus_{k<n} H^k(L(u);\Lambda).
$
\end{cor}
\begin{proof}
We put
$$
\text{\bf e}''_K - c_1 \text{\rm vol}_{L(u)} \in \bigoplus_{k<n}H^k(L(u);\Lambda),
\qquad
\text{\bf e}''_{\{1,\dots, n\}} - c_2 \text{\rm vol}_{L(u)} \in \bigoplus_{k<n}H^k(L(u);\Lambda).
$$
Here $c_2 \equiv 1 \mod \Lambda_+$ so $c_2\ne 0$.
By  Corollary \ref{degreemustn} the left hand side of
\eqref{lemmaindd} is
$$
c_1c_2\sum_{I,J \in 2^{\{1,\dots,n\}}} (-1)^{*}g^{IJ}
\langle \frak m_2^{\frak c,\frak b,b^{\frak c}}(\text{\bf e}''_I,
\text{\rm vol}_{L(u)}),
\frak m_2^{\frak c,\frak b,b^{\frak c}}(\text{\bf e}''_J,\text{\rm vol}_{L(u)})
\rangle_{\text{\rm PD}_{L(u)}}.
$$
By Theorem \ref{21mainformula} this is nonzero unless $c_1 = 0$.
Therefore the corollary follows from Lemma \ref{2418}.
\end{proof}
Now we have
$$\aligned
\langle \text{\bf e}''_I,\text{\bf e}''_J \rangle_{PD_{L(u)}}
= \langle \text{\bf e}''_I,\text{\bf e}''_J\cup^{\frak c, Q} e \rangle_{PD_{L(u)}}
&=  \langle \text{\bf e}''_I\cup^{\frak c, Q} \text{\bf e}''_J,e \rangle_{PD_{L(u)}} \\
&= c \langle \text{\bf e}''_{I\ominus J},e \rangle_{PD_{L(u)}} = 0.
\endaligned$$
Proposition \ref{2417} follows.
\end{proof}
Propositions  \ref{clifford3} and \ref{2417}
imply that an appropriate scalar multiplication of
$\text{\bf e}''_I$ gives cyclic Clifford basis.
\end{proof}

It is likely that the components of a cyclic cohomology in $$
F^{2k+1}HC^*(C)/F^{2k+2}HC^*(C)
$$ for $k>1$
can also be realized as a higher loop analogue of the trace $Z$
attached to the unital and cyclic filtered $A_{\infty}$ algebra.
Namely we expect that a $k$-loop analogue of $Z$ for a unital and cyclic filtered $A_{\infty}$ algebra
can be constructed by taking a sum of appropriate products of the structure constants of
$\frak m_{k'}$ ($k' \le k$) over certain Feynman diagrams into each of whose exterior vertices
the unit is inserted. 
Thus the current situation looks very similar to those appearing in the contexts of
\begin{enumerate}
\item higher residue pairing and primitive form. \cite{Sai83}.\index{Saito theory}
\item $S^1$ equivariant cohomology and gravitational
descendent. \cite{givental2}
\item asymptotic expansion and perturbative Chern-Simons
gauge theory.
\end{enumerate}
We hope to explore the relationship of our story with these
in a sequel of this paper.
\begin{rem}
We however note that our invariant $Z$ and its higher-loop cousins are slightly different
from the invariant of cyclic $A_{\infty}$ algebra described in \cite{konts:feyman}. In fact,
in our work the unit plays a significant role as well as the $A_{\infty}$ operations. 
The $0$-loop (or genus $0$) case of the invariant of \cite{konts:feyman} looks closer to
the invariant defined in \cite{fu091}, where
the Maslov class is assumed to vanish.
\end{rem}
\par

\section{Orientation}
\label{sec:ori}

The aim of this section is to compare the orientations of the moduli spaces
${\mathcal M}_{1,1}(\beta_{a}) ~_{\text{ev}^+} \times_{\text{ev}^+}
{\mathcal M}_{1,1}(\beta_{b})$
and
${\mathcal M}_{3,0} (\beta_{a'}) ~_{(\text{ev}_0,\text{ev}_1)} \times_{(\text{ev}_0,\text{ev}_1)}
{\mathcal M}_{3,0} (\beta_{b'})$ , which
appear as codimension two strata of the moduli space
of ${\mathcal M}^{\text{ann;main}}_{(1,1);0}(\beta)$.
\par
We recall from Lemma \ref{stratifyAn11} that
we have two points
$[\Sigma_1]$ and $[\Sigma_2]$
in the moduli space of annuli
${\mathcal M}_{(1,1);0}^{\text{\rm ann; main}}$.
Let $\gamma$ be a path in ${\mathcal M}_{(1,1),0}^{\text{ann;main}}$ joining $[\Sigma_1]$ and
$[\Sigma_2]$.
Taking the inverse image of $\gamma$ under the forgetful map
${\mathcal M}_{(1,1),0}^{\text{ann;main}}(\beta) \to {\mathcal M}_{(1,1),0}^{\text{ann;main}}$
in
(\ref{forgetfromann}),
we obtain a compact space with oriented Kuranishi structure, which bounds the inverse images of
$[\Sigma_1]$ and $[\Sigma_2]$.
Note that the inverse images of
$[\Sigma_1]$ and $[\Sigma_2]$ correspond to
${\mathcal M}_{1,1}(\beta_{a}) ~_{\text{ev}^+} \times_{\text{ev}^+}
{\mathcal M}_{1,1}(\beta_{b})$
and
${\mathcal M}_{3,0} (\beta_{a'}) ~_{(\text{ev}_0,\text{ev}_1)} \times_{(\text{ev}_0,\text{ev}_1)}
{\mathcal M}_{3,0} (\beta_{b'})$, respectively.
Then our main result of this section is as follows.

\begin{prop}\label{compori}
Let $\gamma$ be a path in ${\mathcal M}_{(1,1);0}^{\text{\rm ann; main}}$ joining
$[\Sigma_1]$ and $[\Sigma_2]$.
Then $\gamma$ intertwines the orientations between
$$
(-1)^{\frac{n(n-1)}{2}}\bigcup_{\beta_a + \beta_b = \beta}
{\mathcal M}_{1,1}(\beta_a) ~_{\text{ev}^+} \times_{\text{ev}^+} {\mathcal M}_{1,1}(\beta_b)
$$
and
$$
\bigcup_{\beta_{a'} + \beta_{b'}=\beta}
{\mathcal M}_{3,0}(\beta_{a'}) ~_{(\text{ev}_0,\text{ev}_1)} \times_{(\text{ev}_0,\text{ev}_1)} {\mathcal M}_{3,0}(\beta_{b'}).
$$
\end{prop}
The orientation of the moduli spaces appearing in Proposition \ref{compori}
will be defined in the subsections below.
To compare the orientations of these moduli spaces,
we follow the argument in \cite{fooobook2} and reduce the problem to the one for linearized operators.
As in the proof of \cite[Proposition 8.1.4]{fooobook2}, we further reduce the problem to the case of product bundle pairs.
We discuss this point in detail below.

\par
\subsection{Operators $D_1$ and
$D_2$}
\label{subsec:opeD_i}
Firstly, we explain the orientation of   ${\mathcal M}^{\text{ann;main}}_{(1,1);0}(\beta)$.
\par
Since $L$ is  a principal homogeneous space under the action of $T^n$,
its tangent bundle $TL$ is trivialized by the fundamental vector fields of the $T^n$-action.
We fix this trivialization in this section.
We reduce the problem of determining orientation to the one in the linearized problem.
\par
Let $\rho > 1$, $ A=A_{\rho}= [1,\rho] \times S^1$ and
$u: A\to X$ a smooth map with
$u(\partial_i A) \subset L$ for $i=1,2$.
Here we denote  $\partial_1 A=\{\rho \} \times S^1$ and $\partial_2 A=\{ 1 \} \times S^1$.
Pick $p>2$ and consider the linearization of the Cauchy-Riemann equation at $u$:
$$D_u \overline{\partial}:W^{1,p}((A; \partial_1 A, \partial_2 A), u^*(TX;TL,TL)) \to
L^p(A; \Lambda^{0,1}A \otimes u^*TX).$$
\par
We adopt the argument in
\cite[Subsection 8.1.1]{fooobook2} as follows. 
We pick a boundary parallel circles $ C_1=\{\rho -\epsilon\} \times S^1 \subset A,
C_2=\{1+ \epsilon\} \times S^1$
with a sufficiently small  $\epsilon > 0$.
Denote by $Z \cong D^2_{(1)} \cup \C P^1 \cup D^2_{(2)}$
the quotient space of $A$ by collapsing $C_1$ and $C_2$ to points.
Here we identify the origin $O_{(1)}$ of the disk $D^2_{(1)}$ with the ``south pole'' $S \in S^2 \cong \C P^1$
and the origin $O_{(2)}$ of the disk $D^2_{(2)}$ with the ``north pole'' $N \in S^2 \cong \C P^1$.
We identify $\partial D^2_{(i)}=\partial_i A$.
Using the trivialization of $TL$, we can descend the vector bundle
$u^*TX \to A$  to $E \to Z$.
Note that $(u\vert_{\partial_i A} )^*TL \subset E\vert_{\partial_i A}$, $i=1,2$, as totally real subbundles.
Then the orientation of the index $D_u \overline{\partial}$ is described as the fiber product
orientation of the index problems on $D^2_{(1)}, \C P^1$ and
$D^2_{(2)}$ as follows.
We identify $E\vert_{D^2_{(i)}}$ with a  trivial bundle $D^2_{(i)} \times \C^n$ such that
the totally real subbundle $(u\vert_{\partial_i A})^* TL$ is identified with $\partial_i A \times
\R^n$.
\par
By deforming the operators in the space of Fredholm operators keeping the boundary condition,
we can assume our
the operators on $D^2_{(1)}, \C P^1$ and $D^2_{(2)}$ are Dolbeault operators:
$$
\aligned
\overline{\partial}_{(i)}=\overline{\partial} ~&:~ W_i
= W^{1,p}(D^2_{(i)}, \partial D^2_{(i)};\underline{\C}^n, \underline{\R}^n) \to L^p (D^2_{(i)}, \Lambda^{0,1} D^2_{(i)} \otimes \underline{\C}^n),
~i=1,2 \\
\overline{\partial}_{\C P^1}=\overline{\partial}~&:~
W_{\C P^1}
=
W^{1,p}(\C P^1, E\vert_{\C P^1}) \to L^p(\C P^1,\Lambda^{0,1} \C P^1 \otimes E).
\endaligned
$$
Here and hereafter we denote by $\underline{\C}^n$ or $\underline{\R}^n$, the trivial
$n$ dimensional,
complex or real vector bundles, respectively.
Denote by ${\rm\bf Dom}$ the triple $(\xi_1,\zeta,\xi_2) \in W_1 \times W_{\C P^1} \times W_2$
such that
$$\xi_1(O_{(1)})=\zeta (S) \in E_{O_{(1)}=S}, \quad
\xi_2(O_{(2)})=\zeta (N) \in E_{O_{(2)}=N}.
$$
Then we find that the index of $D_u \overline{\partial}$ is identified with
the index of the restriction of $\overline{\partial}_{(1)} \oplus  \overline{\partial}_{\C P^1} \oplus
\overline{\partial}_{(2)}$ to ${\rm\bf Dom}$.
\par
We now study the orientation of
${\mathcal M}_{1,1}(\beta_{a}) ~_{\text{ev}^+} \times_{\text{ev}^+} {\mathcal M}_{1,1}(\beta_{b})$ as a fiber product.
We first
give orientations on ${\mathcal M}_{1,1} (\beta_a)$ and
${\mathcal M}_{1,1} (\beta_b)$ by following the argument in
\cite[Subsection 8.1.1]{fooobook2} again. That is,
we collapse boundary parallel circles $C_i$ of disks in each factor and reduce the
orientation problem to the one for the Dolbeault operator on $\C P^1_{(i)}$ with holomorphic
vector bundle $E_{(i)}$ and the one on $(D^2, \partial D^2)$ with
$(\underline{\C}^n, \underline{\R}^n)$.
By homotopy, we may assume that the origin $O_{(i)} \in D^2_{(i)}$ is located between
$C_i$ and $\partial D^2_{(i)}$.
Since the Dolbeault operators on $\C P^1_{(i)}$ are complex Fredholm operators and
the fiber products with their indices are taken over complex vector spaces, the orientation
problem is reduced to the case that $E_{(i)} \cong D^2_{(i)} \times \C ^n$ and the totally
real subbundle over $\partial D^2_{(i)}$ are $\partial D^2_{(i)} \times \R ^n$.
Namely,
it is enough to study
the orientation of the index of the operator $D_1$, which is the restriction of
$\overline{\partial}_{(1)} \oplus \overline{\partial}_{(2)}$ to
$${\rm\bf Dom}_1=\{(\xi_1,\xi_2) \in W_1 \oplus W_2
~\vert~ \xi_1(O_{(1)}) = \xi_2(O_{(2)}) \}.
$$
\par
Next, we explain the orientation of
${\mathcal M}_{3,0} (\beta_a) ~_{(\text{ev}_0,\text{ev}_1)} \times_{(\text{ev}_0,\text{ev}_1)} {\mathcal M}_{3,0} (\beta_b)$.
Denote by $z_0^{(i)}, z_1^{(i)}, z_2^{(i)}$ are boundary marked points on $\partial D^2_{(i)}$
respecting the counter clockwise orientation of the boundary circles.
As in the previous case, the orientation problem is reduced to the orientation of the index of
the operator $D_2$, which is the restriction of the Dolbeault operator $\overline{\partial}_{(1)}
\oplus \overline{\partial}_{(2)}$  to
$$
{\rm\bf Dom}_2=\{(\xi_1,\xi_2) \in W_1 \oplus W_2 ~\vert~ \xi_1(z_i^{(1)})=\xi_2(z_i^{(2)}) \text{ for } i=1,2 \}.
$$

\par
To prove Proposition \ref{compori} it suffices to
compare the orientations of the indices of these operators $D_1$ and $D_2$.
\begin{notation}\label{basis}
\begin{enumerate}
\item In the rest of this section, $\langle v_1, \dots , v_n \rangle$ stands for the $n$ dimensional real vector space
spanned by $v_i$ equipped with the orientation
given by the ordering $v_1 ,\dots , v_n$.
\item We fix an orientation on ${\R}^n$. (In geometric setting, we fix an orientation of the Lagrangian
submanifold $L$.)  Pick an oriented basis
${\bf e}_1, \dots, {\bf e}_n$ of ${\R}^n$ so that
${\R}^n=
\langle {\bf e}_1, \dots, {\bf e}_n \rangle$.
Then
$\langle {\bf e}_1, \sqrt{-1}{\bf e}_1, \dots, {\bf e}_n, \sqrt{-1}{\bf e}_n \rangle$
gives the orientation on ${\C}^n$, which will be
identified with the tangent space of the ambient
toric manifold $X$.
\end{enumerate}
\end{notation}

\subsection{Orientation of $\text{Index } D_1$}
\label{subsec:oriD_1}
For the bundle pair $(\underline{\C}^n,\underline{\R}^n)$,
the Dolbeault operators $\overline{\partial}_{(i)}$
$(i=1,2)$ are surjective. 
To study orientation of ${\rm Index} D_1$, we need to pick finite dimensional reductions of $\overline{\partial}_{(i)}$
such that the coincidence condition for values at $O_{(1)}$ and $O_{(2)}$ is transversal.
\par
Let $A'_{i}=\{z \in \C ~\vert~ 1/4 \leq \vert z \vert \leq 3/4\} \subset D^2_{(i)}$.
Pick and fix a decreasing smooth function $\sigma:[0,1] \to \R$ such that
$\sigma(r) \equiv 1$ for $r \leq 1/4$ and $\sigma(r) \equiv 0$ for $r \geq 3/4$.
We put
$$
f_j^{(i)}=\sqrt{-1} \sigma(\vert z \vert){\bf e}_j, \qquad g_j^{(i)}={\bf e}_j, \qquad i=1,2,  ~j=1,\dots ,n
$$
and define
$2n$-dimensional oriented real vector spaces $U_i$
by
$$
U_i=\langle
f_1^{(i)},\dots ,f_n^{(i)}, g_1^{(i)},\dots ,g_n^{(i)}
\rangle.
$$
We define $E_i$ by
$$E_i = \langle \nu_1^{(i)},\dots ,
\nu_n^{(i)} \rangle
$$
where
$$
\nu_j^{(i)}=\sqrt{-1} \sigma'(\vert z \vert) (dr - \sqrt{-1}r d\theta) \otimes {\bf e}_j, \quad j=1, \dots ,n.
$$
Here $\sigma'$ is the first derivative of $\sigma$ and $(r, \theta)$ is the polar coordinate so that $z=r e^{\sqrt{-1} \theta}$.
The restriction of the Dolbeault operator $\overline{\partial}_{(i)}$ is the following mapping
$s_i:U_i \to \underline{E}_i = U_i \times E_i$ where
\begin{equation}
s_i(f_j^{(i)})=(f_j^{(i)}, \nu_j^{(i)}), \quad s_i(g_j^{(i)})=(g_j^{(i)}, 0). \nonumber
\end{equation}
\par
We denote by $X_i=(s_i;\underline{E}_i \to U_i)$ the Kuranishi
neighborhood defined by $s_i : U_i \to \underline{E}_i$. 
(See  \cite[Definition A1.1]{fooobook2} for the definition of the Kuranishi neighborhood.) 
The orientation on $X_i$ is defined by the orientation on 
$\det \underline{E}_i \otimes \det TU_i$. 
(In general, in case when $s: U \to E$ is transverse to zero, the orientation on $s^{-1}(0)$ is given by 
$$
E_{\ast} \times T_{\ast}s^{-1}(0)=
T_{\ast}U.
$$ See
\cite[Convention 8.2.1 (1)]{fooobook2}.)
The Kuranishi neighborhood of the fiber product $X_1 \times_{{\R}^n}X_2$ of the Kuranishi neighborhoods $X_i$ with respect to the evaluation maps $\text{ev}^+_{(i)}$ at the points
$O_{(i)} \in D_{(i)}^2$
is given by
$$
s=s_1 \oplus s_2: U_1 ~_{\text{ev}^+_{(1)}} \times_{\text{ev}^+_{(2)}} U_2 \to  \underline{E}_1 \oplus \underline{E}_2.
$$ 
Then by \cite[Convention 8.2.1 (4)]{fooobook2}, 
the orientation on the Kuranishi neighborhood of the fiber product
is twisted from the orientation defined by 
$s=s_1 \oplus s_2: U_1 ~_{\text{ev}^+_{(1)}} \times_{\text{ev}^+_{(2)}} U_2 \to  \underline{E}_1 \oplus \underline{E}_2$
as follows:
\begin{equation}\label{-n}
X_1 \times_{{\R}^n}X_2= (-1)^{n} (s=s_1 \oplus s_2; \underline{E}_1 \oplus \underline{E}_2 \to U_1 ~_{\text{ev}^+_{(1)}} \times_{\text{ev}^+_{(2)}} U_2).
\end{equation}

Since the evaluation maps are given by
$\text{ev}^+_{(i)}(f_j^{(i)})=\sqrt{-1} {\bf e}_j$ and $\text{ev}^+_{(i)}(g_j^{(i)})={\bf e}_j$ in this model, we find that
$$U_1 ~_{\text{ev}^+_{(1)}} \times_{\text{ev}^+_{(2)}} U_2=
\langle
{\bf f}_1, \dots , {\bf f}_n, ~{\bf g}_1, \dots ,
{\bf g}_n \rangle
$$
as an oriented vector space. Here
$$
{\bf f}_j=
(f_j^{(1)},f_j^{(2)}),
\quad {\bf g}_j=(g_j^{(1)},g_j^{(2)}),
\quad
j=1, \dots, n.
$$
We note that for each $i=1,2$
the orientation given by
$$
\aligned
& \langle
\text{ev}^+_{(i)}(f_1^{(i)}),\dots ,
\text{ev}^+_{(i)}(f_n^{(i)}),
~\text{ev}^+_{(i)}(g_1^{(i)}), \dots ,
\text{ev}^+_{(i)}(g_n^{(i)})
\rangle  \\
& =
\langle
\sqrt{-1}{\bf e}_1, \dots, \sqrt{-1}{\bf e}_n, {\bf e}_1, \dots, {\bf e}_n
\rangle
\endaligned$$ differs from
the complex orientation of ${\C}^n$ by $(-1)^{n(n+1)/2}$. (See Notation \ref{basis}.2.)
Here
we regard ${\C}^n$ as the tangent space of the ambient toric manifold $X$ which is the target space of the evaluation maps.
\par
Therefore, by combining this with (\ref{-n})
and
canceling $s ({\bf f}_j)=(\nu_j^{(1)},\nu_j^{(2)})$ out from the front of the basis,
we obtain a finite dimensional reduction of the operator $D_1$:
\begin{equation}\label{obstE2}
(-1)^{n(n-1)/2} (0; \underline{E}_2 \to U).
\end{equation}
Here we define $U$ by
$$
U=\langle
{\bf g}_1,\dots ,
{\bf g}_n
\rangle.
$$
We identify $E_2 \times U$
with $\underline{E}_2 = E_2 \times U_2$ in an
obvious way.
This is the obstruction bundle of the above Kuranishi structure. 
\subsection{Orientation of $\text{Index } D_2$}
\label{subsec:oriD_2}
For convenience of description, we use $\R \times [0,1] \subset \C$ instead of $D^2 \subset \C$.
Let $x,y$ be coordinates such that $z=x+\sqrt{-1}y$ on $\C$.
We identify $z_0$ (resp. $z_1$) with the limit $x \to +\infty$ (resp. $x \to -\infty$).
Pick and fix an increasing smooth function $\rho:\R \to \R$ such that $\rho(x) \equiv 0$ for $x \leq -1$
and $\rho(x) \equiv 1$ for $x \geq 1$.
We put
$$
h_j^{(i)}=\rho(x) {\bf e}_j,
\qquad
k_j^{(i)}={\bf e}_j, \qquad  i=1,2,~j= 1, \dots, n
$$
and define the oriented vector spaces
$$
\aligned
V_i & =\langle
h_1^{(i)}, \dots , h_n^{(i)}, ~
k_1^{(i)}, \dots , k_n^{(i)}
\rangle \\
F_i &
=\langle
\chi_1^{(i)}, \dots , \chi_n^{(i)}\rangle,
\endaligned
$$
where
$$
\chi_j^{(i)}=\rho'(x) (dx - \sqrt{-1} dy), \quad j= 1, \dots, n.
$$
Then the Dolbeault operator
$\overline{\partial}_{(i)}$ restricts to a mapping
$t_i: V_i \to \underline{F}_i= V_i \times F_i$ which is defined by
$$
t_i(h_j^{(i)})=(h_j^{(i)}, \chi_j^{(i)}), \quad
t_i(k_j^{(i)})=(k_j^{(i)}, 0).
$$
Namely
$t_i:V_i \to \underline{F}_i$ is a finite dimensional reduction of $\overline{\partial}_{(i)}$.
\par
By \cite[Convention 8.2.1 (4)]{fooobook2} again, the fiber product of the linearized Kuranishi models
$(t_i;\underline{F}_i \to V_i)$ with respect to the evaluation maps $(\text{ev}_0,\text{ev}_1)$ at
$z_0$ and $z_1$ is
given by
\begin{equation}\label{-n2}
(-1)^n (t=t_1 \oplus t_2; \underline{F}_1 \oplus \underline{F}_2 \to
V_1 ~_{(\text{ev}_0,\text{ev}_1)} \times_{(\text{ev}_0,\text{ev}_1)} V_2).
\end{equation}
Since the evaluation maps are given by $\text{ev}_0(h_j^{(i)})=\text{ev}_0(k_j^{(i)})={\bf e}_j$, $\text{ev}_1(h_j^{(i)})=0$ and $\text{ev}_1(k_j^{(i)})={\bf e}_j$,
we find that
$$
V_1 ~_{(\text{ev}_0,\text{ev}_1)} \times_{(\text{ev}_0,\text{ev}_1)} V_2 =
\langle
{\bf h}_1, \dots , {\bf h}_n, ~
{\bf k}_1, \dots , {\bf k}_n
\rangle
$$
where
$$
{\bf h}_j=(h_j^{(1)},h_j^{(2)}), \quad {\bf k}_j=(k_j^{(1)}, k_j^{(2)}), \quad j= 1, \dots, n.
$$
By canceling $t({\bf h}_j)=(\chi_j^{(0)},\chi_j^{(1)})$ out from the front of the basis,
we obtain a finite dimensional reduction of the operator $D_2$:
\begin{equation}\label{obstF2}
(-1)^n (0;\underline{F}_2 \to V:=\langle
{\bf k}_1, \dots ,
{\bf k}_n \rangle).
\end{equation}
Here the obstruction bundle is $\underline{F}_2$, which we identify
with $F_2 \times V$.

\subsection{Continuation of linear Kuranishi models}
\label{subsec:continuation}

In order to compare the orientations of the index bundles of $D_1$ and $D_2$,
we will find a continuous family of linear Kuranishi models joining them.

Firstly, we use of gluing to reduce the problems to the ones on the annulus.
For $\rho >1$ we denote the annulus by
$A_{\rho}=\{ z \in \C ~\vert~ 1 \leq \vert z \vert \leq \rho \}$.
We use the polar coordinates $(r, \theta)$.

\subsubsection{Linear Kuranishi model for $\widehat D_1$}
\label{subsec:kuraD_1}
Gluing the two disks $D^2_{(1)}$ and $D^2_{(2)}$ by $O_{(1)}=O_{(2)}$ and
smoothening the double point, we obtain the annulus $A=A_{\rho_1}$.
$(\rho_1 > 1).$
Here we identify
$$\partial D^2_{(1)}=\{z \in \C ~\vert~ \vert z \vert = \rho_1 \}, \quad
\partial D^2_{(2)}=\{z \in \C ~\vert~ \vert z \vert = 1\}.
$$
(This makes the function $\varphi$ introduced below  increasing.)
Let $\widehat D_1$ be the  operator obtained from $D_1$ after this gluing.
Below we describe a finite dimensional model of $\widehat D_1$.
\par
Note that the elements of the obstruction bundle $\underline{E}_2$ are supported
in the annular region in $D^2_{(2)}$ away from $O_{(2)}$.
Therefore after we glue two disks the support of $\underline{E}_2$
becomes to a subannulus  $A'=\{z \in \C ~\vert~ r_1 \leq \vert z \vert \leq r_2\}$ in $A$.
The elements of $\underline{E}_2$ are also identified with ${\C}^n$-valued $(0,1)$-forms with support
on $A'$.
In other words,  there exist $1 < r_1 < r_2 < \rho_1$ and a increasing smooth function
$\varphi (r)$ with the property that
$\varphi (r) \equiv 0$ for $1 \leq r \leq r_1$ and $\varphi (r) \equiv 1$ for $r_2 \leq r \leq \rho_1$
such that after we glue two disks, the vector space $E_2$ becomes:
$$
E=\left\langle \varphi'(r) \left(\frac{\sqrt{-1}}{r} dr +d\theta \right) \otimes {\bf e}_1, \dots ,
\varphi'(r) \left(\frac{\sqrt{-1}}{r} dr +d\theta \right) \otimes {\bf e}_n
\right\rangle.
$$
Denote by $g_j: A=A_{\rho_1} \to {\C}^n$ the constant map with the value ${\bf e}_j$ and
define an $n$ dimensional oriented real vector space $$\overline V
:= \langle g_1, \dots ,g_n\rangle.$$
Then the linear Kuranishi model $(-1)^{n(n-1)/2}(0;\underline{E}_2 \to U)$ in (\ref{obstE2})
becomes the following
\begin{equation}\label{linearE}
(-1)^{n(n-1)/2} (0;\underline{E} \to \overline V),
\end{equation}
after gluing.
Here $\underline{E} = E\times \overline V$.
This is a linear Kuranishi model for $\widehat D_1$.

\subsubsection{Linear Kuranishi model for $\widehat D_2$}
\label{subsec:kuraD_2}
Recall that we identified $D^2$ with $\R \times [0,1]$
with $z_0$ (resp. $z_1$).
Here $z_0$ (resp. $z_1$) corresponds to the limit $x \to + \infty$ (resp. $x \to -\infty$).
We pick $z_2 =0 \in \R \times [0,1] \subset \C$.
For a sufficiently large $R >0$, we glue two copies of
$([-R, R] \times [0,1])_{(i)}$ by identifying
$(\{R\} \times [0,1])_{(1)}=(\{R\} \times [0,1])_{(2)}$ and
$(\{-R\} \times [0,1])_{(1)}=(\{-R\} \times [0,1])_{(2)}$ to obtain an annulus  $A$.
\par
We pick an identification of $A$ with $A_{\rho_2}$ ($\rho_2 > 1$)
 such that $z_2^{(1)} $ corresponds to a point on  $\partial_1 A_{\rho_2}$ and $z_2^{(2)}$
corresponds to a point on the unit circle $\partial_2 A_{\rho_2}$.
\par
Let $\widehat D_2$ be the operator obtained from $D_2$ after this gluing.
Below we describe a finite dimensional model of
$\widehat D_2$.
\par
Recall that elements of the obstruction bundle
$\underline{F}_2$ are ${\C}^n$-valued $(0,1)$-form
supported in the region $\{-1 \leq x \leq 1 \} \subset \R \times [0,1]$.
Hence, after the gluing, these elements correspond to
${\C}^n$-valued $(0,1)$-forms on the annulus $A_{\rho_2}$ with support in $\{ c_1 \leq \theta \leq c_2 \} \subset A$,
where $\pi/2 < c_1 < c_2 < 3\pi/2$.
Moreover, there exists a decreasing smooth function
$\psi:[0,2\pi ) \to \R$
with the property that
$\psi (\theta) \equiv 1$ for $\theta \leq c_1$ and $\psi (\theta) \equiv 0$ for $\theta \geq c_2$
such that $F_2$ becomes the space
$$
F=\left\langle
\psi'(\theta ) \left(\frac{\sqrt{-1}}{r} dr
+ d\theta \right) \otimes {\bf e}_1, \dots ,
\psi'(\theta ) \left(\frac{\sqrt{-1}}{r} dr
+ d\theta \right) \otimes {\bf e}_n
\right\rangle
$$
after gluing. (The function $\psi$ is decreasing because of the way we identify
$A$ with $A_{\rho_2}$.)
\par
Then the linear Kuranishi model $(-1)^n(0;\underline{F}_2 \to V)$
described in Subsection \ref{subsec:oriD_2}
becomes
\begin{equation}\label{linearF}
(-1)^n (0;\underline{F} \to\overline V
=\langle g_1 ,\dots , g_n \rangle),
\end{equation}
after gluing.
Here $g_j:A=A_{\rho_2} \to {\C}^n$ is the constant map with the value ${\bf e}_j$
and $\underline{F} = F \times \overline V$.
This is our linear Kuranishi model for $\widehat D_2$.

\subsubsection{Continuation of the linear Kuranishi models}
\label{subsec:cont}

Let $A = A_{\rho_1} = A_{\rho_2}$.
We take functions $\varphi (r)$ (resp. $\psi(\theta)$) as in Subsection
\ref{subsec:kuraD_1} (resp. Subsection \ref{subsec:kuraD_2}).
\par
We consider the following families of oriented vector spaces: 
$$
\aligned
E(\tau) = &
\bigg\langle
(\tau + (1-\tau) \varphi'(r))
\Big(\frac{\sqrt{-1}}{r} dr +d\theta \Big)
\otimes {\bf e}_1, \\
& \qquad \dots ,
(\tau + (1-\tau) \varphi'(r))
\Big(\frac{\sqrt{-1}}{r} dr +d\theta \Big)
\otimes {\bf e}_n
\bigg\rangle \\
F(\tau)= &
\bigg\langle(-\tau+(1-\tau)\psi'(\theta))
\Big(\frac{\sqrt{-1}}{r}dr
+d\theta
\Big) \otimes {\bf e}_1, \\
& \qquad \dots ,
(-\tau+(1-\tau)\psi'(\theta))
\Big(\frac{\sqrt{-1}}{r}dr
+d\theta
\Big) \otimes {\bf e}_n
\bigg\rangle,
\endaligned
$$
which are subspaces of $C^{\infty}(A;\Lambda^{0,1}\otimes \underline{\C}^n)$.
We note
\begin{equation}\label{vectspacebdry}
F(0) = F, \quad E(0) = E, \quad  F(1) = (-1)^n E(1).
\end{equation}
\par
Recall that $\varphi (r)$ (resp. $\psi(\theta)$) is an increasing (resp. decreasing) function.
Therefore for any $\tau \in [0,1]$ we have
\begin{equation}\label{posnegineq}
\aligned
\tau + (1-\tau) \varphi'(r) & \ge 0 \\
-\tau + (1-\tau) \psi'(r) & \le 0,
\endaligned
\end{equation}
both of which are not identically zero.
\par
We show the following:

\begin{lem}\label{continuation}
Let $P=\overline{\partial}:W^{1,p}(A,\partial A;{\C}^n, {\R}^n) \to L^p(A;\Lambda^{0,1}A \otimes {\C}^n)$
be the Dolbeault operator on the annuls $A$ with coefficients in the trivial bundle pair
$(\underline{{\C}}^n, \underline{{\R}}^n)$.
Here $p>2$.
For any $\tau \in [0,1]$, we have the following:
\begin{enumerate}
\item $\text{\rm Im } \overline{\partial} \oplus E(\tau) = L^p(A;\Lambda^{0,1}A \otimes {\C}^n).$
\item $\text{\rm Im } \overline{\partial} \oplus F(\tau) = L^p(A;\Lambda^{0,1}A \otimes {\C}^n).$
\end{enumerate}
\end{lem}

\begin{proof}
For the Dolbeault operator $P$ on the annulus
$(A, \partial A)$ with coefficients in the trivial
bundle pair $(\underline{\C}^n,\underline{\R}^n)$, we have
$\text{index }P = 0$.
We also find that $\text{Ker }P\cong{\R}^n$ consisting of constant sections with values in ${\R}^n$.
Thus $\text{Coker }P$ is $n$-dimensional.
In order to prove the lemma, it is enough to show
$\text{\rm Im } \overline{\partial} \cap E(\tau) = 0$ and
$\text{\rm Im } \overline{\partial} \cap F(\tau) =0$, respectively.
Since the bundle pair
$(\underline{{\C}}^n,\underline{{\R}}^n)$ is the direct product of $n$-copies of $(\underline{\C},\underline{\R})$,
we consider the case that $n=1$ from now on.
\par\smallskip
\noindent
Proof of 1:  Suppose that there exists $f:(A,\partial A) \to (\C,\R)$ such that
$$\overline{\partial} f = \left(\tau + (1-\tau) \varphi'(r)\right)
\left(\frac{\sqrt{-1}}{r}dr +d\theta\right).$$
Pick a function $g(r):[1,\rho_1 ] \to \R$ such that $g(1)=0$ and
$$r g'(r)=\tau + (1-\tau) \varphi'(r).$$
Then we have
$$\overline{\partial} (\sqrt{-1} g(r))=(\tau + (1-\tau) \varphi'(r)) \left(\frac{\sqrt{-1}}{r}dr+d\theta\right).
$$
Then $f-\sqrt{-1}g$ is a holomorphic function on $A$ such that
\begin{equation}
\text{Im }(f -\sqrt{-1}g) =
\begin{cases}
0 &\quad \quad \text{on } \partial_1A \\
g(\rho_1 ) &\quad \quad \text{on } \partial_2 A.
\end{cases}
\nonumber
\end{equation}
Here
$\partial_1A = \{ z \in \C \mid \vert z\vert = \rho_1\}$,
$\partial_2A = \{ z \in \C \mid \vert z\vert = 1\}$.
\par
By the Schwarz reflection principle, we extend $f-\sqrt{-1}g$ to a holomorphic function on
$\C \setminus \{0\}$, the real part of which is bounded.
Since the region $\{z \in \C ~\vert~ \text{Re }z < R \}$ for a constant $R>0$ is conformally equivalent
to the unit disk, the removable singularity theorem implies that $f-\sqrt{-1}g$ extends across
the origin in $\C$ and it is a constant function.
On the other hand, by (\ref{posnegineq}),
$\tau + (1-\tau) \varphi'(r)$ is non-negative and positive
somewhere for $\tau \in [0,1]$.
Thus $g(\rho_1 ) > 0$.  This is a contradiction.   Hence we have
$\text{\rm Im } \overline{\partial} \cap E(\tau) = 0$.
\par\smallskip
\noindent
Proof of 2: Let ${\rm pr}:\widetilde{A}=[1,\rho_2 ] \times \R \to A$ be the universal cover of $A$.
By abuse of notation, we denote by $(r, \theta)$ the coordinates on $[1,\rho_2 ] \times \R$.
We define $\widetilde{f}(r,\theta)=\widetilde{f}(\theta)$ by
$$\widetilde{f}(\theta)=\int_0^{\theta}(-\tau + (1-\tau) \psi'(\theta)) d\theta,$$
which satisfies
$$\overline{\partial} \widetilde{f}=(-\tau+(1-\tau)\psi'(\theta))\left(\frac{\sqrt{-1}}{r}dr+d\theta\right).$$
Suppose that there exists $h:(A,\partial A) \to (\C,\R)$ such that
$$\overline{\partial}h =(-\tau + (1-\tau) \psi'(\theta))\left(\frac{\sqrt{-1}}{r}dr+d\theta\right).$$
Then $\widetilde{k}=\widetilde{f}-h \circ {\rm pr}:(\widetilde{A}, \partial \widetilde{A}) \to (\C,\R)$
is a holomorphic function such that
$$\widetilde{k}(r,2\pi)-\widetilde{k}(r,0)=\widetilde{f}(2\pi)-\widetilde{f}(0).$$
By the Schwarz reflection principle, we extend $\widetilde{k}$ to an entire holomorphic function
on $\C$, the imaginary part of which is bounded.
Hence $\widetilde{k}$ is a constant function.

On the other hand, by (\ref{posnegineq}),
$-\tau+(1-\tau)\psi'(\theta)$ is non-positive and negative somewhere
for $\tau \in [0,1]$.
Hence $\widetilde{f}(2\pi)-\widetilde{f}(0) \neq 0$.
This is a contradiction and we have $\text{\rm Im } \overline{\partial} \cap F(\tau) =0$.
\end{proof}

Lemma \ref{continuation} and
(\ref{vectspacebdry}) imply that the linear Kuranishi models
$(0;\underline{E} \to \overline V)$ and $(-1)^n (0;\underline{F} \to \overline V)$ are joined by a one parameter
family of linear Kuranishi models $\{(0;\underline{E}(\tau) \to \overline V)\}_{0 \leq \tau \leq 1}$
and $\{(0;\underline{F}(\tau) \to \overline V)\}_{0 \leq \tau \leq 1}$ and the deformation of
the complex structure on the annuli as $\rho$ varies from $\rho_1$ to $\rho_2$.
Therefore by (\ref{linearE}) and (\ref{linearF})
the orientations of the linear Kuranishi models of the indices of $D_1$ and $D_2$ differ
by the multiplication by $(-1)^{n(n-1)/2}$.
This implies Proposition \ref{compori}.
\qed

\section{Sign in Theorem \ref{annulusmain} and
Proposition \ref{mainprop}}
\label{sec:sign}

Using Proposition \ref{compori}, we check the signs in the formulas in Theorem \ref{annulusmain} and
Proposition \ref{mainprop}.
Throughout this section,
space means one with oriented Kuranishi structure.
(See \cite[Definition A1.17]{fooobook2}).
The dimension of a space with Kuranishi structure
means the virtual dimension.
We use the same notation as in
\cite{fooobook2}
and follow the convention.
Especially see \cite[Section 8.2]{fooobook2}.
Other than the convention of
the pairing we use the same convention as \cite{fooobook2} on orientation and sign. See Remark \ref{rempair}.

\subsection{Some lemmata}
\label{subsec:signpreliminary}
In this subsection we prepare some general lemmata on orientation
of fiber product and Poinca\'e duality which will be used later.
Let $Z$ be a smooth closed oriented manifold.
For the later argument we consider the case either $Z$ is
a toric manifold $X$ of dimension $2n$, or a Lagrangian submanifold $L$ of $X$.
Let
$S$ and $T$ be spaces with oriented Kuranishi structures.
Let $f : S \to Z$ and $g : T \to Z$ be weakly submersive maps.
(See  \cite[Definition A1.13]{fooobook2} for the definition of weakly submersive maps.) Then we can define the fiber product
$S\times _Z T :=S {}_{f}\times_{g}T$ over $Z$
in the sense of Kuranishi structure.
We put
$$
\deg S = \dim Z - \dim S, \quad \deg T = \dim Z -\dim T.
$$
By \cite[Convention 8.2.1 (4)]{fooobook2}, we define the orientation on
$S\times _Z T$.
On the other hand,
we denote by $\Delta$ the diagonal set in $Z \times Z$.
Let $i:\Delta \to Z \times Z$ be an obvious inclusion map.
By \cite[Remark A1.44]{fooobook2},
we can also define the fiber product
$\Delta \times_{Z\times Z} (S\times T)$.
Clearly $S\times _Z T=\Delta \times_{Z\times Z} (S\times T)$ as a set, but as for the orientations we have

\begin{lem}\label{LemmaB}
$S\times _Z T=
(-1)^{\dim Z\deg S}
\Delta \times_{Z\times Z} (S\times T)$.
\end{lem}

\begin{proof}
First we assume that $S, T$ are smooth manifolds and $f_S, f_T$ are submersion.
Put $S=S^{\circ} \times Z$ and $T=Z\times {}^{\circ}T$.
See \cite[Convention 8.2.1]{fooobook2}  for this notation.
Then  \cite[Convention 8.2.1 (3)]{fooobook2} gives the orientation on $S\times _Z T$ by
$S\times _Z T = S^{\circ} \times Z \times {}^{\circ}T
=(-1)^{\dim Z\deg S}Z \times S^{\circ} \times
{}^{\circ}T$.
On the other hand, if we put $S \times T= (Z \times Z) \times {}^{\circ}(S \times T)$,
we find $S^{\circ} \times {}^{\circ}T = {}^{\circ}(S\times T)$.
By identifying $\Delta =Z$ as oriented spaces, we
obtain $S\times _Z T=
(-1)^{\dim Z\deg S}
\Delta \times_{Z\times Z} (S\times T)$.
For general cases we can prove the lemma by using
\cite[Convention 8.2.1 (4)]{fooobook2}.
\end{proof}

Next, let $f_i : S \to Z$ and $g_i : T \to Z$ ($i=1,2$)
be weakly submersive maps.
We can define the fiber product
$S\times_{Z\times Z}T
:=S{}_{(f_1, f_2)}\times_{(g_1, g_2)}T$
over $Z \times Z$.
Then we can show the following.

\begin{lem}\label{LemmaA}
$$
S\times_{Z\times Z}T = (-1)^{\dim Z(\deg S +\dim T)}
\big(
\Delta {}_{i}
\times _{f_1 \times g_1} (S \times T) \big) {}_{f_2 \times g_2}\times _{i}\Delta .
$$
\end{lem}
\begin{rem}
By the associative property of the fiber product (\cite[Lemma 8.2.3]{fooobook2}),
this is equivalent to
$$
S\times_{Z\times Z}T =
(-1)^{\dim Z(\deg S +\dim T)}
\Delta {}_{i}
\times _{f_1 \times g_1} \big( (S \times T)  {}_{f_2 \times g_2}\times _{i}\Delta  \big).
$$
If no confusion can occur, we simply denote the right hand side by
$$
(-1)^{\dim Z(\deg S +\dim T)}
\Delta
\times _{Z \times Z}  (S \times T) \times_{Z \times Z}
\Delta .
$$
\end{rem}
\begin{proof}
We put $Z_i^S =Z_i^T =Z$ ($i=1,2$) and
$$
\aligned
S & =Z_1^S \times {}^{\circ}S^{\circ} \times Z_2^S
= (-1)^{\dim Z \dim S} {}^{\circ}S^{\circ} \times
Z_1^S \times Z_2^S\\
T & =Z_1^T \times {}^{\circ}T^{\circ} \times Z_2^T
=  (-1)^{\dim Z \dim T}
Z_1^T \times Z_2^T \times {}^{\circ}T^{\circ}.
\endaligned
$$
Then by the definition of fiber product we have
$$
S\times_{Z\times Z}T = (-1)^{\dim Z(\dim S + \dim T)}
{}^{\circ}S^{\circ} \times
Z\times Z
\times {}^{\circ}T^{\circ}.
$$
On the other hand, we find that
$$
\aligned
& \Delta
\times _{Z \times Z}  (S \times T) \times_{Z \times Z}
\Delta \\
& = \Delta
\times _{Z \times Z}  (Z_1^S \times {}^{\circ}S^{\circ} \times Z_2^S \times
Z_1^T \times {}^{\circ}T^{\circ} \times Z_2^T)
\times _{Z\times Z} \Delta \\
& = (-1)^{\dim Z(\deg S +\dim T)}
\Delta
\times _{Z \times Z}
(Z_1^S \times Z_1^T \times {}^{\circ}S^{\circ}
\times {}^{\circ}T^{\circ} \times
Z_2^S \times Z_2^T)
\times_{Z\times Z} \Delta\\
& =
(-1)^{\dim Z(\deg S +\dim T)}
\Delta \times {}^{\circ}S^{\circ}
\times {}^{\circ}T^{\circ} \times \Delta \\
& =
(-1)^{\dim Z(\deg S +\dim T)}
S\times_{Z\times Z}T.
\endaligned
$$
\end{proof}

\begin{defn}\label{def:PDpairing}
 We define a pairing on the set of differential forms $\Omega(Z)$ by
$$\langle \eta_1, \eta_2 \rangle_{PD_Z} = \int_Z \eta_1 \wedge \eta_2.$$
\end{defn}

Let $f:\Delta^p \to P \subset Z$ and $g: \Delta^q \to Q \subset Z$ be
smooth singular simplicies in $Z$.
For simplicity of notation, we sometimes denote them by $P$ and $Q$,
respectively.
We denote by $\deg P= \dim Z - \dim P$ and $\deg Q = \dim Z - \dim Q$.
Denote by $PD(P), PD(Q)$ the Poincar\'e dual forms of $P, Q$  
respectively    
satisfying 
\begin{equation}\label{signPD}
\int_{P} v= \int_Z PD(P) \wedge v, \quad
\int_{Q} v= \int_Z PD(Q) \wedge v
\end{equation}
for any differential form $v$. See \cite[Remark 3.5.8 (1)]{fooobook}.
Then for a transversal pair $P,Q$ of complementary dimension we have
$$
\langle PD(P), PD(Q) \rangle_{PD_Z} = \# P\cap Q.
$$
We can regard $P$ and $Q$ as spaces with Kuranishi structures such that
$f$ and $g$ are weakly submersive.
Then we can take their fiber product $P \times_Z Q$.  
If $P$ and $Q$ are
transversal, we have
$$P \cap Q= (-1)^{\deg P \deg Q} P \times_Z Q.$$
Hence, by \cite[Remark 8.4.7]{fooobook2}, we find that
$$\langle PD(P), PD(Q) \rangle_{PD_Z} = (-1)^{\deg P \deg Q} P \times_Z Q$$
in such cases.
In this sense when we calculate signs below, it suffices to 
check them for smooth singular simplicies and we 
use the fiber product description of the pairing.
 
\begin{rem}\label{rempair}
We note that we used a different pairing, which is denoted by $\langle \cdot, \cdot
\rangle_{\text{book}}$ here in
\cite{fooobook2}.
The relation between $\langle \cdot, \cdot \rangle_{PD_Z}$ and $\langle \cdot, \cdot 
\rangle_{\text{book}}$ is given by
$$\langle PD(P),PD(Q) \rangle_{PD_Z}=(-1)^{\deg P \deg Q} \langle P, Q
\rangle_{\text{book}}.$$
This is the only point where the sign convention of this paper is
different from the one in \cite{fooobook},
\cite{fooobook2}.
Hereafter we simply write
$$
\langle P, Q
\rangle_{PD_Z} = \langle PD(P),PD(Q) \rangle_{PD_Z}
$$
for smooth singular chains $P,Q$.
\end{rem}

The following lemma is an immediate consequence of
Lemma \ref{LemmaB}.

\begin{lem}\label{LemmaC}
We have
$$
\langle P, Q\rangle_{PD_Z} =
(-1)^{\deg P\dim Q}
\Delta \times_{Z\times Z}
(P\times Q).
$$
\end{lem}

Next, we decompose the diagonal set $\Delta$
into sum of $\text{\bf e}_I \times \text{\bf e}_J$.
(Here we note that in this section we use homological notation.
That is, $\text{\bf e}_I$
denotes a chain which is the Poincar\'e dual via
\eqref{signPD} to
$\text{\bf e}_I$ used in Section \ref{sec:statements}, Section \ref{sec:annuli} -- Section \ref{sec:ResHess} and
Section \ref{sec:cyclic cohomology}.)
Namely we show the following.

\begin{lem}\label{LemmaE}
$
\Delta =
\sum_{I, J} (-1)^{\vert I \vert \vert J \vert}g^{IJ}
\text{\bf e}_I \times \text{\bf e}_J,
$
where $\vert I \vert = \deg \text{\bf e}_{I}$,
$\vert J \vert = \deg \text{\bf e}_{J}$,
$g_{IJ}=\langle \text{\bf e}_I, \text{\bf e}_J
\rangle_{PD_L}$ and
$g^{IJ}$ is its inverse matrix.
\end{lem}

\begin{rem}
This is a standard fact.
See \cite[Lemma 11.22]{BT}, for example.
Note that the formula in \cite {BT} is written in terms of differential forms, while ours above is written as the equality of homology classes.
To translate the equality of homology classes into
one of differential forms, we need to use
the Poincar\'e duality.
There are different sigh conventions of
the Poincar\'e duality.
In fact,
our convention (\ref{signPD}) is different form
\cite[(6.21)]{BT}.
We denote by $PD^{BT}(P)$ the Poincar\'e dual to $P$
in the sense of Bott-Tu's book.
Then we can easily see that
\begin{equation}\label{BTPDsign}
PD^{BT}(P)= (-1)^{\deg P \dim P}PD(P).
\end{equation}
By noticing this difference, we can see that
Lemma \ref{LemmaE} is the same as \cite[Lemma 11.22]{BT}.
(Anyway, our equality itself is an equality of homology classes. So the sign here does not depend on the choice of conventions.)
\end{rem}

Now when $Z=L$,
the Poincar\'e pairing
$\langle \cdot , \cdot \rangle_{PD_L}$ itself
does not satisfy the cyclic symmetry like
\eqref{cycsymsign}.
However, if we use
$\overline{\frak m}_2 (P,Q)$ for the pairing, it satisfies this
cyclic property.
Using the cohomological notation, we define

\begin{equation}\label{cycpair}
\langle P,Q \rangle_{\text{cyc}} :=
\int_{L} \overline{\frak m}_2 (P,Q).
\end{equation}
We also extend it over $\Lambda$ coefficients naturally and denote it by the same symbol.
The difference between
$\langle P, Q \rangle_{PD_L}$
and $\langle P,Q \rangle_{\text{cyc}}$
is described by the following lemma which is nothing but \cite[Corollary 8.6.4]{fooobook2}.

\begin{lem}\label{Lemmadifference}
$\langle P, Q \rangle_{PD_L}
=(-1)^{\deg P (\deg Q +1)} \langle P,Q \rangle_{\text{\rm cyc}}$.
\end{lem}
We define

\begin{equation}\label{defn:h1}
h_{IJ}:= \langle \text{\bf e}_I,
\text{\bf e}_J
\rangle_{\text{\rm cyc}} = (-1)^{\vert I \vert (\vert J \vert +1)}
\langle \text{\bf e}_I,
\text{\bf e}_J \rangle_{PD_L}
=(-1)^{\vert I \vert (\vert J \vert +1)}g_{IJ},
\end{equation}
with $\vert I \vert = \deg \text{\bf e}_{I}$,
$\vert J \vert = \deg \text{\bf e}_{J}$.
Then we have

\begin{equation}\label{defn:h2}
h^{IJ}=(-1)^{(\vert I \vert +1)\vert J \vert}g^{IJ}.
\end{equation}

\begin{lem}\label{LemmaD}
$$
\aligned
& g^{IJ}\langle \overline{\frak m}_2(\text{\bf e}_I, pt),~
\overline{\frak m}_2(\text{\bf e}_J, pt)
\rangle_{PD_L} \\
= &
\sum_{I,J,A,B,C,D}
(-1)^{\vert A \vert \vert J \vert}
h^{IJ}h^{AB}h^{C0}h^{D0}~
\langle
\text{\bf e}_A \cup^{{\frak b}, b} \text{\bf e}_I,
\text{\bf e}_C \rangle_{\text{\rm cyc}}~
\langle
\text{\bf e}_B \cup^{{\frak b}, b} \text{\bf e}_J,
\text{\bf e}_D \rangle_{\text{\rm cyc}}.
\endaligned
$$
\end{lem}
\begin{proof}
We first note 
the following general formula.
\begin{sublem}\label{sublemmageneral}
$$\langle P,Q \rangle_{PD_L}
=\sum _{A,B}(-1)^{\vert A \vert \vert B\vert}
g^{AB}\langle P,\text{\bf e}_A \rangle_{PD_L}
\langle Q, \text{\bf e}_B \rangle_{PD_L}.
$$
\end{sublem}
\begin{proof}
Put $P=\sum P^C \text{\bf e}_C$ and
$Q=\sum Q^D \text{\bf e}_D$.
Then we have
$$
\aligned
\langle P, \text{\bf e}_A \rangle_{PD_L}
& = \sum P^C \langle \text{\bf e}_C, \text{\bf e}_A
\rangle_{PD_L} = \sum P^C g_{CA}\\
\langle Q, \text{\bf e}_B \rangle_{PD_L}
& = \sum Q^D \langle \text{\bf e}_D, \text{\bf e}_B
\rangle_{PD_L} = \sum Q^D g_{DB}\\
\langle P, Q \rangle_{PD_L}
& = \sum P^C Q^D \langle
\text{\bf e}_C, \text{\bf e}_D
\rangle_{PD_L} = \sum P^C Q^D g_{CD}.
\endaligned
$$
Therefore we have
$$
\aligned
g^{AB} \langle P, \text{\bf e}_A \rangle_{PD_L}
\langle Q , \text{\bf e}_B \rangle_{PD_L}
& =
\sum g^{AB} P^C g_{CA}Q^Dg_{DB} \\
& = \sum
P^C Q^D g^{AB}g_{CA}g_{DB}.
\endaligned
$$
Here we note that $g_{DB}=(-1)^{\vert D \vert \vert B \vert}g_{BD}$. Thus we get
$$
g^{AB}g_{CA}g_{DB}
=(-1)^{\vert D \vert \vert B \vert}
g_{CA}g^{AB}g_{BD}
=(-1)^{\vert A\vert \vert B \vert}g_{CD},
$$
because $\vert A\vert =\vert D \vert$.
Hence we obtain Sublemma \ref{sublemmageneral}.
\end{proof}
By this sublemma, we have
$$
\aligned
& \sum_{I,J} g^{IJ}\langle \overline{\frak m}_2(\text{\bf e}_I, pt),~
\overline{\frak m}_2(\text{\bf e}_J, pt)
\rangle_{PD_L} \\
= &
\sum_{I,J,A,B}(-1)^{\vert A \vert \vert B \vert} g^{IJ}g^{AB}
\langle \overline{\frak m}_2(\text{\bf e}_I, pt),~
\text{\bf e}_{A}
\rangle_{PD_L}
\langle
\overline{\frak m}_2(\text{\bf e}_J, pt),~
\text{\bf e}_{B}
\rangle_{PD_L}.
\endaligned
$$
Now in order to use the cyclic symmetry,
we replace the pairing $\langle \cdot, \cdot \rangle_{PD_L}$ by
$\langle \cdot , \cdot \rangle_{\text{\rm cyc}}$
defined by (\ref{cycpair}).
Then by using
Lemma \ref{Lemmadifference} and (\ref{defn:h2}), we find that the above is equal to
$$
\sum_{I,J,A,B}(-1)^{\vert A \vert \vert B \vert + \gamma_1}
h^{IJ}h^{AB}~
\langle \overline{\frak m}_2(\text{\bf e}_I, pt),~
\text{\bf e}_{A}\rangle_{\text{\rm cyc}} ~
\langle
\overline{\frak m}_2(\text{\bf e}_J, pt),~
\text{\bf e}_{B}\rangle_{\text{\rm cyc}}.
$$
Here
$$
\aligned
\gamma_1 & =
\deg \overline{\frak m}_2(\text{\bf e}_I, pt)
(\vert A \vert +1)
+
\deg \overline{\frak m}_2(\text{\bf e}_J, pt)
(\vert B \vert +1)
+(\vert I \vert +1) \vert J \vert +(\vert A \vert +1)\vert B \vert \\
& = (\vert I \vert + n)(\vert A \vert +1)
+
(\vert J \vert + n)(\vert B \vert +1)
+(\vert I \vert +1) \vert J \vert +(\vert A \vert +1)\vert B \vert,
\endaligned
$$
where $n=\dim L$.
We note that
\begin{equation}\label{ABIJ}
\aligned
\vert A \vert + \vert B \vert & = \vert I \vert + \vert J \vert
=n \\
\vert I \vert + \vert A \vert & = \vert J \vert + \vert B \vert
= 0.
\endaligned
\end{equation}
In particular, we have $\vert A \vert \equiv \vert I \vert$ and
$\vert B \vert \equiv \vert J \vert$.
Therefore it follows that
\begin{equation}\label{gamma1}
\gamma_1 \equiv n \quad \mod 2.
\end{equation}
Then by using the cyclic symmetry
\eqref{cycsymsign}, we have
$$
\sum_{I,J,A,B}(-1)^{\vert A \vert \vert B \vert + \gamma_1 + \gamma_2}
h^{IJ}h^{AB}
~\langle \overline{\frak m}_2(\text{\bf e}_A, \text{\bf e}_I),~
pt \rangle_{\text{\rm cyc}}
~\langle
\overline{\frak m}_2(\text{\bf e}_B, \text{\bf e}_J),~
pt \rangle_{\text{\rm cyc}}
$$
with
\begin{equation}\label{gamma2}
\gamma_2 = (\vert A \vert +1)(\vert I \vert +n)
+ (\vert B \vert +1 )(\vert J \vert +n).
\end{equation}
Furthermore
from \eqref{deformcup}, we obtain
$$
\sum_{I,J,A,B,C,D}(-1)^{\vert A \vert \vert B \vert + \gamma_1 + \gamma_2 + \gamma_3}
h^{IJ}h^{AB}h^{C0}h^{D0}
~\langle\text{\bf e}_A \cup^{\frak b,b}
\text{\bf e}_I,~
\text{\bf e}_C \rangle_{\text{\rm cyc}}
~\langle
\text{\bf e}_B \cup^{\frak b,b}\text{\bf e}_J,~
\text{\bf e}_D \rangle_{\text{\rm cyc}},
$$
where
\begin{equation}\label{gamma3}
\gamma_3 = \vert A \vert (\vert I \vert +1) + \vert B \vert
(\vert J \vert +1).
\end{equation}
By taking (\ref{ABIJ}) into account, (\ref{gamma1}), (\ref{gamma2}) and (\ref{gamma3}) yield
$$
\vert A \vert \vert B \vert
+ \gamma_1 + \gamma_2 + \gamma_3
\equiv
\vert A \vert \vert B \vert
\equiv
\vert A \vert \vert J \vert.
$$
This finishes the proof of Lemma \ref{LemmaD}.
\end{proof}

\subsection{Proofs of signs in Theorem \ref{annulusmain} and
Proposition \ref{mainprop}}
\label{subsec:signproof}
Let $P$ and $Q$ be smooth singular simplicies of $L$.
At the final stage, we will consider the case $P=Q=pt$, a point class of $L$
which corresponds to $\text{vol}_L$.
\par\smallskip
\noindent
{\bf Step 1:}
Firstly we explore orientations on various spaces
of fiber products using the evaluation
maps at the boundary marked points.

\par
We apply Lemma \ref{LemmaA} to the situation that
$$
\text{ev}_{i}^{\beta_1} : \CM^{\text{main}}_{3,0}(\beta_1) \to
L, \quad
\text{ev}_{i}^{\beta_2} : \CM^{\text{main}}_{3,0}(\beta_2) \to
L, \quad (i=0,1).
$$
Note that since the Maslov index $\mu(\beta)$ is even, we have $\dim \CM^{\text{main}}_{3,0}(\beta)= n + \mu(\beta) \equiv n$ and $\deg \CM^{\text{main}}_{3,0}(\beta)\equiv 0 \mod 2$
for any $\beta$.
Thus we have
$$
\aligned
& \Big(
\CM^{\text{main}}_{3,0}(\beta_1)
{}_{(\text{ev}_0^{\beta_1}, \text{ev}_1^{\beta_1})}
\times_{(\text{ev}_0^{\beta_2}, \text{ev}_1^{\beta_2})}
\CM^{\text{main}}_{3,0}(\beta_2)
\Big)
{}_{(\text{ev}_2^{\beta_1}, \text{ev}_2^{\beta_2})}\times
\big(P\times Q \big) \\
& =
(-1)^{\epsilon_1}
\left(
\Delta \times_{L\times L}
\Big( \CM^{\text{main}}_{3,0}(\beta_1) \times
\CM^{\text{main}}_{3,0}(\beta_2)  \Big)
{}_{L\times L}\times \Delta \right)
{}_{(\text{ev}_2^{\beta_1}, \text{ev}_2^{\beta_2})}\times
\big(P\times Q \big),
\endaligned
$$
where
$$
\epsilon_1 \equiv
n(\deg \CM^{\text{main}}_{3,0}(\beta_1) + \dim
\CM^{\text{main}}_{3,0}(\beta_2))
\equiv
n \mod 2.
$$
By the associativity (\cite[Lemma 8.2.3]{fooobook2}) and
Lemma \ref{LemmaE}, we have
\begin{equation}\label{eqdisk1}
\aligned
& \Big(
\CM^{\text{main}}_{3,0}(\beta_1)
{}_{(\text{ev}_0^{\beta_1}, \text{ev}_1^{\beta_1})}
\times_{(\text{ev}_0^{\beta_2}, \text{ev}_1^{\beta_2})}
\CM^{\text{main}}_{3,0}(\beta_2)
\Big)
{}_{(\text{ev}_2^{\beta_1}, \text{ev}_2^{\beta_2})}\times
\big(P\times Q \big) \\
& = \sum_{I,J}
(-1)^{\epsilon_1 +\vert I \vert \vert J \vert}
g^{IJ} \\
& \qquad \Delta \times_{L\times L}
\left(
\Big(
\CM^{\text{main}}_{3,0}(\beta_1) \times
\CM^{\text{main}}_{3,0}(\beta_2)
\Big)
\times_{L\times L}
\big( \text{\bf e}_I \times \text{\bf e}_J \big)
{}_{(\text{ev}_2^{\beta_1}, \text{ev}_2^{\beta_2})}\times
\big(P\times Q \big)\right).
\endaligned
\end{equation}
We show the following lemma.

\begin{lem}\label{Lemmadisk1}
$$
\aligned
& \Delta \times_{L\times L}
\left(
\Big(
\CM^{\text{\rm main}}_{3,0}(\beta_1) \times
\CM^{\text{\rm main}}_{3,0}(\beta_2)
\Big)
\times_{L\times L}
\big( \text{\bf e}_I \times \text{\bf e}_J \big)
{}_{(\text{ev}_2^{\beta_1}, \text{ev}_2^{\beta_2})}\times
\big(P\times Q \big)\right) \\
& =
(-1)^{n + \vert J \vert \deg P}
\Delta \times_{L\times L}
\Big(
\CM^{\text{\rm main}}_{3,0}(\beta_1;\text{\bf e}_I, P) \times
\CM^{\text{\rm main}}_{3,0}(\beta_2;\text{\bf e}_J, Q)
\Big).
\endaligned
$$
\end{lem}

\begin{proof}
We write
the LHS in Lemma \ref{Lemmadisk1} as
$$
\Delta \times_{L\times L}
\left(
\Big(
\CM^{\text{\rm main}}_{3,0}(\beta_1) \times
\CM^{\text{\rm main}}_{3,0}(\beta_2)
\Big)
\times_{L_1\times L_2}
\big( \text{\bf e}_I \times \text{\bf e}_J \big)
\times_{L_3 \times L_4}
\big(P\times Q \big)\right)
$$
and put
\begin{equation}\label{CX}
\aligned
\text{\bf e}_I & = L_1 \times {}^{\circ}\text{\bf e}_{I},
\quad
\text{\bf e}_J = L_2 \times {}^{\circ}\text{\bf e}_{J}, \\
P & = L_3 \times {}^{\circ}P, \quad
Q = L_4 \times {}^{\circ}Q, \\
\CX & =
\Delta \times_{L\times L}
\Big(
\CM^{\text{\rm main}}_{3,0}(\beta_1) \times
\CM^{\text{\rm main}}_{3,0}(\beta_2)
\Big) \\
& = \CX^{\circ} \times L_1 \times L_2 \times L_3 \times L_4
= \CX^{\circ} \times L_3 \times L_4 \times L_1 \times L_2.
\endaligned
\end{equation}
By the associativity again, it is equal to
$$
\left(
\Delta \times_{L\times L}
\Big(
\CM^{\text{\rm main}}_{3,0}(\beta_1) \times
\CM^{\text{\rm main}}_{3,0}(\beta_2)
\Big)\right)
\times_{L_1\times L_2}
\big( \text{\bf e}_I \times \text{\bf e}_J \big)
\times_{L_3 \times L_4}
\big(P\times Q \big).
$$
Using the notation introduced in (\ref{CX}), we can rewrite it
as
\begin{equation}\label{eqepsilon2}
\aligned
& \big(
\CX^{\circ} \times L_3 \times L_4 \times L_1 \times L_2
\big)
\times_{L_1\times L_2}
\big(
L_1 \times {}^{\circ}\text{\bf e}_{I} \times
L_2 \times {}^{\circ}\text{\bf e}_{J}
\big) \\
& \qquad \qquad \times_{L_3\times L_4}
\big(
L_3 \times {}^{\circ}P \times
L_4 \times {}^{\circ}Q
\big) \\
& =
(-1)^{\epsilon_2}
\big(
\CX^{\circ} \times L_3 \times L_4 \times L_1 \times L_2
\big)
\times_{L_1\times L_2}
\big(
L_1 \times L_2 \times {}^{\circ}\text{\bf e}_{I}
\times {}^{\circ}\text{\bf e}_{J}
\big) \\
& \qquad \qquad
\times_{L_3\times L_4}
\big(
L_3 \times L_4 \times {}^{\circ}P
\times {}^{\circ}Q
\big),
\endaligned
\end{equation}
where $\epsilon_2$ is given by
\begin{equation}\label{epsilon2}
\epsilon_2 =
n(\vert  I \vert + \deg P).
\end{equation}
If we put $\text{\bf e}_I \times \text{\bf e}_J
=L_1 \times L_2 \times {}^{\circ}(\text{\bf e}_I \times
\text{\bf e}_J)$ and
$P\times Q = L_3 \times L_4 \times {}^{\circ}(P\times Q)$, we can easily see

\begin{equation}\label{PXQ}
\aligned
{}^{\circ}(\text{\bf e}_I \times
\text{\bf e}_J) & =
(-1)^{n\vert I\vert}{}^{\circ}\text{\bf e}_{I}
\times {}^{\circ}\text{\bf e}_{J}, \\
{}^{\circ}(P\times Q) & =
(-1)^{n\deg P}{}^{\circ}P \times {}^{\circ}Q.
\endaligned
\end{equation}
Therefore (\ref{eqepsilon2}) is equal to
\begin{equation}\label{eqepsilon34}
\aligned
& (-1)^{\epsilon_2 + \epsilon_3 + \epsilon_4}
\Big(
\CX^{\circ} \times L_3 \times L_4 \times L_1 \times L_2
\times
{}^{\circ}(\text{\bf e}_I \times
\text{\bf e}_J)
\Big)
\times_{L_3 \times L_4}
\Big(
L_3 \times L_4 \times {}^{\circ}(P
\times Q)
\Big) \\
& =
(-1)^{\epsilon_2 + \epsilon_3 + \epsilon_4}
\Big(
\CX^{\circ}  \times L_1 \times L_2
\times
{}^{\circ}(\text{\bf e}_I \times
\text{\bf e}_J)
\times L_3 \times L_4
\Big)
\times_{L_3 \times L_4}
\Big(
L_3 \times L_4 \times {}^{\circ}(P
\times Q)
\Big) \\
& =
(-1)^{\epsilon_2 + \epsilon_3 + \epsilon_4}
\CX^{\circ}
\times L_1 \times L_2 \times
{}^{\circ}(\text{\bf e}_I \times
\text{\bf e}_J)
\times
L_3 \times L_4 \times {}^{\circ}(P
\times Q),
\endaligned
\end{equation}
where
$$
\epsilon_3 = n \vert I \vert, \qquad
\epsilon_4 = n \deg P.
$$
Next we study $\CX^{\circ}$. We put
$$
\aligned
\CM^{\text{main}}_{3,0}(\beta_1) & =
\CM^{\text{main}}_{3,0}(\beta_1)^{\circ} \times
L_1 \times L_3, \\
\CM^{\text{main}}_{3,0}(\beta_2) & =
\CM^{\text{main}}_{3,0}(\beta_2)^{\circ} \times
L_2 \times L_4.
\endaligned
$$
Then we find that
$$
\aligned
\CX & =
\Delta \times_{L\times L}
\Big(
\CM^{\text{main}}_{3,0}(\beta_1)^{\circ} \times
L_1 \times L_3 \times
\CM^{\text{main}}_{3,0}(\beta_2)^{\circ} \times
L_2 \times L_4
\Big) \\
& =
(-1)^{n} \Delta \times_{L\times L}
\Big(
\CM^{\text{main}}_{3,0}(\beta_1)^{\circ} \times
\CM^{\text{main}}_{3,0}(\beta_2)^{\circ} \times
L_1 \times L_2 \times L_3 \times L_4
\Big) \\
& =
(-1)^{n} \Delta \times_{L\times L}
\Big(
\CM^{\text{main}}_{3,0}(\beta_1)^{\circ} \times
\CM^{\text{main}}_{3,0}(\beta_2)^{\circ} \Big)
\times
L_1 \times L_2 \times L_3 \times L_4.
\endaligned
$$
Hence we have
\begin{equation}\label{CXcirc}
\CX^{\circ} =
(-1)^n
\Delta \times_{L\times L}
\Big(
\CM^{\text{main}}_{3,0}(\beta_1)^{\circ} \times
\CM^{\text{main}}_{3,0}(\beta_2)^{\circ} \Big).
\end{equation}
We substitute (\ref{CXcirc}) into (\ref{eqepsilon34}). Then
the LHS in Lemma \ref{Lemmadisk1} is equal to
$$
\aligned
& (-1)^{\epsilon_2 + \epsilon_3 + \epsilon_4 + n}
\Delta \times_{L\times L}
\Big(
\CM^{\text{main}}_{3,0}(\beta_1)^{\circ} \times \\
& \qquad \qquad \qquad \CM^{\text{main}}_{3,0}(\beta_2)^{\circ}
\times L_1 \times L_2 \times
{}^{\circ}(\text{\bf e}_I \times
\text{\bf e}_J)
\times L_3 \times L_4
\times {}^{\circ}(P
\times Q)
\Big)\\
& =
(-1)^{\epsilon_2 + n + \epsilon_5}
\Delta \times_{L\times L}
\Big(
\CM^{\text{main}}_{3,0}(\beta_1)^{\circ} \times \\
& \qquad \qquad \qquad
L_1 \times {}^{\circ}\text{\bf e}_I \times
L_3 \times {}^{\circ}P
\times
\CM^{\text{main}}_{3,0}(\beta_2)^{\circ}
\times  L_2 \times
{}^{\circ}\text{\bf e}_J
\times L_4
\times {}^{\circ}Q
\Big).
\endaligned
$$
Here an elementary calculation shows that
\begin{equation}\label{epsilon5}
\epsilon_5 =
(n + \deg P) \dim \text{\bf e}_J.
\end{equation}
By   \cite[Remark 8.2.6]{fooobook2}
it is equal to the following fiber product
$$
\aligned
(-1)^{\epsilon_2 + n + \epsilon_5}
& \Delta\times_{L\times L} \\
& \left(
\Big(
\CM^{\text{main}}_{3,0}(\beta_1) \times_{L_1 \times L_3}
(\text{\bf e}_{I} \times P)
\Big)
\times
\Big(
\CM^{\text{main}}_{3,0}(\beta_2) \times_{L_2 \times L_4}
(\text{\bf e}_{J} \times Q)
\Big)
\right).
\endaligned
$$
Then by \cite[Definition 8.4.1]{fooobook2} we get
$$
(-1)^{\epsilon_2 + n + \epsilon_5 + \epsilon_6}
\Delta\times_{L\times L}
\Big(
\CM^{\text{main}}_{3,0}(\beta_1 ;
\text{\bf e}_{I}, P)
\times
\CM^{\text{main}}_{3,0}(\beta_2 ;
\text{\bf e}_{J}, Q)
\Big),
$$
where
\begin{equation}\label{epsilon6}
\epsilon_6 =
(n+1) (\vert I \vert + \vert J \vert)=(n+1)n \equiv 0 \mod 2.
\end{equation}
Therefore, from (\ref{epsilon2}), (\ref{epsilon5}) and
(\ref{epsilon6}) we obtain Lemma \ref{Lemmadisk1}.
Here we also used (\ref{ABIJ}) again.
\end{proof}
Next task is to rewrite the RHS in Lemma \ref{Lemmadisk1}
by using the Poincar\'e pairing which is introduced in Definition
\ref{def:PDpairing}.

\begin{lem}\label{LemmadiskPD}
$$
\aligned
& \Delta \times_{L\times L}
\Big(
\CM^{\text{\rm main}}_{3,0}(\beta_1 ;
\text{\bf e}_{I}, P)
\times
\CM^{\text{\rm main}}_{3,0}(\beta_2 ;
\text{\bf e}_{J}, Q)
\Big) \\
& =
(-1)^{\epsilon_7}
\langle
{\frak m}_{2,\beta_1}(\text{\bf e}_I, P), ~
{\frak m}_{2,\beta_2}(\text{\bf e}_J, Q)
\rangle_{PD_L},
\endaligned
$$
where
$$
\epsilon_7 =
(\vert I \vert +\deg P)(n+\vert J\vert +\deg Q).
$$
\end{lem}
\begin{proof}
By Lemma \ref{LemmaC},
$\epsilon_7 =
\deg \CM^{\text{\rm main}}_{3,0}(\beta_1 ;
\text{\bf e}_{I}, P)
\dim \CM^{\text{\rm main}}_{3,0}(\beta_2 ;
\text{\bf e}_{J}, Q)
\equiv (\vert I \vert +\deg P)(n+\vert J\vert +\deg Q) \mod 2$.
\end{proof}
An easy calculation shows that
$$
\aligned
\epsilon & :=
\epsilon_1 + \vert I \vert \vert J \vert + n + \vert J \vert \deg P
+ \epsilon_7 \\
& \equiv
(\deg P + \vert I \vert)
(\deg Q + n) \mod 2.
\endaligned
$$
Hence combing it with (\ref{eqdisk1}), Lemmas \ref{Lemmadisk1} and
\ref{LemmadiskPD}, we obtain
the following.
\begin{prop}\label{propdisk}
$$
\aligned
& \Big(
\CM^{\text{\rm main}}_{3,0}(\beta_1)
{}_{(\text{\rm ev}_0^{\beta_1}, \text{\rm ev}_1^{\beta_1})}
\times_{(\text{\rm ev}_0^{\beta_2}, \text{\rm ev}_1^{\beta_2})}
\CM^{\text{\rm main}}_{3,0}(\beta_2)
\Big)
{}_{(\text{\rm ev}_2^{\beta_1}, \text{\rm ev}_2^{\beta_2})}\times
\big(P\times Q \big) \\
& =
\sum_{I,J}
(-1)^{\epsilon}
g^{IJ}
\langle
{\frak m}_{2,\beta_1}(\text{\bf e}_I, P), ~
{\frak m}_{2,\beta_2}(\text{\bf e}_J, Q)
\rangle_{PD_L},
\endaligned
$$
where
$$
\epsilon = (\deg P + \vert I \vert)
(\deg Q+n).
$$
\end{prop}

\par\noindent
{\bf Step 2:}
Next, we explore orientations on various spaces of fiber products using the evaluation maps
at the interior marked points.
In this case we first apply Lemma \ref{LemmaB} to the situation that
$$
\text{ev}^{+}_{\beta'} : \CM_{1,1}(\beta') \to X, \quad
\text{ev}^{+}_{\beta''} : \CM_{1,1}(\beta'') \to X.
$$
Then we have
\begin{equation}\label{eqannuli1}
\aligned
& \Big( \CM_{1,1}(\beta')
{}_{\text{ev}^{+}_{\beta'}}
\times_{\text{ev}^{+}_{\beta''}}
\CM_{1,1}(\beta'') \Big)
\times _{L_1 \times L_2}
(P\times Q) \\
& =
(-1)^{\delta_1}
\Delta_X \times _{X\times X}
\Big( \CM_{1,1}(\beta')
\times
\CM_{1,1}(\beta'') \Big)
\times _{L_1 \times L_2}
(P\times Q)
\endaligned
\end{equation}
with
$$
\delta_1 =
\dim X \deg \CM_{1,1}(\beta')
\equiv 0 \mod 2.
$$
Now we show the following.

\begin{lem}\label{Lemmaannuli1}
$$
\aligned
& \Delta_X \times_{X\times X}
\Big( \CM_{1,1}(\beta')
\times
\CM_{1,1}(\beta'') \Big)
\times _{L_1 \times L_2}
(P\times Q) \\
& =
\Delta_X \times _{X\times X}
\Big( \CM_{1,1}(\beta' ; P)
\times
\CM_{1,1}(\beta''; Q)
\Big).
\endaligned
$$
\end{lem}

\begin{proof}
Put
$$
\aligned
\CM_{1,1}(\beta') & =
X \times {}^{\circ}\CM_{1,1}(\beta')^{\circ}
\times L_1, \\
\CM_{1,1}(\beta'') & =
X \times {}^{\circ}\CM_{1,1}(\beta'')^{\circ}
\times L_2, \\
P & = L_1 \times {}^{\circ}P, \\
Q & = L_2 \times {}^{\circ}Q.
\endaligned
$$
We note that
$\dim {}^{\circ}\CM_{1,1}(\beta')^{\circ}
\equiv
\dim {}^{\circ}\CM_{1,1}(\beta'')^{\circ} \equiv 0 \mod 2$.
Then the LHS in Lemma \ref{Lemmaannuli1} can be written as
$$
\aligned
& \Delta_X \times_{X\times X}
\Big( \CM_{1,1}(\beta')
\times
\CM_{1,1}(\beta'') \Big)
\times _{L_1 \times L_2}
(P\times Q) \\
& =
\Delta_X \times_{X\times X}
\Big(
X \times {}^{\circ}\CM_{1,1}(\beta')^{\circ}
\times L_1
\times
X \times {}^{\circ}\CM_{1,1}(\beta'')^{\circ}
\times L_2
\Big) \\
& \qquad \qquad \qquad
\times_{L_1 \times L_2}
\Big(L_1 \times {}^{\circ}P
\times L_2 \times {}^{\circ}Q \Big) \\
& =
(-1)^{\delta_2}
\Delta_X \times_{X\times X}
\Big(
X \times {}^{\circ}\CM_{1,1}(\beta')^{\circ}
\times
X \times {}^{\circ}\CM_{1,1}(\beta'')^{\circ}
\times L_1 \times L_2
\Big) \\
& \qquad \qquad \qquad
\times_{L_1 \times L_2}
\Big(L_1 \times L_2 \times {}^{\circ}P
\times {}^{\circ}Q \Big) \\
\endaligned
$$
where
$$
\delta_2 =
n \deg P.
$$
It is equal to
$$
\aligned
& =
(-1)^{\delta_2}
\Delta_X \times_{X\times X}
\Big(
X \times {}^{\circ}\CM_{1,1}(\beta')^{\circ}
\times
X \times {}^{\circ}\CM_{1,1}(\beta'')^{\circ} \\
& \qquad \qquad \qquad \times L_1 \times L_2
\times {}^{\circ}P
\times {}^{\circ}Q) \Big) \\
& =
(-1)^{\delta_2}
\Delta_X \times_{X\times X}
\Big(
X \times {}^{\circ}\CM_{1,1}(\beta')^{\circ} \times L_1
\times
X \times {}^{\circ}\CM_{1,1}(\beta'')^{\circ} \\
& \qquad \qquad \qquad  \times L_2
\times {}^{\circ}P
\times {}^{\circ}Q) \Big) \\
& =
(-1)^{\delta_2+\delta_3}
\Delta_X \times_{X\times X}
\Big(
X \times {}^{\circ}\CM_{1,1}(\beta')^{\circ} \times L_1
\times {}^{\circ}P \times
X \times {}^{\circ}\CM_{1,1}(\beta'')^{\circ} \\
& \qquad \qquad \qquad  \times L_2
\times {}^{\circ}Q) \Big), \\
\endaligned
$$
where $\delta_3 = n \deg P$.
Since $\delta_2 + \delta_3 \equiv 0 \mod 2$ and
$$
\aligned
& \Big(
X \times {}^{\circ}\CM_{1,1}(\beta')^{\circ} \times L_1
\times {}^{\circ}P \times
X \times {}^{\circ}\CM_{1,1}(\beta'')^{\circ}
\times L_2
\times {}^{\circ}Q) \Big) \\
& =
\Big( \CM_{1,1}(\beta')\times_{L_1} P \Big)
\times
\Big( \CM_{1,1}(\beta'')\times_{L_2} Q \Big),
\endaligned
$$
we find that the LHS in Lemma \ref{Lemmaannuli1}
is equal to
$$
\Delta_X \times_{X\times X}
\left(\Big( \CM_{1,1}(\beta')\times_{L_1} P \Big)
\times
\Big( \CM_{1,1}(\beta'')\times_{L_2} Q \Big)\right).
$$
Moreover, \cite[Definition 8.10.2]{fooobook2} shows
$$
\Big( \CM_{1,1}(\beta')\times_{L_1} P \Big)
\times
\Big( \CM_{1,1}(\beta'')\times_{L_2} Q \Big)
 =
\CM_{1,1}(\beta' ; P)
\times
\CM_{1,1}(\beta''; Q).
$$
Thus this finishes the proof of Lemma \ref{Lemmaannuli1}.
\end{proof}
By using the Poincar\'e pairing on $X$ (which is defined in
Definition \ref{def:PDpairing} generally),
we find from Lemma \ref{LemmaC} that

\begin{equation}\label{eqannuli2}
\aligned
&
\Delta_X \times _{X\times X}
\Big( \CM_{1,1}(\beta' ; P)
\times
\CM_{1,1}(\beta''; Q)
\Big) \\
=
& (-1)^{\delta_4}
\langle
\CM_{1,1}(\beta' ; P),~
\CM_{1,1}(\beta''; Q) \rangle_{PD_X},
\endaligned
\end{equation}
where
$$
\delta_4 =
\deg \CM_{1,1}(\beta' ; P)
\dim \CM_{1,1}(\beta''; Q)
\equiv
\dim P \dim Q.
$$
Note that $\deg \CM_{1,1}(\beta' ; P)
\equiv \dim \CM_{1,1}(\beta' ; P)$ because
$\dim X$ is even.
Therefore, by (\ref{eqannuli1}), Lemma \ref{Lemmaannuli1} and
(\ref{eqannuli2}), we obtain

\begin{prop}\label{propannuli}
$$
\aligned
& \Big( \CM_{1,1}(\beta')
{}_{\text{ev}^{+}_{\beta'}}
\times_{\text{ev}^{+}_{\beta''}}
\CM_{1,1}(\beta'') \Big)
\times _{L_1 \times L_2}
(P\times Q) \\
& =(-1)^{\delta}
\langle
\CM_{1,1}(\beta' ; P),~
\CM_{1,1}(\beta''; Q) \rangle_{PD_X},
\endaligned
$$
where
$$
\delta = \dim P \dim Q.
$$
\end{prop}

\par\noindent
{\bf Step 3:}
We suppose that
$\beta_1+\beta_2 = \beta' +\beta''$.
We use Proposition \ref{compori}
together with
Proposition \ref{propdisk} and Proposition \ref{propannuli}
to get the following.

\begin{prop}\label{propconclusion}
$$
\aligned
& \langle
{\frak p}_{1,\beta'}(P), ~{\frak p}_{1,\beta''}(Q)
\rangle_{PD_X} \\
& = \sum_{I,J}
(-1)^{\epsilon + \delta + \frac{n(n-1)}{2}}
g^{IJ}
\langle
{\frak m}_{2,\beta_1}(\text{\bf e}_I, P),
~{\frak m}_{2,\beta_2}(\text{\bf e}_J, Q)
\rangle_{PD_L}.
\endaligned
$$
Here
$$
\epsilon = (\deg P + \vert I \vert)
(\deg Q+n), \quad
\delta = \dim P \dim Q,
$$
which are given in
Proposition \ref{propdisk} and Proposition \ref{propannuli} respectively.
\end{prop}
Now, we consider the case
$$
P=Q=pt,
$$
which corresponds to ${\rm vol}_L$.
Then since $\dim P =\dim Q =0$, we have

\begin{equation}\label{eqpreconclusion}
\aligned
& \langle
{\frak p}_{1,\beta'}(pt), ~{\frak p}_{1,\beta''}(pt)
\rangle_{PD_X} \\
& = \sum_{I,J}
(-1)^{\frac{n(n-1)}{2}}
g^{IJ}
\langle
{\frak m}_{2,\beta_1}(\text{\bf e}_I, pt),
~{\frak m}_{2,\beta_2}(\text{\bf e}_J, pt)
\rangle_{PD_L}.
\endaligned
\end{equation}
This finishes verification of sign
in Theorem \ref{annulusmain}.
\par
Combining (\ref{eqpreconclusion}) and Lemma \ref{LemmaD}, we obtain

\begin{equation}\label{eqconclusion}
\aligned
& \langle
{\frak p}_{1,\beta'}(pt), ~{\frak p}_{1,\beta''}(pt)
\rangle_{PD_X} \\
& = \sum_{I,J,A,B,C,D}
(-1)^{\vert A \vert \vert J \vert + \frac{n(n-1)}{2}}
h^{IJ}h^{AB}h^{C0}h^{D0}
\\
& \qquad\qquad\qquad\qquad\qquad
\langle
\text{\bf e}_A \cup^{\frak b,b}\text{\bf e}_I,~\text{\bf e}_C
\rangle_{\text{cyc}}
\langle
\text{\bf e}_B \cup^{\frak b, b}\text{\bf e}_J,~\text{\bf e}_D
\rangle_{\text{cyc}} \\
& = \sum_{I_1,I_2,I_3,J_1,J_2,J_3}
(-1)^{\vert I_1 \vert \vert J_2 \vert + \frac{n(n-1)}{2}}
h^{I_1J_1}h^{I_2J_2}h^{I_30}h^{J_30} \\
& \qquad\qquad\qquad\qquad\qquad
\langle
\text{\bf e}_{I_1} \cup^{\frak b,b}\text{\bf e}_{I_2},~\text{\bf e}_{I_3}
\rangle_{\text{cyc}}
\langle
\text{\bf e}_{J_1} \cup^{\frak b, b}\text{\bf e}_{J_2},
~\text{\bf e}_{J_3}
\rangle_{\text{cyc}}.
\endaligned
\end{equation}
Here $A=I_1, B=J_1, C=I_3, D=J_3, I=I_2$ and $J=J_2$.
By Definition \ref{invariantZ},
the right hand side of \eqref{eqconclusion} is nothing but $Z(\frak b,b)$ and hence we
have verified Proposition \ref{mainprop} with sign.

\begin{rem}\label{remconclude}
If we rewrite the RHS in (\ref{eqconclusion}) in terms of the Poincar\'e pairing
$\langle \cdot, ~ \cdot \rangle_{PD_L}$, we use
Lemma \ref{Lemmadifference} and (\ref{defn:h2}) to obtain

\begin{equation}\label{eqconclusionPD}
\aligned
& \langle
{\frak p}_{1,\beta'}(pt), ~{\frak p}_{1,\beta''}(pt)
\rangle_{PD_X} \\
& = \sum_{I_1,I_2,I_3,J_1,J_2,J_3}
(-1)^{\zeta}
g^{I_1J_1}g^{I_2J_2}g^{I_30}g^{J_30} \\
& \qquad\qquad\qquad\qquad\qquad
\langle
\text{\bf e}_{I_1} \cup^{\frak b,b}\text{\bf e}_{I_2},~\text{\bf e}_{I_3}
\rangle_{PD_L}
\langle
\text{\bf e}_{J_1} \cup^{\frak b, b}\text{\bf e}_{J_2},
~\text{\bf e}_{J_3}
\rangle_{PD_L}
\endaligned
\end{equation}
where
\begin{equation}\label{zeta}
\aligned
\zeta & =
\vert I_1 \vert \vert J_2 \vert + \frac{n(n-1)}{2}
+ (\vert I_1 \vert +1)\vert J_1 \vert
+
(\vert I_2 \vert +1)\vert J_2 \vert
\\
& \quad
+ \vert I_3 \vert + 1+ \vert J_3 \vert +1
+ (\vert I_1 \vert + \vert I_2 \vert )
\vert I_3 \vert
+
(\vert J_1 \vert + \vert J_2 \vert )
\vert J_3 \vert \\
& \equiv
\vert I_1 \vert \vert J_2 \vert + \frac{n(n-1)}{2}  \quad \mod 2
\endaligned
\end{equation}
by taking (\ref{ABIJ}) into account.
\end{rem}

\par

\chapter{Appendix}
\section{Coincidence of the two definitions of
$\delta^{\frak b,b}$}
\label{sec:deltaisthesame}

Let
$$
b = \sum_{i=1}^n \frak x_i {\bf e}_i \in H^1(L(u);\Lambda_0).
$$
Decompose $b = b_0 + b_+$ with $b_0 \in H^1(L(u);\C)$ and $b_+ \in H^1(L(u);\Lambda_+)$.
We start with comparing the two representations
\begin{eqnarray*}
\rho^{b_0}: H_1(L(u);\Z) & \to &\C^*; \quad \gamma \mapsto e^{\gamma \cap b_0} \\
\rho^b: H_1(L(u);\Z) & \to & \Lambda_0 \setminus \Lambda_+; \quad \gamma \mapsto e^{\gamma \cap b}.
\end{eqnarray*}
If we put
$
\frak y_i = e^{\frak x_i} \in \Lambda_0 \setminus \Lambda_+
$
and $\gamma = \sum k_i \text{\bf e}_i^* \in H_1(L(u);\Z)$, then we can write
$
\rho^{b}: H_1(L(u);\Z) \to \Lambda_0 \setminus \Lambda_+
$
as
$$
\rho^{b}(\gamma) = e^{\gamma \cap b}
= \frak y_1^{k_1} \dots \frak y_n^{k_n}
$$
as in Section \ref{sec:frakqreview}. With these notations, we also have
$$
\rho^b = \rho^{b_0} \rho^{b_+}.
$$

Now we unravel the definitions of $\frak m_1^{\frak b, b}$ given in Definition
\ref{defmrho} and $\delta^{\frak b,b}$ in Definition \ref{bdryformdef}. By definition, we have
\begin{equation}\label{eq:m1bb}
\aligned
\frak m_1^{\frak b,b}(h)
& =  \sum_{l_0, l_1 \in \Z_{\ge 0}}\frak m_{1+ l_0 + l_1;\beta}^{\frak b,\rho}(b_+^{l_0},h,b_+^{l_1})
\\
& =  \sum_{l_0, l_1 \in \Z_{\ge 0}} \sum_{\beta}\sum_{\ell=0}^{\infty}
\frac{T^{\beta \cap \omega/2\pi}}{\ell!} \rho^{b_0}(\partial \beta)
\frak q_{\ell;1+l_0 + l_1;\beta}(\frak b^{\ell};b_+^{l_0},h,b_+^{l_1})
\\
& =  \sum_{\beta}\sum_{\ell=0}^{\infty} 
\frac{T^{\beta \cap \omega/2\pi}}{\ell!} \rho^{b_0} (\partial \beta)
\sum_{l_0, l_1 \in \Z_{\ge 0}}  \frak q_{\ell;1+l_0+l_1;\beta}(\frak b^{\ell};b_+^{l_0},h,b_+^{l_1})
\\
&=
\sum_{\beta}\sum_{\ell=0}^{\infty} \frac{T^{\beta \cap \omega/2\pi}}{\ell!} \rho^{b_0} (\partial \beta)
\sum_{k=0}^{\infty} \sum_{l_0+l_1=k}
\frak q_{\ell;1+l_0+l_1;\beta}(\frak b^{\ell};b_+^{l_0},h,b_+^{l_1}).
\endaligned
\end{equation}
Here we used the definition of $\frak m_k^\rho$ given in \eqref{mcyclicdef0}.
On the other hand, $\delta^{\frak b,b}$ is defined as
\begin{equation}\label{eq:deltabb}
\aligned
\delta^{\frak b,b} & = \frak m_{1,0}  + 
\sum_{\beta\ne 0}\sum_{\ell=0}^{\infty}
\frac{T^{\beta\cap \omega/2\pi}}{\ell!} \rho^b(\del \beta)
\delta_{\beta,\ell}^{\frak b} \\
& =  \frak m_{1,0} + \sum_{\beta\ne 0}\sum_{\ell=0}^{\infty}
\frac{T^{\beta\cap \omega/2\pi}}{\ell!}
\rho^{b_0}(\del \beta)
\left(\rho^{b_+}(\del \beta)\delta_{\beta,\ell}^{\frak b}\right)
\endaligned
\end{equation}
in Definition \ref{bdryformdef}. And by the definition \eqref{pushout} of $\delta_{\beta,\ell}^{\frak b}$,
we have
\begin{equation}\label{eq:defdelta}
\delta_{\beta,\ell}^{\frak b}(h) = \frak q_{\ell;1;\beta}(\frak b^{\ell};h).
\end{equation}
By comparing \eqref{eq:m1bb}-\eqref{eq:defdelta}, the proof boils down to proving the following
identity
\begin{equation}\label{eq:identityb}
\sum_{l_0+l_1=k} \frak q_{\ell;1+l_0+l_1;\beta}(\frak b^{\ell};b_+^{l_0},h,b_+^{l_1})
= \frac{\langle\del \beta, b_+\rangle^k}{k!}  \frak q_{\ell;1;\beta}(\frak b^{\ell};h).
\end{equation}
Then this can be proved in the same way as that of  \cite[Lemma 7.1]{fooo09}. 
(Note that because the convention 
of the space $E_{\ell}C$ in the present paper is different from one in \cite{fooo09}, $1/\ell !$ does not appear in the formula above. 
See Remark \ref{quotientsub} and Remark \ref{rem:6.11}.) 
This proves the coincidence of the two definitions.
\qed

\begin{rem}\label{rem:correction11.13}
In the proof of \cite[Lemma 11.8]{fooo08}  which corresponds to \eqref{eq:identityb} with $\frak b =0$,
we wrote down the evaluation maps as  
 in \cite[(11.13) and (11.14)]{fooo08}. 
Even in the case $\frak b =0$, the explicit formulae (11.13) and (11.14) do not hold 
for general class $\beta$, but 
hold only for maps given in 
\cite[Proposition 7.3]{cho-oh}  
with the classes  
$\beta=\beta_1, \dots, \beta_m$. 
However, we do not need such an
explicit description of the evaluation maps 
to prove \cite[Lemma 11.8]{fooo08}  or \eqref{eq:identityb} above. 
It is enough to calculate 
some iterated integral of 
$\text{\rm ev}_1^* b_+ \wedge \dots \wedge \text {\rm ev}_k^* b_+$ 
over a configuration space $\widehat{C}_k$ as follows:
See the first half of the proof of 
\cite[Lemma 11.8]{fooo08}.
Let  $S = \del D$ be the boundary of the unit disk $D=D^2 \subset \C$ and $\beta_D \in H_2(\C,S)$ be the
homology class of the unit disk. We consider the moduli
space $\mathcal M_{k+1}(\C,S;\beta_D)$ and the evaluation map
$
\text{\rm ev} = (\text{\rm ev}_0,\dots,\text{\rm ev}_k) : \mathcal M_{k+1}(\C,S;\beta_D)
\to (S^1)^{k+1}.
$
We fix a point $p_0 \in S \subset \C$ and define 
$$
\widehat C_k := \text{\rm ev}_0^{-1}(p_0) \subset \mathcal
M_{k+1}(\C,S;\beta_D).
$$
We choose the $\frak q$-multisections on
$\mathcal M^{\text{\rm main}}_{k+1,\ell}(\beta;{\bf p})$ that
satisfies  Condition \ref{condmultisq}.
For given ${\bf p}$, consider
the perturbed moduli spaces
$\mathcal M^{\text{\rm main}}_{k+1,\ell}(\beta;{\bf p})^{\frak q}$.
Then we have
\begin{equation}\label{eq:Mk+1}
\mathcal M^{\text{\rm main}}_{k+1,\ell}(\beta;{\bf p})^{\frak q} \cong \mathcal
M^{\text{\rm main}}_{1,\ell}(\beta;{\bf p})^{\frak q} \times \widehat C_k.
\end{equation}
In fact, $\mathcal M^{\text{\rm main}}_{1,\ell}(\beta;{\bf p})^{\frak q}$ consists of finitely many
free $T^n$ orbits (with multiplicity $\in \Q$) and $\mathcal M^{\text{\rm main}}_{1}(\beta;{\bf p})^{\frak q}
= \mathcal M^{\text{\rm main,reg}}_{1}(\beta;{\bf p})^{\frak q}$. 
Hence \eqref{eq:Mk+1} holds after taking the $\frak q$-multisections.
By Condition \ref{condmultisq}.4 we have a map
$\mathcal M^{\text{\rm main,reg}}_{k + 1;\ell}(\beta;{\bf p})^{\frak q}
\to \mathcal M^{\text{\rm main,reg}}_{1;\ell}(\beta;{\bf p})^{\frak q}$.
It is easy to see that the fiber can be identified with $\widehat C_k$.
\par
Now we find that the iterated integral along a loop only depends on
the cohomology class $[\text{\rm ev}_1^* b_+]= \dots =[\text{\rm ev}_k^* b_+]$ using
 in \cite[Proposition 4.1.1 and (4.1.3), (4.1.4)]{Chen}. 
Then we find that for any general class $\beta \in H_2(X,L(u);\Z)$, 
the iterated integral on $\widehat{C}_k$ is 
calculated as 
\begin{equation}\label{iteratedintegral}
\int_{\widehat{C}_k} \text{\rm ev}_1^* b_+ \wedge \dots \wedge \text {\rm ev}_k^* b_+ = \frac{1}{k!} \langle \partial \beta, b_+ \rangle^k.  
\end{equation}
See also \cite[Lemma 7.2]{fooo09}.  
This equality is sufficient for the proof of 
\cite[Lemma 11.8]{fooo08}  or \eqref{eq:identityb} above. 
\end{rem}
\par
\section{Interpolation between Kuranishi structures}
\label{sec:equikuracot}

In this section we clarify the point mentioned in Remark \ref{rem:1014} and complete the proof of Lemma \ref{perturbnbdsigmasec9}. 
Namely we will construct the system of continuous families of multisections $\frak s$.
We first explain the problem in more detail.
We review the notion `component-wise Kuranishi structure' and  `disk component-wise Kuranishi structure'
introduced in \cite{fooo091} for this purpose.
\par
We consider the moduli space $\mathcal M_{k+1;\ell}^{\text{\rm main}}(\beta;\text{\bf p})$
introduced in Section \ref{sec:frakqreview}. 
Here $\text{\bf p} = (P_1,\dots,P_{\ell})$ and
$P_i \in \mathcal A(X)$ are $T^n$ invariant cycles. 
(We allow the case $P_i=P_j$ for some $i\ne j$.)
Let $\text{\bf p}_i = (P_1^i,\dots,P_{\ell_i}^i)$ with $P_j^i \in \mathcal A(X), i=1,\dots ,m$. 
We say that a decomposition 
$\text{\bf p}_1
\cup \dots \cup \text{\bf p}_m = \text{\bf p}$ 
is a {\it disjoint union} if 
$\{ P_1^1, \dots, P_{\ell_1}^1, P_1^2, \dots, P_{\ell_m}^m \} =\text{\bf p}$ 
(with multiplicity)
as an unordered set and 
$\sum_{i=1}^{m}\ell_i =\vert \text{\bf p} \vert$. 
(It may happen that $P_j^i=P_{j'}^{i'}$ in $\mathcal A(X)$.)
We consider the following two kinds of stratifications of $\mathcal M_{k+1;\ell}^{\text{\rm main}}(\beta;\text{\bf p})$.
\begin{enumerate}
\item The stratum of this stratification is a fiber product whose factors are
$$
\mathcal M_{k_i+1;\ell_i}^{\text{\rm main}}(\beta_{(i)};\text{\bf p}_i),
\qquad i=1,\dots,m
$$
where $\sum k_i = k+m$, $\sum \beta_{(i)} = \beta$, $\sum \ell_i = \ell$,
$\text{\bf p}_1
\cup \dots \cup \text{\bf p}_m = \text{\bf p}$ (disjoint union).
Elements of the interior of  this stratum have $m$ disk components.
We include the tree of sphere components to the disk component
on which it is rooted.
\item
 The stratum of this stratification is a fiber product whose factors are
$$
\mathcal M_{k_i+1;\ell_i}^{\text{\rm main}}(\beta_{(i)};\text{\bf p}_i),
\qquad i=1,\dots,m
$$
or
$$
\mathcal M_{\ell'_j}(\alpha_{j};\text{\bf p}_j^{\prime}),
\qquad j=1,\dots,m'
$$
where $\sum k_i = k+m$, $\sum \beta_{(i)} + \sum \alpha_{j} = \beta$, $\sum \ell_i
+ \sum \ell'_j= \ell+m'$,
$\text{\bf p}_1
\cup \dots \cup \text{\bf p}_m
\cup \text{\bf p}'_1
\cup \dots \cup \text{\bf p}'_{m'}= \text{\bf p}$ (disjoint union).
\end{enumerate}
We call the stratification 1  {\it disk stratification}
and 2 {\it disk-sphere stratification}.
See \cite[Subsection 7.1.1]{fooobook2} for the stratification 1.
\begin{defn}\label{exdiskcompt}
We call the union of one disk component and all the trees of sphere components 
rooted on that disk component an 
{\it extended disk component}\index{extended disk component}. 
\end{defn}
\begin{defn}\label{cpwkura}
\begin{enumerate}
\item We say a system of Kuranishi structures of $\mathcal M_{k+1;\ell}^{\text{\rm main}}(\beta;\text{\bf p})$ is
{\it disk component-wise}\index{Kuranishi structure!disk component-wise Kuranishi structure} 
if it is compatible with the fiber product description of
each of the strata of disk stratification.
\item 
The definition of disk component-wise-ness of 
multisection\index{multisection (perturbation)!disk component-wise multisection}  
on disk component-wise Kuranishi structure is similar.
\item We say a system of Kuranishi structures of $\mathcal M_{k+1;\ell}^{\text{\rm main}}(\beta;\text{\bf p})$
and of $\mathcal M_{\ell}(\alpha;\text{\bf p})$ is {\it component-wise}
\index{Kuranishi structure!component-wise Kuranishi structure} if it is compatible with the fiber product description of
each of the strata of disk-sphere stratification.
\item 
The definition of component-wise-ness of multisection 
on component-wise Kuranishi structure is similar.
  \index{multisection (perturbation)!component-wise multisection}  
\end{enumerate}
\end{defn}
We observe that the multisections of $\mathcal M_{k+1;\ell}^{\text{\rm main}}(\beta;\text{\bf p})$ which we use
in this paper and in \cite{fooo08,fooo09} are disk component-wise and $T^n$ equivariant.
In fact, Condition \ref{perturbforq}.3  implies that it is disk component-wise.
The detail of the construction of $T^n$ equivariant Kuranishi structure is given in \cite[Section B]{fooo08}.
(See also Remark \ref{rem457}.1.)
\par
On the other hand, the multisections of $\mathcal M_{k+1;\ell}^{\text{\rm main}}(\beta;\text{\bf p})$, that are the 
multisection $\frak q$ and the continuous family of multisections
$\frak s$, that
we use in this paper and in \cite{fooo08, fooo09} are not necessarily component-wise.
The reason is discussed in \cite[Remark 11.4]{fooo08}.
The difficulty to make our Kuranishi structure both  component-wise and $T^n$ equivariant
lies in the fact that the $T^n$ action for the moduli space of pseudo-holomorphic
{\it spheres} is not necessarily free though $T^n$ action for the moduli space of pseudo-holomorphic
{\it disks} with at least one boundary marked point is free.
\begin{rem}
In \cite{fooo08} we find $T^n$ equivariant multisection as follows.
We first consider the quotient by this $T^n$ action, then find a 
multisection on the quotient and finally lift it.
If the $T^n$ action has isotropy group of positive dimension, 
the quotient space is neither a manifold nor an orbifold.
So this method does not work.
\end{rem}
\par
In the rest of this section, we discuss how we go around this trouble and complete the proof of 
Lemma \ref{perturbnbdsigmasec9}.
\par
Let $\text{\bf p} = (P_1,\dots,P_{\ell})$ be as above.
We consider the fiber product
\begin{equation}
\mathcal M_{k+1;\ell+2}^{\text{\rm main}}(\beta;\text{\bf p}) = 
\mathcal M_{k+1;\ell+2}^{\text{\rm main}}(\beta)
{}_{{\text{\rm ev}}_3^{\rm int},\dots,{\text{\rm ev}_{\ell+2}^{\rm int}}} \times_{X^{\ell}} (P_1 \times \dots \times P_{\ell}).
\end{equation}
The evaluation maps at the boundary marked points give a map
\begin{equation}
{\text{\rm ev}} = ({\text{\rm ev}}_0,\dots,{\text{\rm ev}}_k) :
\mathcal M_{k+1;\ell+2}^{\text{\rm main}}(\beta;\text{\bf p})
\to L^{k+1}
\end{equation}
and the evaluation maps at the 1st and 2nd interior marked points give:
\begin{equation}
{\text{\rm ev}}_{\rm int} = ({\text{\rm ev}}^{\rm int}_1, {\text{\rm ev}}^{\rm int}_2) : \mathcal M_{k+1;\ell+2}^{\text{\rm main}}(\beta;\text{\bf p})
\to X^2.
\end{equation}
We also consider the forgetful map
\begin{equation}
\mathfrak{forget} : \mathcal M_{k+1;\ell+2}^{\text{\rm main}}(\beta;\text{\bf p})
\to \mathcal M_{1;2}^{\text{\rm main}}
\end{equation}
as in (\ref{forgetmap}).
We consider the point $[\Sigma_0] \in \mathcal M_{1;2}^{\text{\rm main}}$
as in Lemma \ref{stratifyD12} and take its small neighborhood $U_0$ in
$\mathcal M_{1;2}^{\text{\rm main}}$.
We consider
$$
\mathcal M_{k+1;\ell+2}^{\text{\rm main}}(\beta;\text{\bf p})
\cap \mathfrak{forget}^{-1}(U_0),
$$
which we write
$$
\mathcal M_{k+1;\ell+2}^{\text{\rm main}}(\beta;\text{\bf p};U_0)
$$
for simplicity.
We can define the notion of disk component-wise
(resp. component-wise) system of Kuranishi structures and multisections on it
by modifying Definition \ref{cpwkura} in an obvious way.
\par
Its codimension one stratum is a union of fiber product (\ref{fipro1})
or (\ref{fipro2})
\begin{equation}\label{fipro1}
\mathcal M_{k_1+1;\ell_1}^{\text{\rm main}}(\beta_1;\text{\bf p}_1)
{}_{\text{\rm ev}_0} \times_{\text{\rm ev}_i}
\mathcal M_{k_2+1;\ell_2+2}^{\text{\rm main}}(\beta;\text{\bf p}_2;U_0)
\end{equation}
where $k_1+k_2 = k+1$, $\ell_1+\ell_2=\ell$,
$\text{\bf p}_1 \cup \text{\bf p}_2 = \text{\bf p}$ (disjoint union).
\begin{equation}\label{fipro2}
\mathcal M_{k_1+1;\ell_1+2}^{\text{\rm main}}(\beta_1;\text{\bf p}_1;U_0)
{}_{\text{\rm ev}_0} \times_{\text{\rm ev}_i}
\mathcal M_{k_2+1;\ell_2}^{\text{\rm main}}(\beta;\text{\bf p}_2)
\end{equation}
where $k_1+k_2 = k+1$, $\ell_1+\ell_2=\ell$,
$\text{\bf p}_1 \cup \text{\bf p}_2 = \text{\bf p}$ (disjoint union).
\par
We recall that during the construction of the operator $\frak q$,
we use a $T^n$ equivariant, disk component-wise
Kuranishi structure and multisection on
$\mathcal M_{k+1;\ell+2}^{\text{\rm main}}(\beta;\text{\bf p})$.
(However, it is not component-wise). We call it {\it $\frak q$-Kuranishi structure}\index{Kuranishi structure!$\frak q$-Kuranishi structure}  
and {\it $\frak q$-multisection}\index{multisection (perturbation)!$\frak q$-multisection}.
\par
The next step is to construct a continuous family of multisections on
$\mathcal M_{k+1;\ell+2}^{\text{\rm main}}(\beta;\text{\bf p};U_0)$
that is compatible
with the fiber product description (\ref{fipro1}), (\ref{fipro2})
of the codimension one stratum. Here we use $\frak q$-multisection
on the first factor of (\ref{fipro1}) and the
second factor of (\ref{fipro2}).
Moreover our Kuranishi structure and multisection are {\it partially component-wise}\index{Kuranishi structure!partially component-wise Kuranishi structure}\index{multisection (perturbation)!partially component-wise multisection}.
We now define this notion precisely.
\par
Let $(\Sigma,\vec z,\vec z^{+},u,(x_1,\dots,x_{\ell}))$
be an element of
$\mathcal M_{k+1;\ell+2}^{\text{\rm main}}(\beta;\text{\bf p};U_0)$.
(Namely $(\Sigma,\vec z,\vec z^{+},u) \in \mathcal M_{k+1;\ell+2}^{\text{\rm main}}(\beta)$
and $x_i \in P_i$.)
We consider a tree of sphere components $\Sigma_a$ of $\Sigma$,
which is rooted at a point $w^a$ contained in a disk component.
We denote the corresponding disk by
$D_{i(a)}$. We define $z^a \in \partial D_{i(a)}$ as follows.
\begin{enumerate}
\item If the $0$-th boundary marked point $z_0$ is on $D_{i(a)}$, then
$z^a = z_0$.
\item Otherwise there exists a unique singular point of $D_{i(a)}$ such that
every path joining $D_{i(a)}$ with $z_0$ contains it. We put this singular point
to be $z^a$.
\end{enumerate}
We next define points $w_i$, $i=1,2$ as follows:
\begin{enumerate}
\item If the $i$-interior marked point $z^+_i$ is on a disk component,
then we put $w_i = z^+_i$.
\item Otherwise, there is a tree of sphere components containing $z^+_i$ and
is rooted on a disk component. We set $w_i$ to be the point on the disk
component where this tree of sphere components is rooted.
\end{enumerate}
\begin{defn}\label{def4232}
Let $\epsilon$ be a sufficiently small positive number,
which we determine later and $\Sigma_a$ a maximal tree of sphere components.
\par
We say $\Sigma_a$ is of {\it Type I} if one of the
following two conditions holds:
\begin{enumerate}
\item $w_1 \notin D_{i(a)}$.
\item Suppose $w_1 \in D_{i(a)}$.
We require $w_1 \ne w^a$.
Then
$(D_{i(a)},z^a,(w_1,w^a))$ defines an element of $\mathcal M_{1;2}^{\text{\rm main} }
\cong D^2$. We then require
\begin{equation}
dist([D_{i(a)},z_a,(w_1,w^a)],\partial\mathcal M_{1;2}^{\text{\rm main} }) < \epsilon.
\end{equation}
\end{enumerate}
\par
We say that $\Sigma_a$ is  of {\it Type II} if both of the
following two conditions hold:
\begin{enumerate}
\item $w_1 \in D_{i(a)}$.
\item  Either $w_1 = w^a$ or
\begin{equation}
dist([D_{i(a)},z^a,(w_1,w^a)],[\Sigma_0]) < \epsilon.
\end{equation}
\end{enumerate}
\par
If  $\Sigma_a$ is neither of Type I nor of Type II, we say it is of {\it Type III}.
\par
We say a sphere component is of Type I (resp. Type II, Type III) if it is
contained in a maximal tree of sphere components of Type I (resp. Type II, Type III). 
See Figures 4.2.1 and 4.2.2.
\end{defn}
\par
\hskip-1.2cm
\epsfbox{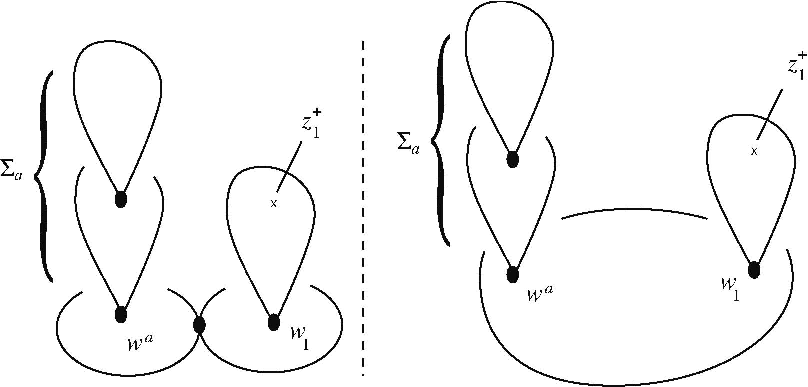}
\par
\centerline{\bf Figure 4.2.1 (Type I)}
\par
\hskip-1.2cm
\epsfbox{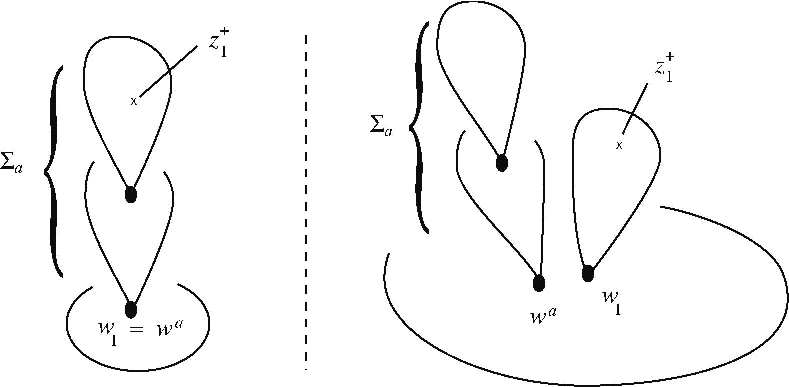}
\par
\centerline{\bf Figure 4.2.2 (Type II)}
\par
\begin{lem}\label{epsmallII}
Let a sufficiently small constant $\epsilon > 0$ be given. Then
we may choose $U_0$ sufficiently small so that
whenever $(\Sigma,\vec z,\vec z^{+},u,(x_1,\dots,x_{\ell}))$
lies in $\mathcal M^{\text{\rm main}}_{k+1;\ell+2}(\beta;\text{\bf p};U_0)$ and
$\Sigma_a$ contains one of the first two interior marked points $z^+_1,  z^+_2$,
then $\Sigma_a$ is of Type II.
\end{lem}
\begin{proof}
If $\Sigma_a$ contains $z_1^+$, then $w_1 \in D_{i(a)}$ and
$w^a = w_1$.
The lemma follows
in this case.
\par
In case $\Sigma_a$ contains $z_2^+$, then $w_2 \in D_{i(a)}$ also.
Moreover by the assumption
$(\Sigma,\vec z,\vec z^{+},u,(x_1,\dots,x_{\ell}))
\in \mathcal M^{\text{\rm main}}_{k+1;\ell+2}(\beta;\text{\bf p};U_0)$, we have
either $[D_{i(a)},z^a,w_1,w^a] = [D_{i(a)},z^a,w_1,w_2] \in U_0$ or $w_1 = w_2$.
Hence we may take $U_0$ small so that $\Sigma_a$ is of Type II.
\end{proof}
\begin{defn}\label{defcompwise}
A system of Kuranishi structures on $\mathcal M_{k+1;\ell+2}^{\text{\rm main}}(\beta;\text{\bf p};U_0)$
is said to be {\it partially component-wise}\index{Kuranishi structure!partially component-wise Kuranishi structure} if the following holds:
\par
We consider the stratification of $\mathcal M_{k+1;\ell+2}^{\text{\rm main}}(\beta;\text{\bf p};U_0)$
such that each of its stratum is a fiber product of
\begin{equation}\label{diskgeneral}
\mathcal M_{k_i+1;\ell_i}^{\text{\rm main}}(\beta_{(i)};\text{\bf p}_i),
\qquad i=1,\dots,m,
\end{equation}
or
\begin{equation}\label{spherreq}
\mathcal M_{\ell'_j}(\alpha_j;\text{\bf p}_j^{\prime}),
\qquad j=1,\dots,m'
\end{equation}
or
\begin{equation}\label{disk0comp}
\mathcal M_{k_0+1;\ell_0+2}^{\text{\rm main}}(\beta_{(0)};\text{\bf p}_0;U_0)
\end{equation}
where we have $\sum k_i = k+m$, $\sum \beta_{(i)} + \sum \alpha_j = \beta$, $\sum \ell_i
+ \sum \ell'_j= \ell+m'$, and
$$
\text{\bf p}_0
\cup \dots \cup \text{\bf p}_m
\cup \text{\bf p}'_1
\cup \dots \cup \text{\bf p}'_{m'}= \text{\bf p},
$$
a disjoint union. (Note that $\beta_{(0)}$ does not stand for $\beta_0 =0$.)
For the factor (\ref{spherreq}), we require the corresponding sphere components
are of type II.
(In other words, we include a tree of sphere components of Type I or III to the
disk component on which it is rooted.)
\par
Then the Kuranishi structure is compatible with this fiber product description of the
stratum of this stratification. Here we use the $\frak q$-Kuranishi structure for the factor (\ref{diskgeneral}),
and the Kuranishi structure in this definition for the factor (\ref{disk0comp}).
\par
The definition of partially component-wise multisection is similar.
\end{defn}

\begin{lem}\label{pcwkura}
There exist a system of Kuranishi structures and 
a continuous family of multisections on $\mathcal M_{k+1;\ell+2}^{\text{\rm main}}(\beta;\text{\bf p};U_0)$
with the following properties:
\begin{enumerate}
\item
They are compatible with the forgetful map
\begin{equation}\label{426forgetcomp}
\mathcal M_{k+1;\ell+2}^{\text{\rm main}}(\beta;\text{\bf p};U_0)
\to \mathcal M_{1;\ell+2}^{\text{\rm main}}(\beta;\text{\bf p};U_0).
\end{equation}
\item
They are invariant under the action of symmetric group of order $\ell!$ which
exchanges the 3rd - $(\ell+2)$nd interior marked points and
exchanges $P_1,\dots,P_{\ell}$ at the same time.
\item
They are compatible with the fiber product description $(\ref{fipro1})$,
$(\ref{fipro2})$ of the codimension one stratum of
$\mathcal M_{k+1;\ell+2}^{\text{\rm main}}(\beta;\text{\bf p};U_0)$.
Here we use the $\frak q$-Kuranishi structure for the
first factor of $(\ref{fipro1})$ and the second factor of $(\ref{fipro2})$.
\item They are partially component-wise in the sense of Definition \ref{defcompwise}.
\item 
The continuous family of multisections is 
transversal to zero. The restriction of the evaluation map $\text{\rm ev}_0$ to 
zero set of multisection  is 
a submersion.
\end{enumerate}
\end{lem}
We remark that if a sphere component of an element of
a codimension one stratum is on the first factor of
$(\ref{fipro1})$ or the second factor of $(\ref{fipro2})$
then this sphere component is of type I.
So we do not require the component-wise-ness
for this sphere component.
Once this point is understood,
we can use the induction on $\beta\cap\omega$  and $k$,
$\ell$.
The proof of Lemma \ref{pcwkura} is in Subsection \ref{subsec:lemma426}.

\begin{proof}[Proof of Lemma \ref{perturbnbdsigmasec9}]
We check that the Kuranishi structure and the family of multisections in Lemma \ref{pcwkura}
have the properties required in Lemma \ref{perturbnbdsigmasec9}.
\par
Lemma \ref{perturbnbdsigmasec9}.1 and 7 are  consequences of Lemma \ref{pcwkura}.5.
\par
Lemma \ref{perturbnbdsigmasec9}.2 and 6 follow from Lemma \ref{pcwkura}.4.
In fact, if $[\Sigma,\vec z,u]$ lies in $\frak{forget}^{-1}([\Sigma_0])$
then there is a unique disk component that contains the tree of sphere components 
containing 1st and 2nd interior marked points. And that disk component is of Type II.
So Lemma \ref{pcwkura}.4 implies the required properties.
\par
Lemma \ref{perturbnbdsigmasec9}.3 is a consequence of Lemma \ref{pcwkura}.3.
\par
Lemma \ref{perturbnbdsigmasec9}.4 is a consequence of Lemma \ref{pcwkura}.1.
\par
Lemma \ref{perturbnbdsigmasec9}.5 is a consequence of Lemma \ref{pcwkura}.2.
\par
The proof of Lemma \ref{perturbnbdsigmasec9} is complete,
assuming Lemma \ref{pcwkura}.
\end{proof}

\section{$T^n$ equivariant and
cyclically symmetric Kuranishi structures}
\label{sec:cyclicKura}
In \cite{fooo08} we gave a construction of a system of $T^n$-equivariant multisections
of moduli space of pseudo-holomorphic disks.
In \cite{fooo09} we gave a construction of a continuous family of multisections
of moduli space of pseudo-holomorphic disks.
In \cite{fooo091} a construction of
a continuous family of multisections with cyclic symmetry is explained.
In this paper we use a continuous family of 
$T^n$-equivariant and cyclically symmetric
multisections that are disk-component-wise.
We describe the construction in this and the next sections.
In this section, we construct a system of $T^n$-equivariant and cyclically symmetric
Kuranishi structures ($\frak c$-Kuranishi structures).
\index{Kuranishi structure!$\frak c$-Kuranishi structure}
Namely we prove the following:
\begin{prop}\label{transconcl}
There exists a system of Kuranishi structures on $\mathcal M_{k+1;\ell}^{\rm main}(\beta)$
($k\ge 0, \ell \ge 0$) with the following properties.
\par
\begin{enumerate}
\item It is compatible with the map
$$
\frak{forget} : \mathcal M_{k+1;\ell}^{\rm main}(\beta) \to \mathcal M_{1;\ell}^{\rm main}(\beta)
$$
forgetting the 1-st \dots $k$-th boundary marked points. (Namely only the $0$-th boundary
marked point remains.) (See \cite[Section 3]{fooo091}  for the precise meaning of
this compatibility.)
\item
It is invariant under the cyclic permutation of the boundary marked
points.
\item It is invariant under the permutation of the interior marked points.
\item
$$
\text{\rm ev}_0 : \mathcal M_{k+1;\ell}^{\rm main}(\beta) \to L
$$
is weakly submersive.
\item
The Kuranishi structures are disk-component-wise.
\item
The Kuranishi structures are $T^n$-equivariant.
\end{enumerate}
\end{prop}
\begin{proof}
The proof occupies the rest of this section.
\par
Let $ (\Sigma,\vec z,\vec z^+)$ be a genus $0$ bordered marked semi-stable curve with 
$\partial\Sigma = S^1$ and
$u : (\Sigma,\partial \Sigma) \to (X,L)$ a pseudo-holomorphic map that is stable.
Here $\vec z = (z_1,\dots,z_k), \vec z^+ = (z^+_1,\dots,z^+_{\ell})$ are boundary and interior marked points, respectively. 
Hereafter we put
$${\bf x}=(\Sigma, \vec z,\vec z^+,u).
$$
\begin{defn}\label{symmetryG}
We define by $G(\text{\bf x})$ the group consisting of the pairs $(g,\varphi)$
of $g \in T^n$ and a biholomorphic map $\varphi : \Sigma \to \Sigma$
such that:
\begin{enumerate}
\item $u\circ\varphi = gu$.
\item $\varphi$ preserves  $\vec z^+$ as a set.
\item $\varphi$ fixes all the points of $\vec z$.
\end{enumerate}
We define a group structure on $G(\text{\bf x})$ by
\begin{equation}\label{Gxgroup}
(g_1,\varphi_1) \cdot (g_2,\varphi_2)
= (g_1g_2,\varphi_1\circ \varphi_2).
\end{equation}
\end{defn}
Here we allow $\varphi$ to permute the interior marked points $\vec z^+$ because
we want our Kuranishi structure to be invariant under the
permutation of the interior marked points. By definition, we have the projections:
$$
\xymatrix{ {} & G(\text{\bf x}) \ar[dl]_{\pi_{T^n}} \ar[dr]^{\pi_{\text{\rm Aut}}} &{}\\
T^n & & \text{\rm Aut}(\Sigma)}
$$
\begin{lem} The natural projection $G(\text{\bf x}) \to \Aut(\Sigma)$
is injective.
\end{lem}
\begin{proof} Let $(g,\varphi) \in G(\text{\bf x})$ with $\varphi = id$.
By Definition \ref{symmetryG}.1, we have
$g u \equiv u$. Since $u$ has its boundary lying on a Lagrangian torus orbit
$L$ on which $T^n$ acts freely, it follows that $g = id$.
\end{proof}
\begin{lem}\label{isotropycompact}
$G(\text{\bf x})$ is compact.
\end{lem}
\begin{proof}
We decompose $\Sigma$ into the union of irreducible components 
$\Sigma_a$, $a \in A$.  
There exists an index finite subgroup $G(\text{\bf x})_0$ of $G(\text{\bf x})$
such that elements of $G(\text{\bf x})_0$ preserve each of $\Sigma_a$.
It suffices to prove that $G(\text{\bf x})_0 \to T^n$ is proper.
\par
We prove it by contradiction.
Suppose that there exists $(g_i,\varphi_i) \in G(\text{\bf x})_0$ 
such that $g_i$ converges to $g_{\infty} \in T^n$ but $\varphi_i$ diverges 
in $\text{\rm Aut}(\Sigma)$.
\par
Since $\varphi_i$ diverges, by taking a subsequence if necessary, 
we may assume that there exists $\Sigma_a$ such that 
$$
\lim_{i\to \infty} \varphi_i(z) = z_{\infty}
$$
holds for $z \in \Sigma_a\setminus \text{\rm a finite set}$. Here $z_{\infty}$ is independent of $z$.
\par
Since $u\varphi_i(z) = g_iu(z)$, we have
$$
u(z) = \lim_{i\to \infty} g_i^{-1}u\varphi_i(z) = g_{\infty}^{-1}u(z_{\infty}).
$$
Namely, $u$ is constant on $\Sigma_a$. Since ${\rm Aut}(\Sigma_a)$ is of infinite 
order, this contradicts  the stability.
\end{proof}
We note that each element of $G(\text{\bf x})$ induces a permutation of
the elements of $\vec z^+$, if we regard $\vec z^+ = (z^+_i)_{i=1,\ldots,\ell}$
as an ordered set.
\par
We use it to define a $G(\text{\bf x})$-action on the set of $\text{\bf x}'
= (\Sigma',\vec z',\vec z^{\prime +},u')$ by
\begin{equation}\label{Gzactiondef}
(g,\varphi) \cdot  (\Sigma',\vec z',\vec z^{\prime +},u')
=
 (\Sigma',\vec z',\rho\vec z^{\prime +},gu')
\end{equation}
where $\rho$ is defined by $\varphi(z_a^+) = z^+_{\rho(a)}$ and
$\rho\vec z^{\prime +}$ is defined by
$(\rho\vec z^{\prime +})_a = z^{\prime +}_{\rho^{-1}(a)}$.
\par
This is a left action. Namely
$$
((g_1,\varphi_1) \cdot (g_2,\varphi_2))\cdot
(\Sigma',\vec z',\vec z^{\prime +},u')
= (g_1,\varphi_1) \cdot ((g_2,\varphi_2)\cdot
(\Sigma',\vec z',\vec z^{\prime +},u')).
$$
We also note that  
$$
(g,\varphi) \cdot \text{\bf x} = \text{\bf x}
$$
by definition.
\par
We also define an action of $T^n \times \frak S_{\ell}$ on the
set of $\text{\bf x}' = (\Sigma',\vec z',\vec z^{\prime +},u')$ by
$$
(g,\rho)\text{\bf x}' = (\Sigma',\vec z',\rho\vec z^{\prime +},gu').
$$
The isotropy group of $\text{\bf x}$ of this action can be identified with
$G(\text{\bf x})$.
\par
We begin with the following lemma. We fix a $T^n$ invariant metric on $X$.
\begin{lem}\label{app2lema}
Let $\text{\bf x} = (\Sigma,\vec z,\vec z^+,u)$  and $G(\text{\bf x})$ as above.
We consider a relatively compact $G(\text{\bf x})$ invariant open subset
$K \subset \Sigma \setminus \vec z^+$  whose closure does not contain
singular points and does not intersect with $\partial \Sigma$. Then there exists
a subspace
$$
E \subset C^{\infty}(\Sigma;u^*TX\otimes \Lambda^{0,1})
$$
such that
\begin{enumerate}
\item
$E$ is $G(\text{\bf x})$-invariant.
\item
Elements of $E$ have compact support in $K$.
\item
$E$ is finite dimensional.
\item
$
D_u\overline{\partial}\left(
\{ \xi \in W^{1,p}(\Sigma,u^*TX)\mid \text{$\xi = 0$ at $\vec z \cup \vec z^+$} \}\right) + E =
L^p(\Sigma;u^*TX\otimes \Lambda^{0,1}).
$
\end{enumerate}
\end{lem}
\begin{proof}
The existence of $E_0$ that satisfies properties 2,3,4 above can be proved by the unique continuation
theorem in a standard way. 
(We also use compactness of the group $G({\bf x})$ (Lemma \ref{isotropycompact}) 
to prove the existence of $K$.)
We will modify the construction of such $E_0$ so that it also 
satisfies property 1.
We define an action of $G(\text{\bf x})$ on
$L^2 (K;u^*TX\otimes \Lambda^{0,1})$ which extends the action
$$
\eta \mapsto (g,\varphi)_*(\eta) = (dg \otimes (d\varphi^{-1})^*)(\eta).
$$
More specifically as a one-form on $\Sigma$, $(g,\varphi)_*(\eta)$ is defined
by
$$
(g,\varphi)_*(\eta)(v_x) = dg(\eta(d\varphi^{-1}(v_x)), \quad v_x \in T_x\Sigma
$$
for a smooth $u^*TX$-valued one-form $\eta \in C^\infty(u^*TX \otimes \Lambda^{0,1})$.
We recall that by definition of $G(\text{\bf x})$ we have the commutative diagram
$$
\xymatrix{ u^*TX \otimes \Lambda^{0,1} \ar[d]\ar[r]^{(g,\varphi)_*} &  u^*TX \otimes \Lambda^{0,1} \ar[d] \\
\Sigma \ar[r]_\varphi & \Sigma}
$$
(Here we recall $g u \circ \varphi^{-1} = u$ from the definition of $G(\text{\bf x})$.)
Using the compactness of $G(\text{\bf x})$, we can
equip the space $L^2(K;u^*TX\otimes \Lambda^{0,1})$
with a $G(\text{\bf x})$-invariant inner product. 
Then we have a decomposition
$$
L^2(K;u^*TX\otimes \Lambda^{0,1})
= \widehat{\bigoplus_{\lambda \in \frak L}} V_{\lambda}
$$
where $V_{\lambda}$ is a finite dimensional representation of $G(\text{\bf x})$
and $\frak L$ is a certain index set, which is countable.
(Here we use compactness of $G(\text{\bf x})$. (Lemma \ref{isotropycompact}.))
We denote by $\widehat{\bigoplus}$ the $L^2$ completion of the direct sum.
Moreover we may choose $V_{\lambda}$ so that each element of $V_{\lambda}$ is smooth.
In fact, we take $G(\text{\bf x})$ invariant metrics of $\Sigma$ and $u^*TX$
and use them to define Laplace operator on
$L^2(K;u^*TX\otimes \Lambda^{0,1})$.
Then each of its eigenspaces (satisfying Dirichlet boundary condition, for example) is a finite dimenisonal vector space consisting of
smooth sections. 
Moreover it is $G(\text{\bf x})$ invariant. 
(Recall that $K$ is $G({\bf x})$ invariant.) 
Therefore we may take $V_{\lambda}$ as
the eigenspaces of this Laplace operator.
\par
Let $e_i$, $i=1,\dots,N$ be an orthonormal basis of $E_0$.
We put
\begin{equation}\label{decompf}
e_i = \sum_{\lambda=1}^{\infty} e_{i,\lambda}
\end{equation}
where $e_{i,\lambda} \in V_{\lambda}$. 
Here we can choose $e_i$'s so that the series 
\eqref{decompf}
converges uniformly and so does in $L^p$.
(This is well known for the eigenspaces of Laplace operator.)
We put
$$
e_{i,k} =  \sum_{\lambda=1}^k e_{i,\lambda}
$$
and let $E_k$ be the vector space generated by $e_{1,k},\ldots,e_{N,k}$.
If $k$ is large, $E_k$ satisfies properties 3 and 4.
(Property 4 can be proved by, for example, open mapping theorem.)
We put
$$
E' = \bigoplus_{\lambda=1}^k V_{\lambda},
$$
which satisfies properties 1, 3, 4.
\par
To promote $E'$ to be a subspace that also satisfies property 2, we proceed as follows.
Let $\chi_i : \Sigma \to [0,1]$ be a sequence of
$G(\text{\bf x})$ invariant
smooth functions with compact support in $K$
such that $\chi_i$ converges to $1$ pointwise
everywhere in $K$.
We put
$$
E(i) = \{ \chi_i s \mid s \in E' \}.
$$
We put $E = E(i)$ for sufficiently large $i$. Then property 4 holds.
Properties 1,2,3 are immediate.
\end{proof}
The main part of the inductive construction of 
$T^n$-equivariant and
cyclically symmetric Kuranishi structure on 
$\mathcal M_{k+1;\ell}(\beta)$ 
is the proof of the following proposition.
During the proof we will use Lemma \ref{app2lema}.
\begin{prop}\label{eqcyckuramain}
There exist systems of Kuranishi structures on $\mathcal M_{0;\ell}(\beta)$
and $\mathcal M_{1;\ell}(\beta)$ with the following properties:
\begin{enumerate}
\item They are compatible with forgetful maps in the sense of \cite[Definition 3.1]{fooo091}.
\item The evaluation map $\text{\rm ev}_0
: \mathcal M_{1;\ell}(\beta) \to L(u)$ at the (unique) boundary marked point
is weakly submersive.
\item The codimension one stratum of $\mathcal M_{0;\ell}(\beta)$
is a union of
$$
\mathcal M_{1;\ell_1}(\beta_1) {}_{\rm ev_0} \times_{\rm ev_0}
\mathcal M_{1;\ell_2}(\beta_2)
$$
with $\ell_1+\ell_2=\ell$ and $\beta_1+\beta_2 = \beta$,
$(\ell_1,\beta_1) \ne (\ell_2,\beta_2)$, 
$$
(\mathcal M_{1;\ell'}(\beta') {}_{\rm ev_0} \times_{\rm ev_0}
\mathcal M_{1;\ell'}(\beta'))/\Z_2
$$
with $2\ell' = \ell$, $2\beta'=\beta$,
(Here $\Z_2$ acts by exchanging the factors.)
and of
\begin{equation}\label{diskspherebubble}
\mathcal M_{\ell+1}(\widetilde\beta) \,{}_{{\rm ev}_0}\times_X L.
\end{equation}
(We explain this boundary component in Remark \ref{pbubblerem} below.)
\par
Our systems of Kuranishi structures are compatible with this fiber product
description.
\item  Our Kuranishi structures are $T^n$ equivariant
in the sense of \cite[Definition B.4]{fooo09}.
\item Our Kuranishi structures are invariant under the permutation of
interior marked points.
\end{enumerate}
\end{prop}
\begin{rem}\label{pbubblerem}
Here is the description of the boundary component (\ref{diskspherebubble}).
We denote by $\widetilde\beta$  an element of $H_2(X;\Z)$
such that its image in $H_2(X,L;\Z)$ is $\beta$.
The space 
$\mathcal M_{\ell+1}(\widetilde\beta)$ is a moduli space of pseudo-holomorphic 
sphere with $\ell+1$ marked points and ${\rm ev}_0$ is the evaluation map
at the $0$-th marked point. We take the fiber product with $L$ 
by this map. 
Such a boundary corresponds to the following configuration: 
We take a constant map $u_0 : D^2 \to L$. 
Let $(\Sigma^{\rm sph},z^+_0\cup \vec z^+,u_1)$ be an element of  
$\mathcal M_{\ell+1}(\widetilde\beta)$. We 
assume $u_1(z^+_0)$ is the point $u_0(D^2)$.
We glue $D^2 = \{x\in \C \mid \vert x\vert \le 1\}$ and $\Sigma^{\rm sph}$ 
at $0\in D^2$ and $z_0^+$ to obtain $\Sigma$.
We define $u : \Sigma \to X$ so that it is $u_0$ on $D^2$ and 
is $u_1$ on $\Sigma^{\rm sph}$.
Such an object $(\Sigma,\vec z^+,u)$ appears as an element of the 
compactification $\mathcal M_{0;\ell}(\beta)$.
(See \cite[Subsection 3.8.3]{fooobook}, \cite[Subsection 7.4.1]{fooobook2}, \cite{Liu}.)
\par
Including $(\Sigma,\vec z^+,u)$ in the compactified moduli space $\mathcal M_{0;\ell}(\beta)$ is slight abuse of 
notation, since the automorphism of $(\Sigma,\vec z^+,u)$ 
contains $S^1$ so it is not a stable map.
(This group $S^1$ is the automorphism of $(D^2,0,u_0)$ which rotates 
$D^2$ around the center.)
We include such an object $(\Sigma,\vec z^+,u)$ into 
the compactified moduli space $\mathcal M_{0;\ell}(\beta)$ and 
point out the appearance of such an object whenever it appears, 
and will explain the point 
we need to modify the necessary arguments.
\end{rem}
The proof of Proposition \ref{eqcyckuramain} occupies most of the 
rest of this section.
The strategy of the proof of Proposition \ref{eqcyckuramain}
is similar to the proof of \cite[Theorem 3.1]{fooo091}.
However we need 
to be careful for the local construction of the
Kuranishi neighborhood since we need to perform the
whole construction in a $T^n$ equivariant way.
We will discuss this point in detail below.
\par
Let $\text{\bf x}$ be as in the proof of Proposition 
\ref{transconcl} and 
$G= G(\text{\bf x})$ as in Lemma \ref{app2lema}.
The following is a $T^n$ equivariant version of
\cite[Lemma 3.1]{fooo091}.
\begin{lem}\label{lem31Fu2ana}
There exists a finite dimensional vector space $E = E_{\text{\bf x}}$ as in Lemma \ref{app2lema}
with the following properties:
\begin{enumerate}
\item Lemma \ref{app2lema}.1-4 hold.
\item We put
\begin{equation}
\text{\rm Ker}(\text{\bf x}) = (D_u\overline{\partial})^{-1}(E_{\text{\bf x}}).
\end{equation}
Then for {\bf any} $z_0 \in \partial \Sigma$, the map
$
Ev_{z_0} :  \text{\rm Ker}(\text{\bf x}) \to T_{u(z_0)}L
$
defined by
\begin{equation}\label{evalatz0}
Ev_{z_0}(v) = v(z_0)
\end{equation}
is surjective.
\end{enumerate}
\end{lem}
\begin{proof}
We apply Lemma \ref{app2lema} to obtain $E_{\text{\bf x}}$
satisfying Item 1. Then Item 2 is an immediate consequence of the fact that
$\text{\rm Ker}(\text{\bf x})$ contains the Lie algebra of $T^n$.
\end{proof}
\begin{rem}
Thus the proof of Lemma $\ref{lem31Fu2ana}$ is easier than that of
\cite[Lemma 3.1]{fooo091}. Namely in our case 2 is automatic 
due to the presence of the free action of $T^n$ in the current circumstance.
\end{rem}
For $g\in T^n$ we put
\begin{equation}
E(g\text{\bf x}) := g_*E_{\text{\bf x}}.
\end{equation}
Here $g\text{\bf x} = (\Sigma,\vec{z}, \vec z^+,gu)$ and
$$
g_* : C^{\infty}(\Sigma,u^*TX \otimes \Lambda^{0,1})
\cong C^{\infty}(\Sigma,(gu)^*TX \otimes \Lambda^{0,1})
$$
is an obvious map.  
When $g$ is the identity, we have $E(\text{\bf x})=E_{\text{\bf x}}$.  
We consider the $T^n \times \frak S_{\ell}$ orbit
$(T^n \times \frak S_{\ell}) \text{\bf x}$ of $\text{\bf x}$.
For each element $(g,\rho)\text{\bf x}$ of it we put
$$
E((g,\rho)\text{\bf x})
= g_*E(\text{\bf x})
\subset C^{\infty}(\Sigma,(gu)^*TX \otimes \Lambda^{0,1}).
$$
\begin{lem}\label{lem:Einvbyrhog}
Suppose that $(g,\rho)\text{\bf x}$ and $\text{\bf x}$ define the 
same element in
$\mathcal M_{k;\ell}^{\text{\rm main}}(\beta)$, that is, 
there exists a biholomorphic map
$\varphi : (\Sigma,\vec z) \to (\Sigma,\vec z) $ such that
$\varphi(z^+_i) = z^+_{\rho^{-1}(i)}$ and
$gu = u\circ \varphi$. 
Then we have 
\begin{equation}\label{Einvbyrhog}
E((g,\rho)\text{\bf x}) = \varphi_*E(\text{\bf x}).
\end{equation}
\end{lem}

\begin{proof} 
This is an immediate consequence of Lemma \ref{app2lema}.
\end{proof}
\par
We will use this choice $E(g\text{\bf x})$ to define a
$T^n$ equivariant Kuranishi neighborhood of the $T^n$ orbit of
$\text{\bf x}$.
(This is the main part of the proof of Proposition \ref{eqcyckuramain}
and occupies the most of this section.)
\par
\subsection{The case of one disk component I}
\label{subsec:onedisk1}
We first note that it suffices to consider the case when $\beta \ne 0$.
Indeed, in the case of $\mathcal M_{k;0}^{\text {main}}(0)$  where $\beta =0$ and no interior 
marked points, 
the moduli space $\mathcal M_{k;0}^{\text {main}}(0)$ is transversal and the evaluation map at one boundary  marked point
is a submersion. So we can put our obstruction bundle $E$ to be trivial, i.e., $E=\mathcal M_{k;0}^{\text {main}}(0) \times \{0\}$.  
(This is not the case when we are studying the $\frak p$-perturbation.
We will discuss it in Subsection \ref{subsec:pkura}.)
\par
In Subsections \ref{subsec:onedisk1}, \ref{subsec:onedisk2}, \ref{subsec:onedisk3},
we study the case when $\Sigma$ has only one disk component. 
We will construct a Kuranishi neighborhood of $T^n \times \frak S_{\ell}$
orbit of ${\bf x} = (\Sigma, \vec z, \vec{z}^+,u)$ which is $T^n$ equivariant and invariant under the permutation of the interior marked points. 
In particular, we will describe the construction of obstruction bundles 
in detail.  
\par
In this subsection  
we start with the case when the following assumption is satisfied.
\begin{assum}\label{trivialisotopystrong}
For each $\frak t \in Lie(T^n)$ there exists a
(disk or sphere) component $\Sigma_a$ of $\Sigma$ and
$p \in \Sigma_a$ such that $u$ is an immersion at $p$ but $\frak t(u(p))$ is not tangent
to $u(\Sigma)$ at $p$. Here we also denote by $\frak t$ the Killing vector field
generated by $\frak t$.
\end{assum}
\par
We observe that Assumption \ref{trivialisotopystrong} implies that
the group $G(\text{\bf x})$ is finite.
\par
We emphasize that $(\Sigma,\vec z^+)$ may not be stable and so 
we stabilize it by adding a finite number of additional interior marked points
$$
\vec w^{+} = (w_a^+)_{a\in A}, \quad m = \#A, 
$$
on $\Sigma$ so that the following conditions are satisfied: 
such a choice, especially the one satisfying Items 3 and 4,  
is possible by Assumption \ref{trivialisotopystrong}.
\begin{conds}[Stabilization]\label{zacondition}
\begin{enumerate}
\item $w_a^+$ is disjoint from $\vec z^+$, the singular point set
and the boundary of $\Sigma$.
\item $(\Sigma,\vec w^{+})$ is stable. Moreover it has no nontrivial automorphism.
Exception: In the situation appearing in Remark \ref{pbubblerem} we do not put any of $\vec w^{+}$ on the 
disk component where $u$ is constant.
\item $u$ is an immersion at $w_a^+$.
\item For any $\frak t \in Lie(T^n)$ there exists $w_a^+$ such that
$\frak t(u(w_a^+))$ is not contained in
$u_*T_{w_a^+}\Sigma$.
\item The subset $\vec w^{+} \subset \Sigma$ is $G(\text{\bf x})$ invariant, i.e.,
$\varphi(w^+_a) \in \vec w^{+}$ for all $a \in A$, $(g,\varphi) \in G(\text{\bf x})$.
\end{enumerate}
\end{conds}
We denote by $\text{\bf v} = (\Sigma,\vec z^+,\vec w^{+})$ the resulting stable curve contained in $\mathcal M_{0;\ell+m}$.
Let $\frak S_{\ell}$, $\frak S_{m}$ be permutation groups of order $\ell!$, $m!$,
respectively.
The group $\frak S_{\ell}$ acts on $\mathcal M_{0;\ell+m}$ as permutation of
$1$-st,\dots,$\ell$-th marked points $\vec z^{+}$ and 
$\frak S_{m}$ acts on it as permutation of
$(\ell+1)$-st,\dots,$(\ell+m)$-th marked points $\vec w^{+}$.
\begin{defn}\label{auto+forv}
We define $\text{\rm Aut}_+(\text{\bf v})$ as the set of all
$\varphi : \Sigma \to \Sigma$ such that:
\begin{enumerate}
\item
$\varphi$ is biholomorphic.
\item
There exists $\rho = (\rho_1,\rho_2) \in \frak S_{\ell} \times \frak S_{m}$
such that
$$
\varphi(z^+_i) = z^+_{\rho_1(i)},
\quad
\varphi(w^+_a) = w^+_{\rho_2(a)}.
$$
\item
In the situation appearing in Remark \ref{pbubblerem} we assume $\varphi$ is identity on the 
disk component where $u$ is constant and which has only one singular point and has no 
marked point. 
\end{enumerate}
\end{defn}
\par
By Condition \ref{zacondition}.2, 
the assignment 
$\varphi \mapsto \rho_2$ gives an
injective homomorphism
$$
\Psi_2 : \text{\rm Aut}_+(\text{\bf v}) \to \frak S_{m}.
$$
We also write the homomorphism $\varphi \mapsto \rho_1$ as
$$
\Psi_1 : \text{\rm Aut}_+(\text{\bf v}) \to \frak S_{\ell}.
$$
We put $\Psi = (\Psi_1,\Psi_2)$.
\par
Let  $\rho = (\rho_1,\rho_2) \in \frak S_{\ell} \times \frak S_{m}$.
Then
$\rho \text{\bf v} = \text{\bf v}$ if and only if $\rho$ is in the image of $\Psi$.
\par
Condition \ref{zacondition}.2 implies that $\Psi : \text{\rm Aut}_+(\text{\bf v})
\to \frak S_{\ell} \times \frak S_m$ is injective. Using this fact, we can
define an action of $\text{\rm Aut}_+(\text{\bf v})$ on a neighborhood of
$\text{\bf v}$ in $\mathcal M_{0;\ell+m}$ as follows:
Let $\text{\bf v}' = (\Sigma,\vec z^{\prime +},\vec w^{\prime +})$
be in a neighborhood of
$\text{\bf v}$ in $\mathcal M_{0;\ell+m}$ and
$\varphi \in \text{\rm Aut}_+(\text{\bf v})$. We put
$\Psi(\varphi) = (\rho_1,\rho_2)$. We then define
\begin{equation}\label{varphiactonv'}
\varphi \cdot \text{\bf v}'
= (\Sigma,\rho_1\vec z^{\prime +},\rho_2\vec w^{\prime +}),
\end{equation}
where
$$
(\rho_1\vec z^{\prime +})_i =z^{\prime +}_{\rho^{-1}_1(i)},
\quad (\rho_2\vec w^{\prime +})_a =w^{\prime +}_{\rho^{-1}_2(a)}.
$$
In other words, the action is by changing the enumeration of the
interior marked points in the same way as $\varphi$ does on $\text{\bf v}$.
This action induces an action of $G(\text{\bf x})$ on a neighborhood of
$\text{\bf v}$ in $\mathcal M_{0;\ell+m}$.
Namely
$$
(g,\varphi) \cdot \text{\bf v}' =
\varphi \cdot \text{\bf v}'.
$$
(See (\ref{Gzactiondef}).)
\par
We take neighborhoods $\frak U(\text{\bf v})$ of $\text{\bf v}$ in $\mathcal M_{0;\ell+m}$
with the following properties.
There exists a neighborhood $\mathcal U(S(\Sigma))$ of the singular point set of
$\Sigma$ such that for any
$\text{\bf v}' = (\Sigma',\vec z^{\prime +},\vec w^{\prime +}) \in \frak U(\text{\bf v})$
we have a smooth embedding
\begin{equation}\label{15apend3}
i_{\text{\bf v}'} : \Sigma \setminus \mathcal U(S(\Sigma))
\to \Sigma'
\end{equation}
such that
\begin{equation}\label{15apend33}
i_{\text{\bf v}'}(z_i^+) = z_{i}^{\prime +},
\quad
i_{\text{\bf v}'}(w_a^+) = w_{a}^{\prime +}.
\end{equation}
We may assume that $i_{\text{\bf v}'}$ is holomorphic on the support of
$E(\text{\bf x})$.
(See \cite[\S 12]{FO}. The map in \cite[(12.10)]{FO}  is a version of
our map $i_{\text{\bf v}'}$ for the case $\partial \Sigma = \emptyset$.)
\par
Furthermore we assume the following condition
 (\ref{equivchangautenum}). 
Let $\text{\bf v}' = (\Sigma',\vec z^{\prime +},\vec w^{\prime +}) \in \frak U(\text{\bf v})$ 
and $\varphi \in \text{\rm Aut}_+(\text{\bf v})$.
We put
$$
\rho =(\rho_1,\rho_2) = \Psi(\varphi)
$$
and require
\begin{equation}\label{equivchangautenum}
i_{\rho\text{\bf v}'} = i_{\text{\bf v}'}\circ \varphi^{-1} :
 \Sigma \setminus \mathcal U(S(\Sigma)) \to \Sigma'.
\end{equation}
\begin{rem}
We note
$$
(i_{\text{\bf v}'}\circ \varphi^{-1})(z^+_i)
= i_{\text{\bf v}'}(z^+_{\rho^{-1}_1(i)})
=i_{\rho\text{\bf v}'}(z^+_i).
$$
\end{rem}
\par
We also choose $i_{\text{\bf v}'}$ 
with the following properties.
Let $(\rho_1,\rho_2) \in \frak S_{\ell} \times \frak S_m$.
Let $\text{\bf v} = (\Sigma,\vec z^+,\vec w^{+})$ 
and  $\text{\bf v}' = (\Sigma,\vec z^{\prime +},\vec w^{\prime +})$.
We put
$$
\text{\bf v}_{(2)} = (\Sigma,\rho_1\vec z^+,\rho_2\vec w^{+}),
\quad
\text{\bf v}'_{(2)} = (\Sigma,\rho_1\vec z^{\prime +},\rho_2\vec w^{\prime +}).
$$
We then require the identity 
$$
i_{\text{\bf v}'_{(2)};\text{\bf v}_{(2)}} = i_{\text{\bf v}';\text{\bf v}}.
$$
(Here we write $i_{\text{\bf v}';\text{\bf v}}$ etc. in place of $i_{\text{\bf v}'}$
to specify $\text{\bf v}$.) 
\par
We also require that $i_{\text{\bf v}'}$ depends continuously on $\text{\bf v}'$
in $C^{\infty}$ sense. 
The existence of $i_{\text{\bf v}'}$ follows from the trivialization of the core part of the universal family 
with coordinate at infinity in the sense of \cite[Definition 16.2]{foootech} (in particular, (3), (5)).  
\par
We fix $\text{\bf x} = (\Sigma,\vec z^+,u)$ and $\text{\bf v} = (\Sigma,\vec z^+,\vec w^+)$.
For each $w_a^+$ we choose a normal slice $N_{w^+_a}$ of $u$ such that:
\begin{conds}[Normal slices]\label{1.8}
\begin{enumerate}
\item $N_{w^+_a}$ is a codimension $2$ smooth submanifold of $X$.
\item $N_{w^+_a}$ is perpendicular to $u(\Sigma)$ at
$u(w_a^+)$.
\item If $(g,\varphi) \in G(\text{\bf x})$, then
$gN_{w^+_a} = N_{\varphi(w^+_a)}$.
\end{enumerate}
\end{conds}
We note that the choices of the above slices $N_{w^+_a}$ depend only on $\text{\bf x}
= (\Sigma,\vec z^+,u)$ and
$\vec w^+=\{w_a^+\}_{a \in A}$
which will be \emph{fixed once and for all} when the latter is given. They also carry the metrics induced from
one given in the ambient space $X$.
\par
We next consider $\text{\bf x}' = (\Sigma',\vec z^{\prime +},u')$ 
and the following condition on it.
Let $\epsilon_1$, $\epsilon_2$ be positive real numbers.
\begin{defn}\label{primeisclose}
Let $h \in T^n$.
We say that $\text{\bf x}'$ is {\it $(\epsilon_1,\epsilon_2)$-close to $h\text{\bf x}$}
if
there exist $\vec w^{\prime +}$ and a neighborhood $\mathcal U(S(\Sigma))$ of the singular point set of 
$\Sigma$ satisfying \eqref{15apend3} and \eqref{15apend33} with the
following properties.
\begin{enumerate}
\item $\vec w^{\prime +}$ is disjoint from $\vec z^{\prime +}$, the singular point set
and the boundary of $\Sigma^{\prime}$.
\item $\text{\bf v}' = (\Sigma',\vec z^{\prime +},\vec w^{\prime +}) \in \frak U(\text{\bf v})$
and $\text{dist}(\text{\bf v}',\text{\bf v}) < \epsilon_1$. Here we use
an appropriate metric on $\mathcal M_{0;\ell+m}$ to define the distance between $\text{\bf v}'$ and $\text{\bf v}$.
\item
The $C^1$ distance between the two maps $u' \circ i_{\text{\bf v}'}$ and $hu$ is smaller than $\epsilon_1$
on $\Sigma \setminus \mathcal U(S(\Sigma))$.
\item $\text{diam}(u'(\mathcal S)) < \epsilon_1$ if
$\mathcal S$ is any connected  component of
$\Sigma' \setminus  i_{\text{\bf v}'}(\Sigma \setminus \mathcal U(S(\Sigma)))$.
\item
$u'(w^{\prime +}_a) \in hN_{w^+_{a}}$ for each $a$
and
$$
\text{dist}(u'(w^{\prime +}_a),hu(w^+_{a}))  < \epsilon_2.
$$
Here $\vec w^{\prime +} _a = (w^{\prime +}_{a})_{a\in A}$.
\end{enumerate}
\par
When $\vec w^+,\vec w^{\prime +},h$ are specified in the above condition, we say that
$\text{\bf x}' = (\Sigma',\vec z^{\prime +},u')$ is {\it $(\epsilon_1,\epsilon_2)$-close to} $(h\text{\bf x},\vec w^+)$
with respect to $\vec w^{\prime +}$.
\end{defn}
Definition \ref{primeisclose} means that $(\Sigma',\vec z^{\prime +},u')$
is sufficiently close to $h\text{\bf x}$.
\par
Later on we will take $\epsilon_1$, $\epsilon_2$ so that $\epsilon_2$ is small and $\epsilon_1$ is smaller
than a number depending on  $\epsilon_2$ (and $\text{\bf x}$).
(See the proof of Lemma \ref{numberlocalmini}.)
\par
\begin{lem}\label{Gxindependence}
Consider $u'$, $\text{\bf v}' = (\Sigma',\vec z^{\prime +},\vec w^{\prime +})$, $h$,
$\vec w^{+}$ given as above and put $\text{\bf x}' = (\Sigma',\vec z^{\prime +},u')$.
Let $\text{\bf g} = (g,\varphi) \in
G(\text{\bf x})$ and $\sigma = (\sigma_1,\sigma_2) = \Psi(\varphi)$.
Suppose $\text{\bf x}' $ is $(\epsilon_1,\epsilon_2)$-close to  $(h\text{\bf x},\vec w^+)$
with respect to $\vec w^{\prime +}$.
We consider $\sigma\text{\bf v}'=(\Sigma',\sigma\vec z^{\prime +},\sigma\vec w^{\prime +})$,
where
$(\sigma z^{\prime +})_{i} = z^{\prime +}_{\sigma_1^{-1}(i)}$,
$(\sigma w^{\prime +})_{a} = w^{\prime +}_{\sigma_2^{-1}(a)}$.
We put $\sigma\text{\bf x}' = (\Sigma',\sigma\vec z^{\prime +},u')$.
\par
Then $\sigma\text{\bf x}' $ is $(\epsilon_1,\epsilon_2)$-close to  $(g^{-1}h\text{\bf x},\vec w^+)$
with respect to $\sigma\vec w^{\prime +}$.
\end{lem}
\begin{proof}
Definitions \ref{primeisclose}.1, 2 and 4 are easy to check.
Let us check  Definition \ref{primeisclose}.3.
We observe 
$$
u'\circ i_{\sigma\text{\bf v}'}
=
u' \circ i_{\text{\bf v}'} \circ \varphi^{-1}
$$
by (\ref{equivchangautenum}).
On the other hand,
$$
h \circ u \circ \varphi^{-1} = hg^{-1} \circ u,
$$
since $(g,\varphi) \in G(\text{\bf x})$.
By assumption
$u' \circ i_{\text{\bf v}'}$ is $C^1$ close
to $h \circ u$. We thus obtain Definition \ref{primeisclose}.3.
(Note $hg^{-1} = g^{-1}h$.) 
Definition \ref{primeisclose}.5 can be checked by
$$
u'(w^{\prime +}_{\sigma^{-1}_2(a)}) \in h N_{w^{+}_{\sigma^{-1}_2(a)}} =
h N_{\varphi^{-1}(w^{+}_{a})} =  hg^{-1} N_{w^+_a}.
$$
\end{proof}
In the next lemma and thereafter we denote by $\frak o(\epsilon)$ 
\index{$\frak o(\epsilon)$} positive numbers 
depending on $\epsilon$ such that $\lim_{\epsilon \to 0}\frak o(\epsilon) = 0$.
\begin{lem}\label{w+primeunique}
Let $\mu_1 \in \frak S_{\ell}$ and let $\epsilon_1,\epsilon_2$ be sufficiently small.
Suppose $\text{\bf x}' $ is $(\epsilon_1,\epsilon_2)$-close to  $(h\text{\bf x},\vec w^+)$
with respect to $\vec w_{(1)}^{\prime +}$
and
$\mu_1\text{\bf x}' $ is $(\epsilon_1,\epsilon_2)$-close to  $(h'\text{\bf x},\vec w^+)$
with respect to $\vec w_{(2)}^{\prime +}$.
\par
Then there exists
$(g,\psi) \in G(\text{\bf x})$ such  that
\begin{equation}\label{wchikai}
d(\vec w_{(2)}^{\prime +},\mu_2 \vec w_{(1)}^{\prime +}) 
< \frak o(\epsilon_1)
\end{equation}
with $\mu = (\mu_1,\mu_2) = \Psi(\psi)$
and 
\begin{equation}\label{form4314}
{\rm dist}(h',gh) < \frak o(\epsilon_1).
\end{equation}
\par
If $h=h'$ in adition then $g=1$ and 
$\vec w_{(2)}^{\prime +} = \mu_2 \vec w_{(1)}^{\prime +}$.
\end{lem}
\begin{rem}
In (\ref{wchikai}) above and  several times later, we use a metric $d$ on the source curve.
When we use it we are always considering certain compact family of source curves  
and we consider the distance among two points which are not singular point and lie 
on the same irreducible component.
Therefore, such metric can be chosen so that the expression such as appearing 
in (\ref{wchikai}) is independent of its choice.
\end{rem}
\begin{proof}
We put $\text{\bf v}'_{(k)} = (\Sigma',\vec z^{\prime +},\vec w^{\prime +}_{(k)})$,
$k = 1,2$.
For each $a \in \{1,\ldots,m\}$, we have ${}_{(2)}w^+_a \in \Sigma$ such that
$$
i_{\text{\bf v}'_{(1)}}({}_{(2)}w^+_a) = w^{\prime +}_{(2),a}.
$$
(We note that $i_{\text{\bf v}'_{(2)}}(w^+_a)=w^{\prime +}_{(2),a}$ 
but $i_{\text{\bf v}'_{(1)}}(w^+_a)$ may not be close to $w^{\prime +}_{(2),a}$.
So ${}_{(2)}w^+_a$ may be far away from $w^+_a$.) 
Then
$(\Sigma,\mu_1\vec z^+,{}_{(2)}\vec w^+)$  is close to
$\text{\bf v}'_{(2)}$ in $\mathcal M_{0; \ell + m}$.
On the other hand, by assumption,
$(\Sigma,\vec z^+,\vec w^+)$ is close to $\text{\bf v}'_{(2)}$ in $\mathcal M_{0;\ell + m}$.
(See Figure 4.3.1.)
\begin{sublem}\label{sublem4320}
There exists a biholomorphic map $\psi : \Sigma \to \Sigma$
such that
$$
\psi(z^+_a) = z^+_{\mu_1(a)} 
$$
and
$$
{\rm dist}(\psi(w^+_a),{}_{(2)}w_a^+) < \frak o(\epsilon_1).
$$
\end{sublem}
\begin{proof}
Using the fact that $(\Sigma,\vec z^+,\vec w^+)$ is $\epsilon_1$-close to $(\Sigma,\mu_1\vec z^+,{}_{(2)}\vec w^+)$, 
we can prove that $(\Sigma,\vec z^+) = (\Sigma,\mu_1\vec z^+)$
in $\mathcal M_{0;\ell}$.  
This is because $\mu_1$ is in a finite set $\frak S_{\ell}$.   
Note that $(\Sigma, \vec{z}^+)$ is the underlying data of ${\bf x}$, which we fix in the argument. 
For each $\mu'_1 \in {\frak S}_{\ell}$, we have either that $(\Sigma, \vec{z}^+) = \mu'_1 (\Sigma, \vec{z}^+)$ 
or that there exists a positive constant $\delta=\delta((\Sigma. \vec{z}^+),\mu'_1)$ such that the distance between 
$(\Sigma, \vec{z}^+)$ and $\mu'_1 (\Sigma, \vec{z}^+)$ is at least $\delta$.   
Hence, if  $\epsilon_1 >0$ is chosen sufficiently small compared with $\delta$, 
we find that $(\Sigma,\vec z^+) = (\Sigma,\mu_1\vec z^+)$
in $\mathcal M_{0;\ell}$.  
Then the sublemma follows easily.
\end{proof}
\par
We next prove the following.
\begin{sublem}\label{sublem4319}
$d_{C^0}(h u \circ \psi,h' u) < \frak o(\epsilon_1)$.
\end{sublem}
\begin{proof}
The proof is by contradiction.
Suppose we have $\epsilon_1^n$ which converge to zero and 
$u_n, \psi_n, h_n, h'_n$, 
${}_n{\bf v} = ({}_n\Sigma, {}_n\vec{z}^+, {}_n\vec{w}^+)$, 
${}_n{\bf v}'_{(k)} = ({}_n\Sigma', {}_n\vec z^{\prime +},{}_n\vec w^{\prime +})$
($k=1,2$), 
${}^n_{(2)}w^+_a$, ${}^nw^+_a$,
with the following properties:
we write $i_{{}_n{\bf v}'_{(k)}} = i_{{{}_n{\bf v}'_{(k)}}, {}_n{\bf v}}$.

\begin{enumerate}
\item
$
i_{{}_n{\bf v}'_{(1)}}({}^n_{(2)}w^{+}_a) = {}^nw^{\prime +}_{(2),a}. 
$
\item
$
i_{{}_n{\bf v}'_{(2)}}({}^nw^+_a) = {}^nw^{\prime +}_{(2),a}.
$
\item
$
d(\psi_n({}^nw^+_a),{}^n_{(2)}w^+_a) < \frak o(\epsilon_1^n) \to 0.
$
\item
$\psi_n({}_nz^+_a) = {}_nz^+_{\mu_1(a)} $.
\item
$
d_{C^0}(h_nu_n\circ \psi_n, h'_n u_n) > \delta > 0.
$
\end{enumerate}
Here $i_{{}_n{\bf v}'_{(k)}} : \Sigma \setminus \mathcal U_n(S(\Sigma)) 
\to {}_n\Sigma'$ is a map such that:
\begin{enumerate}
\item[6.]
$
d_{C^0}(u'_n \circ i_{{}_n{\bf v}'_{(1)}}, h_nu) < \epsilon_1^n \to 0.
$
\item[7.]
$
d_{C^0}(u'_n \circ i_{{}_n{\bf v}'_{(2)}}, h'_nu) < \epsilon_1^n \to 0.
$
\item[8.]
The diameter of the $u$ image of each connected component of 
$\mathcal U_n(S(\Sigma))$ is smaller than $\epsilon_1^n$ that converges to $0$.
\end{enumerate}

Item 8 implies that the intersection $\bigcap_n \mathcal U_n(S(\Sigma))$
is $S(\Sigma)$ the set of singular points of $\Sigma$.
\par
By taking a subsequence, we may assume that $i_{{}_n{\bf v}'_{(k)}}$
converges to $i_{{}_{\infty}{\bf v}'_{(k)}} : \Sigma \to \Sigma$
and $\psi_n$ converges to $\psi_{\infty} : \Sigma \to \Sigma$.
Moreover we may assume ${}^nw^+_a$ 
converges to ${}^{\infty}w^+_a$.
\par
Then Items 1,2,3,4 imply
$$
(i_{{}_{\infty}{\bf v}'_{(1)}} \circ \psi_{\infty})({}^{\infty}w^+_a)
=
\lim_{n\to \infty} {}^nw^{\prime +}_{(2),a}
=
i_{{}_{\infty}{\bf v}'_{(2)}}({}^{\infty}w^+_a).
$$
Note $i_{{}_{\infty}{\bf v}'_{(k)}}$ is a biholomorphic map. Therefore
Condition \ref{zacondition}.2 implies
$$
i_{{}_{\infty}{\bf v}'_{(1)}} \circ \psi_{\infty}
= 
i_{{}_{\infty}{\bf v}'_{(2)}}.
$$
Then Items 6 and 7 imply  
$$
h u \circ \psi_{\infty} = h' u.
$$
This contradicts to Item 5.
The proof of the sublemma is complete.
\end{proof}
Sublemma \ref{sublem4319} and  the finiteness of $G({\bf x})$ imply that 
there exists $g\in T^n$ such that 
$(g,\psi) \in G({\bf x})$
and ${\rm dist}(h',gh) < \frak o(\epsilon_1)$, if we choose $\epsilon_1$ small.
\par
We put ${}_{(3)}\vec w^+ = \psi(\vec w^+)$.
By Condition  \ref{zacondition}.5 it implies that there exists $\mu_2$ such that
${}_{(3)}\vec w^+ = \mu_2\vec w^+$.
For this $\mu_2$ we can prove $\mu = (\mu_1,\mu_2) = \Psi(\psi)$ easily. 
(\ref{wchikai}) follows from Sublemma \ref{sublem4320}.
\par
Suppose $h=h'$ in addition. 
Then (\ref{form4314}) and the finiteness of $G({\bf x})$ 
implies $g=1$.
Moreover since
$$
w^{\prime +}_{(2), a}, w^{\prime +}_{(1), \mu_2^{-1}(a)}
\in hN_{w^+_{\mu_2^{-1}(a)}}
$$
we have $w^{\prime +}_{(2), a} = w^{\prime +}_{(1), \mu_2^{-1}(a)}$.
\end{proof}
Lemma \ref{w+primeunique} and its proof is an equivariant version of 
\cite[Lemma 20.15]{foootech}.
\par\newpage
\epsfbox{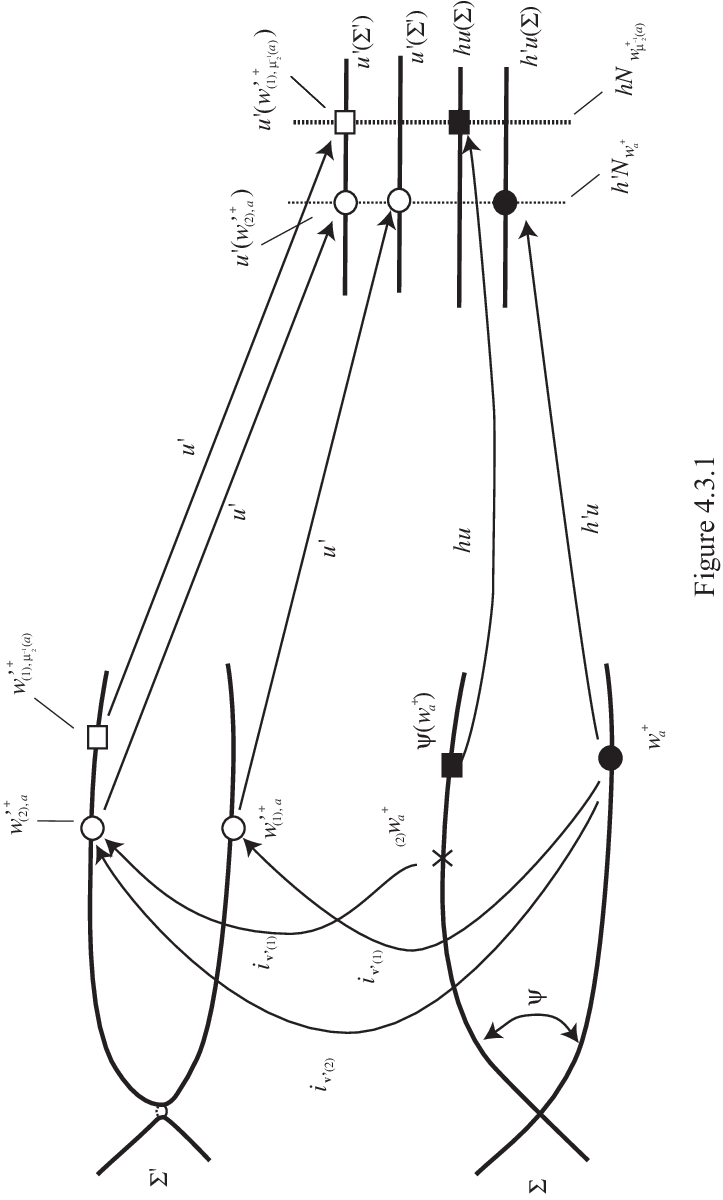}
\par
\par
Suppose $\text{\bf x}' $ is $(\epsilon_1,\epsilon_2)$-close to $(h\text{\bf x},\vec w^+)$
with respect to $\vec w^{\prime +}$.
We define
\begin{equation}\label{deffbyave}
f(h;\text{\bf x}',\vec w^{\prime +})
= \sum_a \text{dist}(u'(w^{\prime +}_a),hu(w^+_{a}))^2.
\end{equation}
Here $\text{dist} = \text{dist}_{hN_{w^+_a}}$is the distance function of the submanifold
$hN_{w^+_a}$ induced from the given ambient K\"ahler metric.
Recall the square of the distance function is strictly convex and smooth in a neighborhood of
the diagonal.
\begin{cor}\label{cortow+primeunique}
Assume that $h=h'$. 
Let $(\text{\bf x}',h,\vec w^{\prime +}_{(1)})$ and
$(\mu_1\text{\bf x}',h,\vec w^{\prime +}_{(2)})$ be
as in Lemma \ref{w+primeunique}. Then
$$
f(h;\text{\bf x}',\vec w^{\prime +}_{(1)})
=
f(h;\mu_1\text{\bf x}',\vec w^{\prime +}_{(2)}).
$$
\end{cor}
This is immediate from  Lemma \ref{w+primeunique} and the definition of $f$.
\par
Once we fix $({\bf x}, \vec{w}^+)$, 
Corollary \ref{cortow+primeunique}, in the case of $\mu_1 = id$, implies that
$f(h;\text{\bf x}',\vec w^{\prime +})$ depends only on
$h,\text{\bf x}'$ and is independent of $\vec w^{\prime +}$ in Definition \ref{primeisclose} 
(It certainly depends on $\vec{w}^+$.)    
Hereafter we write $f(h;\text{\bf x}')$
instead of $f(h;\text{\bf x}',\vec w^{\prime +})$.
\par
For given $\text{\bf x}', \vec w^{+}$, our function $f$
is a function of $h$.
More precisely speaking, we define the subset
$\frak H(\text{\bf x}') \subset T^n$  by
$$
\frak H(\text{\bf x}') =\{h \in T^n \mid \exists  \mu_1 \in \frak S_{\ell},
\text{$\mu_1\text{\bf x}' $ is $(\epsilon_1,\epsilon_2)$-close to  $(h\text{\bf x},\vec w^+)$} \}.
$$
This is an open subset of $T^n$ for each given $\text{\bf x}'$.
Then for each given $\text{\bf x}'$ and sufficiently small $\epsilon_1,\epsilon_2> 0$,
the assignment $h \mapsto f(h;\mu_1 \text{\bf x}')$
defines a function on $\frak H(\text{\bf x}')$.
(Here $\mu_1$ is as in the definition of $\frak H(\text{\bf x}')$.)
We write this function as $f(h;[\text{\bf x}'])$.
(We remark that $\frak H(\text{\bf x}')$ is a disjoint union of small open subsets of 
$T^n$.)
\par
Since $\frak H(\text{\bf x}')$ is an open subset of $T^n$, it
carries a natural affine structure. We use this to define convexity
used in the statement of the following lemma.
\begin{lem}\label{strictconvexf}
If $\epsilon_1,\epsilon_2$ are sufficiently small, then $f(\cdot;[\text{\bf x}']) :
\frak H(\text{\bf x}') \to \R_{\ge 0}$ is strictly convex on each connected component.
\end{lem}
\begin{proof}
This is a consequence of strictly convexity of the square of the Riemannian distance function  and Condition \ref{zacondition}.4.
\end{proof}
\begin{lem}\label{Gxindependence2}
Let $\text{\bf g} = (g,\varphi) \in G(\text{\bf x})$.
Then we have
$$
f(h;[\text{\bf x}'])
= f(h;[\text{\bf g}\text{\bf x}']).
$$
Here $\text{\bf x}' = (\Sigma',\vec z^{\prime +},u')$ and
$\text{\bf g}\text{\bf x}' = (\Sigma',\Psi_1(\varphi)\vec z^{\prime +},gu')$.
\end{lem}
\begin{proof}
Lemma \ref{Gxindependence} imlies
$$
f(h;[\text{\bf x}']) = f(g^{-1}h;[\Psi_1(\varphi)\text{\bf x}']).
$$
On the other hand, it is easy to see that
$$
f(h;[g\text{\bf x}']) = f(g^{-1}h;[\text{\bf x}']).
$$
The lemma follows easily.
\end{proof}
We recall that
$$
\pi_{T^n}(G(\text{\bf x})) =
\{ g \in T^n\mid \exists \varphi,  (g,\varphi) \in G(\text{\bf x})\}.
$$
Let $\frak c = \#\pi_{T^n}(G(\text{\bf x}))$ be its order and
write as
$
\pi_{T^n}(G(\text{\bf x}))= \{g_1=1,\ldots,g_{\frak c}\}.
$
We put
$$
G(\text{\bf x})
= \{\text{\bf g}_{i,j} =  (g_i,\varphi_{i,j}) \mid i = 1,\ldots, \frak c,
 j = 1,\ldots, \frak d\}.
$$
Here
$
\frak d = \# \{ \varphi \mid (1,\varphi) \in G(\text{\bf x})\}
$ and
$\varphi_{1,1} = id.$
\begin{lem}\label{numberlocalmini}
For each sufficiently small $\epsilon_2$ we may choose $\epsilon_1$ small enough
so that the following holds: Whenever $\frak H(\text{\bf x}')$ is nonempty,
$f(\cdot ;[\text{\bf x}']) : \frak H(\text{\bf x}') \to \R_{\ge 0}$ assumes local minima at exactly
$\frak c$ points $h_1,\dots, h_{\frak c}$.
\par
Furthermore the following also holds: 
For each $i$ there exist $\frak d$ choices
$\vec w^{\prime +}_{(i,j)}$ and $\rho_{(i,j),1} \in \frak S_{\ell}$ such that
$\rho_{(i,j),1}\text{\bf x}' $ is $(\epsilon_1,\epsilon_2)$-close to  $(h_i\text{\bf x},\vec w^+)$
with respect to $\vec w_{(i,j)}^{\prime +}$.
\par
If we put $\sigma_{(i,j)} = (\sigma_{(i,j),1},\sigma_{(i,j),2}) = \Psi(\varphi_{i,j})$
for each $\text{\bf g}_{i,j}= (g_i,\varphi_{i,j})$, we have
\begin{equation}\label{formulafor1170}
h_i = h_1g_i^{-1}
\end{equation}
and
\begin{equation}\label{formulafor1172}
 \rho_{(i,j),1} = \sigma_{(i,j),1}\rho_{(1,1),1},
\qquad w^{\prime +}_{(i,j)} = \sigma_{(i,j),2} w^{\prime +}_{(1,1)}.
\end{equation}
\end{lem}
\begin{proof}
We first consider the case $\text{\bf x}'= \text{\bf x}$.
By Lemma \ref{Gxindependence2} we have
\begin{equation}
f(h; \sigma_{(i,j),1}\text{\bf x}') = f(g_i^{-1}h; \text{\bf x}').
\end{equation}
Therefore, $(\Sigma,\sigma_{(i,j),1}\vec z^{+},u)$ is
$(\epsilon_1,\epsilon_2)$-close to
$(g_i^{-1}\text{\bf x},\vec w^+)$ with respect to
$\sigma_{(i,j),2}\vec w^+$.
Moreover, $g_i^{-1}$ is a critical point of 
$f(\cdot;\sigma_{(i,j),1}\text{\bf x})$.
Namely $h_i$ defined by (\ref{formulafor1170}) 
is a critical point and (\ref{formulafor1172}) is satisfied
for $h_1 = 1$, $\rho_{(1,1),1} = id$, $w^{\prime +}_{(1,1)} = \vec w^+$.
\par
We next show that they are all the local minima. Let $h\in \frak H({\bf x}')$. Lemma
\ref{w+primeunique} implies that there exists $g_i$ such that $g_i^{-1}$ is close to $h$.
Since $f(\cdot;[{\bf x}])$ is convex and since it is minimal at $h$, we have $h=g_i^{-1}$.
\par
Let us consider the general $\text{\bf x}' = (\Sigma',\vec z^{\prime +},u')$.
Let $h_0 \in \frak H(\text{\bf x}')$.
Replacing $u'$ by $h_0^{-1}u'$ we may  and will assume $1 \in \frak H(\text{\bf x}')$. 
Let $h \in  \frak H(\text{\bf x}')$. By Lemma \ref{w+primeunique} again there exists ${\bf g}= (g_i,\varphi_{i,j}) \in G(\text{\bf x})$ 
such that $g_i^{-1}$ is close to $h$. Replacing ${\bf x}$ by ${\bf g}{\bf x}$ we may 
and will assume that $h$ is close to $1$. 
Then we can find a short path $\text{\bf x}(t)$ joining
$\text{\bf x}$ to $\text{\bf x}'$.
By the strict convexity of $f$, the number of local minima is
constant on this (short) path.
We may choose $\epsilon_1$ small enough so that
the distance between $u'(w^{\prime +}_{(i,j),t;a})$ and $h_{i,t}u(w^{+}_{a})$
does not go beyond $\epsilon_2$ while $t$ moves from $0$ to $1$.
(Here $h_{i,t}$ and $w^{\prime +}_{a,t}$ are such that
$\rho_{(i,j),1}\text{\bf x}(t)$ is $(\epsilon_1,\epsilon_2)$-close to  $(h_{i,t}\text{\bf x},\vec w^+)$
with respect to $\vec w^{\prime +}_{(i,j),t}$.)
\par
We observe that if (\ref{formulafor1170}), (\ref{formulafor1172}) hold and
$\rho_{(1,1),1}\text{\bf x}' $ is $(\epsilon_1,\epsilon_2)$-close to  $(h_1\text{\bf x},\vec w^+)$
with respect to $\vec w_{(1,1)}^{\prime +}$,
then
$$
u'(\vec w'_{(i,j),a}) = u'(\vec w'_{(1,1),\sigma^{-1}_{(i,j),2}(a)})
\in h_1N_{w^+_{\sigma^{-1}_{(i,j),2}(a)}} = h_1g^{-1}_i N_{w^+_{a}}.
$$
Namely $\rho_{(i,j),1}\text{\bf x}'$ is $(\epsilon_1,\epsilon_2)$-close to
$(h_i\text{\bf x},\vec w^+)$ with respect to $\vec w^{\prime +}_{(i,j)}$.
\par
Using Lemma \ref{Gxindependence2}, we can show that
$f(\cdot;[\text{\bf x}'])$ assumes local minima at $h_i$.
The proof of (\ref{formulafor1170}) of Lemma \ref{numberlocalmini} is complete.
\par
We can prove  (\ref{formulafor1172}) using the fact that it holds for $t=0$ and
the continuity argument.
\end{proof}
Now we use Lemma \ref{numberlocalmini} to send our obstruction vector space $E({\bf x}) = E_{\bf x}$ to ${\bf x}'$ as follows.
Here we assume the following.
\begin{conds}\label{immersionsup}
On the support of each element of $E({\bf x})$ the map $u$ is an immersion.
\end{conds}
Since $\beta \ne 0$ we can choose such $E({\bf x})$ by the unique continuation theorem as in Lemma \ref{app2lema}.
(We remark that we do not need to put obstruction bundle on the irreducible component where 
$u$ is constant, since such an irreducible component is stable and genus $0$.)
\par
Let $h_i = h_1g_i^{-1}, \sigma_{(i,j),1}, \sigma_{(i,j),2}, w^{\prime +}_{(i,j)} = \sigma_{(i,j),2} \vec w^{\prime +}_{(1,1)}$ be as in Lemma \ref{numberlocalmini}. In particular $\sigma_{(i,j),1} {\bf x}'$ is $(\epsilon_1,\epsilon_2)$ close to 
$(h_i{\bf x},\vec w^+)$.
Therefore there exists 
$i_{i,j} : \Sigma \to \Sigma'$ which sends $\vec z^+$ to $\sigma_{(i,j),1}\vec z^{\prime +}$ 
and $\vec w^+$ to $\sigma_{(i,j),2}\vec w^{\prime +}$.
We modify $i_{i,j}$ slightly to define 
\begin{equation}
I_{i,j} : {\rm Supp}(E(\text{\bf x})) \to \Sigma'
\end{equation}
as follows. (Here ${\rm Supp}(E(\text{\bf x}))$ is the union of the  supports of the 
elements of $E(\text{\bf x})$.)
\par
Let $x \in {\rm Supp}(E(\text{\bf x}))$ we require :
\begin{conds}\label{condsIII}
\begin{enumerate}
\item 
${\rm dist}(i_{i,j}(x), I_{i,j}(x)) < \epsilon_3$.
\item
The minimal geodesic $\ell_{x;i,j}$ joining $h_iu(x)$ to $u'(I_{i,j}(x))$ is perpendicular 
to $h_iu(\Sigma)$ at $h_iu(x)$.
\end{enumerate}
\end{conds}
We can choose $\epsilon_3$ sufficiently small depending on 
$u$ and $\vec w^+$ and then choose $\epsilon_1, \epsilon_2$ depending on $\epsilon_3$ etc.
so that there exists unique $I_{i,j}(x)$ satisfying Condition \ref{condsIII}.
This is a consequence of Condition \ref{immersionsup}.
\par
Now we define
\begin{equation}\label{Pprelimi}
\mathcal P_{\text{\bf x}^{\prime},i,j} :
E(\text{\bf x}) \to C^{\infty}(\Sigma',u^{\prime *}TX \otimes \Lambda^{0,1}).
\end{equation}
\par
For $x \in \Sigma$ in the support of $E(\text{\bf x})$,
we consider the minimal geodesic $\ell_{x;i,j}$ joining $hu(x)$ to $u'(I_{i,j}(x))$.
We also 
take the complex linear part 
$$
\text{\rm Pal}_{x;i,j}^c
$$ 
of the parallel transport map
$$
\text{\rm Pal}_{\ell_{x;i,j}} :
T_{h_iu(x)}X   \to T_{u'(I_{i,j}(x))}X.
$$
along the geodesic $\ell_{x;i,j}$.
(Here we use an appropriate connection as in \cite[Section 7.1]{fooobook2}  for which our Lagrangian submanifold 
is totally geodesic and which is $T^n$ invariant.)
Using also the facts that $I_{i,j}$ is $C^1$ close to a 
biholomorphic map and that
$h_i$ is holomorphic, we obtain (\ref{Pprelimi}) as follows: 
We take the successive compositions of the following
target maps
\begin{equation}\label{287}
(u^*TX)_x = T_{u(x)}X
\overset{(h_i)*}{\longrightarrow}  T_{h_iu(x)}X
\overset{\text{\rm Pal}_{x;i,j}^c}{\longrightarrow}
T_{u'(I_{i,j}(x)))}X,
\end{equation}
and then take the tensor product with the map
\begin{equation}\label{288}
\Lambda_x^{0,1}(\Sigma)
\overset{(I_{i,j}^{-1})^*}{\longrightarrow} \Lambda_{I_{i,j}(x)}^{1}(\Sigma')
\to \Lambda_{I_{i,j}(x)}^{0,1}(\Sigma').
\end{equation}
Here the second map in (\ref{288}) is the projection.
\begin{lem}\label{115}
$\mathcal P_{\text{\bf x}^{\prime},i,j} \circ {\bf g}_{i,j *}^{-1}
$
is independent of $i = 1,\dots, \frak c$ and  $j = 1,\dots, \frak d$.
\end{lem}

\begin{proof}
By  (\ref{equivchangautenum}), we have
$
i_{i,j}
= i_{1,1} \circ \varphi_{i,j}^{-1}.
$
Moreover $h_{1} g^{-1}_i u = h_{1} u \varphi_{i,j}^{-1}$. 
Here recall ${\bf g}_{i,j}=(g_i, \varphi_j) \in G({\bf x})$.
They imply
$
I_{i,j}
= I_{1,1} \circ \varphi_{i,j}^{-1}$.
In fact, together with Lemma \ref{numberlocalmini}
they imply
\begin{equation}\label{nfin4328}
h_iu(x) = h_1 u( \varphi_{i,j}^{-1}(x)).
\end{equation}
Therefore
$$
\ell_{x;i,j} = \ell_{\varphi_{i,j}^{-1}(x);1,1}
$$
since they join (\ref{nfin4328}) to a point of $u'(\Sigma')$ 
that is close to $u'(i_{i,j}(x)) = u'(i_{1,1}(\varphi_{i,j}^{-1}(x))$
and they are perpendicular to $h_iu(\Sigma) = h_1u(\Sigma)$.
\par
Since $I_{i,j}(x)$ and $ I_{1,1} \circ \varphi_{i,j}^{-1}(x)$ 
are the end points of $\ell_{x;i,j}$ and $\ell_{\varphi_{i,j}^{-1}(x);1,1}$ 
respectively, we have $
I_{i,j}
= I_{1,1} \circ \varphi_{i,j}^{-1}$. 
We also recall that $h_i=h_1g_i^{-1}$. Then 
the lemma now follows from the definition.
\end{proof}

We now define
$
\mathcal P_{\text{\bf x}'} :
E(\text{\bf x}) \to
C^{\infty}(\Sigma',u^{\prime *}TX \otimes \Lambda^{0,1})
$
by
$$
\mathcal P_{\text{\bf x}'}=
\mathcal P_{\text{\bf x}^{\prime},i,j} \circ {\bf g}_{i,j *}^{-1}.
$$
(When $G({\bf x})$ is a trivial group, we have $\frak c = \frak d =1$.)
\begin{rem}
In the above definition we replace $i_{i,j}$ by $I_{i,j}$ using 
the assumption that $u$ is an immersion on the domain of $I_{i,j}$.
In the case we are studying here, this process is not necessary.
However it becomes useful later. (See Remark \ref{itoIremark}.)
\end{rem}
\begin{lem}\label{invhhh}
For any $h \in T^n$, we have
\begin{equation}\label{equivalenceP2}
\mathcal P_{h\text{\bf x}'} = h_* \mathcal P_{\text{\bf x}'}.
\end{equation}
\end{lem}
\begin{proof}
This is immediate from the construction.
\end{proof}
We now put
$E(\text{\bf x}') = \text{\rm Im} \mathcal P_{\text{\bf x}'}$.
\par
We consider the assignment ${\bf x}' \mapsto E(\text{\bf x}')$ 
and will prove its $G({\bf x})$ equivariance, by which we mean the following: 
Let $\varphi \in {\rm Aut}_+(\Sigma,\vec z^+)$ and $\Psi_1(\varphi) = \rho_1$.
We define
$$
\varphi\cdot {\bf x}' = (\Sigma', \rho_1 \vec z^{\prime +},u')
$$
where ${\bf x}' = (\Sigma', \vec z^{\prime +},u')$, $(\rho_1 \vec z^{\prime +})_a = z^{\prime +}_{\rho_1^{-1}(a)}$.
Clearly
$
E(\text{\bf x}') 
$
and 
$
E(\varphi\cdot \text{\bf x}') 
$
are both subspaces of $C^{\infty}(\Sigma',u^{\prime *}TX \otimes \Lambda^{0,1})$.
If $(g,\varphi) \in G({\bf x})$, then 
$E((g,\varphi)\cdot\text{\bf x}') $ is a subspace of 
$C^{\infty}(\Sigma',(gu')^{ *}TX \otimes \Lambda^{0,1})$ 
where $(g,\varphi)\cdot\text{\bf x}' = (\Sigma', \rho_1 \vec z^{\prime +},gu')$.
Furthermore we have 
$$
g_* : C^{\infty}(\Sigma',u^{\prime *}TX \otimes \Lambda^{0,1}) \to C^{\infty}(\Sigma',(gu')^{*}TX \otimes \Lambda^{0,1}).
$$
\begin{lem}\label{G(x)equivalence1}
$$
g_*(E(\text{\bf x}')) = E((g,\varphi)\cdot\text{\bf x}').
$$
\end{lem}
\begin{proof}
This is a consequence of the definition and Lemma \ref{115}.
\end{proof}
So far we have defined $E(\text{\bf x}')$ for $\text{\bf x}'$ which is 
close to the $T^n \times G({\bf x})$ orbit of $\text{\bf x}$ in the 
sense of Definition \ref{primeisclose}.
Since we want our Kuranishi structure to be invariant under the permutations of the interior marked points 
$\vec z^{\prime +}$, we extend the assignment 
$\text{\bf x}' \mapsto E(\text{\bf x}')$ and associate  $E(\text{\bf x}')$
to $\text{\bf x}'$ which is 
close to the $T^n \times \frak S_{\ell}$-orbit of $\text{\bf x}$, 
as follows.
\par
Let us consider the homomorphism $\Psi_1 : G(\bf x) \to \frak S_{\ell}$.
We take a complete set of the representatives 
$\rho(1)=1, \dots, \rho(\frak k) \in \frak S_{\ell}$ 
of the right coset space $\frak S_{\ell}/\Psi_1(G(\bf x))$.
For each $\rho(k)$ we consider 
$\rho(k)\text{\bf x} = (\Sigma,\rho(k) \vec z^+,u)$
and use the marked points $\vec w$ (which is the same as one used for ${\bf x}$).
Then if $\text{\bf x}$ is close to the $T^n  G(\rho(k)\text{\bf x})$ that is a 
$T^n$ 
orbit of $\rho(k)\text{\bf x} $, we obtain the obstruction space $E({\bf x}')$.
Since
$$
\bigcup_k (T^n  G(\rho(k) \text{\bf x}))  \rho(k) \text{\bf x}
$$
is the $T^n \times \frak S_{\ell}$-orbit of $\text{\bf x}$ 
and the above union is a disjoint union, we have associated $E(\text{\bf x}')$ to each
$\text{\bf x}'$ that is close to the $T^n \times \frak S_{\ell}$-orbit of $\text{\bf x}$.
By construction and Lemmata \ref{equivalenceP2}, \ref{G(x)equivalence1}, $\text{\bf x}'\mapsto E(\text{\bf x}')$ is  $T^n \times \frak S_{\ell}$ equivariant.
\par
We consider the set of the 
triples $\text{\bf x}' = (\Sigma',\vec z^{\prime +},u)$ satisfying
$$
\overline{\partial}u' \in E(\text{\bf x}').
$$
Together with $E$ and
$$
s(\text{\bf x}') = \overline{\partial} u',
$$
this defines a Kuranishi neighborhood of the  $T^n \times \frak S_{\ell}$ orbit of
$\text{\bf x}$.  (\ref{equivalenceP2}) implies that
it is $T^n$ invariant.
By construction, it is invariant under the permutation of the interior marked points
$\vec z^+$.
\par
We have thus constructed a Kuranishi neighborhood of the $T^n \times \frak S_{\ell}$ 
orbit of $\text{\bf x}$ in case the source $\Sigma$ has only one disk component 
and Assumption \ref{trivialisotopystrong} is satisfied.
\begin{rem}
In \cite[Appendix]{FO}, we took the minimum choice of extra marked points which
breaks the symmetry with respect to the automorphism.
The symmetry is recovered in \cite[Appendix]{FO}  later on.
In this section we take a slightly different way and take
more additional marked points so that they are preserved by the
automorphisms. This approach may be more suitable than the one in
\cite{FO} in the situation where a continuous group of symmetry exists.
(It occurs in the situation we handle next.)
\end{rem}
\par
\subsection{The case of one disk component II}
\label{subsec:onedisk2}
We next remove Assumption \ref{trivialisotopystrong}.
\begin{lem}\label{Tnsplits}
Suppose $\text{\bf x} = (\Sigma,\vec z^+,u)$ does not
satisfy Assumption \ref{trivialisotopystrong}. Then
there exists a splitting $T^n = S^1 \times T^{n-1}$ with
the following properties.
\par
\begin{enumerate}
\item If $\frak t \in Lie(T^{n-1})$, then there exists a
(disk or sphere) component $\Sigma_a$ of $\Sigma$ and
$p \in \Sigma_a$ such that $u$ is an immersion at $p$
and $\frak t(u(p))$ is not contained in
$u_*T_{p}\Sigma$. $p$ is not a singular, boundary, or marked point.
Here we regard $\frak t$ as a vector field on $X$.
\item If $\frak t \in Lie(S^1)$ and $w \in \Sigma$ is not
a singular point and $u$ is an immersion at $w$, then
$\frak t(u(w)) \in u_*(T_w \Sigma)$.
\end{enumerate}
\end{lem}
\begin{proof}
Let $\Sigma_{0}$ be the disk component of $\Sigma$.
\par
We first assume that $u$ is nonconstant on $\Sigma_0$.
Since $T^n$ orbits are totally real and the $T^n$ action on $L$ is free, there exists a at
most one dimensional subspace $\frak T
\subset Lie(T^n)$ such that the following holds:
If $\frak t \notin \frak T$ and $p \in \Sigma_{0}$
is a generic point, then $\frak t(u(p))$ is not contained in
$u_*T_{p}\Sigma_0$. Moreover $\frak t \in \frak T$ satisfies Item 2 above on $\Sigma_0$.
\par
Since Assumption \ref{trivialisotopystrong} is not satisfied, this
subspace $\frak T$ is nonempty. Moreover $\frak t \in \frak T$ satisfies Item 2 above on sphere components also.
\par
We note that the Lie group associated to $\frak T$ must be compact
since otherwise its closure in $T^n$ would have dimension $2$ or more and
its Lie algebra would generate the set of Killing vector fields which
spans the tangent space of the image of $u$.
This would give rise to a contradiction
since $u$ is pseudo-holomorphic.
\par
Since all $S^1$ in $T^n$ are direct summands, we obtain the
required splitting.
\par
Suppose $u$ is constant on $\Sigma_{0}$. Since $\beta \ne 0$, 
there exists a sphere  component $\Sigma_a$ with the following properties: 
$\Sigma_a$
can be joined with $\Sigma_0$ by a path $\gamma$ 
which is contained in a union of components where $u$ is constant: 
$u$ is not constant on $\Sigma_a$.
\par
Note $u(\Sigma_a)$ interests $L$ at $u(x)$ where $x$ is a singular point 
of $\Sigma$ contained in $\Sigma_a$.
The intersection of $T^n$ orbit of $u(x)$ with $u(\Sigma_a)$ has positive 
dimension since Assumption \ref{trivialisotopystrong} is not satisfied.
So we can find a subgroup $S^1 \subset T^n$ such that the $S^1$ orbit of $u(x)$ is 
contained in $u(\Sigma_a)$.
Since $T^n$ orbit is totally real, such $S^1$ is unique.
The rest of the proof is similar to the  case when $u$ is nonconstant on $\Sigma_0$.
\end{proof}
Hereafter we write
$$T^n = S^1_0 \times T^{n-1}_0$$
for the splitting
in Lemma \ref{Tnsplits}. (The factor $T^{n-1}_0$ is not unique. 
But the construction below does not depend on its choice.)
A point of (bordered) curve which is singular or 
marked is called {\it special}.
\begin{lem}\label{markedofsymmetry}
Suppose $\text{\bf x} = (\Sigma,\vec z^+,u)$ does not
satisfy Assumption \ref{trivialisotopystrong}.
Then one of the following two alternatives holds.
\begin{enumerate}
\item There exists a special point $z_0$ on
$(\Sigma,\vec z^+)$
such that $u(z_0)$ is not a fixed point of $S^1_0$.
The group $G(\text{\bf x})$ is finite.
Exception: In the situation of Remark \ref{pbubblerem},
if $u$ on the disk is constant and disk has only one special point 
(the interior singular point) then the identity component 
of  $G(\text{\bf x})$
is $S^1$. Its projection to $T^n$ is $\{1\}$ and its action 
on the sphere bubbles is trivial.
\item
There exists a subset $\mathcal C \subset \Sigma$, 
that is the image of an embedding of a one-dimensional CW complex such that the embedding 
is smooth on each open one-dimensional cell,  
and a subgroup $G^+(\text{\bf x}) \subset T^n \times {\rm Diff}(\Sigma\setminus \mathcal C)$
with the following properties.  
\begin{enumerate}
\item
The group $G^+(\text{\bf x})$ consists of  $(g,\varphi)$ where 
$\varphi : \Sigma\setminus \mathcal C \to\Sigma\setminus \mathcal C$ is a biholomorphic 
map and $g \in T^n$ such that
$$
u(\varphi(x)) = g u(x).  
$$
\item
$G^+(\text{\bf x})$ contains $G(\text{\bf x})$.
\item
The projection $G^+(\text{\bf x}) \to T^n$, $(g,\varphi) \mapsto g$ induces a finite to one group homomorphism 
onto $S_0^1$.
\end{enumerate}
\end{enumerate}
\end{lem}
\begin{proof}
Suppose there exists a special point $z_0$ such that $u(z_0)$ is not a fixed point of $S^1_0$.
Then clearly the Lie algebra of $S^1_0$ does not vanish at $u(z_0)$.
It follows easily that $G(\text{\bf x})$ is finite, 
unless the situation is as in the exception.
(See Example \ref{example43333}.)
\par
We next assume that, for any special point $z_0$ of $(\Sigma,\vec z^+)$, the point
$u(z_0)$ is a fixed point of $S^1_0$.
Let $\Sigma_{\rm tri}$ be the union of all irreducible components of $\Sigma$ on which 
$u$ is a constant map.
We put 
$$
\mathcal C_0 = \{ z \in \Sigma \setminus \Sigma_{\rm tri} \mid d_zu = 0, 
\,\, \text{$u(z)$ is not a fixed point of $S^1_0$} \}.
$$
This is a finite set. 
It is easy to see that the vector field generated by $\frak t \in Lie(S^1_0)$
lifts to a holomorphic vector field on $\Sigma \setminus \mathcal C_0$ that vanishes on the special points 
and on $\Sigma_{\rm tri}$. We denote the lift by $\widetilde{\frak t}$. 
We put
$$
\mathcal C  = \{ z \in \Sigma \setminus \Sigma_{\rm tri} \mid 
u(z) \in S^1_0 u(\mathcal C_0)\}.
$$
It consists of the union of finitely many orbits of $\widetilde{\frak t}$ and $\mathcal C_0$.  
By definition $\mathcal C$ is compact and is one dimensional.
The lift $\widetilde{\frak t}$ generates an $S^1$ action on $\Sigma\setminus \mathcal C$.    
It is clear that the action by $G({\bf x})$ preserves $\mathcal C_0$ and $\mathcal C$, respectively.  
We denote the group generated by this $S^1$ action and $G(\text{\bf x})$
by $G^+(\text{\bf x})
\subset T^n \times {\rm Diff}(\Sigma\setminus \mathcal C)$.
\par
Note that the restriction $u\vert_{\partial \Sigma}$ gives a map 
$S^1 \to L$ whose image is a $1$ dimensional subset of an $S^1_0$ orbit in $L$ .
Therefore for any $(g,\varphi) \in G({\bf x})$ we have $g\in S^1_0$.
In other words, the image of the homomorphism $G({\bf x}) \to T^n$ 
is contained in $S^1_0$. 
Using this fact,
it is easy to see that $G^+(\text{\bf x})$ has the properties claimed in Item 2 above.
\par
The proof of
Lemma \ref{markedofsymmetry} is complete.
\end{proof}
\begin{exm}
Let $X=S^2$ and let $L$ be the equator. 
The group $S^1$ acts on $X$ by rotations that fixes the north and 
south poles.
We pick up $u_1 : D^2 \to X$ that is a biholomorphic map to the northern hemisphere 
and a holomorphic map $u_2 : S^2 \to S^2$.
We glue them in two different ways.
\begin{enumerate}
\item
We glue $D^2$ and $S^2$ at $p_1 \in D^2$, $p_2 \in S^2$. Here $u_1(p_1) = u_2(p_2)$ is 
neither on the equator nor the north pole.
We thus obtain $u : \Sigma \to X$. This is an example where 
Lemma \ref{markedofsymmetry}.1 is satisfied.
\item 
We glue $D^2$ and $S^2$ at $p_1 \in D^2$, $p_2 \in S^2$. 
Here we assume $u_1(p_1) = u_2(p_2)$ is the north pole.
We thus obtain $u : \Sigma \to X$.  
Lemma \ref{markedofsymmetry}.2 is satisfied. 
In fact, if $u_2$ is biholomorphic, then 
$G(\Sigma,u) = S^1$ and $\mathcal C$ is an empty set. 
If $u_2$ is not biholomorphic and has a branched locus $B$, then 
$\mathcal C$ is an inverse image of the $S^1$ orbit of $u(B)$.
\end{enumerate}
\end{exm}
\begin{exm}\label{example43333}
Let $X=S^2$ and $L$ the equator.
We glue $D^2$ and $S^2$ at the center $0\in D^2$ and $z_0 \in S^2$ to obtain $\Sigma$.
Let $u_0 : D^2 \to X$ be a constant map to $p \in L$ 
and $u_1 : S^2 \to X$ a biholomorphic map with $u_1(z_0) = p$.
We can glue $u_0$ and $u_1$ to obtain $u : \Sigma \to X$.
In this case Lemma \ref{markedofsymmetry}.1 exception is satisfied.
\end{exm}
\begin{exm}
Let us identify $S^2 = \C \cup \{\infty\}$, 
$L = \R \cup \{\infty\}$ and 
$D^2 =\{ z \in\C \mid {\rm Im} z \ge 0\} \cup \{\infty\}$.
We consider the map $u : D ^2 \to S^2$ defined by $u(z) = z^2$.
In this case we find $\mathcal C = \{z \in D^2 \mid {\rm Re} z = 0 \} \cup \partial D^2$ and 
$G({\bf x})$ is trivial and  
$G^+({\bf x})\cong S^1$.
\par
This example is called a {\it lantern}.\index{lantern} 
We can compose a branched covering $D^2 \to D^2$ with this map and can construct 
an example where $\mathcal C$ is more complicated and/or
$G({\bf x})$ is a nontrivial finite cyclic group.
\end{exm}
We discuss the two cases in Lemma \ref{markedofsymmetry} separately.
In this subsection we consider Case 1. 
We will discuss Case 2 in the next subsection. 
Now recall our situation: 
Let $(\Sigma,\vec z^+)$ be our marked bordered curve of genus zero
with $\partial\Sigma = S^1$ and let $u : (\Sigma,\partial \Sigma) \to (X,L)$
be pseudo-holomorphic.
We choose finitely many points $\vec w^{+}
= (w^+_a)_{a\in A}$ so that the following is satisfied.
\begin{conds}\label{zacondition2}
\begin{enumerate}
\item $w_a^+$ is disjoint from $\vec z^+$, the singular point set
and the boundary of $\Sigma$.
\item $(\Sigma,\vec w^{+})$ is stable. Moreover the group of its automorphisms is trivial.
Exception: In the situation appearing in Remark \ref{pbubblerem} we do not put any of $\vec w^{+}$ on the 
disk component.
\item $u$ is an immersion at $w_a^+$.
\item For any $\frak t \in Lie(T^{n-1}_0)$ there exists $w^+_a$ such that
$\frak t(u(w_a^+))$ is not contained in
$u_*T_{w^+_a}\Sigma$.
\item The set $\vec w^{+}$ is $G(\text{\bf x})$ invariant, i.e.,
$\varphi(w_a^+) \in \vec w^+$ for all $a \in A$,
$(g,\varphi) \in G(\text{\bf x})$.
\item
Each disk or sphere component remains 
stable after removing 2 points 
of $\vec w^+$.
Exception: In the situation appearing Remark \ref{pbubblerem} we do not assume it 
for the disk component.
\end{enumerate}
\end{conds}
This is almost the same as Condition \ref{zacondition}.
The difference is that we take $T^{n-1}_0$ in place of $T^n$ in
Item 4. 
(We also add Item 6 by certain technical reason. See the proof of 
Lemma \ref{convexitylemma}.) 
The existence of such $\vec w^{+}$ follows from Lemma \ref{Tnsplits}.
\par
Let $\vec w^+$ satisfy Condition \ref{zacondition2}.
We put $\text{\bf v} = (\Sigma,\vec z^+ \cup \vec w^{+})$.
We consider its neighborhood $\frak U(\text{\bf v})$ as before.
\par
We next consider $\text{\bf x}' = (\Sigma',\vec z^{\prime +},u')$
which satisfies the conditions in
Definition \ref{primeisclose}.
We define the set $\frak H(\text{\bf x}')$ and a function
$f(\cdot;\text{\bf x}')$ or $f(\cdot;[\text{\bf x}'])$  on it in the same way. 
We note, however, that $f$ may not be strictly convex in the $S^1_0$ direction.
(This is because we replace $T^n$ by $T^{n-1}_0$ in Condition  \ref{zacondition2}.4.) 
So to generalize Lemma \ref{numberlocalmini} to our situation we need to modify $f$ appropriately.
We use the special points on $ (\Sigma,\vec z^+)$  for this purpose as follows.
\par
Let $\text{\bf x}' =  (\Sigma',\vec z^{\prime +},u')$.
Let $z_c$ be a special point of $(\Sigma,\vec z^+)$.
Below, we will associate a real number $f_c(h;\text{\bf x}')$ to $(h;\text{\bf x}')$.
\par
If $z_c$ is an (interior) marked point, then there exists a
corresponding marked point $z'_c$ on  $\Sigma'$. We put
\begin{equation}\label{specialpointdistance}
f_c(h;\text{\bf x}') = \text{\rm dist}(u'(z'_c),hu(z_c))^2.
\end{equation}
Here $\text{\rm dist}$ is a Riemannian metric on $X$ which is invariant
under the $T^n$ action.
\par
Next, we consider the case when $z_c$ is a singular point.
If this point remains singular in $\text{\bf v}'$,  there
exists a corresponding singular point $z'_c$ on  $\Sigma'$.
In that case we define $f_c(h;\text{\bf x}')$ by (\ref{specialpointdistance}).
\par
To study the case when $z_c$ will not remain a singular point in $\Sigma'$, 
we need a digression.
\par
We divide the bordered curve $\Sigma$ 
equipped with the interior marked points $\vec z^+ \cup \vec w^+$ to disk and sphere components and 
obtain $\text{\bf v}_{0} \in \mathcal M_{0;\ell(0)+m(0)+n(0)}$, 
$\text{\bf v}_{e} \in \mathcal M_{\ell(e)+m(e)+n(e)}$, $e=1,\dots,\frak e$.
(We regard all the special points as marked points.) 
Here $\ell(0) = \# \vec z^+\cap \Sigma_0$,  $m(0) = \# \vec w^+\cap \Sigma_0$
and $n(0)$ is the number of singular points on $\Sigma_0$.
We define $\ell(e)$, $m(e)$, $n(e)$ in the same way.
\par
We take the universal families of the moduli space $\mathcal M_{0;\ell(0)+m(0)+n(0)}$ 
of marked disks (with possibly sphere bubbles) and of the moduli spaces $\mathcal M_{\ell(e)+m(e)+n(e)}$ 
of genus 0 marked curves without boundary, in neighborhoods $\mathcal V({0})$ of $\text{\bf v}'_{0}$
and $\mathcal V(e)$ of $\text{\bf v}'_{e}$, respectively.
(We note that automorphism groups of $\text{\bf v}'_{0}$, $\text{\bf v}'_{e}$
are trivial since they are of genus $0$. So 
$\mathcal V({0})$ and $\mathcal V(e)$ are 
not only orbifolds but also manifolds.)
We also fix (families of) coordinates around each marked point.
Using those data, we obtain a map
\begin{equation}\label{Psiclosedmap}
\Psi : D^2(\epsilon)^{\mathcal S} \times \mathcal V({0}) \times \prod_{e =1}^{\frak e} \mathcal V(e) \to \mathcal M_{0,\ell+m}
\end{equation}
where ${\mathcal S} $ is the set of singular points.
The map $\Psi$ is an isomorphism onto an neighborhood of $(\Sigma,\vec z^+ \cup \vec w^+)$.
Note $D^2(\epsilon)^{{\mathcal S} }$ is the gluing parameter. 
\par
The definition of the map $\Psi$ is as follows.
(See \cite[Section 16]{foootech} for detail. In \cite{foootech} the coordinates $T_c, \theta_c$ are used in place of $\xi_c \in D^2$.  
We need to use these coordinates, when we establish smoothness of coordinate changes of Kuranishi neighborhoods.  )
\par
Let $\text{\bf v}'_{0} \in \mathcal V({0})$ and  $\text{\bf v}'_e \in \mathcal V(e)$.
We will define: 
$$
\Psi((\xi_c)_{c},\text{\bf v}'_{0},(\text{\bf v}'_e)_{e=1,\dots,\frak e}) = \text{\bf v}',
$$
where $c \in\mathcal S$.
Let us consider the component $\xi_c$ of $D^2(\epsilon)^{\mathcal S}$
which corresponds to a singular point $\{y_c\} = \text{\bf v}'_{b(1)} \cap \text{\bf v}'_{b(2)}$.
We use it to glue $\text{\bf v}'_{b(1)}$ and $\text{\bf v}'_{b(2)}$, as follows. 
(Here $b(1), b(2) \in \{0,1,\dots,\frak e\}$.) We take the coordinates $z_1$ and $z_2$ 
of a neighborhood of $y_c$ in $\text{\bf v}'_{b(1)}$, $\text{\bf v}'_{b(2)}$, respectively, which we fixed as a part of the data.
Then to glue $\text{\bf v}'_{b(1)}$ and $\text{\bf v}'_{b(2)}$ we put
the relation $z_1z_2 = \xi_c$.
We thus obtain $\text{\bf v}'$.
In case $\xi_c = 0$ the point $y_c$ remains singular in ${\bf v}'$.
This is the definition of $\Psi$.  
\par
When ${\bf v}$ is equipped with the marked points $\vec{w}^+$,  we denote it by ${\bf v}(\vec{w}^+)$.   
We may take  $\Psi$ so that it is invariant under ${\rm Aut}_+({\bf v}(\vec w^+))$ 
action.
\par
We may suppose that our coordinate neighborhoods of $y_c$ is given by
$\{z_1 \mid \vert z_1 \vert < 1\} \subset \text{\bf v}'_{b(1)}$ and $\{z_2 \mid \vert z_2\vert < 1\} \subset 
\text{\bf v}'_{b(2)}$, respectively.
\par
We put $\xi_c = \exp(-2\pi(10T_c+\sqrt{-1}\theta_c))$.
\par
Now for each singular point $y_c$ we define a closed loop $\gamma_c \subset {\bf v}'$ by
$$
\gamma_c = \{ z_1 \mid \vert z_1 \vert = \exp(-10\pi T_c)\} = \{ z_2 \mid \vert z_2 \vert = \exp(-10\pi T_c)\}.
$$
It comes with a parametrization
\begin{equation}\label{defngammac}
t \mapsto \gamma_c(t) = \exp(-2\pi(5T_c+\sqrt{-1}t))
\end{equation}
using the coordinate $z_1$.
\begin{rem}
Let us consider the situation appearing Remark \ref{pbubblerem}, where $u$ is 
constant on the disk component and the disk component has only one singular point and 
no marked point on it.
In this case $\mathcal V(0)$ is one point. Its element has a nontrivial 
automorphism group $S^1$. This $S^1$ kills a part of the component $D^2(\epsilon)$ as follows.
Let $z_c$ be the (interior) singular point which is on the disk component.
We consider the factor $D^2(\epsilon)$ corresponding to this singular point in the 
left hand side of (\ref{Psiclosedmap}). Then the automorphism group of 
the disk component acts on it and the quotient of the gluing parameter $D^2(\epsilon)$ 
by this $S^1$ action becomes $[0,\epsilon)$.
Namely 
in this case the 
left hand side of (\ref{Psiclosedmap}) 
should be replaced by
$[0,\epsilon) \times D^2(\epsilon)^{\# \mathcal S - 1}  \times \prod_{e =1}^{\frak e} \mathcal V(e)$.
The map $\Psi$ still exists after we replace the domain as above.
(See \cite[Proposition 3.8.27]{fooobook}.)
We can define $\gamma_c$ in the same way.
\end{rem}
We finish the digression.
Let us consider ${\bf x}'$ and $h \in  \frak H(\text{\bf x}')$.
We define $w^{\prime +}(h)_a \in \Sigma'$ such that the following holds.
We put $\vec w^{\prime +}(h) = (w^{\prime +}(h)_a)_{a\in A}$.
\begin{conds}\label{condwh}
\begin{enumerate}
\item
$(\Sigma',\vec z^{\prime +} \cup \vec w^{\prime +}(h))$ 
is $\epsilon_1$-close to $(\Sigma,\vec z^{+} \cup \vec w^{+})$.
\item
$
u'(w^{\prime +}(h)_a) \in hN_{w^{+}_a}.
$
\end{enumerate}
\end{conds}
It is easy to see that such $\vec w^{\prime +}(h)$ exists uniquely.
\par
By Condition \ref{condwh}.1,
there exists $((\xi_c)_{c}(h),\text{\bf v}'_{0}(h),(\text{\bf v}'_e(h))_{e=1,\dots,\frak e})$ 
such that
$$
\Psi((\xi_c)_{c}(h),\text{\bf v}'_{0}(h),(\text{\bf v}'_e(h))_{e=1,\dots,\frak e}) = (\Sigma';\vec z^{\prime +} \cup \vec w^{\prime +}(h)).
$$
Thus to each singular point of $\Sigma$ which is no longer singular in $\Sigma'$, 
we can associate a loop $\gamma_c$ defined by (\ref{defngammac}). We write it as $\gamma_c^h$ to clarify the fact that it depends on $h$.

We now put
\begin{equation}\label{specialpointdistance2}
f_c(h;\text{\bf x}') = \int_{t\in [0,1]}\text{\rm dist}(u'(\gamma^h_c(t)),hu(z_c))^2dt.
\end{equation}
\begin{lem}\label{convexitylemma}
If $\epsilon_1$ is sufficiently small, then (\ref{specialpointdistance2}) is convex.
Moreover it is strictly convex on each connected component of intersection of $S^1_0$ orbits and $\frak H({\bf x}')$.
\end{lem}
\begin{proof}
The element $h \in T^n$ appears in the right hand side of (\ref{specialpointdistance2})  twice,
once in  $hu(z_c)$ and once in $\gamma^h_c$. If we fix $h=h_0$ in $\gamma^h_c$ and consider
$$
h \mapsto \int_{t\in [0,1]}\text{\rm dist}(u'(\gamma^{h_0}_c(t)),hu(z_c))^2dt,
$$
its convexity and strictly convexity on $S^1_0$ orbit follows from Lemma \ref{markedofsymmetry}.1 
and the strict convexity of the square of Riemannian distance function.
Therefore it suffices to show the estimate 
\begin{equation}\label{uprimehexdec}
\vert u'(\gamma^{h}_c(t))\vert_{C^2} \le C e^{-\delta T_c}
\end{equation}
where $\delta >0$ and $\exp(-10 T_c \pi) = \vert\xi_c\vert$. Here we regard $h \mapsto u'(\gamma^{h}_c(t))$ as a function of $h$ and 
takes its $C^2$ norm.
\par
Let us prove (\ref{uprimehexdec}).
\par
When $h$ moves on a connected component of $\frak H(\text{\bf x}')$,  
the family $(\Sigma';\vec z^{\prime +} \cup \vec w^{\prime +}(h))$ becomes a compact family in 
the disk moduli space and is contained in a single stratum of it.
Moreover the map $h \mapsto (\Sigma',\vec z^{\prime +} \cup \vec w^{\prime +}(h))$ has 
bounded $C^2$ norm.
\par
Let $\vec w^{\prime +}$ be in a neighborhood of the given $\vec w^{\prime +}(h)$.
Then using $\vec w^{\prime +}$ in place of $\vec w^{\prime +}(h)$ 
we obtain  a loop $\gamma_c^{\vec w^{\prime +}}$. Thus the map
\begin{equation}
\vec w^{\prime +} \mapsto F_c(\vec w^{\prime +}) = \int_{t\in [0,1]}\text{\rm dist}(u'(\gamma_c^{\vec w^{\prime +}}(t)),hu(z_c))^2dt
\end{equation}
is defined. 
The inequality (\ref{uprimehexdec}) will follow from the estimate
\begin{equation}\label{Fnoestimate}
\vert F_c\vert_{C^2} < Ce^{-\delta T_c}.
\end{equation}
Note the domain of this map $F_c$ is contained in a compact set 
and is a manifold. So the $C^2$ norm above is well defined up to 
bounded ratio.
\par
We finally prove (\ref{Fnoestimate}).
We first observe that it suffices to prove the same estimate while we move two of the points 
$w^{\prime +}_{a(1)}, w^{\prime +}_{a(2)}$ but do not move other points 
in $\vec w^{\prime +}$. 
We fix $a(1), a(2)$ and prove the inequality (\ref{Fnoestimate}) 
when we move $w^{\prime +}_{a(1)}, w^{\prime +}_{a(2)}$ only.
\par
We next recall that the loop  $\gamma_c^{\vec w^{\prime +}}$ depends on the choice of
the coordinates at the singular points of $\Sigma$.
Therefore the map $F$ depends on it.
However we can estimate the difference of $F$ between two 
different choices of those coordinates by using exponential decay of $u'$ and 
\cite[Lemma 16.18]{foootech}.  Therefore we can prove that
if exponential decay estimate (\ref{Fnoestimate}) for $F$ holds for one choice of 
coordinates at singular points, then it also holds 
for any other choices of coordinates at singular points.
\par
We now use Condition \ref{zacondition2}.6.
It implies that each irreducible component of 
$(\Sigma',\vec z^{\prime +}\cup \vec w^{\prime +})$ is still stable when 
we forget two 
marked points  
$w^{\prime +}_{a(1)}, w^{\prime +}_{a(2)}$.
We use this fact to find a family of coordinates at  singular points of 
$\Sigma$, 
so that $F_c(\vec w^{\prime +})$ does not change when we move 
$w^{\prime +}_{a(1)}, w^{\prime +}_{a(2)}$.
The proof of Lemma \ref{convexitylemma} is complete.
\end{proof}
Then we define
\begin{equation}\label{form2827}
f_+(h;\text{\bf x}')
= f(h;\text{\bf x}') + \sum_{c}f_c(h;\text{\bf x}').
\end{equation}
Now by Lemma \ref{convexitylemma}, $f_+(\cdot;\text{\bf x}')$ is
{\it strictly} convex.
We can perform the rest of the construction in the same way
using $f_+$ instead of $f$, and
obtain $\mathcal P_{\text{\bf x}'}$ and $E(\text{\bf x}')$.
We have thus obtained a Kuranishi neighborhood of a  neighborhood of the $T^n \times \frak S_{\ell}$ orbit of $\text{\bf x}$ 
in case Lemma \ref{markedofsymmetry}.1 holds.
\par
\subsection{The case of one disk component III}
\label{subsec:onedisk3}
\par
We next consider the case when Lemma \ref{markedofsymmetry}.2 holds.
Let $G^+_0(\text{\bf x}) \cong S^1$ be the identity component of $G^+(\text{\bf x})$.
(It is obtained by integrating the vector field $\tilde{\frak t}$.
See the proof of Lemma \ref{markedofsymmetry}.)
The map $G^+(\text{\bf x}) \to S^1_0 \subset T^n$ is a finite covering. We define
$$
{\rm Aut}_+^+({\bf x}) 
=
{\rm Ker} (G^+(\text{\bf x}) \to S^1_0),
$$
which is a finite group.
The group $G^+(\text{\bf x})$ is generated by $G^+_0(\text{\bf x})$ and ${\rm Aut}_+^+({\bf x})$ 
since the image of $G^+(\text{\bf x}) \to T^n$ is $S^1_0$.
\begin{lem}
Elements of $G^+_0(\text{\bf x})$ commute with elements of ${\rm Aut}_+^+({\bf x})$.
\end{lem}
\begin{proof}
${\rm Aut}_+^+({\bf x})$ is a normal subgroup and has discrete topology.
On the other hand $G^+_0(\text{\bf x})$ is connected. The lemma follows.
\end{proof}
We take a (single) $G^+(\text{\bf x})$ orbit ${\bf S} \subset \Sigma$,
(${\bf S}$ is a finite disjoint union of circles) so that the following holds.
\begin{conds}\label{Scondition}
\begin{enumerate}
\item
${\bf S}$ is on the disk component.
\item
If $x \in {\bf S}$ then $\{g \in T^n \mid gu(x) = u(x)\} = \{1\}$. 
\item
The map $u$ is an immersion at each point in ${\bf S}$.
\item
${\bf S}$ does not intersect with the union of $\vec z^+$, the set of singular points of $\Sigma$,  
the one dimensional complex $\mathcal C$ as in Lemma \ref{markedofsymmetry}.2 and  $\partial \Sigma$.
\end{enumerate}
\end{conds}
We note that $u$ is not constant on the disk component.
In fact, if $u$ is constant on the disk component, then there should be a
sphere bubble, because $\beta \ne 0$. This implies that there exists a special
point which is sent to  a point of $L$ by the map $u$.
Then Lemma \ref{markedofsymmetry}.1 holds. A contradiction.
We can prove the existence of ${\bf S}$ by using this fact.
\par
We next choose a finite number of additional interior marked points
$$
\vec w^{+} = (w_a^+)_{a\in A}, \quad m = \#A
$$
on $\Sigma$ so that the following conditions are satisfied. 
\begin{conds}\label{zacondition3}
\begin{enumerate}
\item $w_a^+$ is disjoint from $\vec z^+$, the singular point set
and the boundary of $\Sigma$.
\item $(\Sigma,\vec w^{+})$ is stable. Moreover it has no nontrivial automorphism.
\item $u$ is an immersion at $w_a^+$.
\item For any $\frak t \in Lie(T^{n-1}_0)$ there exists $w_a^+$ such that
$\frak t(u(w_a^+))$ is not contained in
$u_*T_{w_a^+}\Sigma$.
\item The subset $\vec w^{+} \subset \Sigma$ is ${\rm Aut}_+^+({\bf x})$ invariant, i.e.,
$\varphi(w^+_a) \in \vec w^{+}$ for all $a \in A$, $\varphi \in {\rm Aut}_+^+({\bf x})$.
\item
In case $G({\bf x})$ is of finite order, $\vec w^{+}$ is $G({\bf x})$ invariant.
\end{enumerate}
\end{conds}
Condition \ref{zacondition3} is almost the same as Condition \ref{zacondition}.
In Item 4 we replace $T^n$ by $T^{n-1}_0$.
In Item 5 we replace $G({\bf  x})$ by ${\rm Aut}_+^+({\bf x})$.
\par
Here we note the following
\begin{lem}\label{GGchigai3}
Suppose Lemma \ref{markedofsymmetry}.2 holds.
Then the following three conditions are equivalent.
\begin{enumerate}
\item
$\# G({\bf x}) = \infty$.
\item
$G({\bf x}) = G^+({\bf x})$. 
\item $\mathcal C$ is empty.
\end{enumerate}
\end{lem}
\begin{proof}
We use the notation in the proof of Lemma \ref{markedofsymmetry}.
We note that  for each $(g,\varphi) \in G({\bf x})$  we have
$\varphi(\mathcal C_0) = \mathcal C_0$. 
Therefore using the fact that $\mathcal C_0$ is a finite set, it is easy to see that 
if $\# G({\bf x}) = \infty$ then $\mathcal C_0 = \emptyset$. 
Thus we proved 1 $\Rightarrow$ 3.
3 $\Rightarrow$ 2 $\Rightarrow$ 1 is obvious.
\end{proof}
For each $w_a^+$ we choose a normal slice $N_{w^+_a}$ of $u$ such that:
\begin{conds}\label{1.82}
\begin{enumerate}
\item $N_{w^+_a}$ is a codimension $2$ smooth submanifold of $X$.
\item $N_{w^+_a}$ is perpendicular to $u(\Sigma)$ at
$u(w_a^+)$.
\item If $\varphi \in {\rm Aut}_+^+({\bf x})$, then
$N_{w^+_a} = N_{\varphi(w^+_a)}$.
\item
In case $G({\bf x})$ is of finite order,  
we have
$gN_{w^+_a} = N_{\varphi(w^+_a)}$ for $(g,\varphi) \in G({\bf x})$.
\end{enumerate}
\end{conds}
This is mostly the same as Condition \ref{1.8}.
\par
We next consider $\text{\bf x}' = (\Sigma',\vec z^{\prime +},u')$ 
and the following condition on it.
Let $\epsilon_1$, $\epsilon_2$ be positive real numbers.
\begin{defn}\label{primeisclose2}
Let $h \in T^n$.
We say that $\text{\bf x}'$ is {\it $(\epsilon_1,\epsilon_2)$-close to $h\text{\bf x}$}
if
there exist $\vec w^{\prime +}$ and a neighborhood $\mathcal U(S(\Sigma))$ of the singular point set of 
$\Sigma$ satisfying \eqref{15apend3} and \eqref{15apend33} with the
following properties.
\begin{enumerate}
\item $\vec w^{\prime +}$ is disjoint from $\vec z^{\prime +}$, the singular point set
and the boundary of $\Sigma^{\prime}$.
\item $\text{\bf v}' = (\Sigma',\vec z^{\prime +},\vec w^{\prime +}) \in \frak U(\text{\bf v})$
and $\text{dist}(\text{\bf v}',\text{\bf v}) < \epsilon_1$. Here we use
an appropriate metric on $\mathcal M_{0;\ell+m}$ to define the distance between $\text{\bf v}'$ and $\text{\bf v}$.
\item
The $C^1$ distance between the two maps $u' \circ i_{\text{\bf v}'}$ and $hu$ is smaller than $\epsilon_1$
on $\Sigma \setminus \mathcal U(S(\Sigma))$.
\item $\text{diam}(u'(\mathcal S)) < \epsilon_1$ if
$\mathcal S$ is any connected  component of
$\Sigma' \setminus  i_{\text{\bf v}'}(\Sigma \setminus \mathcal U(S(\Sigma)))$.
\item
$u'(w^{\prime +}_a) \in hN_{w^+_{a}}$ for each $a$
and
$$
\text{dist}(u'(w^{\prime +}_a),hu(w^+_{a}))  < \epsilon_2.
$$
Here $\vec w^{\prime +} = (w^{\prime +}_{a})_{a\in A}$.
\end{enumerate}
\par
When $\vec w^+,\vec w^{\prime +},h$ are specified in the above condition, we say that
$\text{\bf x}' = (\Sigma',\vec z^{\prime +},u')$ is {\it $(\epsilon_1,\epsilon_2)$-close to} $(h\text{\bf x},\vec w^+)$
with respect to $\vec w^{\prime +}$.
\end{defn}
This is the same as Definition \ref{primeisclose}.  
We put 
$$
{\rm Aut}_+({\bf x})={\rm Aut}_+^+({\bf x}) \cap G({\bf x}).
$$  
Namely, ${\rm Aut}_+({\bf x})$ consists of 
$\psi:\Sigma \to \Sigma$ such that $u \circ \psi = u$ and $\psi$ preserves $\vec{z}$ as a set.  
\begin{lem}\label{Gxindependence3}
Consider $u'$, $\text{\bf v}' = (\Sigma',\vec z^{\prime +},\vec w^{\prime +})$, $h$,
$\vec w^{+}$ given as above and put $\text{\bf x}' = (\Sigma',\vec z^{\prime +},u')$.
Let $\varphi \in
 {\rm Aut}_+({\bf x})$ and $\sigma = (\sigma_1,\sigma_2) = \Psi(\varphi)$.
Suppose $\text{\bf x}' $ is $(\epsilon_1,\epsilon_2)$-close to  $(h\text{\bf x},\vec w^+)$
with respect to $\vec w^{\prime +}$.
We consider $\sigma\text{\bf v}'=(\Sigma',\sigma\vec z^{\prime +},\sigma\vec w^{\prime +})$,
where
$(\sigma z^{\prime +})_{i} = z^{\prime +}_{\sigma_1^{-1}(i)}$,
$(\sigma w^{\prime +})_{a} = w^{\prime +}_{\sigma_2^{-1}(a)}$.
We put $\sigma\text{\bf x}' = (\Sigma',\sigma\vec z^{\prime +},u')$.
\par
Then $\sigma\text{\bf x}' $ is $(\epsilon_1,\epsilon_2)$-close to  $(h\text{\bf x},\vec w^+)$
with respect to $\sigma\vec w^{\prime +}$.
\end{lem}
This is mostly the same as Lemma \ref{Gxindependence} and can be proved in the same way.
\begin{lem}\label{lem4349}
If $\#G({\bf x}) < \infty$, Lemma \ref{w+primeunique}
also holds in our situation.
\end{lem}
The proof is the same as  Lemma \ref{w+primeunique}.
We next consider the case $\#G({\bf x}) = \infty$.
\begin{lem}\label{w+primeunique2}
Let $\mu_1 \in \frak S_{\ell}$ and let $\epsilon_1,\epsilon_2$ be sufficiently small.
Suppose $\text{\bf x}' $ is $(\epsilon_1,\epsilon_2)$-close to  $(h\text{\bf x},\vec w^+)$
with respect to $\vec w_{(1)}^{\prime +}$
and
$\mu_1\text{\bf x}' $ is $(\epsilon_1,\epsilon_2)$-close to  $(h\text{\bf x},\vec w^+)$
with respect to $\vec w_{(2)}^{\prime +}$.
\par
Then there exists
$\psi \in  {\rm Aut}_+(\text{\bf x})$ such  that
$
\vec w_{(2)}^{\prime +} =  \mu_2 \vec w_{(1)}^{\prime +}
$
with $\mu = (\mu_1,\mu_2) = \Psi(\psi)$.
\end{lem}
\begin{proof}
The proof is mostly the same as the proof of Lemma \ref{w+primeunique}.
We can prove Sublemma \ref{sublem4319} with $h=h'$ in the same way.
Namely 
$d_{C^0}(u\circ \psi,u) < \frak o(\epsilon_1)$.
Using the finiteness of ${\rm Aut}_+(\text{\bf x})$ we can prove
there exists $\psi' \in {\rm Aut}_+(\text{\bf x})$ 
such that $d_{C^0}(\psi,\psi') < \frak o(\epsilon_1)$.
We have $u\circ \psi' = u$ and $\psi'(\vec z^+) = \mu_1(\vec z^+)$.
Then by Condition  \ref{zacondition}.5 there exists $\mu_2$ such that
$\psi'(\vec w^+) = \mu_2\vec w^+$.
In the same way as the last step of the proof of Lemma \ref{w+primeunique}
we can show that $\psi'$ has the required properties.
\end{proof}
We define
$\frak H(\text{\bf x}';\vec w) \subset T^n$  by
$$
\frak H(\text{\bf x}';\vec w) =\{h \in T^n \mid \exists  \mu_1 \in \frak S_{\ell},
\text{$\mu_1\text{\bf x}' $ is $(\epsilon_1,\epsilon_2)$-close to  $(h\text{\bf x},\vec w^+)$} \}
$$
and will define 
$
f : \frak H(\text{\bf x}';\vec w)  \to \R_{\ge 0}
$
below.
\par
Suppose $\mu_1\text{\bf x}' $ is $(\epsilon_1,\epsilon_2)$-close to  $(h\text{\bf x},\vec w^+)$
with respect to $\vec w^{\prime +}$.
Let 
$i_{\vec w^{\prime +},\vec w^+} : \Sigma \setminus \mathcal U(S(\Sigma)) \to \Sigma'$ 
be the map $i_{{\bf v}',{\bf v}}$
where ${\bf v} = (\Sigma, \vec z^+ \cup \vec w^+)$ 
and ${\bf v}' = (\Sigma', \mu_1\vec z^{\prime +} \cup \vec w^{\prime +})$.
We use it to define a map
\begin{equation}
\frak I^h_{\vec w^{\prime +},\vec w^+} :
{\bf S} \to \Sigma'
\end{equation}
so that the following conditions are satisfied.
\begin{conds}\label{Smapconds}
Let $x \in {\bf S}$.
\begin{enumerate}
\item
$
{\rm dist}(i_{\vec w^{\prime +},\vec w^+}(x),
\frak I^h_{\vec w^{\prime +},\vec w^+}(x))
< \epsilon_3.
$
\item
The shortest geodesic joining $hu(x)$ to $u'(\frak I^h_{\vec w^{\prime +},\vec w^+}(x))$ is perpendicular to 
$hu(\Sigma)$ at $hu(x)$.
\end{enumerate}
\end{conds}
\begin{lem}
For each sufficiently small $\epsilon_3$ there 
exist $\epsilon_1,\epsilon_2$ such that 
the map $\frak I^h_{\vec w^{\prime +},\vec w^+} $ satisfying Condition 
\ref{Smapconds}
exists uniquely.
\end{lem}
\begin{proof}
Note that $u'\circ i_{\vec w^{\prime +},\vec w^+}$ 
is close to $hu$. Moreover $u$ is an immersion at 
${\bf S}$. The lemma follows easily.
\end{proof}
We now put
\begin{equation}\label{formula4331}
f(h;\text{\bf x}',\vec w^{\prime +},\vec w^+)
= 
\int_{{\bf S}} {\rm dist}(hu(x),
u'(\frak I^h_{\vec w^{\prime +},\vec w^+}(x)))^2
dx.
\end{equation}
(Since ${\bf S}$ is a finite union of  $S^1$ orbits, its has unique invariant measure, 
which we use to integrate.)
\begin{lem}
If $\vec w_{(1)}^{\prime +}$, $\vec w_{(2)}^{\prime +}$ are as in Lemma \ref{w+primeunique2}
then 
$$
f(h;\text{\bf x}',\vec w_{(1)}^{\prime +},\vec w^+)
=
f(h;\text{\bf x}',\vec w_{(2)}^{\prime +},\vec w^+).
$$
\end{lem}
The proof is the same as the proof of Corollary 
\ref{cortow+primeunique}. We will write
$f(h;\text{\bf x}',\vec w^+)$ hereafter.
\begin{lem}\label{Sstrictlyconvex}
If $\epsilon_1$, $\epsilon_2$ are sufficiently small, then 
$f(\cdot;\text{\bf x}',\vec w^+)$ is a convex function.
Moreover it is strictly convex on each of the $T^{n-1}_0$ orbit.
\end{lem}
\begin{proof}
The convexity of $f$ is a consequence of the convexity of the 
Riemannian distance function.
Strict convexity on the $T^{n-1}_0$ orbits follows from the fact that 
$\frak t(u(x)) \notin T_{u(x)}u(\Sigma)$ for $x \in {\bf S}$.
(This follows from Condition \ref{Scondition}.2.)
\end{proof}
We next prove that $f(\cdot;\text{\bf x}',\vec w^+)$ is independent of 
$\vec w^+$. For this purpose we first prove:
\begin{lem}\label{4546independentw}
Let $\vec w^+_{(1)}$, $\vec w^+_{(2)}$ both satisfy Condition 
\ref{zacondition3}. Then for each $\epsilon_1$, $\epsilon_2$ 
there exist $\epsilon'_1$, $\epsilon'_2$ such that the following holds.
\par
If ${\bf x}'$ is  $(\epsilon'_1,\epsilon'_2)$ close to $(h\text{\bf x},\vec w^+_{(1)})$
with respect to $\vec w^{\prime +}_{(1)}$, 
then there exists $\vec w^{\prime +}_{(2)}$ such that 
${\bf x}'$ is  $(\epsilon_1,\epsilon_2)$ close to $(h\text{\bf x},\vec w^+_{(2)})$
with respect to $\vec w^{\prime +}_{(2)}$.
Moreover, we have
\begin{equation}\label{formulain4546}
{\rm dist}(i_{{\bf v}'_{(1)},{\bf v}_{(1)}}(w^{+}_{(2), a}),
w^{\prime +}_{(2), a})
< \frak o(\epsilon_1).
\end{equation}
Here ${\bf v}'_{(1)} = (\Sigma', \vec  z^{\prime +}\cup \vec w^{\prime +}_{(1)})$
and 
${\bf v}_{(1)} = (\Sigma, \vec  z^{+}\cup \vec w^{+}_{(1)})$.
\end{lem}
\begin{proof}
We first assume that $\vec w^{+}_{(1)} \cap \vec w^{+}_{(2)}
= \emptyset$. 
We put
$$
\vec w^{\prime +}_{(3)} = i_{{\bf v}'_{(1)},{\bf v}_{(1)}}(\vec w^{+}_{(2)}).
$$
Then 
$(\Sigma, \vec  z^{+}\cup \vec w^{+}_{(1)} \cup \vec w^{+}_{(2)} )$
is $\frak o(\epsilon'_1)$ close to 
$(\Sigma', \vec  z^{\prime +}\cup \vec w^{\prime +}_{(1)} \cup \vec w^{\prime +}_{(3)})$.
Therefore
$(\Sigma, \vec  z^{+}\cup \vec w^{+}_{(2)} )$
is $\frak o(\epsilon'_1)$ close to 
$(\Sigma', \vec  z^{\prime +}\cup \vec w^{\prime +}_{(3)})$.
\par
On the other hand, we can find $\vec w^{\prime +}_{(2)}$ such that
$$
{\rm dist}(
w^{\prime +}_{(2),a},
w^{\prime +}_{(3),a}) 
< \frak o(\epsilon'_1),
\quad
u'(w^{\prime +}_{(2),a}) \in hN_{w^{+}_{(2),a}}.
$$
In particular, $(\Sigma, \vec  z^{+}\cup \vec w^{+}_{(2)} )$
is $\frak o(\epsilon'_1)$ close to 
$(\Sigma', \vec  z^{\prime +}\cup \vec w^{\prime +}_{(2)})$.
It is easy to see that this $\vec w^{\prime +}_{(2)}$ 
has the required properties.
\par
In the case when 
$\vec w^{+}_{(1)} \cap \vec w^{+}_{(2)}
\ne \emptyset$ we take $\vec w^{+}_{(0)}$
satisfying Condition 
\ref{zacondition3} and such that 
$$
\vec w^{+}_{(1)} \cap \vec w^{+}_{(0)}
=
\vec w^{+}_{(2)} \cap \vec w^{+}_{(0)} 
= \emptyset.
$$
We apply the first half of the proof twice.
Once to  $\vec w^{+}_{(1)}$, $\vec w^{+}_{(0)}$
and once to $\vec w^{+}_{(2)}$, $\vec w^{+}_{(0)}$.
The lemma follows.
\end{proof}
\begin{lem}\label{4546independentw3}
In the situation of Lemma \ref{4546independentw}, 
we have
$$
f(h;\text{\bf x}',\vec w_{(1)}^+)
=
f(h;\text{\bf x}',\vec w_{(2)}^+),
$$
provided $\epsilon_1, \epsilon_2$ are small.
\end{lem}
\begin{proof}
An argument similar to the proof of Sublemma 
\ref{sublem4319} shows that 
(\ref{formulain4546}) implies
$$
{\rm dist}(i_{{\bf v}'_{(1)},{\bf v}_{(1)}}(x),
i_{{\bf v}'_{(2)},{\bf v}_{(2)}}(x))
< \frak o(\epsilon_1)
$$
for $x \in {\bf S}$.  Therefore
$$
\frak I^h_{\vec w^{\prime +}_{(1)},\vec w_{(1)}^+}(x)
=
\frak I^h_{\vec w^{\prime +}_{(2)},\vec w_{(2)}^+}(x)
$$
by Condition \ref{Smapconds}.
The lemma follows.
\end{proof}
\begin{lem}\label{S1equivalenceSS}
If $s \in S^1_0$ then
$$
f(sh;\text{\bf x}',\vec w^+)
=
f(h;\text{\bf x}',\vec w^+).
$$
\end{lem}
\begin{proof}
We consider the following two cases separately.
\par\noindent 
Case 1: $\# G^+({\bf x}) = \infty$:
In this case $G^+({\bf x}) = G({\bf x})$ by Lemma \ref{GGchigai3}.
We
take $\widehat s = (s,\varphi) \in G^+_0({\bf x})
\subset G({\bf x})$. Then clearly
$$
f(sh;\text{\bf x}',\vec w^+)
=
f(h;\text{\bf x}',\varphi^{-1}(\vec w^+)).
$$
Therefore Lemma \ref{4546independentw3} implies Lemma \ref{S1equivalenceSS}.
\par\smallskip
\noindent 
Case 2: $\# G^+({\bf x}) < \infty$:
Let $H \subset S^1_0$ be the image of the projection
$G({\bf x}) \to T^n$.  
We can show that $
f(sh;\text{\bf x}',\vec w^+)
=
f(h;\text{\bf x}',\varphi^{-1}(\vec w^+))
=f(h;{\bf x}', \vec{w}^+)
$
for $(s, \varphi) \in G({\bf x})$, by using Condition \ref{zacondition3}.6 and 
Condition \ref{1.82}.4 and Lemma \ref{4546independentw3}.
\par
Lemma \ref{lem4349} implies the following: if $h, hs \in \frak H(\text{\bf x}',\vec w^+)$ and $s \in S^1_0$,
then $s$ is in a small neighborhood of $H$. 
\par
Since ${\bf S}$ is invariant under $G^+_0({\bf x})$ action and 
the image of the projection $G^+_0({\bf x}) \to T^n$ is 
$S^1_0$, the equality
$
f(sh;\text{\bf x}',\vec w^+)
=
f(h;\text{\bf x}',\vec w^+)
$
holds if $s \in S^1_0$  is in a small neighborhood of the identity
Lemma \ref{S1equivalenceSS} is proved.
\end{proof}
We denote by  $\frak H(\text{\bf x}',\vec w^+)/S^1_0$ the image of the following map
$$
\frak H(\text{\bf x}',\vec w^+) \hookrightarrow T^n \to T^n/S^1_0.$$  
(Since $\frak H(\text{\bf x}',\vec w^+) \subset T^n$ is not invariant under the $S^1_0$-action, 
$\frak H(\text{\bf x}',\vec w^+)/S^1_0$ is not the quotient space by the $S^1_0$-action.)  
By Lemma \ref{S1equivalenceSS}, the function
$f(h;\text{\bf x}',\vec w^+)$ induces
\begin{equation}
\overline f(\cdot,\vec w^+) : \frak H(\text{\bf x}',\vec w^+)/S^1_0 \to \R_{\ge 0}.
\end{equation}
\begin{rem}
Lemma \ref{4546independentw3} implies that $\overline f$ does not depend on 
$\vec w^+$, when it is defined. The domain $\frak H(\text{\bf x}',\vec w^+)/S^1_0$ is 
also independent of $\vec w^+$ in the sense of Lemma \ref{4546independentw}.
\end{rem}
\begin{lem}\label{uniqueminimum}
If $(\epsilon_1, \epsilon_2)$ is small then $\overline f$ attains a minimum at a unique element.
\end{lem}
\begin{proof}
We note that $u(\partial\Sigma)$ consists of a single orbit of $S^1_0$.
Therefore if $h,h' \in \frak H(\text{\bf x}',\vec w^+)$ then 
$$
{\rm dist}([h],[h']) \le \frak o(\epsilon_1)
$$
where $[h], [h'] \in \frak H(\text{\bf x}',\vec w^+)/S^1_0$.
(This is a consequence of the definition of $\frak H(\text{\bf x}',\vec w^+)$.)
On the other hand,
by Lemma \ref{Sstrictlyconvex}, $\overline f$ is strictly convex.
The lemma follows.
\end{proof}
Let $[h_0]$ be the  unique element of $\frak H(\text{\bf x}',\vec w^+)/S^1_0$
where $\overline f$ attains its minimum.
We take its lift $h_0 \in \frak H(\text{\bf x}',\vec w^+)$.
We will use $[h_0]$ to move $E(\text{\bf x})$ to a 
space of sections $E(\text{\bf x}') \subset C^{\infty}(\Sigma',(u')^*TX \otimes \Lambda^{0,1})$.
We assume the following:
\begin{conds}\label{Conds4358}
\begin{enumerate}
\item
The support of each element of $E(\text{\bf x})$ is disjoint from $\mathcal C$, where 
$\mathcal C$ is as in Lemma \ref{markedofsymmetry}.2.
\item
$E(\text{\bf x})$ is invariant under the $G^+({\bf x})$  action.
\end{enumerate}
\end{conds}
Note that Item 2 makes sense because of Item 1.
Condition \ref{immersionsup} follows from Item 1 above.
We can find such $E(\text{\bf x})$ in the same way as the proof of Lemma \ref{app2lema}.
We consider the set 
\begin{equation}\label{setofpsiprime}
\{\vec w^{\prime +} \mid \text{$\mu_1\text{\bf x}' $ is $(\epsilon_1,\epsilon_2)$-close to  $(h_0\text{\bf x},\vec w^+)$
with respect to $\vec w^{\prime +}$ for some $\mu_1$} \}.
\end{equation}
By Lemma \ref{w+primeunique2} this set has order 
$\frak d = \#{\rm Aut}_+(\text{\bf x})$. 
We put  ${\rm Aut}_+(\text{\bf x}) = \{\psi_j \mid j=1,\dots, \frak d\}$,  
$\psi_1 = 1$. 
Then we find
$$
(\ref{setofpsiprime}) = \{ \mu_{2,j}\vec w_1^{\prime +} \mid   \mu_{2,j} = \Psi_2(\psi),
\,\, j=1,\dots, \frak d\}.
$$
Moreover,  
$\mu_{1,j}\text{\bf x}' $ is $(\epsilon_1,\epsilon_2)$-close to  $(h_0\text{\bf x},\vec w^+)$
with respect to $\mu_{2,j}\vec w_1^{\prime +}$ with $\mu_{1,j} = \Psi_1(\psi)$.
We put $\vec w_j^{\prime +} = \mu_{2,j}\vec w_1^{\prime +}$.
\par\smallskip
Then the construction described below is similar to (\ref{Pprelimi}): 
For each $j$ there exists 
$i_j : \Sigma \to \Sigma'$ which sends $\vec z^+$ to $\mu_{1,j}\vec z^{\prime +}$ and 
$\vec w^+$ to $\vec w_j^{\prime +}$. 
We modify $i_{j}$ slightly to define 
\begin{equation}
I_{j} : {\rm Supp}(E(\text{\bf x})) \to \Sigma'
\end{equation}
as follows.
Let $x \in {\rm Supp}(E(\text{\bf x}))$. We require:
\begin{conds}\label{condsIII3}
\begin{enumerate}
\item 
${\rm dist}(i_{j}(x), I_{j}(x)) < \epsilon_3$.
\item
The minimal geodesic $\ell_{x;j}$ joining $h_0u(x)$ to $u'(I_{j}(x))$ is perpendicular 
to $h_0u(\Sigma)$ at $h_0u(x)$.
\end{enumerate}
\end{conds}
We may assume that there exists unique $I_{j}(x)$ satisfying Condition \ref{condsIII3}.
Now we define
\begin{equation}\label{Pprelimi22}
\mathcal P^{\vec w^+}_{\text{\bf x}^{\prime},j} :
E(\text{\bf x}) \to C^{\infty}(\Sigma',u^{\prime *}TX \otimes \Lambda^{0,1})
\end{equation}
in the same way as (\ref{Pprelimi}).
\begin{lem}\label{11533}
$\mathcal P^{\vec w^+}_{\text{\bf x}^{\prime},j} \circ (\psi_{j}^{-1})_*
$
is independent of $j = 1,\dots, \frak d$.
\end{lem}
The proof is the same as that of Lemma \ref{115}.
We define
$$
\mathcal P^{\vec w^+}_{\text{\bf x}^{\prime}} = \mathcal P^{\vec w^+}_{\text{\bf x}^{\prime},j} \circ (\psi_{j}^{-1})_*.
$$
\begin{lem}\label{PisindepedentW}
$\mathcal P^{\vec w^+}_{\text{\bf x}^{\prime}}$ is independent of $\vec w^+$.
\end{lem}
The proof is similar to the proof of Lemma \ref{4546independentw3}.
\begin{lem}\label{parallelsindepe}
$\mathcal P^{\vec w^+}_{\text{\bf x}^{\prime}}$ is independent of the 
representative $h_0$ of $[h_0] \in T^n/S^1_0$.
\end{lem}
\begin{proof}
The proof is similar to the proof of Lemma \ref{S1equivalenceSS}.
Let $s \in S^1_0$. We take $(s,\varphi) \in G_0({\bf x})$.
Using Condition \ref{Conds4358}.2 it is easy to see that
$
\mathcal P^{\vec w^+,sh_0}_{\text{\bf x}^{\prime}} = \mathcal P^{\varphi(\vec w^+),h_0}_{\text{\bf x}^{\prime}}$.
(Here we include $h_0$ or $sh_0$ in the above notation to clarify the choice of the representative of $[h_0]$
we use to define the embedding.)
The lemma then follows from Lemma \ref{PisindepedentW}, 
in Case 1 in the proof of Lemma \ref{S1equivalenceSS}.
\par
In Case 2 in the proof of Lemma \ref{S1equivalenceSS},
invariance under $H$ the image of $G({\bf x}) \to T^n$
follows as above.
Then we only need to consider the case when $s \in S^1_0$
is close to the unit.
We can use Condition \ref{Conds4358}.2 to prove the required 
$s$ invariance.
\end{proof}
\begin{rem}\label{itoIremark}
During the construction of the map (\ref{Pprelimi22}) (that is similar to 
(\ref{Pprelimi})) we first take the map $i_j$ (which is $i_{i,j}$ in the situation of (\ref{Pprelimi}))
and replace it by $I_j$ (which is $I_{i,j}$ in the situation of   (\ref{Pprelimi})) that is close to $i_j$, by using 
Condition \ref{condsIII3}. 
For the purpose of defining  (\ref{Pprelimi}) this step is actually unnecessary 
and we can use $i_{i,j}$ instead of $I_{i,j}$.
\par
On the other hand, here in the definition of (\ref{Pprelimi22}) we need to replace $i_j$ by $I_j$.
In fact, Lemmata \ref{PisindepedentW}, \ref{parallelsindepe} would not hold, if we use $i_j$ in place of $I_j$.
\end{rem}
Hereafter we write $\mathcal P_{\text{\bf x}^{\prime}}$ in place of $\mathcal P^{\vec w^+}_{\text{\bf x}^{\prime}}$.
Now in the same way as Lemma \ref{invhhh} 
we can prove
$$
h_*\mathcal P_{\text{\bf x}^{\prime}} = \mathcal P_{h\text{\bf x}^{\prime}}
$$
where $h\in T^n$.
Let $E({\bf x'})$ be the image of $\mathcal P_{\text{\bf x}^{\prime}}$. 
Then in the same way as Lemma \ref{G(x)equivalence1} we can prove
$$
g_*(E(\text{\bf x}')) = E((g,\varphi)\cdot\text{\bf x}')
$$
for $(g,\varphi) \in G({\bf x})$.
The rest of the argument is the same as the last part of 
Subsection \ref{subsec:onedisk1}.
\par
We have thus completed the construction of our Kuranishi neighborhood
of the $T^n \times \frak S_{\ell}$ orbit of ${\bf x}$,  
in the case when $\text{\bf x}$ has only one disk component.
\par
\subsection{The case of more than one disk components}
\label{subsec:1+}
We finally consider the case when the domain may have $2$ or more disk components.
\par
Let $\text{\bf x} = (\Sigma,\vec z^+,u) \in \mathcal M_{0;\ell}(\beta)$.
We decompose $\Sigma$ into extended disk components:
\begin{equation}\label{decompxxxx}
\Sigma = \bigcup_{\alpha \in \frak A} \Sigma_\alpha.
\end{equation}
Each extended disk component $\Sigma_{\alpha}$ together with the restriction of $u$ and $\vec z^+$
defines an element
\begin{equation}\label{xalpha}
\text{\bf x}_\alpha = (\Sigma_\alpha,\vec z^+_{\alpha},u_{\alpha})
\in \mathcal M_{0;\ell_{\alpha}}(\beta_\alpha).
\end{equation}
We may regard boundary singular points on $\Sigma({\alpha})$ as 
marked points. Then we obtain 
 \begin{equation}\label{xalphakk}
\text{\bf x}^+_\alpha = (\Sigma_\alpha,\vec z_{\alpha}.\vec z^+_{\alpha},u_{\alpha})
\in \mathcal M_{k_{\alpha};\ell_{\alpha}}^{\rm main}(\beta_\alpha)/\Z_{k_{\alpha}},
\end{equation}
where $\Z_{k_{\alpha}}$ action is defined by the cyclic permutation of the boundary singular points.
Each of $\text{\bf x}_\alpha$ satisfies one of the
following 4 conditions.
\par
\begin{enumerate}
\item[O.]
$\beta_\alpha =0$. 
\item[A.] $\beta_\alpha \ne 0$, $\text{\bf x}_\alpha$ satisfies Assumption \ref{trivialisotopystrong}.
\item[B.] $\beta_\alpha \ne 0$, $\text{\bf x}_\alpha$ satisfies the conclusion of Lemma \ref{markedofsymmetry}.1.
\item[C.] $\beta_\alpha \ne 0$, $\text{\bf x}_\alpha$ satisfies the conclusion of Lemma \ref{markedofsymmetry}.2.
\end{enumerate}
\begin{rem}
In case $\beta_\alpha =0$ and $\ell_{\alpha} \le 2$ the right hand side of 
(\ref{xalphakk}) is not defined. 
This does not matter because, in case O, we do not put obstruction bundle on the component $\Sigma_\alpha$, 
and hence we do not need to study such a component.
\end{rem}
We say that $\text{\bf x}_\alpha$ is {\it of Type A, B, C} if it satisfies A,B,C above,
respectively.
Sometimes we say $\alpha$ is of Type A,B,C also.
\par
In case $\alpha$ is of Type C, we take and fix a splitting of $T^n$ as in Lemma \ref{Tnsplits}.
We write this splitting as $S^1_{\alpha} \times T^{n-1}_{\alpha}$.
(This splitting depends on $\alpha$ so we include
$\alpha$ in the notation.)
We also choose ${\bf S}_{\alpha}$ satisfying 
Condtion \ref{Scondition}.
\par
For each $\text{\bf x}_{\alpha} = (\Sigma_{\alpha},\vec z^{+}_{\alpha},u_{\alpha})$ of Type A,B,C, we take additional marked points
$w^+_{\alpha}$ so that Condition \ref{zacondition}, \ref{zacondition2}, or \ref{zacondition3} is satisfied,
respectively.
We also choose $N_{w}$ as in Condition \ref{1.8} or \ref{1.82}.
We put 
$$
\text{\bf x}_{\alpha}(\vec w^+) = (\Sigma_{\alpha},\vec z^{+}_{\alpha}\cup \vec w_{\alpha}^+,u_{\alpha}), \quad
\text{\bf v}_{\alpha}(\vec w^+) = (\Sigma_{\alpha},\vec z^{+}_{\alpha}\cup \vec w_{\alpha}^+). 
$$
Let $\vec w^+ = \bigcup_{\alpha} \vec w_{\alpha}^+$ and
$$
\text{\bf x}(\vec w^+) = (\Sigma, \vec z^+ \cup \vec w^+,u), \quad
\text{\bf v}(\vec w^+) = (\Sigma, \vec z^+ \cup \vec w^+).
$$
\par
We defined the group $G(\text{\bf x})$ as in Definition \ref{symmetryG}.
We note that in our situation where the domain $\Sigma$ has
at least two disk components, $G(\text{\bf x})$ is necessarily a finite group.
This is because the isotropy group of a boundary singular point is trivial.
(See the proof of Lemma \ref{markedofsymmetry}.)
\begin{rem}\label{rem146}
\begin{enumerate}
\item
An element $\text{\bf g} \in G(\text{\bf x})$ may exchange components $\text{\bf x}_{\alpha}$.
We choose $\vec w^{+}_{\alpha}$ (the case of Type A or B component),
$\text{\bf S}_{\alpha}$ (the case of Type C component)  so that it is preserved under
this action. In fact, we are proving Proposition \ref{eqcyckuramain} inductively.
(Namely by induction on $E(\beta)$.)
In other words, we are making the choice of $\vec w^{+}_{\alpha}$ or $\text{\bf S}_{\alpha}$
inductively on $E(\beta_{\alpha})$.
In our situation for which $\beta$ is decomposed into $\beta = \sum \beta_{\alpha}$,
the choice for each $\beta_{\alpha}$ had already been made.
Of course, it depends only on the isomorphism classes. So the action
of $\text{\bf g}$ automatically preserves it.
\item
The situation is different in the case of the additional marked points $\vec w^{+}_{\alpha}$ on  
Type C components. Actually the group $G^+({\bf x}_{\alpha})$ is 1 dimensional in that case
and we can not find $\vec w^{+}_{\alpha}$ which is $G^+({\bf x}_{\alpha})$ invariant.
In Subsection \ref{subsec:onedisk3}, we only assume that 
it is invariant under ${\rm Aut}_+^+({\bf x}_{\alpha})$ action, 
in the case of Type C component.
Actually we used $\vec w^{+}_{\alpha}$ during the construction but 
finally proved that it is independent of $\vec w^{+}_{\alpha}$.
(See Lemmata \ref{4546independentw3}, \ref{PisindepedentW}.)
Here we choose $\vec w^{+}_{\alpha}$ so that its union $\vec z^+ \cup \vec w^+$
is invariant under the $G({\bf x})$ action.   
\end{enumerate}
\end{rem}
\par
We consider a neighborhood $\frak U(\text{\bf v}(\vec w^+))$ of $\text{\bf v}(\vec w^+)$.
Let $\text{\bf v}' \in \frak U(\text{\bf v}(w^+))$
and
$\text{\bf v}'  = (\Sigma',\vec z^{\prime +},\vec w^{\prime +})$.
Then we have a neighborhood $\mathcal  U(S(\Sigma))$ 
of the singular point set of $\Sigma$ and
$$
i_{\text{\bf v}'} : \Sigma \setminus \mathcal  U(S(\Sigma)) \to \Sigma'
$$
such that the following holds:
\begin{conds}\label{ivcompatibility}
\begin{enumerate}
\item $i_{\text{\bf v}'}(z_i^+) = z_{i}^{\prime +}$,
$i_{\text{\bf v}'}(w_a^+) = w_{a}^{\prime +}$.
\item $i_{\text{\bf v}'}$ is a diffeomorphism to its image.
\item $i_{\text{\bf v}'}$ is biholomorphic on 
the support of $E(\text{\bf x}_{\alpha})$ for any $\alpha$.
\item If $\sigma = (\sigma_1,\sigma_2) \in \Psi(\psi)$,
$\psi \in \text{\rm Aut}_+(\text{\bf v}(\vec w^+))$, then
$$
i_{\sigma\text{\bf v}'} = i_{\text{\bf v}'}\circ \psi^{-1} :
 \Sigma \setminus \mathcal U(S(\Sigma)) \to \Sigma'.
 $$
\item When $\Sigma$ has 2 or more disk components, $i_{\text{\bf v}'}$ is compatible with
the maps inductively defined in each of the disk components.
\item
$i_{\text{\bf v}'}$ depends continuously on $\text{\bf v}'$ in $C^{\infty}$
topology.
\end{enumerate}
\end{conds}
\begin{rem}
Our construction of the Kuranishi neighborhoods given below is designed so that it does not 
depend on the choice of $i_{{\bf v}'}$ much. Namely if $i'_{{\bf v}'}$ is another choice 
such that $d(i_{{\bf v}'}(x),i'_{{\bf v}'}(x))$ is small, then the Kuranishi neighborhoods 
obtained by using $i_{{\bf v}'}$ and by using $i'_{{\bf v}'}$ coincide.
\end{rem}
\par
Let $\epsilon_1,\epsilon_2$ be small positive numbers.
We take $(\Sigma',\vec z^{\prime +})$ and $u' : (\Sigma',\partial \Sigma')
\to (X,L)$.
\begin{defn}\label{primeisclosemodmod}
For $\text{\bf x}' = (\Sigma',\vec z^{\prime +},u')$ and $\vec w^+$, we say that $\text{\bf x}' = (\Sigma',\vec z^{\prime +},u')$ is
{\it $(\epsilon_1,\epsilon_2)$-close to} $(h\text{\bf x},\vec w^+)$
with respect to $\vec w^{\prime +}$ if the following holds.
\begin{enumerate}
\item
$\text{\bf v}' = (\Sigma',\vec z^{\prime +},\vec w^{\prime +}) \in \frak U(\text{\bf v}(\vec w^+))$
and $\text{dist}(\text{\bf v}',\text{\bf v}(\vec w^+)) < \epsilon_1$. Here we use
an appropriate metric on $\mathcal M_{0;\ell+m}$ to define the distance between $\text{\bf v}'$ and $\text{\bf v}(\vec w^+)$.
\item
The $C^1$ distance between two maps $u' \circ i_{\text{\bf v}'}$ and $hu$ is smaller than $\epsilon_1$
on $\Sigma \setminus \mathcal U(S(\Sigma))$.
\item$\text{diam}(u'(\mathcal S)) < \epsilon_1$ if
$\mathcal S$ is any connected  component of 
$\Sigma' \setminus  i_{\text{\bf v}'}(\Sigma \setminus \mathcal U(S(\Sigma)))$.
\item  $u'(w^{\prime +}_a) \in hN_{w^+_{a}}$ for each $a$
and
$$
\text{dist}(u'(w^{\prime +}_a),hu(w^+_{a}))  < \epsilon_2.
$$
Here $\vec w^{\prime +} _a = (w^{\prime +}_a)_{a\in A}$.
\end{enumerate}
\end{defn}
We note that this definition is exactly the same as Definition \ref{primeisclose}.
\par
\begin{lem}\label{w+primeunique33}
Let $\mu_1 \in \frak S_{\ell}$ and let $\epsilon_1,\epsilon_2$ be sufficiently small.
Suppose $\text{\bf x}' $ is $(\epsilon_1,\epsilon_2)$-close to  $(h\text{\bf x},\vec w^+)$
with respect to $\vec w_{(1)}^{\prime +}$
and
$\mu_1\text{\bf x}' $ is $(\epsilon_1,\epsilon_2)$-close to  $(h'\text{\bf x},\vec w^+)$
with respect to $\vec w_{(2)}^{\prime +}$.
\par
Then there exists
$(g,\psi) \in G(\text{\bf x})$ such  that
$
\vec w_{(2)}^{\prime +}$ is close to  $ \mu_2 \vec w_{(1)}^{\prime +}
$
with $\mu = (\mu_1,\mu_2) = \Psi(\psi)$
and 
$$
{\rm dist}(h',gh) < \frak o(\epsilon_1).
$$
\par
If $h=h'$ in addition, then $g=1$ and $\vec w^{\prime +}_{(2)}=\mu_2 \vec w^{\prime +}_{(1)}$.
\end{lem}
This is the same as Lemma \ref{w+primeunique}. The proof is also the same.
\par
We can proceed in the same way as in Subsection \ref{subsec:onedisk1},
replacing Definition \ref{primeisclose} by
Definition \ref{primeisclosemodmod}.
However, the resulting Kuranishi neighborhood will not be disk component-wise. (See Definition \ref{cpwkura}.)
(See however Remark \ref{colermethod}.)
To obtain a disk component-wise Kuranishi neighborhood we need to consider the
case where the element $h \in T^n$ above may vary depending on the disk component to
which $h$ is applied.
We will describe this point in this subsection.
\par
We take the decomposition
\begin{equation}
\Sigma' = \bigcup_{\widehat\alpha\in {\frak A}'} \Sigma'_{\widehat\alpha}
\end{equation}
into the extended disk components.
We consider the situation of Definition \ref{primeisclosemodmod}.
So we fix $h$ and $\vec w^{\prime +}$ for a while.
\par
We put
\begin{equation}
{\frak A}(\widehat\alpha) = 
\{
\alpha \in {\frak A}\mid i_{{\bf v}'}(\Sigma_{\alpha} \setminus \mathcal U(S(\Sigma)))  \subset \Sigma'_{\widehat\alpha}.
\}
\end{equation}
We note that 
$$
{\frak A} = \bigcup_{\widehat\alpha\in {\frak A}'} {\frak A}(\widehat\alpha),
$$
and the right hand side is the disjoint union. 
Let $\alpha \in \frak A(\widehat\alpha)$.
\begin{defn}
$$
\vec z^{\prime +}_{\alpha} = i_{{\bf v}'}(\vec z^+_{\alpha}) \subset \Sigma'_{\widehat\alpha} \cap \vec z^{\prime +},
\quad
\vec w^{\prime +}_{\alpha} = i_{{\bf v}'}(\vec w^+_{\alpha}) \subset \Sigma'_{\widehat\alpha} \cap \vec w^{\prime +}.
 $$
\par
We denote by $\Sigma'_{\widehat\alpha;\alpha}$ the set $\Sigma'_{\widehat\alpha}$
minus all the sphere bubbles which are not on  $ i_{{\bf v}'}(\Sigma_{\alpha})$. 
We put
$$
{\bf x}'_{\widehat \alpha, \alpha}(\vec w^{\prime +}) = (\Sigma'_{\widehat \alpha, \alpha}, \vec z^{\prime +}_{\alpha}  \cup \vec w^{\prime +}_{\alpha} ,u'),
\quad
{\bf v}'_{\widehat \alpha, \alpha}(\vec w^{\prime +}) = (\Sigma'_{\widehat \alpha, \alpha}, \vec z^{\prime +}_{\alpha}  \cup \vec w^{\prime +}_{\alpha} ).
$$
\end{defn}
\begin{rem}
Here we consider only the marked points on $\Sigma'_{\widehat\alpha}$ which 
correspond to the marked points on $\Sigma_{\alpha}$.
There may be other elements $\alpha'$ of $\frak A(\widehat\alpha)$. However we do not put the marked points corresponding 
to those marked points on $\Sigma_{\alpha'}$.
\end{rem}
\begin{lem}\label{ifremvethenclose}
${\rm dist}({\bf v}'_{\widehat \alpha, \alpha}(\vec w^{\prime +}),
{\bf v}_{\alpha}(\vec w^{+})) < \frak o(\epsilon_1)$.
\end{lem}
\begin{proof}
$\Sigma_{\widehat \alpha, \alpha}$ is $\Sigma_{\alpha}$ with disk components glued.
We do not put marked points to those disk components.
Hence the lemma. (See the proof of Lemma \ref{idonotchangemuch} below.)
\end{proof}
We next compare two maps $i_{{\bf v}'_{\widehat \alpha, \alpha}(\vec w^{\prime +}_{\widehat\alpha})}$ and $i_{{\bf v}'(\vec w^{\prime +})}$ on $\Sigma_{\alpha}$. (See Lemma \ref{idonotchangemuch}.)
For this purpose, we make certain digression: 
We decompose
$$
\Sigma_{\alpha} = \Sigma_{\alpha,0} \cup \bigcup_{a \in A_{\alpha}} \Sigma_{\alpha,a}.
$$
Here $\Sigma_{\alpha,0}$ is a disk and $\Sigma_{\alpha,a}, a \in A_{\alpha}$,  are spheres in the trees that are rooted on $\Sigma_{\alpha,0}$.  
The component $\Sigma_{\alpha,0}$ (resp. $\Sigma_{\alpha,a}$) together with special points on them defines an  element 
$
{\bf v}^+_{\alpha,0}(\vec w^{+}) \in \mathcal M_{k_{\alpha};\ell_{\alpha}(0)+m_{\alpha}(0) + n_{\alpha}(0)}
$
(resp. $
{\bf v}_{\alpha,a}(\vec w^{+}) \in \mathcal M_{\ell_{\alpha}(a)+m_{\alpha}(a) + n_{\alpha}(a)}
$).
(Here $k_{\alpha}$ is the number of boundary singular points on $\Sigma(\alpha)$.
$\ell_{\alpha}(0) = \# \vec z^+ \cap \Sigma_{\alpha,0}$, 
$m_{\alpha}(0) = \# \vec w^+ \cap \Sigma_{\alpha,0}$ 
and $n_{\alpha}(0)$ is the number of interior singular points on $\Sigma_{\alpha,0}$.)
Let $\mathcal V(\alpha,0)$ (resp. $\mathcal V(\alpha,a)$) be a 
neighborhood of ${\bf v}_{\alpha,0}(\vec w^{+})$ (resp. ${\bf v}_{\alpha,a}(\vec w^{+})$) 
in $\mathcal M_{k_{\alpha};\ell_{\alpha}(0)+m_{\alpha}(0) + n_{\alpha}(0)}$ (resp. 
$\mathcal M_{\ell_{\alpha}(a)+m_{\alpha}(a) + n_{\alpha}(a)}$.)
We fix coordinates at each singular point.
\par
Let $\mathcal S_{c}$ (resp. $\mathcal S_d$) be the set of the interior (resp. boundary) 
singular points of  $\Sigma$.
\par
We will define a map
\begin{equation}\label{Psiborered}
\Psi : [0,\epsilon)^{\mathcal S_d} \times D^2(\epsilon)^{\mathcal S_c}
\times \prod_{\alpha}  \mathcal V(\alpha,0)
\times \prod_{\alpha,a} \mathcal M_{\ell_{\alpha}(a)+m_{\alpha}(a) + n_{\alpha}(a)}
\to
\mathcal M_{0,\ell+m}.
\end{equation}
The construction of (\ref{Psiborered}) is similar to the construction of the map
(\ref{Psiclosedmap}).
So we describe the way how we glue the boundary singular points only.
Let $\{y_c\} \in \Sigma_{\alpha_1} \cap \Sigma_{\alpha_2}$ be a boundary singular
point.  We consider the coordinates $z_i$ $(i=1,2)$ of $\Sigma_{\alpha_i}$ in a neighborhood 
of $y_c$. Here a coordinate neighborhood can be identified with 
$\{ z_i \mid \vert z_i\vert <1, \,\, {\rm Im}\, z_i \ge 0\}$.
Let $\zeta_c \in [0,\epsilon)$ be the corresponding component of the factor 
$[0,\epsilon)^{\mathcal S_d}$. 
We glue $\Sigma_{\alpha_1}$ and $\Sigma_{\alpha_2}$ by 
identifying $z_1$ and $z_2$ if $z_1z_2 = - \zeta_c$.
The other part of the construction of (\ref{Psiborered}) is the same as
(\ref{Psiclosedmap}).
We finish digression. 
\par
Now we go back to the situation of Lemma \ref{ifremvethenclose}.
Among the factors of the left hand side of (\ref{Psiborered}), we only take those for 
$\alpha \in \frak A(\widehat\alpha)$. 
Among the factors of $[0,\epsilon)^{\mathcal S_d}$, where each $[0,\epsilon)$ corresponds to the boundary singular point,  
we only take those corresponding to the boundary singular points $y_c$'s  such that both of the disk components 
containing $y_c$ are in $\frak A(\widehat\alpha)$.
We then obtain
\begin{equation}\label{Psiboreredalpha}
\aligned
\Psi_{\widehat\alpha} : &[0,\epsilon)^{\mathcal S_d(\widehat\alpha)} \times D^2(\epsilon)^{\mathcal S_c(\widehat\alpha)}
\\ &\times \prod_{\alpha\in \frak A(\widehat\alpha)}  \mathcal V(\alpha,0)
\times \prod_{\alpha\in \frak A(\widehat\alpha) \atop a\in A_{\alpha}} \mathcal M_{\ell_{\alpha}(a)+m_{\alpha}(a) + n_{\alpha}(a)}
\to
\mathcal M_{0,\ell(\widehat\alpha)+m(\widehat\alpha)}.
\endaligned
\end{equation}
Suppose
$$
{\bf v}'_{\widehat\alpha}(\vec w^{\prime +})
=
\Psi_{\widehat\alpha}((\zeta_c)_{c\in {\mathcal S_d(\widehat\alpha)}},(\xi_c)_{c\in {\mathcal S_c(\widehat\alpha)}},(\text{\bf v}'_{\alpha,0}(\vec w^{\prime +}))_{\alpha\in \frak A(\widehat\alpha)},
(\text{\bf v}'_{\alpha,a}(\vec w^{\prime +}))_{\alpha\in \frak A(\widehat\alpha),a}).
$$
Here, to simplify the notation we denote just by $c$ the boundary singular 
point $y_c$.  
Note that by forgetting the marked points of ${\bf v}'_{\widehat\alpha}(\vec w^{\prime +})$ other than those coming from $\Sigma_{\alpha}$, 
we obtain $ {\bf v}'_{\widehat \alpha, \alpha}(\vec w^{\prime +})$.
\begin{lem}\label{idonotchangemuch}
If $\vert\zeta_c\vert < \epsilon_4$ for all $c \in \mathcal S_d(\widehat\alpha)$ 
and $\vert\xi_c\vert < \epsilon_4$ for all $c \in \mathcal S_c(\widehat\alpha)$,
then we have the following. 
There exists a neighborhood $U_{\alpha}$ of the set of boundary singular points in $\Sigma_{\alpha}$ 
such that
\begin{equation}\label{4563mainestimate}
{\rm dist}(i_{{\bf v}'(\vec w^{\prime +})}(x),i_{{\bf v}'_{\widehat \alpha, \alpha}(\vec w^{\prime +})}(x)) < \frak o(\epsilon_4)
\end{equation}
holds for $x \in  \Sigma_{\alpha}\setminus \mathcal U(S(\Sigma_{\alpha})) \setminus U_{\alpha}$.
\end{lem}
Note $i_{{\bf v}'(\vec w^{\prime +})} = i_{{\bf v}'_{\widehat\alpha}(\vec w^{\prime +})}$ on $\Sigma_{\alpha}\setminus \mathcal U(S(\Sigma_{\alpha})) \setminus U_{\alpha}$.
\begin{proof}
The map $\Psi_{\widehat\alpha}$ depends on the choice of
coordinates around the singular points.
However, we can estimate the difference between the two choices by using 
\cite[Lemma 16.18]{foootech}. Therefore we can prove that 
if (\ref{4563mainestimate}) holds for some choice of coordinates around the singular points
then it holds also for other choices of coordinates around the singular points.
Therefore to prove the lemma it suffices to find a choice of the 
coordinates around the singular points so that 
$i_{{\bf v}'(\vec w^{\prime +})}(x) = i_{{\bf v}'_{\widehat \alpha, \alpha}(\vec w^{\prime +})}(x)$.
\par
We give such a choice below.
Let $y_c \in \Sigma_{\alpha,a}$ be a singular point contained in a sphere 
component $\Sigma_{\alpha,a}$.
We identify 
$\Sigma_{\alpha,a} = \C \cup {\infty}$ and $y_c = 0$. Then 
we take the standard coordinate of $\C$.
\par
Let $y_c \in \Sigma_{\alpha,0}$ be an interior singular point
contained in the disk component $\Sigma_{\alpha,0}$.
We identify $\Sigma_{\alpha,0}$ with 
$D^2(2) = \{ z\in \C \mid \vert z\vert < 2\}$ and $y_c = 0$.
Then 
we take the restriction of the standard coordinate of $\C$ to $D^2(2)$.
\par
Let $y_c \in \Sigma_{\alpha,0}$ be a boundary singular point.
We identify
$\Sigma_{\alpha,0}$ with the union of upper half plane  $\frak h= \{z \in \C \mid {\rm Im}\,z \ge 0\}$  and ${\infty}$
and $y_c = 0$.
Then we take the restriction of the standard coordinate of $\C$ to 
the upper half plane.  
\par
We claim that if we take these choices of coordinates at the 
singular points, then $i_{{\bf v}'(\vec w^{\prime +})}(x) = i_{{\bf v}'_{\widehat \alpha, \alpha}(\vec w^{\prime +})}(x)$.
We prove it in the following case 
for simplicity.
Suppose $\frak A(\widehat\alpha) = \{\alpha, \alpha'\}$.
${\bf v}_{\alpha}(\vec w^{+}) = (\frak h \cup \{\infty\},z_1,z_2)$, ${\bf v}_{\alpha'}(\vec w^{+}) = (\frak h \cup \{\infty\},w_1,w_2)$
 (disks without sphere component and with 
two marked points).
Suppose they are glued at $0$. 
The coordinate we took is $z$ on $\Sigma(\alpha)$ and $w$ on $\Sigma(\alpha')$.
We identify them by $zw = \zeta$ and 
obtain $\Sigma(\widehat\alpha)$.
It comes with four marked points $z_1,z_2,w_1,w_2$. 
This is ${\bf v}'_{\widehat\alpha}(\vec w^{\prime +})$.
We forget $w_1,w_2$ and obtain ${\bf v}'_{\widehat\alpha,\alpha}(\vec w^{\prime +})$.
\par
We take $U_{\alpha} = \{z \in \frak h \mid \vert z\vert < \epsilon_3\}$.
Then the equality
$i_{{\bf v}'_{\alpha}(\vec w^{\prime +})}(x) = i_{{\bf v}'_{\widehat \alpha, \alpha}(\vec w^{\prime +})}(x)$ is obvious in this case.
The general case is similar. The proof of Lemma \ref{idonotchangemuch} is complete.
\end{proof}
Now, for any $\alpha \in \frak A(\widehat\alpha)$ we will define a map 
\begin{equation}\label{Palphamanycomp}
\mathcal P_{\alpha} : E({\bf x}_{\alpha}) \to C^{\infty}(\Sigma'(\widehat\alpha),(u')^*TX \otimes \Lambda^{0,1}),
\end{equation}
according to the type of $\alpha$ as follows:
We note that we still fix $h$ and $\vec w^{\prime +}$.
\par\smallskip
\noindent{\bf Case A: $\alpha$ is of Type A.}
\par
Let $h_{\alpha} \in T^n$ be an element which is $\epsilon_4$ close to $h$.
For each $a \in A_{\alpha}$, we take 
additional interior marked points 
$\vec w^{\prime +}_{\alpha}(h_{\alpha})
:= (w^{\prime +}_{\alpha}(h_{\alpha})_b)_{b}
\subset \Sigma'(\widehat\alpha)$, 
which depend on $h_{\alpha}$ and satisfy the following conditions. 
\begin{conds}\label{wprimeha}
\begin{enumerate}
\item
${\rm dist}(w^{\prime +}_{\alpha}(h_{\alpha})_b,w^{\prime +}_{\alpha,b}) < \epsilon_3$.
\item
$u'(w^{\prime +}_{\alpha,b}) \in h_{\alpha}N_{h_{\alpha}u(w^{+}_{\alpha,b})}$.
\end{enumerate}
\end{conds}
We may take $\epsilon_3$ smaller than a number depending only on $X,L$, and $N_{h_{\alpha}u(w^{+}_{\alpha,b})}$
such that if $\epsilon_1,\epsilon_2,\epsilon_4$ is smaller than a number depending $\epsilon_3$ then
there exists a unique $\vec w^{\prime +}_{\alpha}(h_{\alpha})$ satisfying Condition \ref{wprimeha}.
\begin{rem}
Since
$i_{{\bf v}'(w^{\prime +})}(x) \ne i_{{\bf v}'_{\widehat \alpha, \alpha}(w^{\prime +})}(x)$ in general, 
$w^{\prime +}_{\alpha}(h_{\alpha})_b = w^{\prime +}_{\alpha,b}$ may not hold.
\end{rem}
\par
We define
\begin{equation}\label{fincase2components}
f_{\alpha}(h_{\alpha}) = \sum_{b \in  A_{\alpha}} {\rm dist}(u'(w^{\prime +}_{\alpha,b}),h_{\alpha}u(w^{+}_{\alpha,b}))^2.
\end{equation}
The definition is very much similar to the definition of $f$ in \eqref{deffbyave}.
Indeed, this coincides with one in case $\Sigma'(\alpha)$ has only one disk component.
\begin{lem}
$f_{\alpha}$ does not change if we change $u'$ on $\Sigma'(\widehat\alpha')$ for 
$\widehat\alpha' \ne \widehat\alpha$.
\end{lem}
This is obvious from the definition.
\begin{lem}
$f_{\alpha}$ attains its minimum at a unique point, denoted by $h^0_{\alpha}$, in a neighborhood of $h$.
\end{lem}
This is a consequence of strict convexity of $f_{\alpha}$.
\par
We next define a map 
\begin{equation}
I_{\widehat\alpha,\alpha} :  {\rm Supp}(E(\text{\bf x}_{\alpha})) \to \Sigma'(\widehat\alpha)
\end{equation}
as follows. 
Let $x \in {\rm Supp}(E(\text{\bf x}_{\alpha}))$ we require:
\begin{conds}\label{condsIato}
\begin{enumerate}
\item 
${\rm dist}(i_{{\bf v}'_{\widehat\alpha}(\vec w^{\prime +})}(x), I_{\widehat\alpha,\alpha}(x)) < \epsilon_3$.
\item
The minimal geodesic $\ell_{x;\widehat\alpha,\alpha}$ joining $h^0_{\alpha}u(x)$ to $u'(I_{\widehat\alpha,\alpha}(x))$ is perpendicular 
to $h_\alpha^0 u(\Sigma(\alpha))$ at $h_\alpha^0 u(x)$.
\end{enumerate}
\end{conds}
We may choose $\epsilon_4$ small so that $U_{\alpha} \cap {\rm Supp}(E(\text{\bf x}_{\alpha})) 
= \emptyset$. (Here $U_{\alpha}$ is as in Lemma \ref{idonotchangemuch}.)
Then Lemma \ref{idonotchangemuch} and Definition \ref{primeisclosemodmod} imply the 
unique existence of such $I_{\widehat\alpha,\alpha}$.
\par
Using $I_{\widehat\alpha,\alpha}$ and the parallel transport along the geodesic $\ell_{x;\widehat\alpha,\alpha}$, 
we can define (\ref{Palphamanycomp}) in the same way as 
(\ref{Pprelimi}).
\par\smallskip
\noindent{\bf Case B: $\alpha$ is of Type B.}
\par
We define $f_{\alpha}$ in the same way as 
(\ref{fincase2components}).
This function is convex in the $T^{n-1}_{\alpha}$ orbit.
We again use special point of $\Sigma({\alpha})$
to obtain a convex function in the same way as in 
Subsection \ref{subsec:onedisk2}, as follows.
\par
Let $z_c$ be a special point of $\Sigma(\alpha)$.
Then as in Subsection \ref{subsec:onedisk2} 
we can find either a special point $z'_c$ on 
 $\Sigma'(\widehat\alpha)$ or a 
loop $\gamma_c : S^1 \to \Sigma'(\widehat\alpha)$.
We then put
\begin{equation}
f_{c,\alpha}(h_{\alpha})
=
\begin{cases}
{\rm dist}(h_{\alpha}u(z_c),u'(z'_c))^2
&\text{if $z'_c$ exists,}\\
\displaystyle
\int_{S^1} {\rm dist}(h_{\alpha}u(z_c),u'(\gamma_c(t)))^2 dt
&\text{if $\gamma_c$ exists.}
\end{cases}
\end{equation}
\begin{rem}
We note that we can take $\gamma_c$ in a way 
depending on the coordinates around the interior singular points 
on the component $\Sigma_{\alpha}$ but is 
{\it independent} of the coordinates around the 
boundary singular points of $\Sigma$ on $\Sigma_{\alpha}$,
as follows.
We consider the map
\begin{equation}\label{Psiclosedmap2}
\Psi : D^2(\epsilon)^{\mathcal S_{\alpha}} \times \mathcal V(\alpha,{0}) \times \prod_{a\in A_{\alpha}} \mathcal V(\alpha,a) \to \mathcal M_{0,\ell_{\alpha}+m_{\alpha}}
\end{equation}
which is defined in a similar way as (\ref{Psiclosedmap}).
Here $\mathcal S_{\alpha}$ is the set of interior singular points 
of $\Sigma(\alpha)$ and other notations are the same as
(\ref{Psiboreredalpha}).
\par
By Lemma \ref{ifremvethenclose}, 
${\bf v}'_{\widehat \alpha, \alpha}(\vec w^+)$ is in the image of this 
map.
We use this fact to define $\gamma_c$ in the same way as in 
Subsection \ref{subsec:onedisk2}.
Note that (\ref{Psiclosedmap2}) uses the coordinate around the 
interior singular points of $\Sigma_{\alpha}$ but does 
not use the coordinates around the 
boundary singular points of $\Sigma$ on $\Sigma_{\alpha}$.
\par
On the other hand, the map (\ref{Psiboreredalpha}) 
uses the coordinates around the 
boundary singular points of $\Sigma$ on $\Sigma_{\alpha}$.
We however use (\ref{Psiboreredalpha})  only to {\it prove}
Lemma \ref{ifremvethenclose}.
\end{rem}
We now put
$$
f_{+,\alpha} = f_{\alpha} + \sum_c f_{c,\alpha},
$$
where $f_{\alpha}$ is given by \eqref{fincase2components}. 
The rest of the argument is the same as Case A.
\par\smallskip
\noindent{\bf Case C: $\alpha$ is of Type C.}
\par
Let $h_{\alpha}$ be close to $h$.
By Lemma \ref{idonotchangemuch} we may choose 
$\epsilon_1,\epsilon_2$ small so that
$U_{\alpha} \cap {\bf S}_{\alpha} = \emptyset$.
We use it to define a map
\begin{equation}
\frak I^{h_{\alpha}}_{\widehat\alpha,\alpha} :
{\bf S} \to \Sigma'(\widehat\alpha)
\end{equation}
so that the following conditions are satisfied.
\begin{conds}\label{Smapcondsmulticomp}
Let $x \in {\bf S}_{\alpha}$.
\begin{enumerate}
\item
$
{\rm dist}(i_{{\bf v}'_{\widehat \alpha, \alpha}(\vec w)}(x),
\frak I^{h_{\alpha}}_{\widehat\alpha,\alpha}(x))
< \epsilon_3.
$
\item
The shortest geodesic joining $h_{\alpha}u(x)$ to 
$u'(\frak I^{h_{\alpha}}_{\widehat\alpha,\alpha}(x))$ is perpendicular to 
$h_{\alpha}u(\Sigma(\alpha))$ at $h_{\alpha}u(x)$.
\end{enumerate}
\end{conds}
Using Lemma \ref{idonotchangemuch}, we can prove the unique existence 
of such $\frak I^{h_{\alpha}}_{\widehat\alpha,\alpha}$.
We then define
\begin{equation}
f_{\alpha}(h_{\alpha})
= 
\int_{{\bf S}} {\rm dist}(h_{\alpha}u(x),
u'(\frak I^{h_{\alpha}}_{\widehat\alpha,\alpha}(x)))^2
dx.
\end{equation}
In the same way as Lemma \ref{S1equivalenceSS}, we can prove
$
f_{\alpha}(sh_{\alpha}) = f_{\alpha}(h_{\alpha}) 
$
for $s \in S^1_{\alpha}$ when both sides are defined.
We thus obtain a function $\overline f_{\alpha}$ 
on an open set of $T^n/S^1_{\alpha}$.
We can prove that it is strictly convex and attains its minimum 
at a unique point $[h^0_{\alpha}]$, 
in the same way as Lemma \ref{uniqueminimum}.
\par
We take its lift $h^0_{\alpha}$ and use it
to define 
\begin{equation}\label{PalphamanycomptypeIII}
\mathcal P_{\alpha} : E({\bf x}_{\alpha}) \to C^{\infty}(\Sigma'(\widehat\alpha),(u')^*TX \otimes \Lambda^{0,1}).
\end{equation}
We can prove that it is independent of the choice of the representative 
$h^0_{\alpha}$ but depends only on $[h^0_{\alpha}] 
\in T^n/S^1_{\alpha}$.
Since we have chosen $[h^0_{\alpha}]$ using $\overline f_{\alpha}$ 
already, the map (\ref{PalphamanycomptypeIII}) is defined. 
Moreover, in the same way as Lemma \ref{PisindepedentW}, 
we can prove that the image of $\mathcal P_{\alpha}$ 
is independent of the choice of $\vec w^+_{\alpha}$.
\begin{rem}\label{rem4571}
In case of Type A or Type B components, 
we have chosen $\vec w^+_{\alpha}$
inductively, at each stage of the construction. 
However we do not fix a particular choice of 
$\vec w^+_{\alpha}$ for Type C component.
Actually in case of Type A or Type B, our marked points $\vec w^+_{\alpha}$
play two different roles simultaneously. In case of Type C, 
one of the roles of $\vec w^+_{\alpha}$
in Type A or B cases is played by the circles ${\bf S}_{\alpha}$.
So the function $f_{\alpha}$ depends on the choice of 
$\vec w^+_{\alpha}$ in either Type A or B case but is independent of it 
in Type C case. 
($f_{\alpha}$ depends on ${\bf S}_{\alpha}$ in Type C case.)
\end{rem}
\medskip
We have thus defined $\mathcal P_{\alpha}$ for each $\alpha$.
We put
\begin{equation}
E({\bf x}') = \bigoplus_{\alpha\in \frak A}
{\rm Im} \mathcal P_{\alpha}.
\end{equation}
Note that so far we fixed $\vec w^{\prime +}$ and $h$ as in 
Definition \ref{primeisclosemodmod}. But we can show:
\begin{lem}
$E({\bf x}')$ is independent of the choice of $\vec w^{\prime +}$ and $h$ as in 
Definition \ref{primeisclosemodmod}.
Moreover it is independent of $\vec w^+_{\alpha}$ 
if $\alpha$ is of Type C.
\end{lem}
\begin{proof}
By construction, $E({\bf x}')$ does not change if we move $h$ 
to $h'$ that is in a sufficiently small neighborhood of $h$.
By Lemma \ref{w+primeunique33},
if 
$(h,\vec w^{\prime +}_{(1)})$,
$(h',\vec w^{\prime +}_{(2)})$ are two choices, 
then $h$ is transformed to a neighborhood of $h'$ by the $T^n$ component of an element of 
$G({\bf x})$. 
The first half of the lemma then follows from the construction.
The second half can be proved in the same way as 
Lemma \ref{PisindepedentW} (whose proof is similar 
to the proof of Lemma \ref{4546independentw3}.)
\end{proof}
We can prove that the assignment
${\bf x}' \mapsto E({\bf x}')$ is 
$T^n \times G({\bf x})$-equivariant in the same way as in Subsections
\ref{subsec:onedisk1} - \ref{subsec:onedisk3}.
We can then extend 
${\bf x}' \mapsto E({\bf x}')$ to the case when ${\bf x}'$
is close to a $T^n \times \frak S_{\ell}$ orbit of ${\bf x}'$,
in the same way as in the last part of Subsection \ref{subsec:onedisk1}.
\begin{lem}\label{dcwieness}
The Kuranishi neighborhood we defined above is disk component-wise.
\end{lem}
\begin{proof}
The choice of our map $\mathcal P_{\alpha}$ is made so that it depends only on ${\bf x}_{\alpha}$
and ${\bf x}'$. Especially it is independent of ${\bf x}_{\alpha'}$ for $\alpha' \ne \alpha$.
The lemma is immediate from this fact.
\end{proof}
\par
\subsection{Taking the sum of obstruction spaces so that 
coordinate change exists}
\label{subsec:coordinatechange}

We have thus constructed a Kuranishi neighborhood of each point 
of $\mathcal M_{0;\ell}(\beta)$.
To obtain a Kuranishi {\it structure} of $\mathcal M_{0;\ell}(\beta)$,
we need to define Kuranishi neighborhoods such that 
coordinate changes among them are defined.
For this purpose, we take a finite subset 
\index{$\frak P_{\ell}(\beta)$} $\frak P_{\ell}(\beta) \subset \mathcal M_{0;\ell}(\beta)/(T^n\times \frak S_{\ell})$ and 
associate data to determine obstruction bundle to each point of $\frak P_{\ell}(\beta)$, 
and then for each ${\bf x} \in \mathcal M_{0;\ell}(\beta)$ we take sum of the 
obstruction bundles obtained from all the nearby elements of $\frak P_{\ell}(\beta)$.
(See \cite[p.1003]{FO}, \cite[pp.423-424]{fooobook2} and \cite[Section 18]{foootech} for detail.)
In this subsection, we explain the way how we perform this process in a way 
compatible with $T^n$ action and disk-component-wise-ness.
\par
The construction is by induction on $\beta\cap \omega$.
In the case $\beta\cap \omega = 0$, the obstruction bundle is trivial and 
there is nothing to show.
Suppose $\beta \cap \omega$ is minimal among the nonzero $\beta$'s for which
$\mathcal M_{0;\ell}(\beta)$ is nonempty.
We take a sufficiently dense subset $\frak P_{\ell}(\beta) \subset \mathcal M_{0;\ell}(\beta)/(T^n\times \frak S_{\ell})$.
(More precisely, we require that (\ref{Pisdense}) is satisfied.)
\par
Let $\frak c \in \frak P_{\ell}(\beta)$. We take its representative 
$(\Sigma_{\frak c},\vec z^+_{\frak c},u_{\frak c}) \in \mathcal M_{0;\ell}(\beta)$.
By the minimality of $\beta$ we may assume that $\Sigma_{\frak c}$ contains 
only one disk component.

We take a $T^n \times \frak S_{\ell}$ invariant closed neighborhood $\frak U(\frak c)$ of $\frak c$ in $\mathcal M_{0;\ell}(\beta)$.
We fix $\vec w^+_{\frak c}$ if $\frak c$ is of type A or B as in 
Subsections \ref{subsec:onedisk1} - \ref{subsec:onedisk2}, 
and ${\bf S}_{\frak c}$ if $\frak c$ is of Type C as in Subsection \ref{subsec:onedisk3}.
We also fix $N_{w^+_{\frak c,b}}$ as in Condition \ref{1.8} and a finite dimensional 
subspace $E_{\frak c} \subset C^{\infty}(\Sigma_{\frak c};u_{\frak c}^*TX\otimes \Lambda^{0,1})$.
We assume that $\frak U(\frak c)$ is in the Kuranishi neighborhood constructed in 
Subsections \ref{subsec:onedisk1} - \ref{subsec:onedisk3}.
We take $\frak P_{\ell}(\beta)$ so that
\begin{equation}\label{Pisdense}
\bigcup_{\frak c \in \frak P_{\ell}(\beta)}\frak U(\frak c)
=\mathcal M_{0;\ell}(\beta).
\end{equation}
For ${\bf x} = (\Sigma_{\bf x},\vec z^+_{\bf x},u_{\bf x}) \in \mathcal M_{0;\ell}(\beta)$, 
we put
\begin{equation}\label{defofPd}
\frak C({\bf x}) = \{\frak c \in \frak P_{\ell}(\beta)
\mid {\bf x} \in \frak U(\frak c)\}.
\end{equation}
We choose $(\epsilon'_1,\epsilon'_2)$, 
such that if ${\bf y} = (\Sigma_{\bf y},\vec z^+_{\bf y},u_{\bf y})$ 
is $(\epsilon'_1,\epsilon'_2)$-close to $h {\bf x}$ for some $h \in T^n$ in the sense of 
Definition \ref{primeisclose} and $\frak c \in \frak C({\bf x})$, then 
${\bf y}$ is $(\epsilon_1,\epsilon_2)$-close to $h' {\frak c}$ for some $h' \in T^n$ 
in the sense of Definition \ref{primeisclose}.
\begin{rem}
More precisely, we specify $\vec w_{{\bf x},{\frak c}}$ as follows.
We assume the condition that ${\bf x}$ is $(\epsilon_1/2,\epsilon_2/2)$-close to 
an $T^n$ orbit of ${\frak c}$ means
that there exist $\vec w_{{\bf x},{\frak c}}$ 
(additional marked points on $\Sigma_{{\bf x}}$ and $h
\in T^n,
\mu_1 \in \frak S_{\ell}$ such that
$\mu_1{\bf x}$ is $(\epsilon_1,\epsilon_2)$-close 
to $(h{\bf c},\vec w_{\frak c})$ with respect to 
$\vec w_{{\bf x},{\frak c}}$.
\par
Then we take $(\epsilon'_1,\epsilon'_2)$ such that if
${\bf y}$ is $(\epsilon'_1,\epsilon'_2)$-close to $({\bf x},\vec w_{{\bf x},{\frak c}})$
then
${\bf y}$ is $(\epsilon_1,\epsilon_2)$-close to $(\frak c,\vec w_{{\frak c}})$.
\end{rem}
Then we can send $E_{\frak c}$ to 
$$
E_{\frak c}({\bf y})
=
C^{\infty}(\Sigma_{\bf y},u^*_{\bf y}TX\otimes \Lambda^{0,1})
$$
in the way described in Subsections \ref{subsec:onedisk1} - \ref{subsec:onedisk3}.
We now define
\begin{equation}\label{sumupdisk1}
E_{\bf x}({\bf y})
=
\bigoplus_{\frak c \in {\frak C}({\bf x})} E_{\frak c}({\bf y}).
\end{equation}
We can choose $E_{\frak c}$ so that 
the right hand side of (\ref{sumupdisk1}) is a direct sum.
(See \cite[Section 27]{foootech}.)
Now we put
$$
V_{\bf x} 
= \{{\bf y} \mid \overline\partial u_{\bf y} \in E_{\bf x}(\bf y)\}.
$$
Then ${\bf y} \mapsto E_{\bf x}(\bf y)$ defines a smooth vector bundle on it 
and $s_{\bf x}({\bf x}) = \overline\partial u_{\bf y}$ gives its smooth 
section.
We thus obtain a Kuranishi neighborhood of ${\bf x}$ in 
$\mathcal M_{0;\ell}(\beta)$.
We can define a coordinate change among this 
Kuranishi neighborhoods.
(See \cite[Section 22]{foootech} for detail.)
This completes the case when $\beta \cap \omega$ is smallest positive.
\begin{rem}
An important point of the construction above is that 
the space $E_{\frak c}({\bf y})$ depends only on ${\bf y}$ and 
${\frak c}$ and is independent of $\bf x$ as long as 
$\frak c \in \frak C({\bf x})$.
This is important for the coordinate change to be well-defined.
\end{rem}
\par
Now we assume as the induction hypothesis that 
we have chosen $\frak P_{\ell}(\beta')$
and $\frak U(\frak c)$ for $\frak c \in \frak P_{\ell}(\beta')$ when 
$\beta' \cap \omega < \beta \cap\omega$.
Moreover, we assume that we defined 
a Kuranishi structure on $\mathcal M_{0;\ell}(\beta')$ 
for $\beta' \cap \omega < \beta \cap\omega$ using the data on
$\frak P_{\ell}(\beta')$ etc..
We will construct $\frak P_{\ell}(\beta)$ and a Kuranishi structure on 
$\mathcal M_{0;\ell}(\beta)$.
\par
For $\frak d=2,3,\dots$ we denote by $S_{\frak d}\mathcal M_{0;\ell}(\beta)$ 
the subset of $\mathcal M_{0;\ell}(\beta)$ whose element 
has $\frak d$ or more disk components.
It is a closed subset.
We will construct a Kuranishi structure of a neighborhood of 
$S_{\frak d}\mathcal M_{0;\ell}(\beta)$ by a downward induction on $\frak d$.
\par
We consider $\frak d_0$ that is the maximum among those with 
$S_{\frak d_0}\mathcal M_{0;\ell}(\beta) \ne \emptyset$.
Let ${\bf x} = (\Sigma_{\bf x},\vec z^+_{\bf x},u_{\bf x}) \in S_{\frak d_0}\mathcal M_{0;\ell}(\beta)$.
We will construct a Kuranishi neighborhood of 
a $T^n\times \frak S_{\ell}$ neighborhood of ${\bf x}$.
We decompose $\Sigma_{\bf x}$ into disk components 
and obtain 
${\bf x}_{\alpha} = (\Sigma_{\bf x_{\alpha}},\vec z^+_{\bf x_{\alpha}},u_{\bf x_{\alpha}})
\in \mathcal M_{0;\ell_{\alpha}}(\beta_{\alpha})$
for $\alpha \in \frak A$.
($\#\frak A = \frak d_0$.)
Note $\Sigma_{\bf x_{\alpha}}$ has only one disk component.
We may assume $\beta_{\alpha}\omega < \beta\cap \omega$.
Therefore we have $\frak P_{\ell_{\alpha}}(\beta_{\alpha})$ and 
$\frak U(\frak c_{\alpha})$ for $\frak c_{\alpha} \in \frak P_{\ell_{\alpha}}(\beta_{\alpha})$
by induction hypothesis.
Note
$$
\frak C({\bf x}_{\alpha}) = \{\frak c_{\alpha} \mid {\bf x}_{\alpha} \in \frak U(\frak c_{\alpha})\}.
$$
We fix one $\alpha \in \frak A$ and 
$\frak c_{\alpha} \in \frak C({\bf x}_{\alpha})$.
Then there exist $h^0$ and $\mu_1$ such that
$\mu_1{\bf x}_{\alpha}$ is $(\epsilon_1/2,\epsilon_2/2)$ close 
to $h^0\frak c_{\alpha}$.
For simplicity of notation we assume $\mu_1 = 1$, $h^0 = 1$.
\par
Now we glue $\frak c_{\alpha}$ with ${\bf x}_{\alpha'}$ ($\alpha' \ne \alpha$)
to obtain ${\bf x}(\alpha)$ that is close to ${\bf x}$.  
Denote by ${\bf v}(\frak c_{\alpha})$, (resp. 
${\bf v}_{\alpha'}$) the domains of stable maps 
$\frak c_{\alpha}$, (resp. ${\bf x}_{\alpha'}$). 
Note that $\frak c_{\alpha}$ and ${\bf x}_{\alpha'}$ may not satisfy the incidence condition.  
We use cut-off functions for gluing them.  Since we only use the maps away from the gluing 
part, i.e., the supports of elements in $E(\frak c_{\alpha}), E({\bf x}_{\alpha'})$, additional interior 
marked points added on $\frak c_{\alpha}$ and ${\bf x}_{\alpha'}$, the gluing using cut-off function is 
enough for our purpose.  

Suppose ${\bf y}$ is $(\epsilon'_1,\epsilon'_2)$ close to 
${\bf x}$. Then it is $(\epsilon_1,\epsilon_2)$ close to
${\bf x}(\alpha)$.
Thus we can repeat the argument of Subsection \ref{subsec:1+} 
to obtain $h_{\alpha}$ and 
$$
\mathcal P_{c_{\alpha}} :
E_{\frak c_{\alpha}} \to h_{\alpha *}E_{\frak c_{\alpha}}
\to C^{\infty}(\Sigma_{\bf y},u_{\bf y}^*TX\otimes \Lambda^{0,1}).
$$
Denote by $E_{\frak c_{\alpha}}({\bf y})$ the image of 
$\mathcal P_{\frak c_{\alpha}}$.
We then define
\begin{equation}
E({\bf y})
=
\bigoplus_{\alpha \in \frak A} \bigoplus_{\frak c_{\alpha} \in \frak C({\bf x}_{\alpha})}
E_{\frak c_{\alpha}}({\bf y}).
\end{equation}
\begin{rem}
We note that the point $\in \partial \Sigma_{\frak c_{\alpha}}$ where we glue 
$\frak c_{\alpha_0}$ with ${\bf x}_{\alpha'}$ is not well-defined.
In fact, it is well-defined only up to small perturbations.
However, the construction of Subsection \ref{subsec:1+} 
is carefully designed so that this small ambiguity does 
not affect the element $h_{\alpha}$ and the map $\mathcal P_{\alpha}$ at all.
\end{rem}
Using $E({\bf y})$, we can define a Kuranishi neighborhood of ${\bf x}$ 
in the same way as the case of $\beta$ with the smallest $\beta \cap \omega$.
\par
We note that the space $E_{\frak c_{\alpha}}({\bf y})$ is independent of $\bf x$ as far as 
$\frak c_{\alpha} \in \frak C({\bf x}_{\alpha})$,
by construction. Therefore we can define coordinate change
among those Kuranishi neighborhoods to obtain a 
Kuranishi structure on a neighborhood of $S_{\frak d_0}\mathcal M_{0;\ell}(\beta)$.
\par
We next consider 
$S_{\frak d_0-1}\mathcal M_{0;\ell}(\beta)$.
Let
$\frak V(S_{\frak d_0}\mathcal M_{0;\ell}(\beta))$ be the 
intersection of $\mathcal M_{0;\ell}(\beta)$ and  the Kuranishi neighborhood of $S_{\frak d_0}\mathcal M_{0;\ell}(\beta)$
we constructed above.
We consider
$$
\mathcal K_{\frak d_0-1}\mathcal M_{0;\ell}(\beta) 
=
S_{\frak d_0-1}\mathcal M_{0;\ell}(\beta) 
\setminus  \frak V(S_{\frak d_0}\mathcal M_{0;\ell}(\beta)).
$$
We can repeat the same construction as 
the case of $S_{\frak d_0}\mathcal M_{0;\ell}(\beta)$
to obtain a Kuranishi structure on a neighborhood of 
$\mathcal K_{\frak d_0-1}\mathcal M_{0;\ell}(\beta)$.
Let
$\mathcal V(\mathcal K_{\frak d_0-1}\mathcal M_{0;\ell}(\beta))$ be 
the intersection of this neighborhood and 
$\mathcal M_{0;\ell}(\beta)$.
We glue these two Kuranishi structures as follows.
\begin{enumerate}
\item
If ${\bf x} \in \mathcal K_{\frak d_0-1}\mathcal M_{0;\ell}(\beta)$, 
we use the Kuranishi neighborhood of $\mathcal K_{\frak d_0-1}\mathcal M_{0;\ell}(\beta)$
constructed above to be the Kuranishi neighborhood of ${\bf x}$.
\item For ${\bf x}\in \frak V(S_{\frak d_0}\mathcal M_{0;\ell}(\beta))
\setminus \mathcal V(\mathcal K_{\frak d_0-1}\mathcal M_{0;\ell}(\beta))$, 
we use the Kuranishi neighborhood of ${\bf x}$ 
that is induced by the Kuranishi structure of $S_{\frak d_0}\mathcal M_{0;\ell}(\beta)$.
\item Let 
${\bf x} \in \frak V(S_{\frak d_0}\mathcal M_{0;\ell}(\beta)) \cap 
\mathcal V(\mathcal K_{\frak d_0-1}\mathcal M_{0;\ell}(\beta))$. 
For any ${\bf y}$ which is sufficiently close to ${\bf x}$,
we take sum of the two obstruction spaces to obtain $E({\bf y})$. Here one of the 
obstruction spaces is defined from 
$\mathcal K_{\frak d_0-1}\mathcal M_{0;\ell}(\beta)$
and the other is defined from $S_{\frak d_0}\mathcal M_{0;\ell}(\beta)$.
If necessary, we can perturb the obstruction spaces 
to define our Kuranishi structure on 
$\mathcal K_{\frak d_0-1}\mathcal M_{0;\ell}(\beta)$
and can assume that this sum is a direct sum.
We use this $E({\bf y})$ (that is the sum) 
to define a Kuranishi neighborhood of $E({\bf x})$.
\end{enumerate}
By this choice of Kuranishi neighborhoods, we have  the coordinate changes. Thus we obtain 
a Kuranishi structure on a neighborhood of 
$S_{\frak d_0-1}\mathcal M_{0;\ell}(\beta)$.
\par
Now we repeat the same construction inductively by the 
downward induction on $\frak d$.
We then obtain a Kuranishi structure on a 
neighborhood of 
$
\partial \mathcal M_{0;\ell}(\beta) 
= S_{2}\mathcal M_{0;\ell}(\beta)
$.
Let  $\frak V(S_{2}\mathcal M_{0;\ell}(\beta))$
be the intersection of this neighborhood and 
$\mathcal M_{0;\ell}(\beta)$.
We then take a sufficiently dense subset 
$\frak P_{\ell}(\beta)$
of 
$\mathcal K_{1}\mathcal M_{0;\ell}(\beta)
= \mathcal M_{0;\ell}(\beta) \setminus \frak V(S_{2}\mathcal M_{0;\ell}(\beta))$
to obtain a Kuranishi structure of an neighborhood of 
$\mathcal K_{1}\mathcal M_{0;\ell}(\beta)$.
We then glue two Kuranishi structures in the same way 
to obtain the required Kuranishi structure on $\mathcal M_{0;\ell}(\beta)$. 
The proof of Proposition \ref{eqcyckuramain} is now complete.
\end{proof}
\begin{rem}\label{colermethod}
So far, in this section, we have described the construction of a system of Kuranishi structures 
on $\mathcal M_{0;\ell}(\beta)$ which is $T^n$ equivariant and is disk component-wise.
\par
In this remark we sketch another method, in which we first construct a  
$T^n$ equivariant Kuranishi structure for each $\mathcal M_{0;\ell}(\beta)$ 
and modify it later near the boundary and corner so that 
compatibility with fiber product (something similar to disk component-wise-ness)
will be achieved. 
Note as we mentioned right after Lemma \ref{w+primeunique33}, 
the construction of $T^n$ equivariant Kuranishi neighborhood centered at 
${\bf x}$ with 2 or more disk components is easier than the argument 
of the second half of this subsection when we do not require disk-component-wise-ness.
\par
Suppose that we have  $T^n$ equivariant (but not necessary disk-component-wise)
Kuranishi structures on $\mathcal M_{0;\ell}(\beta)$.
We first describe the way how to deform 
them in a neighborhood of 
$$
\mathcal M_{1;\ell_1}(\beta_1) \times_L \mathcal M_{1;\ell_2}(\beta_2)
\subset \partial \mathcal M_{0;\ell}(\beta),
$$
where $\beta_1 + \beta_2 = \beta$, $\ell_1 + \ell_2 = \ell$.
\par
We pull back Kuranishi structures of $\mathcal M_{0;\ell}(\beta)$ to 
$\mathcal M_{k;\ell}(\beta)$ via the forgetting map of 
the boundary marked points.
\par
We have two Kuranishi structures on 
$\mathcal M_{1;\ell_1}(\beta_1) \times_L \mathcal M_{1;\ell_2}(\beta_2)$,
one is the fiber product Kuranishi structure and
the other is the restriction of the Kuranishi structure on 
$\mathcal M_{0;\ell}(\beta)$.
It is easy to construct a cobordism between them.
The cobordism is a ($T^n$-equivariant) Kuranishi structure on 
$(\mathcal M_{1;\ell_1}(\beta_1) \times_L \mathcal M_{1;\ell_2}(\beta_2)) \times [0,1]$.
Let $t$ be the coordinate of $[0,1]$.
Then at $t=0$ it coincides the restriction and at $t=1$ it coincides with 
the fiber product.
Therefore we glue it with $\mathcal M_{0;\ell}(\beta)$ at $t=0$ to obtain a Kuranishi structure on 
$\mathcal M_{0;\ell}(\beta) \cup ((\mathcal M_{1;\ell_1}(\beta_1) \times_L \mathcal M_{1;\ell_2}(\beta_2)) \times [0,1])$
so that its boundary $\mathcal M_{1;\ell_1}(\beta_1) \times_L \mathcal M_{1;\ell_2}(\beta_2)$
has fiber product Kuranishi structure.
\par
We next consider the corner point of $\mathcal M_{0;\ell}(\beta)$.
Typically it is described as
\begin{equation}\label{3fiber}
(\mathcal M_{1;\ell_1}(\beta_1) \times \mathcal M_{1;\ell_2}(\beta_2)) \times_{L^2} \mathcal M_{2;\ell_3}(\beta_3),
\end{equation}
where $\beta_1 + \beta_2 + \beta_3 = \beta$, $\ell_1 + \ell_2 + \ell_3 = \ell$.
The fiber product Kuranishi structure does not coincide with the restriction 
of the Kuranishi structure on $\mathcal M_{0;\ell}(\beta)$.
We attach an appropriate collar to compensate this disagreement as follows  
(see Figure 4.3.2).
\par
We put $\beta' = \beta_1 + \beta_3$, $\beta'' = \beta_2 + \beta_3$,
$\ell' = \ell_1 + \ell_3$, $\ell'' = \ell_2 + \ell_3$. 
By applying the first part of the argument, we obtain Kuranishi structures on 
\begin{equation}\label{idomenofiber}
\aligned
&\mathcal M_{0;\ell'}(\beta') \cup (\partial\mathcal M_{0;\ell'}(\beta') \times [0,1]_t), \\
&\mathcal M_{0;\ell''}(\beta'') \cup (\partial\mathcal M_{0;\ell''}(\beta'') \times [0,1]_s),
\endaligned
\end{equation}
so that the boundaries are given by the fiber product Kuranishi structures
$$
\mathcal M_{1;\ell_1}(\beta_1) \times_L \mathcal M_{1;\ell_3}(\beta_3),
\quad
\mathcal M_{1;\ell_2}(\beta_2) \times_L \mathcal M_{1;\ell_3}(\beta_3).
$$
(Here $t$, $s$ in the indices show that the parameters we use for 
$[0,1]$.)
\par
We add one boundary marked point to (\ref{idomenofiber}),
pull back the Kuranishi structure and then 
take a fiber product with  $\mathcal M_{1;\ell_2}(\beta_2)$ and $\mathcal M_{1;\ell_1}(\beta_1)$ respectively.
We then obtain Kuranishi structures on  
\begin{equation}\label{idomenofiber2}
\aligned
&(\mathcal M_{1;\ell'}(\beta') \times_L \mathcal M_{1;\ell_2}(\beta_2)) \cup (\partial\mathcal M_{1;\ell'}(\beta') \times_L \mathcal M_{1;\ell_2}(\beta_2)\times [0,1]_t \times \{1\}_s), \\
&(\mathcal M_{1;\ell''}(\beta'') \times_L \mathcal M_{1;\ell_1}(\beta_1)) \cup (\partial\mathcal M_{1;\ell''}(\beta'') \times_L \mathcal M_{1;\ell_1}(\beta_1)\times \{1\}_t \times [0,1]_s).
\endaligned
\end{equation}
The boundary of them (that is the part where $t=1$ or $s=1$ respectively) 
is the fiber product Kuranishi structure (\ref{3fiber}).
\par
We next take the cobordism between 
$\mathcal M_{1;\ell'}(\beta') \times_L \mathcal M_{1;\ell_2}(\beta_2)$ 
(resp. 
$\mathcal M_{1;\ell''}(\beta'') \times_L \mathcal M_{1;\ell_1}(\beta_1)$)
with fiber product Kuranishi structures 
and 
their Kuranishi structures which are restrictions of the 
Kuranishi structure of  $\mathcal M_{0;\ell}(\beta)$.
We note that  
$$
\partial\mathcal M_{0;\ell}(\beta)
\supset (\mathcal M_{1;\ell'}(\beta') \times_L \mathcal M_{1;\ell_2}(\beta_2))
\cup 
(\mathcal M_{1;\ell''}(\beta'') \times_L \mathcal M_{1;\ell_1}(\beta_1)).
$$
These cobordisms are Kuranishi structures on 
\begin{equation}\label{idomenofiber3}
\aligned
&(\mathcal M_{1;\ell'}(\beta') \times_L \mathcal M_{1;\ell_2}(\beta_2)) \times \{0\}_t \times [0,1]_s, \\
&(\mathcal M_{1;\ell''}(\beta'') \times_L \mathcal M_{1;\ell_1}(\beta_1))  \times [0,1]_t \times \{0\}_s,
\endaligned
\end{equation}
respectively.
We can glue Kuranishi structures of (\ref{idomenofiber2}), (\ref{idomenofiber3}) and of 
$\mathcal M_{0;\ell}(\beta)$ at $s=0$ or $t=0$.
We finally extend the Kuranishi structures we gave on
$$
(\mathcal M_{1;\ell_1}(\beta_1) \times \mathcal M_{1;\ell_2}(\beta_2)) \times_{L^2} \mathcal M_{2;\ell_3}(\beta_3)
\times \partial ([0,1]_t \times [0,1]_s)
$$
(that are restriction of those on (\ref{idomenofiber2}), (\ref{idomenofiber3})) 
to a Kuranishi structure on 
$$
(\mathcal M_{1;\ell_1}(\beta_1) \times \mathcal M_{1;\ell_2}(\beta_2)) \times_{L^2} \mathcal M_{2;\ell_3}(\beta_3)
\times ([0,1]_t \times [0,1]_s).
$$
They altogether give the desired Kuranishi structure on $\mathcal M_{0;\ell}(\beta)$ plus collar.
(See Figure 4.3.2.)
\par
We can continue in the same way in the case of higher codimensional corners 
to obtain a desired system of Kuranishi structures.
\par 
We remark that this method can be used in the other situations such as those appearing in 
Section \ref{sec:compwiseplusTnequiv}.
\end{rem}
\par
\hskip0.3cm
\epsfbox{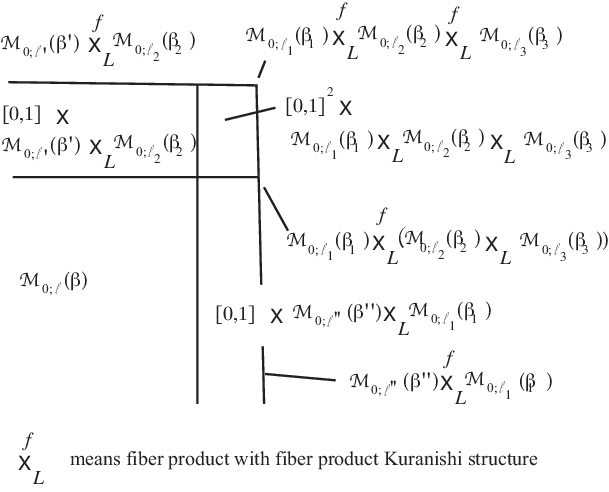}
\par
\centerline{\bf Figure 4.3.2}
\par
\par
\subsection{Taking fiber product with $T^n$ invariant cycles}
\label{subsec:fiberproduct}

We can use Proposition \ref{eqcyckuramain} to prove
Proposition \ref{transconcl}
in a similar way as in \cite[Corollary 3.1]{fooo091}.
We however need to modify the argument slightly by the following reason:
The evaluation map
$$
{\rm ev}^{\rm int} : \mathcal M_{k+1;\ell}^{\rm main}(\beta) \to X^{\ell}
$$
at the interior marked points is not weakly submersive with respect 
to the Kuranishi structure we constructed in Proposition \ref{eqcyckuramain}.
This is because we do not put any obstruction space on the disk or sphere component 
where the map is constant.
So the definition of the fiber product given in \cite[Section A1.2]{fooobook2}
does not apply here.
\par
One way to go around this trouble is to use de Rham representative 
of the cohomology group of $X$ rather than taking submanifolds 
(that are $T^n$ invariant submanifolds in the toric divisor).
\par
In this subsection we take a different way and use an alternative definition
of the fiber product to produce a required Kuranishi structure on 
$\mathcal M_{k+1;\ell}^{\rm main}(\beta;\bf p)$.
It then implies Proposition \ref{existcycpert} by using Lemma \ref{append3mainlemma} and 
the argument of  \cite[Corollary 3.1]{fooo091}.
\par
Let ${\bf p} = ({\bf p}(1),\dots,{\bf p}({\ell}))$.
We take a $T^n$ equivariant tubular neighborhood $U({\bf p}(i))$ of ${\bf p}(i)$.
We consider the normal bundle $N_{{\bf p}(i)}X$ and its $\epsilon$ ball bundle 
$N^{\epsilon}_{{\bf p}(i)}X$. We have a $T^n$ equivariant diffeomorphism 
\begin{equation}\label{tublariso}
I_{{\bf p}(i)} : N^{\epsilon}_{{\bf p}(i)}X \cong U({\bf p}(i)).
\end{equation}
Let ${\bf x} \in \mathcal M_{k+1;\ell}^{\rm main}(\beta;\bf p)$.
Namely  ${\bf x} \in \mathcal M_{k+1;\ell}^{\rm main}(\beta)$ 
such that 
${\rm ev}^{\rm int} _i({\bf x}) \subset {\bf p}(i)$. 
(Since $ {\bf p}(i)$ are submanifolds of $X$, 
the fiber product  
$$
\mathcal M_{k+1;\ell}^{\rm main}(\beta;{\bf p})
=
\mathcal M_{k+1;\ell}^{\rm main}(\beta) \,\,{}_{{\rm ev}^{\rm int} }\times \prod_{i=1}^{\ell} {\bf p}(i)
$$ 
is a subspace of
$\mathcal M_{k+1;\ell}^{\rm main}(\beta)$.) 
Let 
$(V_{\bf x},E_{\bf x},\Gamma_{\bf x},s_{\bf x},\psi_{\bf x})$ be a 
Kuranishi neighborhood of $\bf x$ in 
$\mathcal M_{k+1;\ell}^{\rm main}(\beta)$.
We will define a 
Kuranishi neighborhood
$(V_{\bf x}^{\bf p},E_{\bf x}^{\bf p},\Gamma_{\bf x}^{\bf p},s_{\bf x}^{\bf p},\psi_{\bf x}^{\bf p})$
of $\bf x$ in $\mathcal M_{k+1;\ell}^{\rm main}(\beta;\bf p)$.
\par
We put
$$
V_{\bf x}^{\bf p} = \{{\bf y} \in V_{\bf x} \mid {\rm ev}^{\rm int}_i({\bf y}) \in U({\bf p}(i)), 
\text{ for $i=1,\dots,\ell$}\}.
$$
This is a $T^n$ and $\Gamma_{\bf x}$ invariant open subset of $V_{\bf x}$.
We put $\Gamma_{\bf x}^{\bf p} = \Gamma_{\bf x}$.
We next define
$$
E_{\bf x}^{\bf p} = 
E_{\bf x}\vert_{V_{\bf x}^{\bf p}} \oplus \bigoplus_{i=1}^{\ell} (\pi\circ I_{{\bf p}(i)}^{-1}\circ  {\rm ev}_i^{\rm int} )^*N_{{\bf p}(i)}X.
$$
Here we use $\pi : N^{\epsilon}_{{\bf p}(i)}X \to {\bf p}(i)$.
We next put
$$
s_{\bf x}^{\bf p}({\bf y}) =  s_{\bf x}({\bf y}) \oplus (I_{{\bf p}(i)}^{-1}\circ {\rm ev}_i^{\rm int} ({\bf y}))_{i=1,\dots,\ell}.
$$
Note $I_{{\bf p}(i)}^{-1}\circ {\rm ev}_i^{\rm int} ({\bf y})$ is an element of the (total space of) the normal bundle
$N_{{\bf p}(i)}X$. We can regard it as an element of the bundle 
$(\pi\circ I_{{\bf p}(i)}^{-1}\circ  {\rm ev}_i^{\rm int} )^*N_{{\bf p}(i)}X$ in an obvious way. 
It is easy to see from the definition that
$
s_{\bf x}^{\bf p}({\bf y}) = 0
$
holds if and only if 
$
s_{\bf x}({\bf y}) = 0
$
and 
${\rm ev}_i^{\rm int} ({\bf y}) \in {\bf p}(i)$ for $i=1,\dots,\ell$.
Therefore the restriction of $\psi_{\bf x}$ gives $\psi_{\bf x}^{\bf p}$.
\par
Using the coordinate change between Kuranishi neighborhoods $(V_{\bf x},E_{\bf x},\Gamma_{\bf x},s_{\bf x},\psi_{\bf x})$,  
we can easily construct the coordinate change between 
Kuranishi neighborhoods $(V_{\bf x}^{\bf p},E_{\bf x}^{\bf p},\Gamma_{\bf x}^{\bf p},s_{\bf x}^{\bf p},\psi_{\bf x}^{\bf p})$.
We have thus constructed the required Kuranishi structure on $\mathcal M_{k+1;\ell}^{\rm main}(\beta;\bf p)$.
\qed 

\section{Continuous family of $T^n$ equivariant and cyclic symmetric multisections}
\label{sec:equimulticot}

In this section we construct a continuous family of multisections with cyclic symmetry
for the Kuranishi structure produced in the last section.
Most of the construction is a straightforward combination of those
employed in \cite{fooo08, fooo09, fooo091}.
We only need to prove the following Lemma \ref{append3mainlemma}.
\par
We consider $(U,E,\Gamma,s)$ and $f$, $N$ with the following properties:
\begin{enumerate}
\item
$U$ is a smooth manifold with a $T^n$ action.
\item
$E$ is a $T^n$ equivariant vector bundle on $U$.
\item
$\Gamma$ is a finite group which has {\it right} action on $U$ and
$E$. Moreover $E$  is a $\Gamma$ equivariant vector bundle.
The action of $\Gamma$ on $U$ is effective.
\item The $\Gamma$ action and $T^n$ action commute and so defines a
$T^n \times \Gamma$ action on $E$.
\item
$N$ is a smooth manifold and
$f: U\to N$ is a smooth submersion.
\item $s$ is a $T^n \times \Gamma$ equivariant section of $E$.
\end{enumerate}
\begin{lem}\label{append3mainlemma}
There exist a manifold $W$ with a $T^n$ action,
and a multisection $\frak s$ of $E/\Gamma \times W \to U_0/\Gamma \times W$
that satisfy the following properties:
\begin{enumerate}
\item $U_0$ is a $T^n\times \Gamma$ invariant neighborhood of $s^{-1}(0)$.
\item $\frak s$ is smooth and is transversal to $0$.
\item
The restriction of $f$ to $\frak s^{-1}(0)$ is a submersion to $N$.
\item Each branch of $\frak s$ can be chosen to be as close to $s$ as we want in $C^0$
sense.
\item $\frak s$ is $T^n$ equivariant.
\end{enumerate}
\end{lem}
\begin{proof}
Let $p \in U$. We put
$I_p = \{ g \in T^n \mid gp=p\}$, $\Gamma_p = \{\gamma \in \Gamma
\mid  p = p\gamma\}$, $G_{p,+} = \{(g,\gamma)
\in T^n \times \Gamma \mid gp =  p\gamma\}$,
which contains $I_p \times \Gamma_p$.
\par
We define a {\it left} action of $\Gamma$ by $\gamma p = p {\gamma^{-1}}$.
Hereafter we only use this left action.
The group $G_{p,+} $ acts on $X$ and $E$ as a subgroup of $T^n \times \Gamma$.
By this action, $p$ is a fixed point of $G_{p,+}$.
So $G_{p,+} $ acts on the fiber $E_p$.
\par
We choose a submanifold $X_p$ containing $p$ which is
a local transversal to the $T^n$ orbit of $p$.
We may choose $X_p$ so that it is $G_{p,+} $ invariant.
We take a $T^n \times \Gamma$ invariant metric of $U$ and a
$T^n$ and $\Gamma$ invariant connection of $E$.
\par
We put $W_p = E_p$ and define a section $s^{(1)}_p$ of $E$ on $X_p \times W_p$ by
\begin{equation}\label{spfirst}
s^{(1)}_p(x,w) = s(x) + \Pi_p^x(w).
\end{equation}
This is an $I_p$-invariant section.  
Here $\Pi_p^x$ is the parallel transport along the minimal geodesic joining $p$ and $x$. The section
$s^{(1)}_p$ is transversal to $0$. 
\par
Let $I_{p,0}$ be the connected component of $I_p$.
(In the situation we apply Lemma \ref{append3mainlemma},
the group $I_{p,0}$ is $S^1$.) 
We take a closed subgroup $J_p \subset T^n$ such that 
$I_{p,0} \times J_p \cong T^n$. 
We note that $I_{p}$ acts on $W_p$ and $X_p$, and 
the section $s^{(1)}_p$ is $I_{p}$ equivariant.
We put $U_p = T^nX_p = J_pX_p$.
\par
For any $x$ there exists $g \in J_p$ unique up to $I_p \cap J_p$  such that 
$g^{-1}x \in X_p$. We extend $s^{(1)}_p$ to $U_p$ by
$$
s^{(2)}_p(gx,w) = g_* s^{(1)}_p(x,w)  
$$
for $x \in X_p$ and $g \in J_p$.
\par
We define $T^n$ action on $W_p$ by 
using $I_{p,0}$ on $W_p$ and the splitting 
$I_{p,0} \times J_p \cong T^n$,
where $J_p$ action on $W_p$ is trivial.
(We remark that for $g \in I_p \setminus I_{p,0}$ 
the $g$ action on $E_p = W_p$ may be different 
from the above $g \in T^n$ action on $W_p$.
But this does not matter.) 
Then the section $s^{(2)}_p$ on $U_p \times W_p$ is 
obviously $T^n$ equivariant.
\par
We next take a $W_p$ parametrized family of multisections 
on $U_p$ such that it is also $\Gamma_{p,+}$
invariant, where
$$
\Gamma_{p,+} = \{ \gamma \in \Gamma \mid \gamma U_p \cap U_p
\ne \emptyset \}.
$$
We may choose $X_p$ so small that $\gamma \in \Gamma_{p,+}$
holds only when $\gamma p \in T^np$.
Then we have 
\begin{equation}\label{Upisgammainvariant}
\gamma U_p = U_p
\end{equation}
for each $\gamma \in \Gamma_{p,+}$.
Let $\ell$ be the order of $\Gamma_{p,+}$.
We define $\ell$-multisection $s^{(3)}_p$ on 
$U_p \times W_p$ by
 by
\begin{equation}\label{sp3ban}
s_p^{(3)}(x,w) = \{ \gamma s_p^{(2)}(\gamma^{-1}x,w) \mid \gamma 
\in \Gamma_{p,+}\}.
\end{equation}
Since $\Gamma$ action commutes with $T^n$ action and 
$s^{(2)}_p$ is $T^n$ equivariant, 
it follows that $s_p^{(3)}$ is $T^n$ equivariant.
Moreover $s_p^{(3)}$  is $\Gamma_{p,+}$ invariant. 
\par
We now extend $s_p^{(3)}$ to $U_p^+ = \Gamma U_p$ as follows.
Let $(y,w) \in U_p^+ \times W_p$. We have
$y = \gamma x$ with $\gamma \in \Gamma$, $x \in U_p$.
We put
\begin{equation}\label{sp4ban}
s_p^{(4)}(y,w) = \gamma_* s_p^{(3)}(x,w).
\end{equation}
The multisection $s_p^{(4)}$ is well-defined since $s_p^{(3)}$ is $\Gamma_{p,+}$ invariant.
\par
It is easy to see that $s_p^{(4)}$ is $\Gamma \times T^n$-equivariant.
(\ref{spfirst}) implies that $s^{(1)}_p$ satisfies Lemma \ref{append3mainlemma}.2.
Therefore $s_p^{(4)}$ satisfies Lemma \ref{append3mainlemma}.2
also.
\par
The multisection $s_p^{(4)}$ also satisfies Lemma \ref{append3mainlemma}.3.
In fact, by (\ref{spfirst}) the map 
$(s^{(1)}_p)^{-1}(0) \cap (X_p \times W_p) \to X_p$ is a submersion.
Since $s^{(2)}_p$ is $T^n$ equivariant it implies that
$(s_p^{(2)})^{-1}(0)  \to U$ is a submersion. 
By construction each branch of 
$s_p^{(4)}$ has the same property.
Lemma \ref{append3mainlemma}.3 follows.
\par
Thus
$s_p^{(4)}$ gives a required $W_p$ parametrized family of multisections on $U_p^+$.
\par
Then by gluing $s_{p_i}^{(4)}$ for various $p_i$ by an appropriate partitions of unity, we obtain
the required family of multisections.
\end{proof}
Using Lemma \ref{append3mainlemma}, the rest of the construction is
the same as in \cite{fooo08}, \cite{fooo09} and \cite{fooo091} and so is omitted.
\par
\section{Construction of various other Kuranishi structures.}
\label{sec:compwiseplusTnequiv}

\par
\subsection{Construction of $\frak p$-Kuranishi structure}
\label{subsec:pkura}
\par
In this subsection we prove Lemma \ref{pmodulikura}.\index{Kuranishi structure!$\frak p$-Kuranishi structure}
The proof is mostly the same as the proof of 
Proposition  \ref{transconcl}, which is the construction of $\frak c$-Kuranishi structure.
The difference is Lemma \ref{pmodulikura}.5.
Namely we need to construct a Kuranishi structure 
so that the evaluation map
\begin{equation}
{\rm ev}^{\rm int}_{0} : \mathcal M^{\rm main}_{k;\ell+1}(\beta) \to X
\end{equation}
at the $0$-th interior marked point is weakly submersive.
On the other hand, the Kuranishi structure constructed in Lemma \ref{pmodulikura} may not 
be invariant under the permutaiton between $0$-th 
{\it interior} marked point 
and 1-st, \dots, $\ell$-th {\it interior} marked points.
\par
We first observe that 
when $\beta \ne 0$, we can choose our obstruction 
space $E$ as in Lemma \ref{app2lema} so that the map
$$
(D_u\overline\partial)^{-1}(E) \to T_{u(z^+_0)}X
$$
is a submersion. Here $u \in \mathcal M^{\rm main}_{0;\ell+1}(\beta)$ 
and $z^+_0$ is the $0$-th (interior) marked point and 
the map is evaluation at $z_0^+$.
\par
Then most part of the proof of Proposition  \ref{transconcl}
goes without change.
(Note we use the $\frak c$-Kuranishi structure,
which is one obtained in Proposition \ref{eqcyckuramain},  
on the extended disk component if  $z^+_0$ is not on it.)
\par
The main difference between the construction of $\frak p$-Kuranishi structure and of $\frak c$-Kuranishi 
structure, comes from the fact that  we need to study also the case 
when $\beta=\beta_0 = 0$ to construct $\frak p$-Kuranishi structure.
In fact, we need to put some obstruction space 
on the component where $u$ is constant.
This is because
the image of 
${\rm ev}^{\rm int}_{0} :  \mathcal M^{\rm main}_{k;\ell+1}(\beta_0) \to X$
is in $L$ so it is not a submersion if we do not put 
obstruction bundle on it.
We can certainly obtain an obstruction space on this component 
so that ${\rm ev}_{0}^{\rm int}$ becomes a submersion on the Kuranishi neighborhood.
The new point we need to discuss is the way how we send this obstruction 
space to $\Sigma'$ that is the source of the nearby object.
We will discuss this point below.
\par
Let
$$
{\bf x} = (\Sigma,z_0^+\cup \vec  z^+,u) 
\in \mathcal M^{\rm main}_{0;\ell+1}(\beta).
$$
We assume that there exists only one disk component in $\Sigma$.
Let $\Sigma_0 \subset \Sigma$ be the disk component.
Then $\Sigma_0$ and $u\vert_{\Sigma_0}$ together with special points on $\Sigma_0$
define an object which we denote ${\bf x}_0$.
We define $z_0^{++} \in \Sigma_0$ as follows.
If $z_0^{+} \in \Sigma_0$, then $z_0^{++} = z_0^+$.
If $z_0^{+} \notin \Sigma_0$, 
then there exists a (maximal) tree of sphere components containing 
$z_0^+$. The point $z_0^{++}$ is the root of this tree of 
sphere components.
We assume that $u\vert_{\Sigma_0}$ 
is a constant map to a point of $L$.
We take an obstruction space
$$
E({\bf x}) \subset C^{\infty}(\Sigma;u^*TX \otimes \Lambda^{0,1})
$$
so that the following holds.
$E({\bf x})$ is a direct sum of $E({\bf x}_0)$ and  $E({\bf x}_a)'s$.
(The latter's are supported on the sphere bubbles.)
We require that they satisfy the following Conditions \ref{submersivez0} and \ref{submersiveatsing}. 
\begin{conds}\label{submersivez0}
\begin{enumerate} 
\item
The evaluation map
$$
{\rm ev}^{\rm int}_0 : (D_u\overline\partial)^{-1}(E({\bf x}_0)) \to T_{u(z^{++}_0)}X
$$
at $z^{++}_0$ is a submersion. 
\item
We do not put obstruction spaces on the {\it sphere} bubbles where the map $u$ is constant.
\item
If $\Sigma_0$ contains more than one special points, then 
${\bf x}_0$ is stable. We choose $E({\bf x}_0)$ so that it is invariant under 
${\rm Aut}_+({\bf x}_{0})$ action.
\item 
If $\Sigma_0$ contains only one special point
(which is $z^{++}_0$) and $u$ is constant on $\Sigma_0$, then 
the group of automorphisms of ${\bf x}_0$ is $S^1$ and so is compact. We take 
$E({\bf x}_0)$ so that 
it is invariant under the $S^1$ action.
\item
$E({\bf x})$ is invariant under ${\rm Aut}_+({\bf x})$ action.
\end{enumerate}
\end{conds}
To find a Kuranishi structure that satisfies the 
required property in the case when $z^+_0$ is on the sphere bubble 
also, we need some additional properties of the obstruction bundle on the sphere bubble.
Let $\Sigma_a$ be any sphere component on which $u$ is nonconstant.
Let $\vec z_a^{+,{\rm sing}}$ be the set of interior singular points 
contained in $\Sigma_a$.
Let 
$$
E(\Sigma_a)\subset C^{\infty}(\Sigma_a;u^*TX \otimes \Lambda^{0,1}).
$$
We assume:
\begin{conds}\label{submersiveatsing}  
If $u$ is nonconstant on $\Sigma_a$, 
the evaluation map
$$
{\rm ev}_{\rm sing} : 
(D_u\overline\partial)^{-1}(E(\Sigma_a))
\to \prod_{z_{a,i}^+ \in \vec z_a^{+,{\rm sing}}}
T_{u(z_{a,i}^+)}X
$$
at the points of $\vec z_a^{+,{\rm sing}}$ is surjective.
\end{conds}
We can show that ${\rm ev}_0^{\text{\rm int}}:{\mathcal M}_{k;\ell +1}(\beta) \to X$ is weakly subversive if
Conditions \ref{submersivez0} and \ref{submersiveatsing} are satisfied by the induction over the number of
irreducible components and the energy of $\beta$.
We consider the weak submersivity at
$
{\bf x} = (\Sigma,z_0^+\cup {\vec z}^{\,+},u)
\in \mathcal M^{\rm main}_{0;\ell+1}(\beta).
$
\par
We argue by the induction over the number of irreducible components of the domain of $u$.
If $z_0^{++} = z_0^+$
then it is a consequence of Condition \ref{submersivez0} (1).
If $u$ is nonconstant on the sphere component
which contains $z_0^+$, then it is a consequence of
Condition \ref{submersiveatsing}.
\par
Otherwise we take minimal union of sphere components
$\bigcup_{a \in A}S^2_a$ which is
connected and contains $z_0^{++}$ and $z_0^+$.
Suppose that the claim holds for $u$ with $N$ irreducible components
and consider the case that the domain of $u$ consists of $N+1$ irreducible components.
Let $a(z_0^+) \in A $ be the component such that $z_0^+ \in S^2_{a(z_0^+)}$.
\par
Suppose that there are no sphere components other than those in $A$.
Then we remove $S^2_{a(z_0^+)}$ from the domain of the stable map $u$ and
leave the unique nodal point of $S^2_{a(z^+_0)}$ as a marked point marked point
denoted by $z^{+++}_0$ of the other component to obtain a stable map $u'$.
Since the number of irreducible components of $u'$ is $N$, the evaluation map at $z^{+++}_0$
is weakly submersive.
Then the claim follows from Condition \ref{submersivez0} (2), resp. 
Condition  \ref{submersiveatsing}
if $u$ is constant on $S^2_{a(z^+_0)}$, resp. nonconstant on $S^2_{a(z^+_0)}$.

If there are sphere components other than in $A$, then we pick a sphere component $S^2_{a_1}$
($a_1 \notin A$) which contains exactly one nodal point.
Then the claim again follows from Condition \ref{submersivez0} (2), resp. 
Condition  \ref{submersiveatsing}
if $u$ is constant on $S^2_{a_1}$, resp. non-constant on $S^2_{a_1}$.
\par\smallskip
Now we describe the way how we send the obstruction spaces to nearby objects.
For the obstruction bundle which is supported on the component where 
$u$ is nonconstant, the construction is the same as given in  
Section \ref{sec:cyclicKura}.
So we only discuss the case of the disk component on which 
$u$ is constant.
Let ${\bf x} 
= (\Sigma,z_0^+\cup \vec  z^+,u)  \in \mathcal M^{\rm main}_{0;\ell+1}(\beta)$.
We decompose it into extended disk components 
$$
\Sigma = \Sigma_{\alpha_0} \cup \bigcup_{\alpha \in \frak A}\Sigma_{\alpha}
$$ 
as in (\ref{decompxxxx}), where $\Sigma_{\alpha_0}$ contains the $0$-th marked point.
$\Sigma_{\alpha_0}$ together with restriction of $u$ and special points on it
define an object which we write ${\bf x}_{\alpha_0}$.
The object ${\bf x}_{\alpha}$ is obtained from $\Sigma_{\alpha}$ in the same way.
We decompose $\Sigma_{\alpha_0}$ furthermore as 
$$
\Sigma_{\alpha_0} = \Sigma_{\alpha_0,0} \cup \bigcup_{a \in \frak A_{\alpha}}\Sigma_{\alpha_0,a}
$$
where $\Sigma_{\alpha_0,0}$ is a disk and $\Sigma_{\alpha_0,a}$ are spheres.
Then $\Sigma_{\alpha_0,0}$ together with restriction of $u$ and special points on it
define an object, which we write ${\bf x}_{\alpha_0,0}$.
The object ${\bf x}_{\alpha_0,a}$ is obtained from $\Sigma_{\alpha_0,a}$ in the same way.
\par
Let $\text{\bf x}' =  (\Sigma',z_0^{\prime +}\cup \vec  z^{\prime +},u')$.
We consider $\text{\bf x}$ together with $\vec w^+$ which are additional marked points, and 
$\text{\bf x}'$ together with 
$\vec w^{\prime +}$,
such that
$\text{\bf x}'$ is
$(\epsilon_1,\epsilon_2)$-close to $(h\text{\bf x},\vec w^+)$
with respect to $\vec w^{\prime +}$,
in the sense of Definition \ref{primeisclosemodmod}.
We will explain the way how we send
$E({\bf x}_{\alpha_0})$ to $\Sigma'$.
Let $\Sigma'_{\widehat\alpha_0}$ be the extended disk component 
of $\Sigma'$ that contains 
$z^{\prime +}_0$.
We assume that the restriction of $u$ to $\Sigma_{\alpha_0,0}$ is a constant map.
The construction is divided into two steps.
\par\smallskip
\noindent {\bf Step 1}:
We find $h_{\alpha_0}$ in a neighborhood $h$ such that $h_{\alpha_0}$ 
is determined disk component wise. (Namely 
it is independent of the disk components of 
$\Sigma$ other than $\Sigma_{\alpha_0}$ and of the disk components of
$\Sigma'$ other than $\Sigma'_{\widehat\alpha_0}$.)
\par\smallskip
\noindent {\bf Step 2}:
We use $h_{\alpha_0}$ to move $E({\bf x}_{\alpha_0})$ to a subspace of 
$C^{\infty}(\Sigma_{\alpha_0};(h_{\alpha_0}u)^*TX \otimes \Lambda^{0,1})$.
We then move it to 
$C^{\infty}(\Sigma_{\widehat\alpha_0};(u')^*TX \otimes \Lambda^{0,1})$.
\par\smallskip
We first discuss Step 1.
We consider the two cases separately.
\par\smallskip
\noindent
{\bf Case 1}: 
We assume that the restriction of $u$ to $\Sigma_{\alpha_0}$ is  
nonconstant.
(It means that $u$ is nonconstant on certain sphere component $\Sigma_{\alpha_0,a}$.)
In this case we can fix $h_{\alpha_0}$ in the same way as in Section 
\ref{sec:cyclicKura}.
(Note we are either in Case A or Case B.)
\par\smallskip
\noindent
{\bf Case 2}: 
We assume that $u$ is constant on $\Sigma_{\alpha_0}$.  
Then $u(z^+_0) \in L$.
We consider 
\begin{equation}\label{defofhfin45}
f(h') =  {\rm dist}(h'u(z^+_0),u'(z^{\prime +}_0))^2.
\end{equation}
It attains a unique minimum in a neighborhood of $h$. 
This is because $u(z^+_0) \in L$ and $T^n$ action is free on $L$.
Let  $h_{\alpha_0}$ 
be the point where this minimum is attained.  
For Step 2, we consider the following two subcases of Case 2 separately.  

\noindent
{\bf Case 2-1}: There are at least two special points on $\Sigma_{\alpha_0,0}$. 

\noindent
{\bf Case 2-2}: There is only one special point on $\Sigma_{\alpha_0,0}$.  
\par\smallskip
We next discuss Step 2.  
First we consider Case 1 and Case 2-1.  
In these cases, we find that $G({\bf x}_{\alpha_0})$ is finite.  
In a similar way as in Subsection \ref{subsec:1+}, 
we will define
\begin{equation}\label{defiforp}
i_{{\bf v}_{\widehat\alpha,\alpha}(\vec w^+)} : \Sigma_{\alpha_0} 
\setminus \mathcal U(S(\Sigma_{\alpha_0}))
\to 
\Sigma'_{\widehat\alpha_0}
\end{equation}
as follows.
We consider $z^+_{0}$, $\vec z^+_{\alpha_0}$, $\vec w^+_{\alpha_0}$ that are interior 
marked points on $\Sigma_{\alpha_0}$. 
Note at each point $w^+_{\alpha_0,b}$ of $\vec w^+_{\alpha_0}$ we assume $u$
to be an immersion and we fix codimension two submanifolds $N_{w^+_{\alpha_0,b}}$
perpendicular to $u(\Sigma_{\alpha_0})$ at $u(w^+_{\alpha_0,b})$.
We moreover assume that those data are invariant under $G(\bf x_{\alpha_0})$ action.
We use it to define $w^{\prime +}_{\alpha_0,b} \in \Sigma'_{\widehat\alpha_0}$.
(Namely we require $u'(w^{\prime +}_{\alpha_0,b}) \in h_{\alpha_0}N_{w^+_{\alpha_0,b}}$.)
\par
We also note that we have $z^{\prime +}_{0}$, $\vec z^{\prime +}_{\alpha_0}$ 
as a part of the object $\text{\bf x}' = (\Sigma',z_0^{\prime +}\cup\vec z^{\prime +},u')$.
We then put
$$
{\bf v}_{\widehat\alpha_0,\alpha_0}(\vec w^{\prime +})
=
(\Sigma'_{\widehat\alpha_0},z^{\prime +}_{0}\cup\vec z^{\prime +}_{\alpha_0}
\cup \vec w^{\prime +}_{\alpha_0},u').
$$
Let ${\bf v}_{\alpha_0}(\vec w^+)$ be $\Sigma_{\alpha_0}$ together with the marked points $z^{+}_{0}$,  $\vec z^{+}_{\alpha_0}$,
$\vec w^{+}_{\alpha_0}$ on it.
In the same way as Lemma \ref{ifremvethenclose} we can prove
\begin{equation}\label{453form}
{\rm dist}({\bf v}_{\widehat\alpha_0,\alpha_0}(\vec w^{\prime +}),{\bf v}_{\alpha_0}(\vec w^+)) < \frak o(\epsilon_1).
\end{equation}
\par
Thanks to (\ref{453form}) we can use a trivialization of the universal family of a neighborhood of 
${\bf v}_{\alpha_0}(\vec w^+)$ 
in the moduli space of 
marked disks (possibly with sphere bubbles) to obtain a map (\ref{defiforp}).
(See \cite[Section 16]{foootech} for detail.)
\par
We can prove a similar lemma as Lemma \ref{idonotchangemuch}.
Namely 
\begin{equation}\label{4563revvv}
{\rm dist}(i_{{\bf v}'(\vec w^{\prime +})}(x),i_{{\bf v}'_{\widehat \alpha_0, \alpha_0}(\vec w^{\prime +})}(x)) < \frak o(\epsilon_4).
\end{equation}
Here 
${\bf v}'(\vec w^{\prime +})$ 
is $\Sigma'$ together with the marked points $z^{\prime +}_{0}$,  $\vec z^{\prime +}$,
$\vec w^{\prime +}$ on it.
And $i_{{\bf v}'(\vec w^{\prime +})}$ is obtained from ${\bf v}'(\vec w^{\prime +})$   and 
${\bf v}(\vec w^+)$, where ${\bf v}(\vec w^+)$ is $\Sigma$ together with 
$\vec z^+$ and $\vec w^+$.
(\ref{4563revvv}) holds on 
$\Sigma_{\alpha_0} \setminus \mathcal U(S(\Sigma_{\alpha_0})) \setminus U_{\alpha_0}$ where
$U_{\alpha_0}$  is a neighborhood of the boundary singular points.
\par
(\ref{4563revvv}) implies that we have
$$
{\rm dist}(u'(i_{{\bf v}'_{\widehat \alpha_0, \alpha_0}}(x)),h_{\alpha_0}u(x)) < 
\frak o(\epsilon_4)
$$
if $x$ is in the support of $E({\bf x}_{\alpha_0,0})$.
Let $\ell_{x}$ be the shortest geodesic joining $h_{\alpha_0}u(x)$ to $u'(i_{{\bf v}'_{\widehat \alpha_0, \alpha_0}}(x))$.
Now using parallel transport along this geodesic 
and the diffeomorphism $i_{{\bf v}'_{\widehat \alpha_0, \alpha_0}}$ defined in a neighborhood of $x$
we define
\begin{equation}\label{456formpara}
\mathcal P_{\alpha_0,0} :
(h_{\alpha_0})_*E({\bf x}_{\alpha_0,0}) \to C^{\infty}(\Sigma'_{\widehat\alpha_0};
(u')^*TX \times \Lambda^{0,1}).
\end{equation}
\begin{rem}\label{Iandiremimp}
We note that the construction here is different from (\ref{Pprelimi}) in the following important point:
To define (\ref{Pprelimi}) we modify $i_{i,j}$ (that corresponds to $i_{{\bf v}'_{\widehat \alpha_0, \alpha_0}}$)
to another map $I_{i,j}$ by requiring Condition \ref{condsIII}.
But here 
we use $i_{{\bf v}'_{\widehat \alpha_0, \alpha_0}}$ directly. 
Note that we can not require a condition similar to Condition \ref{condsIII} 
in the current situation, since the map $u$ is constant on the support of $E({\bf x}_{\alpha_0,0})$.
\end{rem}
In spite of Remark \ref{Iandiremimp}, the map $\mathcal P_{\alpha_0,0}$ has all the required properties.
In fact, as we mentioned in Remark \ref{rem4571}, we need to take $I$'s in place of $i$'s 
only when the component in question is of Type C. 
In our case where $u$ is constant on the disk component, the object ${\bf x}_{\alpha_0}$ is never 
of Type C.
\par
Now we consider Case 2-2.  
Note because of $S^1$ symmetry of $\Sigma_{\alpha_0,0}$, the map $i_{{\bf v}'_{\widehat \alpha_0, \alpha_0}}$
and hence the linear embedding $\mathcal P_{\alpha_0,0}$ 
is well defined only up to this $S^1$ action of the source.
However we assumed that $E({\bf x}_{\alpha_0,0})$ is invariant under this $S^1$ action.
Thus we obtain the obstruction space as the image of $\mathcal P_{\alpha_0,0}$.
\par
The other part of the proof of Lemma \ref{pmodulikura}
is  the same as the proof of 
Proposition  \ref{transconcl} and so is omitted.
The proof of Lemma \ref{pmodulikura} is completed.
\qed

\par
\subsection{Construction of $\frak t$-Kuranishi structure}
\label{subsec:tkuranishi}
\par
In this subsection we prove Lemma \ref{existstkura}. 
\index{Kuranishi structure!$\frak t$-Kuranishi structure} 
The proof is mostly the same as the 
proof of Proposition  \ref{transconcl}  (construction of the $\frak c$-Kuranishi structure)
and the proof of Lemma \ref{pmodulikura} (construction of the $\frak p$-Kuranishi structure).
\par
We will construct a Kuranishi structure on 
$\mathcal M^\text{\rm main}_{0,0,0;\ell}(\beta)$
by induction on $\beta \cap \omega$.
Note in the description of the boundary 
of $\mathcal M^\text{\rm main}_{0,0,0;\ell}(\beta)$ the moduli space 
$\mathcal M^\text{\rm main}_{k_1,k_2,k_3;\ell'}(\beta')$
with $(k_1,k_2,k_3) = (1,0,0), (0,1,0), (0,0,1)$ appears.
The $\frak t$-Kuranishi structure on it is 
obtained from the case $(k_1,k_2,k_3) = (0,0,0)$ by pull back.
The case $(k_1,k_2,k_3) = (0,0,0)$ is given by induction hypothesis.
\par
Let ${\bf x} = (\Sigma,(z^0_0,z^0_1,z^0_2),\vec z^+,u)$ be an element of 
$\mathcal M^\text{\rm main}_{0,0,0;\ell}(\beta)$. Here $z^0_0,z^0_1,z^0_2$ are unforgetable marked points, 
and $\vec z^+$ are interior marked points.
We decompose $\Sigma$ into extended disk components:
\begin{equation}\label{Sgmadecfort}
\Sigma = \bigcup_{\alpha\in \frak A} \Sigma_{\alpha}.
\end{equation}
We first describe the way how we determine the component $\Sigma_{\alpha_0}$
to which we apply $\frak t$-perturbation.
(To the other components we apply $\frak c$-perturbation.)
\par
For each $\alpha$ and $i=0,1,2$  we define a special point $z^0_{\alpha,i} \in \Sigma_{\alpha}$ as follows.
If $z^0_i \in \Sigma_{\alpha}$, then  $z^0_{\alpha,i} = z^0_i$. If not, 
there is a maximal tree of disk components $\subset \Sigma \setminus \Sigma_{\alpha}$ 
such that $z^0_i$ is contained in it. Then $z^0_{\alpha,i}$ is the root of this tree that is 
a singular point contained in $\Sigma_{\alpha}$.
\begin{lem}\label{alpha0choceunie}
There exists a unique $\alpha_0$ such that $z^0_{\alpha_0,i}$, $i=0,1,2$ are  different to each other.
\end{lem}
The proof is easy and is left to the reader.
We call the component $\Sigma_{\alpha_0}$ (or $\alpha_0$) 
the {\it ${\frak t}$-component} and all the other components 
{\it $\frak c$-component}.
\par
Let $\alpha$ be a $\frak c$-component.
We denote by ${\bf v}_{\alpha}$ the component $\Sigma_{\alpha}$ together with the following marked points:
\begin{enumerate}
\item
$\{ z^0_{\alpha,i} \mid i=0,1,2\}$. Note this set has order $1$ or $2$. We do not double count 
the elements of this set when two of them coincide.
\item $\vec z^+ \cap \Sigma_{\alpha}$.
\end{enumerate}
We denote ${\bf v}_{\alpha}$ together with $u$ by ${\bf x}_{\alpha}$.
\par
When we defined $\frak c$-multisection, we defined a finite dimensional
subspace $E({\bf x}_{\alpha};\frak c)$ of 
$C^{\infty}(\Sigma_{\alpha};u^*TX \otimes \Lambda^{0,1})$.
(They are independent of the boundary marked points.)
\par
Let $\alpha_0$ be the $\frak t$-component.
We decompose $\Sigma_{\alpha_0}$ as 
$$
\Sigma_{\alpha_0} 
= \Sigma_{\alpha_0,0} \cup \bigcup_{a\in \frak A(\alpha_0)}\Sigma_{\alpha_0,a}, 
$$
where $\Sigma_{\alpha_0,0}$ is a disk and 
$\Sigma_{\alpha_0,a}$ are spheres.
Note in our situation $\Sigma_{\alpha_0,0}$ is given 
three boundary singular points $z^0_{\alpha_0,i}$ ($i=0,1,2$) and so 
is always stable. Its automorphism group is trivial.
We take
$E({\bf x}_{\alpha_0,0};\frak t) \subset C^{\infty}(\Sigma_{\alpha_0,0};u^*TX \otimes \Lambda^{0,1})$
such that
the evaluation map
$$
({\rm ev}^0_{\alpha_0,0},{\rm ev}^0_{\alpha_0,1},
{\rm ev}^0_{\alpha_0,2}), : (D_u\overline\partial)^{-1}(E({\bf x}_{\alpha_0,0};\frak t)) \to \bigoplus_{i=0}^2T_{u(z^0_{\alpha_0,i})}L(u)
$$
at the unforgetable marked points is surjective.
We remark that we need to put nontrivial $E({\bf x}_{\alpha_0,0};\frak t)$
in case when $u$ is constant on $\Sigma_{\alpha_0,0}$, also.
\par
On the sphere component $\Sigma_{\alpha_0,a}$ we put 
nontrivial obstruction bundle only in case $u$ is nonconstant there.
We put some additional marked points $\vec w^+$ on the sphere components
which are unstable.
Let ${\bf v}(\vec w^{+})$ and ${\bf v}_{\alpha_0}(\vec w^{+})$
be ${\bf v}$ (resp. ${\bf v}_{\alpha_0}$) plus additional marked points $\vec w^+$ (resp. $\vec w^+\cap \Sigma_{\alpha_0}$).
\par
We have thus described the obstruction space we put on ${\bf x}$.
\par
Let 
${\bf x}' = (\Sigma',(z^{\prime 0}_0,z^{\prime 0}_1,z^{\prime 0}_2),\vec z^{\prime +},u)$
such that $\mu_1{\bf x}'$ is sufficiently close to $h{\bf x}$ 
with respect to $\vec w^{\prime +}$ in the sense of 
Definition \ref{primeisclosemodmod}.
(Here $\mu_1 \in \frak S_{\ell}$.)
For simplicity of notation we assume without loss of generality that $\mu_1 = 1$.
We will describe the way how we send $E({\bf x}_{\alpha_0,0};\frak t)$
to $\Sigma'$.
\par
We decompose $\Sigma'$ as
$$
\Sigma' = \bigcup_{\widehat\alpha \in \frak A'} \Sigma_{\widehat\alpha}.
$$
We can define $\widehat{\alpha}_0$ in the same way as $\alpha_0$ 
using $(z^{\prime 0}_0,z^{\prime 0}_1,z^{\prime 0}_2)$.
We define ${\bf v}'(\vec w^{\prime +})$ and ${\bf v}'_{\widehat\alpha_0}( \vec w^{\prime +})$
in an obvious way. Then we have
$$
i_{{\bf v}'(\vec w^{\prime +}),{\bf v}(\vec w^{+})}(\Sigma_{\alpha_0} \setminus \mathcal U(S(\Sigma_{\alpha_0})) \setminus U_{\alpha_0})
\subset \Sigma'_{\widehat\alpha_0},
$$
where $U_{\alpha_0}$ is an appropriate neighborhood of the 
boundary singular point sets.
\par
Let $\Sigma'_{\widehat\alpha_0, 0}$ be the disk component of $\Sigma'_{\widehat\alpha_0}$
which is defined in a similar way for the case of $\Sigma_{\alpha_0}$. 
There exists a unique biholomorphic map
$$
i_{\widehat\alpha_0,\alpha_0;0} :  \Sigma_{\alpha_0,0}  \to \Sigma'_{\widehat\alpha_0,0}
$$
which sends $(z^{0}_0,z^{0}_1,z^{0}_2)$ to $(z^{\prime 0}_0,z^{\prime 0}_1,z^{\prime 0}_2)$.
The map $i_{\widehat\alpha,\alpha;0}$ is close to $i_{{\bf v}'(\vec w^{\prime +}),{\bf v}(\vec w^{+})}$ on 
$\Sigma_{\alpha_0,0} \setminus U_{\alpha_0}$.
(We can prove it in the same way as Lemma \ref{idonotchangemuch}.)
\par
We first define $h_{\alpha_0}$ which is close to $h$.
In the case when $u$ is nonconstant on $\Sigma_{\alpha_0}$,  
we can define $h_{\alpha_0}$ in the same way as in Subsection \ref{subsec:1+}.
(Namely, according to whether $\Sigma_{\alpha_0}$ is of Type A, B, or C, we use
additional marked points, additional marked points and loops $\gamma_c^h$,
or a circle ${\bf S}_{\alpha_0}$, respectively.)
\par
We assume that $u$ is constant on $\Sigma_{\alpha_0}$.
We take an additional interior marked point $w_0$ on $\Sigma_{\alpha_0,0}$.
For example, we identify $\Sigma_{\alpha_0,0} = D^2$ such that
$z^0_{\alpha,i} = \exp(2\pi i\sqrt{-1}/3)$ for $i=0,1,2$, and we put $w_0 = 0$.
(In case this $w_0$ happens to be a singular point, then we move it to a nearby point.
The important point here is $w_0$ depends only on ${\bf x}_{\alpha_0}$ and not 
on any other components. It should be independent of the boundary singular 
point also. It can depend on the interior singular point since the position of the interior singular point is a part of 
the data consisting ${\bf x}_0$.)
\par
Then since $u(w_0) \in L(u)$ the function
$$
f(h') = {\rm dist}(h'u(w_0),u'(i_{\widehat\alpha_0,\alpha_0;0}(w_0))^2
$$
is strictly convex and attains its minimum at a unique point in a neighborhood of $h$.
We define $h_{\alpha_0}$ to be the point where $f(h')$ attains the minimum.
\begin{rem}
Note 
$ {\rm dist}(hu(w_0),u'(i_{\widehat\alpha_0,\alpha_0;0}(w_0))$ is small since 
$w_0 \in \Sigma_{\alpha_0} \setminus U_{\alpha_0}$.
We cannot use $z^0_{\alpha_0,i}$ in place of $w_0$, since 
$z^0_{\alpha_0,i} \in U_{\alpha_0}$ and
so ${\rm dist}(hu(w_0),u'(i_{\widehat\alpha_0,\alpha_0;0}(w_0))$
may not be small.
\end{rem}
We define $h_{\alpha}$ for $\alpha \ne \alpha_0$ in the same way as Subsection \ref{subsec:1+}.
\par
Now we use $h_{\alpha_0}$ (resp. $h_{\alpha}$) to send the obstruction bundle on the 
components of $\Sigma_{\alpha_0}$ (resp. $\Sigma_{\alpha}$), in case where $u$ is nonconstant on $\Sigma_{\alpha_0}$
(resp. $\Sigma_{\alpha}$).
\par
Namely we first use $i_{{\bf v}'(\vec w^{\prime +}),{\bf v}(\vec w^{+})}$ 
and modify a bit by using a condition similar to Condition \ref{Smapcondsmulticomp}
on the support of the obstruction bundle.
(We assume that $u$ is an immersion there.)
We then use it to send the obstruction bundle.
\par
The situation is different from Subsection \ref{subsec:1+} in the case of the component $\Sigma_{\alpha_0,0}$ on which $u$ is constant.
In this case, we use $i_{\widehat\alpha_0,\alpha_0;0}$ as follows.
Let $x$ be in the support of $E({\bf x}_{\alpha_0,0};\frak t)$.
We take the minimum geodesic $\ell_x$ joining 
$h_{\alpha_0}u(x)$ to $u'(i_{\widehat\alpha_0,\alpha_0;0}(x))$.
Then using parallel transport along $\ell_x$ and 
the biholomorphic map $i_{\widehat\alpha_0,\alpha_0;0}$
we can send $E({\bf x}_{\alpha_0,0};\frak t)$
to $C^{\infty}(\Sigma'_{\widehat\alpha_0,0},(u')^*TX \otimes \Lambda^{0,1})$.
\par
The rest of the proof of Lemma \ref{existstkura} is the same as the proof of Proposition  \ref{transconcl}.
The proof of Lemma \ref{existstkura} is completed.
\qed
\begin{rem}
We do not claim cyclic symmetry of boundary marked points in 
Lemma \ref{form316}. Actually the cyclic symmetry in the naive sense 
does not hold since we distinguish forgetable marked points and 
unforgetable marked points. However we can prove a 
cyclic symmetry among the unforgetable marked points 
in the following sense.
We consider the map
$$
{\rm cyc} : \mathcal M^\text{\rm main}_{k_0,k_1,k_2;\ell}(\beta) \to \mathcal M^\text{\rm main}_{k_1,k_2,k_0;\ell}(\beta)
$$
such that ${\rm cyc} \circ {\rm ev}^0_i = {\rm ev}^0_{i+1}$. Then ${\rm cyc}$ induces an isomorphism of 
$\frak t$-Kuranishi structure.
We do not use this fact in this paper.
\end{rem}
\par
\subsection{Proof of Lemma \ref{pcwkura}}
\label{subsec:lemma426}
\par

In this subsection, we prove Lemma \ref{pcwkura}.
We first remark that we use $\frak q$-Kuranishi structure
on  the disk bubbles. 
\begin{rem}\label{rem457}
\begin{enumerate}
\item
We use the Kuranishi structure produced in Proposition  \ref{transconcl}
to construct $\frak q$-multisection.
The Kuranishi structure constructed in \cite[Section B, Proposition B.7]{fooo08} has 
most of the properties claimed there. However it is not 
disk-component-wise.
\item
Note we need to break cyclic symmetry to construct $\frak q$-{\it multisection}.
However the {\it Kuranishi structure} in Proposition  \ref{transconcl} 
has all the properties 
we need to define $\frak q$-multisection.
(The Kuranishi structure in Proposition  \ref{transconcl} 
is cyclically symmetric moreover. However we do not use this 
extra property to construct $\frak q$-multisection.)
Therefore we can use the $\frak c$-Kuranishi structure 
as $\frak q$-Kuranishi structure also.
\par
Alternatively we can use $\frak q$-Kuranishi structure produced 
in \cite[Section B, Proposition B.7]{fooo08} and then use
the method of Remark \ref{colermethod} to recover the 
compatibility with fiber product.
\end{enumerate}
\end{rem}

\begin{proof}[Proof of Lemma \ref{pcwkura}]
We construct the Kuranishi structure in Lemma \ref{pcwkura}.
Our construction is by induction on $\beta \cap \omega$.
We use the notation and terminology of Section \ref{sec:equikuracot}.
Let
$$
{\bf x} = (\Sigma,\vec z^+,u) \in \mathcal M_{0;\ell+2}^{\text{\rm main}}(\beta;U_0).
$$
We assume:
\begin{assum}\label{asssuYpteIIII}
All the Type I sphere components are in disk bubbles.
\end{assum}
We  decompose
\begin{equation}\label{Sgmatype12doc}
\Sigma = \bigcup_{a\in A} \Sigma_a
\end{equation}
such that $\Sigma_a$ satisfies one of the following conditions 1,2,3:
\begin{conds}\label{decompconds}
\begin{enumerate}
\item
$\Sigma_a$ is a disk component together with tree of sphere components 
of Type I on it. 
The first and the second interior marked points $z_1^+, z_2^+$ are neither on $\Sigma_a$ nor on  
sphere components rooted on $\Sigma_a$.
\item
$\Sigma_a$ is a sphere component of Type II.
\item
$\Sigma_a$ together with all trees of sphere components of Type II is an element of 
$\mathcal M_{1;\ell'+2}^{\text{\rm main}}(\beta';\text{\bf p};U_0)$ for some $\beta'$, $\ell'$.
(Note we do not include sphere components of Type II in $\Sigma_a$.
There may be certain sphere components of Type III included in $\Sigma_a$.)
\end{enumerate}
\end{conds}
We note that Conditions \ref{decompconds}.1, 2, 3 correspond to 
(\ref{diskgeneral}), 
(\ref{spherreq}),
(\ref{disk0comp}),
respectively.
Observe that the sphere bubble on (\ref{diskgeneral}) is necessary of Type I by 
the definition of Type I.
By Assumption \ref{asssuYpteIIII} there is no component of Type I if $\Sigma_a$ does not satisfy Conditions \ref{decompconds}.1, 2.
Therefore we have such a decomposition.
\par
We now describe the way how we put the obstruction spaces in Case 1,2,3 below.
We put ${\bf x}_a = (\Sigma_a, \Sigma_a \cap \vec z^+,u\vert_{\Sigma_a}) 
= (\Sigma_a, \vec z_a^+,u_a) $, 
${\bf v}_a = (\Sigma_a, \vec z_a^+)$.
\par\smallskip
\noindent {\bf Case 1}:
$\Sigma_a$ satisfies Condition \ref{decompconds}.1.
\par
During the construction of $\frak c$-Kuranishi structure in Section \ref{sec:cyclicKura},
we defined
$$
E_a({\bf x}_a) \subset C^{\infty}(\Sigma_a;u^*TX\otimes \Lambda^{0,1}).
$$
We also fixed additional marked points $\vec w^+_a$ 
together with codimension two submanifolds $N_{w^+_{a,b}}$ (in case ${\bf x}_a$ is of Type A or B),
circles ${\bf S}_a$ (in case  ${\bf x}_a$ is of Type C).
\par\smallskip
\noindent {\bf Case 2}:
$\Sigma_a$ satisfies Condition \ref{decompconds}.2.
\par
If $u$ is constant on $\Sigma_a$, then we put $E_a({\bf x}_a) = 0$.
\par
Suppose $u$ is nonconstant on $\Sigma_a$.
We fix $\vec w^+_a \subset \Sigma_a$ such that
${\bf v}_a$ together with $\vec w^+_a$,
(which we denote by ${\bf v}_a(\vec w^+_a)$) is stable.
\par
We also take codimension two submanifold $N_{w^+_{a,b}}$ for each 
$w^+_{a,b} \in \vec w^+_a$.
\par
We will take them so that they are invariant under 
the action of the automorphism group in an obvious sense. Note we do 
not need to choose it so that it is invariant under $T^n$ action, since we do not claim the Kuranishi structure we are constructing 
in this subsection to be invariant under the $T^n$ action.
\par
Let ${\bf v}_a^+(\vec w^+_a)$ (resp. ${\bf v}_a^+$) be  ${\bf v}_a(\vec w^+_a)$ (resp. ${\bf v}_a$) together 
with singular points on $\Sigma_a$ regarded as marked points.
We denote the set of the singular points on $\Sigma_a$ by
$$
\{z_c \mid c \in {\rm Sing}_a\}.
$$
We now take 
$$
E_a({\bf x}_a) \subset C^{\infty}(\Sigma_a;u^*TX \otimes \Lambda^{0,1}),
$$ 
such that the evaluation map
\begin{equation}
{\rm Eval}^{\rm sing} : (D_{u_a}\overline\partial)^{-1}(E_a({\bf x}_a)) 
\to 
\bigoplus_{c \in {\rm Sing}_a} T_{u(z_c)} X
\end{equation}
at the singular points on $\Sigma_a$ is surjective.
\par
We choose $E_a({\bf x}_a)$ so that it is invariant under 
the action of the automorphism group.
We also assume that $u$ is an immersion on a neighborhood of the support of 
$E_a({\bf x}_a)$.
\par\smallskip
\noindent {\bf Case 3}:
$\Sigma_a$ satisfies Condition \ref{decompconds}.3.
\par
We decompose $\Sigma_a$ into irreducible components:
$$
\Sigma_a = \Sigma_{a,0} \cup \bigcup_{\frak a} \Sigma_{a,\frak a}.
$$
Here $\Sigma_{a,0}$ is a disk and $\Sigma_{a,\frak a}$ are spheres.
We define ${\bf x}_{a,0}$, ${\bf x}_{a,\frak a}$, ${\bf v}_{a,0}$, ${\bf v}_{a,\frak a}$
in an obvious way.
\par
We will define $E_{a,0}({\bf x}_{a,0})$ and  $E_{a,\frak a}({\bf x}_{a,\frak a})$.
\par
If $u$ is constant on $\Sigma_{a,0}$ (resp. $\Sigma_{a,\frak a}$),  
we put $E_{a,0}({\bf x}_{a,0})=0$ (resp.   $E_{a,\frak a}({\bf x}_{a,\frak a})=0$).
\par
If $u$ is non-constant on $\Sigma_{a,0}$ (resp. $\Sigma_{a,\frak a}$),
we take additional marked points $\vec w^+_{a,0}$ (resp. $\vec w^+_{a,\frak a
}$)
so that the source becomes stable. We also fix codimension $2$ submanifolds 
$N_{a,\frak a,b}$
(resp. $N_{a,0,b}$).
\par
We then take $E_{a,0}({\bf x}_{a,0})$  so that 
the the image of $D_{u\vert_{a,0}}\overline\partial$ 
and $E_{a,0}({\bf x}_{a,0})$ generate 
$C^{\infty}(\Sigma_{a,0};u^*TX \otimes \Lambda^{0,1})$.
\par
When $u$ is non-constant on the component ${\bf x}_{a,0}$,  $E_{a,0}({\bf x}_{a,0})$ is chosen such 
that the evaluation map at interior singular points is submersive.  
\par
In the case of sphere component we assume 
that the evaluation map
\begin{equation}\label{495formula}
{\rm Eval}^{\rm sing} : (D_{u\vert_{a,\frak a}}\overline\partial)^{-1}(E_{a,\frak a}({\bf x}_{a,\frak a})) 
\to 
\bigoplus_{c \in {\rm Sing}_{a,\frak a}} T_{u(z_c)} X
\end{equation}
at the singular points on $\Sigma_{a,\frak a}$ is surjective.
\par\smallskip
We observe that the evaluation map at the interior singular points 
of non-constant sphere bubbles is surjective as we 
require in Case 3.
We also observe that the evaluation map at each boundary singular point  
of any disk bubbles is surjective by the $T^n$ equivariance.
Therefore we can show the surjectivity of the linearized operator after 
taking the fiber products of moduli spaces associated to irreducible components.
\par\smallskip
Let ${\bf x}' = (\Sigma',\vec z^{\prime +},u')$ be an object. We assume that it is close to 
$({\bf x},\vec w^+)$ with respect to $\vec w^{\prime +}$ in the following sense.
\par
We took $\vec w^+_a$ for various components of $\Sigma$ 
in (\ref{Sgmatype12doc}).
In case $\Sigma_{a}$ is of Type I and is of Type C, we add $\vec w^+_a$ also.
(Later on it turns out that 
the construction is independent of the choice of 
$\vec w^+_a$ in case $\Sigma_{a}$ is of Type I and is of Type C.
But we temporally fix this choice. See Remark \ref{rem4571}.)
The additional marked point $\vec w^+$ is the union of all of them.
\par
Note each member $w^+_a$ of $\vec w^+$ comes with 
a codimension $2$ submanifold $N_{w^+_a}$ of $X$.
\par
For each member $w^+_a$ of $\vec w^+$ we are supposed to have a 
corresponding member $w^{\prime +}_a$ of $\vec w^{\prime +}$.
We require:
\begin{conds}\label{xprimeclosex453}
\begin{enumerate}
\item
$
u'(w^{\prime +}_a) \in N_{w^+_a}.
$
\item
Let ${\bf v}'(w^{\prime +})$ (resp. ${\bf v}(w^{+})$) be $(\Sigma',\vec z^{\prime +})$ 
together with additional marked points $w^{\prime +}$ 
(resp. ${\bf v}$ together with $\vec w^+$).
We require ${\bf v}'(w^{\prime +})$ to be close to ${\bf v}(w^{+})$ in the moduli space of
marked bordered curves.
\item
By Item 2 above, there exists a map
$$
i_{{\bf v}'(w^{\prime +});{\bf v}(w^{+})} : \Sigma\setminus \mathcal U(S(\Sigma))
\to \Sigma',
$$
where $\mathcal U(S(\Sigma))$ is a neighborhood of the singular point set of 
$\Sigma$. This map is obtained by a local trivialization of the universal family 
of the moduli space of marked (bordered) curves and a choice of the coordinates around the 
singular points of $\Sigma$. (See (\ref{15apend3}).)
\par
We require that the composition $u' \circ i_{{\bf v}'(w^{\prime +});{\bf v}(w^{+})}$ is 
$C^1$ close to $u$ on $\Sigma\setminus \mathcal U(S(\Sigma))$.
\item
The diameter of the image by $u'$ of each connected component of $\mathcal U(S(\Sigma))$ is 
small.
\end{enumerate}
\end{conds}
\begin{rem}
The above definition is similar to 
Definition \ref{primeisclosemodmod}.
We however do not move ${\bf x}$ by $h \in T^n$.
Note that 
we do not claim the Kuranishi structure we are constructing on $\mathcal M_{0;\ell'+2}^{\text{\rm main}}(\beta';U_0)$
has $T^n$ action.
\par
On the other hand, we do claim that the Kuranishi structure we are constructing on $\mathcal M_{0;\ell'+2}^{\text{\rm main}}(\beta';U_0)$
is invariant under the permutation of the interior marked points.
In the situation of Section \ref{sec:cyclicKura}, the $\frak S_{\ell}$ action is mixed up with the $T^n$ action,
which makes it harder to study. Since in our case we only has $\frak S_{\ell}$ action, it is easier to handle.
So we omit the argument about $\frak S_{\ell}$ invariance.
\end{rem}
We now describe the way how we send the obstruction space we fixed on 
${\bf x}$ to $\Sigma'$.
We do so for each component $\Sigma_a$ as in Condition \ref{decompconds}.
\par\smallskip
\noindent
{\bf Case 1}:
$\Sigma_a$ satisfies Condition \ref{decompconds}.1.
\par
Note, for the consistency with the $\frak c$-Kuranishi structure, we first need to 
find $h_{a} \in T^n$ in a neighborhood of the unit. Using the data $\vec w_a^+$ or ${\bf S}_a$ we are given, 
we can find such $h_{a}$ in the same way as in Subsection \ref{subsec:1+}.
\par
Then we can send $E_{a}({\bf x}_a)$ to $C^{\infty}(\Sigma';(u')^*TX\otimes \Lambda^{0,1})$ 
in the same way as (\ref{Pprelimi}).
Namely we first modify $i_{{\bf v}'(w^{\prime +});{\bf v}(w^{+})}$ to $I_a$ on the support of $E_{a}({\bf x}_a)$ 
by requiring a condition similar to Condition \ref{condsIII} then use parallel transport along the minimal 
geodesic to send  $E_{a}({\bf x}_a)$.
\par\smallskip
\noindent
{\bf Case 2}:
$\Sigma_a$ satisfies Condition \ref{decompconds}.2.
\par
In this case we do not need to find $h_a$. We replace $i_{{\bf v}'(w^{\prime +});{\bf v}(w^{+})}$ by $I_a$ 
on the support of $E_{a}({\bf x}_a)$ and use it to send $E_{a}({\bf x}_a)$ in the same way as (\ref{Pprelimi}).
\par\smallskip
\noindent
{\bf Case 3}:
$\Sigma_a$ satisfies Condition \ref{decompconds}.3.
\par
This is the same as Case 2.
\par\smallskip
We have thus constructed the required Kuranishi neighborhood of ${\bf x}$.
\begin{rem}
We remark that the evaluation map
$$
({\rm ev}^+_1,{\rm ev}^+_2) :
V_{\bf x} \to X^2
$$
at the two distinguished interior marked points is not a submersion on our 
Kuranishi neighborhood $V_{\bf x}$.
(Note $V_{\bf x}$ is the set of ${\bf x}'$ such that $\overline\partial u' \equiv 0
\mod E({\bf x}')$. Here $E({\bf x}')$ is the sum of the images of various $E({\bf x}_a)$.)
\par
We however can define a fiber product of it with the $T^n$ equivariant cycles ${\bf p}(i)$ 
in the same way as in Subsection \ref{subsec:fiberproduct}.
\end{rem}
\par\smallskip
To construct a Kuranishi structure (which is a system of Kuranishi neighborhoods and coordinate changes) 
we proceed in the same way as \cite[end of 1003-begining 1004]{FO} (See \cite[Section 18]{foootech} for detail.) 
as follows.
Our argument is also similar to Subsection \ref{subsec:coordinatechange}.
So our discussion below is rather sketchy.
\par
We first note that while we construct the $\frak c$-Kuranishi structure
we take a finite subset $\frak P_{\ell}(\beta) \subset \mathcal M_{0;\ell}^{\text{\rm main}}(\beta)/(T^n \times \frak S_{\ell})$ 
and, for each $\frak c \in \frak P_{\ell}(\beta)$, a data we need to put obstruction bundle on it, (that is, $\vec w_{\frak c}^+$, $N_{w^{\frak c,a}}$, 
${\bf S}_{\frak c}$, $E_\frak c$ etc.) and a representative ${\bf x}_{\frak c}$.
\par
We also take a finite subset $\frak P_{\ell}^{\rm sph}(\alpha)$
of the moduli space of pseudo-holomorphic sphere
$\mathcal M_{\ell}(\alpha)/\frak S_{\ell}$ and
data to define obstruction bundle for each $\frak c \in \frak P_{\ell}^{\rm sph}(\alpha)$.
We require that $E_{\frak c}$ is $0$ if $\alpha =0$.
We also require (\ref{495formula}) to be surjective.
We use it to define a system of component-wise Kuranishi structures on $\frak P_{\ell}^{\rm sph}(\alpha)$
in the same way as in Subsection \ref{subsec:coordinatechange}.
\par
Now
we will take a finite subset $\frak P_{\ell}(\beta,U_0)$  
of $\mathcal M_{0;\ell+2}^{\text{\rm main}}(\beta;U_0)$
by induction on $\beta\cap \omega$ and take various data on it 
to define a Kuranishi structure on $\mathcal M_{0;\ell+2}^{\text{\rm main}}(\beta;U_0)$, below.
\begin{rem}
Those finite sets $\frak P_{\ell}^{\rm sph}(\alpha)$, $\frak P_{\ell}(\beta,U_0)$ 
are the subsets of  the moduli space itself 
(divided by the symmetric group $\frak S_{\ell}$ only)
and not of its quotient by 
$T^n$ action. (The Kuranishi structures on those moduli spaces 
are not required to carry $T^n$ action.)
\end{rem}
\par
We consider 
$\text{\bf x} = (\Sigma,\vec z^+,u) \in \mathcal M_{0;\ell+2}^{\text{\rm main}}(\beta;U_0)$, 
for which in the decomposition 
(\ref{Sgmatype12doc}) there exists at least one component satisfying 
Condition \ref{decompconds}.1.
(Namely we assume that there exists at least one disk bubble.) Suppose also that Assumption \ref{asssuYpteIIII} is satisfied.
Let ${\bf x}' = (\Sigma',\vec z^{\prime +},u')$ be an object which is close to $\text{\bf x}$ in the sense of Condition \ref{xprimeclosex453}.
\par
Let ${\bf x}_a$ be an object which we obtain by restricting the data on ${\bf x}$ to $\Sigma_a$.
\par\smallskip
\noindent
{\bf Case 1}:
$\Sigma_a$ satisfies Condition \ref{decompconds}.1.
\par
Let $\frak c_a \in \frak C({\bf x}_a)$. This means that 
the element ${\bf x}_a \in \mathcal M_{0;\ell_a}^{\text{\rm main}}(\beta_a)$
is in a small neighborhood $\frak U(\frak c_a)$ of 
$T^n \times \frak S_{\ell_a}$ orbit of $\frak c_a \in \frak P_{\ell_a}(\beta_a)$.
So there exists $h_{\frak c_a} \in T^n$, $\mu_{\frak c_a} \in \frak S_{\ell_a}$ such that $\mu_{\frak c_a}{\bf x}_a$ is close to $h_{\frak c_a}\frak c_a$.
We can then send $E_{\frak c_a}$ to $\Sigma'$
in the same way as above, by 
using ${\bf x}(\frak c_a)$ where we replace ${\bf x}_a$ by $h_{\frak c_a}\frak c_a$ in ${\bf x}$.
We thus obtain
$$
E_{\frak c_a}({\bf x}') \subset C^{\infty}(\Sigma',(u')^*TX \otimes \Lambda^{0,1}).
$$
\par\smallskip
\noindent
{\bf Case 2}:
$\Sigma_a$ satisfies Condition \ref{decompconds}.2.
\par
Let us consider $\frak c_a \in \frak P_{\ell_a}^{\rm sph}(\alpha_a)$ such that 
${\bf x}_a \in \frak U(\frak c_a)$.
In other words $\frak c_a \in \frak C({\bf x}_a)$.
We can then send $E_{\frak c_a}$ in the same way as above to obtain
$$
E_{\frak c_a}({\bf x}') \subset C^{\infty}(\Sigma',(u')^*TX \otimes \Lambda^{0,1}).
$$
\par\smallskip
\noindent
{\bf Case 3}:
$\Sigma_a$ satisfies Condition \ref{decompconds}.3.
\par
We consider $\Sigma_a$ together Type II sphere bubbles on it.
We thus obtain an element 
$\widehat{\bf x}_a \in \mathcal M_{0;\ell_a+2}^{\text{\rm main}}(\beta_a;U_0)$.
\par
Since there exists at least one disk bubble, we have
$$
\beta_a \cap \omega < \beta \cap \omega.
$$
Therefore
by induction hypothesis, we have $\frak P_{\ell_a}(\beta_a,U_0)$ and 
data to define obstruction bundles on each element of it.
(On the Type II sphere components, those data are the same 
as we described in Case 2. This is a part of induction hypothesis.)
\par
Let $\frak c_a \in \frak P_{\ell_a}(\beta_a,U_0)$ 
such that $\widehat{\bf x}_a \in \frak U(\frak c_a)$.
In other words, $\frak c_a \in \frak C(\widehat{\bf x}_a)$.
We have defined $E_{a,\frak a_0}$ or $E_{a,\frak a}$ for each 
irreducible component of ${\bf x}_a$ and various data to 
fix obstruction bundle.
We send them to ${\bf x}'$ in the same way as before.
We thus obtain
$$
E_{\frak c_a}({\bf x}') \subset C^{\infty}(\Sigma',(u')^*TX \otimes \Lambda^{0,1}).
$$
\par\medskip
Summing up these three cases, we define
$E_{\bf x}({\bf x}')$ by
$$
E_{\bf x}({\bf x}') 
=
\bigoplus_{a : \text{Case 1}} \bigoplus_{\frak c_a \in \frak C({\bf x}_a)}E_{\frak c_a}({\bf x}') 
\oplus
\bigoplus_{a : \text{Case 2}} \bigoplus_{\frak c_a \in \frak C({\bf x}_a)}E_{\frak c_a}({\bf x}') 
\oplus
\bigoplus_{a : \text{Case 3}} \bigoplus_{\frak c_a \in \frak C({\bf x}_a)}E_{\frak c_a}({\bf x}').
$$
We use it to define a Kuranishi neighborhood of ${\bf x}$.
We have a coordinate change among these Kuranishi neighborhoods.
\par
By choosing $\epsilon$ in Definition \ref{def4232} sufficiently small, we may assume that 
the union of them consists a Kuranishi neighborhood of $\partial\mathcal M_{0;\ell+2}^{\text{\rm main}}(\beta;U_0)$.
In fact, $\partial\mathcal M_{0;\ell+2}^{\text{\rm main}}(\beta;U_0)$ consists of elements ${\bf x}$ which 
have at least one disk bubble. We may choose $\epsilon$ in Definition \ref{def4232} 
small so that all the 
Type I sphere bubbles are obtained by gluing from the sphere bubbles which are in the tree of sphere bubbles 
rooted on a disk bubbles. This implies that we can cover $\partial\mathcal M_{0;\ell+2}^{\text{\rm main}}(\beta;U_0)$
by taking the Kuranishi neighborhoods of ${\bf x}$ which satisfies Assumption \ref{asssuYpteIIII}.
\par
We have thus constructed a Kuranish structure on a neighborhood of $\partial\mathcal M_{0;\ell+2}^{\text{\rm main}}(\beta;U_0)$.
We can then take $\frak P_{\ell}(\beta,U_0)$ and various data defining obstruction bundles, 
to obtain a Kuranishi structure of $\mathcal M_{0;\ell+2}^{\text{\rm main}}(\beta;U_0)$
outside a small  neighborhood of $\partial\mathcal M_{0;\ell+2}^{\text{\rm main}}(\beta;U_0)$.
We glue these two in the same way as in Subsection \ref{subsec:coordinatechange} to obtain 
a Kuranishi structure on $\mathcal M_{0;\ell+2}^{\text{\rm main}}(\beta;U_0)$.
\par
Therefore we have constructed a system of Kuranishi structures on $\mathcal M_{0;\ell+2}^{\text{\rm main}}(\beta;U_0)$
that satisfies the required properties.
The proof of Lemma \ref{pcwkura} 
is now complete.
\end{proof}

\par
\subsection{Kuranishi structure on the moduli space 
of pseudo-holomorphic annuli I}
\label{subsec:kuraanul}
\index{Kuranishi structure!Kuranishi structure on moduli of holomorphic annuli}
In this and the next two subsections, we construct a Kuranishi structure 
on the moduli space of pseudo-holomorphic annuli.
In this subsection we prove Lemma \ref{bdryMannu}. 
\par
Let ${\bf x} 
= (\Sigma,(z_1,z_2),\vec z^+,u) \in \mathcal M_{(1,1);\ell}^{{\rm annu};{\rm main}}(\beta)$.
We will define a Kuranishi neighborhood of ${\bf x}$.
In this subsection, 
we study the case 
$\frak{forget}({\bf x}) \in \mathcal M_{(1,1);0}^{{\rm annu};{\rm main}}$
is in the stratum 1 in Lemma \ref{stratifyAn11}.
(This is the case when $\Sigma$ is a (nonsingular) annulus together 
with sphere bubbles.)
\par
We decompose $\Sigma$ as
\begin{equation}\label{decomp454}
\Sigma = \Sigma_0 \cup \bigcup_{a} \Sigma_a,
\end{equation}
where $\Sigma_0$ is an annulus together with 
sphere bubbles rooted on it, 
and $\Sigma_a$ are extended disk components.
Then $\Sigma_a$ together with marked points on it 
and the restriction of $u$ define an element 
${\bf x}_a \in \mathcal M_{0;\ell_a}(\beta_a)$.
\par
We define $z^0_1, z^0_2 \in \partial \Sigma_0$ as follows.
If $z_i \in \Sigma_0$, then $z^0_i =z_i$.
If not, $z^0_i$ is the root of the tree of disk components 
containing $z_i$.
We define 
${\bf x}_0 \in \mathcal M_{(1,1);\ell_0}^{{\rm annu};{\rm main}}(\beta(0))$ 
to be
$(\Sigma_0,(z^0_1,z^0_2),\vec z^+ \cap \Sigma_0,u\vert_{\Sigma_0})$.
We also put
${\bf v}_0 = (\Sigma_0,(z^0_1,z^0_2),\vec z^+ \cap \Sigma_0)$.
\par
We take $\frak c(a) = (\Sigma_{\frak c(a)},\vec z^+_{\frak c(a)},u_{\frak c(a)})
\in \frak P_{\ell_a}(\beta_a) \subset \mathcal M_{0;\ell_a}(\beta_a)/(T^n \times \frak S_{\ell_a})$
such that ${\bf x}_a \in \frak U(\frak c(a))$.
\par
While we constructed $\frak c$-Kuranishi structure, 
we chose additional marked points $\vec w^+_{\frak c(a)}$,
and $N_{w_{\frak c(a),b}}$ (in case $\frak c(a)$ is of Type A or B),
${\bf S}_{\frak c(a)}$ (in case $\frak c(a)$ is of Type C).
We  chose $\vec w^+_{\frak c(a)}$
and $N_{w_{\frak c(a),b}}$ 
in Type C case also temporary. (The 
resulting Kuranishi neighborhood we will construct will be independent of the 
choice of them in Type C case.)
\par
By assumption, $\mu_{\frak c(a)}{\bf x}_a$ is close to $(h_{\frak c(a)}
{\frak c(a)},\vec w^+_{\frak c(a)})$ 
with respect to $\vec w^{{\bf x},+}_{\frak c(a)}$, 
where $\mu_{\frak c(a)} \in \frak S_{\ell_a}$,
$h_{\frak c(a)} \in T^n$. 
We also fixed 
$$
E_{\frak c(a)} \subset C^{\infty}(\Sigma_{\frak c(a)},u_{\frak c(a)}^*TX \otimes 
\Lambda^{0,1}).
$$
We next describe the data we need to define an obstruction bundle on 
$\Sigma_0$. 
We will define a finite set 
$\frak P_{\ell}^{{\rm annu}}(\beta)
\subset \mathcal M_{(1,1);\ell}^{{\rm annu};{\rm main}}(\beta)/\frak S_{\ell}$
by an induction on $\beta\cap \omega$.
\par
Let 
$\frak c(0)
= 
(\Sigma_{\frak c(0)},(z^{\frak c(0)}_1,z^{\frak c(0)}_2),\vec z_{\frak c(0)}^+, u_{\frak c(0)}) 
\in \frak P_{\ell_0}^{{\rm annu}}(\beta(0))$ 
which is close to ${\bf x}_0$.
By assumption 
the decomposition (\ref{decomp454}) is trivial, i.e., only one component, for $\frak c(0)$.
\par
Let 
$$
\Sigma_{\frak c(0)} 
= \Sigma_{\frak c(0),0} \cup \bigcup_{\frak a} \Sigma_{\frak c(0),\frak a}
$$
be the decomposition into the irreducible components.
(Here $\Sigma_{\frak c(0),0}$ is an annulus and 
$\Sigma_{\frak c(0),\frak a}$ are spheres.)
Then $\Sigma_{\frak c(0),0}$ (resp. $\Sigma_{\frak c(0),\frak a}$)
together with data induced from $\frak c(0)$ gives an object 
which we denote by $\frak c(0)_0$ (resp. $\frak c(0)_{\frak a}$).
\par
We take additional interior marked points
$$
\vec w^+_{\frak c(0)} 
= 
\vec w^+_{\frak c(0),0} \cup \bigcup_{\frak a} \vec w^+_{\frak c(0),\frak a},
$$
where $\vec w^+_{\frak c(0),0} \subset \Sigma_{\frak c(0),0}$,
$\vec w^+_{\frak c(0),{\frak a}} \subset \Sigma_{\frak c(0),{\frak a}}$.
We assume that they have the following properties:
\begin{enumerate}
\item
The source becomes stable after we add $\vec w^+_{\frak c(0)}$.
Exception: If $u_{\frak c(0)}$ is constant on the annulus component 
$\Sigma_{\frak c(0),0}$,
then  $\vec w^+_{\frak c(0)} \cap \Sigma_{\frak c(0),0} = \emptyset$.
\item
$\vec w^+_{\frak c(0),\frak a}$ is an empty set if $u_{\frak c(0)}$ is 
constant on $\Sigma_{\frak c(0),\frak a}$.
\item
$u_{\frak c(0)}$ is an immersion in a neighborhood of $\vec w^+_{\frak c(0)}$.
\end{enumerate}
We also choose 
a codimention $2$ submanifold 
$N_{w^+_{{\frak c}(0),0, b}}$
(resp.  $N_{w^+_{{\frak c}(0),\frak a, b}}$) 
of $X$ which intersects transversally 
to $u_{\frak c(0)}(\Sigma_{\frak c(0),0})$
(resp. $u_{\frak c(0)}(\Sigma_{\frak c(0),a})$)
at $u_{\frak c(0)}( w^+_{\frak c(0),0,b})$
(resp. $u_{\frak c(0)}( w^+_{\frak c(0),\frak a,b})$).
\par
We define a finite dimensional subspace 
\begin{equation}\label{defnofE454}
E_{\frak c(0)} 
= E_{\frak c(0),0} \oplus \bigoplus_{\frak a} E_{\frak c(0),\frak a} 
\subset C^{\infty}(\Sigma_{\frak c(0)},u_{\frak c(0)}^*TX \otimes 
\Lambda^{0,1})
\end{equation}
as follows:
\begin{enumerate}
\item
We assume $E_{\frak c(0)}$ and the image of 
$D_{u_{\frak c(0)}}\overline\partial$ generate  
$C^{\infty}(\Sigma_{\frak c(0)},u_{\frak c(0)}^*TX \otimes 
\Lambda^{0,1})$.
\item
If $u_{\frak c(0)}\vert_{\Sigma_{\frak c(0),0}}$ 
(resp. $u_{\frak c(0)}\vert_{\Sigma_{\frak c(0),\frak a}}$) is nonconstant, 
we assume that $u_{\frak c(0)}$ is an immersion on the support of 
$E_{\frak c(0),0}$ (resp. $E_{\frak c(0),\frak a}$).
\item
If $u_{\frak c(0)}\vert_{\Sigma_{\frak c(0),\frak a}}$ is constant,  
we assume that $E_{\frak c(0),\frak a} = 0$.
\item
$E_{\frak c(0)}$ is invariant under the action of 
${\rm Aut}_+(\frak c(0))$.
Here ${\rm Aut}_+(\frak c(0))$ is the set of biholomorphic maps $\varphi : \Sigma_{\frak c(0)} \to \Sigma_{\frak c(0)}$ such that 
$u_{\frak c(0)} \circ \varphi = u_{\frak c(0)}$ and that 
$\varphi$ preserves $\vec z_{\frak c(0)}^+$ as a set.
\item
$E_{\frak c(0),0}$ 
(resp. $E_{\frak c(0),\frak a}$)  is invariant under the action of 
${\rm Aut}_+(\frak c(0)_{0})$ (resp. ${\rm Aut}_+(\frak c(0)_{\frak a})$).
\end{enumerate}
\par
By Conditions 1 and 5, we find that $E_{\frak c(0)}$ is $S^1$-invariant, when $u$ is constant on the annulus component 
$\Sigma_{\frak c(0),0}$.  
We assume that $\mu_{\frak c(0)}{\bf x}_0$ is close to $(\frak c(0),\vec w_{\frak c(0)})$ 
with respect to $w^{{\bf x} +}_{\frak c(0)}$ in the same sense as 
Definition \ref{primeisclose}.
\par
We consider ${\bf y} = (\Sigma',(z^{\prime}_1,z^{\prime}_2),\vec z^{\prime +},u')$
which is close to $({\bf x},\vec w^{{\bf x}+})$ 
with respect to $\vec w^{\prime +}$.
Here $\vec w^{{\bf x}+} = \vec w^{{\bf x}+}_{\frak c(0)}
\cup \bigcup_a \vec w^{{\bf x}+}_{\frak c(a)}$.
Note the choice of $\vec w^{{\bf x}+}$ depends on $\frak c$.
So this definition of closed-ness of ${\bf y}$ to ${\bf x}$ depends on $\frak c$.
However we have the following lemma. 
Here $\vec w^{{\bf x}+}_{(i)}$, $\vec w^{\prime +}_{(i)}$
are  $\vec w^{{\bf x}+}$, $\vec w^{\prime +}$ when we consider 
$\frak c^{(i)}(\frak a)$, $\frak c^{(i)}(0)$
in place of $\frak c(\frak a)$, $\frak c(0)$.
\par
\begin{lem}
For each $\epsilon > 0$, $\frak c^{(1)}(\frak a)$, $\frak c^{(1)}(0)$, 
$\frak c^{(2)}(\frak a)$, $\frak c^{(2)}(0)$,  there exists $\epsilon'>0$ 
with the following properties. 
If ${\bf y}$ is $\epsilon'$ close to $({\bf x},\vec w^{{\bf x}+}_{(1)})$ 
with respect to $\vec w^{\prime +}_{(1)}$, then 
there exists  $\vec w^{\prime +}_{(2)}$ such that 
${\bf y}$ is $\epsilon$ close to $({\bf x},\vec w^{{\bf x}+}_{(2)})$ 
with respect to $\vec w^{\prime +}_{(2)}$.
\end{lem}
The proof of the lemma is easy and is omitted.
We will not mention a similar lemma in later subsections.
\par
We fix $\frak c(a)$ and $\frak c(0)$ for a while.
As in the previous sections, we denote by 
${\bf x}(\vec w^{{\bf x}+})$ 
the object obtained by adding $\vec w^{{\bf x}+}$ to 
${\bf x}$ as data of interior marked points.
We define ${\bf y}(\vec w^{\prime +})$ similarly.
\par
We glue $h_{\frak c(a)}\frak c(a)$ and $\frak c(0)$ and obtain 
${\bf x}_{\frak c} = (\Sigma_{\frak c},(z^{\frak c}_1,z^{\frak c}_2),
\vec z^+_{\frak c},u_{\frak c})$
such that
$\mu_{\frak c}{\bf y}$ is close to $({\bf x}_{\frak c} ,\vec w^+_{\frak c})$
with respect to $\vec w^{\prime +}$. 
(Here $\mu_{\frak c} \in \frak S_{\ell}$.) 
\begin{rem}\label{smallambrem}
We note that the way how we glue $\frak c(a)$ and $\frak c(0)$ 
to obtain ${\bf x}_{\frak c}$ is {\it not} canonical.
Namely we glue two of the components at their 
boundaries at the point which is close to the point where 
the corresponding components are glued in ${\bf x}$.  
We use cut-off functions in this process.  
But there is no canonical way to determine the precise position
where we glue  $\frak c(a)$ with $\frak c(a')$ or $\frak c(0)$.
\par
By this reason there is an ambiguity of the definition of ${\bf x}_{\frak c}$.
However the difference between the two choices of ${\bf x}_{\frak c}$
is small and the argument below is carefully designed so that 
this ambiguity will be removed during the construction.
\end{rem}
We have a map
\begin{equation}
i_{{\bf y}(\vec w^{\prime +})} :
\Sigma_{\frak c}\setminus \mathcal U(S(\Sigma_{\frak c})) \to \Sigma'.
\end{equation}
\begin{rem}\label{smallambrem2}
By the reason explained in Remark \ref{smallambrem},
$i_{{\bf y}(\vec w^{\prime +})}$ is well-defined only up to small ambiguity.
\end{rem}
We now describe the way how we send $E_{\frak c}$ to ${\bf y}$.
\par\smallskip
\noindent
{\bf Case 1}: 
Suppose $u_{\frak c(0)}$ is nonconstant on $\Sigma_{\frak c(0),0}$.
We replace $i_{{\bf y}(\vec w^{\prime +})}$ by a map $I$ 
on the support of $E_{\frak c(0),0}$ in the same way as Condition 
\ref{condsIato}
and send $E_{\frak c(0),0}$ to ${\bf y}$ in the same way as 
(\ref{Pprelimi}).
\par
The construction is the same for $E_{\frak c(0),\frak a}$.
(Note we assume that $E_{\frak c(0),\frak a}$  is zero if $u$ is constant there.)
\par
We note that the ambiguity we mentioned in Remarks \ref{smallambrem},
\ref{smallambrem2} is removed during the process to go from $i$ to $I$.
We also note that the map $I$ depends only on ${\bf x}_{\frak c(0)}$ 
and ${\bf y}$ and is independent of ${\bf x}$.
\par\smallskip
\noindent
{\bf Case 2}: Suppose $u_{\frak c(0)}$ is constant on $\Sigma_{\frak c(0),0}$.
We consider the universal family
\begin{equation}\label{annunivfam}
\mathfrak M^{{\rm ann}}_{(1,1),0} \to \mathcal M^{{\rm ann}}_{(1,1),0}
\end{equation}
of annuli with one marked point on each boundary component.
In a neighborhood of $(\Sigma_{\frak c(0),0},(z^{\frak c(0)}_1,z^{\frak c(0)}_2))$
we take a trivialization of (\ref{annunivfam}) as a smooth fiber bundle.
It induces a diffeomorphism 
\begin{equation}\label{anndiffeobyunivfac}
i_{\Sigma'_0} : \Sigma_{\frak c(0),0} \to \Sigma'_0.
\end{equation}
Here $\Sigma'_0$ is a component of $\Sigma'$ which is diffeomorphic to an annulus.
\par
We use $i_{\Sigma'_0}$ together with the parallel transport along the geodesic joining 
$u(x)$ to $u'(i_{\Sigma'_0}(x))$ to send $E_{\frak c(0),0}$ to ${\bf y}$.
\par
We note that the local trivialization of (\ref{annunivfam}) is a part of the data 
which we fix for each $\frak c(0)$. Therefore the map (\ref{anndiffeobyunivfac}) 
depends only on $\frak c(0)$ and ${\bf y}$.
\par\smallskip
\noindent
{\bf Case 3}: Suppose $\frak c(a)$ is of Type A.
Note we assumed that $\mu_{\frak c(a)}{\bf x}_a$ is close to $h_{\frak c(a)}\frak c(a)$.
For the simplicity of notation we assume without loss of generality 
that $\mu_{\frak c(a)} = 1$.
\par
Note the choice of this $h_{\frak c(a)}$ is not canonical. 
We need to replace $h_{\frak c(a)}$ by $h_a$ that is determined uniquely by 
$\frak c(a)$ and ${\bf y}$.
\par
If $h$ is close to $h_{\frak c(a)}$, then we find $w^{h +}_{a,b}$ for each $b$ by 
requiring  $u'(w^{h +}_{a,b}) \in hN_{a,b}$.
\par
We also assume that $w^{h +}_{a,b}$ is close 
to $i_{{\bf y}(\vec w^{\prime +})}(w^{+}_{a,b})$.
Note the choice of $w^{h +}_{a,b}$ is 
independent of the ambiguity mentioned in
Remarks \ref{smallambrem}, \ref{smallambrem2}.
\par
We define the function $f(h)$ by  (\ref{deffbyave}).
Then we define $h_a$ so that $f$ attains the unique minimum at 
$h_a$.
\par
We next replace $i_{{\bf y}(\vec w^{\prime +})}$ by $I$ 
in the same way as Condition 
\ref{condsIato}
and send $E_{\frak c(a)}$ to ${\bf y}$ in the same way as 
(\ref{Pprelimi}).
\par\smallskip
\noindent
{\bf Case 4}: Suppose $\frak c(a)$ is of Type C.
We first find $h_a$. We fixed ${\bf S}_a \subset \Sigma_{\frak c(a),0}$ 
as a part of the data we took for $\frak c(a)$.
Then, on ${\bf S}_a$, we replace $i_{{\bf y}(\vec w^{\prime +})}$ by $\frak I^{h}$ 
by requiring Condition \ref{Smapconds}.
We then use it to define $f(h)$ by (\ref{formula4331}).
Thus we obtain $h_a$ by requiring that $f$ attains its unique minimum there.
The rest of the construction is the same as Case 3.
\par\smallskip
\noindent
{\bf Case 5}: Suppose $\frak c(a)$ is of Type B.
Let $\Sigma'_{\widehat{a}}$ be the component of $\Sigma'$ 
containing 
$i_{{\bf y}(\vec w^{\prime +})}(\Sigma_{\frak c(a)} \setminus \mathcal U(S(\Sigma_{\frak c(a)})))$.
\par\smallskip
\noindent
{\bf Case 5-1}:
Suppose $\Sigma'_{\widehat{a}}$ is an extended disk component.
We define $f_+$ in the same way as Case 4 and Subsection
\ref{subsec:onedisk2}.
We use it to define $h_a$. The rest is the same as Case 3.
\par\smallskip
\noindent
{\bf Case 5-2}:
Suppose $\Sigma'_{\widehat{a}}$ is an annulus plus sphere bubbles.
Most of the proof is the same as Case 5-1. However,  
we need to modify the way how we associate a loop $\gamma^h_c 
\in \Sigma'_{\widehat a}$ to each interior singular point 
$z_c \in \Sigma_a$ (which does not have corresponding singular 
point in $\Sigma'$).
In fact, the construction of $\gamma^h_c$ in 
Subsections \ref{subsec:onedisk2}
and \ref{subsec:1+} uses the fact that $\Sigma'_{\widehat a}$ 
is a disk.
The target of the map (\ref{Psiborered}) is the 
moduli space of marked disks.
In our situation $\Sigma'_{\widehat a}$ is an 
annulus.
\begin{rem}
The choice of $\gamma^h_c$ which we provide below looks a bit cumbersome.
There is a simpler choice as follows.
We take a metric of constant negative curvature $-1$ on 
$\Sigma'_{\widehat a} \setminus (\vec z^{\prime +} \cup 
\vec w^{\prime +}(h))$ that is complete and has finite volume.
On the neck region corresponding to the singular point $z_c$, 
there is a unique geodesic of minimal length, 
by Margulis' lemma.
We may take it as our $\gamma^h_c$.
\par
In order to use this choice for our purpose we
need to prove a similar estimate as (\ref{uprimehexdec})
for this choice of $\gamma_c^h$.
The authors have no doubt that it is correct.
However, they do not find a reference which they can directly quote 
to prove (\ref{uprimehexdec}) for this choice.
\par
For the choice of $\gamma_c^h$ we made in Subsection \ref{subsec:onedisk2}
and we take below, we can use \cite[Lemma 16.18]{foootech},
from which we can prove (\ref{uprimehexdec}) easily as we did in 
Subsection \ref{subsec:onedisk2}.
This is the reason why we take the current choice of ${\gamma}^h_c$.
\end{rem}
To explain this proof in detail, we begin with a digression.
\par
Let $\frak c(a)_{0}$ be the disk $\Sigma_{a,0}$ together 
with data induced from $\frak c(a)$.
Note
$$
\frak c(a)_{0} \in \mathcal M_{0;\ell_{a,0}} (\beta_{a,0}).
$$
Let ${\bf v}_{a,0}(\vec w^+)$ be $\Sigma_{a,0}$ together 
with $\vec z^+\cap \Sigma_{a,0}$, 
$\vec w^+\cap \Sigma_{a,0}$, and the interior singular points of 
$\Sigma$ on $\Sigma_{a,0}$. Note
$$
{\bf v}_{a,0}(\vec w^+) \in \mathcal M_{0;\ell_{a,0}+m_{a,0}+n_{a,0}}.
$$
We define
$$
{\bf v}_{a,\frak a}(\vec w^+) \in \mathcal M_{\ell_{a,\frak a}+m_{a,\frak a}+n_{a,\frak a}}
$$
in the same way by including marked points to $\Sigma_{a,\frak a}$.
Here the right hand side is a moduli space of 
spheres with $\ell_{a,\frak a}+m_{a,\frak a}+n_{a,\frak a}$ marked points.
\par
We consider the universal family
\begin{equation}\label{disksnivfam}
\pi : \mathfrak M_{0,l} \to \mathcal M_{0,l}
\end{equation}
of disks with $l$ interior marked points.
We take and fix a trivialization of this fiber bundle in a neighborhood 
$\frak V_{{\bf v}_{a,0}(\vec w^+)}$
of  ${\bf v}_{a,0}(\vec w^+)$.
The trivialization of (\ref{disksnivfam}) induces a choice of the 
coordinate at each point of the boundary of the member of
our family as follows.
We take and fix a section $s$ of (\ref{disksnivfam}) so that its 
image is not on the boundary of the disk. 
For each $\frak x \in \mathcal M_{0,l}$ and 
$\frak z_0 \in \partial(\pi^{-1}(\frak x))$ we take a
biholomorphic map
$$
\phi : \pi^{-1}(\frak x) \cong \frak h\cup \{\infty\}
$$
where $\frak h = \{z\in \C \mid {\rm Im}\, z \ge 0\}$, such that
$$
\phi(s(\frak x)) = \sqrt{-1},\qquad \phi(\frak z_0) = 0.
$$
We define our coordinate
\begin{equation}\label{fdrmula4515}
\varphi_{\frak z_0}^{\frak x} :
D^2 \cap \frak h \to \pi^{-1}(\frak x)
\end{equation}
as the restriction of the inverse of $\phi$.
Note this coordinate is determined only by the data
we took for $\frak c(a)$.
\par
We next consider a family of annuli.
We put
$$
\Sigma_r = ([0,2] \times \R)/\Z
$$
where $\Z$ action is defined by $1\cdot [x,y] = [x,y+r]$.
(We put the complex structure such that $x+\sqrt{-1}y$ is a 
complex coordinate.)
\par
For each $r$ and $\frak z_0 = [0,0]$ we take a 
coordinate of its neighborhood by
\begin{equation}\label{coordinatasigma3}
\varphi_{\frak z_0}^r(x+\sqrt{-1}y) = [y, - x],
\end{equation}
where
$
\varphi_{\frak z_0}^r : D^2 \cap \frak h \to \Sigma_r.
$
Let ${\bf v} = (\Sigma_{\bf v},z_1^{\rm sing}) \in \mathcal M_{(1;0)}^{\rm ann}$.
(Here $z^{{\rm sing}}_1$ will become the point 
of $\partial \Sigma_{\frak c(0),0}$ where the root of the 
tree of disk components containing $\frak c(a)$ will be attached.)
We identify $\Sigma_{\bf v} \cong \Sigma_r$ by biholomorphic map
which sends $z^{{\rm sing}}_1$ to $\frak z_0\in \Sigma_r$.
Such $r$ and the biholomorphic map exist uniquely. 
Then (\ref{coordinatasigma3}) determines a coordinate at $z_1^{\rm sing}$.
\par
Now we go back to our situation.
Note that $\frak c(a)$ is attached to one of the two boundary components 
$\partial_1\Sigma_{\frak c(0)}$ or $\partial_2\Sigma_{\frak c(0)}$.
We may assume without loss of generality that
it is attached to $\partial_1\Sigma_{\frak c(0)}$.
\par
Let $\frak V^+_{a,0}$ (resp. $\frak V_{a,\frak a}$) be a
neighborhood of $\frak c(a)^+_0$ (resp. $\frak c(a)_{\frak a}$)
in $\mathcal M_{1,\ell_{a,0}+m_{a,0}+n_{a,0}}$
(resp. $\mathcal M_{\ell_{a,\frak a}+m_{a,\frak a}+n_{a,\frak a}}$).
Here to obtain $\frak c(a)^+_0$ from $\frak c(a)_0$ we add one 
boundary marked point as follows.
There exists a unique connected component of 
$\Sigma_{\frak c} \setminus \Sigma_{\frak c(a)}$ 
which contains an annulus. The boundary marked point we add is 
 the singular point 
on $\Sigma_{\frak c(a)}$ which is contained in the closure of this 
connected component.
\par
Let $\frak V_{(1;0)}^{{\rm main}}$
be a neighborhood of $(\Sigma_{\frak c(0),0},z^{\rm sing}_1)$
in $\mathcal M_{(1,0);0}^{\rm ann}$.
\par
As a part of the data consisting $\frak c(a)$, the coordinate 
of the interior marked points of $\Sigma_{\frak c(a)}$ is 
included. 
We glue an element of $\frak V^+_{a,0}$ with 
an element of $\frak V_{(1;0)}^{{\rm main}}$ at
the unique boundary marked point of the former 
and the boundary marked point $z_1^{\rm sing}$ of the latter.
Note we fixed the coordinates around these two boundary marked points.
\par
Thus we obtain the following map:
\begin{equation}\label{Phistaanus}
\Psi : D^2(\epsilon)^{\mathcal S} \times [0,\epsilon) 
\times
\frak V^+_{a,0} \times \prod_{\frak a}\frak V_{a,\frak a}
\times \frak V_{(1,0)}^{\rm main} 
\to \mathcal M^{\rm ann}_{(0,0);\ell_a+m_a},
\end{equation}
which is injective and is a diffeomorphism to an open set.
Here $\mathcal S$ is the set of interior singular points in 
$\frak c({a})$ and
$\ell_a = \# \vec z^+_a$, $m_a = \# \vec w^+_a$.
The definition of $\Psi$ is similar to (\ref{Psiclosedmap}).
\begin{rem}
In the situation of Remark \ref{pbubblerem}, where $u_{\frak c(a)}$ is constant 
on $\Sigma_{\frak c(a),0}$  and $\Sigma_{\frak c(a),0}$ contains 
only one special point (an interior singular point), 
we replace $D^2(\epsilon)^{\mathcal S}$ by
$D^2(\epsilon)^{\mathcal S - 1} \times [0,\epsilon)$.
\end{rem}
For each $h$ which is close to $h_{\frak c(a)}$, we obtain $\vec w_a^{\prime +}(h) 
\subset \Sigma'_{\widehat a}$ by requiring
$w_{a,b}^{\prime +}(h) \in hN_{w_{a,b}}$.
By definition $(\Sigma'_0,z^0_2,\vec z^{\prime +}_a\cup \vec w_a^{\prime +}(h))$
is contained in the image of (\ref{Phistaanus}).
Using this fact, we can define a loop $\gamma^h_c$ in the same way as 
in Subsection \ref{subsec:onedisk2}.
\par
We use it to define $f_+$ by (\ref{specialpointdistance2}).
Then $h_a$ is the point where $f_+$ attains its minimum.
The rest of the construction is the same as the other cases.
\par\medskip
Thus we have defined linear embeddings $E_{\frak c(0)} \to C^{\infty}(\Sigma',(u')^*TX \otimes 
\Lambda^{0,1})$ and 
$E_{\frak c(a)} \to C^{\infty}(\Sigma',(u')^*TX \otimes 
\Lambda^{0,1})$.
We denote by $E_{\frak c(0)}({\bf y})$, $E_{\frak c(a)}({\bf y})$
their images.
We put
\begin{equation}
E({\bf y})
=
\bigoplus_{\frak c(0) \in \frak C({\bf x}_{0})}
E_{\frak c(0)}({\bf y})
\oplus
\bigoplus_{a} \bigoplus_{\frak c(a) \in \frak C({\bf x}_{a})}
E_{\frak c(a)}({\bf y}).
\end{equation}
Using this choice of $E({\bf y})$, we define a Kuranishi neighborhood 
of ${\bf x}$. Our construction is designed so that 
there is a coordinate change among those Kuranishi neighborhoods 
such that the Kuranishi structure obtained has required properties.
\par
\subsection{Kuranishi structure on the moduli space 
of pseudo-holomorphic annuli II}
\label{subsec:kuraanul2}
\par
We next construct a Kuranishi structure on a neighborhood of
$\frak{forget}^{-1}([\Sigma_2])$.
Let
${\bf x} = (\Sigma,(z_1,z_2),\vec z^+,u) \in 
\frak{forget}^{-1}([\Sigma_2]) \cap \mathcal M_{(1,1);\ell}(\beta).$
We decompose $\Sigma$ into extended disk components as
\begin{equation}
\Sigma = \Sigma_1^{\frak t} \cup \Sigma_2^{\frak t} \cup 
\bigcup_a \Sigma_a.
\end{equation}
Here $\Sigma_1^{\frak t}$ and $\Sigma_2^{\frak t}$ are 
the extended disk components to which we apply $\frak t$-perturbation 
and $\Sigma_a$ are extended 
disk components to which we apply $\frak c$-perturbation.
We determine $\Sigma_1^{\frak t}$, $\Sigma_2^{\frak t}$ as follows.
(We will define $z_i^0 \in \Sigma_i^{\frak t}$ at the same time.)
\par
We take $\Sigma_{\rm main} \subset \Sigma$
that is a union of extended disk components such that 
$\Sigma_{\rm main}$ is not simply connected but 
if we remove any one of the extended disk components 
from $\Sigma_{\rm main}$ then it becomes simply connected.
The unique existence of such $\Sigma_{\rm main}$ is obvious.
\par
If $z_i \in \Sigma_{\rm main}$, then $\Sigma_i^{\frak t}$ 
is the unique extended disk component in
$\Sigma_{\rm main}$ which contains $z_i$, by definition.
We also put $z_i  = z_i^0$.
\par
Suppose $z_i \notin \Sigma_{\rm main}$.
Then there exists a maximal tree of disk components of 
$\Sigma \setminus \Sigma_{\rm main}$ that contains $z_i$.
The extended disk component $\Sigma_i^{\frak t}$ is 
by definition the unique extended disk component of 
$\Sigma_{\rm main}$ where this tree of  disk components 
is rooted. The point $z_i^0$ is by definition the root of  
this tree of  disk components.
\par
For each extended disk component $\Sigma_{\alpha}$ 
in $\Sigma_{\rm main}$ there are two boundary singular 
points $z_{\alpha}^1, z_{\alpha}^2$ where $\Sigma_{\alpha}$ 
are attached to other extended disk components of $\Sigma_{\rm main}$.
We use orientation of $\Sigma_{\rm main}$  so that
the enumeration $z_{\alpha}^1, z_{\alpha}^2$ respects the cyclic 
order.
(See Figure 4.5.1.)
\par
\epsfbox{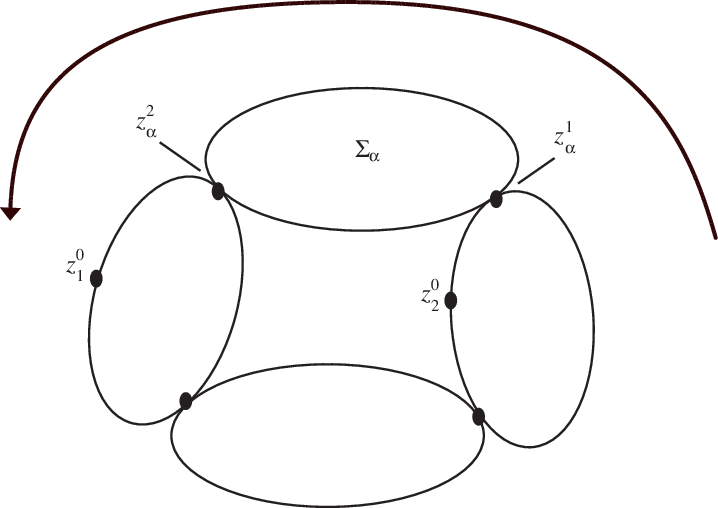}
\par
\centerline{\bf Figure 4.5.1}
\par
\par
We put
\begin{equation}
\aligned
&(\Sigma^{\frak t}_1,(z^0_1,z_{\frak t,1}^1,z_{\frak t,1}^2),\vec z^+\cap \Sigma_1^{\frak t},u)
= {\bf x}^{\frak t}_1 \in \mathcal M^{\rm main}_{3,\ell_1}(\beta_1) = \mathcal M_{(0,0,0),\ell_1}(\beta_1),
\\
&(\Sigma^{\frak t}_2,(z^0_2,z_{\frak t,2}^1,z_{\frak t,2}^2),\vec z^+\cap \Sigma_2^{\frak t},u)
= {\bf x}^{\frak t}_2 \in \mathcal M^{\rm main}_{3,\ell_2}(\beta_2)= \mathcal M_{(0,0,0),\ell_2}(\beta_2),
\\
&(\Sigma_a,\vec z^+\cap \Sigma_a,u)
= {\bf x}_a \in \mathcal M_{\ell_a}(\beta_a).
\endaligned
\end{equation}
Here $z_{\frak t,1}^1$ is the point $z_{\alpha}^1$ where 
$\Sigma_{\alpha} = \Sigma_1^{\frak t}$.
The boundary marked points
$z_{\frak t,1}^2$, $z_{\frak t,2}^1$, $z_{\frak t,2}^2$ are defined 
in the same way.
\par
Let $\frak c_i^{\frak t} \in \frak P_{3,\ell_i}(\beta_i)$ be an element 
such that $[{\bf x}^{\frak t}_i] \in \frak U(\frak c_i^{\frak t})$.
Here $\frak P_{3,\ell}(\beta)$ is a finite subset of 
$\mathcal M^{\rm main}_{3,\ell}(\beta)/(T^n \times \frak S_{\ell})$, which 
we defined while we constructed $\frak t$-Kuranishi structure on 
$\mathcal M^{\rm main}_{3,\ell}(\beta)$.
\par
Let $\frak c_a \in \frak P_{\ell_a}(\beta_a)$ be an element 
such that ${\bf x}_a \in \frak U(\frak c_a)$.
Here $\frak P_{\ell}(\beta)$ is a finite subset of 
$\mathcal M_{\ell}(\beta)/(T^n \times \frak S_{\ell})$, which 
we defined while we constructed the $\frak c$-Kuranishi structure on 
$\mathcal M_{\ell}(\beta)$.
\par
We took additional marked points $\vec w^{\frak t +}_i$  of 
$\frak c_i^{\frak t}$, codimension $2$ submanifolds 
$N_{i,b}^{\frak t}$ of $X$ and obstruction spaces 
$E_{\frak c^{\frak t}_i}$. 
We also took additional marked points $\vec w^+_a$ of 
$\frak c_a$ and codimension $2$ submanifolds 
$N_{a,b}$  if $\frak c_a$ is either of Type A or B.
In case  $\frak c_a$ is of Type C we 
temporary take them.
In case  $\frak c_a$ is  of Type C
we took circles ${\bf S}_a$. 
We also took the obstruction space 
$E_{\frak c_a}$ for each $a$.
\par
By assumption, $\mu^{{\bf x}^{\frak t}_i}{\bf x}^{\frak t}_i$ is 
close to 
$(h_i^{0 \frak t}\frak c_i^{\frak t},\vec w^{\frak t +}_i)$
with respect to $\vec w^{{\bf x} \frak t +}_i$
and
$\mu^{{\bf x}_a}{\bf x}_a$ is 
close to 
$(h_a^0\frak c_a,\vec w^{+}_a)$
with respect to $\vec w^{{\bf x} +}_a$.
(Here $\mu^{{\bf x}^{\frak t}_i} \in \frak S_{\ell_i}$,
$\mu^{{\bf x}_a} \in \frak S_{\ell_a}$, 
$h_i^{0 \frak t}, h_a^0 \in T^n$.)
\par
We fix $\frak c^{\frak t}_i$ and $\frak c_a$ for a while.
\par
By gluing $\frak c^{\frak t}_i, \frak c_a$, we obtain 
$\frak c = (\Sigma_{\frak c},(z^{\frak c}_1,z^{\frak c}_2),
\vec z^+_{\frak c},u_{\frak c})$.
(There is a similar ambiguity as mentioned in Remark \ref{smallambrem}.
We will remove it during the construction below.)
On $\Sigma_{\frak c}$, we have the set of additional marked points $\vec w_{\frak c}^{+}$
that is a union of $\vec w^{\frak t +}_i$ and $\vec w^{+}_a$.
Let  $\vec w^{{\bf x} +}$ be the union of $\vec w^{{\bf x} \frak t +}_i$ and
$\vec w^{{\bf x} +}_a$.
\par\medskip
Let ${\bf y} =(\Sigma',(z'_1,z'_2),\vec z^{\prime +},u')$
be an object. Suppose $\mu_1{\bf y}$ is close to 
$(h{\bf x},\vec w^{{\bf x} +})$ with respect to 
$\vec w^{\prime +}$.
We have a map
\begin{equation}\label{iembsubsec545}
i_{{\bf y}(\vec w^{\prime +})} :
\Sigma_{\frak c}\setminus \mathcal U(S(\Sigma_{\frak c})) \to \Sigma'.
\end{equation}
Below we define the way how we send $E_{\frak c^{\frak t}_i}$,
$E_{\frak c_a}$ to $\Sigma'$.
For simplicity of notation we assume $\mu_1 = 1$.
\par\smallskip
We first describe the way how we send $E_{\frak c^{\frak t}_1}$.
The case of $E_{\frak c^{\frak t}_2}$ is entirely similar.
\par
Note 
${\bf x}^{\frak t}_1$ is 
close to 
$(h_1^{0 \frak t}\frak c_1^{\frak t},\vec w^{\frak t +}_1)$
with respect to $\vec w^{{\bf x} \frak t +}_1$.
However, the choice of $h_1^{0 \frak t}$ is not canonical.
We first replace it by $h_1^{\frak t}$ which is determined by  
${\bf y}$ and $\frak c^{\frak t}_{1}$ only without ambiguity.
\par
We start with a digression.
Since 
$
(\Sigma^{\frak t}_{1,0},(z^0_1,z_{\frak t,1}^1,z_{\frak t,1}^2)) 
$
is a disk with three boundary marked points, 
there exists uniquely a biholomorphic map
\begin{equation}\label{disk3tenstandar}
\varphi : \Sigma^{\frak t}_{1,0} 
\to \{z \in \C \mid {\rm Re}\, z \in [0,1]\} \cup \{ \pm \sqrt{-1}\infty \},
\end{equation}
such that $\varphi(z^0_1) = 1$, $\varphi(z_{\frak t,1}^1) = + \sqrt{-1}\infty$,
$\varphi(z_{\frak t,2}^1) = - \sqrt{-1}\infty$.
Let $\Sigma'_0$ be the irreducible component of $\Sigma'$ that is an 
annulus. There exists a unique $T>0$ such that $\Sigma'$ is 
biholomorphic to
$$
\{ z \in \C \mid {\rm Re}\, z \in [0,1],\,\, {\rm Im}\,\, z \in [-T,T] \}/\sim
$$
where $-T\sqrt{-1} \sim T\sqrt{-1}$.
We define a biholomorphic map 
$$
\Phi : \Sigma' \to \{ z \in \C \mid {\rm Re}\, z \in [0,1],\,\, {\rm Im}\,\, z \in [-T,T] \}/\sim
$$
by requiring $\Phi(z^0_1) = 1$.
Here $z^0_1 \in \partial \Sigma'_0$ is defined as follows.
If $z_1 \in \Sigma'_0$, then $z_1 = z^0_1$. If not, there exists a tree of 
disk components of $\Sigma'$ which contains $z_1$ and which is rooted on $\Sigma'_0$.
Then $z^0_1$ is by definition the root of this tree of disk components.
\par\smallskip
We now describe the way how to define $h^{\frak t}_1$.
\par
In the case when $u_{\frak c^{\frak t}_1}$ is nonconstant, 
we can  find such $h_1^{\frak t}$ in the same way as in Subsection \ref{subsec:1+}.
Note that in case when $\frak c^{\frak t}_1$ is of Type B, we consider the map
\begin{equation}\label{Phiinsec456}
\Psi : D^2(\epsilon)^{\mathcal S} \times (T_0,\infty] \times \frak V_{\frak c^{\frak t}_{1,0}} \times \prod_{\frak a}\frak V_{\frak c^{\frak t}_{1,\frak a}}
\to \mathcal M_{(0,1);\ell_{1}+m_1}^{\rm ann}.
\end{equation}
Here $\frak c^{\frak t}_{1,0}$ (resp. $\frak c^{\frak t}_{1,\frak a}$) is obtained by restricting 
the data of $\frak c^{\frak t}_1$ to the disk component (resp. a sphere component) and regarding 
interior singular points as marked points.
$\frak V_{\frak c^{\frak t}_{1,0}}$ (resp. $\frak V_{\frak c^{\frak t}_{1,\frak a}}$)
is an neighborhood of $\frak c^{\frak t}_{1,0}$ (resp. $\frak c^{\frak t}_{1,\frak a}$) in the
moduli space of marked disks (resp. spheres). 
$m_1 = \# \vec w^{\frak t +}_{1}$. The factor $D^2(\epsilon)^{\mathcal S} $ parametrizes the way how to 
smoothen the interior singular points. The factor $(T_0,\infty]$ parametrizes the way 
how to obtain $
\{ z \in \C \mid {\rm Re}\, z \in [0,1],\,\, {\rm Im}\,\, z \in [-T,T] \}/\sim
$
from $\{z \in \C \mid {\rm Re}\, z \in [0,1]\} \cup \{ \pm \sqrt{-1}\infty \}$.
The map $\Psi$ is  defined in the same way as (\ref{Psiclosedmap}).
Using $\Psi$ we define $\gamma^h_c$ in the same way as in Subsection \ref{subsec:onedisk2}. We use it to define a function $f^+$ and
then an element $h_t^{\frak t}$.
Using $h_t^{\frak t}$ we send $E_{\frak c^{\frak t}_1}$ to $\Sigma'$ in the same way as 
(\ref{Pprelimi}).
\par
Suppose $u_{\frak c^{\frak t}_i}$ is constant.
While we defined $\frak t$-Kuranishi structure, we chose $w_0 
\in {\rm Int}\,\,\Sigma^{\frak t}_{i,0}$ where $\Sigma^{\frak t}_{i,0}$ is the 
disk component.
We use $\varphi$ in (\ref{disk3tenstandar}) to obtain
$\varphi(w_0)$.
Assuming ${\bf y}$ is sufficiently close to ${\bf x}$, we may assume $T$ large so that
$\vert{\rm Im}\,\varphi(w_0)\vert < T$.
Then we may regard $\Phi^{-1}(\varphi(w_0)) \in \Sigma'$.
We put
$$
f(h') = {\rm dist}(u'(\Phi^{-1}(\varphi(w_0))), h'u(w_0))^2.
$$
This is a strictly convex function. We define $h_1^{\frak t}$ so that 
$f$ attains its minimum at $h_1^{\frak t}$.
\par
We may also assume that 
$$
\sup \{ \vert{\rm Im}\,\varphi(z)\vert \mid z \in {\rm Support}\, E_{\frak c^{\frak t}_{1,0}}\} <T.
$$
Therefore $\Phi^{-1} \circ \varphi$ defines an embedding $: {\rm Support}\, E_{\frak c^{\frak t}_{1,0}} 
\to \Sigma'$.
We use it and the parallel transport along the minimal geodesic in $X$ to send  $E_{\frak c^{\frak t}_{1,0}}$ 
to $\Sigma'$.
\par\smallskip
We next describe the way how we send $E_{\frak c_a}$ to $\Sigma'$.
Let $\Sigma'_{\widehat a}$ be a component of $\Sigma'$ such that
$$
i_{\bf y(\vec w^{\prime +})}(\Sigma_{\frak c_a} \setminus \mathcal U(S(\Sigma_{\frak c}))) 
\subset \Sigma'_{\widehat a}
$$
where $i_{\bf y(\vec w^{\prime +})}$ is as in  (\ref{iembsubsec545}).
\par
In case $\Sigma'_{\widehat a}$ is an extended disk component we can 
send $E_{\frak c_a}$ to $\Sigma'$ in the same way as in Subsection \ref{subsec:1+}.
\par
In case $\frak c_a$ is of Type A or C, the way how we send $E_{\frak c_a}$ to $\Sigma'$ is again the same as in Subsection \ref{subsec:1+}.
\par
We finally study the case when $\frak c_a$ is of Type B.
We consider $(\Sigma_{\frak c}, \vec z^{+}_{a} \cup \vec w^{+}_{a})$.
(Namely we consider only the marked points on $\Sigma_{\frak c_a}$ and forget all the 
other marked points.)
We shrink the disk or sphere components of $\Sigma_{\frak c}$ which becomes unstable, inductively,
and denote by $\overline\Sigma$ the curve we obtain in this way.
\begin{lem}\label{selfintlemma}
$\overline\Sigma$ satisfies one of the following two alternatives.
\begin{enumerate}
\item
$\overline\Sigma$ has only one extended disk component which intersects with itself 
at one point.
\item
$\overline\Sigma$ has two extended disk components.
The first component intersects with itself. The second component intersects with
the first component at one point. 
\end{enumerate}
\end{lem}
\begin{proof}
If $\Sigma_{a} \subset \Sigma$ is contained in   $\Sigma_{\rm main}$, then 1 occurs. 
Otherwise 2 occurs.
(See Figure 4.5.2.)
\end{proof}

\par\medskip
\hskip-0.3cm\epsfbox{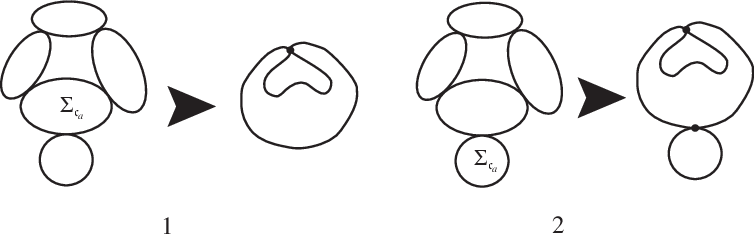}
\par
\centerline{\bf Figure 4.5.2}
\par

We note that all the marked points are on the second component in case of Lemma \ref{selfintlemma}.2.
\par
\noindent
{\bf Case 1}:
Suppose $\overline\Sigma$ satisfies Lemma \ref{selfintlemma}.1.
We decompose $\Sigma_{\frak c_a}$ into irreducible components.
Its disk component (resp. sphere components) gives $\frak c_{a,0}$
(resp. $\frak c_{a,\frak a}$).
Let $\ell_{a,0}, m_{a,0}, n_{a,0}$ (resp. $\ell_{a,\frak a}, m_{a,\frak a}, n_{a,\frak a}$)
be the number of elements of $\vec z^+_a$, $\vec w^+_a$ and the singular points 
on $\Sigma_{\frak c_{a,0}}$ (resp. on $\Sigma_{\frak c_{a,\frak a}}$).
We consider the universal family 
\begin{equation}\label{diskunivfam3plus}
\pi: \mathfrak M_{0;\ell_{a,0}+m_{a,0}+n_{a,0}} \to \mathcal M_{0;\ell_{a,0}+m_{a,0}+n_{a,0}}
\end{equation}
of marked disks. Let $\frak V_{\frak c_{a,0}}$ be a neighborhood of $\frak c_{a,0}$ 
(with the map $u_{\frak c_{a,0}}$ being forgotten) in  $\mathcal M_{0;\ell_{a,0}+m_{a,0}+n_{a,0}}$.
As a part of the data we fix for $\frak c_a$ we took a trivialization 
(as a smooth manifold) of (\ref{diskunivfam3plus}) 
on $\frak V_{\frak c_{a,0}}$.
We use it to define a coordinate at each point of 
$\partial\pi^{-1}(\frak V_{\frak c_{a,0}})$ as in (\ref{fdrmula4515}).
\par
We use this coordinate (and the coordinate at the interior marked points) to define a map
\begin{equation}\label{4525casemform}
\Psi : D^2(\epsilon)^{\mathcal S} \times [0,\epsilon) \times C_1 \times C_2 \times \frak V_{\frak c_{a,0}} \times \prod_{\frak a}\frak V_{\frak c_{a,\frak a}}
\to \mathcal M_{(0,1);\ell_{a}+m_a}^{\rm ann}.
\end{equation}
Here 
$\frak V_{\frak c_{a,\frak a}}$ is a neighborhood of 
$\frak c^{\frak t}_{a,\frak a}$ (with map being forgotten) in the
moduli space of marked  spheres, and 
$m_a = \# \vec w_{a}$. The factor $D^2(\epsilon)^{\mathcal S} $ parametrizes the way how to 
smoothen the interior singular points. The factor $[0,\epsilon)$ parametrizes the way 
how to smoothen the self intersection of the disk component.
The parameter spaces $C_1$, $C_2$ are diffeomorphic to the arcs and parametrize the 
position of the self intersection of the disk component.
The map $\Psi$ is  defined in the same way as (\ref{Psiclosedmap}) and (\ref{Phiinsec456}).
\begin{rem}
In the situation of Remark \ref{pbubblerem}, where $u_{\frak c(a)}$ is contant 
on $\Sigma_{\frak c(a),0}$  and $\Sigma_{\frak c(a),0}$ contains 
only one special point (an interior singular point), 
we replace $D^2(\epsilon)^{\mathcal S}$ by
$D^2(\epsilon)^{\mathcal S - 1} \times [0,\epsilon)$.
\end{rem}
It is easy to see that $(\Sigma'_0,\vec z^{\prime +}_{a} \cup \vec w^{\prime +}_a)$ is in the image of
(\ref{4525casemform}). We can then define $\gamma_c^h$ in the same way as 
in Subsection \ref{subsec:onedisk2}. We use it to define a function $f_+$ and
then an element $h_a$.
We use $h_a$ to send $E_{\frak c_a}$ to $\Sigma'$ in the same way as 
(\ref{Pprelimi}).
\par\smallskip
\noindent
{\bf Case 2}:
The case when $\overline\Sigma$ satisfies Lemma \ref{selfintlemma}.2.
We identify the first comonent there with 
\begin{equation}\label{domain111}
\{ z \in \C \mid {\rm Re}\, z \in [0,2]\} \cup \{\infty\}.
\end{equation}
Here $\infty$ is the boundary point (regarded as a point in the Riemann sphere) 
where the disk intersects with itself at the boundary.
We define its coordinate at $\frak z_0 = (0,0)$ by (\ref{coordinatasigma3}).
We then obtain a map
\begin{equation}\label{4525casemform222}
\Psi : D^2(\epsilon)^{\mathcal S} \times [0,\epsilon) \times (T_0,\infty] \times C \times \frak V_{\frak c_{a,0}} \times \prod_{\frak a}\frak V_{\frak c_{a,\frak a}}
\to \mathcal M_{(0,1);\ell_{a}+m_a}^{\rm ann}.
\end{equation}
Here we use the factor $T \in (T_0,\infty]$ to deform the singularity of (\ref{domain111}) to
\begin{equation}\label{domain112}
\{ z \in \C \mid {\rm Re}\, z \in [0,2], \,\,  {\rm Im}\, z \in [-T,T]\} /\sim,
\end{equation}
where $z \sim 2T\sqrt{-1} + z$.
$[0,\epsilon)$ is the parameter to deform the singularity which is the intersection point of the first and the second 
components. $C$ is an open set of an arc which parameterizes the position of the boundary singular point of the second component.
\par
We use (\ref{4525casemform222}) to define the loop $\gamma_c^h$ in the same way as Case 1.
The rest of the proof is the same as Case 1.
\par\smallskip
We have thus described the way how we send obstruction spaces to $\Sigma'$.
We can then define $E({\bf y})$ by
\begin{equation}
E({\bf y})
=
\bigoplus_{i=1,2}\bigoplus_{\frak c^{\frak t}_i \in \frak C({\bf x}^{\frak t}_i)}
E_{\frak c^{\frak t}_i}({\bf y})
\oplus
\bigoplus_{a} \bigoplus_{\frak c_a \in \frak C({\bf x}_{a})}
E_{\frak c_a}({\bf y}).
\end{equation}
We can use it to define a required Kuranishi structure on a neighborhood of 
$\frak{forget}^{-1}([\Sigma_2])$.
\par
\subsection{Kuranishi structure on the moduli space 
of pseudo-holomorphic annuli III}
\label{subsec:kuraanul3}
\par
We next construct a Kuranishi structure on a neighborhood of
$\frak{forget}^{-1}([\Sigma_1])$.
Namely we prove Lemma \ref{perturbnbdsigma1}. 
The construction is similar to one given in Subsection \ref{subsec:lemma426}.
Let
${\bf x} = (\Sigma,(z_1,z_2),\vec z^+.u) \in 
\frak{forget}^{-1}([\Sigma_1]) \cap \mathcal M_{(1,1);\ell}(\beta).$
\par
Since $\frak{forget}({\bf x}) = [\Sigma_1]$, there exists a canonical map
$\pi : \Sigma \to \Sigma_1$. 
Recall that $\Sigma_1$ is a union of two disks 
glued at  interior points. So we have exactly two 
disk components $\Sigma^{\frak m}_1$, $\Sigma^{\frak m}_2$
of $\Sigma$ on which $\pi$ is nonconstant.
We enumerate them so that $\partial\Sigma^{\frak m}_i \subset 
\partial_i\Sigma$.
\par 
We take the (unique) minimal union of sphere components 
$\Sigma({\rm II},\frak a)$, $\frak a=1,\dots,N$, such that
\begin{equation}\label{mainin546}
\Sigma_{\rm main}
= \Sigma^{\frak m}_1 \cup \Sigma^{\frak m}_2
\cup \bigcup_{\frak a=1,\dots,N}\Sigma({\rm II},\frak a)
\end{equation}
is connected. We take $z^0_i \in {\rm Int}\, \Sigma^{\frak m}_i$
such that
$$
\{z^0_i\} = \overline{(\Sigma_{\rm main} \setminus \Sigma^{\frak m}_i)}
\cap \Sigma^{\frak m}_i.
$$
The disk components which are not $\Sigma^{\frak m}_1$, $\Sigma^{\frak m}_2$
are denoted by $\Sigma^{\frak c}_{i,\frak a}$,
$i=1,2$, $\frak a=1,\dots,N_i$, where $\partial \Sigma^{\frak c}_{i,\frak a}
\subset \partial_i\Sigma$.
\par
A sphere component other than $\Sigma({\rm II},\frak a)$ 
is contained in a tree of sphere components which satisfy 
one of the following:
\begin{conds}\label{spherecomp2jyuusannkaku}
\begin{enumerate}
\item
The tree of sphere components is rooted on $\Sigma^{\frak m}_1$.
\item
The tree of sphere components is rooted on $\Sigma^{\frak m}_2$.
\item
The tree of sphere components is rooted on $\Sigma({\rm II},\frak a)$.
\item
The tree of sphere components is rooted on $\Sigma^{\frak c}_{1,\frak a}$.
\item
The tree of sphere components is rooted on $\Sigma^{\frak c}_{2,\frak a}$.
\end{enumerate}
\end{conds}
We identify $(\Sigma^{\frak m}_i,z^0_i) \cong (D^2,0)$
by a biholomorphic map and choose $\epsilon > 0$ sufficiently small.
\begin{defn}
\begin{enumerate}
\item
A sphere component is of Type I if it is contained in the 
tree of sphere components satisfying Condition \ref{spherecomp2jyuusannkaku}.4 or 5.
\item
A sphere component is of Type I if it is contained in the 
tree of sphere components satisfying 
Condition \ref{spherecomp2jyuusannkaku}.1 or 2 and 
${\rm dist}(w,\partial \Sigma^{\frak m}_i) < \epsilon$,
where $w$ is the root of this tree of sphere components.
\item
$\Sigma({\rm II},a)$ is a sphere component of Type II.
\item
A sphere component is of Type II if it is contained in the 
tree of sphere components satisfying Condition \ref{spherecomp2jyuusannkaku}.3.
\item
A sphere component is of Type II if it is contained in the 
tree of sphere components satisfying 
Condition \ref{spherecomp2jyuusannkaku}.1 or 2 and 
${\rm dist}(w,z^0_i) < \epsilon$,
where $w$ is the root of this tree of sphere components.
\item
All the other sphere components are of Type III.
\end{enumerate}
\end{defn}
\begin{defn}\label{defcompwise2}
A system of Kuranishi structures of $
\mathcal M_{(k_1+1,k_2+1);\ell}^{{\rm ann};\text{\rm main}}(\beta)$
in a neighborhood of its intersection with  $\frak{forget}^{-1}([\Sigma_1])$
is said to be {\it partially component-wise}\index{Kuranishi structure!partially component-wise Kuranishi structure} on 
$\frak{forget}^{-1}([\Sigma_1])$ if the following holds:
\par
We consider the stratification of $
\mathcal M_{(k_1+1,k_2+1);\ell}^{{\rm ann};\text{\rm main}}(\beta)
\cap \frak{forget}^{-1}([\Sigma_1])$
such that each stratum is a fiber product of
\begin{equation}\label{diskgeneral2}
\mathcal M_{k_{i,j}+1;\ell_{i,j}}^{\text{\rm main}}(\beta_{(i,j)}),
\qquad i=1,2, \,\, j=1,\dots,m_i,
\end{equation}
or
\begin{equation}\label{spherreq2}
\mathcal M_{\ell'_j}(\alpha_j),
\qquad j=1,\dots,m'
\end{equation}
or
\begin{equation}\label{disk0comp2}
\mathcal M_{k'_1+1;\ell''_1}^{\text{\rm main}}(\beta'_1),
\end{equation}
or
\begin{equation}\label{disk0comp3}
\mathcal M_{k'_2+1;\ell''_2}^{\text{\rm main}}(\beta'_2),
\end{equation}
where we have $\sum_j k_{i,j} = k_i+m_i$, $\sum_{i,j}\beta_{(i,j)} + \sum_j \alpha_j + \beta'_1 + \beta'_2  = \beta$, $\sum_{i,j} \ell_{i,j}
+ \sum_j \ell'_j + \ell''_1 + \ell''_2= \ell+m'+1$. 
\par
We require that there is exactly one component each as in (\ref{disk0comp2}), 
(\ref{disk0comp3}). They play the role of $\Sigma^{\frak m}_1$ 
and $\Sigma^{\frak m}_2$ in (\ref{mainin546}), respectively.
\par
For the factor (\ref{spherreq2}), we require the corresponding sphere components
are of type II.
(In other words, we include a tree of sphere components of Type I or III to the
disk component on which it is rooted.)
\par
Then the Kuranishi structure on 
$\frak{forget}^{-1}([\Sigma_1])$ is compatible with this fiber product description of the
stratum of this stratification. Here we use the $\frak c$-Kuranishi structure for the factor (\ref{diskgeneral2}).
For the factor (\ref{spherreq2}),  
we use the Kuranishi structure which was used in Subsection \ref{subsec:lemma426} for the factor (\ref{spherreq}).
For the factors (\ref{disk0comp2}), (\ref{disk0comp3}), 
we use a Kuranishi structure similar to one we used in Subsection \ref{subsec:lemma426} for the factor (\ref{disk0comp}).
(Namely the case of Condition \ref{decompconds}.3). We however require the following in this case:
The exponential map at $z^0_i \in {\rm Int}\, \Sigma^{\frak m}_i$ gives a {\it surjective} map:
\begin{equation}\label{surjeval4523}
{\rm Eval}: (D_{u_{\frak c^{\frak t}_i}} \overline\partial)^{-1}(E_{\frak c^{\frak t}_i}) \to
T_{u_{\frak c^{\frak t}_i}(z^0_i)}X.
\end{equation}
\end{defn}
\begin{rem}
\begin{enumerate}
\item
Since we assume (\ref{surjeval4523}) to be surjective, we need to put nontrivial
obstruction space on the disk component in case $u_{\frak c^{\frak t}_i}$ is constant on the disk component, also.
\item
We use surjectivity of  (\ref{surjeval4523}) to show that the fiber product which corresponds to the 
gluing at $D^{\frak m}_1 \cap D^{\frak m}_2$ (and a similar gluing in case there is a line of 
Type II sphere bubbles between $D^{\frak m}_1$ and $D^{\frak m}_2$) is transversal.
\end{enumerate}
\end{rem}
In the rest of this subsection, we will construct 
such partially component-wise Kuranishi structure
which also has the properties formulated in Lemma \ref{bdryMannu}.
Since our Kuranishi strucrture is compatible 
with forgetful map
$$
\mathcal M_{(k_1+1,k_2+1);\ell}^{{\rm ann};\text{\rm main}}(\beta)
\to
\mathcal M_{(1,1);\ell}^{{\rm ann};\text{\rm main}}(\beta)
$$
it suffices to consider the case $k_1=k_2=0$.
The obstruction spaces we put to each of the factors 
are specified in Definition \ref{defcompwise2}.
So it suffices to describe the way how we send them to nearby objects.
We assume:
\begin{assum}\label{assump4524}
All the Type I sphere components are on the disk bubble.
\end{assum}
This is the same as Assumption \ref{asssuYpteIIII}.
By the same reason as in Subsection \ref{subsec:lemma426}
we only need to consider the case 
when the center of our Kuranishi neighborhood ${\bf x}$ 
satisfies Assumption \ref{assump4524}.
\par
Let
${\bf x} = (\Sigma,(z_1,z_2),\vec z^+.u) \in 
\frak{forget}^{-1}([\Sigma_1]) \cap \mathcal M_{(1,1);\ell}(\beta).$
We decompose 
$$
\Sigma = \bigcup_{\alpha\in A} \Sigma_{a}
$$
such that each $\Sigma_a$ satisfies one of the following:
\begin{conds}\label{conds45251}
\begin{enumerate}
\item
$\Sigma_{a}$ is one of $\Sigma^{\frak c}_{\frak a}$ together with 
 trees of sphere components of type I.
\item
$\Sigma_{a}$ is a sphere component of Type II.
\item
$\Sigma_a$ is one of $\Sigma_{i}^{\frak w}$ together 
with trees of sphere components of type III.
(In other words it is  one of (\ref{disk0comp2}) or (\ref{disk0comp3})
where $k'_i = 0$.)
\end{enumerate}
\end{conds} 
Thanks to Assumption \ref{assump4524}, such a decomposition exists and is unique.
\par
We denote by ${\bf x}_{a}$ the component $\Sigma_{a}$ together with data induced from ${\bf x}$. 
We take $\frak c(a) = (\Sigma_{\frak c(a)},\vec z^+_{\frak c(a)},u_{\frak c(a)})$ such that:
\begin{enumerate}
\item[(1)]
The case of Condition \ref{conds45251}.1. Here
$\frak c(a) \in \frak P_{\ell_{a}}(\beta_{a})$,
${\bf x}_{a} \in \frak U(\frak c(a))$ and 
${\bf x}_{a} \in \mathcal M_{0;\ell_{a}}(\beta_{a})$.
\item[(2)]
The case of Condition \ref{conds45251}.2.
Here $\frak c(a) \in \frak P_{\ell_{a}}^{\rm sph}(\alpha_{a})$,
${\bf x}_{a} \in \frak U(\frak c(a))$ and 
${\bf x}_{a} \in \mathcal M_{\ell_{a}}(\alpha_{a})$.
\item[(3)]
The case of Condition \ref{conds45251}.3.
Here $\frak c(a) \in \frak P_{\ell_{a}}(\beta_{a})$,
${\bf x}_{a} \in \frak U(\frak c(a))$ and 
${\bf x}_{a} \in \mathcal M_{0;\ell_{a}}(\beta_{a})$.
\end{enumerate}
\par
Let ${\bf y} = (\Sigma',\vec z^{\prime +},u')$ be close to ${\bf x}$.
We will describe the way how we send $E_{\frak c(a)}$ (which is the obstruction space that we defined 
already as is mentioned in Definition \ref{defcompwise2}) to $\Sigma'$.
\par
We note that in Case (1) we need to determine $h_a$ also.
In the other two cases, we do not need to determine $h_a$.
(Note that we do not claim our Kuranishi structure to be $T^n$ equivariant.)
\par
In Case (2) the way how we send $E_{\frak c(a)}$ to $\Sigma'$ 
is the same as (\ref{Pprelimi}).
(Note in this case we assumed that $u_{\frak c(a)}$ is an immersion 
on the support of $E_{\frak c(a)}$.)
\par
Case (3):
We decompose
$$
\Sigma_{\frak c(a)}
= 
\Sigma_{\frak c(a),0} \cup \bigcup_{\frak a}\Sigma_{\frak c(a),\frak a}
$$
where $\Sigma_{\frak c(a),0}$ is a disk and $\Sigma_{\frak c(a),\frak a}$ 
are spheres. 
(Note $\Sigma_{\frak c(a),0}$ is either $\Sigma^{\frak m}_1$ or $\Sigma^{\frak m}_2$.)
The obstruction spaces are decomposed as 
$$
E_{\frak c(a)} = E_{\frak c(a),0} \oplus \bigoplus_{\frak a}E_{\frak c(a),\frak a}.
$$
Note we assumed that $u_{\frak c(a)}$ is an immersion on the support of 
$E_{\frak c(a),\frak a}$. So the way how we send $E_{\frak c(a),\frak a}$ to $\Sigma'$ 
is the same as  (\ref{Pprelimi}).
\par
In case $u_{\frak c(a)}$ is nonconstant on $\Sigma_{\frak c(a),0}$, we can send
$E_{\frak c(a),0}$ to $\Sigma'$ 
in the same way as (\ref{Pprelimi}).
\par
Let us discuss the case when $u_{\frak c(a)}$ is constant on $\Sigma_{\frak c(a),0}$.
We consider marked points $z^+_{\frak c(a)}$ on each component 
$\Sigma_{\frak c(a),0}$, $\Sigma_{\frak c(a),\frak a}$.
We also took additional marked points $\vec w^{+}_{\frak c(a)}$. 
(When $u\vert_{\Sigma_{\frak c(a),a}}$ is constant, we do not take additional marked points on $\Sigma_{\frak c(a),a}$.)
Then $\Sigma_{\frak c(a),0}$ and  $\Sigma_{\frak c(a),\frak a}$
together with them and singular points on them define elements  
${\bf v}(c(a),0)(\vec w^+) \in \mathcal M_{0;\ell_{\frak c(a),0}+m_{\frak c(a),0}+n_{\frak c(a),0}}$,  ${\bf v}(c(a),\frak a)(\vec w^+) \in  \mathcal M_{\ell_{\frak c(a),\frak a}+m_{\frak c(a),\frak a}+n_{\frak c(a),\frak a}}$, respectively.
Let $\frak V_{\frak c(a),0}$, $\frak V_{\frak c(a),\frak a}$ be their neighborhoods.
\par
We define a map
\begin{equation}\label{mapPhi4532}
\Psi : D^2(\epsilon)^{\mathcal S} \times [0,\epsilon) \times 
\frak V_{\frak c(a),0} \times \prod_{\frak a}\frak V_{\frak c(a),\frak a}
\to \mathcal M_{(0,0);\ell_a+m_{a}}^{\rm ann}.
\end{equation}
Here ${\mathcal S}$ is the set of interior singular points on $\frak c(a)$ and 
$D^2(\epsilon)^{\mathcal S}$ is used as the gluing parameter of the 
interior singular points. The parameter $[0,\epsilon)$ is used as follows.
Suppose $\Sigma_a = \Sigma_{1}^{\frak w}$.
Then $[0,\epsilon)$ is the gluing parameter to glue 
$\Sigma_{1}^{\frak w}$ with $\Sigma_{2}^{\frak w}$.
Namely, for $\delta \in [0,\epsilon)$ we identify 
$x \in \Sigma_{1}^{\frak w}$ with $y \in \Sigma_{2}^{\frak w}$
if $xy = \delta$.
(Note we identify $(\Sigma^{\frak m}_i,z^0_i) \cong (D^2,0)$.)
(This gluing parameter is $[0,\epsilon)$ and is not $D^2(\epsilon)$,
since we forget all the marked points on $\Sigma_2^{\frak m}$.)
The other parts of the construction of (\ref{mapPhi4532})
is the same as (\ref{Psiclosedmap}) and (\ref{Phiinsec456}).
\par
When we took $\vec w^+_{\frak c(a)}$, we took 
codimension $2$ submanifolds $N_{w_{\frak c(a),b}}$ at the same time.
We require $u'(w^{\prime +}_{\frak c(a),b}) \in N_{w_{\frak c(a),b}}$
to define $\vec w^{\prime +}_{\frak c(a)} \subset \Sigma'$.
\par
Among the interior marked points $\vec z^{\prime +}$ we take those corresponding to 
$\vec z^{+}_{\frak c(a)}$ and denote them as $\vec z^{\prime +}_{\frak c(a)}$.
\par
Let $\Sigma'_0$ be the irreducible component of $\Sigma'$ that is an annulus.
(Here we consider the case when $\frak{forget}({\bf y}) \ne [\Sigma_1]$.
The case $\frak{forget}({\bf y}) = [\Sigma_1]$ is similar and simpler and so is omitted.) 
Note $\vec w^{\prime +}_{\frak c(a)}, \vec z^{\prime +}_{\frak c(a)} \subset \Sigma'_0$.
Moreover, by construction, 
$(\Sigma'_0,\vec w^{\prime +}_{\frak c(a)} \cup \vec z^{\prime +}_{\frak c(a)} )$
is in the image of (\ref{mapPhi4532}).
\par
Using this fact and a trivialization of the universal family on $\frak V_{\frak c(a),0}$, we obtain a smooth open embedding from the support of $E_{\frak c(a),0}$ to
$\Sigma'_0$. We use this map together with the parallel transport along the minimal geodesic on
$X$ to send $E_{\frak c(a),0}$ to $\Sigma'$ as (\ref{456formpara}).
\par
We finally discuss Case (1).
If $\frak c(a)$ is either of Type A or C, the construction is the same as in Subsection \ref{subsec:1+}.
(Namely we first define $h_a$ by using a function similar to $f$ (that is defined 
by using additional marked points (Type A) or circles (Type C)). 
Then we use it to send  $E_{\frak c(a)}$ to
$\Sigma'_0$.
\par
Suppose $\frak c(a)$ is of Type B. 
We may assume without loss of generality that $\frak c(a)$ is contained in a tree of extended 
disk components which is attached to $\Sigma^{\frak m}_1$.
We identify
\begin{equation}
\aligned
\Sigma^{\frak m}_1 &=
\{ z \in \C \mid {\rm Re}\, z \in [0,2], \,\,  {\rm Im}\,z \ge 0\} /\sim \\
\Sigma^{\frak m}_2 &=
\{ z \in \C \mid {\rm Re}\, z \in [0,2], \,\,  {\rm Im}\,z \le 0\} /\sim,
\endaligned
\end{equation}
where $z \sim 2 + z$.
For $T$ large, we glue them by identifying $z \in \Sigma^{\frak m}_1$ and  
$z' \in \Sigma^{\frak m}_2$ if $z - 2T\sqrt{-1} = z' + 2T\sqrt{-1}$.
(We then obtain an annulus 
$
\cong \{ z \in \C \mid {\rm Re}\, z \in [0,2], \,\,  {\rm Im}\,z \in  [-T,T]\} /\sim
$
where $z \sim 2 + z$.
\par
We take $\frak z_0 = a \in \partial  \Sigma^{\frak m}_1$ and 
define a coordinate $z \mapsto z+1$, $D^1(1) \cap \frak h \to \Sigma^{\frak m}_1$.
We use this coordinate to define a map
\begin{equation}\label{mapPhi4532222}
\Psi : D^2(\epsilon)^{\mathcal S} \times (T_0,\infty] \times C \times
\frak V_{\frak c(a),0} \times \prod_{\frak a}\frak V_{\frak c(a),\frak a}
\to \mathcal M_{(0,0);\ell_a+m_{a}}^{\rm ann}
\end{equation}
in a similar way as (\ref{mapPhi4532}).
Here $T \in (T_0,\infty]$ is the parameter as above and $C$ is a an open set of an arc
which parametrizes the point of the boundary of $\frak c(a)$ where we attach it to  $\Sigma^{\frak m}_1$ at $\frak z_0$.
We use (\ref{mapPhi4532222})
to obtain a loop $\gamma_c^h$. We use $\gamma_c^h$ to define a function similar to $f_+$ and use it 
to take $h_a$. Then we use $h_a$ to send  $E_{\frak c(a)}$ to
$\Sigma'_0$. 
\par\smallskip
Therefore we have constructed a Kuranishi structure of a 
neighborhood of $
\mathcal M_{(k_1+1,k_2+1);\ell}^{{\rm ann};\text{\rm main}}(\beta)
\cap \frak{forget}^{-1}([\Sigma_1])$.
\par
\medskip
Now, we glue three Kuranishi structures obtained in Subsections \ref{subsec:kuraanul},
\ref{subsec:kuraanul2}, \ref{subsec:kuraanul3}
by taking the sum of the obstruction bundles on the `annulus part'.
Here `annulus part' means the components contained in $\Sigma_{\rm main}$ in the 
situation of Subsection \ref{subsec:kuraanul2}
and the cases (\ref{spherreq2}), (\ref{disk0comp2}), (\ref{disk0comp3}) in the situation of 
Subsection \ref{subsec:kuraanul3}.
\par
We can take the obstruction spaces of the `annulus part' so that the above sum is a 
direct sum, as follows.
Our construction is done by induction of the stratum and on the dimension of the stratum.
So the actual order to perform the construction is: we first 
construct Kuranishi neighborhoods in the situation of  Subsections 
\ref{subsec:kuraanul2}, \ref{subsec:kuraanul3} and then 
in the situation of Subsection \ref{subsec:kuraanul}.
Therefore while we define the obstruction space (\ref{defnofE454}) we may choose it 
so that it is linearly independent to those we have chosen in the situation of Subsections 
\ref{subsec:kuraanul2}, \ref{subsec:kuraanul3}.
\par
We also note that the way how we send the obstruction bundles supported on the extended disk 
components on the disk bubble (that is the part other than `annulus part' in the above sense), 
coincides for the constructions in Subsections \ref{subsec:kuraanul},
\ref{subsec:kuraanul2}, \ref{subsec:kuraanul3}.
This fact is mostly obvious. The only point to check is that the coordinates we use on the boundary marked points 
on the annuli coincide (up to constant multiplication) for the three constructions.
This is immediate from the explicit choice of such coordinates we gave.
\par
Therefore the proofs of Lemmata \ref{bdryMannu} and   
\ref{perturbnbdsigma1} are now complete. \qed
\par
\section{Proof of Lemma \ref{sublemma:Sigma3limit}}
\label{sec:cornerestimate}.
\par
In this section we prove Lemma \ref{sublemma:Sigma3limit}.
Let us restate it below.
We continue to use the notation of Section \ref{sec:ringhomo}.

\begin{prop}\label{sublemma:Sigma3limitprop}
\begin{equation}\label{3moonaji2}
\lim_{i\to \infty}(\text{\rm ev}_0)_* \left(\frak{forget}^{-1}([\Sigma_{3,i}])\right)
= (\text{\rm ev}_0)_* \left(\frak{forget}^{-1}([\Sigma_{0}])\right).
\end{equation}
\end{prop}
\begin{proof}
Let
$$
\overset{\circ}{\mathcal M}_{1}(\alpha;\beta';\overline{\text{\bf f}}_{a(1)},
{\overline{\text{\bf f}}}_{a(2)};\ell_1,\ell_2;\Delta)
$$
be the subset of 
$\mathcal M_{1}(\alpha;\beta';\overline{\text{\bf f}}_{a(1)},
{\overline{\text{\bf f}}}_{a(2)};\ell_1,\ell_2;\Delta)$
consisting of elements that do not satisfy any of
Condition \ref{1tai1cond}.
(Here $\beta' + \alpha = \beta$. We note that 
$\alpha$ is determined by $\beta$, $\beta'$ under this condition.)

We take an increasing sequence of compact subsets
$$
\mathcal N^{K}_{1}(\alpha;\beta';\overline{\text{\bf f}}_{a(1)},
{\overline{\text{\bf f}}}_{a(2)};\ell_1,\ell_2;\Delta)
$$
of 
$\overset{\circ}{\mathcal M}_{1}(\alpha;\beta';\overline{\text{\bf f}}_{a(1)},
{\overline{\text{\bf f}}}_{a(2)};\ell_1,\ell_2;\Delta)$
indexed by positive integers $K$ such that
$$
\bigcup_K
\mathcal N^{K}_{1}(\alpha;\beta';\overline{\text{\bf f}}_{a(1)},
{\overline{\text{\bf f}}}_{a(2)};\ell_1,\ell_2;\Delta)
=
\overset{\circ}{\mathcal M}_{1}(\alpha;\beta';\overline{\text{\bf f}}_{a(1)},
{\overline{\text{\bf f}}}_{a(2)};\ell_1,\ell_2;\Delta).
$$
We first prove the following:
\begin{lem}\label{claim2627}
\begin{equation}\label{equality2619}
\aligned
&\sum_{\alpha \# \beta' = \beta}\sum_{\ell_1+\ell_2=\ell}\lim_{K\to \infty}
\int_{\mathcal N^{K}_{1}(\alpha;\beta';\overline{\text{\bf f}}_{a(1)},
{\overline{\text{\bf f}}}_{a(2)};\ell_1,\ell_2;\Delta)}
{\rm ev}_0^*(\text{\rm vol}_{L(u)}) \\
&=
(\text{\rm ev}_0)_* \left(\frak{forget}^{-1}([\Sigma_{0}])\right).
\endaligned
\end{equation}
Here $\text{\rm vol}_{L(u)}$ is a smooth differential $n$-form on $L(u)$ such that
$\int_{L(u)}\text{\rm vol}_{L(u)} = 1$
and $\frak{forget}$ in the right hand side is the map given by 
\eqref{forgetmap}:
$$
\frak{forget} : \mathcal M^{\text{\rm main}}_{1;\ell+2}
(\beta;\overline{\text{\bf f}}_{a(1)} \otimes
{\overline{\text{\bf f}}}_{a(2)} \otimes \frak b_{\rm high}^{\otimes \ell})
\to  \mathcal M_{1;2}^{\text{\rm main}}.
$$
\end{lem}
We will define the left hand side of (\ref{equality2619}) precisely below during the proof.
\begin{proof}
Note that we use a family of multisections parametrized by a certain space with a smooth 
probability measure to define the integration along the fiber $(\text{\rm ev}_0)_*$
in \cite[Section 12]{fooo09}.
We first review this definition.
We cover the moduli space 
$ \mathcal M^{\text{main}}_{1;\ell+2}
(\beta;\overline{\text{\bf f}}_{a(1)} \otimes
{\overline{\text{\bf f}}}_{a(2)} \otimes \frak b_{\rm high}^{\otimes \ell})
\cap \frak{forget}^{-1}([\Sigma_0])$
by finitely many Kuranishi neighborhoods
$(V_{\frak p},E_{\frak p},\Gamma_{\frak p},s_{\frak p},\psi_{\frak p})$, 
$\frak p \in \frak P$ of 
$\mathcal M^{\text{main}}_{1;\ell+2}
(\beta;\overline{\text{\bf f}}_{a(1)} \otimes
{\overline{\text{\bf f}}}_{a(2)} \otimes \frak b_{\rm high}^{\otimes \ell})$
so that they consist a part of the good coordinate system 
(see \cite[Lemma A1.11]{fooobook2}) of 
$\mathcal M^{\text{main}}_{1;\ell+2}
(\beta;\overline{\text{\bf f}}_{a(1)} \otimes
{\overline{\text{\bf f}}}_{a(2)} \otimes \frak b_{\rm high}^{\otimes \ell})$.
(Here $E_{\frak p}$ is the obstruction bundle, $\Gamma_{\frak p}$ is a finite group 
acting on $V_{\frak p}$, $s_{\frak p}$ is the Kuranishi map, and $\psi_{\frak p}$
is a homeomorphism from $s_{\frak p}^{-1}(0)/\Gamma_{\frak p}$ 
to an open set of our moduli space 
$\mathcal M_{1}(\alpha;\beta';\overline{\text{\bf f}}_{a(1)},
{\overline{\text{\bf f}}}_{a(2)};\ell_1,\ell_2;\Delta)$.)
\par
By Lemma \ref{perturbnbdsigmasec9}.6 and 
Lemma \ref{normalcrossing}, the
Kuranishi neighborhood of 
elements of 
$\mathcal M_{1}(\alpha;\beta';\overline{\text{\bf f}}_{a(1)},
{\overline{\text{\bf f}}}_{a(2)};\ell_1,\ell_2;\Delta)$
is obtained from $(V_{\frak p},E_{\frak p},\Gamma_{\frak p},s_{\frak p},\psi_{\frak p})$ 
as follows.
\par
Let $\frak x \in \mathcal M_{1}(\alpha;\beta';\overline{\text{\bf f}}_{a(1)},
{\overline{\text{\bf f}}}_{a(2)};\ell_1,\ell_2;\Delta)$
and $\text{\rm Glue}(\frak x) = \psi_{\frak p}([\frak y])$ with 
$\frak y \in  s_{\frak p}^{-1}(0) \subset V_{\frak p}$.
Then a Kuranishi neighborhood of $\frak x$ is 
$(V^-_{\frak p}(\beta'),(\Gamma_{\frak p})_{\frak y},s_{\frak p},\text{\rm Glue}^{-1}\circ\psi_{\frak p})$,
where $V^-_{\frak p}(\beta')$ is a codimension $2$ submanifold of $V_{\frak p}$, 
$(\Gamma_{\frak p})_{\frak y} = \{\gamma \in \Gamma_{\frak p} \mid \gamma\frak y = \frak y\}$
and we denote the restrictions of $s_{\frak p}$, $\psi_{\frak p}$ by the same symbol.
Note the forgetful map 
$$
\frak{forget} : \mathcal M_{1}(\alpha;\beta';\overline{\text{\bf f}}_{a(1)},
{\overline{\text{\bf f}}}_{a(2)};\ell_1,\ell_2;\Delta) \to \mathcal M_{1;2}^{\text{\rm main}}
$$
is extended to the Kuranishi neighborhoods and 
we may take
$$
V^-_{\frak p}(\beta') \subseteq V_{\frak p} \cap \frak{forget}^{-1}([\Sigma_0]).
$$
Moreover, $V^-_{\frak p}(\beta')$ contains an open neighborhood of $\frak y$ in
$V_{\frak p} \cap \frak{forget}^{-1}([\Sigma_0])$ if 
$\frak x \in \overset{\circ}{\mathcal M}_{1}(\alpha;\beta';\overline{\text{\bf f}}_{a(1)},
{\overline{\text{\bf f}}}_{a(2)};\ell_1,\ell_2;\Delta)
$.
(Otherwise there are finitely many codimension two submanifolds 
which contain $\frak x$ and intersect transversally so that 
their union contains a neighborhood of $\frak x$ in 
$V_{\frak p} \cap \frak{forget}^{-1}([\Sigma_0])$.)
\par
We can choose a good coordinate system of 
$\mathcal M_{1}(\alpha;\beta';\overline{\text{\bf f}}_{a(1)},
{\overline{\text{\bf f}}}_{a(2)};\ell_1,\ell_2;\Delta)$ and  
a good coordinate system 
of $\mathcal M^{\text{main}}_{1;\ell+2}
(\beta;\overline{\text{\bf f}}_{a(1)} \otimes
{\overline{\text{\bf f}}}_{a(2)} \otimes \frak b_{\rm high}^{\otimes \ell})$
on a neighborhood of  $\frak{forget}^{-1}([\Sigma_0])$
so that the index sets are both $\frak P$ and the members of the good coordinate system  
of $\mathcal M_{1}(\alpha;\beta';\overline{\text{\bf f}}_{a(1)},
{\overline{\text{\bf f}}}_{a(2)};\ell_1,\ell_2;\Delta)$ are
$(V^-_{\frak p}(\beta'),(\Gamma_{\frak p})_{\frak y},s_{\frak p},\text{\rm Glue}^{-1}\circ\psi_{\frak p})$.
(This is a consequence of Lemma \ref{perturbnbdsigmasec9}.6 and 
the construction of good coordinate system.)
Moreover, we may assume 
$(\Gamma_{\frak p})_{\frak y} = \Gamma_{\frak p}$.
\par
To define $({\rm ev}_0)_*$ we take a finite dimensional manifold
$W_{\frak p}$ for each $\frak p$ together with $\omega_{\frak p}$
a top degree form with compact support on $W_{\frak p}$.
We also take a partition of unity $\chi_{\frak p}$ associated to our good coordinate system.
(See \cite[Definition 12.10]{fooo09}.)
\par
Our family of multisections, $\frak s$, gives a multisection $\frak s_{\frak p}$ on $V_{\frak p}/\Gamma_{\frak p}
\times W_{\frak p}$ for each $\frak p$.
They satisfy appropriate compatibility conditions, as formulated in \cite[Definition 12.8]{fooo09}.
(We may also regard $s_{\frak p}$ as a $\Gamma_{\frak p}$ equivariant 
map $V_{\frak p} \times W_{\frak p} \to S^{l_{\frak p}}(E_{\frak p})$,
where $S^{l_{\frak p}}(E_{\frak p})$ is the $l_{\frak p}$-th symmetric power of 
the fiber of our vector bundle. We use the same symbol $s_{\frak p}$ 
by abuse of notation.)
\par
The zero set $\frak s_{\frak p}^{-1}(0) \cap (V_{\frak p} \times W_{\frak p})$  is a finite union of smooth submanifolds. 
(Note $\frak s_{\frak p}^{-1}(0)$ is the union of the zero sets of the branches. 
The zero set of each branch is a smooth submanifold by transversality.)
Moreover, its zero set is transversal to $\text{\rm Glue}$ in the sense of 
Lemma \ref{perturbnbdsigmasec9}.6. Namely 
$\frak s_{\frak p}^{-1}(0) \cap V_{\frak p}$ is transversal to 
$V_{\frak p}^{-}(\beta)$.
\par
Now by definition we have
\begin{equation}\label{form2620}
\aligned
&(\text{\rm ev}_0)_* \left(\frak{forget}^{-1}([\Sigma_{0}])\right)\\
&= 
\sum_{\frak p \in \frak P} 
\frac{1}{\#\Gamma_{\frak p} l_{\frak p}}
\int_{\frak s_{\frak p}^{-1}(0) \cap V_{\frak p} \cap \frak{forget}^{-1}([\Sigma_0])}
\chi_{\frak p} (\text{\rm ev}_0)^{*}({\rm vol}_{L(u)}) \wedge \omega_{\frak p}.
\endaligned
\end{equation}
Here $l_{\frak p}$ is the number of branches.
\begin{rem}
In \cite[Section 12]{fooo09} we cover each Kuranishi neighborhood 
by finitely many open sets so that on each open set the number of 
branches does not change. (See \cite[Definition 12.5]{fooo09}.)
We omit this process here for simplicity. Actually we can subdivide the Kuranishi neighborhood 
$V_{\frak p}$ furthermore so that the number of branches does not change 
on each $V_{\frak p}$ already  without taking finer cover.
\end{rem}
We note that the intersection $(V_{\frak p} \cap \frak{forget}^{-1}([\Sigma_0])) \times 
W_{\frak p}$ is {\it not} a 
smooth submanifold in general. However, it is a finite union of smooth submanifolds
$V_{\frak p}^{-}(\beta') \times W_{\frak p}$ for various $\beta'$.
Moreover, the zero set of each branch of $\frak s_{\frak p}$ 
is transversal to these submanifolds by Lemma \ref{perturbnbdsigmasec9}.6.
Therefore the integration in the right hand side of (\ref{form2620}) is well defined.
\par
We next discuss the left hand side of (\ref{equality2619}).
We note that an element of $V_{\frak p}$ is regarded as $((\Sigma,\vec z,w),(\frak y_1,\frak y_2))$, 
where 
$(\Sigma,\vec z,w)$ is a pair of marked bordered semi-stable curve
$(\Sigma,\vec z)$ and a map $w : \Sigma \to X$, and
$(\frak y_1,\frak y_2) \in \overline{\text{\bf f}}_{a(1)} \times \overline{\text{\bf f}}_{a(2)}$.
(The map $w$ may not be pseudo-holomorphic.) 
So if $(\Sigma,\vec z,w) \in \frak{forget}^{-1}([\Sigma_0])$, 
we may say whether it satisfies Condition \ref{1tai1cond} or not.
We define  $\overset{\circ}V_{\frak p}^-(\beta')$ 
as the set of the elements of $V_{\frak p}^-(\beta')$ such that none of 
Condition \ref{1tai1cond} is satisfied.
They consist a Kuranishi neighborhood of 
$
\overset{\circ}{\mathcal M}_{1}(\alpha;\beta';\overline{\text{\bf f}}_{a(1)},
{\overline{\text{\bf f}}}_{a(2)};\ell_1,\ell_2;\Delta).
$
We choose  $V_{\frak p,K}$ for each $K=1,2,\dots$ with the following 
properties.
\begin{conds}\label{cond2629}
\begin{enumerate}
\item
$V_{\frak p,K} \subset V_{\frak p,K+1}$.
\item 
$V_{\frak p,K} \cap V_{\frak p}^-(\beta') \subset \overset{\circ}V_{\frak p}^-(\beta')$.
\item
If $\overset{\circ}V_{\frak p}^-(\beta') = V_{\frak p}^-(\beta')$, then 
$V_{\frak p,K} \supseteq V_{\frak p}^-(\beta')$.
\item
$$
\bigcup_{K=1}^{\infty}(V_{\frak p,K} \cap V_{\frak p}^-(\beta'))
=  \overset{\circ}V_{\frak p}^-(\beta').
$$
\item
$$
\aligned
&\text{Glue}\circ\psi_{\frak p}((s_{\frak p}^{-1}(0) \cap V_{\frak p,K} \cap V_{\frak p}^-(\beta'))/\Gamma_{\frak p})
\\
&=
\text{Glue}\circ\psi_{\frak p}((s_{\frak p}^{-1}(0) \cap V_{\frak p}^-(\beta'))/\Gamma_{\frak p})
\cap
\mathcal N^{K}_{1}(\alpha;\beta';\overline{\text{\bf f}}_{a(1)},
{\overline{\text{\bf f}}}_{a(2)};\ell_1,\ell_2;\Delta).
\endaligned
$$
Here we recall that $s_{\frak p}$ is the Kuranishi map of the Kuranishi structure of $\mathcal M_{1}(\alpha;\beta';\overline{\text{\bf f}}_{a(1)},
{\overline{\text{\bf f}}}_{a(2)};\ell_1,\ell_2;\Delta)$, which is not necessarily transversal to $0$.
We remark that $s_{\frak p} \ne \frak s_{\frak p}$.
\end{enumerate}
\end{conds}
Now we define
\begin{equation}\label{form2621}
\aligned
&\int_{\mathcal N^{K}_{1}(\alpha;\beta';\overline{\text{\bf f}}_{a(1)},
{\overline{\text{\bf f}}}_{a(2)};\ell_1,\ell_2;\Delta)}
{\rm ev}_0^*(\text{\rm vol}_{L(u)})\\
&=
\sum_{\frak p \in \frak P} \frac{1}{l_{\frak p}\#\Gamma_{\frak p}}\int_{\frak s_{\frak p}^{-1}(0) \cap   ((V_{\frak p,K}
\cap V_{\frak p}^-(\beta'))\times W_{\frak p})}
\chi_{\frak p} (\text{\rm ev}_0)^{*}(\text{\rm vol}_{L(u)}) \wedge \omega_{\frak p}.
\endaligned
\end{equation}
Condition \ref{cond2629}.5 above justifies the notation in the left hand side.
Now Lemma \ref{claim2627} follows 
from (\ref{form2620}), (\ref{form2621}) and Condition \ref{cond2629}.4 above.
\end{proof}
We next prove:
\begin{lem}\label{claim2628}
$$
\aligned
&\lim_{i\to \infty}
\int_{\frak s_{\frak p}^{-1}(0) \cap   ((V_{\frak p,K} \cap \frak{forget}^{-1}([\Sigma_{i,3}])\times W_{\frak p})
}
\chi_{\frak p} (\text{\rm ev}_0)^{*}(\text{\rm vol}_{L(u)}) \wedge \omega_{\frak p}
\\
&=
\sum_{\beta'+\alpha = \beta}
\int_{\frak s_{\frak p}^{-1}(0) \cap   ((V_{\frak p,K}
\cap V_{\frak p}^-(\beta')) \times W_{\frak p})}
\chi_{\frak p} (\text{\rm ev}_0)^{*}(\text{\rm vol}_{L(u)}) \wedge \omega_{\frak p}.
\endaligned
$$
\end{lem}
\begin{proof}
Note the union of $\frak s_{\frak p}^{-1}(0) \cap   ((V_{\frak p,K}
\cap V_{\frak p}^-(\beta')) \times W_{\frak p})$
over $\beta'$ is a disjoint union.
Moreover,  
$\frak s_{\frak p}^{-1}(0) \cap   ((V_{\frak p,K}\cap \frak{forget}^{-1}([\Sigma_{i,3}]))\times W_{\frak p}) $ 
converges to 
$$
\bigcup_{\beta'+\alpha = \beta}
\frak s_{\frak p}^{-1}(0) \cap   \overline{((V_{\frak p,K}
\cap V_{\frak p}^-(\beta')) \times W_{\frak p})}
$$
as $i$ goes to infinity, in $C^1$ topology of submanifolds.
(The proof of this $C^1$ convergence uses the fiber product 
description of the Kuranishi neighborhood that is 
similar to and easier than the discussion appearing in the 
proof of Lemma \ref{claim26299} below. So we omit the detail.)
The lemma follows.
\end{proof}
The next lemma is the most important part in the proof of Proposition \ref{sublemma:Sigma3limitprop}.
\begin{lem}\label{claim26299}
For any positive number $\epsilon$ there exists $K_0$ 
and for any $K$ there exists $I_0(\epsilon,K)$ such that, if $K > K_0$,  $i>I_0(\epsilon,K)$,  
$\beta'$, $\frak p \in \frak P$, then we have an inequality:
\begin{equation}\label{eq2622}
\left\vert
\int_{\frak s_{\frak p}^{-1}(0) \cap   (((V_{\frak p}\setminus V_{\frak p,K}) \cap \frak{forget}^{-1}([\Sigma_{i,3}]))\times W_{\frak p})}
\chi_{\frak p} (\text{\rm ev}_0)^{*}(\text{\rm vol}_{L(u)}) \wedge \omega_{\frak p}\right\vert
< \epsilon.
\end{equation}
\end{lem}
\begin{proof}
For simplicity of notation, we prove the case when $\ell_1 = \ell_2 = 0$. 
(Namely we consider the case when we do not take bulk deformation.)
The general case is similar.
\par
We consider the Kuranishi neighborhood $(V_{\frak p},E_{\frak p},\Gamma_{\frak p},s_{\frak p},\psi_{\frak p})$.
We have $\widehat V_{\frak p}$ together with evaluation maps $\widehat V_{\frak p} \to X^2$ such that
$$
V_{\frak p} = \widehat V_{\frak p} \times_{X^2} (\overline{\text{\bf f}}_{a(1)} \times
{\overline{\text{\bf f}}}_{a(2)}).
$$
We have a corresponding Kuranishi neighborhood, which we write 
$(\widehat V_{\frak p},E_{\frak p},\Gamma_{\frak p},s_{\frak p},\psi_{\frak p})$.
This Kuranishi neighborhood 
is centered at a certain element $(\Sigma_{\frak p},\vec z_{\frak p},w_{\frak p})$.
Here $(\Sigma_{\frak p},\vec z_{\frak p})$ is a bordered semi-stable curve with two interior and 
one exterior marked points and of genus $0$ with $1$ boundary component and
$w_{\frak p} : (\Sigma_{\frak p},\partial\Sigma_{\frak p}) \to (X,L(u))$ is 
pseudo-holomorphic with homology class $\beta$.
Moreover $\frak{forget}([\Sigma_{\frak p},\vec z_{\frak p}]) = [\Sigma_0]$.
\par
If $(\Sigma_{\frak p},\vec z_{\frak p},w_{\frak p})$ satisfies one of the Condition 
\ref{1tai1cond}, 
then we have $V_{\frak p} = V_{\frak p,K}$ for sufficiently large $K$.
Therefore there is nothing to prove in this case.
\par
We assume that 
$(\Sigma_{\frak p},\vec z_{\frak p},w_{\frak p})$ does not satisfy Condition 
\ref{1tai1cond}.
For simplicity of notation we discuss the case when this element is of the 
form of Example \ref{doublepointexam}. The general case is 
similar.
We also consider the case of $\Gamma_{\frak p} = \{1\}$ for simplicity.
We use the notation of Example \ref{doublepointexam}.
We have
$$
\Sigma_{\frak p} = D^2 \cup S^2_1 \cup S^2_2
$$ 
and $\vec z_{\frak p} = (z_0;z_1^{\rm int},z_2^{\rm int})$ 
where $z_0$ is a boundary marked point and 
$z_1^{\rm int},z_2^{\rm int}$ are interior marked points so that 
$z_0 \in D^2$, $z_1^{\rm int},z_2^{\rm int} \in S^2_2$.
We put $D^2 \cap S^2_1 = \{p_1\}$, 
$S^2_1 \cap S^2_2 = \{p_2\}$.
We also put
\begin{equation}\label{frakxinMM}
\aligned
&\frak x_0 = (D^2,(z_0;p_1),w\vert_{D^2}) \in \mathcal M_{1;1}(\beta_{(0)}),\\
&\frak x_1 = (S^2_1,(p_1,p_2),w\vert_{S^2_1}) \in \mathcal M_{2}(\alpha_1),\\
&\frak x_2 = (S^2_2,(p_2,z_1^{\rm int},z_2^{\rm int}),w\vert_{S^2_2}) \in \mathcal M_{3}(\alpha_2).
\endaligned
\end{equation}
Here $\mathcal M_{2}(\alpha_1)$ (resp. $\mathcal M_{3}(\alpha_2)$) is the 
moduli space of marked stable maps of genus $0$ in $X$ with $2$ (resp. 3) marked points
and of homology class $\alpha_1$ (resp. $\alpha_2$).
Note $\beta_{(0)} + \alpha_1 + \alpha_2 = \beta$.
\par
We describe Kuranishi neighborhoods of $\frak x_i$ of the Kuranishi structures of 
the moduli spaces in the right hand side of (\ref{frakxinMM}).
\par
We start with $\frak x_0$. We take 
$$
E_0(w) \subset C^{\infty}(D^2,w^*TX \otimes \Lambda^{0,1}).
$$
This is a finite dimensional linear space of smooth sections.
(\cite[(12.7.4)]{FO}.)
We assume that the support of its elements are away from 
the marked points and does not intersect with boundary 
$\partial D^2$.
\par
Let $w'_0 : (D^2,\partial D^2) \to (X,L(u))$ be a smooth map 
that is $C^{1}$ close to $w\vert_{D^2}$.
We use parallel transport along the minimal geodesic 
joining $w(z)$ and $w'_0(z)$ to send $E_0(w)$ to 
$$
E_0(w'_0) \subset C^{\infty}(D^2,(w'_0)^*TX \otimes \Lambda^{0,1}).
$$
Now let $\widehat V_0$ be the moduli space of maps $w'_0$ that satisfy
\begin{equation}
\overline{\partial} w'_0 \equiv 0 \mod E_0(w'_0)
\end{equation}
and are $C^{1}$ close to $w\vert_{D^2}$. The family of vector spaces $E_0(w'_0)$ 
parametrized by $w'_0$ defines a smooth vector bundle on $\widehat V_0$, which we denote by $E_0$.
Then the map $w'_0\mapsto \overline\partial w'_0$ defines a 
section of this bundle $E_0$, which is the Kuranishi map.
We have evaluation maps
$$
({\rm ev}_0,{\rm ev}^p_1) : \widehat V_0 \to L(u) \times X,
$$
at $z_0$ and $p_1$.
We may choose $E_0(w)$ so that $\widehat V_0$ is a smooth manifold and 
$({\rm ev}_0,{\rm ev}^p_1)$ is a submersion.
\par
We next study $\frak x_1$. Since the source is $S^2$ with two marked points 
and is 
not stable, we use the method of \cite[appendix]{FO} to stabilize the source 
by adding a marked point as follows.
We take and fix a point $p_3 \in S^2_1$ such that $p_3$ is different from $p_1,p_2$ and 
that $w$ is an immersion at $p_3$. 
We also fix a codimension $2$ smooth submanifold $\mathcal D$ of $X$ with the following 
property. Let $U$ be a sufficiently small neighborhood of $p_3$ in $S^2_1$. Then 
$\mathcal D$ intersects with $w(U)$ transversally at one point $w(p_3)$.
\par
We take 
$$
E_1(w) \subset C^{\infty}(S^2_1,w^*TX \otimes \Lambda^{0,1}).
$$
This is a finite dimensional linear space of smooth sections.
We assume that the support of its elements are away from 
the marked points $p_1,p_2,p_3$.
\par
Let $w'_1 : S^2_1 \to X$ be a smooth map 
that is $C^{1}$ close to $w\vert_{S^2_1}$.
We use parallel transport along the minimal geodesic 
joining $w(z)$ and $w'_1(z)$ to send $E_1(w)$ to 
$$
E_1(w'_1) \subset C^{\infty}(S^2_1,(w'_1)^*TX \otimes \Lambda^{0,1}).
$$
Now let $\widehat V_1^{+}$ be the moduli space of maps $w'_1$ that satisfy
\begin{equation}
\overline{\partial} w'_1 \equiv 0 \mod E_1(w'_1).
\end{equation}
We have evaluation maps
$$
({\rm ev}^p_1,{\rm ev}^p_2,{\rm ev}^p_3) : \widehat V_1^+ \to  X^3,
$$
at $p_1,p_2,p_3$.
We require an extra condition that
\begin{equation}\label{transcons}
{\rm ev}^p_3(w'_1) = w'_1(p_3) \in \mathcal D.
\end{equation}
We denote by $\widehat V_1$ the subset of $\widehat V_1^+$ consisting the elements 
$w'_1$ satisfying (\ref{transcons}).
\par
We may choose $E_1(w)$ so that $\widehat V_1^+$ is a smooth manifold 
and $\widehat V_1$ is a codimension $2$ smooth submanifold.
Moreover we may choose $E_1(w)$ so that
$$
({\rm ev}_2,{\rm ev}_3) : \widehat V_1 \to  X^2
$$
is a submersion.
\par
We next consider $\frak x_3$. Since the source has three marked points and so is stable, 
we do not need to add marked points in this case.
We take
$$
E_2(w) \subset C^{\infty}(S^2_2,w^*TX \otimes \Lambda^{0,1}).
$$
This is a finite dimensional linear space of smooth sections whose supports 
are away from $p_2, z_1^{\rm int},z_2^{\rm int}$.
\par
If $w'_2$ is $C^{1}$ close to $w\vert_{S^2_2}$, we use parallel 
transport along minimal geodesic to obtain 
$$
E_2(w'_2) \subset C^{\infty}(S^2_2,(w'_2)^*TX \otimes \Lambda^{0,1}).
$$
Now let $\widehat V_2$ be the moduli space of maps $w'_2$ that satisfy
\begin{equation}
\overline{\partial} w'_2 \equiv 0 \mod E_2(w'_2).
\end{equation}
We have evaluation maps
$$
({\rm ev}^p_2,{\rm ev}_1,{\rm ev}_2) : \widehat V_2 \to  X^3,
$$
at $p_2,z_1^{\rm int},z_2^{\rm int}$.
We may choose $E_2(w)$ so that $\widehat V_2$ is a smooth manifold 
and $({\rm ev}^p_2,{\rm ev}_1,{\rm ev}_2)$ is a submersion.

By the submersivity of the evaluation maps, the following 
fiber product is transversal.
$$
\widehat V_{\infty}^+ 
=
\widehat V_0 \,\,{}_{{\rm ev}_1^p} \times_{{\rm ev}_1^p}
\widehat V_1^+  {}_{{\rm ev}_2^p} \times_{{\rm ev}_2^p}
\widehat V_2.
$$
Moreover, the evaluation map
$$
({\rm ev}_0,{\rm ev}_1,{\rm ev}_2) :
\widehat V_{\infty}^+ 
\to 
L(u) \times X^2
$$
at $z_0, z_1^{\rm int}, z_2^{\rm int}$ is a submersion.
\par
The space $\widehat V_{\infty}^+$ parametrizes a family of marked stable maps 
$(\Sigma,\partial \Sigma) \to (X,L(u))$ such that 
the source $\Sigma$ has two singular points.
We next perform the gluing construction at these two singular points.
(Below we follow the argument given in \cite[Section A1.4]{fooobook2}.
See also \cite[Section 21]{foootech} for detail.)
\par
We consider the space 
$((T_0,\infty] \times S^1)/\sim$ where 
$(T,\theta) \sim (T',\theta')$ 
if and only if either $(T,\theta) = (T',\theta')$ 
or $T=T'=\infty$.
(Here $\theta,\theta' \in S^1 = \R/\Z$.)
We note that the quotient space 
$((T_0,\infty] \times S^1)/\sim$ is 
homeomorphic to $D^2(1/T_0) 
= \{ z \in \C \mid \vert z\vert < 1/T_0\}$ by the map 
which sends $(T,\theta)$ to $\exp(2\pi \sqrt{-1} \theta)/T$ and 
$[(\infty,\theta)]$ to $0$.
(This defines the smooth coordinate of our Kuranishi structure 
that was introduced in \cite[Section A1.4]{fooobook2}. 
See \cite[Section 21]{foootech} for detail.)
\par
We identify 
\begin{equation}\label{coordinateforglue}
\aligned
D^2 \setminus \{p_1\}
&= [0,\infty) \times S^1, \\
S^2_1 \setminus \{p_1, p_2\}
&= \R \times S^1, \\
S^2_2 \setminus \{p_2, z_1^{\rm int}\}
&= \R \times S^1,
\endaligned
\end{equation}
where $(\tau_0,t_0)$, $(\tau_1,t_1)$, and $(\tau_2,t_2)$
are the coordinates of the first, second, and third lines of 
the right hand side, respectively.
We take the identification (\ref{coordinateforglue})
so that $p_1$ corresponds to either $\tau_0 = + \infty$ or 
$\tau_1 = - \infty$, 
$p_2$ corresponds to either $\tau_1 = + \infty$ or 
$\tau_2 = - \infty$, 
and $z_1^{\rm int}$ corresponds to $\tau_2 = + \infty$.
Moreover, $z_0$ is $(\tau_0,t_0) = (0,0)$, 
$z_2^{\rm int}$ is $(\tau_2,t_2) = (0,0)$,
and $p_3$ is $(\tau_1,t_1) = (0,0)$.
(These conditions determine the 
conformal and orientation preserving diffeomorphism (\ref{coordinateforglue}) 
uniquely.)
\par 
Let $T_1,T_2 > T_0$ and $\theta_1,\theta_2 \in S^1 = \R/\Z$.
We put
\begin{equation}
\aligned
&\tau_0 = \tau_1 + 10T_1,\qquad
\tau_1 = \tau_2 + 10T_2, \\
&t_0 = t_1 - \theta_1,\qquad
t_1 = t_2  - \theta_2.
\endaligned
\end{equation}
We glue $D^2, S^2_1, S^2_2$ by this identification and obtain 
$$
\Sigma(T_1,T_2;\theta_1,\theta_2)
=
([0,\infty) \times S^1) \cup \{z_1^{\rm int}\} \cong D^2.
$$
More precisely, $\Sigma(T_1,T_2;\theta_1,\theta_2)$ 
is a union of three subsets
$$
[0,5T_1) \times S^1 \subset D^2, 
\quad
[-5T_1,5T_2] \times S^1 \subset S^2_1,
\quad
([-5T_2,\infty) \times S^1) \cup \{z_1^{\rm int}\} \subset S^2_2.
$$
In case either $T_1 =\infty$, $T_2 \ne \infty$ or 
$T_1 \ne \infty$, $T_2 =\infty$, 
the space $\Sigma(T_1,T_2;\theta_1,\theta_2)$ has one nodal singularity.
We put $\rho = (T_1,T_2;\theta_1,\theta_2)$ and 
write $\Sigma_{\rho} = \Sigma(T_1,T_2;\theta_1,\theta_2)$.
\par
Let $C_i$ ($i=0,1,2$) be sufficiently large positive  numbers. We denote
$$
K_0 = [0,C_0] \times S^1 \subset D^2, 
\quad
K_1 = [-C_0,C_1] \times S^1 \subset S^2_1, 
\quad
K_2 = [-C_1,C_2] \times S^1 \subset S^2_2.
$$
We may also regard $K_i$ as a subset of  $\Sigma_{\rho}$.
We put $K = K_1 \cup K_2 \cup K_3 \subset \Sigma_{\rho}$ 
and call it the {\it core}\index{core (of a stable curve)}.
We define Riemannian metics on $K_i$ by restricting the standard metric of $D^2$ or $S^2_i$.
We note that $K_i$ and $K$ are independent of $\rho$.
\par
We call the complement $\Sigma_{\rho} \setminus K$ 
the {\it neck region}\index{neck region (of a stable curve)}.
\par
Let $\frak w = (w'_0,w'_1,w'_2)$ be an element of  
$
\widehat V_{\infty}^+ 
=
\widehat V_0 \,\,{}_{{\rm ev}_1^p} \times_{{\rm ev}_1^p}
\widehat V_1^+  {}_{{\rm ev}_2^p} \times_{{\rm ev}_2^p}
\widehat V_2.
$
We define a map
$$
w^1_{\frak w,\rho} : (\Sigma_{\rho},\partial \Sigma_{\rho}) \to (X,L(u))
$$
so that $w^1_{\frak w,\rho} = w'_i$ on $K_i$ for $i=0,1,2$.
We use an appropriate bump function to define it.
(See \cite[Subsection 7.1.3]{fooobook2} and \cite[Section 19]{foootech}.)
\par
We take $C_i$ large so that the supports of the elements of 
$E_i$ are contained in the core $K_i$.
Let $w' :  (\Sigma_{\rho},\partial \Sigma_{\rho}) \to (X,L(u))$
be a smooth map that is $C^{1}$ close to 
$w^1_{\frak w,\rho}$ for some $\frak w$.
Using parallel transport along the minimal geodesic, 
we send $E_i$ to
$$
E_i(w') \subset C^{\infty}(\Sigma_{\rho},(w')^*TX \otimes \Lambda^{0,1}).
$$
We put
$$
E(w') = E_1(w') \oplus E_2(w') \oplus E_3(w').
$$
Now by construction we have $\xi$ such that
$$
\overline\partial w^1_{\frak w,\rho} \equiv \xi \mod E(w')
$$
and
$$
\Vert\xi\Vert_{L^2_m} \le C e^{-c \min\{ T_1,T_2\}}.
$$
Here $C,c$ are positive numbers independent of $T_1,T_2$.
\par
Then if $T_0$ is sufficiently large, we can use Newton's iteration 
together with alternating method to obtain 
$$
w^2_{\frak w,\rho} :  (\Sigma_{\rho},\partial \Sigma_{\rho}) \to (X,L(u))
$$
such that
\begin{equation}\label{equationOKw2}
\overline\partial w^2_{\frak w,\rho} \equiv 0 \mod E(w^2_{\frak w,\rho}).
\end{equation}
It also enjoys the following decay estimate of exponential order:
\begin{equation}\label{decayequationOKw2}
\left\Vert
\frac{\partial  w^2_{\frak w,\rho}}{\partial T_i}
\right\Vert_{C^m(K)} \le C e^{-c T_i}.
\end{equation}
The construction of $w^2_{\frak w,\rho}$ satisfying 
(\ref{equationOKw2}) and (\ref{decayequationOKw2}) 
is written in detail in \cite[Subsection A1.4]{fooobook2} and \cite[Section 19]{foootech}.
Note when $\rho = \rho_{\infty} = (\infty,\infty)$ (in other words when $T_1=T_2=\infty$) the maps $w^2_{\frak w,\rho_{\infty}} = w^1_{\frak w,\rho_{\infty}}$ 
are $w'_0$, $w'_1$, $w'_2$ on $D^2$, $S^2_1$, $S^2_2$, respectively.
\par
We have thus defined a family of maps 
$
w^2_{\frak w,\rho} :  (\Sigma_{\rho},\partial \Sigma_{\rho}) \to (X,L(u))
$
parametrized by 
$$
(\frak w,\rho) \in \widehat V_{\infty}^+ \times ((T_0,\infty] \times S^1)/\sim)^2.
$$ 
Let $\widehat{\frak V}$ be the subset of $\widehat V^+ \times ((T_0,\infty] \times S^1)/\sim)^2$ 
defined by the equations
\begin{equation}\label{transconstraintD}
w^2_{\frak w,\rho}(p_3) \in \mathcal D.
\end{equation}
We put
$$
\widehat V_{\infty} 
=
\widehat V_0 \,\,{}_{{\rm ev}_1^p} \times_{{\rm ev}_1^p}
\widehat V_1  {}_{{\rm ev}_2^p} \times_{{\rm ev}_2^p}
\widehat V_2.
$$ 
Here $\widehat V_{\infty}$ is a codimension $2$ submanifold of  $\widehat V_{\infty}^+$ 
and is the intersection of $\widehat{\frak V}$ with $T_1=T_2=\infty$. 
We then define
$$
\frak V = \widehat{\frak V} \,\,{}_{({\rm ev}_1,{\rm ev}_2)} \times_{X^2}  (\overline{\text{\bf f}}_{a(1)} \times
{\overline{\text{\bf f}}}_{a(2)})
$$
and 
$$
V_{\infty}  = \widehat V_{\infty} \,\,{}_{({\rm ev}_1,{\rm ev}_2)} \times_{X^2}  (\overline{\text{\bf f}}_{a(1)} \times
{\overline{\text{\bf f}}}_{a(2)}).
$$
The intersection of ${\frak V}$ with $T_1=T_2=\infty$ is $V_{\infty}$.
\par
Using the inequality (\ref{decayequationOKw2}), we can find a map
$$
\Phi :  V_{\infty} \times ((T_0,\infty] \times S^1)/\sim)^2
\to \widehat V_{\infty}^+ \times  (\overline{\text{\bf f}}_{a(1)} \times
{\overline{\text{\bf f}}}_{a(2)})
$$
with the following properties:
\begin{enumerate}
\item[(A)]
$\frak V$ is the graph of $\Phi$. Namely
$$
\frak V = \{ (\Phi(\frak w,\rho),\rho) \mid (\frak w,\rho) \in V_{\infty} \times ((T_0,\infty] \times S^1)/\sim)^2\}.
$$
(We note that both the right hand side and the left hand side  are subsets of 
$
\widehat V_{\infty}^+ \times ((T_0,\infty] \times S^1)/\sim)^2 
\times (\overline{\text{\bf f}}_{a(1)} \times
{\overline{\text{\bf f}}}_{a(2)}).
$)
\item[(B)]
If $\rho = \rho_{\infty} = (\infty,\infty)$, then $\Phi(\frak w,\rho_{\infty}) = \frak w$.
\item[(C)]
$$
\left\Vert
\frac{\partial  \Phi}{\partial T_i}
\right\Vert_{C^m} \le C e^{-c T_i}.
$$
\end{enumerate}
The existence of such $\Phi$ is immediate from implicit function theorem and (\ref{decayequationOKw2}).
(It is a special case of \cite[Proposition 21.1]{foootech}.)
\par
We put
\begin{equation}
w_{\frak w,\rho} = w^2_{\Phi_1(\frak w,\rho),\rho}.
\end{equation}
Here $\frak w \in V_{\infty}$ and 
$\Phi_1(\frak w,\rho)$ is the $\widehat V_{\infty}^+$ factor of $\Phi(\frak w,\rho) \in \widehat V_{\infty}^+ \times  (\overline{\text{\bf f}}_{a(1)} \times
{\overline{\text{\bf f}}}_{a(2)})$. 
We have
\begin{equation}\label{equationOKw332}
\overline\partial w_{\frak w,\rho} \equiv 0 \mod E(w_{\frak w,\rho})
\end{equation}
and
\begin{equation}\label{decayequationOKw233}
\left\Vert
\frac{\partial  w_{\frak w,\rho}}{\partial T_i}
\right\Vert_{C^m(K)} \le C e^{-c T_i}.
\end{equation}
We now put
\begin{equation}
V_{\frak p} = V_{\infty} \times ((T_0,\infty] \times S^1)/\sim)^2.
\end{equation}
We define a vector bundle $E_{\frak p}$ on it so that
its fiber at $(\frak w,\rho)$ is $E(w_{\frak w,\rho})$.
We define a section $s_{\frak p}$ of $E_{\frak p}$ by
\begin{equation}
s_{\frak p}(\frak w,\rho) = \overline\partial w_{\frak w,\rho}.
\end{equation}
(\ref{decayequationOKw233}) implies 
\begin{equation}\label{decayequationOKw233222}
\left\Vert
\frac{\partial  s_{\frak p}}{\partial T_i}
\right\Vert_{C^m} \le C e^{-c T_i}.
\end{equation}
We next define $\psi_{\frak p} : s_{\frak p}^{-1}(0) \to 
\mathcal M^{\text{main}}_{1;2}
(\beta;\overline{\text{\bf f}}_{a(1)} \otimes
{\overline{\text{\bf f}}}_{a(2)})$ by
\begin{equation}
\psi_{\frak p}(\frak w,\rho) = ((\Sigma_{\rho},z_0,z_1^{\rm int},z_2^{\rm int}, w_{\frak w,\rho}),\Phi_2(\frak w,\rho)).
\end{equation}
Here $\Phi_2(\frak w,\rho)$ is the $ (\overline{\text{\bf f}}_{a(1)} \times
{\overline{\text{\bf f}}}_{a(2)})$ factor of $\Phi(\frak w,\rho) \in \widehat V_{\infty}^+ \times  (\overline{\text{\bf f}}_{a(1)} \times
{\overline{\text{\bf f}}}_{a(2)})$. 
\par
We put $\Gamma_{\frak p} = \{1\}$.
(When we discuss case 
$(S^2_1,(p_1,p_2),w\vert_{S^2_1})$ has an automorphism, 
that is, $\Gamma_{\frak p} \ne \{1\}$,
we need to take more marked points than $p_3$ so that the 
construction is symmetric with respect to this automorphisms.
(See \cite[Definition 17.5]{foootech}.)
For simplicity we assumed that $(S^2_1,(p_1,p_2),w\vert_{S^2_1})$ has 
no automorphism.)
\par
Now it is easy to see that
$(V_{\frak p},E_{\frak p},\Gamma_{\frak p},s_{\frak p},\psi_{\frak p})$ 
is a Kuranishi neighborhood of
$((\Sigma_{\frak p},\vec z_{\frak p},w_{\frak p}),(\frak y_1,\frak y_2))$.
\par
The evaluation map ${\rm ev}_0$ extends to $V_{\frak p}$.
Again by (\ref{decayequationOKw233}) we have:
\begin{equation}\label{decayequationOKw23324422}
\left\Vert
\frac{\partial  {\rm ev}_0}{\partial T_i}
\right\Vert_{C^m} \le C e^{-c T_i}.
\end{equation}
\par
We next construct a family of multisections perturbing $s_{\frak p}$.
(Actually we construct a family of single valued sections.
This is because in our case $\Gamma_{\frak p}$ is trivial.)
We first study the part where $T_1=T_2=\infty$.
We defined the manifold $\widehat V_i$ ($i=0,1,2$) and 
the vector bundle $E_i$ on it. The section $s_i$ of $E_i$ is defined as follows.
An element of $\widehat V_i$  is represented by 
either a map $w'_0 : D^2 \to X$ or $w'_i : S^2 \to X$ ($i=1,2$).
By assumption
$\overline\partial w'_i \in E(w'_i)$. 
We put
$$
s_i(w'_i) = \overline\partial w'_i.
$$
This is the Kuranishi map.
They are not necessarily transversal to $0$.
We take a finite dimensional manifold $W_i$. 
(We may take it so that it is an open neighborhood of 
the origin in an Euclidean space of sufficiently large dimension.) 
Then for any positive number $\delta$ there exist a section  
$\frak s_i$ of the bundle $E_i \times W_i \to \widehat V_i \times W_i$
so that it satisfies the following conditions.
\begin{conds}
\begin{enumerate}
\item
For each $w \in W_i$ we regard $\frak s_i(\cdot,w)$ as a section 
of $E_i$ on $\widehat V_i$. Then
$$
\vert \frak s_i(\cdot,w) - s_i\vert < \delta.
$$
Here the norm is $C^0$ norm.
\item
The section $\frak s_i$ is transversal to $0$.
(In particular, the zero set
$$
\frak s_i^{-1}(0) = \{ (\frak x,w) \in  \widehat V_i \times W_i \mid \frak s_i(\frak x,w) = 0\}
$$
is a smooth submanifold.)
\item
We consider the composition of the evaluation map and the projection to the first factor 
$: \widehat V_i \times W_i \to \widehat V_i$.
Then the following maps, which are their restrictions to $\frak s_i^{-1}(0)$, are submersions:
$$\aligned
({\rm ev}_0,{\rm ev}^p_1) &: \frak s_0^{-1}(0) \to L(u) \times X, \\
({\rm ev}^p_1,{\rm ev}^p_2,{\rm ev}^p_3) &: \frak s_1^{-1}(0) \to  X^3, \\
({\rm ev}^p_2,{\rm ev}_1,{\rm ev}_2) &: \frak s_2^{-1}(0) \to  X^3.
\endaligned$$
\end{enumerate}
\end{conds}
We note that we can always take a perturbation $\frak s_i$ satisfying Items 1, 2 above
without introducing $W_i$. (In other words we may take $W_i =$ point.)
However, in general, it is impossible to find $\frak s_i$ such that Item 3 is satisfied as well, 
if $W_i$ is one point.
\par
We put $W_{\frak p} = W_1 \times W_2 \times W_3$.
Note
$
\widehat V_{\infty} 
=
\widehat V_0 \,\,{}_{{\rm ev}_1^p} \times_{{\rm ev}_1^p}
\widehat V_1  {}_{{\rm ev}_2^p} \times_{{\rm ev}_2^p}
\widehat V_2.
$
We pull back $E_i$ to $\widehat V_{\infty}$ and take Whitney sum 
for $i=0,1,2$. We thus obtain a vector bundle $E_{\infty}$ on $\widehat V_{\infty}$.
The sections $\frak s_0,\frak s_1,\frak s_2$ induce a section of 
$E_{\infty}\times W \to \widehat V_{\infty} \times W$ in an obvious way.
We denote it by $\frak s_{\infty}$.
Then for any $w \in W$, we have
\begin{equation}\label{deltaclosefamimulti}
\vert \frak s_{\infty}(\cdot,w) - s_{\infty}\vert < \delta,
\end{equation}
where $s_{\infty}$ is the Kuranish map.
The section $\frak s_{\infty}$ is transversal to $0$. Moreover the map
\begin{equation}\label{famimultizerosubmer}
({\rm ev}_0,{\rm ev}_1,{\rm ev}_2) : \frak s_{\infty}^{-1}(0) \to L(u) \times X^2
\end{equation}
is a submersion.
Note
$
V_{\infty}  = \widehat V_{\infty} \,\,{}_{(ev_1,ev_2)} \times_{X^2}  (\overline{\text{\bf f}}_{a(1)} \times
{\overline{\text{\bf f}}}_{a(2)}).
$
We pull back $E_{\infty}$ to it and denote it by the same symbol.
The section $\frak s_{\infty}$ induces a section of this bundle on $V_{\infty} \times W$
which we denote by the same symbol. 
\par
The inequality (\ref{deltaclosefamimulti}) holds for this section 
$\frak s_{\infty}$ on $V_{\infty} \times W$.
Using the submersivity to the second and the third factors 
of (\ref{famimultizerosubmer}), the section $\frak s_{\infty}$ on $V_{\infty} \times W$
is transversal to $0$. Moreover the evaluation map
$$
{\rm ev}_0 : \frak s_{\infty}^{-1}(0) \cap (V_{\infty} \times W) \to L(u)
$$
is a submersion.
\par
We have thus described our family of multisections on the part $T_1=T_2=\infty$.
We next extend it.
We put $\rho_{\infty} = (\infty,\infty)$.
By construction, the fiber of the obstruction bundle $E(w_{\frak w,\rho})$ 
is induced by the parallel transport from $E(w_{\frak w,\rho_{\infty}})$.
Therefore the obstruction bundle on 
$V_{\frak p} = V_{\infty} \times ((T_0,\infty] \times S^1)/\sim)^2$ 
is canonically isomorphic to 
$E_{\infty} \times ((T_0,\infty] \times S^1)/\sim)^2$.
We use this isomorphism to extend $\frak s_{\infty}$ 
to $V_{\frak p} \times W$ so that it is constant in $\rho$ direction.
This is our family of multisections $\frak s_{\frak p}$.
\par
By taking $T_0$ sufficiently large,  (\ref{deltaclosefamimulti}) implies that
$$
\vert \frak s_{\frak p}(\cdot,w) - s_{\frak p} \vert < 2\delta.
$$ 
Moreover $\frak s_{\frak p}$ is transversal to zero.
Furthermore we have a canonical diffeomorphism
\begin{equation}\label{dirctprodzerosets}
\frak s_{\frak p}^{-1}(0) 
\cong
(\frak s_{\infty}^{-1}(0) \cap (V_{\infty} \times W)) \times ((T_0,\infty] \times S^1)/\sim)^2.
\end{equation}
\begin{sublem}\label{claim2633}
On $\frak s_{\frak p}^{-1}(0)$ we have the following inequality:
\begin{equation}\label{eval0expdecay}
\left\vert
\frac{\partial  {\rm ev}_0}{\partial T_i}
\right\vert \le C e^{-c T_i}.
\end{equation}
\end{sublem}
Note the partial derivative $\frac{\partial }{\partial T_i}$ in (\ref{eval0expdecay}) is defined by using the diffeomorphism 
(\ref{dirctprodzerosets}).
\begin{proof}
This is a consequence of (\ref{decayequationOKw233}).
\end{proof}
We have now described our family of multisections.
We also have the next two sublemmata.
\begin{sublem}\label{claim2634}
There exists a sequence of positive numbers $\frak T_i$ converging to $\infty$ and a sequence $\theta(i) \in \R/\Z$ such that 
$\Sigma_{\rho} = \Sigma_{3,i}$ if and only if
$\rho = ((T_1,\theta_1),(T_2,\theta_2))$ with 
$T_1+T_2 = \frak T_i$,  $\theta_1+\theta_2 = \theta(i)$.
\end{sublem}
\begin{proof}
We recall that $\Sigma_{3,i} \in \mathcal M_{1;2}^{\text{\rm main}}$ is a sequence converging to $\Sigma_0$.
Sublemma \ref{claim2634} then follows from the definition of $\Sigma_{\rho}$.
\end{proof}
\begin{sublem}\label{claim2635}
There exists a sequence of positive numbers $\mathcal T_K$ converging to $\infty$  and a sequence of positive numbers $I_K$ such that the following
holds for $i>I_K$.
\par
If $(\Phi(\frak w,\rho),\rho)$ is an element of the complement of $V_{\frak p,K}$ then $\rho = ((T_1,\theta_1),(T_2,\theta_2))$ 
with 
$T_1 > \mathcal T_K$, $T_2 > \mathcal T_K$.
\end{sublem}
\begin{proof}
This is an immediate consequence of Condition \ref{cond2629}.
\end{proof}
Now using Sublemmata \ref{claim2633}, \ref{claim2634}, \ref{claim2635},  the left hand side of (\ref{eq2622}) 
is estimated by
$$
C \int_{\mathcal T_K}^{\frak T_i-\mathcal T_K} \exp (- c \min\{ T, 
\frak T_i  -T\}) dT
<
2C \int_{\mathcal T_K}^{\infty} e^{-cT} dT < C' e^{-c\mathcal T_K}
$$
for $i>I_K$.
The proof of Lemma \ref{claim26299} is complete.
\end{proof}
Proposition \ref{sublemma:Sigma3limitprop} is a consequence of 
Lemmata \ref{claim2627}, \ref{claim2628} and \ref{claim26299}.
\end{proof}

\section{Hochschild and quantum cohomologies}
\label{sec:QMHOch}\index{Hochschild cohomology}
In this paper we study the operators $\frak p$, $\frak q$ and the cyclically
symmetric version $\frak q^{\frak c}$ of $\frak q$, in the case of toric manifold and its
Lagrangian submanifold that is an orbit of the $T^n$ action.
We use it to prove our main result, Theorem \ref{mtintro}.
We use various properties of these operators for this purpose.
We have proved them only in our toric case.
\par
However, those operators are defined in greater generality and most of
the properties we use here indeed hold in greater generality.
There are various literatures which
conjecture/claim/prove/use them in various
level of generality and of rigor.
\par
The purpose of this section is to briefly summarize the properties of operators
$\frak p$, $\frak q$ etc. which are expected to hold, in order
to put the discussion of this paper in proper perspective.
We also mention various earlier results. We might overlook some other existing literature
because our knowledge are certainly incomplete.
\par
We note that the operators $\frak p$, $\frak q$ etc. are defined in
\cite[Chapter 3]{fooo06}  and \cite[Section 3.8]{fooobook}  for any pair $(X,L)$
where $X$ is a symplectic manifold and $L$ is its relatively spin closed Lagrangian submanifold.
We first review its relation to Hochschild cohomology here.
\par
Let $C$ be a filtered $\Lambda_{0,{\rm nov}}$ module. We assume it is free.
We consider the Hochschild complex:
\begin{equation}\label{Hochcomplex}
CH(C[1],C[1]) = \prod_{k=0}^{\infty} Hom_{\Lambda_{0,{\rm nov}}}
(B_k(C[1]),C[1]).
\end{equation}
More precisely, for $\varphi_0 \in Hom_{\Lambda_{0,{\rm nov}}}
(B_0(C[1]),C[1]) = C[1]$, we assume
$\varphi_0(1) \equiv 0 \mod \Lambda_{0,{\rm nov}}^+$, in addition.
Its boundary operator is defined by
\begin{equation}\label{Hochboundary}
\aligned
(\delta\varphi)(x_1,\ldots,x_k)
=
&\sum_{k_1+k_2=k+1}\sum_{j=1}^{k_1} (-1)^{*_1}
\frak m_{k_1}(x_1,\ldots,\varphi_{k_2}(x_j,\ldots),\ldots,x_{k})\\
&+\sum_{k_1+k_2=k+1}\sum_{j=1}^{k_1} (-1)^{*_2}
\varphi_{k_1}(x_1,\ldots,\frak m_{k_2}(x_j,\ldots),\ldots,x_{k}).
\endaligned
\end{equation}
Here $\varphi = (\varphi_k)_{k=0}^{\infty}$ and
$$
*_1 = \deg\varphi_{k_2}(\deg' x_1 + \ldots + \deg' x_{j-1}),
\quad
*_2= \deg' x_1 + \ldots + \deg' x_{j-1}.
$$
(See \cite[Subsection 4.4.5]{fooobook}.)
Let us consider the case when $C = H(L;\Lambda_{0,{\rm nov}})$ as in
\cite[Theorem A]{fooobook}. The operators $\frak q_{\ell;k}$ defines a map
$$
\widehat{\frak q} : E(H(X;\Lambda_{0,{\rm nov}})[2])
\to 
{CH}(H(L;\Lambda_{0,{\rm nov}})[1],H(L;\Lambda_{0,{\rm nov}})[1])
$$
by
\begin{equation}\label{hatq}
\widehat{\frak q}(\text{\bf p}_1,\ldots,\text{\bf p}_{\ell})(x_1,\ldots,x_k)
=
\frak q_{\ell;k}(\text{\bf p}_1,\ldots,\text{\bf p}_{\ell};x_1,\ldots,x_k).
\end{equation}
It is easy to see that the image of this map consists of cycles.
See \cite[Section 7.4]{fooobook2}.
\par
We next define product structure on $CH(C[1],C[1])$.
Let $\varphi = (\varphi_k)_{k=0}^{\infty}$,  $\psi = (\psi_k)_{k=0}^{\infty}$.
We put
\begin{equation}
\aligned
\frak m_2(\varphi,\psi)(x_1,\ldots,x_k)& \\
= \sum_{k_1+k_2+\ell_0 + \ell_1 + \ell_2 =k}(-1)^{*}&\frak m_{\ell_0 + \ell_1 + \ell_2 + 2}(x_1, \dots, x_{\ell_0}, \varphi_{k_1}(x_{\ell_0 +1},\ldots,x_{\ell_0 + k_1}), \\
& x_{\ell_0 + k_1 + 1}, \dots, x_{\ell_0+ k_1 + \ell_1},
\psi_{k_2}(x_{\ell_0 + k_1+ \ell_1 +1},\\
&\ldots, x_{\ell_0 + k_1 + \ell_1 + k_2}), x_{\ell_0 + k_1 + \ell_1 + k_2 + 1}, \dots, x_k)
\endaligned
\end{equation}
where
$$
* = \deg \psi_{k_2} (\deg' x_1 + \cdots + \deg'x_{k_1}).
$$
Then
$$
\varphi \cup \psi = (-1)^{\deg\varphi(\deg \psi+1)} \frak m_2(\varphi,\psi)
$$
defines an associative product.
\begin{rem}
\begin{enumerate}
\item
We can define
$
\frak m_k(\varphi_1,\ldots,\varphi_k)
$
in a similar way. It defines a structure of $A_{\infty}$ algebra on $CH(C[1],C[1])$.
\item
This $A_{\infty}$ structure is a particular case of one defined in \cite[Theorem-Definition 7.55]{fooo099}.
Namely we consider the filtered $A_{\infty}$ category with one object
whose endomorphism algebra is the  $A_{\infty}$ algebra $C$. Then
$CH(C[1],C[1])$ is identified with the set of all
prenatural transformations from identity functor to identity functor.
In \cite[Theorem-Definition 7.55]{fooo099}, the filtered $A_{\infty}$ structure is constructed on it. If we apply it to our particular case,
it coincides with the above $A_{\infty}$ structure on $CH(C[1],C[1])$.
\end{enumerate}
\end{rem}
\begin{clm}\label{ringhomogeneral}
The map $\widehat{\frak q}$ defines a ring homomorphism
\begin{equation}\label{hatqmap}
QH(X;\Lambda_{0,{\rm nov}}) \to HH(H(L;\Lambda_{0,{\rm nov}})[1],H(L;\Lambda_{0,{\rm nov}})[1]).
\end{equation}
Here we put quantum cohomology ring structure on the left hand side.
The right hand side is $\delta$-cohomology of
$CH(H(L;\Lambda_{0,{\rm nov}})[1],H(L;\Lambda_{0,{\rm nov}})[1])$, that is the
Hochschild cohomology.
\end{clm}
\begin{rem}\label{rem293}
\begin{enumerate}
\item
When we consider several Lagrangian submanifolds instead of
one Lagrangian submanifold, it forms a filtered $A_{\infty}$ category
called Fukaya category. Then we can
modify the statement so that the right hand side is a
Hochschild cohomology of Fukaya category.
\item
We can include bulk deformation by an element
$\frak b \in H^{ev}(L;\Lambda_{0,{\rm nov}}^+)$ then Claim
\ref{ringhomogeneral} still holds. Namely
we put $QH(X;\Lambda_{0,{\rm nov}})$ the $\frak b$-deformed quantum
cup product. (See (\ref{defcup2}).) We use also  the $\frak b$-deformed
$A_{\infty}$ structure on Floer cohomology of the right hand side of
(\ref{ringhomogeneral}).
Then (\ref{ringhomogeneral}) becomes a ring homomorphism in this generality. See \cite{AFOOO1}.
\end{enumerate}
\end{rem}
\begin{conj}\label{comg294}
We use quantum (higher) Massey product on $QH(X;\Lambda_{0,{\rm nov}})$ to obtain
a filtered $A_{\infty}$ structure on it. (See \cite[Corollary 1.10]{fooo092}.)
Then $(\ref{hatqmap})$ is extended to a filtered $A_{\infty}$ homomorphism.
\end{conj}
Remark \ref{rem293}.1 and 2 also apply to Conjecture \ref{comg294}.
\par
Theorem \ref{multiplicative} is closely related to Claim \ref{ringhomogeneral},
but it looks slightly different. We will explain its relation now.
(For the simplicity of notation we consider the case when $\frak b = 0$, namely the case
without bulk deformation. Including bulk deformation is fairly straightforward.)
\par
We first consider the case when the potential function $\frak{PO}$ is a Morse function.
In that case there are finitely many $(L(u_i),b_i)$ $(i=1,\ldots,\text{rank} H(X))$, such that
$$
HF((L(u_i),b_i),(L(u_i),b_i)) \cong H(T^n).
$$
(See Lemma \ref{HHvanish}.)
Moreover $HF((L(u_i),b_i),(L(u_i),b_i))$ is a (nondegenerate) Clifford algebra as a ring.
So we have
$$
HH(H(L(u_i)),H(L(u_i));\Lambda) \cong \Lambda.
$$
(Here we use the $A_{\infty}$ structure $\frak m_k^{b_i}$ to define the Hochschild cohomology in the left hand side.)
Thus Claim \ref{ringhomogeneral}  implies that there exists a ring homomorphism
\begin{equation}\label{ringhomoifactor}
QH(X;\Lambda) \to \prod_{i}  \Lambda.
\end{equation}
Using the identification (\ref{factorizeJac}), it is easy to identify
$\prod_{i}  \Lambda$ with $\text{Jac}(\frak{PO}) \otimes_{\Lambda_0}\Lambda$.
Then the ring homomorphism (\ref{hatqmap}) induces a ring homomorphism to
$\text{Jac}(\frak{PO}) \otimes_{\Lambda_0}\Lambda$.
It coincides with the Kodaira-Spencer map appearing in Theorem \ref{multiplicative}.
\par
We next consider the case when $\frak{PO}$ is not
necessary a Morse function.
We consider the basis $\text{\bf e}$, $\text{\bf e}_i$
($i=1,\ldots,n$), $\text{\bf e}_I$ ($I \subseteq \{1,\ldots,n\}$) of
$H^*(T^n;\Q)$. Here $\text{deg}\text{\bf e} = 0$,
$\text{deg}\text{\bf e}_i = 1$, $\text{deg}\text{\bf e}_I \ge 2$.
\par
Let $\varphi = (\varphi_k)_{k=0}^{\infty}$ be an element of
$CH(H(L;\Lambda_0)[1],H(L;\Lambda_0)[1]) \otimes_{\Lambda_0} \Lambda$.
We put
$$
f_{j_1,\ldots,j_k} = \text{the coefficient of } {\bf e} \text{ in } \varphi_k(\text{\bf e}_{j_1},\ldots,\text{\bf e}_{j_k})
\in \Lambda.
$$
We define
\begin{equation}\label{297map}
\frak k(\varphi)
= \sum_{k=0}^{\infty}\sum_{j_1,\ldots,j_k = 1}^n \frac{1}{k!} f_{j_1,\ldots,j_k}x_{j_1}\cdots x_{j_k}
\in \Lambda_0[[x_1,\ldots,x_n]] \otimes_{\Lambda_0} \Lambda.
\end{equation}
\begin{lem}\label{HHandJac}
$(\ref{297map})$ induces a $\Lambda$ linear map
$$
HH(H(L;\Lambda_0)[1],H(L;\Lambda_0)[1]) \otimes_{\Lambda_0} \Lambda
\to
\prod_{j} \frac
{\Lambda_0[[x_1,\ldots,x_n]] \otimes_{\Lambda_0} \Lambda}
{\big(\frac{\partial\frak{PO}^{u_j}(x+b_j)}{\partial x_i};i=1,\ldots,n\big)}.
$$
Here we use $\frak m_k^b$ to define Hochschild boundary operator in the left hand side.
\end{lem}
We will write the boundary operator $\delta$ in (\ref{Hochboundary})
as $\delta_{\frak b}$ in case we  use $\frak m^{\frak b}_*$.
\begin{proof}
First of all, we note that Hochschild cohomology is isomorphic to the cohomology of
its normalized complex consisting of Hochschild cochains which vanish on elements
$(a_1, \dots, a_n)$ such that one of $a_i's$  is the unit ${\bf e}$.
So we consider the normalized complex.

It suffices to show that the image of the Hochschild boundary operator maps to the
Jacobian ideal in the right hand side by the map $(\ref{297map})$. Let
$\psi = (\psi_k)_{k=0}^{\infty}$ be a Hochschild cochain.
Then we find that $ k! {\frak k}(\delta_{\frak b}(\psi))$ is the coefficient of ${\bf e}$ in
$\delta_{\frak b}(b, \dots, b)$, where $b=\sum_{i=1}^n x_i {\bf e}_i$.
Since $b \in H^1(L(u_j);\Lambda_0)$, we find that ${\frak m}_{\frak b, k}(b, \dots, b)$
is a multiple of ${\bf e}$.
Hence the second term of the right hand side in \eqref{Hochboundary} with $\delta \phi = \delta_{\frak b} \psi$
vanishes, since $\psi$ is a normalized cochain.
The first term of the right hand side of \eqref{Hochboundary} with $\delta \phi = \delta_{\frak b} \psi$
is equal to ${\frak m}_1^{\frak b, b}(\psi(b, \dots, b))$.
By Proposition \ref{generatedede}, we find that it belong to the Jacobian ideal multiplied to
the module $H(L(u_j);\Lambda_0)$.
Hence the lemma.
\end{proof}
We can prove
$$
\text{\rm Jac}(\frak{PO}) \otimes_{\Lambda_0} \Lambda
\cong
\prod_j\frac
{
\Lambda_0[[x_1,\ldots,x_n]] \otimes_{\Lambda_0} \Lambda}
{\big(\frac{\partial\frak{PO}^{u_j}(x+b_j)}{\partial x_i};i=1,\ldots,n\big)}.
$$
We remark that the left hand side is a quotient of the
strictly convergent power series ring, while
the right hand side is a quotient of the
formal power series ring.
We can however prove that they are isomorphic in the present situation
in a way similar to Lemma \ref{polyapprox}.
See \cite{AFOOO2}.
\par
Therefore, by Lemma \ref{HHandJac}, we obtain a $\Lambda$-linear map
\begin{equation}\label{HHtoJacmap}
HH(H(L;\Lambda_0)[1],H(L;\Lambda_0)[1]) \otimes_{\Lambda_0} \Lambda
\to
\text{\rm Jac}(\frak{PO}) \otimes_{\Lambda_0} \Lambda.
\end{equation}
Composition of (\ref{HHtoJacmap}) and (\ref{hatqmap}) is the Kodaira-Spencer map
appearing in Theorem \ref{multiplicative}.
\par
We can prove Claim \ref{ringhomogeneral} in a way similar to Theorem \ref{multiplicative}.
We do not try to prove this here since we do not use it in this paper.
\par
The homomorphism $\widehat{q}$ in (\ref{hatqmap}) is related to the homomorphism
$i^{\ast}_{\text{\rm qm},(b,u)}$ given in \eqref{eq:iqmbu} as follows.
Let $\varphi = (\varphi_k)_{k=0}^{\infty}$ be a Hochschild cochain.
We put
$$
\overline{\varphi} = \varphi_0(1) \in C.
$$
It is easy to see that $\overline{\varphi}$ is an $\frak m_1$-cocycle if $\varphi$ is a
Hochschild cocycle. Thus we have
\begin{equation}\label{HHtoHF}
HH(H(L;\Lambda_0)[1],H(L;\Lambda_0)[1]) \otimes_{\Lambda_0} \Lambda
\to
HF((L(u),b),(L(u),b);\Lambda).
\end{equation}
$i^{\ast}_{\text{\rm qm},(b,u)}$ is a composition of (\ref{HHtoHF}) and
$\widehat{q}$. Since (\ref{HHtoHF})
and $\widehat{q}$ are both ring homomorphisms\index{$i^{\ast}_{\text{\rm qm},(\frak b,b,u)}$}\index{operator $\frak q$}
\begin{equation}\label{iqinduce}
i^{\ast}_{\text{\rm qm},(b,u)} : QH(X;\Lambda) \to HF((L(u),b),(L(u),b);\Lambda)
\end{equation}
is also a ring homomorphism.
\begin{rem}
\begin{enumerate}
\item
We note that $i^{\ast}_{\text{\rm qm},(b,u)}$  is constructed in
\cite{fooo06} and   \cite[Theorem 3.8.62]{fooobook}.
However the fact that it is a {\it ring} homomorphism is neither
stated nor proved there.
\item
For the case of monotone $L$, the fact that $i^{\ast}_{\text{\rm qm},(b,u)}$
is a ring homomorphism is closely related to a statement proved in
\cite[Theorem A]{bircor2}.
\par
Here are some more precise statements on their relationship.
\begin{enumerate}
\item
The algebra which is denoted by $QH_*(L;\mathcal R)$ in \cite{bircor2}
is isomorphic to the Floer homology. Namely
$$
QH_*(L;\mathcal R) \otimes_{\mathcal R} \Lambda_{0,{\rm nov}}
\cong
HF((L,0),(L,0);\Lambda_{0,{\rm nov}}),
$$
as rings.
(Here it is assumed that $L$ is monotone.)
This fact had been proved in the year 2000 version \cite[pages 287--291]{fooo00},
\emph{without the monotonicity assumption}, except the compatibility with ring structure.
Note we may take $\mathcal R =\Q[T,T^{-1}]$ as the coefficient ring in the monotone
case as explained in \cite[Chapter 2]{fooobook}.
\item
The map
\begin{equation}\label{BCmap}
QH(X;\mathcal R) \otimes QH(L;\mathcal R) \to QH(L;\mathcal R)
\end{equation}
appears in \cite[Theorem A (ii)]{bircor2}. Via the isomorphism
mentioned above, this map coincides with the map
\begin{equation}\label{iqinduce2}
(\text{\bf p},x) \mapsto \frak q_{1;1}(\text{\bf p},x).
\end{equation}
It is proved in \cite{bircor2} that (\ref{BCmap}) determines a module structure on
$QH(L;\mathcal R)$ over the quantum cohomology ring $QH(X;\mathcal R)$.
In view of the fact that $\frak m_0 =0$ in monotone case,
this statement is a special case of the Claim \ref{ringhomogeneral}.
(See also \cite{seidel:graded}.)
\par
We note that (\ref{iqinduce2}) is slightly different from
(\ref{iqinduce}). The latter map corresponds to
$
\text{\bf p} \mapsto \frak q_{1;0}(\text{\bf p}).
$
\end{enumerate}
\item
The homomorphism (\ref{hatqmap}) had been known much earlier.
It seems that the first reference where it is mentioned is
\cite{konts:hms}. (In \cite{konts:hms}, the right hand is
replaced by Hochschild homology of Fukaya category which is a more general
statement than the present context. There may also be several physics literature
containing statements related to this homomorphism.)
The category version of the map (\ref{hatqmap}) is believed to be an isomorphism under suitable assumptions.
This is stated explicitly as a conjecture by P. Seidel in  \cite[Conjecture 1]{seidel:HH}
for the exact case. (In noncompact case such as the exact case,
one needs to replace quantum cohomology by symplectic homology.)
There are also earlier references such as \cite{seidel:vm} but
\cite{seidel:HH} contains the statement in the most transparent way.
Seidel also constructed the map (\ref{hatqmap}) and
proved that it is a ring homomorphism in the exact case in  \cite[Section 6 (6B)]{seidel:vm}.
(`The technicalities' mentioned in the 4-th line of  \cite[Section 6 (6B)]{seidel:vm}
can be easily taken care of by now, that is 10 years after \cite{seidel:vm} was written.)
\par
Recently this ring homomorphism  (\ref{hatqmap}) is familiar to
many experts and becomes a kind of `folklore' in the field.
However, as far as the authors know at the stage of September 2010, rigorous proof of
the ring homomorphism property is written in detail only
for the cases where $L$ is either exact (Seidel), monotone (Biran-Cornea).
See \cite{AFOOO1} for general cases.
\end{enumerate}
We also remark that this kinds of application of  `open-closed Gromov-Witten invariant'
in Floer  homology
first appeared in the literature in \cite{floer:cup}.
\end{rem}
There is a cyclic analogue of the map  (\ref{hatqmap}).
Namely we can use the cyclically symmetric version $\frak q^{\frak c}$ of $\frak q$ to define a map
\begin{equation}\label{hatqcmap}
\widehat{\frak q}^{\frak c} : QH(X;\Lambda_{0,{\rm nov}}) \to HC(H(L;\Lambda_{0,{\rm nov}})[1];\Lambda_{{\rm nov}})
\end{equation}
by
$$
\widehat{\frak q}^{\frak c}(\text{\bf p}_1,\ldots,\text{\bf p}_{\ell})(x_1,\ldots,x_k,x_0)
=
\left\langle\frak q^{\frak c}_{\ell;k}(\text{\bf p}_1,\ldots,\text{\bf p}_{\ell};x_1,\ldots,x_k),
x_0\right\rangle.
$$
Here the right hand side of (\ref{hatqcmap}) 
is the cyclic cohomology.\index{cyclic cohomology}
\index{operator $\frak q^{\frak c}$}
This homomorphism is expected to play an important role
to pursue the line of ideas explained in Section
\ref{sec:cyclic cohomology}.
\par
The operator $\frak p$ is constructed in \cite[Section 3.8]{fooobook}
for any pair $(X,L)$. It induces a map
\begin{equation}\label{hatp}
\widehat{\frak p} : HC(H(L;\Lambda_{0,{\rm nov}})) \to H(X;\Lambda_{0,{\rm nov}})
\end{equation}
by
$$
[x_1 \otimes \ldots \otimes x_k]
\mapsto
{\frak p}_k(x_1, \ldots, x_k).
$$
Here the left hand side of (\ref{hatp}) is cyclic homology.\index{cyclic homology}\index{operator $\frak p$}
\begin{rem}
We can use $\frak p_{\ell;k}$ of \cite[Section 3.8.9]{fooobook}  to extend (\ref{hatp})
to a map
$$
EH(X;\Lambda_{0,{\rm nov}})\otimes HC(H(L;\Lambda_{0,{\rm nov}})) \to H(X;\Lambda_{0,{\rm nov}})
$$
in a similar way.
\end{rem}
If $b$ is a bounding cochain (namely if $\frak m_0^b =0$), then
we have
$$
HF((L,b);\Lambda_{0,{\rm nov}}) \to HC(H(L;\Lambda_{0,{\rm nov}}))
$$
by
$$
[x] \mapsto [x] \in B_1^{\text{cyc}}H(L;\Lambda_{0,{\rm nov}})[1].
$$
(We use the $A_{\infty}$ structure $\frak m^b$ to define the above cyclic homology.)
Combined with (\ref{hatp}), we obtain
\begin{equation}\label{LMinduce}
i_{L,b} : HF((L,b);\Lambda_{0,{\rm nov}}) \to H(X;\Lambda_{0,{\rm nov}}).
\end{equation}
It coincides with the map $i_{\ast,\text{\rm qm},(b,u)}$ in (\ref{qsharpmap}).
\par
Theorem \ref{annulusmain} is expected to hold in that generality.
Namely we expect
\begin{equation}\label{21mainformulagen}
\aligned
&\langle i_{L,b}(\frak v), i_{L,b}(\frak w)
\rangle_{\text{\rm PD}_{X}} \\
&=
\sum_{I,J} \pm g^{IJ}
\langle \frak m_2^{\frak c,b}(\text{\bf e}_I,\frak v),
\frak m_2^{\frak c,b}(\text{\bf e}_J,\frak w)
\rangle_{\text{\rm PD}_{L}}.
\endaligned\end{equation}
Here $\text{\bf e}_I$ is a basis of the cohomology $H(X;\Lambda_{\rm nov})$.
See Proposition \ref{propconclusion} for sign $\pm$.
Other notations are the same as in Section \ref{sec:annuli}.
\par
Many parts of the proof in Section \ref{sec:annuli} can be generalized to the non-toric
case. However we take short cut in several places in Section \ref{sec:annuli}.
So several parts of the proof should be rewritten in the general case.
\begin{rem}
\begin{enumerate}
\item
In the monotone case the map $i_{L,0}$ coincides with the
map $i_L$ appearing in \cite[Theorem A (iii)]{bircor2}.
The authors of \cite{bircor2} attributes those maps to the paper \cite{albers}.
Indeed, the map $\xi_k$ and $\tau_k$ which appear in
\cite[Theorem 1.5]{albers} (defined in the monotone case) coincide with
(\ref{iqinduce}) and (\ref{LMinduce}), respectively.
\par
One of the main applications of  \cite[Theorem 1.5]{albers} is to prove
nontriviality of Lagrangian Floer cohomology under certain conditions.
The fundamental chain of $L$ is used for this purpose there.
In the proof of \cite[Theorem 6.13 (6.14.3)]{fooo00}  a similar argument had been
used for the same purpose. The operators $\frak p$ was introduced in
\cite{fooo06} in order to make this argument more systematic.
\item
We used the moduli space of pseudo-holomorphic annuli to prove
(\ref{21mainformulagen}) in the toric case.
The moduli space of pseudo-holomorphic annuli is also used in
\cite{bircor} and in \cite{abouz}.
\end{enumerate}
\end{rem}
\par

\bibliographystyle{amsalpha}

\printindex
\end{document}